\pdfoutput=1
\documentclass[11pt]{amsart}
\usepackage{amsmath,amssymb,amsthm}
\usepackage{mathtools}
\usepackage{mathrsfs}
\usepackage[margin=1.2in]{geometry}
\usepackage{booktabs}
\usepackage{array}
\usepackage{enumitem}
\usepackage{xcolor}
\usepackage{tikz}
\usepackage[colorlinks=true,linkcolor=blue,citecolor=blue,urlcolor=blue]{hyperref}
\hypersetup{pdfauthor={Bernd Johannes Wuebben},
  pdftitle={Cobordism maps and Atiyah–Floer isomorphisms for framed instanton homology I},
  pdfsubject={Naturality under cobordisms of mixed isomorphisms between framed instanton homology and symplectic framed Floer homology},
  pdfkeywords={Atiyah-Floer conjecture, framed instanton homology, Lagrangian Floer homology,
  mixed equation, cobordism maps, adiabatic limits}}

\newtheorem{theorem}{Theorem}[section]
\newtheorem{proposition}[theorem]{Proposition}
\newtheorem{lemma}[theorem]{Lemma}
\newtheorem{corollary}[theorem]{Corollary}
\theoremstyle{definition}
\newtheorem{definition}[theorem]{Definition}
\newtheorem{example}[theorem]{Example}
\newtheorem{convention}[theorem]{Convention}
\newtheorem{hypothesis}[theorem]{Hypothesis}
\theoremstyle{remark}
\newtheorem{remark}[theorem]{Remark}

\theoremstyle{plain}
\newtheorem{thmletter}{Theorem}

\numberwithin{equation}{section}

\newcommand{\len}{\operatorname{len}}
\newcommand{\ZZ}{\mathbb{Z}}
\newcommand{\RR}{\mathbb{R}}
\newcommand{\CC}{\mathbb{C}}
\newcommand{\FF}{\mathbb{F}_2}
\newcommand{\Hb}{\mathbb{H}}
\newcommand{\cB}{\mathcal{B}}

\newcommand{\cF}{\mathcal{F}}
\newcommand{\cG}{\mathcal{G}}

\newcommand{\cM}{\mathcal{M}}
\newcommand{\cP}{\mathcal{P}}

\newcommand{\bB}{\mathbf{B}}
\newcommand{\ad}{\operatorname{ad}}
\newcommand{\ind}{\operatorname{ind}}
\newcommand{\coker}{\operatorname{coker}}
\newcommand{\im}{\operatorname{im}}
\newcommand{\tr}{\operatorname{tr}}

\newcommand{\SO}{\mathrm{SO}}
\newcommand{\SU}{\mathrm{SU}}
\newcommand{\U}{\mathrm{U}}
\newcommand{\Ish}{I^{\#}}
\newcommand{\HF}{\mathit{HF}}
\newcommand{\CF}{\mathit{CF}}
\newcommand{\Msf}{\mathsf{M}}
\newcommand{\eps}{\varepsilon}

\newcommand{\dist}{\operatorname{dist}}

\title[Atiyah--Floer cobordism maps I]{Cobordism maps
and Atiyah--Floer isomorphisms\\ for framed instanton homology I}
\author{Bernd Johannes Wuebben}

\date{First version: 25 September 2026.\newline This version: 3 November 2026}

\begin{document}

\begin{abstract}
We prove that isomorphisms of the type constructed by Daemi, Fukaya and Lipyanskiy, between
framed instanton homology and the Lagrangian Floer homology of a Heegaard splitting,
intertwine the cobordism maps of Kronheimer and Mrowka with maps defined by holomorphic
curves alone, over the field of two elements. This holds for every connected cobordism
between connected three-manifolds carrying a framed arc between its ends and an $\SO(3)$
bundle trivialized along the ends and the arc. The symplectic maps count holomorphic triangles in
moduli spaces of flat connections on surfaces, with the correspondence of a compression body
as boundary condition for one- and three-handles. Consequently these maps are independent of the handle decomposition and of the
Heegaard data, and form a functor isomorphic to framed instanton homology; where both are
defined, our isomorphisms agree with those of Daemi, Fukaya and Lipyanskiy modulo two. The
proof rests on two analytic results for the mixed anti-self-duality and holomorphic curve
equation: an adiabatic limit in which a compression body collapses next to a holomorphic strip
and becomes the Lagrangian boundary condition of its flat connections, and an excision theorem
along a torus.
\end{abstract}

\maketitle

\setcounter{tocdepth}{1}
\newpage
\tableofcontents
\newpage

\section{Introduction}\label{sec:intro}

This is the first of three papers on the compatibility of Atiyah--Floer
isomorphisms with cobordism maps. The present paper proves the comparison
over $\FF$ and establishes the mixed-equation analysis used throughout
the series. Paper~II~\cite{Wuebben-II} develops the integral theory:
orientations are transported from the almost-complex gauge convention,
and the comparison holds up to chain homotopy and one overall sign for
each cobordism. Paper~III~\cite{Wuebben-III} identifies this integral
comparison with the usual singular Hopf-link instanton maps and, at
common regular data, with the Daemi--Fukaya--Lipyanskiy mixed count in
a specified coherent orientation system. It also constructs, for every
closed connected oriented three-manifold, a comparison with the
normalization in which attachment of a generating sphere preserves
orientation. That construction uses finite integer combinations of
products of geometric counts and retains the stated surgery and
oriented gluing inputs. Thus the series proceeds from the comparison
modulo two, through integral naturality, to identification and
normalization of the integral maps.

\subsection{The Atiyah--Floer isomorphism}

Let $M$ be a closed, connected, oriented three-manifold. Its framed
instanton homology $\Ish(M)$ is the instanton Floer homology of the
connected sum $M\#T^3$ with the $\SO(3)$ bundle that is trivial over $M$
and nontrivial over $T^3$, with the gauge group enlarged by the automorphism of the
torus summand \cite[Section~7]{Scaduto}; Kronheimer and Mrowka define it equivalently
through a Hopf link \cite{KM-unknot}. A Heegaard
splitting $M=H_0\cup_\Sigma H_1$ of genus $g$ gives a second invariant.
Take the boundary sum of each handlebody with $[0,1]\times T^2$. The
two resulting compression bodies $C_0,C_1$ have boundary
$\Sigma_h\sqcup T^2$ with $h=g+1$. Let $\Msf_h$ denote the moduli space
of flat connections on the nontrivial $\SO(3)$ bundle over $\Sigma_h$.
The flat connections that extend over $C_i$ form an embedded Lagrangian
$L(C_i)\subset\Msf_h$. The Lagrangian Floer homology
$\HF(L(C_0),L(C_1))$ is the symplectic framed Floer homology of \cite{DFL},
introduced by Wehrheim and Woodward as torus-summed Floer homology
\cite[Definition~4.4.1]{WW-fft}; they show by Cerf theory and geometric composition
that it is an invariant of $M$. Atiyah conjectured that, for an integral homology
sphere, instanton Floer homology \cite{Floer} is recovered in this way from the
$\SU(2)$ character varieties of a Heegaard splitting \cite{Atiyah}; for the version
for admissible $\SO(3)$ bundles adopted here see \cite[Section~10]{DF}. Daemi, Fukaya and Lipyanskiy proved that the two groups are
isomorphic \cite[Theorem~3]{DFL}. Their isomorphism counts solutions of a
\emph{mixed equation}, introduced by Lipyanskiy \cite{Lipyanskiy}: the
anti-self-duality equation on part of a four-dimensional domain, coupled along
a matching surface to the holomorphic curve equation in $\Msf_h$ on the rest
\cite{DFL, DFL-mixed}.

Both groups carry cobordism maps. On the gauge side a cobordism $W$ with
a path between its ends induces a map
$\Ish(W)\colon\Ish(M_0)\to\Ish(M_1)$, and these maps form a functor
\cite[Section~5]{KM-unknot}, \cite[Section~7]{Scaduto}. On the
symplectic side a handle decomposition of $W$ can be assigned a composite of
counts of holomorphic polygons and quilts, using the quilted Floer theory of
Wehrheim and Woodward \cite{WW-quilt, WW-composition}. For the symplectic instanton
homology of Manolescu and Woodward \cite{MW}, such handle maps were defined and
shown to be independent of the decomposition by Cazassus \cite{Caz-nat}, using
Cerf theory, and Horton proved the corresponding functoriality for his
theory built from traceless character varieties \cite{Hor1},
\cite[Theorems~1.1--1.2]{Hor2}. For the
symplectic framed theory used here no maps for four-dimensional cobordisms were available, and in
none of these settings had the symplectic maps been compared with the gauge
cobordism maps. The comparison closest in spirit is the open--closed
isomorphism of Salamon and Wehrheim \cite[Theorem~5.8]{SW-openclosed}, which intertwines
relative invariants of four-manifolds, with gauge theory on both sides. This
paper compares the two families of maps.

\subsection{Main result}

We work over $\FF$ throughout. Objects of the connected framed
cobordism category $\bB$ are closed, connected, oriented three-manifolds
$M$ with a parametrized ball, which fixes where the torus is attached
(Definition~\ref{def:category}). A morphism $(W,E,\nu)\colon M_0\to
M_1$ is an isomorphism class of connected oriented cobordisms $W$ with a
framed path tube $\nu$ joining the two balls, and an $\SO(3)$ bundle $E$
on $W$ trivialized over $M_0$, $M_1$ and along the tube. The bundle is
arbitrary in the interior; its class relative to the ends lies in
$H^2(W,\partial W;\FF)$. Composition is gluing.

Framed instanton homology is a functor
$\Ish\colon\bB\to\mathrm{Vect}_{\FF}$
(Proposition~\ref{prop:gauge-functor}). For each object $M$ and each
admissible presentation $p$ (a Heegaard splitting adapted to the ball,
together with the analytic data of Section~\ref{sec:models}) there are
the gauge complex $C_G(M,p)$, the symplectic complex
$C_S(M,p)=\CF(L(C_0),L(C_1))$ and a mixed map
$N_p\colon C_G(M,p)\to C_S(M,p)$. It is a quasi-isomorphism by the
triangularity argument of \cite[Theorem~3.29]{DFL}
(Theorem~\ref{thm:mixed-iso}). The Lagrangians $L(C_i)$ are moved by \emph{primary}
perturbations, holonomy perturbations supported inside $C_i$; the mixed equation is
regularized by \emph{secondary} perturbations in the interior of its domain, and the gauge
complex by an \emph{end} perturbation. All of them are taken invariant under the
\emph{lifting group} $\phi_M\cong\ZZ/2$, the extra automorphism of the torus summand
\cite[Section~5.1]{KM-unknot}.

For a morphism we define symplectic maps by counting holomorphic curves
in the Heegaard moduli spaces, without reference to gauge theory:
\begin{itemize}[leftmargin=2em]
\item the trace of surgery on a framed link is sent to the triangle map
  $\mu_S(\,\cdot\otimes\theta_S)$ of a subordinate Heegaard triple
  $(\Sigma;H_0,H_1,H_2)$, where $\theta_S$ is the generator of maximal relative grading (the maximum of a height function) of the Floer complex of the Lagrangians of $H_1$ and $H_2$ (Section~\ref{sec:two-handles});
\item a $k$-fold one-handle is sent to a compression cup, defined by the
  triangle map of $(A_0,L_C,A_1)$ in $\Msf_h^-\times\Msf_{h+k}$.
  Here $L_C$ is the correspondence of the compression body, and $A_i$ is the
  product of a reversed handlebody Lagrangian of $M$ and a handlebody
  Lagrangian of $M^+$. The inputs are the generators of maximal grading
  for $(A_0,L_C)$ and $(L_C,A_1)$, the maxima of explicit Morse functions
  (Section~\ref{subsec:folded}). Unfolding gives the quilted triangle with
  seam $L_C$ (Remark~\ref{rem:quilted-cup},
  Figure~\ref{fig:folded-triangle}). A three-handle is sent to the transposed
  compression cap (Section~\ref{sec:one-three});
\item a change of trivialization $u\colon M\to\SO(3)$ of the bundle (a
  \emph{frame change}) is sent to the pushforward by the induced
  symplectomorphism of $\Msf_h$;
\item a change of presentation is sent to the corresponding word of
  continuations, cups, caps and triangles (Section~\ref{sec:presentations}).
\end{itemize}
An \emph{admissible decomposition} $d$ of $(W,E,\nu)$ is a factorization into
presentation changes, frame changes, single one- and three-handles and one
two-handle cobordism, from $(M_0,p_0)$ to $(M_1,p_1)$
(Section~\ref{subsec:generation}); it gives a composite map $S_d(W,E,\nu)$ from
$\HF(M_0,p_0)=H(C_S(M_0,p_0))$ to $\HF(M_1,p_1)$. A presentation is \emph{admissible} if its primary holonomy perturbations lie in a fixed
small open neighbourhood $\mathcal U(F)$ of zero in the span of a finite family $F$ of
holonomy perturbations, depending only on the Heegaard splitting, the metric data and $F$
(Definitions~\ref{def:presentation} and~\ref{def:admissible-nbhd}). Within $\mathcal U(F)$ we
prove that the mixed maps are natural under changes of primary data and that the transpose
of the reflected mixed map inverts $N_p$ (Propositions~\ref{prop:small-primary}
and~\ref{prop:duality}); this is used in the one-handle comparison.

\begin{thmletter}\label{thm:A}
For every morphism $(W,E,\nu)\colon M_0\to M_1$ of $\bB$, all admissible
presentations $p_0$ of $M_0$ and $p_1$ of $M_1$, and every admissible
decomposition $d$,
\[
H(N_{p_1})\circ\Ish(W,E,\nu)=S_d(W,E,\nu)\circ H(N_{p_0}).
\]
Consequently $S_d(W,E,\nu)$ depends only on the morphism and the end
presentations. The groups $\HF(M,p)$ form a transitive system over the
presentations of $M$, and with its colimit $\HF^{\#}(M)$ the maps
$S(W,E,\nu)$ form a functor $\HF^{\#}\colon\bB\to\mathrm{Vect}_{\FF}$.
The maps $H(N_p)$ descend to a natural isomorphism
$\Ish\Rightarrow\HF^{\#}$. The same statements hold in the homotopy
category of finite chain complexes over $\FF$, with the colimit replaced by
a choice of representative of the transitive system of homotopy
equivalences.
\end{thmletter}

The symplectic side of Theorem~\ref{thm:A} contains no gauge theory. In
particular, independence of the decomposition and functoriality of the
symplectic maps are consequences, not hypotheses. The argument uses no
Cerf theory on the symplectic side, in contrast with the direct constructions
of \cite{Caz-nat, Hor2} for other symplectic models, and no geometric composition
isomorphism: the embedded compositions of the compression correspondence
with handlebody Lagrangians enter only as sets, to identify generators and
gradings. The quilted description of the compression cup, which rests on
the composition theorems of Wehrheim--Woodward and Lekili--Lipyanskiy, is
recorded in Remark~\ref{rem:quilted-cup} and is not used.

Our mixed maps $N_p$ are built by the method of \cite{DFL}, with one change:
Chern--Simons trajectories are regularized by an additional holonomy perturbation $f$ of
$M^\#$ rather than by the primary perturbations (Remark~\ref{rem:dfl-map}). When the primary
perturbations already regularize them, as in \cite[Lemma~2.56]{DFL}, one may take $f=0$;
with the almost complex structures of \cite{DFL} and secondary perturbations invariant
under the lifting group, $N_p$ is then the isomorphism of \cite[Theorem~3.29]{DFL} reduced
modulo two (Proposition~\ref{prop:dfl-comparison}). Such primary data exist at every Heegaard
splitting and all metric data, without the rescaling of the metric used in
\cite[Proposition~6.8]{DFL} (Theorem~\ref{thm:dfl-regular}, proved in
Appendix~\ref{app:dfl-regular}). Hence Theorem~\ref{thm:A} holds for the isomorphism of
Daemi, Fukaya and Lipyanskiy itself, reduced modulo two, which on homology is independent
of the choices (Corollary~\ref{cor:dfl-natural}), and the mixed map at any admissible
presentation is, on homology, that isomorphism composed with gauge and symplectic
continuation maps (Section~\ref{sec:assembly}). Two points at which we depart from
\cite{DFL} are recorded in the text: the inverse function statement
\cite[Lemma~4.1(iii)]{DFL} is used with an $L^q$-ball, $q>2$, since it fails with an
$L^2$-ball (Remark~\ref{rem:dfl-radius}); and for the inductive step of
\cite[Proposition~6.8]{DFL} we could not deduce the independence of a cokernel element from
the translation vector on one piece of the splitting, which Appendix~\ref{app:dfl-regular}
supplies for separated base data (Remark~\ref{rem:dfl-regular}).

\subsection{The collapse theorem}

Every elementary comparison in Theorem~\ref{thm:A} reduces to one
degeneration, stated here informally. Let $X$ be a domain for the mixed
equation containing a \emph{compression rectangle}: a product
$I\times C$ of an interval with a compression body, attached to a
holomorphic strip along the matching surface $I\times\Sigma$, where $\Sigma$ is the
boundary component of $C$ of larger genus (Figure~\ref{fig:collapse}). Scale
the metric of $C$ by $\eps^2$. After a conformal change each endpoint of
$I$ becomes a cylinder over $(C,g_C)$ of length $\log(1/\eps)$, and
$\eps$ enters only through explicit weights. As $\eps\to0$ the rectangle
disappears. In the limit the strip meets $I$ along the Lagrangian
boundary condition $L(C)$, and two mixed ends appear at the endpoints.

Let $X_0$ denote the limiting domain, in which the compression rectangle is
replaced by the boundary condition $L(C)$, and $X_\eps$ the domain at scale $\eps$.

\begin{thmletter}\label{thm:B}
Assume Hypothesis~\ref{hyp:collapse} of Section~\ref{sec:collapse} (near the rectangle the
compression body is unperturbed and the almost complex structure on the holomorphic part is
$*_\Sigma$; all
other perturbations are supported away from it), that all data on $X_\eps$ outside a fixed neighbourhood of the rectangle
are independent of $\eps$, and that the moduli spaces on $X_0$ of index at most
zero are regular. Then the following hold.
\begin{enumerate}[label=(\roman*)]
\item If on $X_0$ the Fredholm index and the topological energy $\kappa$, normalized as in \cite[Proposition~3.20]{DFL}, satisfy $\ind=c+8\kappa$ for all configurations with fixed asymptotic labels, then the same relation, with the same constant $c$, holds on $X_\eps$ for every $\eps>0$.
\item For every sequence of solutions of index $d$ on $X_{\eps_\nu}$, with
  $\eps_\nu\to0$, a subsequence converges to a solution on $X_0$ of index $d_0$
  and broken trajectories of indices $d_j$, and $d-d_0-\sum_jd_j=4n$ for an
  integer $n\ge0$.
\item There is $\eps_0>0$ such that for $0<\eps<\eps_0$ gluing the constant flat
  solution on the rectangle defines a bijection
  $\cG_\eps\colon\cM^0(X_0)\to\cM^0(X_\eps)$ between the index-zero moduli spaces,
  and the solutions at scale $\eps$ are regular.
\end{enumerate}
\end{thmletter}

The precise statement is Theorem~\ref{thm:collapse}. The limit is of the type
proposed in the programme of Salamon and Wehrheim \cite[Conjecture~3.4]{Salamon-ICM},
\cite[Section~4]{Wehrheim-lagrangian}, in which the metric on a handlebody
degenerates onto its spine and the handlebody is replaced by the Lagrangian boundary
condition of its flat connections. (The second degeneration of that programme,
shrinking $\Sigma$ in $[0,1]\times\Sigma$, is the adiabatic limit of Dostoglou and
Salamon \cite{DS}; Fukaya \cite{Fukaya-degenerate} coupled the anti-self-duality
equation to holomorphic curves through a degenerate metric.) Here the collapsing
body is a compression body whose boundary is matched, along part of it, to a
holomorphic strip that does not collapse, and two mixed ends appear. To our
knowledge the handlebody degeneration of \cite{Salamon-ICM, Wehrheim-lagrangian}
has not been carried out; Salamon and Wehrheim prove their open--closed
isomorphism by a different route \cite{SW-openclosed}. After a conformal change the
collapse becomes neck stretching along cylinders over $C$. In the compactness
theorem of Xu \cite[Theorem~1.4]{Xu}, for instantons on $\RR\times M$ with $M$ a
three-manifold with a cylindrical end over $\Sigma$, energy escaping along the end
produces holomorphic curves in the moduli space of $\Sigma$ with the Lagrangian
boundary condition of $M$. Theorem~\ref{thm:B} adds to compactness the
index--energy law at every scale and the bijection of index-zero counts. The
symplectic analogue is strip shrinking in quilted Floer theory
\cite{WW-composition, BW-squiggly}. Three features of
the proof are specific to this situation. First, the model linear
estimates (Section~\ref{sec:collapse-linear}) lose no power of
$\log(1/\eps)$; the patched inverse on the whole domain loses exactly one
factor $1+\log(1/\eps)$, which is sharp, and the surjectivity argument is arranged so that this factor is always
multiplied by a positive power of $\eps$. The reason for the first point is
that the slow harmonic directions of $C$
are tied, pointwise in time, to the boundary value of the strip and are
controlled by the strip's unweighted norm. This depends on the
orientation convention~\ref{conv:orientation}; with the other sign the
estimate is false. Second, compactness requires no classification of
bubbles. Energy loss is excluded by the affine law~(i), which holds at
every scale, together with nonnegativity of the limiting indices.
Third, surjectivity of the gluing map is proved without entering a logarithmically shrinking
uniqueness ball directly. A solution is first brought within a power of
$\eps$ of the glued solution, and only then compared with it.

\subsection{Excision along a torus}

The one-handle comparison needs a second new result. Let $M$ have a Heegaard
splitting and let $M^+=M\#^k(S^1\times S^2)$ have the splitting obtained by
adding $k$ trivial handles. Cutting $(-M)^\#$ and $(M^+)^\#$ along torus fibres and regluing gives the closed three-manifold $Z\cong(-M)\#M^+\#T^3$. The pair of compression bodies
$(D_0,D_1)$ splitting $Z$ defines a Lagrangian pair $(A_0,A_1)$ in $\Msf_h^-\times\Msf_{h+k}$,
the \emph{folded pair} (Section~\ref{subsec:folded}). Let $N_Z$ be its mixed map, $N_-$ the
mixed map of the reflected presentation of $-M$ (valued in the dual Floer complex), $N_+$
that of $M^+$, $\mathsf E$ Floer's excision map
\cite{BD, KM-sutures}, \cite[Theorem~5.6]{KM-unknot} from the instanton complex
of $Z$ to the tensor product of those of $(-M)^\#$ and $(M^+)^\#$, and $\mathcal K$ the
folding isomorphism of Lagrangian Floer complexes for a split almost complex
structure.

\begin{thmletter}\label{thm:C}
With $N_Z$, $N_\pm$, $\mathsf E$ and $\mathcal K$ as above,
$[\mathcal K][N_Z]=([N_-]\otimes[N_+])[\mathsf E]$ on homology over $\FF$.
\end{thmletter}

This is Theorem~\ref{thm:excision}. It is not an instance of neck stretching along a closed
hypersurface: the torus region runs into the mixed ends, so no closed three-manifold separates
the two sides. Near the torus the domain is a product of $T^2$ with a planar region
(Figure~\ref{fig:torus-domains}). We insert a flat region $T^2\times[-\lambda,\lambda]\times\RR$
and, on its boundary, impose only flatness of the restriction to each torus fibre; since the
moduli space of flat connections on the torus with $w_2\ne0$ (the \emph{odd torus}) is a
point, this is the whole matching condition there. Two
facts drive the proof. On the plane, the half-plane and straight strips, the linearized
operator at the odd flat torus has an inverse bounded by $\mu_T^{-1}$, independently of the
strip width (Proposition~\ref{prop:flat-fibre}); and this boundary condition admits
small-energy regularity (Theorem~\ref{thm:flat-boundary-regularity}). The proof compares the
two sides at a fixed width with Floer's excision, varies the width, and lets it tend to
infinity (Section~\ref{sec:torus}).

\subsection{Strategy}

Theorem~\ref{thm:A} is deduced in Section~\ref{sec:assembly} from four
inputs by a formal argument (Lemma~\ref{lem:descent}):
\begin{enumerate}[label=(\arabic*)]
\item $\Ish$ is a functor on $\bB$ (Section~\ref{sec:category});
\item $N_p$ is an isomorphism on homology, and its homology class is
  independent of the end and secondary data (Theorem~\ref{thm:mixed-iso}(ii)) and,
  within the admissible neighbourhoods, of the primary data
  (Proposition~\ref{prop:small-primary});
\item every morphism is a composite of elementary arrows: presentation
  moves, frame changes, one-handles, a single two-handle cobordism and
  three-handles (Section~\ref{sec:handles});
\item every elementary arrow $e$ satisfies the square
  $H(N)\Ish([e])=S_eH(N)$, where $N$ denotes the mixed map $N_p$ of the
  chosen presentation at each end
  (Sections~\ref{sec:two-handles}--\ref{sec:presentations}).
\end{enumerate}
The maps $H(N_p)$ are invertible. Hence (3) and (4) determine every
symplectic word as $H(N)\,\Ish(W)\,H(N)^{-1}$. Independence of the word and
functoriality are then inherited from the gauge side.

The squares for two-handles are proved as follows. The gauge two-handle
map is a triangle map with an auxiliary input
(Section~\ref{sec:two-handles}). A mixed triangle degenerates in two ways.
Stretching one family of necks produces the gauge triangle followed by
$N$. Collapsing a compression rectangle by Theorem~\ref{thm:B} and then
stretching ordinary necks produces $N\otimes N$ followed by the
symplectic triangle. A compact family of mixed triangles joins the two
(Section~\ref{sec:families}).

For a $k$-fold one-handle the symplectic cup is a triangle count in the
product $\Msf_h^-\times\Msf_{h+k}$ of two moduli spaces, with the correspondence
of the compression body as middle Lagrangian and top classes of two
auxiliary pairs as inputs (Section~\ref{sec:one-three}). Folding the
triangle along the diagonal and applying the same two degenerations
identifies the cup with the mixed map of a folded pair composed with a gauge
triangle, and Theorem~\ref{thm:C} splits off the extra three-torus summand.
The result is the \emph{cup formula}
\[
U_C=N_{C\circ p}\circ\Ish(W_k)\circ N^\dagger_p ,
\]
where $W_k$ is the one-handle cobordism and $N^\dagger_p$ is the transpose of the mixed
map $N_{\bar p}$ of the reflected presentation. It remains to show that
$N^\dagger_p=N_p^{-1}$ on homology, and this is done without any composition theorem for
Lagrangian correspondences. For $k=0$ the cup is a continuation map. At a common auxiliary
presentation $A$, the cup formula for the reflections of two presentations with primary
data $P$ and $Q$ reads (indexing mixed maps by primary data, and using
$N^\dagger_{\bar P}=N_P^\vee$)
\[
\mathrm{cont}_S(\bar P,A)=N_A\circ\mathrm{cont}_G(\bar P,A)\circ N_P^\vee,\qquad
\mathrm{cont}_S(\bar Q,A)=N_A\circ\mathrm{cont}_G(\bar Q,A)\circ N_Q^\vee .
\]
Eliminating the invertible $N_A$ gives the naturality of $N$ under small changes of
primary data (Proposition~\ref{prop:small-primary}). Applying the $k=0$ cup
formula once more, with $P$ itself as the lower presentation, and comparing with
that naturality gives $N_P\circ N^\dagger_P=1$ (Proposition~\ref{prop:duality}). The
three-handle squares follow by reflection.

\subsection{Scope}

The theorem is stated over $\FF$. Over $\ZZ$ the symplectic counts
require coherent orientations of the polygon and quilt moduli spaces,
together with a comparison of those orientations with the homology
orientations of the gauge side. Neither is addressed here. Integer
coefficients are treated in Paper~II~\cite{Wuebben-II}, which proves the
integral comparison with orientations transported from the almost-complex
gauge convention, up to one sign for each cobordism. Paper~III~\cite{Wuebben-III}
compares this convention with the usual singular Hopf-link maps under its
stated inputs. The objects
carry the trivial bundle on the original three-manifold. Disconnected
objects, the units and counits of a field theory, and closed pairings are
not treated. The mixed maps $N_p$ are constructed by the method of
\cite{DFL}, with the invariant perturbations of Section~\ref{sec:models}.
Their agreement with the isomorphism of \cite{DFL} modulo two, at the common
data described above, is Proposition~\ref{prop:dfl-comparison}; existence of
such data at fixed metric data is Theorem~\ref{thm:dfl-regular}. The maps of the transitive system
of Theorem~\ref{thm:A} are symplectic words; their independence of the word is proved through
the gauge side, and we do not compare them with the Cerf-theoretic identifications
of \cite{WW-fft}.

\subsection*{Organization}

Section~\ref{sec:category} defines $\bB$ and the gauge functor.
Section~\ref{sec:models} recalls the Heegaard moduli spaces and the mixed
isomorphism and fixes conventions. Section~\ref{sec:handles} contains the
topology of handle decompositions and bundles. Section~\ref{sec:families}
treats neck-stretching, families and regularity in the invariant
perturbation class. Sections~\ref{sec:collapse}--\ref{sec:collapse-coverage}
prove Theorem~\ref{thm:B}. Section~\ref{sec:two-handles} proves the two-handle square.
Section~\ref{sec:torus} proves the torus excision theorem (Theorem~\ref{thm:C}),
which the one-handle comparison requires. Section~\ref{sec:one-three} proves
the one- and three-handle squares, and Section~\ref{sec:presentations} treats
changes of presentation and frame. Section~\ref{sec:assembly} proves
Theorem~\ref{thm:A}. Appendix~\ref{app:dfl-regular} proves that primary
perturbations regularizing the gauge complex, as required in \cite{DFL},
exist at the metric data of a presentation. A table of notation follows.

\subsection*{Notation}

\begin{center}
\small
\begin{tabular}{>{$}l<{$}p{0.7\linewidth}}
\toprule
M^\#,\ E^\#_M & torus-summed object and its bundle (Section~\ref{sec:category})\\
\phi_M & lifting group, the subgroup $\phi$ of \cite[Section~5.1]{KM-unknot}\\
\Sigma_h,\ \Msf_h & Heegaard surface of genus $h$ and its moduli space of odd flat connections\\
C_0,\ C_1,\ L(C) & compression bodies of a presentation and their Lagrangians;
 in Sections~\ref{sec:collapse}--\ref{sec:collapse-coverage} also $P=L(C)$\\
p,\ N_p & presentation and its mixed map\\
P,\ N_P,\ N^\dagger_P & primary data of $p$, $N_P=N_p$, and $N^\dagger_P=N_{\bar P}^\vee$\\
C_G(M,p),\ C_S(M,p) & gauge and symplectic complexes\\
\mathfrak h_i,\ f,\ \eta & primary perturbations, end perturbation, secondary perturbation\\
\kappa,\ \mathfrak c & topological energy and the index constants of \cite[Proposition~2.53]{DFL}\\
\mathrm{cont}_S,\ \mathrm{cont}_G & symplectic and gauge continuation maps\\
X_0,\ X_\eps & collapsing family (Section~\ref{subsec:insertion})\\
I,\ N_I & collapsing arc and its fixed neighbourhood\\
\mathfrak s,\ \mathsf A & compression scale, and $\mathsf A=1+\log(1/\eps)$\\
W_{\mathsf R},\ W_\infty,\ W^J_\lambda & collared, completed and joined domains (Definition~\ref{def:torus-domains})\\
\mu_T & spectral gap of the odd flat connection on the torus\\
R_i,\ \widetilde P_i,\ L_C & height Lagrangians, lower Lagrangians, correspondence of $C_{h,g}$\\
\mathbb X,\ A_0,\ A_1 & product target $\Msf_h^-\times\Msf_g$ and the folded pair in it\\
\mathcal K,\ \mathsf E & folding isomorphism, Floer's excision map\\
U_0,\ U_1 & in Sections~\ref{sec:torus}--\ref{sec:one-three}, compression bodies of the folded splitting\\
P_{\mathsf L}(V),\ \mathcal G_{\mathsf L}(V) & preglued configuration and glued solution at torus length $\mathsf L$\\
\mathcal U(F) & admissible neighbourhood of primary data (Definition~\ref{def:admissible-nbhd})\\
U_C,\ V_C & compression cup and cap (Definition~\ref{def:cup})\\
\bottomrule
\end{tabular}
\end{center}

\subsection*{Acknowledgements}
The author used LLMs during the development of this work for literature
exploration, editorial and coding assistance and as an adversarial reader.

\section{The framed cobordism category and the gauge functor}\label{sec:category}

\subsection{The torus attachment}

Write $T^3=S^1\times T^2$ and let $P_T\to T^3$ be the pullback of the
$\SO(3)$ bundle on $T^2$ with $w_2\ne0$. Fix once and for all a ball
$B_T\subset T^3$ disjoint from a fibre $\{\theta_0\}\times T^2$, a
parametrization of $B_T$ by the unit ball $B^3$, a trivialization of
$P_T|_{B_T}$ and a metric on $T^3$ that is standard near $B_T$ and flat on a
neighbourhood of the fibre $\{\theta_0\}\times T^2$; the flat metric it induces
on $T^2$ is used for every torus collar below, in all presentations. Put
$T^\circ=T^3\setminus\operatorname{int}B_T$. Let
\begin{equation}\label{eq:alphaT}
\alpha_T=\operatorname{pr}_{S^1}^*(\text{generator of }H^1(S^1;\FF))
\in H^1(T^3;\FF).
\end{equation}
Then $\langle w_2(P_T),[T^2]\rangle=1$, $\alpha_T|_{T^2}=0$ and
$\langle\alpha_T\smile w_2(P_T),[T^3]\rangle=1$.

\begin{definition}\label{def:category}
An \emph{object} of $\bB$ is a closed, connected, oriented three-manifold
$M$ together with an orientation-preserving embedding $b\colon B^3\to M$.
The \emph{torus-summed object} is
\[
M^\#=(M\setminus b(\operatorname{int}B^3))\cup_{S^2}T^\circ,
\qquad
E_M^\#=\underline{\RR}^3\cup P_T|_{T^\circ},
\]
glued along the boundary spheres by the fixed parametrizations and
trivializations.

A \emph{decorated cobordism} $(W,E,\nu)\colon(M_0,b_0)\to(M_1,b_1)$
consists of:
\begin{itemize}[leftmargin=2em]
\item a compact, connected, oriented four-manifold $W$ with
  $\partial W=-M_0\sqcup M_1$ and product collars at both ends;
\item an embedding $\nu\colon B^3\times[0,1]\to W$ with
  $\nu(\cdot,i)=b_i$, a product in the collars (the \emph{tube});
\item an $\SO(3)$ bundle $E\to W$ with trivializations over $M_0$,
  $M_1$ and over the image of $\nu$, which agree on the overlaps.
\end{itemize}
An \emph{isomorphism} of decorated cobordisms is an
orientation-preserving diffeomorphism covered by a bundle isomorphism,
both respecting the boundary identifications, the tubes and the
trivializations. The morphisms of $\bB$ are isomorphism classes of
decorated cobordisms. Composition glues along the ends using the
collars, tubes and trivializations; the identity of $(M,b)$ is the
product $[0,1]\times M$ with the vertical tube and the product bundle.
\end{definition}

Associativity and the unit laws hold on isomorphism classes, by the
usual collar-reparametrization isomorphisms. A product
$[0,1]\times M$ whose bundle identifications at the two ends differ by a
nontrivial automorphism is not the identity. Its image under the functor
below is the corresponding bundle automorphism map.

For a decorated cobordism put $W^\circ=W\setminus\nu(\operatorname{int}
B^3\times[0,1])$. Then
\begin{equation}\label{eq:Wsharp}
W^\#=W^\circ\cup_{S^2\times[0,1]}(T^\circ\times[0,1]),\qquad
E^\#=E|_{W^\circ}\cup\bigl(P_T|_{T^\circ}\times[0,1]\bigr),
\end{equation}
glued by the tube parametrization and the trivializations. Its ends are
$-M_0^\#\sqcup M_1^\#$ with the bundles $E^\#_{M_i}$.

\begin{lemma}\label{lem:H1-sharp}
Restriction gives an isomorphism
$H^1(W^\#;\FF)\to H^1(W;\FF)\oplus H^1(T^3;\FF)$.
\end{lemma}

\begin{proof}
$W^\circ$ is connected, since a path in the four-manifold $W$ can be
displaced off the core arc of the tube and then off its tubular
neighbourhood. The pieces $W^\circ$ and $T^\circ\times[0,1]$ are
connected, their intersection $S^2\times[0,1]$ is connected and has
vanishing $H^1$, so Mayer--Vietoris gives
$H^1(W^\#)=H^1(W^\circ)\oplus H^1(T^\circ)$. The same argument applied to
$W=W^\circ\cup(B^3\times[0,1])$ and $T^3=T^\circ\cup B_T$ identifies
$H^1(W^\circ)=H^1(W)$ and $H^1(T^\circ)=H^1(T^3)$.
\end{proof}

Let $\alpha_W\in H^1(W^\#;\FF)$ correspond to $(0,\alpha_T)$ and put
$\phi_W=\{0,\alpha_W\}$. On each end $M_i^\#$ the class $\alpha_W$
restricts to the class $\alpha_{M_i}$ corresponding to $(0,\alpha_T)$,
and $\phi_{M_i}=\{0,\alpha_{M_i}\}$. The class $\alpha_W$ is realized by
a bundle automorphism of $E^\#$: take the standard automorphism of $P_T$
with lifting class $\alpha_T$ (the gauge transformation $g$ of
\cite[Section~7.1]{Scaduto}), deform it to be the identity over $B_T$,
use it constantly along $T^\circ\times[0,1]$, and extend by the identity
over $W^\circ$.

\subsection{The gauge functor}

Recall from \cite[Section~5.1]{KM-unknot} that for an $\SO(3)$ bundle on
a three-manifold the quotient of the automorphism group by the image of
the determinant-one gauge group is $H^1(\,\cdot\,;\FF)$, so a subgroup
$\phi\subset H^1$ acts on the configuration space of connections modulo
determinant-one gauge transformations. The class in $H^1(\,\cdot\,;\FF)$ of an automorphism $g$ is
its \emph{lifting class} $c(g)$, the obstruction to lifting $g$ to a
determinant-one gauge transformation; $g$ is \emph{even} if its lifting class vanishes. A subgroup $\phi$ of lifting classes used in this way is the
\emph{lifting group} of the problem; it is the subgroup $\phi$ of \cite[Section~5.1]{KM-unknot}. A triple satisfies the
\emph{$\phi$-non-integral condition} if it contains a surface on which
$w_2$ is odd and $\phi$ vanishes \cite[Definition~5.1]{KM-unknot}. For
$M^\#$ the torus $\{\theta_0\}\times T^2$ is such a surface for
$\phi_M$. Under this condition $\phi$ acts freely on the critical set,
the instanton homology of the quotient $I(M^\#,E^\#_M)^{\phi_M}$ is
defined, and a cobordism with a subgroup $\phi\subset H^1(W;\FF)$
restricting to the given subgroups at the ends induces a map
$I(W^\#,E^\#)^{\phi}$ by counting points in the quotient moduli spaces
\cite[pp.~154--155]{KM-unknot}. Over $\FF$ no orientation data is
needed.

\begin{definition}\label{def:Ish}
For an object $(M,b)$ put $\Ish(M)=I(M^\#,E^\#_M;\FF)^{\phi_M}$. For a
decorated cobordism put
$\Ish(W,E,\nu)=I(W^\#,E^\#;\FF)^{\phi_W}\colon\Ish(M_0)\to\Ish(M_1)$.
\end{definition}

Here $I(M^\#,E^\#_M;\FF)^{\phi_M}$ denotes a single vector space: the Floer
groups for different metrics and perturbations are identified by continuation
maps, which form a transitive system, and $\Ish(M)$ is its limit, that is, any
one of the groups together with these identifications.
The group $\Ish(M)$ is the framed instanton homology of
\cite{KM-unknot, Scaduto} with $\FF$ coefficients. Since
$\langle\alpha_T\smile w_2(P_T),[T^3]\rangle=1$, the nonzero element of
$\phi_M$ shifts the grading of the determinant-one (unquotiented) theory by four
(the proof of \cite[Proposition~5.2(a)]{KM-unknot}), and the determinant-one Floer
homology of $(M^\#,E^\#_M)$ is two copies of $\Ish(M)$.

\begin{proposition}\label{prop:gauge-functor}
$\Ish$ is a functor $\bB\to\mathrm{Vect}_{\FF}$.
\end{proposition}

\begin{proof}
\emph{Well defined.} Isomorphic decorated cobordisms give isomorphic
data $(W^\#,E^\#,\phi_W)$, hence equal maps by the usual continuation
argument for metrics and perturbations \cite[Section~3]{KM-unknot},
\cite[Section~5.3]{Donaldson-book}.

\emph{Composition.} Let $W_{01}\colon M_0\to M_1$ and
$W_{12}\colon M_1\to M_2$ be composable, with composite $W_{02}$. Then
$W_{02}^\#=W_{01}^\#\cup_{M_1^\#}W_{12}^\#$ as bundles. By
\cite[Proposition~5.4]{KM-unknot},
$I(W_{02}^\#)^{\phi}=I(W_{12}^\#)^{\phi_{12}}\circ
I(W_{01}^\#)^{\phi_{01}}$ provided that $\phi$ contains the image $\psi$
of the Mayer--Vietoris connecting map
$H^0(M_1^\#;\FF)\to H^1(W_{02}^\#;\FF)$, and $\phi$ is the largest
subgroup restricting to $\phi_{01}$ and $\phi_{12}$ on the pieces. Since
$M_1^\#$ is connected, the restriction
$H^0(W_{01}^\#)\oplus H^0(W_{12}^\#)\to H^0(M_1^\#)$ is onto, so
$\psi=0$. The class $\alpha_{W_{02}}$ restricts to $\alpha_{W_{01}}$ and
$\alpha_{W_{12}}$. A class restricting to $0$ or to $\alpha$ on both
pieces is determined by its restrictions, because
$H^1(W_{02}^\#)\to H^1(W_{01}^\#)\oplus H^1(W_{12}^\#)$ is injective when
$\psi=0$. Hence the largest such subgroup is $\phi_{W_{02}}$.

\emph{Unit.} For the product $[0,1]\times M$ with translation-invariant
regular data the index-zero solutions are the constant trajectories.
The group $\phi_M$ acts freely on the critical set. A trajectory from a
critical point $\beta$ to $\alpha_M\cdot\beta$ has index in the determinant-one theory
congruent to $4$ modulo $8$ by the proof of \cite[Proposition~5.2(a)]{KM-unknot}, so
it does not contribute in index zero. The map is therefore the identity
on the chain level for such data, and on homology for all data.
\end{proof}

\begin{remark}
The functor uses the enlarged gauge group of \cite[Section~5]{KM-unknot}
with the subgroup $\phi_W$. It is not the projective functor of
\cite[Section~4.1]{KM-unknot}, which uses determinant-one gauge
transformations and yields two copies of $\Ish$. Scaduto enlarges the
gauge group on $W^\#$ by the transformations that restrict to $g$ on the
torus \cite[Section~7.2]{Scaduto}; restricted to classes compatible with
the ends, this is $\phi_W$. The subgroup $\phi_W$ vanishes on the
original piece $W^\circ$. Classes of $H^1(W)$ that vanish on $\partial W$
are not added.
\end{remark}

\begin{remark}\label{rem:grading}
With the $\ZZ/4$ grading of \cite[Proposition~7.1]{Scaduto},
$\Ish(W,E,\nu)$ is homogeneous of degree
$-\tfrac32(\chi(W)+\sigma(W))+\tfrac12(b_1(M_1)-b_1(M_0))+2P(E)$ modulo
$4$, where $P(E)$ is the mod-two self-intersection of a surface in the
interior of $W$ dual to $w_2(E)$ relative to the trivialized ends.
\end{remark}

\section{Heegaard data, the two complexes and the mixed map}\label{sec:models}

\subsection{The moduli space of a surface}

For $h\ge1$ let $\Sigma_h$ be a closed oriented surface of genus $h$ and
$F_h\to\Sigma_h$ the $\SO(3)$ bundle with $w_2\ne0$. Let $\Msf_h$ be the
moduli space of flat connections on $F_h$ modulo determinant-one
automorphisms. It is a closed, simply connected, monotone symplectic
manifold of dimension $6h-6$. For $h\ge2$ it has
$\pi_2(\Msf_h)=\ZZ$ and minimal Chern number $2$; for $h=1$ it is
a point \cite{AB, Newstead}. Its tangent space at $[\alpha]$ is
$H^1(\Sigma_h;\ad\alpha)$.

\begin{convention}[orientation]\label{conv:orientation}
Let $C$ be an oriented three-manifold and $\Sigma\subset\partial C$ a
boundary component with the outward-normal-first orientation. On the
moduli space of $\Sigma$ viewed as a boundary component of $C$ we use
\[
\omega(v,w)=\int_\Sigma\langle v\wedge w\rangle,\qquad
J=*_\Sigma,
\]
where $\langle\cdot,\cdot\rangle$ is the inner product $\langle a,b\rangle=
-\tr(ab)$ on $\ad F\cong\mathfrak{su}(2)$, used for all $L^2$ norms. In the notation of
\cite{DFL}, whose bilinear form $\tr$ is $-\tfrac12$ of the Euclidean inner product on $\RR^3\cong\mathfrak{so}(3)$,
the matching surface is $\Sigma_{\mathrm D}=-\partial C$, with
$\Omega_{\mathrm D}=-\omega$ and $J_{\mathrm D}=-J$; the two
descriptions define the same mixed equation. Here $\Sigma$ is the fibre of
$U_-\times\Sigma$ in the four-manifold of \cite{DFL}: its upper arm is
$[2,\infty)_\theta\times Y$ with $\partial Y$ at $s=0$ and outward normal
$+\partial_s$, so the orientation is
$d\theta\wedge ds\wedge\mathrm{vol}_{\partial Y}=-ds\wedge d\theta\wedge\mathrm{vol}_{\partial Y}$
and the fibre is $-\partial Y$.
\end{convention}

Convention~\ref{conv:orientation} is not a matter of taste. It is the
sign for which the Green identity of the compression body cancels the
Green identity of the adjacent strip (Lemma~\ref{lem:green}). With the
opposite sign the uniform estimate of
Proposition~\ref{prop:garding} fails, already for constant scale
(Remark~\ref{rem:wrong-sign}).

\subsection{Compression bodies and their Lagrangians}\label{subsec:compressions}

A \emph{torus compression body} of genus $h$ is a compact oriented
three-manifold $C$ with $\partial C=\Sigma_h\sqcup T^2$ obtained from a
handlebody $H$ of genus $h-1$ by boundary sum with $[0,1]\times T^2$,
together with the $\SO(3)$ bundle that is pulled back from the odd bundle
on $T^2$. More generally we consider compression bodies $C$ from
$\Sigma_h$ to $\Sigma_{h'}$ with the odd bundle. Every flat connection on
$C$ is irreducible with $H^0(C;\ad)=H^2(C;\ad)=0$, the flat moduli space
$\mathcal R(C)$ is a compact manifold, and restriction to the boundary is a
Lagrangian embedding into the product of the moduli spaces of the boundary
components with the orientation of $\partial C$ \cite[Proposition~2.21]{DFL}.
For a torus compression body the moduli space of the torus boundary is a
point, so restriction to $\Sigma_h$ is a Lagrangian embedding
\[
L(C)\subset\Msf_h ,
\]
with $L(C)\cong(S^3)^{h-1}$, simply connected with $\pi_2=0$. The pair
$(L(C_0),L(C_1))$ of two torus compression bodies is monotone with minimal
Maslov number $4$ \cite[Theorem~1]{DFL}. For a compression body from
$\Sigma_h$ to $\Sigma_{h'}$ with $h'\ge2$ the Lagrangian lies in
$\Msf_h^-\times\Msf_{h'}$ (Definition~\ref{def:clean-datum} and
Section~\ref{subsec:folded}).

The properties used in the collapse theorem are the following.

\begin{definition}\label{def:clean-datum}
Let $C$ be a compact oriented three-manifold
with a metric $g_C$ that is a product near $\partial C$, and an
$\SO(3)$ bundle $\mathbb E\to C$. We impose the following conditions:
\begin{enumerate}[label=(C\arabic*)]
\item $w_2(\mathbb E)$ is nonzero on every component of $\partial C$;
\item every flat connection $B$ on $\mathbb E$ has
  $H^0(C;\ad B)=H^2(C;\ad B)=0$;
\item the restriction $H^1(C;\ad B)\to H^1(\partial C;\ad B)$ is
  injective, and the flat moduli space of $C$ restricts to an embedded
  Lagrangian $L(C)$ in the product of the moduli spaces of the
  boundary components with the signs of
  Convention~\ref{conv:orientation}.
\end{enumerate}
\end{definition}

The compression bodies of Sections~\ref{subsec:compressions}
and~\ref{subsec:clean-data} satisfy these conditions. Beyond (C1)--(C3), Sections~\ref{sec:collapse}--\ref{sec:collapse-coverage}
use only that $L(C)$ is compact and connected, not the structure of a compression
body.

For $v\in T_qL(C)$ let $h_v$ denote the harmonic one-form on $C$ with
values in $\ad B_q$, with vanishing normal component and boundary
cohomology class $v$. The map $v\mapsto h_v$ is well defined because
$H^1(C;\ad)\to H^1(\Sigma_h;\ad)$ is injective.

\emph{Primary perturbations.} Following \cite[Section~2.3]{DFL}, a
\emph{primary perturbation} of $C$ is a holonomy perturbation $\mathfrak h$
supported in the interior of $C^\circ$, the part of $C$ outside the unit collar on which
\cite[(2.1) and p.~9]{DFL} require the perturbations to vanish (the special quintuple
and the torus analysis of Section~\ref{sec:torus} use this collar unperturbed), small
enough that the perturbed flat
connections still form a compact manifold whose restriction
$L_{\mathfrak h}(C)\subset\Msf_h$ is an embedded Lagrangian. It is Hamiltonian isotopic
to $L(C)$: for $\mathfrak h$ small the implicit function theorem makes $L_{\mathfrak h}(C)$
$C^1$-close to $L(C)$, hence the graph of a closed one-form in a Weinstein
neighbourhood, and the one-form is exact because $L(C)$ is simply connected.
(\cite{DFL} call these perturbations cylinder functions. In their notation $C_0$ and $C_1$
are the halves $Y_1=Y_0\cup[0,1]\times\Sigma$ and $Y_1'$ of $Y_\#$ cut along $\{0\}\times\Sigma$
\cite[p.~79]{DFL}, and $C_i^\circ$ is $Y_0$, $Y_0'$; Section~\ref{sec:torus} writes $(Y_\pm)_0$.) The universal restriction map from pairs (perturbation, flat
connection) to $\Msf_h$ is a submersion \cite[Section~2.3]{DFL}, so
finitely many primary perturbations make any finite collection of such
Lagrangians pairwise transverse with empty triple intersections.

\emph{Trace functions.} The determinant connection is fixed, so the holonomy of $B$ around a
loop $\ell$ is $\lambda_\ell\,g_\ell(B)$ with $g_\ell(B)\in\SU(2)$ and $\lambda_\ell\in U(1)$
independent of $B$ \cite[Section~5.6, p.~147]{Donaldson-book}; $\lambda_\ell$ is determined up to
sign. Over the loops of one embedding below we choose it continuously; this is possible because
they are parametrized by a disk, and the choice is unique up to one sign. For an embedding
$\gamma\colon S^1\times D^2\to\operatorname{int}C^\circ$ and a two-form $\mu$ supported in the
interior of $D^2$ with $\int\mu=1$, the \emph{averaged trace function} is
$\tau_{\gamma}(B)=\int_{D^2}\tau_z(B)\,\mu$, where $\tau_z(B)=\tr g_\ell(B)$ for the loop
$\ell=\gamma(S^1\times\{z\})$; a cylinder function is
$G(\tau_{\gamma_1},\dots,\tau_{\gamma_n})$ with $G$ smooth, the $\gamma_j$ being allowed to be
immersions \cite[Section~6.1]{DFL}. Daemi, Fukaya and Lipyanskiy use the trace of the holonomy
itself; since $\tr g_\ell=\lambda_\ell^{-1}\tr\mathrm{Hol}_\ell$ with $\lambda_\ell$ independent of
$B$, the results of \cite[Section~6]{DFL} used below hold for these functions with the same
proofs. On flat connections $\tau_\gamma$ is $\tr g_\ell$ for the core
$\ell=\gamma(S^1\times\{0\})$, and, with $\lambda$ transported along a free homotopy, it depends
only on the free homotopy class of the core. The transport is ambiguous only along a
self-homotopy sweeping a torus on which $w_2$ of the bundle is nonzero, and there
$\tr g_\ell$ vanishes on flat connections, because the lifts of the two generators of that torus
anticommute. A finite family $F$ of primary perturbations of $(C_0,C_1)$ is \emph{separated} if
there are finitely many embeddings $\gamma_\alpha\colon S^1\times D^2\to\operatorname{int}C_0^\circ
\sqcup\operatorname{int}C_1^\circ$ with pairwise disjoint images such that every member of $F$ is
a real affine combination of averaged trace functions of the $\gamma_\alpha$. The height
functions of Sections~\ref{sec:two-handles} and~\ref{sec:torus} enter the arguments only through
their restrictions to flat connections. We fix, once for each Heegaard splitting, a realization of
them as affine combinations of averaged trace functions of embedded solid tori, all the solid
tori used in one compression body being pairwise disjoint, with parallel copies where two cores
are freely homotopic; for instance $h_Q=\sum_j(1-\tfrac12\tau_{\gamma_j})$ in
Section~\ref{sec:two-handles}, where the signs of the $\lambda$ are chosen so that
$\tau_{\gamma_j}=\tr b_j$ on flat connections, with the lifts $b_j$ fixed there.
We also fix, once for each Heegaard splitting, a finite separated family $F_*$ of primary
perturbations of $(C_0,C_1)$ containing the height functions of Section~\ref{sec:torus}, such
that the universal
restriction map of $\operatorname{span}F_*$ is a submersion at every flat connection of $C_0$ and
of $C_1$; this is possible by compactness and \cite[Section~2.3]{DFL}, and
Lemma~\ref{lem:app-separated}(i) provides a separated one. Every finite family of primary
perturbations considered in Sections~\ref{sec:one-three}--\ref{sec:assembly} contains $F_*$,
and the generic perturbations of the auxiliary height Lagrangians of
Section~\ref{subsec:one-three-squares} are taken from $\operatorname{span}F_*$, so that they do
not depend on the family. The height functions of Section~\ref{sec:two-handles} belong to the
families of the presentations compared there, which may be enlarged
(Definition~\ref{def:admissible-nbhd}).

\subsection{Admissible presentations}\label{subsec:admissible}

\begin{definition}\label{def:presentation}
Let $(M,b)$ be an object of $\bB$. \emph{Presentation data} $p$
for $M$ consist of:
\begin{enumerate}[label=(\roman*)]
\item a Heegaard splitting $M=H_0\cup_\Sigma H_1$ of genus $g\ge1$
  adapted to $b$;
\item the torus compression bodies $C_i=H_i\natural([0,1]\times T^2)$,
  the boundary sums taken inside $b(B^3)$ through the equatorial disk,
  so that $C_0\cup C_1$ is the torus-summed object $M^\#$ of
  Definition~\ref{def:category}, with $h=g+1\ge2$;
\item identifications of collars of $\partial C_0$, $\partial C_1$ as in \cite[Section~2.1]{DFL},
  and metrics on $C_0$, $C_1$ that are products on these collars and
  agree on $\Sigma_h$ and on the collar of $T^2$, which carries the flat metric fixed
  in Section~\ref{sec:category} (the same for all presentations) (\emph{metric data});
  the product collar $I\times\Sigma_h$ of $\Sigma_h$ in $M^\#$ given by these identifications
  is compatible with the gluing sphere $S=b(\partial B^3)$ of Definition~\ref{def:category}:
  it meets $S$ in $I\times\partial D$, where $D$ is the equatorial disk of (ii), so that,
  with $T^\circ=T^3\setminus\operatorname{int}B_T$ as in Section~\ref{sec:category},
  every collar line through a point of $\Sigma_h\cap T^\circ$ lies in $T^\circ$ and every
  collar line through a point of $\Sigma_h\setminus T^\circ$ lies in
  $M\setminus b(\operatorname{int}B^3)$;
\item primary perturbations $\mathfrak h_0$, $\mathfrak h_1$ with
  $L_{\mathfrak h_0}(C_0)\pitchfork L_{\mathfrak h_1}(C_1)$ (\emph{primary data});
\item the regularizing data of Section~\ref{subsec:complexes}, with the
  end perturbation $f$ and the secondary perturbation below the
  thresholds of Theorem~\ref{thm:mixed-iso}.
\end{enumerate}
They form an \emph{admissible presentation} if, for some finite family $F\supseteq F_*$ of
primary perturbations, with $F_*$ the family fixed above, the primary data lie in a fixed open
neighbourhood $\mathcal U(F)$ of zero in the span of $F$, depending only on the splitting, the
metric data and $F$ (Definition~\ref{def:admissible-nbhd}): the primary data are joined to zero
by a path in $\operatorname{span}F$ along which the Lagrangians of the primary data, and those
of their reflections, stay in a $C^1$-neighbourhood of the unperturbed pair that is small
relative to the heights of the auxiliary height functions of Section~\ref{subsec:folded},
which are fixed with the splitting.
\end{definition}

Theorem~\ref{thm:mixed-iso} and Sections~\ref{sec:families}--\ref{sec:torus}
apply to all presentation data, and more generally to the pairs of
Remark~\ref{rem:clean-pairs}; admissibility is used only in
Sections~\ref{sec:one-three}--\ref{sec:assembly}.

\subsection{The two complexes}\label{subsec:complexes}

Let $p$ be presentation data and put $P_i=L_{\mathfrak h_i}(C_i)$.

\emph{The symplectic complex.} $C_S(M,p)=\CF(P_0,P_1)$ is the Floer
complex over $\FF$ generated by $P_0\cap P_1$, with the differential
counting index-one holomorphic strips for a generic $t$-dependent
compatible almost complex structure equal near the boundary to $J_{\mathrm D}=-J$ in the
strip coordinates $(s,\theta)$, that is, to $J$ in the coordinates $(t,r)=(\theta,s)$ of
Convention~\ref{conv:orientation} \cite{Oh, DFL}. Monotonicity and minimal Maslov number $4$ exclude disk
bubbling in dimensions at most one.

\emph{The gauge complex.} Glue $C_0$ and $-C_1$ along $\Sigma_h\sqcup
T^2$ to obtain $M^\#$ with its bundle. The primary perturbations define a
perturbed Chern--Simons functional $\mathrm{CS}_{\mathfrak h}$ on $M^\#$, whose critical
orbits correspond bijectively to $P_0\cap P_1$ and are nondegenerate
\cite[Lemmas~2.45, 2.47]{DFL}. The functional $\mathrm{CS}_{\mathfrak h}$ is invariant
under $\alpha_M$: since $\alpha_M$ vanishes on each $C_i$, the $\SU(2)$ traces of
loops in $C_i^\circ$ are unchanged \cite[p.~19]{DFL}. The nonzero element of
$\phi_M$ acts freely on the critical set. Choose an \emph{end perturbation}
$f$: a finite holonomy perturbation of the Chern--Simons functional on $M^\#$, invariant under
$\phi_M$, vanishing on a fixed neighbourhood of the critical set, and
such that all nonconstant trajectories of $\mathrm{CS}_{\mathfrak h}+f$ of index at most seven in the determinant-one theory are regular. Such $f$ exist, may be taken
arbitrarily small in the Banach norm of holonomy perturbations, and do not change the critical set or the Hessians
(Proposition~\ref{prop:input-functions}). Regularity is achieved by $f$ rather than
by the primary perturbations because in the two-handle comparison the primary
perturbation of the middle compression body must vanish
(Remark~\ref{rem:dfl-map}). If the trajectories of $\mathrm{CS}_{\mathfrak h}$ of index at most seven
are already regular, $f=0$ is allowed; we then call the primary data
\emph{DFL-regular}. These are the primary data required in \cite[Lemma~2.56]{DFL} and
Theorem~3.29 there; they exist at the metric data of every presentation
(Theorem~\ref{thm:dfl-regular}). The gauge complex
$C_G(M,p)$ is the Floer complex of $\mathrm{CS}_{\mathfrak h}+f$ modulo $\phi_M$. Its
homology is canonically $\Ish(M)$.

Both complexes have the same generator set $P_0\cap P_1$. This is an
identification of generators, not of differentials.

\begin{proposition}\label{prop:input-functions}
For every choice of primary data there are end perturbations $f$ as above,
arbitrarily small in the Banach norm of holonomy perturbations
\cite[Section~5.5]{Donaldson-book} and vanishing on a fixed neighbourhood of the
critical set.
\end{proposition}

\begin{proof}
This is the standard transversality argument for holonomy perturbations
\cite[Section~5.5]{Donaldson-book}, carried out with perturbations that vanish
near the critical set; we indicate the three points to check.

\emph{The perturbations reach every nonconstant trajectory.} The critical set is
finite and nondegenerate. Let $\gamma$ be a nonconstant trajectory of index at most seven in the determinant-one theory. Since $\gamma$ is nonconstant, there is a time
$s_0$ at which $\gamma(s_0)$ lies outside a fixed neighbourhood of the critical
set, and there the perturbation $f$ may be nonzero.

\emph{The perturbations detect cokernel elements.} Let $\psi$ be a nonzero element
of the cokernel of the linearization at $\gamma$. By unique continuation for the
trajectory equation, forward and backward in time
\cite[Proposition~8.6 and Theorem~C.1]{SW-lagrangian}, both $\psi$ and
$\partial_s\gamma$ are nonzero at $s_0$, which we may choose generic. They are
moreover linearly independent there: $\partial_s\gamma$ lies in the kernel of the
linearization and $\psi$ in the kernel of its adjoint, so
$\frac{d}{ds}\langle\psi(s),\partial_s\gamma(s)\rangle=0$ because the Hessian is
symmetric, and the pairing vanishes by decay at the ends. This rules out a
cokernel element proportional to $\partial_s\gamma$ near $s_0$, which no
perturbation localized near $\gamma(s_0)$ could detect. Holonomy
functions separate $\gamma(s_0)$ from the rest of $\gamma$, from its translate
$\alpha_M\cdot\gamma$ and from the critical set, and their gradients at $\gamma(s_0)$
span a dense subspace of the tangent space of the configuration space (compare
the finite embedding argument of \cite[Section~5.5, Lemma~5.16]{Donaldson-book}).
So some holonomy direction pairs nontrivially with $\psi$, and the Sard--Smale
theorem makes all such trajectories regular for generic $f$.

\emph{Invariance and finiteness.} Averaging over $\phi_M$ preserves the
construction: a cylinder function whose gradient is supported, in the configuration
space, near $\gamma(s_0)$ has its $\alpha_M$-translate supported near $\alpha_M\cdot\gamma(s_0)$,
and these neighbourhoods are disjoint because $\phi_M$ acts freely. There are finitely many energy levels below the
index-seven bound, and the choice is made by induction on the energy.
\end{proof}

\emph{Existence of DFL-regular primary data.} Let $\mathfrak h$ be primary data with transverse
Lagrangians. For $i\in\{0,1\}$ let $\cB(C_i^\circ)$ be the space of connections on $C_i^\circ$
modulo all determinant-one gauge transformations of $C_i^\circ$, with no boundary condition, and
let $X_i\subset\cB(C_i^\circ)$ be the finite set of restrictions of the critical points of
$\mathrm{CS}_{\mathfrak h}$. A \emph{split addition} is a finite sum $g=g_0+g_1$ of cylinder
functions, $g_i$ with loops in $\operatorname{int}C_i^\circ$ and vanishing on a neighbourhood of
$X_i$. For an irreducible connection $B$ on $C_i^\circ$ let $\pi^i_B$ be the $L^2$-orthogonal
projection of $L^2(C_i^\circ;\Lambda^1\otimes\ad)$ onto the orthogonal complement of the closed
subspace $d_BL^2_1(C_i^\circ)$. One-forms on $M^\#$ are \emph{linearly independent modulo gauge
on $C_i^\circ$} at $B$ if the images of their restrictions under $\pi^i_B$ are linearly
independent. Trajectories are taken in temporal gauge, in the determinant-one theory, between
lifts of critical points; the cokernel of the linearization at a trajectory $\gamma=(A_s)_{s\in\RR}$
is represented by the decaying solutions $\psi$ of the adjoint equation, and $\psi(s)$ is a
one-form on $M^\#$ (Lemma~\ref{lem:app-trajectories}).

\begin{lemma}\label{lem:dfl-detect}
Let $\mathfrak h$ be primary data in the span of a separated family, with transverse
Lagrangians, and let $g$ be a split addition so small that $\mathrm{CS}_{\mathfrak h+g}$ and
$\mathrm{CS}_{\mathfrak h}$ have the same critical points. Let $\gamma=(A_s)_{s\in\RR}$ be a
nonconstant trajectory of $\mathrm{CS}_{\mathfrak h+g}$, of any index and in any path class, and
$\psi\ne0$ an element of the cokernel of its linearization. Then for each $i\in\{0,1\}$:
\begin{enumerate}[label=(\roman*)]
\item the path $s\mapsto[A_s|_{C_i^\circ}]$ in $\cB(C_i^\circ)$ is injective and does not meet
  $X_i$;
\item there is $T$ such that for every $s_0<-T$ the restriction $A_{s_0}|_{C_i^\circ}$ is
  irreducible, and $\partial_s\gamma(s_0)$ and $\psi(s_0)$ are linearly independent modulo gauge
  on $C_i^\circ$;
\item there is a cylinder function $g_i'$ with loops in $\operatorname{int}C_i^\circ$, vanishing
  on a neighbourhood of $X_i$, such that
  $\int_\RR\langle\nabla g_i'(A_s),\psi(s)\rangle\,ds\ne0$.
\end{enumerate}
\end{lemma}

\begin{theorem}\label{thm:dfl-regular}
Fix a Heegaard splitting and metric data as in Definition~\ref{def:presentation}(i)--(iii).
There are families $F$ as in Definition~\ref{def:presentation} that are separated, and for each of them
generic primary data in a neighbourhood of zero in $\operatorname{span}F$ have transverse
Lagrangians. Let $F$ be a family as in Definition~\ref{def:presentation} that is separated, and
$\mathfrak h\in\mathcal U(F)$ primary data with transverse Lagrangians. Then there are
arbitrarily small split additions $g$, each lying in the span of a finite family $G$ of
cylinder functions with loops in
$\operatorname{int}C_0^\circ\sqcup\operatorname{int}C_1^\circ$, such that
$\mathrm{CS}_{\mathfrak h+g}$ has the same critical
points as $\mathrm{CS}_{\mathfrak h}$, the primary data $\mathfrak h+g$ are DFL-regular, and
$\mathfrak h+g\in\mathcal U(F\cup G)$; presentation data with primary data $\mathfrak h+g$ are
therefore admissible. The family $G$ may depend on the required smallness of $g$, and its loops
may be taken in $\operatorname{int}C_0^\circ$ alone. The metric data are not changed.
\end{theorem}

Lemma~\ref{lem:dfl-detect} and Theorem~\ref{thm:dfl-regular} are proved in
Appendix~\ref{app:dfl-regular}. Near $s=-\infty$ both $\partial_s\gamma$ and $\psi$ are governed by
eigenvectors of the Hessian at the limiting critical point, with eigenvalues of opposite signs.
The Hessian acts on each piece separately, so eigenvectors with distinct eigenvalues stay
independent modulo gauge on each piece, provided no eigenvector vanishes on the unperturbed part
of the neck. That last property is where separatedness is used
(Proposition~\ref{prop:nontrapping}).

\begin{remark}\label{rem:dfl-regular}
Daemi, Fukaya and Lipyanskiy obtain DFL-regular primary data after a constant rescaling of the
metric of $M^\#$ \cite[Lemma~2.56 and Proposition~6.8]{DFL}. The cylinder functions added in the
inductive step of their proof are split, with loops in the interior of one of the two pieces,
and the step on \cite[p.~80]{DFL} passes from the linear independence of $\psi$ and
$\partial_s\gamma$ on $M^\#\setminus\Sigma$ to split cylinder functions that pair nontrivially
with $\psi$. For split functions the condition needed concerns one piece: independence of
$\psi$ and $\partial_s\gamma$ modulo gauge there, together with injectivity and irreducibility
of the restricted path. We could not deduce independence on one piece from the
$L^2$-orthogonality of $\psi$ and $\partial_s\gamma$ on $M^\#$: $\psi=c_i\,\partial_s\gamma$
modulo gauge on $C_i$, with constants $c_0\ne c_1$, is compatible with it, and for such $\psi$
every cylinder function with loops in one piece pairs to zero. The isomorphism of
\cite[Theorem~3.29]{DFL} is not affected, since Theorem~\ref{thm:dfl-regular} supplies the
regular data it requires.
Lemma~\ref{lem:dfl-detect} provides these conditions on each piece when the base data lie in a
separated family. Theorem~\ref{thm:dfl-regular} is therefore an alternative proof of the existence
statement of \cite[Proposition~6.8]{DFL} for base data in a separated family, with the added
perturbations supported in $\operatorname{int}C_0^\circ\sqcup\operatorname{int}C_1^\circ$ and without
rescaling. For a general base perturbation in the Banach space of \cite[Section~6.1]{DFL} the
argument would need the conclusion of Proposition~\ref{prop:nontrapping}, which we prove only for
separated data. The rescaling preserves the regularity of the trajectories, which it only
reparametrizes in time, and the Lagrangians, since $*_3\nabla h$ does not change under a constant
rescaling of the metric; it changes the flat metric of the torus collar, which
Definition~\ref{def:presentation}(iii) fixes.
\end{remark}

\subsection{The mixed map}\label{subsec:mixed-iso}

Let $q_s$ be the special quintuple of \cite[Section~3.1]{DFL} associated
to the compression bodies $C_0$, $C_1$: a four-manifold $X$ containing
$\RR\times C_0^\circ$, $\RR\times(-C_1^\circ)$ and $\mathbb U_-\times\Sigma_h$,
with a gauge-theoretic end $(-\infty,-3]\times M^\#$, and a surface
$S=\mathbb U_+$ with a symplectic end $[3,\infty)\times[-1,1]$ and two boundary
components labelled by $P_0$ and $P_1$. Here $\CC$ is cut by the matching line
$\mathbb U_\partial=i\RR$ into $\mathbb U_-$ and $\mathbb U_+$; the planar region $\mathbb U_-$ contains the input
arm $(-\infty,-3]\times[-1,1]$, two arms $[-2,0]\times[2,\infty)$ and
$[-2,0]\times(-\infty,-2]$ and the square $[-1,0]\times[-2,2]$, and $\RR\times C_0^\circ$
and $\RR\times(-C_1^\circ)$ are glued to $\mathbb U_-\times\Sigma_h$ along the two boundary
curves of $\mathbb U_-$ other than $\mathbb U_\partial$ \cite[Figure~1]{DFL}. The ends of the two
arms are the mixed ends of the quintuple. We call a domain of this type, for
any pair covered by the results below, a \emph{DFL domain}. A \emph{DFL disk} is a DFL domain obtained by cutting a larger
domain along ordinary strip ends.

\begin{figure}[ht]
\centering
\begin{tikzpicture}[scale=0.62,font=\small]
\fill[orange!18] (-6,-1) -- (-3,-1) .. controls (-2.3,-1) and (-2,-1.4) .. (-2,-2) -- (-2,-4) -- (0,-4) -- (0,4) -- (-2,4) -- (-2,2) .. controls (-2,1.4) and (-2.3,1) .. (-3,1) -- (-6,1) -- cycle;
\fill[gray!12] (6,-1) -- (3,-1) .. controls (2.3,-1) and (2,-1.4) .. (2,-2) -- (2,-4) -- (0,-4) -- (0,4) -- (2,4) -- (2,2) .. controls (2,1.4) and (2.3,1) .. (3,1) -- (6,1) -- cycle;
\draw[thick] (-6,1) -- (-3,1) .. controls (-2.3,1) and (-2,1.4) .. (-2,2) -- (-2,4);
\draw[thick] (-6,-1) -- (-3,-1) .. controls (-2.3,-1) and (-2,-1.4) .. (-2,-2) -- (-2,-4);
\draw[thick,blue!70!black] (6,1) -- (3,1) .. controls (2.3,1) and (2,1.4) .. (2,2) -- (2,4);
\draw[thick,blue!70!black] (6,-1) -- (3,-1) .. controls (2.3,-1) and (2,-1.4) .. (2,-2) -- (2,-4);
\draw[dashed] (0,-4) -- (0,4);
\node at (-1.2,0.35) {$\mathbb U_-\times\Sigma_h$};
\node at (1.1,0) {$\mathbb U_+$};
\node[left] at (-6,0) {$M^\#$};
\node[right] at (6,0) {$P_0\cap P_1$};
\node[above] at (-4.5,1) {$\RR\times C_0$};
\node[below] at (-4.5,-1) {$\RR\times(-C_1)$};
\node[blue!70!black,above] at (4.5,1) {$P_0$};
\node[blue!70!black,below] at (4.5,-1) {$P_1$};
\node[above] at (0,4) {mixed end of $C_0$};
\node[below] at (0,-4) {mixed end of $C_1$};
\node[right] at (0.05,-3) {$\mathbb U_\partial$};
\end{tikzpicture}
\caption{The planar base of the special quintuple of \cite[Section~3.1]{DFL}, after
\cite[Figure~1]{DFL}. Over the region $\mathbb U_-$ (shaded) the gauge domain is
$\mathbb U_-\times\Sigma_h$, with $\RR\times C_0$ and $\RR\times(-C_1)$ glued along the two boundary
curves; the section is defined over $\mathbb U_+$, with boundary conditions $P_0$ and
$P_1$; the two are matched along the line $\mathbb U_\partial$. The ends are the gauge
input $M^\#$, the strip output at $P_0\cap P_1$, and two mixed ends along the matching
line.}
\label{fig:quintuple}
\end{figure}
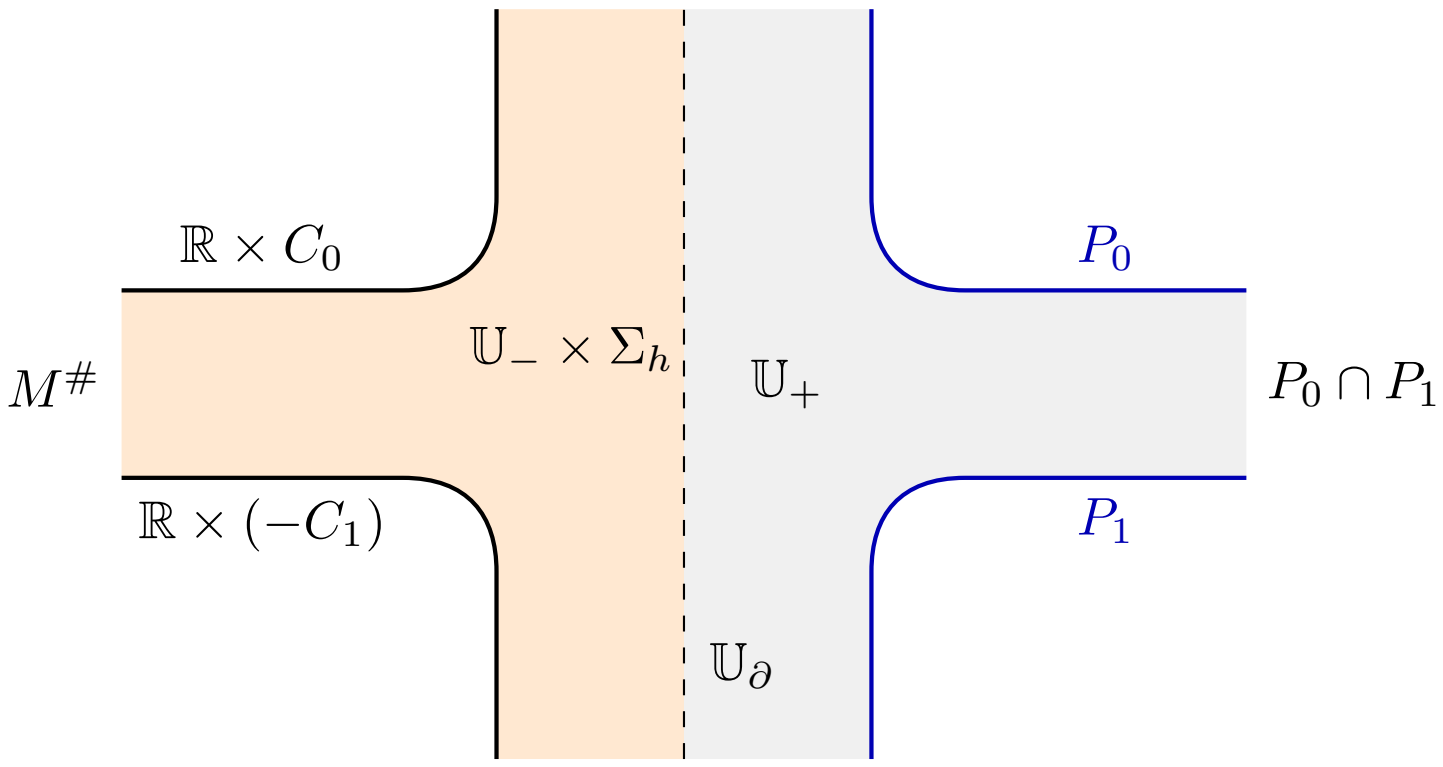
 The torus $T^2$ and its point
moduli space are carried along; the matching condition on
$\mathbb U_\partial\times T^2$ imposes flatness of the torus restriction. On the
gauge-theoretic end we add $\chi f$, where $\chi$ is a cutoff equal to
$1$ near $-\infty$ and supported before the matching region, and we
choose a secondary perturbation $\eta$: a finite holonomy perturbation
supported in a compact subset of $X$ away from the matching line and
invariant under $\phi_M$, together with a compact deformation of the
almost complex structure away from the matching line. The mixed equation
is
\begin{equation}\label{eq:mixed}
F_A^++(*_3\nabla\mathfrak h)^++\chi(*_3\nabla f)^++\eta(A)=0,\qquad
\partial_su+J_{s,\theta}(u)\partial_\theta u=0,
\end{equation}
with the matching condition along $\mathbb U_\partial\times\Sigma_h$ and the
Lagrangian boundary conditions on $S$. Here $J_{s,\theta}$ is a family of compatible almost
complex structures equal to $J_{\mathrm D}=-J$ near the matching line and near $\theta=\pm1$, in
the coordinates $(s,\theta)$ of \cite[Section~3.1]{DFL}. For generators $a$, $x$ let
$\cM^d(a,x)$ be the space of finite-energy solutions of index $d$
asymptotic to $a$ at the gauge end and to $x$ at the symplectic end,
modulo the gauge group enlarged by $\phi_M$.

\begin{theorem}\label{thm:mixed-iso}
\begin{enumerate}[label=(\roman*)]
\item For $f$ sufficiently small and invariant secondary perturbations $\eta$
  that are sufficiently small, including the deformation of the almost complex
  structure, and generic, every $\cM^d(a,x)$ with $d\le1$ is regular, the
  zero-dimensional spaces are finite, and
  \[
  N_p(a)=\sum_x\#\cM^0(a,x)\,x
  \]
  is a chain map $C_G(M,p)\to C_S(M,p)$. With respect to the filtration
  by the constants $\mathfrak c$ of \cite[Proposition~2.53]{DFL}, which
  normalize the index formula,
  $N_p(a)=a+(\text{terms of lower filtration})$; in particular $N_p$ is a
  chain isomorphism.
\item The chain homotopy class of $N_p$, composed with the gauge
  continuation map between two choices of $f$, is independent of $f$,
  $\chi$ and $\eta$.
\end{enumerate}
\end{theorem}

We call $N_p$ the \emph{mixed map} at $p$. In Section~\ref{sec:presentations} the
presentations are the vertices of a graph whose edges are changes of presentation.

\begin{proof}
\emph{Energy and index.} Let $\kappa$ be the topological energy normalized as in
\cite[Proposition~3.20]{DFL} (Section~\ref{subsec:energy-convention}). For a solution
of~\eqref{eq:mixed} in the homotopy class of $(A,u)$, the analytic energy equals
\[
8\pi^2\kappa(A,u)+2\|\eta(A)\|^2_{L^2(X)}+2\int\chi'(s)f(A_s)\,ds ,
\]
as in \cite[(3.24)]{DFL} with the additional end term, whose absolute value is at most
$2\sup|f|$, since $\int|\chi'|=1$ and $f$ vanishes at the critical point at the
gauge end. Indeed, if $K=F_A+*_3\nabla\mathfrak h$ and
$P=\chi*_3\nabla f$, then $P\wedge P=0$ and, with the inner product
$-\tr$, the cross term is
$2\int\tr(K\wedge P)=-2\int\chi\,\frac{d}{ds}f(A_s)\,ds
=2\int\chi'f(A_s)\,ds$. The index is
$\ind D_{(A,u)}=8\kappa(A,u)+\mathfrak c(x)-\mathfrak c(a)$; the \emph{asymptotic labels} of a configuration are the
critical points, intersection points or components of critical manifolds to which it
converges at its ends. For two configurations with the same labels,
$2(\kappa-\kappa')\in\ZZ$ \cite[Lemma~3.15]{DFL}, and $\mathfrak c$ is the function on
generators of \cite[Proposition~2.53]{DFL}, determined up to an additive constant
by this formula.

\emph{Regularity.} Transversality for secondary perturbations is obtained
by DFL's induction over finite-dimensional families of perturbations and the
Sard--Smale theorem \cite[Section~6.2]{DFL}. The
perturbations can be chosen invariant under $\phi_M$ by
Lemma~\ref{lem:invariant-span}: holonomies along loops with vanishing lifting class span the relevant test fields. Regularity in index at most one
is achieved by a finite choice, proceeding upward from the lowest index
allowed by the energy bound.

\emph{Chain map.} The one-dimensional spaces compactify with boundary
the broken configurations of \cite[Proposition~3.26]{DFL}: a gauge
trajectory followed by an index-zero mixed solution, or an index-zero
mixed solution followed by a holomorphic strip. Sphere and disk bubbles
and energy concentration are excluded in dimension at most one by the
index-energy relation and monotonicity \cite[Section~5.4]{DFL}. Hence
$\partial_SN_p=N_p\partial_G$.

\emph{Triangularity.} Fix $\kappa_+>0$ below every positive value of
$|\mathfrak c(a)-\mathfrak c(x)|/8$. An index-zero solution from $a$ to $x$ has
$\kappa=(\mathfrak c(a)-\mathfrak c(x))/8$ by the index formula, and since its analytic
energy is nonnegative,
$8\pi^2\kappa\ge-2\sup|f|-2\sup_A\|\eta(A)\|^2_{L^2(X)}$. So if
$2\sup|f|+2\sup_A\|\eta(A)\|^2_{L^2(X)}<8\pi^2\kappa_+$, then $\kappa\ge0$, hence
$\mathfrak c(a)\ge\mathfrak c(x)$, with equality only if $\kappa=0$. The constant pairs of
\cite[Section~3.2]{DFL} are regular by \cite[Proposition~3.19]{DFL}, so by the implicit
function theorem, for small $(f,\eta)$, a fixed neighbourhood of each constant pair contains
exactly one solution, of index zero and with $\kappa=0$. The second threshold, below which
these are all the index-zero solutions with $\kappa=0$, is obtained by contradiction. Suppose
that for data $(f_n,\eta_n)\to0$ there are index-zero solutions $c_n$ with $\kappa=0$ outside
these neighbourhoods. Their analytic energies are at most
$2\sup|f_n|+2\sup_A\|\eta_n(A)\|^2_{L^2(X)}\to0$. A nonconstant trajectory of
$\mathrm{CS}_{\mathfrak h}+f_n$ between critical points, or a nonconstant holomorphic strip, has
energy bounded below independently of $n$, because $f_n$ vanishes near the critical set; so no
energy escapes along the ends, and a subsequence of $(c_n)$ converges to a configuration of zero
energy with the same asymptotic labels, that is, to a constant pair; breaking
along the two mixed ends is excluded by the energy quantum and exponential decay
there \cite[Proposition~3.6 and Section~5.4]{DFL}. The decay of the $c_n$ along
the ends is uniform in $n$, again because $f_n$ vanishes near the critical set, so the
convergence holds in the topology in which the neighbourhoods above are open, and $c_n$ lies in
one of them for large $n$, a contradiction. Below both thresholds, therefore, the index-zero
solutions with $\kappa=0$ are exactly one for each generator, from the generator to itself. This
is the argument of \cite[Theorem~3.29]{DFL}.

\emph{Independence.} At fixed $f$, a path of choices $(\chi_t,\eta_t)$ has
parameter derivative supported in a compact core, and defines a
parametrized moduli space whose index $-1$ elements occur at isolated
parameters and are glued using the parametrized operator
$(\xi,t)\mapsto D\xi+t\,\partial_\lambda F$
(Theorem~\ref{thm:family}). Its one-dimensional part has boundary
$N^1-N^0$ and the broken configurations, giving a chain homotopy
$N^1\simeq N^0$. A change of $f$ changes the equation on the whole gauge end
and is not a compactly supported variation, so it is treated separately:
insert a long gauge cylinder carrying the continuation from $f_0$ to $f_1$ at
the gauge end. By Theorem~\ref{thm:neck} applied at the closed cut $M^\#$, for
long necks the count on the glued domain is $N_{f_1}\circ \mathrm{cont}_G(f_0,f_1)$, while
the glued domain is a compactly supported deformation of the domain for
$f_0$, so the previous argument identifies the count with $N_{f_0}$ up to
chain homotopy.
\end{proof}

\begin{remark}\label{rem:clean-pairs}
Theorem~\ref{thm:mixed-iso} and its proof hold verbatim for any pair of
three-manifolds satisfying Definition~\ref{def:clean-datum}, with common boundary,
possibly disconnected, product target and transverse primary Lagrangians,
including boundary markings twisted by a frame transition: in particular for
the pairs $(C_0,C_j)$ with the $\rho_\zeta$-marking and $(C_1,C_2)$ of
Section~\ref{sec:two-handles} (whose critical sets are described by
Lemma~\ref{lem:twisted-critical}), and for the pairs $(D_i,C)$ and $(D_0,D_1)$
of Section~\ref{subsec:folded}, whose matching surface is
$\Sigma_h\sqcup\Sigma_g$ and target $\Msf_h^-\times\Msf_g$. The inputs
\cite[Lemma~3.15, Proposition~3.20]{DFL} and the triangularity argument use
only embeddedness and monotonicity of the Lagrangians, and
\cite{DFL-mixed} allows a disconnected odd matching surface. The two further
inputs of the triangularity step also hold: regularity of the constant pairs
\cite[Proposition~3.19]{DFL} reduces to transversality of the pair at the
corresponding intersection point, and each listed pair is an admissible
splitting in the sense of \cite{DFL}, with two odd boundary components, so
\cite[Propositions~2.53 and~3.20]{DFL} apply to it directly.
\end{remark}

\begin{remark}\label{rem:dfl-map}
The map $N_p$ is defined by the construction of \cite[Section~3.2]{DFL}, with two
changes: the gauge complex may be regularized by an end perturbation instead of by
the primary perturbations, and the secondary perturbations are invariant under
$\phi_M$. The first is needed in the two-handle comparison, where the primary
perturbation of the middle compression body must vanish, so that
\cite[Lemma~2.56]{DFL} is not available there. The second makes the counts
independent of the lift of a generator to the determinant-one configuration space:
DFL count solutions asymptotic to a fixed representative of each generator
\cite[p.~27]{DFL}, and pass to the quotient by the involution $\iota$ of
\cite[p.~18]{DFL}, which is $\alpha_M$. Where both constructions apply they agree
(Proposition~\ref{prop:dfl-comparison}).
\end{remark}

\begin{proposition}[comparison with the DFL isomorphism]\label{prop:dfl-comparison}
Let $p$ be presentation data, or more generally a pair of torus compression bodies as in
Remark~\ref{rem:clean-pairs} with collars and metrics as in
Definition~\ref{def:presentation}(iii) but an arbitrary flat metric on the torus collar, with
DFL-regular primary data. Take $f=0$; a
family of almost complex structures as in \cite[Section~3.1 and Lemma~2.35]{DFL},
equal near $\theta=\pm1$ to DFL's standard structure, which is
$J_{\mathrm D}=-J$ in DFL's coordinates $(s,\theta)$, so that in the coordinates
$(t,r)=(\theta,s)$ of Sections~\ref{sec:collapse}--\ref{sec:collapse-linear} the
equation is the $J$-holomorphic equation; and a $\phi_M$-invariant secondary perturbation
$\eta$, sufficiently small and supported, with the deformation of the almost complex
structure, in a compact set away from the matching line as in
\cite[Proposition~3.23(i)]{DFL}, for which the mixed moduli spaces of expected dimension
at most three are regular. Such $\eta$ exist, the isomorphism $N^{\rm DFL}$ of
\cite[Theorem~3.29]{DFL} is defined at these data, and
\[
N_p=N^{\rm DFL}\otimes\FF\colon C_G(M,p)\to C_S(M,p)
\]
as chain maps.
\end{proposition}

\begin{proof}
\emph{The complexes.} DFL's gauge complex is generated by the critical points of
the perturbed functional modulo determinant-one gauge transformations and the
involution $\iota$ given by the class Poincaré dual to a component of the matching
surface \cite[pp.~18--19, Lemmas~2.43 and~2.45]{DFL}. In the framed splitting this
component may be taken to be the torus, so $\iota=\alpha_M$ and DFL's generators are
ours. With $f=0$ the gauge equations agree, so the differentials agree mod two; the
strip complexes agree because the Lagrangians and almost complex structures do.
Regular families equal to the standard structure near $\theta=\pm1$ exist, since the
transversality perturbations for nonconstant strips can be supported in the interior of
the strip.

\emph{The equations.} With $f=0$, equation~\eqref{eq:mixed} is \cite[(3.22)]{DFL}
on the same special quintuple; Convention~\ref{conv:orientation} records that the
two sign conventions define the same equation.

\emph{The counts.} DFL's zero-dimensional space $\cM_\eta(\tilde a,x)_0$ consists of
solutions asymptotic to a fixed determinant-one representative $\tilde a$ of $a$,
modulo determinant-one gauge transformations. The automorphism $\alpha_M$ extends over
the special quintuple \cite[Remark~3.7]{DFL}, it maps solutions for $\eta$ to solutions
for $\alpha_M^*\eta$, and $\alpha_M^*\eta=\eta$ by invariance. So it identifies the spaces
for the two representatives $\tilde a$, $\alpha_M\tilde a$, and $\cM_\eta(\tilde a,x)_0$ maps
bijectively onto our $\cM^0(a,x)$, the space of solutions modulo the gauge group
enlarged by $\phi_M$ (which acts freely by Theorem~\ref{thm:regularity}(i)). Hence the
counts agree mod two.

\emph{Existence of $\eta$ and DFL's hypotheses.} By Theorem~\ref{thm:regularity}(ii) the
linearized operator together with the $\phi_M$-invariant holonomy directions supported
in a compact set away from the matching line is onto at every solution, and the
Sard--Smale theorem, applied finitely many dimensions at a time, gives $\eta$, small
and compactly supported away from the matching line, with all moduli spaces of
expected dimension at most three regular. These are the requirements of
\cite[Propositions~3.23 and~3.26]{DFL}, whose remaining conclusions (nonnegative
energy, the unique solution in $\cM_\eta(a,a)_0$, compactness and the boundary
description) follow from smallness and compact support as in
\cite[Section~5]{DFL}.
\end{proof}

By Theorem~\ref{thm:dfl-regular}, at every Heegaard splitting and all metric data there are
admissible presentations with DFL-regular primary data $\mathfrak h+g$, where $\mathfrak h\in\mathcal
U(F)$, for a family $F$ as in Definition~\ref{def:presentation} that is separated, is any primary data with
transverse Lagrangians and $g$ is an arbitrarily small split addition;
Proposition~\ref{prop:dfl-comparison} applies at these presentations.

\subsection{Energy, gradings and filtrations}\label{subsec:energy-convention}

We use DFL's normalization throughout: the topological energy $\kappa$
of a mixed pair is the Chern--Weil integral of the gauge part plus the
symplectic area of the section, divided so that
$\ind=8\kappa+\text{const}$ \cite[Proposition~3.20]{DFL}, and the analytic
energy of a solution is $8\pi^2\kappa$ plus the secondary terms
\cite[(3.24)]{DFL}. For two configurations with the same asymptotic labels
$2(\kappa-\kappa')\in\ZZ$ \cite[Lemma~3.15]{DFL}; we write $\kappa_0$ for the
resulting coset representative.
Thus every nonzero index difference between classes with the same
labels has absolute value at least $4$. We use only this divisibility,
without asserting that every allowed set of classes realizes the
smallest possible period. In this normalization a change
$\Delta\kappa=\tfrac12$ changes $8\pi^2\kappa$ by $4\pi^2$.
For a pure strip the topological energy is
$\kappa(u)=\frac1{4\pi^2}\int u^*\omega$, since its full Dirichlet
energy is twice its symplectic area
\cite[Proposition~4.69]{DFL}. In particular, a positive generator
sphere, when $h\ge2$, contributes $\Delta\kappa=\tfrac12$ and index $4$ already
in the determinant-one surface moduli space. No enlargement of
the gauge group is needed for this sphere contribution.

\section{Handle decompositions and bundles}\label{sec:handles}

This section contains the topology needed to write every morphism of
$\bB$ as a composite of elementary pieces. All handle decompositions are
relative to the tube: handles are attached away from $\nu$, and the tube
is a product throughout.

\subsection{Decompositions relative to the tube}

\begin{proposition}\label{prop:handles}
Let $(W,E,\nu)\colon M_0\to M_1$ be a morphism of $\bB$. Then $W$ has a
handle decomposition relative to $M_0\times[0,1]\cup\nu(B^3\times[0,1])$
with handles of index $1$, $2$ and $3$ only, attached in this order and
away from the tube, such that every regular level is connected. All
two-handles can be attached along one framed link in a single level.
Consequently
\[
W=W_3\circ W_2\circ W_1,
\]
where $W_1\colon M_0\to Z_0$ consists of one-handles, $W_2\colon Z_0\to
Z_1$ is the trace of surgery on a framed link, and $W_3\colon Z_1\to M_1$
consists of three-handles, with $Z_0,Z_1$ connected.
\end{proposition}

\begin{proof}
Let $W^\circ$ be the complement of the open tube; it is connected
(proof of Lemma~\ref{lem:H1-sharp}). Choose a Morse function
$f\colon W^\circ\to[0,1]$ that equals the $[0,1]$ coordinate on a collar
$S^2\times[0,\delta)\times[0,1]$ of the vertical boundary and near the
ends, and a gradient-like vector field equal to $\partial_t$ on the
collar. The collar is a union of flow lines, so every stable and
unstable manifold of an interior critical point, and every modification
below, lies outside it.

\emph{Order of indices.} If a $k$-handle is attached after a $j$-handle
with $k\le j$, the attaching sphere of dimension $k-1$ and the belt
sphere of dimension $3-j$ in a three-dimensional level satisfy
$(k-1)+(3-j)<3$. An isotopy supported off the collar separates them, and
the critical values can be interchanged
\cite[Section~4]{Milnor-h-cob}, \cite[Propositions~4.2.7 and~4.2.13]{GS}.

\emph{No zero- or four-handles.} After ordering, the union of the incoming
collar with the zero- and one-handles is connected, because later handles
have connected attaching regions and cannot join components, while
$W^\circ$ is connected. Choose a spanning tree of the graph whose
vertices are the zero-handles and the incoming collar and whose edges are
the one-handles. A leaf zero-handle meets its tree one-handle in one
point of the attaching sphere, and the pair cancels
\cite[Theorem~5.4]{Milnor-h-cob}. Proceeding from the leaves removes all
zero-handles. Four-handles are removed by the same argument applied to
$1-f$.

\emph{Connected levels.} Attaching a one-handle to a connected level
gives a connected level. Surgery on a circle in a connected
three-manifold gives a connected manifold, since the complement of a
circle is connected. Surgery on a two-sphere changes the number of
components by $0$ or $+1$. Starting from the connected $M_0$, all levels
before the three-handles are connected. Since there are no four-handles,
the number of components does not decrease from the first three-handle
on, and it equals one at $M_1$. Hence all levels are connected.

\emph{One level of two-handles.} A two-handle can be attached below the
belt circle of an earlier two-handle, since $1+1<3$. Pushing all
attaching circles into the level just above the one-handles gives a
single framed link.
\end{proof}

\begin{proposition}\label{prop:outer-bundles}
In the decomposition of Proposition~\ref{prop:handles}, the
trivialization of $E$ over $M_0$ and the tube extends over $W_1$, and
the trivialization over $M_1$ and the tube extends over $W_3$. In
particular $E|_{Z_0}$ and $E|_{Z_1}$ are trivial, and $(W,E,\nu)$ is
the composite of $(W_1,\underline{\RR}^3)$, $(W_2,E|_{W_2})$ with the
induced frames at $Z_0$, $Z_1$, and $(W_3,\underline{\RR}^3)$.
\end{proposition}

\begin{proof}
Up to homotopy equivalence relative to $M_0\cup\nu$, the pair
$(W_1,M_0\cup\nu)$ is $(M_0\times[0,1],M_0\times0\cup b_0\times[0,1])$
with relative one-cells attached along pairs of points. The
trivialization extends over the product. Over a one-handle
$D^1\times D^3$ the two transition functions on the attaching balls are
null-homotopic maps $D^3\to\SO(3)$, and a path in the connected group
$\SO(3)$ joining their values extends the frame across the handle. The
statement for $W_3$ is the same with $f$ replaced by $1-f$. The
decorated pieces glue to $(W,E,\nu)$ by construction.
\end{proof}

The extension in Proposition~\ref{prop:outer-bundles} is not unique:
frames on $Z_0$ may differ by an element of $\pi_1(\SO(3))$ for each
one-handle. This ambiguity is absorbed by the frame-change arrows of
Section~\ref{subsec:frame-change}, and it is invisible in the final
statement, because the composite morphism does not depend on it.

\subsection{One- and three-handles and adapted presentations}

A Heegaard splitting $M=H_0\cup_\Sigma H_1$ is \emph{adapted to} $b$ if
$\Sigma\cap b(B^3)$ is the equatorial disk $D$ and $\Sigma$ is standard
on a collar of $b(B^3)$. The attaching balls of a one-handle can be moved,
by an isotopy of the level fixed near $b(B^3)$, into a neighbourhood of
a Heegaard surface meeting $b(B^3)$ in $D$, because the complement of a
ball in a connected three-manifold is connected. The isotopy extends to
the cobordism by a diffeomorphism fixed on the tube, so the decorated
morphism is unchanged. Reading $W_3$ from $M_1$, each three-handle
becomes a one-handle attached to $-M_1$, and its attaching sphere is the
belt sphere of that one-handle; by Proposition~\ref{prop:handles} it is
nonseparating in its level.

\subsection{Subordinate Heegaard triples and the tube}\label{subsec:triples}

Let $L\subset Z$ be a framed link of $n$ components in a connected
three-manifold with a ball $b$ disjoint from $L$, and let $W(Z,L)$ be the
trace of surgery. By \cite[Lemma~4.5]{OS-triangles} there is a Heegaard
triple $(\Sigma;H_0,H_1,H_2)$ of genus $g$ subordinate to a bouquet for
$L$, in the sense of \cite[Definition~4.2]{OS-triangles}. After an
isotopy in the complement of the bouquet, $\Sigma$ meets $b(B^3)$ in the
equatorial disk. Then $Y_{01}=H_0\cup H_1=Z$,
$Y_{02}=H_0\cup H_2=Z(L)$ and $Y_{12}=H_1\cup H_2\cong\#^r(S^1\times S^2)$
with $r=g-n$. Let $X$ be the four-manifold
\[
X=(\Delta\times\Sigma)\cup(e_0\times H_0)\cup(e_1\times H_1)\cup
(e_2\times H_2)
\]
of \cite[Section~2.2]{OS-triangles}, where $\Delta$ is a triangle with
edges $e_0,e_1,e_2$, and let
$H_r=\natural^r(S^1\times D^3)$.

\begin{proposition}\label{prop:tube}
There is a diffeomorphism
$\Phi\colon X\cup_{Y_{12}}H_r\to W(Z,L)$ that is the identity on $Z$ and
on $Z(L)$, and that restricts to the identity from a tube
$D_z\times[-1,1]\times[-1,1]$ onto $b(B^3)\times[-1,1]\subset
Z\times[-1,1]$, where $D_z\subset\Sigma$ is a disk around the basepoint
missing all attaching curves. The core of the tube is carried to the
product arc with its product normal framing. The same holds for any
other identification of the filling with $H_r$, after an isotopy
relative to the boundary.
\end{proposition}

\begin{proof}
Let $V$ be a regular neighbourhood of the bouquet, a handlebody of genus
$n$ containing $L$. By \cite[Definition~4.2]{OS-triangles},
$Z\setminus\operatorname{int}V=H_0\cup_\Sigma\mathcal K$, where
$\mathcal K$ is the compression body obtained from $\Sigma\times[0,1]$
by attaching two-handles along the curves $\beta_{n+1},\dots,\beta_g$;
moreover $H_1=\mathcal K\cup V$ and $H_2=\mathcal K\cup V(L)$, where
$V(L)$ is $V$ surgered along $L$. Fix a bicollar
$\Sigma\times[-1,1]_s\subset H_0\cup\mathcal K$ and put
$B=D_z\times[-1,1]_s$, a ball with $B\cap\Sigma=D_z$, disjoint from $L$; we
choose the bicollar so that $B=b(B^3)$ as parametrized balls.

\emph{A product model of $X$.} For a planar hexagon $Q$ whose sides
alternate between edge sides $\partial_0,\partial_1,\partial_2$ and
vertex sides, set $X(Q)=(\Sigma\times Q)\cup\bigcup_i H_i^0\times
\partial_iQ$, where $H_i^0$ is $H_i$ minus a collar, with corners
smoothed. A diffeomorphism of hexagons matching labelled sides, and a
product near the edge sides, induces a diffeomorphism of the
corresponding spaces. Take
$Q'=[-1,1]_s\times[-1,1]_t\setminus\bigl((0,1]\times(-\epsilon,\epsilon)
\bigr)$ with $\partial_0=\{s=-1\}$, $\partial_1=\{s=1,t\le-\epsilon\}$
and $\partial_2=\{s=1,t\ge\epsilon\}$. Then
\[
X\cong Z\times[-1,-\epsilon]\ \cup\ H_0\times[-\epsilon,\epsilon]\ \cup\
Z(L)\times[\epsilon,1],
\]
and near the basepoint $X$ contains the product $D_z\times Q'$.

\emph{The filling.} The complement of $X$ in
$(Z\setminus V)\times[-1,1]$ is $\mathcal K^{s\ge0}\times
(-\epsilon,\epsilon)$. Let $Q_V$ be a four-ball bounding the three-sphere
$V\times\{-\epsilon\}\cup\partial V\times[-\epsilon,\epsilon]\cup
V(L)\times\{\epsilon\}$, a Heegaard splitting of $S^3$ of genus $n$.
Put $H=(\mathcal K^{s\ge0}\times[-\epsilon,\epsilon])\cup Q_V$. Since
$\mathcal K\times I$ is $\partial V\times I\times I$ with $g-n$
four-dimensional one-handles, $H\cong\natural^r(S^1\times D^3)$ with
$\partial H=Y_{12}$. Near the basepoint $H$ contains the half-ball
$Q_z=D_z\times[0,1]_s\times[-\epsilon,\epsilon]$.

\emph{The diffeomorphism.} Now
$X\cup H=(Z\setminus V)\times[-1,1]\cup W'_V$ with
$W'_V=V\times[-1,-\epsilon]\cup Q_V\cup V(L)\times[\epsilon,1]$, while
$W(Z,L)=(Z\setminus V)\times[-1,1]\cup W_V$ with $W_V$ the union of
$V\times[-1,1]$ and the two-handles along $L\times\{1\}$. Both $W_V$ and
$W'_V$ are four-balls: in $W_V$ each two-handle cancels a distinct
one-handle of $V\times I$. Their boundaries are identified by the
identity on $V\times\{-1\}\cup\partial V\times[-1,1]$ and the
tautological identification of the surgered handlebody at the top. By
Cerf's theorem \cite{Cerf} this boundary diffeomorphism extends over the
balls. Define $\Phi$ to be the identity on $(Z\setminus V)\times[-1,1]$
and this extension on $W'_V$.

\emph{The tube.} The tube $(D_z\times Q')\cup Q_z=B\times[-1,1]$ lies in
$(Z\setminus V)\times[-1,1]$, where $\Phi$ is the identity. Its core and
its normal frame $(\partial_x,\partial_y,\partial_s)$ are carried to
themselves. For another filling $H'_r$, a diffeomorphism $H'_r\to H$
can be chosen to be the identity on the boundary by Laudenbach--Po\'enaru
\cite{LP}. Since $Q_z$ is a regular neighbourhood of a boundary ball
inside a collar, and collars are unique up to isotopy, the diffeomorphism
may be isotoped relative to the boundary to preserve $Q_z$.
\end{proof}

\subsection{Bundles on a two-handle cobordism}\label{subsec:two-handle-bundles}

Let $V=H^1(\Sigma;\FF)$, and let $\Lambda_i\subset V$ be the image of $H^1(H_i;\FF)$, a
Lagrangian subspace of dimension $g$. Put $W_2=X\cup_{Y_{12}}H_r$.

\begin{proposition}\label{prop:bundle-normal-form}
Let $\delta\colon V^3\to H^2(X;\FF)$ be the Mayer--Vietoris coboundary
for the cover of $X$ by $\Delta\times\Sigma$ and $e_i\times H_i$. Then
$x\mapsto\delta(x,0,0)$ induces an isomorphism
\[
V\big/\bigl(\Lambda_0+(\Lambda_1\cap \Lambda_2)\bigr)\xrightarrow{\ \cong\ }H^2(W_2;\FF).
\]
The class $\delta(x,0,0)$ restricts to $x$ modulo $\Lambda_0+\Lambda_1$ on $Y_{01}$
and to $x$ modulo $\Lambda_0+\Lambda_2$ on $Y_{02}$. In particular every $\SO(3)$
bundle on $W_2$ that is trivial over both ends has
$w_2=\delta(x,0,0)$ with $x\in(\Lambda_0+\Lambda_1)\cap(\Lambda_0+\Lambda_2)$.
\end{proposition}

\begin{proof}
The pieces $e_i\times H_i$ are pairwise disjoint and meet
$\Delta\times\Sigma$ in copies of $e_i\times\Sigma$. Mayer--Vietoris
reads
\[
V\oplus\textstyle\bigoplus_iH^1(H_i)\xrightarrow{\rho}V^3
\xrightarrow{\delta}H^2(X)\to H^2(\Delta\times\Sigma)\to
H^2(\Sigma)^3,
\]
with $\rho(v,u_0,u_1,u_2)=(v-u_0,v-u_1,v-u_2)$. The last map is the
diagonal $\FF\to\FF^3$ and is injective, so
$H^2(X)=(V/\Lambda_0\oplus V/\Lambda_1\oplus V/\Lambda_2)/\Delta(V)$. The same cover of
$Y_{ij}$ gives $H^2(Y_{ij})=V/(\Lambda_i+\Lambda_j)$, with restriction
$[x_0,x_1,x_2]\mapsto x_i-x_j$. Since $H^1(H_r)\to H^1(Y_{12})$ is onto
and $H^2(H_r)=0$, Mayer--Vietoris for $W_2=X\cup H_r$ identifies
$H^2(W_2)$ with the kernel of $H^2(X)\to H^2(Y_{12})$. For a class in
this kernel, $x_1-x_2\in \Lambda_1+\Lambda_2$. Subtracting $\Delta(x_2)$ and then
$\Delta(a_2)$ for suitable $a_2\in \Lambda_2$ brings it to the form
$[x,0,0]$. Two such classes agree exactly when $x-x'\in \Lambda_0+(\Lambda_1\cap
\Lambda_2)$.
\end{proof}

\begin{example}
For the $(-1)$-framed unknot in $S^3$ one may take $g=1$,
$\alpha=\lambda$, $\beta=\mu$, $\gamma=\lambda-\mu$. Then $r=0$,
$\Lambda_1\cap \Lambda_2=0$ and $H^2(W_2;\FF)=\FF^2/\langle\lambda\rangle\cong\FF$,
the class of the nontrivial bundle on $\overline{\CC P}^2$ minus two
balls. Its restrictions to both ends vanish.
\end{example}

\begin{definition}\label{def:normal-form}
Let $E$ be a bundle on $W_2$ trivial over both ends, with
$w_2(E)=\delta(\zeta,0,0)$. A \emph{filling frame} for $E$ consists of a
trivialization of $E$ over the central product $\Delta\times\Sigma$, over
$e_1\times H_1$, $e_2\times H_2$ and $H_r$ with identity transitions, and
over $e_0\times H_0$ with transition a map $\rho_\zeta\colon\Sigma\to\SO(3)$ of
lifting class $\zeta$, equal to the identity near $D_z$; the frame agrees with
the given trivialization of $E$ along the tube.
\end{definition}

Filling frames exist, since $[\Sigma/D_z,\SO(3)]=H^1(\Sigma;\FF)$ (the map
$\SO(3)\to\RR P^\infty$ is $3$-connected), and the agreement along the tube can be
arranged because the tube is contractible. A filling frame is not a
trivialization of $E$ over $Z_0=H_0\cup(-H_1)$: its two handlebody frames
differ along $\Sigma$ by $\rho_\zeta$. In the surface coordinates of the frame
on $\Delta\times\Sigma$, the boundary conditions of the triangle are the three
Lagrangians $(\tau_\zeta L(C_0),L(C_1),L(C_2))$ of
Section~\ref{subsec:triple-data}, where $\tau_\zeta$ is the symplectomorphism
of the surface moduli space induced by $\rho_\zeta$.

\emph{Global frames from the filling.} Since $E|_{Z_0}$ is trivial,
$\zeta\in \Lambda_0+\Lambda_1$ (Proposition~\ref{prop:bundle-normal-form}); write
$\zeta=a_0+a_1$ with $a_i\in \Lambda_i$. Choose $w_1\colon H_1\to\SO(3)$, equal to $1$
near the ball, with lifting class $a_1$. The map $w_1\rho_\zeta^{-1}|_\Sigma$ has
lifting class $a_0\in \Lambda_0$, hence extends, after a homotopy absorbed in a collar,
to $w_0\colon H_0\to\SO(3)$. Then $\phi_0^{\mathrm{gl}}=(w_0\cdot\text{frame on }
H_0)\cup(w_1\cdot\text{frame on }H_1)$ is a trivialization of $E|_{Z_0}$; in the
surface coordinates it acts by $\tau_{a_1}$, which carries
$(\tau_\zeta L(C_0),L(C_1))$ to $(L(C_0),L(C_1))$, because
$\tau_{a_1}\tau_\zeta=\tau_{a_0}$ and $\tau_{a_i}$ preserves $L(C_i)$. For a
perturbed Lagrangian $L_{h_i}(C_i)$ the frame carries $h_i$ to the transported
perturbation $w_i^*h_i$, since trace functions change sign on loops with
$a_i(\ell)=1$; this is a change of primary data, and the presentations
$p_i^\phi$ of Section~\ref{sec:two-handles} carry the transported data. The same
construction at $Z_1$, with $\zeta=a_0'+a_2'$, gives $\phi_1^{\mathrm{gl}}$.

\subsection{Frame changes}\label{subsec:frame-change}

A \emph{frame change} of an object $M$ is a map $u\colon M\to\SO(3)$ equal to
$1$ near $b(B^3)$. It defines a morphism $[u]\colon M\to M$ of $\bB$: the
product $[0,1]\times M$ with the vertical tube and the product bundle, whose
trivialization at the outgoing end is changed by $u$. Frame changes compose
by pointwise multiplication, $[v]\circ[u]=[vu]$, and $[u]$ depends only on
the homotopy class of $u$ relative to the ball.

\begin{lemma}\label{lem:two-handle-frames}
Let $(W_2,E;\sigma_0,\sigma_1)\colon Z_0\to Z_1$ be the two-handle cobordism with
the trivializations $\sigma_i$ induced by Proposition~\ref{prop:outer-bundles},
and put $u_i=\sigma_i(\phi_i^{\mathrm{gl}})^{-1}$. Then in $\bB$
\begin{equation}\label{eq:W2-normal-form}
(W_2,E;\sigma_0,\sigma_1)=[u_1]\circ(W_2,E;\phi_0^{\mathrm{gl}},\phi_1^{\mathrm{gl}})
\circ[u_0]^{-1}.
\end{equation}
If $W_2$ has no two-handles, $E|_{W_2}$ is trivial and the morphism $(W_2,E;\sigma_0,\sigma_1)$ is the
single frame change $[\sigma_1\tilde\sigma_1^{-1}]$, where $\tilde\sigma$ extends
$\sigma_0$ as a product.
\end{lemma}

\begin{proof}
Both sides have the same underlying cobordism and bundle; the
trivializations at the ends agree after composing with the product
cobordisms $[u_i]$, since $u_i=1$ on the ball.
\end{proof}

\section{Necks, families and invariant regularity}\label{sec:families}

This section collects three statements used repeatedly: gluing along
necks with invertible cross-sections, one-parameter families, and
regularity within the class of perturbations invariant under the
lifting group. They are standard in structure, following
\cite[Sections~3--5]{Donaldson-book} and \cite[Section~4]{DFL}; we state them in the form
needed. The points specific to the present setting are the extended weighted spaces at the
mixed ends, with the weight below the stationary gap, and invariance under the lifting
group; the proofs are written out mainly for these.

\subsection{Pieces, ends and hypotheses}

A \emph{piece} is a pure gauge cobordism with cylindrical ends, a DFL
mixed domain \cite[Section~3.1]{DFL}, possibly with several gauge inputs,
or a holomorphic polygon domain with strip-like ends on pairwise
transverse Lagrangians. Ends are \emph{closed} ($\RR_+\times Y$ with a
nondegenerate irreducible critical set), \emph{ordinary} (strip ends at
transverse intersections) or \emph{mixed} (DFL mixed ends, whose stationary solutions form a
Morse--Bott family, treated in extended weighted spaces). Only closed and ordinary ends are joined. We assume:
\begin{enumerate}[label=(H\arabic*)]
\item\label{H:lattice} \emph{Lattice.} For fixed labels every
  finite-energy configuration on a piece, a joined problem or a family
  fibre satisfies $\ind=c+8\kappa$ with $2\kappa\in\ZZ+2\kappa_0$.
\item\label{H:regular} \emph{Component regularity.} Each piece is regular
  in index at most zero (at most one when chain identities are needed),
  and every nonconstant trajectory on a joined end has index at least
  one.
\item\label{H:support} \emph{Supports.} The autonomous primary and end
  perturbations are fixed on the ends and on the corresponding product
  necks. The additional secondary perturbations, the nonautonomous
  cutoff changes, and the parameter derivatives have their output and
  entire loop-dependence regions in a common compact set away from
  the inserted necks. This compact set is fixed after the finite
  selection of secondary directions and before neck lengths vary.
\item\label{H:energy} \emph{Energy.} On solutions, analytic energy equals
  topological energy plus terms bounded independently of neck lengths
  and family parameters.
\item\label{H:framing} \emph{Lifting group.} A finite group $\Phi$ of lifting classes vanishing on every matching surface, satisfying the
  $\phi$-non-integral condition of \cite[Definition~5.1]{KM-unknot} on every
  component of every cut (a surface with $w_2\ne0$ on which $\Phi$ vanishes),
  acts freely on the solution spaces; cuts are connected, so the group of the
  joined problem is the fibre product of the groups of the pieces
  \cite[Proposition~5.4]{KM-unknot}.
\end{enumerate}

\subsection{Neck gluing}

\emph{Function spaces.} Fix an integer $l\ge3$. On a piece, a joined problem
or a family fibre, connections and sections are of class $L^2_{l,\mathrm{loc}}$
and gauge transformations of class $L^2_{l+1,\mathrm{loc}}$. On closed and
ordinary ends, and in particular on every neck, the norms are unweighted:
$L^2_l$ on configurations and tangent vectors, $L^2_{l-1}$ on the target of
the equations. On the mixed ends we use the spaces of
\cite[Definitions~3.8--3.9, (4.5) and (4.8)]{DFL}. The weight is
$e^{\delta\theta}$ in the end coordinate $\theta$ of a mixed end and $1$
elsewhere; a tangent vector is the cut-off pull-back of a tangent vector to
the limiting Lagrangian plus a remainder in $L^2_{l,\delta}$; the target is
$L^2_{l-1,\delta}$. Configurations are compared in the corresponding
distance: the $L^2_l$ distance off the mixed ends, the $L^2_{l,\delta}$
distance of the decaying parts on the mixed ends, and the distance of the
limits. The weight $\delta>0$ is smaller than the constant $\delta_0$ of
\cite[(4.24)]{DFL}, so that the linearized operators are Fredholm, and
smaller than the exponential rate of \cite[Proposition~5.29]{DFL}, so that
solutions lie in these spaces. Moreover
\begin{equation}\label{eq:neck-weight}
\delta<\mu_*,
\end{equation}
where $\mu_*$ is defined as follows. For a mixed end with three-manifold $Y$
and a connection $B$ on $Y$ representing a point of the limiting Lagrangian,
let $\lambda_1(B)$ be the lowest eigenvalue of $d_B^*d_B$ on
$\Omega^0(Y;\ad E)$ with the Neumann condition at $\partial Y$. It is
positive: an element of the kernel is parallel, and $B|_{\partial Y}$ is flat
on a bundle with $w_2\ne0$ on every component of $\partial Y$, so it has no
nonzero parallel section. The eigenvalue depends continuously on $B$ and only
on its gauge class, and the limiting Lagrangians are compact
\cite[paragraph after Proposition~2.16]{DFL}. Hence $\mu_*$, the minimum of
$\lambda_1(B)^{1/2}$ over the mixed ends and the points of their limiting
Lagrangians, is positive. All three conditions on $\delta$ are upper bounds.
Theorem~\ref{thm:neck} is stated and proved in these spaces, and regularity
means surjectivity of the linearized operator between them.

\begin{theorem}[neck gluing]\label{thm:neck}
Let pieces $X_1,\dots,X_k$ be joined along a tree of closed or ordinary
ends with lengths $R=(R_e)$, under
\ref{H:lattice}--\ref{H:framing}, in the function spaces just fixed. There
is $R_0$ such that for $\min_eR_e\ge R_0$:
\begin{enumerate}[label=(\roman*)]
\item the joined problem has no solutions of negative index;
\item gluing is a $\Phi$-equivariant bijection from the union, over
  matched intermediate labels, of products of index-zero component
  spaces onto the index-zero space of the joined problem, whose
  elements are regular;
\item the joined count is the composite of the component counts over
  $\FF$;
\item the same holds with one component an index-one trajectory on an end of the joined domain, modulo translation, and with one component an index $-1$ piece of a
  family together with the parametrized operator.
\end{enumerate}
\end{theorem}

\emph{The chart at a preglued configuration.} Let $p=p_R$ be the preglued
configuration of a matched tuple of component solutions. Its connection lives
on the four-dimensional part $X_R$ of the joined domain and its section on the
two-dimensional part. The boundary of $X_R$ is the union $S$ of the matching
surfaces $\mathbb U_\partial\times\Sigma$ of the pieces, with collars
$[-2,0]_s\times\RR_\theta\times\Sigma$ and outward normal $\partial_s$; we write
$\iota_n$ for the normal component at $S$. Let $B_p$ be the space of
\cite[Section~4.1]{DFL} formed on the joined domain: pairs $w=(\zeta,\nu)$ in
the function spaces above that satisfy the linearized matching,
class-matching and Lagrangian conditions (iii)--(iv) there, with the norm
$|w|$ of \cite[(4.5)]{DFL}. Let $E_p\subset B_p$ be the subspace cut out by
$\iota_n\zeta|_S=0$, and $F_p$ the target \cite[(4.8)]{DFL}. On a piece of
DFL type the chart $P_p$ is that of \cite[Section~4.1]{DFL}. On the joined
domain it is its analogue: $P_p(w)$ is obtained from $p+w$ by three
corrections, all supported in the pieces. On the collar of every component
of $S$, including a component of genus one, where the moduli space of flat
connections on that component is a point and the matching condition is the flatness
of the slices, the slice $A_{s,\theta}+\zeta_{s,\theta}$ is replaced by
\[
A_{s,\theta}+\zeta_{s,\theta}+\rho(s)\bigl(E_{\tau(0,\theta)}(A_{s,\theta},\zeta_{s,\theta})
-A_{s,\theta}-\zeta_{s,\theta}\bigr),
\]
with the maps $E_t$ of \cite[Lemma~4.1 and the paragraph after
Lemma~4.4]{DFL} and DFL's cut-off functions $\rho$ and $\tau$. On the mixed
ends, DFL's finite-dimensional term moves the limits onto fixed families
$B_b$, $B'_{b'}$ of connections representing points of the limiting
Lagrangians. The section is $e_\tau(\nu)$, fibrewise over $u_p$, with the maps
$e_t$ of \cite[Lemma~4.4]{DFL}. The families are
normalized by $\frac{d}{db}B_b|_{b=0}=b$ and by the vanishing of the normal
component of $B_b-B$ on $\partial Y$, and likewise for $B'_{b'}$. The chart is
used with end gauge $1$. Put
\[
\mathcal G_p(w):=\bigl(\mathcal F(P_p(w)),\,d_p^*\zeta\bigr),
\]
where $\mathcal F$ is the perturbed mixed equation, with its Cauchy--Riemann
component regarded as a map into $\RR^N$ for a fixed embedding of the target
symplectic manifold into $\RR^N$. The zeros of $\mathcal G_p$ on $E_p$ are the
points $w$ of the slice \cite[(4.7)]{DFL} for which $P_p(w)$ is a solution.

\begin{remark}[the radius in {\cite[Lemma~4.1]{DFL}}]\label{rem:dfl-radius}
The chart applies $E_t$ to the slices $\zeta_{s,\theta}$, and we use
\cite[Lemma~4.1(iii)]{DFL} with its radius measured in $L^q$ for a fixed
$q>2$. Let $\cB:=L^2_{l-1}(\Sigma;\Lambda^1\otimes\ad F)$ and
$\cB^{(q)}_{<\eps}:=\{c\in\cB:\|c\|_{L^q}<\eps\}$. The statement used is the
following: there is $\eps_q>0$ such that for every connection $\sigma$ of class
$L^2_{l-1}$ on $F$, every $t\in[-1,1]$ and every $s\in[0,1]$, the map
$c\mapsto c+s\bigl(E_t(\sigma,c)-\sigma-c\bigr)$ is a diffeomorphism of
$\cB^{(q)}_{<\eps_q}$ onto an open neighbourhood of $0$ in $\cB$, and, for $\sigma$
of class $L^2_k$ with $k\ge l-1$, restricts to a diffeomorphism at the Sobolev
level $k$. For $s=1$ this is (iii) of \cite[Lemma~4.1]{DFL} with the $L^2$-ball
$\cB_{<\eps}$ replaced by $\cB^{(q)}_{<\eps_q}$; this ball is gauge invariant, so (i) of that
lemma keeps its meaning. The proof is DFL's, made quantitative in $L^q$ as follows. Write
$c=c_0+d_\alpha\xi+{*}d_\alpha\xi'$ with respect to the flat connection $\alpha$ of
\cite[Lemma~4.2]{DFL}. Then $L^q$ elliptic estimates, Kato's inequality and Morrey's
inequality give $\sup_\Sigma|\xi|+\|d_\alpha\xi\|_{L^q}\le C\|c\|_{L^q}$, with $C$ independent of
$\alpha$ by gauge invariance and the compactness of the moduli space. So on a ball of fixed radius
the derivative of the map differs from the identity by at most $\frac12$ in the $L^q$ operator
norm; the map is injective there, and its derivative is injective at every level $k\ge l-1$ and,
being a compact perturbation of an invertible operator, an isomorphism.

We use the statement with an $L^q$-ball, $q>2$; the following example shows that an
$L^2$-ball does not suffice. By (ii) and (iv) of
\cite[Lemma~4.1]{DFL}, $E(\alpha,d_\alpha\xi)$ depends on $\xi$ only through
$\exp\xi$. Let $\psi$ be a cut-off function equal to $1$ on a disc, and let
$e\ne e'$ be sections of $\ad F$ with eigenvalues $\pm i$ that agree where
$\psi<1$, with $e$ parallel on $\operatorname{supp}\psi$. Then $2\pi\psi e$ and $2\pi\psi e'$
have the same exponential, while $d_\alpha(2\pi\psi e)\ne d_\alpha(2\pi\psi e')$.
Choose $e'=\operatorname{Ad}_{\exp\eta}e$ with a nonzero bump $\eta$ supported where $\psi=1$
and with arbitrarily small $\|d\eta\|_{L^2}$; then
$\|d_\alpha(2\pi\psi e')\|_{L^2}\le2\pi\|d\psi\|_{L^2}+C\|d\eta\|_{L^2}$, and a logarithmic cutoff
makes the first term small as well. So the two distinct forms lie in any $L^2$ ball. By the
Sobolev embedding, an $L^2_{l-1}$-ball of fixed radius lies in $\cB^{(q)}_{<\eps_q}$.
\end{remark}

\begin{lemma}[the chart at a preglued configuration]\label{lem:neck-chart}
There are $R_0$, $\kappa_*>0$, $r_*>0$ and $C$, independent of $R$ with
$\min_eR_e\ge R_0$ and of the matched tuple, such that:
\begin{enumerate}[label=(\alph*)]
\item $P_p$ is defined and injective on $\{w\in B_p:|w|<\kappa_*\}$ and takes
  values in the mixed configuration space;
\item on that ball $\mathcal G_p$ is twice continuously differentiable with
  $\|d^2\mathcal G_p\|\le C$, and the same holds in a family for
  $(t,w)\mapsto\bigl(\mathcal F^{\lambda_*+t}(P_p(w)),d_p^*\zeta\bigr)$ with
  $|t|<\kappa_*$, where $\lambda_*$ is the base parameter;
\item every configuration $x$ with $\|x-p\|<r_*$ whose limits on the mixed
  ends lie in the families $B_b$, $B'_{b'}$ is $P_p(w)$ for some $w$ with
  $|w|<\kappa_*$, and $|w|\le C\|x-p\|$.
\end{enumerate}
\end{lemma}

\begin{proof}
The necks lie in closed and ordinary ends, so they meet neither $S$, nor its
collar, nor the mixed ends. Hence the collars, the slices of $p$ on them and
the mixed ends are those of the finitely many component solutions and do not
depend on $R$, and every unit patch of the joined domain is isometric to one of
finitely many models, on which $p$ is bounded in every $C^k$ norm.

(a) By the trace theorem on unit patches,
$\sup_{s,\theta}\|\zeta_{s,\theta}\|_{L^q(\Sigma)}+\sup|\nu|\le C|w|$. For $|w|$
small the slices therefore lie in the ball of Remark~\ref{rem:dfl-radius}, and
$\nu$ lies in a fixed neighbourhood of the zero section on which the maps $e_t$
are fibrewise diffeomorphisms. $P_p(w)$ satisfies the matching,
class-matching and Lagrangian conditions by \cite[Lemma~4.1(iv), Lemma~4.4(iv)
and the paragraph after it]{DFL}. It determines $w$: the end parameters
through its limits, $\nu$ through the injectivity of the $e_t$, and $\zeta$
slice by slice through Remark~\ref{rem:dfl-radius}.

(b) The difference $P_p(w)-p-w$ is the sum of the collar terms $\mathcal N_i(\zeta)$
of the pieces, which depend only on $\zeta$ on the collar of the piece $i$ and
on its slices; the finite-dimensional end term; and a term on the section given
pointwise by the family $e_t$. The last two are smooth, with bounds uniform
over unit patches. Each $\mathcal N_i$ is smooth, with $\mathcal N_i(0)=0$ and
$d\mathcal N_i(0)=0$, on a ball of the space of forms on the collar that are
sums of cut-off pull-backs of profiles in $L^2_l([-2,0]\times\Sigma)$ at
$\theta\to\pm\infty$ and a remainder in $L^2_{l,\delta}$. Indeed, where the slice
has small curvature, DFL's construction expresses $E_t(\sigma,c)$ through the
flat connection $\alpha$ of \cite[Lemma~4.2]{DFL}, the Hodge decomposition of
$c$ with respect to $\alpha$, the exponential of the potential of its exact
part, and a horizontal lift of the path $r\mapsto e_t([\alpha],rc_0)$ of its
harmonic part $c_0$. These ingredients depend smoothly on $\alpha$ in every
Sobolev norm on $\Sigma$, by the inverse function theorem applied at each
Sobolev level to the map $(\alpha,\eta)\mapsto\alpha+{*}d_\alpha\eta$ of
\cite[Lemma~4.2]{DFL} and by Hodge theory for the irreducible flat connection
$\alpha$. Along the collar, $\alpha=\alpha_{s,\theta}$ is smooth in $(s,\theta)$,
with bounded derivatives that converge exponentially along the mixed ends. So
$\mathcal N_i$ is assembled from three kinds of terms: slice-wise linear operators
with bounded $(s,\theta)$-derivatives, which act boundedly on $L^2_l$ of the
collar; entire functions of $\ad\xi$, which are smooth because $L^2_l$ of a
four-manifold of bounded geometry is an algebra for $l\ge3$, also with the
weight on one factor; and a finite-rank term, smooth into every Sobolev space
on $\Sigma$, to which the composition estimates of Moser type apply
\cite[Chapter~13, Proposition~3.9]{Taylor-PDE3}; their proof holds verbatim for
maps from a finite-dimensional space into a Hilbert space. The equation
$\mathcal F$ is quadratic in the connection up to the perturbation terms, and
acts pointwise on the section, with bounds uniform over unit patches; the
fixed autonomous primary and end terms have uniform bounds on unit
patches. The additional secondary terms and nonautonomous parameter
variations have compact dependence regions as in \ref{H:support}. Hence the
second derivative of $\mathcal G_p$ is bounded independently of $R$, the norms
adding over a partition into pieces and necks.

(c) The inverse is constructed slice by slice. The end parameters are read
from the limits of $x$, the families $B_b$ being injective immersions near
$b=0$ by the normalization. The section coordinate $\nu$ is $e_\tau^{-1}$ of the
section of $x$, pointwise. On each collar, $\zeta$ solves
$(\mathrm{id}+\mathcal N_i)(\zeta)=a$, where $a$ is $x-p$ minus the end term; since
$d\mathcal N_i(0)=0$, the inverse function theorem solves this uniquely on a ball
of fixed radius. The resulting $w$ lies in $B_p$. At $s=0$ the slice of $x$ is
$E_t(A_\theta,\zeta_\theta)$ with $A_\theta$ flat, and it is flat; since
$E_t(A_\theta,\cdot)$ maps the small $d_{A_\theta}$-closed forms onto the flat
connections near $A_\theta$ and is injective, $d_{A_\theta}\zeta_\theta=0$. The class
matching follows from $[E_t(\alpha,c)]=e_t([\alpha],[c])$ for $d_\alpha$-closed $c$
and the injectivity of $e_t$, and the Lagrangian condition from the choice of
$e_1$ as the exponential map of a metric in which the Lagrangian is totally
geodesic \cite[Lemma~4.4 and its proof]{DFL}. The constants depend only on the
pieces.
\end{proof}

\begin{proof}[Proof of Theorem~\ref{thm:neck}]
The proof is the standard neck-stretching argument of \cite[Chapters~3--5]{Donaldson-book}.
The two points that differ for the mixed equation are the chart in which the contraction is run
and the gauge slice into which nearby solutions are brought; Section~\ref{subsec:torus-completion}
uses the same steps. Monotonicity is used only through compactness.

\emph{Compactness.} By \ref{H:energy} and \ref{H:lattice} the energy is bounded at fixed
index. By the compactness theorems for the pieces \cite[Theorems~5.1--5.2]{DFL},
\cite[Theorem~3]{DFL-mixed}, \cite[Theorem~5.4]{Donaldson-book} and Gromov compactness, a
sequence with $\min_eR_e\to\infty$ converges to broken configurations on the pieces. The energy
lost in the limit is nonnegative, and the lattice turns it into an index loss $4N$ with $N\ge0$;
with \ref{H:regular} this excludes negative indices and, at index zero, all breaking except at
the joined necks.

\emph{Uniform inverse.} Each neck has an invertible cross-section: on closed necks the
extended Hessian has a spectral gap, the critical points being nondegenerate and irreducible,
and on ordinary necks the intersections are transverse. Hence the linearized operator
$D_p\colon E_p\to F_p$ at the preglued configuration $p$ (on pieces of DFL type, the operator of
\cite[(4.8)]{DFL} with the gauge-fixing component $d_p^*$) has a right inverse $Q_p$ with
$\|Q_p\|\le C_Q$ independent of the neck lengths, by the patching argument of
\cite[Proposition~3.9]{Donaldson-book}, which takes place on the unweighted necks.

\emph{Existence and local uniqueness.} The matching, class-matching and Lagrangian conditions
are nonlinear, and an affine contraction does not preserve them. We therefore run the
contraction for the map $\mathcal G_p$ on $E_p$, in the chart $P_p$ of
Lemma~\ref{lem:neck-chart} with end gauge $1$ (on a piece of DFL type, the chart of
\cite[Section~4.1]{DFL} on the Coulomb slice of \cite[Proposition~4.6]{DFL}). Its
Cauchy--Riemann component is a section along the section of $P_p(w)$, not along $u_p$; let
$\widehat{\mathcal G}_p$ be obtained by projecting it pointwise onto the tangent spaces along
$u_p$ and, in the setting of (iv) with a family, the self-dual component onto the self-dual
forms at the base parameter. These projections are injective nearby, so near $0$ the two maps
have the same zeros. By the normalization $\frac{d}{db}B_b|_{b=0}=b$, the derivative of
$\widehat{\mathcal G}_p$ at $0$ differs from $D_p$ by a zeroth-order term bounded by the error of
$p$, so for $\min_eR_e$ large it has a right inverse $Q'_p$ with $\|Q'_p\|\le2C_Q$ and
$\im Q'_p=\im Q_p$; the normalization of the normal component of $B_b-B$ makes the chart
preserve the normal condition at $S$. The component solutions converge exponentially to
nondegenerate limits, so $\|\widehat{\mathcal G}_p(0)\|_{F_p}\le C_ee^{-\delta_e\min_eR_e}$ for
some $\delta_e>0$, and by Lemma~\ref{lem:neck-chart} the second derivative of
$\widehat{\mathcal G}_p$ is bounded by some $C_2$ on a ball whose radius does not depend on the
neck lengths. Write $\widehat{\mathcal G}_p(z)=\widehat{\mathcal G}_p(0)+D\widehat{\mathcal G}_p(0)z
+\mathcal R_p(z)$, so that $\|\mathcal R_p(z)-\mathcal R_p(z')\|\le C_2(|z|+|z'|)|z-z'|$. The zeros of
$\widehat{\mathcal G}_p$ in $\im Q_p$ are the fixed points of
$z\mapsto-Q'_p\bigl(\widehat{\mathcal G}_p(0)+\mathcal R_p(z)\bigr)$, which is a contraction with constant
at most $\frac12$ on the ball of $\im Q_p$ of radius $\rho_0=(8C_QC_2)^{-1}$ (reduced, if
necessary, to the radius of the chart and to a radius on which the projections are injective). For $\min_eR_e$ large it maps the ball of radius
$4C_QC_ee^{-\delta_e\min_eR_e}$ into itself, so the contraction mapping principle gives a zero
$z_R\in\im Q_p$ with $|z_R|\le4C_QC_ee^{-\delta_e\min_eR_e}$, the only zero in the ball of
radius $\rho_0$; $P_p(z_R)$ is the glued solution, and at index zero $\im Q_p=E_p$. The derivative
of $\widehat{\mathcal G}_p$ at $z_R$ differs from that at $0$ by at most $C_2|z_R|$, so it is onto;
it is the linearization at $P_p(z_R)$ in the chart followed by the projections, which are
isomorphisms there, and the slice is transverse to the gauge orbits, so $P_p(z_R)$ is regular.
In the setting of (iv) with a family the contraction is run in $(t,z)\in\RR\oplus E_p$ with the
parametrized operator $(t,z)\mapsto D_pz+t\,\partial_\lambda\mathcal F$.

\emph{Surjectivity.} Let $y_n$ be solutions with $\min_eR_{e,n}\to\infty$ in a fixed class of
index zero. By the compactness step, after passing to a subsequence and applying gauge
transformations, $y_n$ converges on compact subsets of the pieces to a matched tuple $c$ of
index-zero component solutions, with no energy lost and no trajectory breaking off, so its
energy on the necks and on the ends outside a fixed compact set tends to zero. Exponential decay for solutions of small energy
then controls $y_n$ there, with constants independent of the neck lengths: on closed necks and
gauge ends by \cite[Propositions~4.3--4.4 and Section~5.5]{Donaldson-book}; on ordinary necks
and ends by \cite[Proposition~2.3.1]{Simcevic}, after the Hamiltonian term is removed by its
flow, with all derivatives controlled by \cite[Lemma~C.3]{RS-strips} and the linear equation for
the time derivative; and on the mixed ends by \cite[Proposition~5.29]{DFL}. Patching gauge
transformations over bands on which the limits are irreducible gives $g_n$ with
$\|g_n^*y_n-p_n\|\to0$, where $p_n$ is the preglued configuration of $c$.

To bring $g_n^*y_n$ into the slice at $p_n$ with control independent of the neck lengths, let
$\Delta_p\xi=\bigl(d_p^*d_p\xi,\ \iota_nd_p\xi|_S\bigr)$ act on infinitesimal gauge
transformations of class $L^2_{l+1}$ off the mixed ends that are asymptotic on each mixed end to
the pull-back of a section of $\ad E$ over $Y$, with a remainder in $L^2_{l+1,\delta}$. Its
asymptotic operators are positive: on closed necks and gauge ends because the critical points
are irreducible, and on the mixed ends because $\lambda_1(B)\ge\mu_*^2$. Consequently
$\|\xi\|_{L^2}\le C_P\|d_p\xi\|_{L^2}$ for every $\xi\in L^2_1(X_R)$, with no boundary condition
and with $C_P$ independent of the neck lengths. If this failed, a sequence with
$\|\xi_n\|_{L^2}=1$ and $\|d_{p_n}\xi_n\|_{L^2}\to0$ would, by the Poincar\'e inequalities on the
cross-sections (on $Y$ with no boundary condition, by the variational characterization of
$\lambda_1(B)$), have $L^2$ norm tending to zero on the necks and ends; on the compact parts of
the pieces it would converge to a parallel section vanishing on an open set, hence to zero,
every component of a piece having an end. In addition to \eqref{eq:neck-weight}, take
$\delta C_P\le\frac12$, and let $\tau$ be a smoothing of the end coordinate of the mixed ends,
extended by zero, with $|d\tau|\le1$ and depending only on $\theta$ near $S$. For
$U=e^{\delta\tau}\xi$ the weak Neumann form becomes
$\langle d_pU-\delta\,d\tau\,U,\ d_pV+\delta\,d\tau\,V\rangle_{L^2}$, whose diagonal
$\|d_pU\|^2-\delta^2\|d\tau\,U\|^2\ge\frac34\|d_pU\|^2$ is coercive by the Poincar\'e
inequality. Since $\partial_n\tau=0$, the Lax--Milgram lemma and local Neumann regularity on unit
patches invert $\Delta_p$ on the weighted spaces and, with the Neumann problem on $Y$ for the
asymptotic values, on the spaces with asymptotic values, with bounds independent of the neck
lengths. After the mixed-end limits are brought into the families $B_b$, $B'_{b'}$ by a gauge
transformation on $Y$, $Y'$ (Coulomb gauge with the Neumann condition at $B$ and
\cite[Proposition~2.10]{DFL}), the inverse function theorem for the action of the decaying gauge
group in the chart, whose linearization $(\xi,w)\mapsto w-d_p\xi$ is inverted by means of
$\Delta_p^{-1}$, writes every configuration within a distance $r_C$ of $p$ as $g^*P_p(w)$, with
$w$ in the slice \cite[(4.7)]{DFL} and $|w|$ bounded by a fixed multiple of that distance; by
Lemma~\ref{lem:neck-chart}, $r_C$ and the multiple do not depend on the neck lengths. This is
\cite[Proposition~4.6]{DFL} at $p$ with uniform radii.

Hence $y_n$ is gauge equivalent to $P_{p_n}(w_n)$, with $w_n$ a zero of
$\widehat{\mathcal G}_{p_n}$ in $E_{p_n}$ and $|w_n|\to0$. At index zero $|w_n|<\rho_0$ for $n$
large, and by uniqueness $P_{p_n}(w_n)$ is the glued solution. At augmented index zero in a
family the same argument applies to $(t,w)$, with the compactness step of
Theorem~\ref{thm:family} and the slice at the base parameter. With one index-one trajectory, the
zeros of $\widehat{\mathcal G}_p$ in a ball of radius independent of the neck length form an arc
through the glued point, a graph over $\ker D_p$ with values in $\im Q_p$, along which the glued
point moves with speed bounded below as the gluing parameter varies; so every nearby solution is
glued. The gluing map is continuous, locally injective, and leaves every compact subset of the
one-dimensional moduli space as the parameter tends to infinity; its image lies in a noncompact
component, an interval, on which such a map is monotone, hence injective.

\emph{Injectivity.} Distinct component tuples give distinct glued solutions for large neck
lengths: if the glued solutions of $c$ and $c'$ were gauge equivalent along a sequence of neck
lengths tending to infinity, both would converge in $L^2_l$ on compact subsets of the pieces, to
$c_i$ and $c'_i$; the relating gauge transformations would be bounded in $L^2_{l+1}$ there, a
diagonal subsequence would converge, and its limit would relate $c_i$ to $c'_i$
\cite[Section~4.4]{Donaldson-book}. Equivariance and the fibre-product description of the framed
quotient follow from \ref{H:framing}.
\end{proof}

\subsection{Families}

\begin{theorem}[families]\label{thm:family}
Let $F_\lambda$, $\lambda\in[0,1]$, be a smooth family of data on one
domain, constant near $\lambda=0,1$, with ends and critical data
independent of $\lambda$, $\partial_\lambda F$ supported in a compact
region, all scales in a compact positive range, and $\Phi$-invariant.
Assume both endpoint fibres are regular in index at most zero. After adding a generic
small $\Phi$-invariant perturbation from a finite-dimensional space, vanishing near the
endpoints, with supports and loop regions in a compact set
(Theorem~\ref{thm:regularity}(iii)):
\begin{enumerate}[label=(\roman*)]
\item there are no solutions of fibre index less than $-1$; the index
  $-1$ solutions occur at finitely many parameters, and at each the
  augmented operator $(\xi,t)\mapsto D\xi+t\,\partial_\lambda F$ is
  invertible;
\item for each such solution and each index-one trajectory at an
  closed or ordinary end there is exactly one boundary branch of the
  parametrized index-zero space;
\item the compactified parametrized index-zero space is a compact
  one-manifold whose boundary consists of the two endpoint index-zero
  spaces and the branches of (ii); breaking at the mixed ends does not occur,
  since a nonconstant mixed-end trajectory has $\kappa\in\frac12\ZZ_{>0}$ and
  hence index at least four;
\item consequently the two endpoint counts $P_0,P_1$ satisfy
  $P_0+P_1=dK+Kd$ over $\FF$, where $K$ counts the index $-1$ solutions.
\end{enumerate}
\end{theorem}

\begin{proof}
\emph{Compactness in the family.} The compactness theorem
\cite[Theorem~3]{DFL-mixed} is stated for the fixed collar geometry
$ds^2+d\theta^2+g_\Sigma$ of \cite[(1.1)]{DFL-mixed}, and only its interior
estimate \cite[Lemma~4.50]{DFL-mixed} is formulated for a $C^k$-compact
family of metrics. We explain why it holds uniformly over the family. Near
the seams, the metrics of the family are of the form
$e^{2f_\lambda}(ds^2+d\theta^2)+c_\lambda^2g_\Sigma$ with $f_\lambda$ in a $C^k$-compact
set and scales $c_\lambda$ in a compact positive range. In every family used in
this paper (the path of Lemma~\ref{lem:family-path}, the stages in the proof
of Theorem~\ref{thm:folded-triangle}, and the families of
Lemma~\ref{lem:vertex-identification}, a class closed under convex combinations) one has $e^{-f_\lambda}c_\lambda\equiv1$, and we
use the statement only in that case: the conformal change below then gives
exactly the product metric \cite[(1.1)]{DFL-mixed}, to which
\cite[Theorem~4.2 and Lemma~4.25]{DFL-mixed} apply verbatim. For completeness we
indicate the argument for general warped families. The four-dimensional
conformal change by $e^{-2f_\lambda}$ turns this into the warped product
$ds^2+d\theta^2+(e^{-f_\lambda}c_\lambda)^2g_\Sigma$, whose fibre scale has
logarithmic derivatives bounded uniformly in $\lambda$. On a unit patch, after
rescaling $(s,\theta)$ by the fibre scale at its centre, the metric is a
$C^k$-small perturbation of the fixed product $ds^2+d\theta^2+g_\Sigma$, with
smallness controlled by the patch size. The anti-self-duality equation and the
curvature energy are invariant under four-dimensional conformal changes, the
Cauchy--Riemann equation and the energy of the section depend only on the
conformal structure of the two-dimensional factor, the matching condition
does not involve the metric, and the holonomy perturbations vanish on the
collars. The mean-value inequality \cite[Theorem~4.2]{DFL-mixed} and the seam
derivative bounds \cite[Lemma~4.25]{DFL-mixed} are proved by elliptic
estimates whose constants depend continuously on the metric in $C^k$, so they
hold on these patches with uniform constants, and so does seam regularity. Away from the seams, Uhlenbeck compactness and
\cite[Lemma~4.50]{DFL-mixed} apply to compact families of metrics, and the
perturbations depend smoothly on $\lambda$.

Consequently a sequence of solutions $U_\nu$ of fibre index at most one and
bounded energy at parameters $\lambda_\nu\to\lambda$ has a subsequence converging,
after gauge transformations, in $C^\infty$ on compact subsets of the
complement of a finite set $Z$ and modulo breaking at the ends, to a
solution at $\lambda$, and exactly one of the following holds: $Z=\emptyset$ and
no energy is lost, in which case the convergence is $C^\infty$ on compact sets
up to breaking; or $Z\neq\emptyset$ and at least one quantum of energy is lost
at $Z$, in which case convergence fails in $C^\infty$ near $Z$ and holds only
in $L^p_1$ on compact subsets of the complement. In the second case
\ref{H:lattice} converts the energy loss into an index loss of at least four,
so the limit has fibre index at most $-3$, which the regularity of the
family (\ref{H:regular} and (i)) excludes. Hence in index at most one only
the first alternative occurs. Regularity of the parametrized problem is Theorem~\ref{thm:regularity}
applied with the parameter. Gluing at an index $-1$ piece uses the
augmented operator \cite[Section~5.3, (5.6)]{Donaldson-book}. The
remaining statements are the standard boundary count.
\end{proof}

\subsection{Invariant regularity}

We use holonomy perturbations in the sense of
\cite[Section~6.2, (6.13)--(6.14)]{DFL}: finite sums of terms
$\pi_{\mathfrak{su}(2)}\mathrm{Hol}_\gamma(\tilde A)\otimes\omega$, where
$\gamma$ is a family of loops, $\tilde A$ a local unitary lift and
$\omega$ a compactly supported form, with every loop having zero value
under every element of $\Phi$.

\begin{lemma}\label{lem:invariant-span}
Let $\Sigma$ be a closed oriented surface of genus at least one,
$U\to\Sigma$ a rank-two Hermitian bundle of odd degree with fixed
determinant connection, and $\tilde B$ a unitary connection inducing it
whose associated $\mathrm{PU}(2)$ connection is flat. Let $\mathcal H\subset
\U(2)$ be its based holonomy group. Then $\mathcal H$ acts irreducibly on
$\CC^2$, $\mathcal H$ contains $\U(1)\cdot\mathrm{Id}$, and the real span of
$\mathcal H$ is $M_2(\CC)$. The same holds for holonomies defined only up
to sign.
\end{lemma}

\begin{proof}
An invariant line would carry half the curvature of the determinant,
hence have half-integer degree. Small discs where the central part of
the curvature is nonzero give a continuous arc of scalar holonomies,
which generates $\U(1)$. By Burnside's theorem the complex span of
$\mathcal H$ is $M_2(\CC)$, and since $i\mathcal H\subset\mathcal H$ its
real span is the same. A sign change does not change a real linear span.
\end{proof}

For a holonomy perturbation term $\pi\,\mathrm{Hol}_\gamma\otimes\omega$ call
$\operatorname{supp}\omega$ its \emph{support} and the union of the images of its loops its
\emph{loop region}. On each gauge component choose an open set $O$ containing the supports
of all fixed perturbation derivatives; it may extend down original ends. Choose an open
set $\Lambda$ such that every component of $O$ is joined through $\Lambda$ to an original
end whose limiting connection is flat on a collar of an odd surface, and on which the
lifting classes of loops vanish. On a joined torus domain take the union of these sets on
the two sides; the loops need not cross the inserted torus bridge. The compact set
containing the supports and loop regions of the finitely many perturbations eventually
chosen is fixed afterwards.

\begin{theorem}[invariant regularity]\label{thm:regularity}
\begin{enumerate}[label=(\roman*)]
\item Let $Y$ be a closed connected three-manifold containing a surface
  $\Sigma$ with $w_2(E)|_\Sigma\ne0$, and $g$ an automorphism of the
  $\SO(3)$ bundle $E$ whose lifting class vanishes on $\Sigma$. If
  $g$ fixes a connection whose restriction to $\Sigma$ is flat, then
  $g=1$. In particular $\Phi$ acts freely on critical points and on
  finite-energy solutions.
\item Let $D$ be the linearized operator at a solution of a gauge or mixed problem,
  with $O$ and $\Lambda$ as above. Then $D$, augmented by the derivatives of the
  $\Phi$-invariant holonomy perturbations with support in $O$ and loop region in
  $O\cup\Lambda$, is onto. On purely symplectic components the same holds with the usual
  almost complex and Hamiltonian variations in place of holonomy perturbations.
\item Under the compactness and gluing hypotheses of
  Theorems~\ref{thm:neck} and~\ref{thm:family}, generic small $\Phi$-invariant
  perturbations from a finite-dimensional space of such directions, with supports and loop
  regions in a compact set fixed before the neck lengths are varied, make all spaces of
  index at most one regular, relative to strata that are already regular.
\end{enumerate}
\end{theorem}

\begin{proof}
(i) A fixed automorphism commutes with the holonomy of the connection
on $\Sigma$, which is irreducible by Lemma~\ref{lem:invariant-span}, the
bundle being odd on $\Sigma$; so $g$ is the identity on $\Sigma$, and hence everywhere by unique
continuation of parallel sections along paths, $Y$ being connected. The
hypotheses cannot be dropped, as the example of $T^3$ with holonomy the
Klein four-group shows.

(ii) Write the full gauge-fixed derivative at a solution as
\[
D(a,v)=\bigl(d_A^*a,(d_A^++K_A)a,D_uv\bigr),
\]
where $K_A$ includes all primary and secondary holonomy derivatives.
Let $(\mu,\xi,z)$ annihilate the range of this operator and of all
perturbation directions as in (ii). We show $\mu=0$, then $\xi=0$, then $z=0$. A Neumann gauge variation $(d_Af,0)$ is an admissible
source; its surface harmonic class is zero. Gauge covariance at a
solution gives $(d_A^++K_A)d_Af=0$, so
$\langle\mu,d_A^*d_Af\rangle=0$.
The scalar Neumann operator is onto at each such solution. Indeed
its end operators have positive scalar gaps; on a completed torus
end use the odd-fibre gap. The compact-remainder argument therefore
makes it Fredholm at weight zero. Integration by parts identifies
its kernel with decaying parallel sections, which vanish by
irreducibility. Self-adjointness gives index zero, and sufficiently
small weights preserve that index and introduce no kernel. Hence $\mu=0$.

At the limiting odd reference surface, Lemma~\ref{lem:invariant-span}
provides finitely many holonomies spanning the adjoint fibre over
$\mathbb R$. Move their loops slightly into the unperturbed flat
collar and then to a sufficiently far, finite slice of the original
end. Convergence to the limiting connection preserves this span.
Conjugate the loops by a common path to any ball in $O$. The path and
its reverse cancel in every lifting class, so these loops still have
zero lifting class. Thicken them smoothly to
holonomy families. Their self-dual output forms can be chosen
arbitrarily in the ball; the real coefficients expressing a target
vector in the spanning holonomies are computed at this solution
and fixed as part of the perturbation. Pairing with all such
directions therefore gives $\xi|_O=0$.

All fixed nonlocal derivatives have output in $O$. Consequently,
for every test variation $a$,
$\langle K_A^*\xi,a\rangle=\langle\xi,K_Aa\rangle=0$.
Thus the interior adjoint equation is local:
$d_A^{+,*}\xi=0$. Applying $d_A^+$ gives an elliptic second-order
equation with scalar principal symbol. Unique continuation
\cite{Aronszajn} implies $\xi=0$ on each connected gauge component.

On an open matching arc choose an arbitrary section trace, represent
its components by harmonic forms on the matching surfaces, and
extend them independently into the gauge collars with zero normal
component. These are admissible coupled variations. Since both
gauge cokernel components vanish, their Green pairing forces the
section adjoint to have zero boundary trace on that arc. Mixed
adjoint regularity supplies this trace. Extension by zero across a
smaller arc, followed by unique continuation for the real
Cauchy--Riemann equation, gives $z=0$. Purely symplectic components
are regularized separately by their usual $J$ and Hamiltonian
variations \cite[Theorem~3.9]{WW-relative} and use no holonomy
perturbations. The pure gauge pieces have framed reference ends
with the same unperturbed odd collar. For the closed gauge factors
in the excision theorem we use the regularity results of
\cite[Section~5.3]{KM-unknot}.

(iii) At a fixed solution finitely many directions suffice because the
cokernel is finite dimensional. At each index, remove neighbourhoods
of the already regular broken strata and of fixed regular parameter
faces. Compactness of the remaining locus gives a finite collection
of directions by a finite cover. Parameter cutoffs implement relative
choices. Sard--Smale then gives a generic small regular choice.
Only now take a compact set containing the output supports and the
entire loop regions of these directions, and place the original-end
truncations farther out. This set is independent of subsequently
inserted lengths. Each direction is invariant, so no
averaging or loss of span occurs.
\end{proof}

\section{The collapse theorem: geometry, index and compactness}\label{sec:collapse}

This section and the next two prove Theorem~\ref{thm:B}. The present
section fixes the geometry, states the theorem, proves the index-energy
law at every scale and proves compactness with the loss formula.
Section~\ref{sec:collapse-linear} contains the uniform linear estimates, and Section~\ref{sec:collapse-coverage} the nonlinear gluing and the surjectivity of the gluing map.

\subsection{Geometric hypotheses}\label{subsec:clean-data}

Let $(C,g_C,\mathbb E)$ satisfy Definition~\ref{def:clean-datum},
and write $n=\dim L(C)$ and $P=L(C)$; $\Sigma$ denotes the union of the boundary
components matched to the section and $\Msf$ their moduli space.

The two instances used in this paper are the torus compression body
$C=C_1$ with $\partial C=\Sigma_h\sqcup T^2$, where the torus has a
one-point moduli space, $L(C)\subset\Msf_h$ and $n=3h-3$; and the
compression body $C_{h,g}$ from $\Sigma_h$ to $\Sigma_g$ used in the
folded geometry of Section~\ref{sec:torus}, with
$L(C_{h,g})\subset\Msf_h^-\times\Msf_g$ and $n=3(h+g-2)$. Conditions
(C1)--(C3) hold for both: irreducibility, forced by the odd bundle on a
boundary component, gives $H^0=0$; since $H^0(\partial C;\ad)=0$, duality
gives $H^2(C;\ad)\cong H^1(C,\partial C;\ad)^*$, which vanishes exactly when
$H^1(C;\ad)\to H^1(\partial C;\ad)$ is injective. Injectivity holds because
the inclusion of the large boundary surface $\partial_+C$ is $\pi_1$-surjective:
first cohomology with local coefficients is group cohomology of $\pi_1$, and
inflation $H^1(\pi_1C;\ad)\to H^1(\pi_1\partial_+C;\ad)$ is injective for a
surjective homomorphism. The same surjectivity makes the image an embedded
Lagrangian \cite[Proposition~2.21]{DFL}.

For $v\in T_qL(C)$ let $h_v$ be the harmonic $\ad B_q$-valued one-form on
$C$ with vanishing normal component and boundary class $v$; it is unique
by (C3).

\begin{lemma}[distance to the Lagrangian]\label{lem:distance}
There is $c_C$ such that every smooth connection $B$ on $\mathbb E$ with
flat boundary restriction satisfies
$\dist([B|_{\partial C}],L(C))\le c_C\|F_B\|_{L^2(C)}$.
\end{lemma}

\begin{proof}
Since the target is compact, it suffices to treat small curvature. If
$\|F_{B_i}\|_{L^2}\to0$, weak Uhlenbeck compactness on the compact
three-manifold with boundary \cite[Theorem~A]{Wehrheim-Uhlenbeck} gives gauges in which $B_i\rightharpoonup
B_\infty$ in $W^{1,2}$ and $B_i\to B_\infty$ in $L^4$, with $B_\infty$ flat;
the relative Coulomb slice, with vanishing normal component, and the
normal Hodge estimate
$\|a\|_{W^{1,2}}\le C(\|d_{B_\infty}a\|+\|d^*_{B_\infty}a\|+\|a\|)$ upgrade
this to strong $W^{1,2}$ convergence. In a flat chart about $B_q$, fix the
gauge so that $a=B-B_q$ is in Coulomb gauge with vanishing normal component
and $L^2$-orthogonal to the harmonic fields, adjusting $q$; this is possible
by the implicit function theorem because the gauge and flat directions are
transverse. Then $\|a\|_{W^{1,2}}\le C\|d_{B_q}a\|\le C(\|F_B\|+\|a\|_{L^4}^2)$,
so $\|a\|_{W^{1,2}}\le C\|F_B\|$, and the trace inequality, read in the
$H^{1/2}$ surface Coulomb chart of the flat connections (proof of
Lemma~\ref{lem:seam-regularity}), bounds the distance of the boundary class
from $q\in L(C)$.
\end{proof}

\subsection{The collapsing family $X_\eps$}\label{subsec:insertion}

\begin{definition}\label{def:insertion}
Let $X_0$ be a domain for the mixed equation of \cite{DFL}: a
four-manifold with a bundle, a surface $S_0$ with boundary arcs carrying
Lagrangians (the surface on which the holomorphic part $u$ of a configuration,
called the \emph{section} as in \cite{DFL}, is defined), matching seams, and ends of
three kinds: gauge-theoretic cylindrical ends, strip ends, and mixed ends of
Morse--Bott type. Assume that $S_0$ has a boundary arc $I$ carrying $L(C)$
between two boundary punctures $z_\pm$, and that at $z_\pm$ the domain has
the standard mixed end of $C$: in polar coordinates $z=z_++xe^{i\theta}$ at
$z_+$ and $z=z_--xe^{-i\theta}$ at $z_-$, $0<x<x_0$, $0\le\theta\le\pi$, with
$x_0<|I|/2$ so that the two ends do not overlap, the gauge piece $\{0<x<x_0\}\times C$ is
matched to the section along the ray $\theta=0$, and the ray $\theta=\pi$
lies in $I$. The ends of $X_0$ other than the two mixed ends at $z_\pm$ are its
\emph{original ends}. For $\eps\in(0,x_0/4)$ the domain $X_\eps$ is obtained as
follows.
\begin{enumerate}[label=(\roman*)]
\item The section surface of $X_\eps$ is $S_0$ with $I$ and $z_\pm$ restored:
  near $I$ it is the upper half-plane chart with coordinate $z$, and $I$
  becomes part of a single matching seam $(-2,2)\subset\RR$. Its complex
  structure $j_\eps$ is the push-forward under $F_\eps$ of the conformal
  structure of the collar, extended scale-covariantly so that it is
  self-similar under $z\mapsto(z-z_\pm)/\eps$; it agrees with that of $S_0$
  outside two half-disks of radius $c_{\rm core}\eps$ about $z_\pm$, for a constant
  $c_{\rm core}\ge2$ fixed with the insertion.
\item The gauge piece over the seam is $(-2,2)\times C$ with the metric
  $g_\eps=dt^2+\lambda_\eps(t)^2g_C$, where $\lambda_\eps=\eps\,\ell(x/\eps)$
  near the endpoints ($x$ the distance to the nearer endpoint, positive
  outside $I$), $\lambda_\eps=\eps$ on $I$, and $\lambda_\eps=1$ near
  $t=\pm2$; here $\ell\colon\RR\to[1,\infty)$ is smooth with $\ell(s)=1$ for
  $s\le1$ and $\ell(s)=s$ for $s\ge2$. For $x\ge2\eps$, where
  $\lambda_\eps=x$, this piece is the exterior gauge piece of $X_0$ at
  $z_\pm$.
\item The matching seam is $t\mapsto t$ along $(-2,2)$, the boundary of the
  section collar embedded by
  $F_\eps(t,\theta)=t+\lambda_\eps(t)\bigl[\sigma b_\eps(t)(\cos\theta-1)
  +i\bigl((1-b_\eps(t))\theta+b_\eps(t)\sin\theta\bigr)\bigr]$ for
  $0\le\theta<\theta_0$, a fixed small angle, with $\sigma=\pm1$ at the right
  and left endpoints and $b_\eps$ a cutoff equal to $0$ for $x\le\eps$ and
  $1$ for $x\ge2\eps$. On $I$ it is $t+i\eps\theta$; for $x\ge2\eps$ it is
  the polar collar $z_\pm\pm xe^{\pm i\theta}$. For $\theta$ in this range
  $F_\eps$ is an embedding.
\item Outside a fixed neighbourhood $N_I$ of $I$ and its endpoints, $X_\eps$
  coincides with $X_0$; $N_I$ is taken to contain the patches $\mathcal H_{r_0}(b)$ of
  Lemma~\ref{lem:seam-regularity}, for the radius $r_0$ fixed there, about every seam point $b$
  of $I$ and of the endpoint necks $\{x\le x_0\}$.
\end{enumerate}
\end{definition}

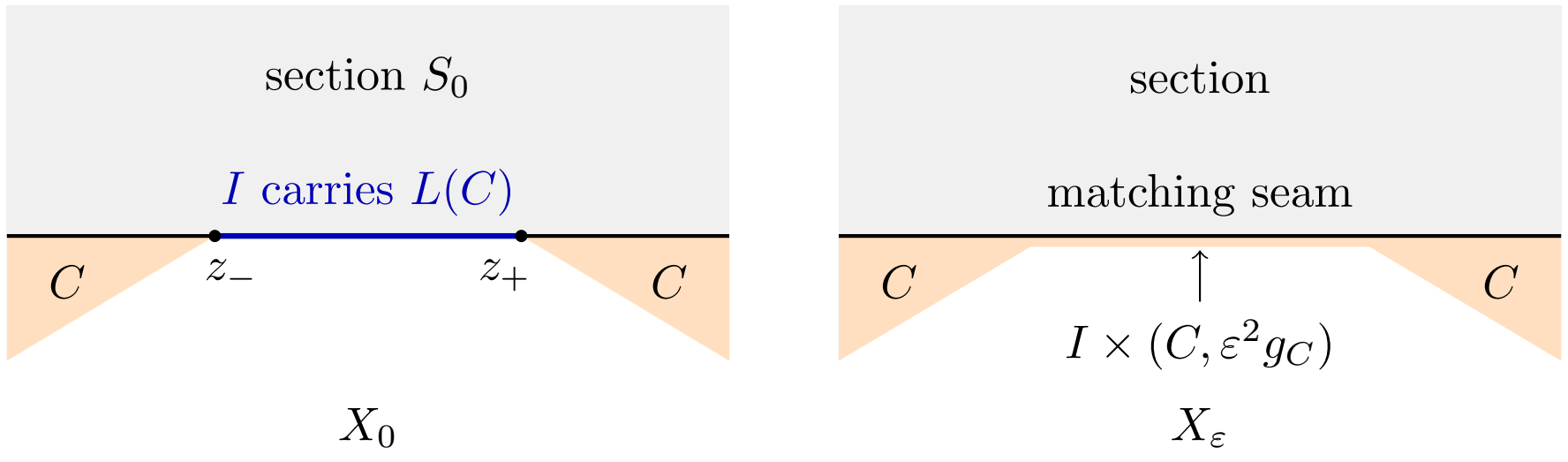
\begin{figure}[ht]
\centering
\begin{tikzpicture}[scale=0.82,font=\small]
\begin{scope}
\fill[gray!12] (-3.3,0) rectangle (3.3,2.1);
\node at (0,1.45) {section $S_0$};
\fill[orange!25] (-3.3,0) -- (-1.4,0) -- (-3.3,-1.14) -- cycle;
\fill[orange!25] (3.3,0) -- (1.4,0) -- (3.3,-1.14) -- cycle;
\draw[thick] (-3.3,0) -- (-1.4,0);
\draw[thick] (1.4,0) -- (3.3,0);
\draw[very thick,blue!70!black] (-1.4,0) -- (1.4,0);
\node[blue!70!black,above] at (0,0.02) {$I$ carries $L(C)$};
\fill (-1.4,0) circle (1.6pt); \fill (1.4,0) circle (1.6pt);
\node[below] at (-1.25,-0.05) {$z_-$}; \node[below] at (1.25,-0.05) {$z_+$};
\node at (-2.75,-0.42) {$C$};
\node at (2.75,-0.42) {$C$};
\node at (0,-1.75) {$X_0$};
\end{scope}
\begin{scope}[xshift=7.6cm]
\fill[gray!12] (-3.3,0) rectangle (3.3,2.1);
\node at (0,1.45) {section};
\fill[orange!25] (-3.3,0) -- (3.3,0) -- (3.3,-1.14) -- (1.55,-0.1) -- (-1.55,-0.1) -- (-3.3,-1.14) -- cycle;
\draw[thick] (-3.3,0) -- (3.3,0);
\node[above] at (0,0.02) {matching seam};
\draw[<-] (0,-0.13) -- (0,-0.6) node[below] {$I\times(C,\eps^2g_C)$};
\node at (-2.75,-0.42) {$C$};
\node at (2.75,-0.42) {$C$};
\node at (0,-1.75) {$X_\eps$};
\end{scope}
\end{tikzpicture}
\caption{The collapsing family of Definition~\ref{def:insertion} near the arc $I$. The
shaded region below the seam is the gauge piece, drawn with thickness equal to
the scale of the metric on $C$: it is proportional to the distance to $z_\pm$ outside
$I$, and equal to $\eps$ over $I$ on $X_\eps$. On $X_0$ (left) the arc $I$ carries the
Lagrangian boundary condition $L(C)$, and after the conformal change of
Lemma~\ref{lem:conformal-models} the two points $z_\pm$ become mixed ends. On
$X_\eps$ (right) the whole segment is a matching seam, and the thin part is the
compression rectangle.}
\label{fig:collapse}
\end{figure}

\begin{lemma}\label{lem:conformal-models}
\begin{enumerate}[label=(\roman*)]
\item On the collar, $\lambda_\eps^{-2}g_\eps=ds^2+g_C$ with
  $ds=dt/\lambda_\eps$, and the section metric is conformal to
  $ds^2+d\theta^2$. The matching conditions on $\partial C$ therefore have
  the product form of \cite[Section~3.1]{DFL} at every scale, with the fixed
  surface metric.
\item Near each endpoint, for $c_{\rm core}\eps<x<x_0$, $s=-\log x$ and polar
  coordinates on $S_0$, the domain is conformal to a cylinder of length
  $\log(x_0/2\eps)$ with slice
  \[
  (C,g_C)\ \sqcup\ [0,\pi]\ \sqcup\ (-C,{\mathfrak s}^2g_C),\qquad
  {\mathfrak s}=\eps e^{s}=\eps/x,
  \]
  in which the first copy is the exterior gauge piece at the point
  $z_\pm\pm x$, matched at $\theta=0$, and the last is the rectangle at the
  point $z_\pm\mp x\in I$, matched at $\theta=\pi$. The two copies carry opposite
  orientations, so that the two seam terms of Lemma~\ref{lem:green} cancel: with the
  time chosen so that (time, $\theta$) is positively oriented, namely $\rho=L-s$ at $z_+$
  and $s$ at $z_-$, where $\theta$ runs clockwise, the copy at $\theta=\pi$ is the reversed
  one. At $z_-$ the tangential operator is conjugate to $-A_{\mathfrak s}$, and every
  estimate below is invariant under this sign change.
\item Along $I$, at distance $r\ge2\eps$ from the nearer endpoint, the
  conformal change by $r^{-2}$ gives the compression the metric
  $dT^2+(\eps/r)^2g_C$ with $dT=dt/r$; so ${\mathfrak s}=\eps/r$ along $I$.
\end{enumerate}
The anti-self-duality equation, the Cauchy--Riemann equation and their
energies are conformally invariant.
\end{lemma}

\begin{proof}
(i) and (iii) are direct computations. For (ii), $dx^2+x^2g_C=
x^2(ds^2+g_C)$, and the collar embedding is the polar map. Conformal
invariance holds in dimensions four and two respectively.
\end{proof}

The gauge-fixing (Coulomb) condition is not conformally invariant. In
the conformal coordinates of Lemma~\ref{lem:conformal-models} it
acquires a zeroth-order term $2\varsigma\,\Pi_{\rm f}$ on the collapsing component,
where $\varsigma=\partial_s\log{\mathfrak s}$ and $\Pi_{\rm f}$ is the projection onto the
scalar component.

\begin{definition}\label{def:admissible-config}
A \emph{configuration} on $X_\eps$ ($\eps\ge0$) with given
asymptotic labels is a pair $(A,u)$ of a connection on the bundle and a
section map, of class $L^2_{k,\mathrm{loc}}$ for a fixed large $k$, such
that: the restriction of $A$ to each matching surface is flat and
represents the value of $u$ on the seam; $u$ takes values in the given
Lagrangians on the Lagrangian arcs (in $L(C)$ on $I$ when $\eps=0$); and at
every end $(A,u)$ converges exponentially, at the fixed weights below the
gaps, to the labelled critical data, with free limits in the critical
manifolds at mixed ends. Configurations are considered up to the
determinant-one gauge group enlarged by the lifting group, and carry
their relative homotopy class and bundle marking.
\end{definition}

\subsection{Data and hypotheses}\label{subsec:collapse-hypotheses}

\begin{hypothesis}\label{hyp:collapse}
The data on $X_\eps$, $\eps\ge0$, satisfy:
\begin{enumerate}[label=(X\arabic*)]
\item\label{hyp:X1} \emph{Structures near the collapse.} On $N_I$ the
  primary perturbation of $C$ vanishes, the section almost complex
  structure is the complex structure $*_\Sigma$ of Convention~\ref{conv:orientation}
  with no inhomogeneous term, and all metrics are products near the
  seams.
\item\label{hyp:X2} \emph{Perturbations.} The end perturbations are the fixed
  autonomous perturbations of the original gauge ends, as in
  Proposition~\ref{prop:input-functions}, vanishing near the critical set and switched on by
  monotone cutoffs along the end coordinates; their outputs (the supports of their forms),
  their dependence regions (the unions of their loops) and their cutoff bands lie in the
  original ends, outside $N_I$. All secondary perturbations, including the nonautonomous
  holonomy terms $\chi_gg$ of the regularizations, have outputs and dependence regions in a
  compact set disjoint from $N_I$, with $\sup\|\eta\|_{L^2}\le b_\eta$ and
  $\sum\|g\|_{C^0}\operatorname{TV}(\chi_g)\le b_\eta$.
\item\label{hyp:X3} \emph{Ends.} The original ends are nondegenerate
  gauge ends, transverse strip ends or mixed ends with fixed
  weights below their gaps. Every nonconstant autonomous trajectory on an
  original end that can occur has index at least one, where at mixed ends
  the index is taken with fixed incoming and free outgoing limit.
\item\label{hyp:X4} \emph{Regularity at scale zero.} For every choice of
  asymptotic labels, there are no solutions on $X_0$ of index less than
  zero, and the index-zero solutions form a finite set of regular
  points, the index being computed in the extended weighted spaces at the two new
  ends.
\item\label{hyp:X5} \emph{Index-energy law at scale zero.} For fixed
  labels, $\ind D_V-8\kappa(V)$ is the same constant $c$ for all configurations $V$ on $X_0$, and $2\kappa(V)-2\kappa_0\in\ZZ$ for a fixed
  $\kappa_0$.
\end{enumerate}
The data outside $N_I$, including the neck lengths of any ordinary necks,
are chosen before $\eps$.
\end{hypothesis}

Hypothesis~\ref{hyp:X5} is verified for the domains used in this paper in
Lemma~\ref{lem:zero-scale-law}.

\subsection{Statement}

\begin{theorem}[collapse theorem]\label{thm:collapse}
Let $(C,g_C,\mathbb E)$ satisfy Definition~\ref{def:clean-datum} and
$\{X_\eps\}_{\eps\ge0}$ the domains of Definition~\ref{def:insertion},
with data satisfying Hypothesis~\ref{hyp:collapse}. Fix asymptotic
labels.
\begin{enumerate}[label=(\roman*)]
\item \emph{Index-energy law.} Every configuration $U$ with
  these labels on $X_\eps$, for any $\eps>0$, has $\ind D_U=c+8\kappa(U)$
  with $2\kappa(U)-2\kappa_0\in\ZZ$, with the constants of~\ref{hyp:X5}.
  Every solution of index $d$ has analytic energy at most
  $\pi^2(d-c)+B$ in the normalization of
  Section~\ref{subsec:energy-convention}, with $B=B(b_\eta,\sum\|f\|_{C^0})$ independent of
  $\eps$.
\item \emph{Loss formula.} For every sequence $\eps_\nu\to0$ and
  solutions $U_\nu$ on $X_{\eps_\nu}$ of index $d$, a subsequence
  converges in the sense of Proposition~\ref{prop:weak-limit} to a solution
  $V$ on $X_0$ and broken trajectories $\Gamma_j$ on the original ends, with
  \[
  d-\ind(V)-\sum_j\ind(\Gamma_j)=4n,\qquad n\in\ZZ_{\ge0},
  \]
  the indices of $\Gamma_j$ taken as in~\ref{hyp:X3}. If $d\le3$,
  $\ind(V)\ge0$ and every $\ind(\Gamma_j)\ge1$, then $n=0$ and no energy is lost.
\item \emph{Counts.} There is $\eps_0>0$ such that for $0<\eps<\eps_0$
  there are no solutions on $X_\eps$ of negative index, and gluing (the
  construction of Section~\ref{sec:collapse-coverage}) defines a bijection
  \[
  \cG_\eps\colon\cM^0(X_0)\longrightarrow\cM^0(X_\eps)
  \]
  onto a set of regular points, preserving asymptotic labels and bundle
  markings.
\end{enumerate}
\end{theorem}

\begin{lemma}\label{lem:zero-scale-law}
Hypothesis~\ref{hyp:X5} holds when $X_0$ is obtained from DFL domains of
compression pairs and a polygon domain with pairwise transverse
Lagrangians by joining along ordinary strip ends, as in
Sections~\ref{sec:two-handles} and~\ref{sec:one-three}, in the target $\Msf_h$
or in the product target $\Msf_h^-\times\Msf_g$ with labels among products of
handlebody Lagrangians and the correspondence $L_C$.
\end{lemma}

\begin{proof}
For a DFL domain, $\ind=8\kappa+\mathfrak c(x)-\mathfrak c(a)$
\cite[Proposition~3.20]{DFL}, and two configurations with the same labels
have $2(\kappa-\kappa')\in\ZZ$ \cite[Lemma~3.15]{DFL}. For a polygon domain
with simply connected Lagrangians with $\pi_2=0$ in the simply connected
$\Msf_h$, relative classes with the same corners form a torsor over
$\pi_2(\Msf_h)=\ZZ$, on which the positive generator changes $\kappa$ by
$\frac12$ in the normalization of Section~\ref{subsec:energy-convention}, and
the index by $2c_1=4$; this
gives the law with its own constant. In the product target
$\Msf_h^-\times\Msf_g$ of Section~\ref{sec:one-three}, with Lagrangians among
products of simply connected Lagrangians with $\pi_2=0$ and the correspondence
$L_C$, relative classes with the same corners form a torsor over
$\pi_2(\Msf_h^-\times\Msf_g)=\ZZ^2$ modulo the image of $\pi_2(L_C)=\ZZ$, which is
generated by spheres of zero area and zero Chern number and maps to
$\ZZ(1,1)$. The quotient is $\ZZ$, and its generator changes $\kappa$ by $\frac12$ and
the index by $4$, so the law holds with its own constant. Along an
ordinary strip end with
transverse intersection the cross-section is invertible, so the index is
additive \cite[Proposition~3.9]{Donaldson-book} and the topological
energy is additive by locality. The constants add.
\end{proof}

\subsection{Fredholm theory and excision at a fixed scale}

\begin{lemma}[Green identity]\label{lem:green}
Let $B$ be a flat connection on $C$, and let $x=(p,b)$, $y=(\pi,c)$ be
pairs of an $\ad B$-valued function and one-form on $C$ with
$b(n)=c(n)=0$ on $\partial C$. Put $\mathcal B(p,b)=(d_B^*b,
-*d_Bb+d_Bp)$. Then
\[
\langle\mathcal Bx,y\rangle_{L^2(C)}-\langle x,\mathcal By\rangle_{L^2(C)}
=-\int_{\partial C}\langle b\wedge c\rangle .
\]
If the traces $b|_{\partial C}$, $c|_{\partial C}$ are closed, the
right side equals $-\omega(\tau b,\tau c)$, where $\tau$ denotes the
harmonic part of the trace and $\omega$ is as in
Convention~\ref{conv:orientation}. For a section variation $\nu$ on a
strip $[0,1]_r$ with $J$-compatible metric,
$\langle-J\partial_r\nu,\mu\rangle-\langle\nu,-J\partial_r\mu\rangle
=-[\omega(\nu,\mu)]_{r=0}^{r=1}$.
\end{lemma}

\begin{proof}
The scalar terms cancel because the normal components vanish. For the
one-form terms, $\langle-*db,c\rangle-\langle b,-*dc\rangle
=-\int_C(db\wedge c-b\wedge dc)=-\int_{\partial C}b\wedge c$ by Stokes.
Exact parts of closed traces pair to zero with closed forms, and on a
torus boundary component every closed trace is exact. The strip identity
is integration by parts in $r$.
\end{proof}

With the matching condition and Convention~\ref{conv:orientation}, the
seam contributions of Lemma~\ref{lem:green} from the compression and from
the section cancel. This makes the mixed tangential operators symmetric,
and it is the mechanism behind every uniform estimate below.

At a mixed end the stationary solutions form a Morse--Bott family, and we use
\emph{extended weighted spaces}, as in \cite[Definitions~3.8--3.9]{DFL}: the
exponentially weighted Sobolev space of decaying fields plus the span of the
stationary fields $c_e(v)$, $v$ in the tangent space of the critical family,
cut off to the end. The coefficient $v$ is the \emph{asymptotic value} of the
field. This is the extended $L^2$ space of Atiyah, Patodi and Singer
\cite{APS} in the weighted form.

\begin{proposition}\label{prop:fredholm}
For every $\eps>0$ and every configuration $U$ on $X_\eps$, the linearized
mixed operator $D_U$, with its gauge-fixing component and on the extended
weighted spaces at the mixed ends, is Fredholm \cite[Proposition~3.18]{DFL}. The same holds on $X_0$, where the two new
ends are mixed ends with asymptotic values in $T_{q_\pm}L(C)$. The index is
locally constant along continuous paths of configurations
with fixed asymptotic data.
\end{proposition}

\begin{proof}
At fixed $\eps>0$ the region $N_I$ is compact with a smooth metric, and by
Lemma~\ref{lem:conformal-models}(i) its seam is a product mixed seam after
a conformal change that does not affect the Fredholm property. The local
estimates on the compact part are \cite[Theorem~4.17 and Remark~4.19]{DFL};
the end parametrices are those of \cite[Section~4.2.1]{DFL}. On $X_0$ the
new ends are mixed ends of the same type. Paths of configurations are identified with paths of operators on fixed spaces by
local matching charts, and the index is continuous.
\end{proof}

\begin{lemma}[excision]\label{lem:excision}
Let $D_1,D_2,D_3,D_4$ be mixed operators, Fredholm on weighted spaces, on
domains $X_1=P\cup F$, $X_2=P'\cup F'$, $X_3=P\cup F'$, $X_4=P'\cup F$,
glued along a common collar $K=[-1,1]\times\Gamma$, where $\Gamma$ is a
slice of the form $(\pm C)\sqcup[0,\pi]\sqcup(\pm C)$ with the boundary
conditions of Definition~\ref{def:admissible-config}. Suppose the
operators agree on the pieces $P$, $P'$, $F$, $F'$ as indicated and on
$K$. Then $\ind D_1+\ind D_2=\ind D_3+\ind D_4$.
\end{lemma}

\begin{proof}
Let $\phi,\psi$ be functions of the collar coordinate with
$\phi^2+\psi^2=1$, $\phi=1$ on the $P$ and $P'$ sides and $\psi=1$ on the
$F$ and $F'$ sides. The map
$R(x_1,x_2)=(\phi x_1+\psi x_2,\ -\psi x_1+\phi x_2)$, where $x_1$ is read
on $X_1$ and $x_2$ on $X_2$ and the products are transplanted to $X_3$ and
$X_4$ through the common collar, is an isomorphism between the domains of
$D_1\oplus D_2$ and $D_3\oplus D_4$: multiplication by a function of the
collar coordinate alone preserves the vanishing of the normal component,
the closedness of the traces and the harmonic matching, because these
conditions are linear and pointwise in the collar coordinate. Then
$(D_3\oplus D_4)R-R(D_1\oplus D_2)$ is multiplication by
$(\partial\phi,\partial\psi)$ on the compact collar, which is compact from
the first-order domain to $L^2$ by the Rellich lemma on the compact
domain $K$, whose corners, where the seams meet the cut, are those of a
product of an interval with a manifold with boundary. Hence the two
Fredholm operators have the same index (compare the addition property of \cite[Section~3.3]{Donaldson-book}).
\end{proof}

\subsection{The index--energy law at positive scale}\label{subsec:index-law}

\begin{lemma}[flattening]\label{lem:flattening}
Let $U$ be a configuration on $X_\eps$, $\eps>0$, and
$q_*\in L(C)$. There is a continuous path of configurations,
with fixed ends and fixed data outside a compact neighbourhood $K_2$ of
the rectangle, from $U$ to a configuration $U^\flat$ equal to the constant
flat pair $(B_{q_*},q_*)$ on a smaller neighbourhood $K_0$ of the
rectangle. The topological energy is constant along the path, and
$\kappa(U)=\kappa_0(V)$, where $V$ is the configuration on $X_0$ obtained
from $U^\flat$ by removing the rectangle and completing the two new ends by
$(B_{q_*},q_*)$.
\end{lemma}

\begin{proof}
\emph{The section.} In a star-shaped half-disk chart about the rectangle,
contract the section to its value $u(0)$ near the rectangle by
$u_s=u\circ((1-s\chi)z)$, with $\chi$ a cutoff; then move $u(0)$ to $q_*$
along a path $c$ in the connected moduli space, $u_{1+s}=c(s\rho)$, with
$\rho$ a smaller cutoff. The seam carries no Lagrangian condition at
positive scale, so the path may leave $L(C)$.

\emph{Flat surface traces.} Lift the seam values to actual flat
connections on the matching surfaces by path lifting in the principal
bundle of flat connections over the moduli space, with structure group the
determinant-one gauge group modulo its centre, keeping the lift constant
where the base path is constant. At the end the trace over a central time
interval differs from $B_{q_*}$ by a path $h$ in the gauge group. The
correction needed is a path of surface gauge transformations over the
interval, equal to $h$ in the middle and to the identity at both ends; it
exists because the determinant-one gauge group of the surface is connected.
No loop of gauge transformations arises, so the correction introduces no lifting class. Do the same on the other boundary components.

\emph{Connections.} Let $\beta_s$ be the resulting boundary connections.
Choose an auxiliary connection equal to $B_{q_*}$ near the rectangle and
to $A$ outside $K_2$, interpolate affinely, and add a collar extension of
the difference between $\beta_s$ and the trace of the interpolant. The
result has the prescribed flat traces at every parameter and equals the
flat product near the rectangle at the end.

\emph{Energy.} The variation of the topological energy along a compactly
supported path is a boundary integral; on the matching surfaces the gauge
term cancels the section term by the flat matching condition and
Convention~\ref{conv:orientation}, exact parts integrate to zero, and on a
torus face the variations are exact because $H^1(T^2;\ad)=0$. So the
energy is constant. The removed region is flat with constant section, and
the added tails have zero energy density, so $\kappa(U^\flat)=\kappa_0(V)$.
\end{proof}

Let $F_\eps(q)$ denote the linearized operator of the \emph{flat cap at
scale} $\eps$: the section is a half-disk bounded by the seam over
$(-2,2)$, the gauge piece is $(-2,2)\times C$ with the metric $g_\eps$,
the two lateral interfaces are completed by the cylinder with slice
$(-C)\sqcup[0,\pi]\sqcup C$, and the configuration is the constant flat pair
$(B_q,q)$. Let $F_0(q)$ be the zero-scale operator: the section is
$S_0=\overline{\Hb}\setminus\{\pm1\}$ with boundary arcs
$\ell_-=(-\infty,-1)$, $I=(-1,1)$, $\ell_+=(1,\infty)$, the gauge pieces are
$Z_\pm=\ell_\pm\times C$ matched along $\ell_\pm$, the condition $L(C)$ is
imposed on $I$, and there are three ends: the mixed ends $E_{\pm1}$ with
slice $C\sqcup[0,\pi]$ and the end $E_\infty$ with slice
$C\sqcup[0,\pi]\sqcup C$. All operators act on extended weighted spaces, with one asymptotic value at
each end, represented by the cut-off stationary fields $c_e(v)=(0,h_v,v)$ at $E_{\pm1}$ and $(h_v,v,h_v)$ at $E_\infty$,
for $v\in T_qL(C)$.

\begin{proposition}\label{prop:caps}
For every $q\in L(C)$, $\ker F_\eps(q)\cong\ker F_0(q)\cong T_qL(C)$ and
$\coker F_\eps(q)=\coker F_0(q)=0$. In particular
$\ind F_\eps(q)=\ind F_0(q)=n$.
\end{proposition}

\begin{proof}
We give the proof for $F_0(q)$. The positive-scale case is the same with one
end and no Lagrangian arc; there the kernel lies in $T_qL(C)$ because the
constant section value must be the boundary class of a closed form on all of
$C$, namely the harmonic field $h_v$. Use the weight $\delta$ below all end gaps
and below the Coulomb gap of the Neumann Laplacian on $C$.

\emph{Energy identity.} For $(\zeta,\nu)$, $(\zeta',\nu')$ in the domain,
\[
\int_{Z_+\sqcup Z_-}\tr(d\zeta\wedge d\zeta')+\int_{S_0}\omega(d\nu\wedge d\nu')=0 .
\]
Indeed both integrands are exact because $B_q$ is flat. On
$\ell_\pm\times T^2$ the boundary terms vanish since closed traces are
exact; on $\ell_\pm\times\Sigma$ the gauge term is
$\pm\int\omega(\nu,\partial_t\nu')dt$ and cancels the seam term of the
section by Convention~\ref{conv:orientation}; on $I$ the term vanishes
because $T_qL(C)$ is Lagrangian; at the ends the cross terms between a stationary field and a decaying remainder form the Green pairing of the
self-adjoint end operator and vanish, and the remaining terms tend to
zero.

\emph{Kernel.} For $(\zeta,\nu)$ in the kernel the identity gives
$\|d\zeta\|^2+\|d\nu\|^2=0$, so $\nu\equiv v$ with $v\in T_qL(C)$ by the
condition on $I$. On $Z_+=\RR\times C$ write $\zeta=b_t+\varphi_tdt$;
$d\zeta=0$ gives $d_Cb_t=0$ and $\partial_tb_t=d_C\varphi_t$, so the class of
$b_t$ is constant with boundary value $v$, and $b_t=h_v+d_C\psi_t$ with
$\partial_t\psi_t=\varphi_t$. Hence $\zeta=\pi^*h_v+d\psi$. The pullback
$\pi^*h_v$ is closed, coclosed and has vanishing normal component for any
metric conformal to $dt^2+g_C$ by a factor depending on $t$, so the gauge-fixing equation gives $\Delta\psi=0$ with Neumann condition and decay, and $\psi=0$ by
the maximum principle. So the kernel is $\{(\pi^*h_v,v)\}\cong T_qL(C)$.

\emph{Adjoint domain.} Let $w=(\xi,\mu,\upsilon)$ in the dual weighted space
annihilate the range. Testing with $(d\chi,0)$, $\chi$ with Neumann
condition, and using that the Neumann Laplacian is onto on the weighted
spaces, gives $\xi=0$. Interior tests give $d\mu=0$ with $\mu$ self-dual
and $d\upsilon=0$ for the section $(0,1)$-form. Near $\ell\times\partial C$ write
$\mu=dt\wedge m+*_Cm$. Boundary tests give $m(n)=0$, closed restrictions
of $m$ to $\Sigma$ and to $T^2$, and, by Convention~\ref{conv:orientation},
$[m|_\Sigma]=\upsilon(\partial_t)$ on the seams; on $I$, $\upsilon(\tau)\in T_qL(C)$,
since a Lagrangian subspace is its own symplectic complement. These are
the conditions of the original domain under $\mu\mapsto m$.

\emph{Decay.} On each end $w$ solves the translation-invariant adjoint
equation, so $w=w_0+O(e^{-\gamma\rho})$ with $w_0$ in the kernel of the end
operator. Pair with $F_0(q)(\beta_e(\cdot-R)c_e(v))=\beta_e'(\cdot-R)c_e(v)$, where
$\beta_e$ is a cutoff; letting $R\to\infty$ removes the contribution of the
decaying part and gives $\langle c_e(v),w_0\rangle=0$ for all $v$, so $w_0=0$ and $w$
decays at a rate exceeding $\delta$.

\emph{Vanishing.} The section surface $S_0$ is contractible and $\upsilon$ is
closed and decaying, so $\upsilon=d\nu'$ with $\nu'$ converging at each puncture;
normalize $\nu'(1)=0$; on $I$, $\nu'\in T_qL(C)$. On $Z_+$ put
$\zeta'_+(t)=\pi^*h_{\nu'(1)}+\int_{-\infty}^tm\,dt'$; then $d\zeta'_+=\mu$,
$\zeta'_+(n)=0$, its traces are closed, and by the seam condition its
boundary class is $\nu'(t)$. At the other end $\zeta'_+$ converges to
$h_{\nu'(1)}+M$ with $M=\int m$ closed, coclosed and normal-free, hence
harmonic with boundary class $\nu'(\infty)-\nu'(1)$; so
$\nu'(\infty)\in T_qL(C)$ and $\zeta'_+\to h_{\nu'(\infty)}$. The same on
$Z_-$. Thus $(\zeta',\nu')$ lies in the domain, with asymptotic values at the three ends, and $0=\langle F_0(q)(\zeta',\nu'),w\rangle=\|\mu\|^2+\|\upsilon\|^2$.

The primitive $\zeta'_\pm$ need not decay, since $M$ may be nonzero; it is
asymptotically constant, and this is exactly what the extended spaces allow. In the positive-scale case normalize instead $\nu'(\infty)=0$ and put
$\zeta'=\int_{-\infty}^tm\,dt'$: both seam ends tend to the same end, so $M=\int m$ is
harmonic with zero boundary class and vanishes by~(C3).
The index is therefore $\dim\ker=n$.
\end{proof}

\begin{proof}[Proof of Theorem~\ref{thm:collapse}(i)]
Let $U$ be a configuration on $X_\eps$. By Lemma~\ref{lem:flattening} and
Proposition~\ref{prop:fredholm}, $\ind D_U=\ind D_{U^\flat}$ and
$\kappa(U)=\kappa_0(V)$. The four operators $D_{U^\flat}$, $F_0(q_*)$,
$D_V$ and $F_\eps(q_*)$ agree on a fixed collar around the cut with slice
$(-C)\sqcup[0,\pi]\sqcup C$: on the $U^\flat$ and $F_\eps$ side it is the
rectangle at positive scale near its lateral interfaces, and on the $V$
and $F_0$ side the corresponding part of the new mixed ends, and all four
configurations are the flat pair there. Lemma~\ref{lem:excision} gives
$\ind D_{U^\flat}+\ind F_0(q_*)=\ind D_V+\ind F_\eps(q_*)$. By
Proposition~\ref{prop:caps} the cap indices are equal, so
$\ind D_U=\ind D_V=c+8\kappa_0(V)=c+8\kappa(U)$ by~\ref{hyp:X5}.

For a solution, the analytic energy equals $8\pi^2\kappa$ plus
$2\|\eta\|^2\le2b_\eta^2$ plus the terms $2\int\chi'f$ of the end perturbations,
bounded by twice the sum of their $C^0$ norms, plus the analogous terms of the nonautonomous
holonomy perturbations, bounded by $2b_\eta$; the first two terms are
\cite[(3.24)]{DFL}, and the end terms are estimated as in the proof
of Theorem~\ref{thm:mixed-iso}.
The bound follows.
\end{proof}

\subsection{Regularity along the collapsing seam}\label{subsec:seam-regularity}

We will use the following quantitative consequence of fixed-geometry
mixed regularity. The energy is measured on a larger patch than
the one on which derivatives are estimated.

\begin{lemma}[Quantitative local seam estimate]\label{lem:quantitative-seam}
Consider the unperturbed mixed equation in the product geometry of
\cite[(1.1)]{DFL-mixed}, with fixed surface metric and its induced
complex structure. Let $Q_2$ be a gauge half-disk times the matching
surface, together with the section half-disk, and let $Q_1$ be the
concentric patch of half the radius. There are $e_*>0$ and $C_k$
such that every solution of energy $e\le e_*$ on $Q_2$ satisfies
\[
 \sup_{Q_1}(|\nabla_A^kF_A|+|\nabla^kdu|)\le C_k\sqrt e.
\]
The constants are uniform on compact families of these product data.
\end{lemma}
\begin{proof}
Choose nested intermediate patches. Mixed small-energy regularity
and compactness give gauges in which solutions of energy below a
fixed concentration threshold have bounded smooth coefficients on
these patches \cite[Theorem~4.2, Lemma~4.25 and pp.~38--39]{DFL-mixed}.
As the energy tends to zero their limits are flat connections and
constant sections. The possible flat surface classes form a compact
space, so finitely many gauge and target charts suffice. These are
qualitative coefficient bounds; the energy rate is obtained next.

Let $t$ be the coordinate tangent to the seam and set
$\zeta=\iota_{\partial_t}F_A$, $\xi=\partial_tu$.
Covariant differentiation of the unperturbed mixed equation in
this product coordinate shows that $(\zeta,\xi)$ satisfies its
linearized equation. The gauge row is zero by the Yang--Mills
identity $d_A^*F_A=0$. Its normal component on the seam is zero,
since anti-self-duality identifies $F_{ts}$ with the curvature of
the flat surface restriction. Differentiating that restriction and
its class gives the closed-trace and harmonic-matching conditions.
Thus $(\zeta,\xi)$ lies in the local mixed boundary domain,
including the gauge row, as in the translation-field argument of
\cite[Proposition~5.29]{DFL}. Use matched scalar cutoffs on the
gauge and section patches, with identical seam values and supports
away from the artificial edges. Their commutators are zeroth order.
The local all-order boundary estimates of
\cite[Theorems~5 and 5.9, Remarks~5.73 and 5.76]{DFL-mixed},
with the uniformity of \cite[Remark~4.19]{DFL}, give
\[
 \|\zeta\|_{H^{j+1}(Q_1)}+\|\xi\|_{H^{j+1}(Q_1)}
 \le C_j\bigl(\|\zeta\|_{L^2(Q')}+\|\xi\|_{L^2(Q')}\bigr)
 \le C'_j\sqrt e,
\]
where $Q_1\Subset Q'\Subset Q_2$. One can equivalently trivialize the moving matching domain by
its smooth chart $Q$ and commute $Q\partial_tQ^{-1}$, rather than
raw tangential differentiation, before recovering normal derivatives
from the equation. Use successively smaller patches. The
qualitative smooth bounds
control every coefficient and its commutators at each fixed order.
They also give uniform constants over the finite chart cover.

For an anti-self-dual two-form on this product, contraction with
$\partial_t$ determines all its components and
$|F_A|^2=2|\zeta|^2$. The Cauchy--Riemann equation similarly
determines $du$ from $\xi$. These identities and Sobolev embedding
prove the claimed bounds. The threshold is chosen once below
concentration; only $C_k$ depends on the derivative order.
\end{proof}

Away from the seam we use the corresponding interior estimates
\cite[Lemma~4.50]{DFL-mixed}.

The following lemma is the source of all uniform local control near the
collapse. For a seam point $b$ (on $I$, or on an endpoint neck) let
${\mathfrak s}_b$ be the local compression scale (${\mathfrak s}_b=\eps$ on $I$, and the
value of Lemma~\ref{lem:conformal-models}(ii) on the necks), and let
$\mathcal H_r(b)$ be the half-disk of radius $r$ in the section about $b$ together
with the gauge region over $[b-r,b+r]$. Let $E(b,r)$ be the energy of a
solution on $\mathcal H_r(b)$, and $r_0>0$ fixed.

In the following lemma $\Sigma\subset\partial C$ is the union of the boundary components
matched to the section and $\Msf$ their moduli space, the remaining components carrying the
flatness condition, as in Section~\ref{subsec:model}. The topological energy $\mathcal Q$ of
a configuration is that of Section~\ref{subsec:energy-convention}: the Chern--Weil integral
of the gauge part plus the symplectic area of the section, normalized so that it equals the
analytic energy of a solution of the unperturbed equations.

\begin{lemma}[quantization of capped configurations]\label{lem:cap-quantization}
There is $\hbar>0$, depending only on $(C,g_C,\mathbb E)$, with the following
property. Let $x\in P$, and let $(A,u)$ be a continuous, piecewise smooth configuration on
$(\RR\times C)\cup\overline\Hb$: $A$ is a connection on $\RR\times\mathbb E$ whose
restriction to every slice $\{t\}\times\partial C$ is flat, and $u\colon\overline\Hb\to\Msf$
satisfies $u(t)=[A|_{\{t\}\times\Sigma}]$ for $t\in\RR=\partial\Hb$. Suppose that $(A,u)$
equals the constant pair $(B_x,x)$ outside a compact set. Then $\mathcal Q(A,u)\in\hbar\ZZ$.
\end{lemma}

\begin{proof}
(1) \emph{Invariance.} $\mathcal Q$ is unchanged by gauge transformations of
$\RR\times\mathbb E$, whose action on the integrands is pointwise, and by homotopies through
configurations of the same kind that are constant outside a fixed compact set. For the
latter, Bianchi's identity and the closedness of $\omega$ reduce the variation of $\mathcal Q$
to boundary integrals over $\RR\times\partial C$ and $\partial\Hb$, which cancel by the
matching condition and Convention~\ref{conv:orientation} exactly as in the energy paragraph of
the proof of Lemma~\ref{lem:flattening}; on the torus faces the variations are exact. The slice
contributions at the interfaces of the pieces cancel by continuity.

(2) \emph{Removing the section.} Choose $T_0$ such that $(A,u)=(B_x,x)$ for $|t|\ge T_0$ and
outside the half-disk of radius $T_0$. For $s\in[0,1]$ put $u_s(z)=u(z+isT_0)$ and
$q_s(t)=u_s(t)$ for $t\in\RR$; then $u_1\equiv x$ and $q_s(\pm T_0)=x$. By (C1) the flat
connections on $\partial C$ form a principal bundle over $\Msf$ whose structure group is the
determinant-one gauge group of $\partial C$ modulo its centre, as in the proof of
Lemma~\ref{lem:flattening}. The map $(s,t)\mapsto q_s(t)$ on $[0,1]\times[-T_0,T_0]$ has the
prescribed lift $\alpha(t)=A|_{\{t\}\times\partial C}$ on $\{0\}\times[-T_0,T_0]$ and the
constant lift $B_x|_{\partial C}$ on $[0,1]\times\{\pm T_0\}$; the union of these three edges
is a deformation retract of the square, so the lift extends to flat connections $\alpha_s(t)$
with $[\alpha_s(t)|_\Sigma]=q_s(t)$. With a bounded extension operator $\mathcal E$ from
tangential one-forms on $\partial C$ to one-forms on $C$, supported in a collar, the
connections $A_s=B(t)+\mathcal E(\alpha_s(t)-\alpha(t))+\Psi\,dt$, where $A=B(t)+\Psi\,dt$,
together with $u_s$, form a homotopy as in (1). Hence $\mathcal Q(A,u)=\mathcal Q(A_1,x)$, and
the latter is the Chern--Weil term alone.

(3) \emph{Reduction to a gauge transformation.} The restrictions of $A_1$ to the slices of
$\RR\times\partial C$ are flat of class $x$, so $\alpha_1(t)=h(t)^*(B_x|_{\partial C})$ for a
path $h$ in the determinant-one gauge group of $\partial C$ with $h(t)=1$ for $t\le-T_0$; for
$t\ge T_0$ the path is constant and central, the stabilizer of an irreducible flat connection
being the centre. Restriction to $\partial C$ is a fibration of gauge groups, by the homotopy
extension property of $(C,\partial C)$, so $h$ lifts to a path $\hat h$ of determinant-one gauge
transformations of $C$ with $\hat h(t)=1$ for $t\le-T_0$ and $\hat h$ constant for $t\ge T_0$.
The gauge transformation $(t,y)\mapsto\hat h(t)(y)^{-1}$ carries $A_1$ to a connection $A_2$
whose restriction to every slice of $\RR\times\partial C$ is $B_x|_{\partial C}$, equal to $B_x$
for $t\le-T_0$ and to $g^*B_x$ for $t\ge T_0$, where $g=\hat h(T_0)^{-1}$ is central on
$\partial C$. By (1), $\mathcal Q(A,u)=\mathcal Q(A_2,x)$.

(4) \emph{The period.} Let $\mathcal G_\partial$ be the group of determinant-one gauge
transformations of $C$ that are central on $\partial C$. For $g\in\mathcal G_\partial$ let
$\Xi(g)=\mathcal Q(A',x)$ for any connection $A'$ as in (3), with slice restrictions
$B_x|_{\partial C}$, equal to $B_x$ for $t\ll0$ and to $g^*B_x$ for $t\gg0$. This does not
depend on $A'$: two choices are joined by the affine homotopy, which keeps the slice
restrictions fixed, and (1) applies. Concatenation and gauge invariance give
$\Xi(g_1g_2)=\Xi(g_1)+\Xi(g_2)$, and $\Xi$ is constant on path components: a path
$\lambda\mapsto g_\lambda$ in $\mathcal G_\partial$ gives the connection
$g_{\lambda(t)}^*B_x$, which is flat on every slice, so its Chern--Weil integrand
$F\wedge F$ vanishes, and whose slice restrictions are $B_x|_{\partial C}$. Thus $\Xi$ is a
homomorphism from $\pi_0(\mathcal G_\partial)$ to $\RR$. Let $\mathcal G_1\subset\mathcal G_\partial$ be the normal subgroup of transformations equal
to $1$ on $\partial C$. Its index $N=[\mathcal G_\partial:\mathcal G_1]$ is finite, since
$\mathcal G_\partial/\mathcal G_1$ embeds in the locally constant central-valued maps on
$\partial C$, and its group of components is $H^3(C,\partial C;\ZZ)\cong\ZZ$ by obstruction
theory, the fibre $S^3$ of the group bundle being $2$-connected. For $g\in\mathcal G_1$ we
may, after a homotopy in $\mathcal G_\partial$, take $g=1$ near $\partial C$, and then
$\Xi(g)$ is the Chern--Simons difference of $g^*B_x$ and $B_x$, whose boundary term vanishes,
so $\Xi(g)=c_\Xi\deg g$ with $c_\Xi$ independent of $x$. For $g\in\mathcal G_\partial$ the
power $g^N$ lies in $\mathcal G_1$, so $N\,\Xi(g)=\Xi(g^N)\in c_\Xi\ZZ$. All values of $\Xi$,
for every $x\in P$ and on every component of $P$, therefore lie in $(c_\Xi/N)\ZZ$. Take
$\hbar=|c_\Xi|/N$ if $c_\Xi\ne0$, and $\hbar=1$ if $c_\Xi=0$, in which case $\Xi$ vanishes
identically. Thus $\hbar$ depends only on $C$.
\end{proof}

\begin{lemma}[seam regularity]\label{lem:seam-regularity}
There are $e_0,\alpha,C>0$, independent of $\eps$ and of the seam point, such that
for every solution $U=(A,u)$ on $X_\eps$ and every seam point $b$ with
$\dist(b,z_\pm)\ge c_{\rm core}\eps$ and $E(b,r_0)\le e_0$, with $c_{\rm core}$ as in
Definition~\ref{def:insertion}(i):
\begin{enumerate}[label=(\roman*)]
\item $E(b,r)\le C\,E(b,r_0)\,(r/r_0)^\alpha$ for ${\mathfrak s}_b\le r\le r_0$;
\item $|du(z)|\le C\sqrt{E(b,r_0)}\,(y+{\mathfrak s}_b)^{\alpha/2-1}$ at distance
  $y\le r_0/2$ from the seam, and
  $d(u(z),u(z'))\le C\sqrt{E(b,r_0)}\,|z-z'|^{\alpha/2}$ on $\mathcal H_{r_0/2}(b)$;
\item $\dist(u(b),L(C))\le C\sqrt{E(b,r_0)}\,{\mathfrak s}_b^{\alpha/2}$; in the
  rescaled coordinates on the patch of size ${\mathfrak s}_b$ about $b$ the
  curvature, the section derivative and all their covariant derivatives are
  bounded by $C_k\sqrt{E(b,r_0)}\,{\mathfrak s}_b^{\alpha/2}$; and in gauges on the
  compression slices over $[b-{\mathfrak s}_b,b+{\mathfrak s}_b]$,
  $\|B(t)-B_{q(t)}\|_{C^k(C)}\le C_k\sqrt{E(b,r_0)}\,{\mathfrak s}_b^{\alpha/2}$ with
  $q(t)\in L(C)$.
\end{enumerate}
\end{lemma}

\begin{proof}
\emph{Rescaled geometry.} In the coordinate $w\mapsto(w-b)/{\mathfrak s}_b$ (on
$I$) or $\Phi_b$ with $\Phi_b'=1/{\mathfrak s}$ (on the necks), the collapsing
compression and the section have exactly the product mixed geometry of
\cite[(1.1)]{DFL-mixed} with the fixed surface metric, and unit patches
correspond to regions of size $C{\mathfrak s}_b$. Energy is conformally
invariant. So Lemma~\ref{lem:quantitative-seam} in that fixed product
geometry and the interior estimate \cite[Lemma~4.50]{DFL-mixed},
which allows a $C^k$-compact family of metrics, apply on every unit patch with constants independent of $\eps$
and $b$: on a patch of energy at most $e_0$ the curvature, the section
derivative and all their covariant derivatives are bounded by $C\sqrt e$
in the rescaled coordinates.

\emph{Decay above the local scale.} Fix $r$ with ${\mathfrak s}_b\le r\le r_0$ and write
$E(b,r)=G(b,r)+S(b,r)$, with $G$ the gauge part and $S$ the section part; the gauge region of
$\mathcal H_r(b)$ is $[b-r,b+r]\times C$, and its section region is a half-disk $D_r$ with diameter
$[b-r,b+r]$ on the seam.

\emph{The endpoint slices.} The gauge energy density on the compression is
${\mathfrak s}^{-1}\|F_B\|^2_{L^2(C)}$ per unit time, so the connections $B_\pm$ on the slices
$t=b\pm r$ satisfy $\|F_{B_\pm}\|^2\le C\lambda_\eps(b\pm r)\partial_rG(b,r)\le Cr\partial_rG(b,r)$,
where $\lambda_\eps$ is the compression scale. When
this is small, the proof of Lemma~\ref{lem:distance} gives gauge transformations $g_\pm$ and
$q_\pm\in P$ with $B_\pm=g_\pm^*(B_{q_\pm}+\beta_\pm)$, where $\beta_\pm$ lies in the Neumann slice
orthogonal to the harmonic fields and $\|\beta_\pm\|_{H^1(C)}\le C\|F_{B_\pm}\|$; the
representatives $B_q$, $q$ near $x$, are taken from a smooth local section of flat connections
through $B_x$. The image $\gamma$ of the semicircle $\partial D_r\setminus[b-r,b+r]$ satisfies
$\len(\gamma)^2\le Cr\partial_rS(b,r)$ by Cauchy--Schwarz. Its endpoints $u(b\pm r)=[B_\pm|_\Sigma]$
lie within $C\|\beta_\pm\|_{H^1}$ of $q_\pm$, so with $x=q_-$,
\[
d(q_+,x)\le C\delta_b,\qquad \delta_b=\len(\gamma)+\|\beta_-\|_{H^1}+\|\beta_+\|_{H^1}.
\]

\emph{Paths with flat boundary.} Let $B$ be a connection on $C$ with flat boundary restriction
$\alpha$, and $\|B-B_x\|_{H^1(C)}=\delta$ small; write $\alpha_x=B_x|_{\partial C}$. Every
component of $\partial C$ carries only irreducible flat connections, so the surface Coulomb
slice, in which the gauge parameter lies in $H^{3/2}(\partial C)$ (an algebra acting on
$H^{1/2}$), and the harmonic chart of the flat connections give
$\alpha=\exp(\zeta)^*(\alpha_x+h_\nu+\sigma(\nu))$, where $h_\nu$ is the harmonic representative
of $\nu\in H^1(\Sigma;\ad\alpha_x)$, $\sigma$ the quadratic correction of the chart, and
$\|\zeta\|_{H^{3/2}}+|\nu|\le C\delta$. Put
$\alpha^\lambda=\exp(\lambda\zeta)^*(\alpha_x+h_{\lambda\nu}+\sigma(\lambda\nu))$ and
\[
B^\lambda=B_x+\lambda(B-B_x)+\mathcal E\bigl(\alpha^\lambda-\alpha_x-\lambda(\alpha-\alpha_x)\bigr),
\qquad0\le\lambda\le1,
\]
with $\mathcal E\colon H^{1/2}(\partial C)\to H^1(C)$ a bounded tangential extension. Then
$B^0=B_x$, $B^1=B$, every $B^\lambda$ has the flat boundary restriction $\alpha^\lambda$,
$\sup_\lambda(\|B^\lambda-B_x\|_{H^1}+\|\partial_\lambda B^\lambda\|_{H^1})\le C\delta$, and the
boundary classes of $B^\lambda$ form a path of length at most $C\delta$. By the Sobolev
embedding $H^1\subset L^4$ on $C$, $\|F_{B^\lambda}\|_{L^2}\le C\delta$, so the path, regarded as
a connection on $[0,1]\times C$ with zero temporal component, has Chern--Weil integral at most
$C\int_0^1\|\partial_\lambda B^\lambda\|\,\|F_{B^\lambda}\|\,d\lambda\le C\delta^2$ in absolute
value. A temporal component $\Psi$ does not change this integral, since
$\int_C\langle d_B\Psi\wedge F_B\rangle=\int_{\partial C}\langle\Psi,F_{B|_{\partial C}}\rangle=0$
by Bianchi's identity and the flatness of the boundary restriction; so $\Psi$ can be chosen to
match the temporal component of the solution at the junctions. We use this with
$B=B_{q_-}+\beta_-$, for which $\delta=\|\beta_-\|_{H^1}$, and with $B=B_{q_+}+\beta_+$, for which
$\delta\le C\delta_b$, the second path traversed backwards.

\emph{The cap.} Composing with a fixed representative of the component class in
$H^1(C;\ZZ/2)$, which acts on flat connections of $C$ and preserves $L(C)$, we may take
$g_\pm$ of determinant one. The determinant-one gauge group of $C$ is connected,
because $H^3(C;\ZZ)=0$ for $\partial C\ne\varnothing$ and the fibre $S^3$ is $2$-connected;
choose a path $g(t)$, $t\in[b-r,b+r]$, from $g_-$ to $g_+$. The gauge transformation $(t,y)\mapsto
g(t)(y)^{-1}$ changes neither the topological energy of the solution on $[b-r,b+r]\times C$ nor
its boundary classes, keeps all its slice restrictions flat, the torus faces included, and makes
its end slices equal to $B_{q_\pm}+\beta_\pm$. Concatenate in $t$: $B_x$ on
$(-\infty,b-r-1]$, the path from $B_x$ to $B_{q_-}+\beta_-$ on $[b-r-1,b-r]$, the transformed
solution on $[b-r,b+r]$, the path from $B_{q_+}+\beta_+$ to $B_x$ on $[b+r,b+r+1]$, and $B_x$ on
$[b+r+1,\infty)$. On the section keep $u$ on $D_r$; on the half-annulus $\{r\le|z-b|\le
r+1,\ \operatorname{Im}z\ge0\}$, whose boundary loop consists of $\gamma$, of the boundary class
paths of the two flat-boundary paths on the seam intervals $[b-r-1,b-r]$ and $[b+r,b+r+1]$, and
of the outer semicircle mapped to $x$, use the geodesic cone from $x$ over this loop; and put
$u=x$ elsewhere. The loop has length at most $C\delta_b$ and lies in the ball of radius
$C\delta_b$ about $x$, so the cone has area at most $C\delta_b^2$. The result is a configuration
as in Lemma~\ref{lem:cap-quantization}. Its topological energy is $E(b,r)+Q_{\rm cap}$, since
the solution satisfies the unperturbed equations near the seam by
Hypothesis~\ref{hyp:collapse}\ref{hyp:X1}--\ref{hyp:X2}, and
\[
|Q_{\rm cap}|\le C\delta_b^2\le C\bigl(\len(\gamma)^2+\|\beta_-\|_{H^1}^2+\|\beta_+\|^2_{H^1}\bigr)
\le Cr\partial_rE(b,r).
\]
By Lemma~\ref{lem:cap-quantization}, $E(b,r)+Q_{\rm cap}\in\hbar\ZZ$.

\emph{Conclusion.} Take $e_0<\hbar/2$. There is $c_2>0$ such that $r\partial_rE(b,r)\le c_2$
ensures the smallness required above and $|Q_{\rm cap}|<\hbar/2$; then $E(b,r)+Q_{\rm cap}$
lies in $\hbar\ZZ\cap(-\hbar,\hbar)=\{0\}$, and $E(b,r)\le|Q_{\rm cap}|\le Cr\partial_rE(b,r)$.
If $r\partial_rE(b,r)>c_2$, then $E(b,r)\le e_0\le(e_0/c_2)\,r\partial_rE(b,r)$. Hence
$E\le C_0r\partial_rE$ for ${\mathfrak s}_b\le r\le r_0$, and integration gives (i) with
$\alpha=1/C_0$.

\emph{Derivative and Hölder bounds.} At distance $y\le c{\mathfrak s}_b$ from the
seam apply the rescaled estimates on patches of radius $C{\mathfrak s}_b$, whose
energy is $O(E{\mathfrak s}_b^\alpha)$ by (i); for $y\ge{\mathfrak s}_b$ apply the interior
estimate for holomorphic curves on disks of radius $y/2$, whose energy is
$O(Ey^\alpha)$. Since $\alpha>0$ the bound is integrable in $y$, which gives
the Hölder bound. For (iii), part (i) gives energy at most
$CE(b,r_0){\mathfrak s}_b^\alpha$ on the patch of size $R{\mathfrak s}_b$ about $b$, for a fixed
$R$; the rescaled small-energy estimates on that patch give the bounds with
the factor ${\mathfrak s}_b^{\alpha/2}$, and the spatial Neumann slice about the nearest
flat connection, with Lemma~\ref{lem:distance}, gives the gauge statement.
\end{proof}

\subsection{Compactness and the loss formula}\label{subsec:compactness}

\begin{proposition}\label{prop:weak-limit}
Let $\eps_\nu\to0$ and $U_\nu=(A_\nu,u_\nu)$ be solutions on
$X_{\eps_\nu}$ with bounded energy. After passing to a subsequence and
gauge transformations there are a finite set $Z\subset X_0$ and a solution
$V$ on $X_0$ such that:
\begin{enumerate}[label=(\roman*)]
\item $U_\nu\to V$ in $C^\infty$ on compact subsets of $X_0\setminus(Z\cup
  I\cup\{z_\pm\})$; here $X_{\eps_\nu}$ agrees with $X_0$ on every compact
  subset of the complement of $I\cup\{z_\pm\}$ for large $\nu$;
\item on compact subsets of $\overline{S_0}$ away from $Z$ and $z_\pm$,
  including points of $I$, the sections converge in $C^0$ and are
  uniformly Hölder, and the limit takes values in $L(C)$ on $I$;
\item $V$ is smooth up to $I$, the points of $Z$ are removable, and $V$
  extends to a finite-energy solution on $X_0$ whose tails at the new ends
  converge exponentially to points $q_\pm\in L(C)$;
\item the original ends break into finitely many trajectories $\Gamma_j$.
\end{enumerate}
\end{proposition}

\begin{proof}
(i) On the complement of $N_I$ and on the part of $N_I$ at distance at least
$c_{\rm core}\eps_\nu$ from $I$ and its endpoints the geometry is fixed, the secondary perturbations have
compact support away from $N_I$, and the end perturbations are autonomous beyond a fixed cutoff
band on the original ends, so this is the compactness theorem
for the mixed equation \cite[Theorem~5.2]{DFL},
\cite[Theorem~3 and Theorem~4.44]{DFL-mixed}, together with Uhlenbeck and
Gromov compactness, with a finite concentration set $Z$ outside of which
energy densities converge. Holonomy perturbation terms converge because
the loop families avoid the finite concentration set for almost every
output point.

(ii) Energy concentration along $I$ occurs at finitely many points, since
each concentration point carries at least a fixed quantum: at points of $I$
the threshold $e_0$ of Lemma~\ref{lem:seam-regularity}, elsewhere the
thresholds of the fixed-geometry theory. Away from
them, $E(b,r_0)\le e_0$ for $\nu$ large and $r_0$ small, and
Lemma~\ref{lem:seam-regularity} gives uniform Hölder bounds for $u_\nu$ up
to the seam and $\dist(u_\nu(b),L(C))\le C\eps_\nu^{\alpha/2}\to0$ on $I$.

(iii) By (ii), near a point of $I$ away from $Z$ the limit is a continuous
$J$-holomorphic map, $J$ being the complex structure $*_\Sigma$ there by~\ref{hyp:X1},
with boundary values in the embedded Lagrangian $L(C)$ and finite energy.
It is of class $W^{1,2}$ and continuous up to $I$. In coordinates in which
$L(C)=\RR^n$ and $J|_{L(C)}=J_0$, reflection across $L(C)$ and the small
oscillation of $J(u)$ give, by the Calder\'on--Zygmund inequality, $W^{1,p}$
regularity for every $p<\infty$, and boundary regularity for
pseudo-holomorphic curves with Lagrangian boundary conditions
\cite[Theorem~B.4.1]{MS} gives smoothness.
Removal of singularities at interior, seam and boundary points of $Z$
follows from \cite[Theorem~4.44]{DFL-mixed}, Uhlenbeck's theorem and
\cite[Theorem~4.1.2]{MS}. At the new ends the energy of $V$ is finite; after
the conformal rescaling of the end to the normalization of \cite{DFL}, in which
$g_C$ is scaled together with the strip width, and passing to temporal
gauge, \cite[Proposition~5.29]{DFL} gives exponential decay.

(iv) The chains on the original ends are extracted as in
\cite[Section~5.4]{DFL}.
\end{proof}

\begin{proof}[Proof of Theorem~\ref{thm:collapse}(ii)]
The analytic energy of $U_\nu$ is $8\pi^2\kappa(U_\nu)$ plus terms that
converge along the subsequence. Energy is lower semicontinuous under the
convergence of Proposition~\ref{prop:weak-limit}, including at the points
of $Z$ and along the necks of the original ends, where the action is
constant on each connected critical family. Hence the topological energy
lost is nonnegative. Glue $V$ and the $\Gamma_j$ into a configuration
with the labels of $U_\nu$; at a mixed end the incoming limit of $\Gamma_j$ equals
the outgoing limit of the previous piece, and with the convention
of~\ref{hyp:X3} (fixed incoming, free outgoing limit) the index is additive.
The topological energy is additive by locality. By part (i) applied on
$X_{\eps_\nu}$ and~\ref{hyp:X5} on $X_0$, both with the constant $c$, the lost
topological energy lies in $\tfrac12\ZZ_{\ge0}$ and the lost index is eight
times it. This is $d-\ind(V)-\sum_j\ind(\Gamma_j)=4n$ with $n\ge0$. If $d\le3$,
$\ind(V)\ge0$ and $\ind(\Gamma_j)\ge1$, then $4n\le3$ forces $n=0$; energy is
then preserved and there is no concentration.
\end{proof}

At index zero, Hypothesis~\ref{hyp:collapse}\ref{hyp:X3}--\ref{hyp:X4}
and the loss formula give $\ind(V)=0$, no broken trajectories, $n=0$, and
energy exhaustion: $Z=\varnothing$, and every compact part of the rescaled
endpoint regions and of the rescaled neighbourhood of $I$ carries energy
tending to zero.

\begin{corollary}[convergence on the middle compression]\label{cor:seam-convergence}
In this situation fix $r_0>0$ so small that the energy of $V$ on every half-disk
$\mathcal H_{r_0}(b)$ with $b\in I$ is less than $e_0/2$. Then for $\nu$ large,
uniformly for $b\in I$ at distance at least $r_0$ from the endpoints and on
rescaled unit patches of size $\eps_\nu$ about $b$: the curvature and all its
derivatives tend to zero; in the spatial Neumann slice the compression
connection tends in $C^N$ to the flat lift $B_{q_V(t)}$ of the boundary values
$q_V$ of $V$ on $I$; and the rescaled time derivatives of the harmonic coefficient $v$
(the component along the harmonic fields $h_v$, $v\in T_qL(C)$) tend to zero.
The same smallness holds on the necks where $x/\eps_\nu\to\infty$ and on the
rescaled cores at $z_\pm$.
\end{corollary}

\begin{proof}
By energy exhaustion $E(b,r_0)\le e_0$ for $\nu$ large, uniformly in $b$. By
Lemma~\ref{lem:seam-regularity}(iii) the rescaled quantities are bounded by
$C_N\sqrt{e_0}\,\eps_\nu^{\alpha/2}\to0$, and $q(t)$ is within
$C\eps_\nu^{\alpha/2}$ of the section value, which converges uniformly to $q_V$
by Proposition~\ref{prop:weak-limit}(ii). Near the endpoints, on the
necks where $x/\eps_\nu\to\infty$, Lemma~\ref{lem:entry} gives smallness
$C\sqrt{e_0}(\eps_\nu/x)^{\alpha/2}\to0$, and on the rescaled cores the convergence to
zero comes from the convergence of the rescaled transitions to the flat pairs at
$q_\pm$ (Section~\ref{subsec:necks}).
\end{proof}

\section{The collapse theorem: uniform linear estimates}\label{sec:collapse-linear}

This section proves the linear estimates on which the gluing and surjectivity
arguments of Section~\ref{sec:collapse-coverage} rest. The basic
estimates are on a model: a family of compressions at varying small scale
coupled to a holomorphic strip along the whole of one boundary line
(Sections~\ref{subsec:model}--\ref{subsec:uniform-inverse}). The endpoint
transitions are treated in Section~\ref{subsec:transition}, and the
patched inverse on the actual domain in Section~\ref{subsec:patching}.

The mechanism is the following. On the collapsing compression the operator splits into a
fast part, controlled by the Hodge theory of the compression body at scale $\mathfrak s$, and
a slow part along the harmonic fields, which is coupled to the strip through the matching
condition. The uniform Gårding estimate (Proposition~\ref{prop:garding}) controls both with
no loss in $\mathfrak s$: the seam terms of the gauge and section parts cancel by
Convention~\ref{conv:orientation}, and the slow coefficient is controlled by the boundary
trace (Lemma~\ref{lem:tie}). The
only loss is the factor $\mathsf A=1+\log(1/\eps)$, which comes from the stationary
coefficients on the two endpoint necks of length $\log(1/\eps)$
(Proposition~\ref{prop:transition-inverse}). In the gluing argument of
Section~\ref{sec:collapse-coverage} it always multiplies a positive power of $\eps$, so it
does no harm.

\subsection{The model collapsing strip}\label{subsec:model}

Let $(C,g_C,\mathbb E)$ satisfy Definition~\ref{def:clean-datum}, with a boundary
component $\Sigma\subset\partial C$ carrying the section, and let $P=L(C)$
(for simplicity of notation we suppress the other boundary components,
whose moduli spaces are points or enter as fixed factors). Fix a
compact Lagrangian $K$ in the moduli space $\Msf$ of $\Sigma$. Fix:
\begin{itemize}[leftmargin=2em]
\item a path $q\colon\RR\to P$ with $Q=\sup|\dot q|<\infty$, flat
  connections $B(t)$ representing $q(t)$, and $\Psi(t)$ with
  $\partial_tB-d_B\Psi=h_{\dot q}$;
\item a map $u\colon S=\RR\times[0,1]\to\Msf$ with $u(t,0)=q(t)$,
  $u(t,1)\in K$, bounded in $C^1$, and a compatible $J(t,r)$;
\item a scale ${\mathfrak s}\colon\RR\to(0,\bar{\mathfrak s}]$ with $|{\mathfrak s}'/{\mathfrak s}|\le H$;
  write $\varsigma={\mathfrak s}'/{\mathfrak s}$.
\end{itemize}
For $U=(p,b)$, a pair of $\ad B(t)$-valued function and one-form on $C$
depending on $t$, put $\mathsf D_t=\partial_t+[\Psi,\cdot]$,
$\mathcal B_t(p,b)=(d_B^*b,-*d_Bb+d_Bp)$ and $\Pi_{\rm f}(p,b)=(p,0)$, and
define
\[
D_{\mathfrak s}(U,\nu)=\bigl(\mathsf D_tU-{\mathfrak s}^{-1}\mathcal B_tU+2\varsigma\Pi_{\rm f}U,\
D_u\nu\bigr),
\]
where $D_u$ is the linearized Cauchy--Riemann operator. The domain
$\mathcal D$ consists of $(U,\nu)$ with $b(n)=0$ on $\partial C$, closed
traces on every component of $\partial C$, harmonic part of
$b|_\Sigma$ equal to $\nu(t,0)$, and $\nu(t,1)\in T_{u(t,1)}K$. Norms:
\begin{gather*}
\|(U,\nu)\|_X^2=\int\bigl({\mathfrak s}\|U\|^2+{\mathfrak s}\|\mathsf D_tU\|^2
+{\mathfrak s}^{-1}\|\mathcal B_tU\|^2\bigr)dt+\|\nu\|^2_{H^1(S)},\\
\|(f,g)\|_Y^2=\int{\mathfrak s}\|f\|^2dt+\|g\|^2_{L^2(S)},
\end{gather*}
all norms on $C$ in the fixed metric $g_C$. This operator is the
linearization of the mixed equation at a configuration that is flat on
the compression, in the conformal coordinates of
Lemma~\ref{lem:conformal-models}; the weight ${\mathfrak s}$ is the Jacobian of
the conformal change.

\subsection{Hodge theory on a compression body}

Let $\mathbb H_t=\{h_v:v\in T_{q(t)}P\}$ be the space of absolute
harmonic fields and $\Pi_t$ the $L^2(C)$-projection onto
$\{(0,h):h\in\mathbb H_t\}$. Constants depend only on $(C,g_C,\mathbb E)$ and on
$P$, which is compact.

\begin{lemma}\label{lem:hodge}
For $U=(p,b)$ in the static domain:
\begin{enumerate}[label=(\roman*)]
\item $\|\mathcal B_tU\|^2=\|d^*b\|^2+\|db\|^2+\|dp\|^2$;
\item $\|U-\Pi_tU\|_{H^1(C)}\le c\,\|\mathcal B_tU\|$;
\item $\ker\mathcal B_t=\{(0,h_v):v\in T_qP\}$.
\end{enumerate}
\end{lemma}

\begin{proof}
(i) $\langle*db,dp\rangle=\int_Cdp\wedge db=\int_{\partial C}p\,d(b|_{\partial C})=0$
since the traces are closed and $B$ is flat. (ii) is the
Gaffney--Friedrichs inequality for one-forms with vanishing normal
component \cite[Corollary~2.1.6 and Theorem~2.1.7]{Schwarz}, whose kernel is the space of
absolute harmonic fields, together with the Poincar\'e inequality for
$p$, which holds because $H^0(C;\ad)=0$. (iii) follows from (i), (ii)
and (C3).
\end{proof}

We call the harmonic fields $h_v$ the \emph{slow directions}: they span the kernel
of the fast operator ${\mathfrak s}^{-1}\mathcal B_t$, and their only dynamics is
inherited from the boundary value of the strip. The corresponding coefficient
is the \emph{slow coefficient}.

\begin{lemma}[slow directions are tied to the trace]\label{lem:tie}
For $U$ in the static domain write $\Pi_tU=(0,h_w)$ and let $\nu_0$ be
the harmonic part of $b|_\Sigma$. Then
\[
|\nu_0-w|\le c\,\|\mathcal B_tU\|,\qquad
\|U\|_{L^2}+\|b\|_{H^1(C)}\le c\,(\|\mathcal B_tU\|+|\nu_0|).
\]
\end{lemma}

\begin{proof}
The harmonic part of $h_w|_\Sigma$ is $w$, and the harmonic projection
of the trace of $b-h_w$ is bounded by $c\|b-h_w\|_{H^1}\le
c\|\mathcal B_tU\|$ by Lemma~\ref{lem:hodge}(ii). The second
inequality follows, since $|w|\le|\nu_0|+c\|\mathcal B_tU\|$ and
$\|h_w\|\le c|w|$.
\end{proof}

\begin{lemma}[the seam term]\label{lem:seam}
For $U(t)$ smooth in $t$ in the static domains, write $b|_\Sigma=\beta
+d\xi$ with $\beta$ harmonic, and $b|_{T}=d\xi_T$ on the other boundary
components. Then
\[
\int_{\partial C}\langle b\wedge\mathsf D_tb\rangle
=\int_\Sigma\langle\beta\wedge\mathsf D_t\beta\rangle+R,\qquad
|R|\le c\,|\dot q|\,\|b\|^2_{H^1(C)} .
\]
\end{lemma}

\begin{proof}
On $\Sigma$, $\mathsf D_t(d\xi)=d(\mathsf D_t\xi)+[a_\Sigma,\xi]$ with
$a=h_{\dot q}$. Expanding, the terms $\int\beta\wedge d(\cdot)$ and
$\int d\xi\wedge d\mathsf D\xi$ vanish by Stokes and flatness, and
$\int d\xi\wedge\mathsf D\beta=\int\xi[a_\Sigma\wedge\beta]$. The
remaining terms are bounded by
$c|a|_\infty(|\beta|+\|d\xi\|)\|\xi\|$, and $\|\xi\|\le c\|d\xi\|$
because $H^0(\Sigma;\ad)=0$. The other components have
$H^0=H^1=0$ and are treated the same way. The trace $H^1(C)\to
L^2(\partial C)$ is bounded, and $|a|\le c|\dot q|$.
\end{proof}

\subsection{The Gårding estimate}\label{subsec:garding}

Fix symplectic frames along $u$ as follows: transport a symplectic frame
along $t\mapsto u(t,0)$ and then along each segment $r\mapsto u(t,r)$ by a
symplectic connection, and write $\nu=\Phi\xi$ with
$\xi\colon S\to\RR^{2m}$, $\xi_0=\xi(\cdot,0)$. Then $D_u\nu=\Phi(\partial_t\xi
+\tilde J\partial_r\xi+\tilde S\xi)$ with $\tilde J=\Phi^{-1}J\Phi$ and $\tilde S$ bounded by the $C^1$ norms of $u$
and $J$; the derivatives of $\tilde J$ contribute zeroth-order terms of the
same size. Let $e_i(t)$ be the harmonic representatives
of $\Phi(t,0)\varepsilon_i$. Under Convention~\ref{conv:orientation},
$\omega(e_i,e_j)=\omega_0(\varepsilon_i,\varepsilon_j)$, and
\begin{equation}\label{eq:seam-frame}
\int_\Sigma\langle\beta\wedge\mathsf D_t\beta\rangle=
\omega_0(\xi_0,\partial_t\xi_0)+R',\qquad |R'|\le c(K_u+Q)|\xi_0|^2,
\end{equation}
where $K_u$ bounds the connection coefficients of the frame. At $r=1$,
$|\omega_0(\xi_1,\partial_t\xi_1)|\le cK_u|\xi_1|^2$ because $K$ is
Lagrangian.

\begin{proposition}[uniform Gårding estimate]\label{prop:garding}
There are ${\mathfrak s}_1>0$ and $C$, depending on $(C,g_C,\mathbb E)$, on $Q$, $H$, on
the Lagrangian $K$ and on the $C^1$ norms of $u$ and $J$, but not on ${\mathfrak s}$, $\inf{\mathfrak s}$ or the
support of $z$, such that for $\sup{\mathfrak s}\le{\mathfrak s}_1$ and every smooth
compactly supported $z=(U,\nu)\in\mathcal D$,
\[
\|z\|_X\le C\bigl(\|D_{\mathfrak s} z\|_Y+\|\nu\|_{L^2(S)}\bigr).
\]
The same holds for the formal adjoint in the measure ${\mathfrak s}\,dt\,dV_C$,
whose gauge-fixing component is $-\mathsf D_t-{\mathfrak s}^{-1}\mathcal B_t-\varsigma+2\varsigma\Pi_{\rm f}$, on the
same domain.
\end{proposition}

\begin{proof}
\emph{The energy identity.} Let $\mathcal A_t=[\mathsf D_t,\mathcal B_t]$,
a zeroth-order operator with $|\mathcal A_tU|\le c|\dot q||U|$. By
Lemma~\ref{lem:green} (without the closedness of the traces),
\[
\tfrac{d}{dt}\langle U,\mathcal BU\rangle
=2\langle\mathsf DU,\mathcal BU\rangle+\langle U,\mathcal AU\rangle
+\int_{\partial C}\langle b\wedge\mathsf D_tb\rangle .
\]
For a gauge-fixing component $f=\mathsf DU-{\mathfrak s}^{-1}\mathcal BU+ZU$ with a pointwise
endomorphism $|Z_t|\le z_0H$ (here $Z=2\varsigma\Pi_{\rm f}$), expanding ${\mathfrak s}|f|^2$ and
integrating the total derivative gives
\begin{equation}\label{eq:gauge-identity}
\int{\mathfrak s}\|f\|^2
=\int\bigl({\mathfrak s}\|\mathsf DU+ZU\|^2+{\mathfrak s}^{-1}\|\mathcal BU\|^2\bigr)
-2\!\int\!\langle ZU,\mathcal BU\rangle+\int\!\langle U,\mathcal AU\rangle
+\int\!\!\int_{\partial C}\!\langle b\wedge\mathsf D_tb\rangle .
\end{equation}
No factor ${\mathfrak s}$ multiplies the last three terms, and no derivative of
${\mathfrak s}$ appears: the weight ${\mathfrak s}$ cancels the ${\mathfrak s}^{-1}$ of
$\mathcal B$ before any integration by parts in $t$.

For the strip, with $D_0=\partial_t+\tilde J\partial_r$,
\begin{equation}\label{eq:strip-identity}
\|D_0\xi\|^2=\|\partial_t\xi\|^2+\|\partial_r\xi\|^2
+\int\omega_0(\xi,\partial_t\xi)\big|_{r=1}dt
-\int\omega_0(\xi,\partial_t\xi)\big|_{r=0}dt ,
\end{equation}
from $-2\omega_0(\partial_t\xi,\partial_r\xi)=\partial_r\omega_0(\xi,\partial_t\xi)
-\partial_t\omega_0(\xi,\partial_r\xi)$, and
$\|D_u\nu\|^2\ge\|D_0\xi\|^2-2\|\tilde S\|\|D_0\xi\|\|\xi\|$.

\emph{Cancellation.} By Lemma~\ref{lem:seam} and~\eqref{eq:seam-frame},
the last term of~\eqref{eq:gauge-identity} is
$\int\omega_0(\xi_0,\partial_t\xi_0)\,dt$ up to an error bounded by
$cQ\int\|\mathcal BU\|^2+c(Q+K_u)\int(|\xi_0|^2+|\xi_1|^2)$; the coefficient $K_u$
multiplies only the trace term. Adding
\eqref{eq:gauge-identity} and~\eqref{eq:strip-identity}, the two terms
$\pm\int\omega_0(\xi_0,\partial_t\xi_0)$ cancel identically. This is the
only point where Convention~\ref{conv:orientation} is used.

\emph{Absorption.} Put $\mathcal E=\int\|\mathcal BU\|^2\le{\mathfrak s}_1\int
{\mathfrak s}^{-1}\|\mathcal BU\|^2$ and $\mathcal J=\int(|\xi_0|^2+|\xi_1|^2)$.
By Lemma~\ref{lem:tie}, $\|U\|$ and the slow coefficient $w$ are bounded
pointwise in $t$ by $\|\mathcal BU\|+|\nu_0|$, so
$|\langle U,\mathcal AU\rangle|\le cQ(\|\mathcal BU\|^2+|\nu_0|^2)$,
$2|\langle ZU,\mathcal BU\rangle|\le\tfrac14{\mathfrak s}^{-1}\|\mathcal BU\|^2
+cH^2{\mathfrak s}(\|\mathcal BU\|^2+|\nu_0|^2)$, and ${\mathfrak s}\|U\|^2\le
c{\mathfrak s}_1(\|\mathcal BU\|^2+|\nu_0|^2)$. The one-dimensional trace
inequality $|g(0)|^2\le\theta\|g'\|^2+(1+\theta^{-1})\|g\|^2$ on $[0,1]$
bounds $\mathcal J$ by a small multiple of $\|\partial_r\xi\|^2$ plus a
multiple of $\|\xi\|^2$. Choosing ${\mathfrak s}_1=c(1+Q+H)^{-2}$ gives
\[
\|z\|_X^2\le C\bigl(\|D_{\mathfrak s} z\|_Y^2+\|\nu\|_{L^2(S)}^2\bigr),\qquad
C=c(1+Q+K_u+H)^4 .
\]
\emph{Adjoint.} By Lemma~\ref{lem:green} and the bilinear form of
\eqref{eq:strip-identity}, $\langle D_{\mathfrak s} z,z'\rangle_Y-\langle
z,D_{\mathfrak s}^*z'\rangle_Y=\int\omega(\nu_0,\mu_0)-\int\omega(\nu_0,\mu_0)=0$
for $z,z'\in\mathcal D$, the term at $r=1$ vanishing because $K$ is
Lagrangian. So the adjoint has the same domain. The reflection
$t\mapsto-t$ carries $D_{\mathfrak s}^*$ to an operator of the same form with
$Z=\mp(2\varsigma\Pi_{\rm f}-\varsigma)$, and the same argument applies.
\end{proof}

\begin{remark}\label{rem:wrong-sign}
With the opposite sign in Convention~\ref{conv:orientation} the
estimate fails already for constant ${\mathfrak s}$ and constant $q$. For a
unit vector $v$ and $\zeta=v-iJv$, the pair
$b=\chi\,\mathrm{Re}(e^{i\tau t+{\mathfrak s}\tau(s-1)}\zeta)$ near the seam, where $s$ is the
normal coordinate of $C$ with the boundary at $s=1$,
$\nu=\chi\,\mathrm{Re}(e^{i\tau t-\tau r}\zeta)$ on the strip, with
cutoffs $\chi$, lies in the domain and decays away from the seam on both
sides. Its ratio $\|z\|_X^2/(\|D_{\mathfrak s} z\|_Y^2+\|\nu\|^2)$ grows like
$\tau^2$. This is a boundary layer at the seam, where a
$\bar\partial$-operator would be glued to a $\partial$-operator; it is
unrelated to the slow directions.
\end{remark}

\begin{remark}[why no logarithm appears]\label{rem:no-log}
The slow component $w(t)$ is never estimated through the gauge part. A
bound through $\int{\mathfrak s}(|\dot w|^2+|w|^2)$ would require a weighted
Hardy inequality, which for ${\mathfrak s}=\bar{\mathfrak s} e^{-H|t|}$ costs a factor
$\log(1/\inf{\mathfrak s})$. Here $w$ is controlled pointwise in $t$ by
$\nu_0$ through Lemma~\ref{lem:tie}, and $\nu_0$ by the unweighted strip
norm.
\end{remark}

\subsection{The uniform inverse at a regular limit}\label{subsec:uniform-inverse}

Assume now that $u$ is $J$-holomorphic, converges exponentially to
points $x_\pm\in P\pitchfork K$, and that $J$ is independent of $t$ near
the ends. Assume
$D_u\colon\{\nu\in H^1:\nu(\cdot,0)\in TP,\ \nu(\cdot,1)\in TK\}\to L^2$
is onto with kernel $\mathcal K$ of dimension $d$, and fix $d$ smooth
functionals $\Gamma$ of $\nu$ with compact support such that
$\Gamma|_{\mathcal K}$ is an isomorphism.

\begin{lemma}[ends]\label{lem:ends}
There are $T_0$, ${\mathfrak s}_1$, $C$ such that for $\sup{\mathfrak s}\le{\mathfrak s}_1$ and
$z\in\mathcal D$ supported in $\{|t|\ge T_0\}$,
$\|z\|_X\le C\|D_{\mathfrak s} z\|_Y$, and similarly for $D_{\mathfrak s}^*$.
\end{lemma}

\begin{proof}
For $|t|\ge T_0$ the Lagrangian frames $\Phi^{-1}TP$ and $\Phi^{-1}TK$
at $r=0,1$ are close to transverse limits, so
$\sup_r|\xi(t,r)|\le C(\|\partial_r\xi(t,\cdot)\|+\|\mathcal B_tU\|)$,
using Lemma~\ref{lem:tie} for the distance of $\nu_0$ from $TP$. Squaring and
integrating over $t$,
\[
\mathcal J+\|\xi\|^2\le C_P\left(\|\partial_r\xi\|^2+\int\|\mathcal B_tU\|^2dt\right),
\]
where $C_P$ is fixed by the transverse limiting pair. It need not be small.
Return to the energy identities \eqref{eq:gauge-identity} and
\eqref{eq:strip-identity}, before replacing the boundary terms by an
unweighted $L^2$ remainder. Put
\[
\mathcal E_0=\|\partial_t\xi\|^2+\|\partial_r\xi\|^2+
\int\bigl({\mathfrak s}\|\mathsf D_tU\|^2+
{\mathfrak s}^{-1}\|\mathcal B_tU\|^2\bigr).
\]
The coefficients involving $\dot q$, the derivatives of the strip, and
the moving Lagrangian frames tend to zero on the ends. Their contributions,
including the lower-order term of the strip operator, are at most
$o_{T_0}(1)(\mathcal E_0+\mathcal J+\|\xi\|^2)$ after Young's inequality.
The bounded scale derivative contributes at most a fixed multiple of
$H^2\int{\mathfrak s}\|U\|^2$, after absorbing one quarter of the
$\mathcal B_tU$ term. Lemma~\ref{lem:tie} and the preceding Poincar\'e bound give
\[
\int{\mathfrak s}\|U\|^2+
\mathcal J+\|\xi\|^2\le C(1+C_P)\mathcal E_0,
\qquad
\int{\mathfrak s}\|U\|^2\le C{\mathfrak s}_1(1+C_P)\mathcal E_0.
\]
Consequently
\[
\mathcal E_0\le C\|D_{\mathfrak s}z\|_Y^2+
\bigl(o_{T_0}(1)(1+C_P)+C H^2{\mathfrak s}_1(1+C_P)\bigr)\mathcal E_0.
\]
Choose $T_0$ first and then ${\mathfrak s}_1$ to absorb the last term.
The complementary spatial estimate and the same Poincar\'e bound recover
the full $X$ norm from $\mathcal E_0$. For the adjoint the extra scalar
endomorphism is bounded by $CH$ and is treated by the identical small-scale
estimate. This proves both assertions without a small Poincar\'e constant.
\end{proof}

Combining Proposition~\ref{prop:garding} on $[-T_0-1,T_0+1]$ with
Lemma~\ref{lem:ends} by a partition of unity in $t$ gives
\begin{equation}\label{eq:localized-garding}
\|z\|_X\le C\bigl(\|D_{\mathfrak s} z\|_Y+\|\nu\|_{L^2([-T_1,T_1]\times[0,1])}\bigr),
\end{equation}
with $C$, $T_1$, ${\mathfrak s}_1$ independent of ${\mathfrak s}$, and the same for
$D_{\mathfrak s}^*$.

\begin{proposition}[uniform inverse without logarithmic loss]\label{prop:uniform-inverse}
There are ${\mathfrak s}_1$, $C$ such that for every scale with
$\sup{\mathfrak s}\le{\mathfrak s}_1$ and $|{\mathfrak s}'/{\mathfrak s}|\le H$, the operator
$(D_{\mathfrak s},\Gamma)\colon X\to Y\oplus\RR^d$ is an isomorphism with
$\|(D_{\mathfrak s},\Gamma)^{-1}\|\le C$. Moreover, if $J_{\mathfrak s} k=((0,h_{k(t,0)}),k)$
for $k\in\mathcal K$ and $G_{\mathfrak s} k$ is the kernel element of $D_{\mathfrak s}$
with $\Gamma G_{\mathfrak s} k=\Gamma k$, then $\|G_{\mathfrak s} k-J_{\mathfrak s} k\|_X\le
C\sqrt{\sup{\mathfrak s}}\,\|k\|$.
\end{proposition}

\begin{proof}
\emph{Injectivity estimate.} Suppose $\|z_j\|_X=1$ and
$\|D_{{\mathfrak s}_j}z_j\|_Y+|\Gamma z_j|\to0$ with $\sup{\mathfrak s}_j\to0$. By
\eqref{eq:localized-garding}, $\liminf\|\nu_j\|_{L^2([-T_1,T_1]\times
[0,1])}>0$. A subsequence of $\nu_j$ converges weakly in $H^1$ and
strongly in $L^2_{\mathrm{loc}}$ to $\nu_\infty$ with $D_u\nu_\infty=0$,
$\nu_\infty(\cdot,1)\in TK$, and $\nu_\infty(\cdot,0)\in TP$, the last
because $\int\dist(\nu_j(t,0),TP)^2\le c\int\|\mathcal BU_j\|^2\le
c\sup{\mathfrak s}_j$ by Lemma~\ref{lem:tie}. Thus $\nu_\infty\in\mathcal K$ and
$\Gamma\nu_\infty=0$, so $\nu_\infty=0$, a contradiction. For the adjoint
the limit satisfies $D_u^*\mu_\infty=0$ with Lagrangian boundary values,
hence is orthogonal to the image of $D_u$, which is everything. The
scales ${\mathfrak s}_j$ play no role in the limit, which lives on the strip
alone: no convergent subsequence of scales and no lower bound on
${\mathfrak s}_j$ is used.

\emph{Surjectivity.} The same contradiction argument, applied to $D_{\mathfrak s}^*$
(whose limit object is an element of the kernel of $D_u^*$ with Lagrangian
boundary values, hence zero), gives $\|y\|_X\le C\|D^*_{\mathfrak s} y\|_Y$. Coercivity of
$\langle D^*_{\mathfrak s} y,D^*_{\mathfrak s} y'\rangle$ and the Lax--Milgram lemma then give, for $f\in Y$, a weak solution
$z=D^*_{\mathfrak s} y$ of $D_{\mathfrak s} z=f$ with $\|z\|_Y\le C\|f\|_Y$; write $R_{\mathfrak s}f=z$, a right
inverse of $D_{\mathfrak s}$. At fixed
${\mathfrak s}$ weak solutions are strong: the slice operator
$({\mathfrak s}^{-1}\mathcal B_t,-J\partial_r)$ with the domain of
Section~\ref{subsec:model} is symmetric by Lemma~\ref{lem:green} and the seam
cancellation. It is self-adjoint: an element of the domain of the adjoint is
a weak solution, with $L^2$ right-hand side, of the same mixed boundary
problem (vanishing normal component, closed traces, harmonic matching), and
the fixed-geometry regularity for exactly this mixed boundary condition
\cite[Theorem~4.17 and Remark~4.19]{DFL} places it in the domain. That theorem
is stated for the four-dimensional cylinder with one compression body; it
applies to a slice operator after extending an element $w$ of the adjoint domain
as $\beta(t)w$ on an interval times the slice, to slices with two gauge branches
after localizing near each seam and reflecting $\theta\mapsto\pi-\theta$ at the
collapsing seam, and the finite-rank harmonic matching does not affect
regularity. After a bounded trivialization of the $t$-dependent domain, regularity of weak
solutions is that of \cite{Robbin-Salamon}. Applying
\eqref{eq:localized-garding} to cut-off versions of $z$ gives $z\in X$
with $\|z\|_X\le C\|f\|_Y$.

\emph{Kernel.} $J_{\mathfrak s} k$ lies in the domain, its strip part is in
$\ker D_u$, and $D_{\mathfrak s} J_{\mathfrak s} k=((0,\mathsf D_th_{k(t,0)}),0)$ has
$Y$-norm at most $C\sqrt{\sup{\mathfrak s}}\|k\|$ since $k$ decays
exponentially. Hence $z_k=J_{\mathfrak s} k-R_{\mathfrak s} D_{\mathfrak s} J_{\mathfrak s} k$ lies in
the kernel, $\Gamma$ is onto on the kernel for small $\sup{\mathfrak s}$, and
$(D_{\mathfrak s},\Gamma)$ is an isomorphism. The estimate for $G_{\mathfrak s}
k-J_{\mathfrak s} k$ follows from the uniform bound.
\end{proof}

\subsection{The endpoint transition}\label{subsec:transition}

In this subsection we work with the scale ${\mathfrak s}$ of the collapsing branch of a
single slice; on the endpoint neck of $X_\eps$ it is ${\mathfrak s}=\eps e^s$. Near an
endpoint, with $L=\log(1/\eps)$, in the coordinate $\rho=L-s=\log(1/{\mathfrak s})$ (with $s$ as in Lemma~\ref{lem:conformal-models}(ii)), the
slice at scale ${\mathfrak s}=e^{-\rho}$
is $(C,g_C)\sqcup[0,\pi]\sqcup(-C,{\mathfrak s}^2g_C)$, the copy at $\theta=\pi$ carrying the
reversed orientation in the conventions of Lemma~\ref{lem:conformal-models}(ii). On the slice Hilbert space
$\mathcal H_{\mathfrak s}$ with norm $\|z\|_{\mathfrak s}^2=\|U_f\|^2+{\mathfrak s}\|U_c\|^2+\|\nu\|^2$, the
tangential operator is
\[
A_{\mathfrak s}(U_f,U_c,\nu)=(\mathcal B U_f,\ {\mathfrak s}^{-1}\mathcal BU_c,\
-J\partial_\theta\nu)
\]
on the ${\mathfrak s}$-independent domain $\mathcal W$ given by vanishing normal
components, closed traces and harmonic matching at both ends of $[0,\pi]$.
The quadratic form $\langle z,A_{\mathfrak s} z\rangle_{\mathfrak s}$ does not depend on ${\mathfrak s}$.
By Lemma~\ref{lem:green} and Convention~\ref{conv:orientation} the seam
terms cancel, and $A_{\mathfrak s}$ is symmetric. It is self-adjoint by the same
argument: with $I_{\mathfrak s}=\operatorname{diag}(1,{\mathfrak s},1)$, an element of the domain of
the adjoint of $I_{\mathfrak s} A_{\mathfrak s}$ is a weak solution of the mixed boundary problem,
hence in the domain by \cite[Theorem~4.17 and Remark~4.19]{DFL}.

\begin{proposition}\label{prop:transition-gap}
There are $\gamma>0$ and ${\mathfrak s}_0>0$ such that for ${\mathfrak s}\le{\mathfrak s}_0$ and
$q\in P$, $\ker A_{\mathfrak s}=\{k_v=(h_v,h_v,v):v\in T_qP\}$ and the spectrum
of $A_{\mathfrak s}$ meets $(-\gamma,\gamma)$ only in $0$. Moreover
$\|U_c-h_{\Pi\nu(\pi)}\|_{H^1}\le C\sqrt{\mathfrak s}\,\|A_{\mathfrak s} z\|_{\mathfrak s}$, where
$\Pi$ is the projection to $T_qP$.
\end{proposition}

\begin{proof}
\emph{Kernel.} If $A_{\mathfrak s} z=0$, then $\mathcal BU_f=\mathcal BU_c=0$, so by
Lemma~\ref{lem:hodge} $U_f=(0,h_{v_f})$ and $U_c=(0,h_{v_c})$; and $\nu$ is
constant, so the matching gives $v_f=\nu=v_c=:v$, with $v\in T_qP$ because
$h_v$ is defined only there.

\emph{The complementary operator.} The orthogonal complement of $T_qP$ in $T_q\Msf$ is
$J\,T_qP$, since $T_qP$ is Lagrangian and $J$ is compatible, so $1-\Pi$ projects onto
$J\,T_qP$. Let $\mathcal W^\perp_q$ be the space of $U=(p,b)$ on $C$ with $b(n)=0$,
closed traces, and harmonic part of $b|_\Sigma$ in $J\,T_qP$, and let
$\mathcal B^\perp$ be $\mathcal B$ with domain $\mathcal W^\perp_q$. By
Lemma~\ref{lem:green} it is symmetric, $J\,T_qP$ being Lagrangian. It satisfies
\begin{equation}\label{eq:complementary-perp}
\|U\|_{H^1(C)}\le c_1\|\mathcal BU\|\qquad(U\in\mathcal W^\perp_q):
\end{equation}
if $w\in T_qP$ is the bulk harmonic coefficient of $U$ and $\nu_0\in J\,T_qP$ the
harmonic part of its trace, Lemma~\ref{lem:tie} gives $|w|\le|\nu_0-w|\le
c\|\mathcal BU\|$, because $\nu_0\perp w$, and Lemma~\ref{lem:hodge}(ii) gives the rest.
Hence $\mathcal B^\perp$ is injective with closed range; it is self-adjoint by the local
regularity used above for $A_{\mathfrak s}$, so it is an isomorphism $\mathcal W^\perp_q\to
L^2(C)$ with inverse bounded by $c_1$ into $H^1$ (this is the complementary problem of
\eqref{eq:complementary}), and $\mathcal B^\perp-\mu$ is invertible with
$\|(\mathcal B^\perp-\mu)^{-1}\|_{L^2\to H^1}\le2c_1$ for $|\mu|\le1/(2c_1)$.

\emph{Elimination of the collapsed component.} Let $A_0$ be the tangential operator of
the mixed end, on $\mathcal H_0=L^2(C)\oplus L^2([0,\pi])$ with the domain $\mathcal W_0$
of pairs $(U_f,\nu)$ with $b_f(n)=0$, closed traces, harmonic part of $b_f|_\Sigma$ equal
to $\nu(0)$, and $\nu(\pi)\in T_qP$. It is self-adjoint and its spectrum meets
$(-2\gamma,2\gamma)$ only in $0$, uniformly in $q$ \cite[(4.20)--(4.24)]{DFL}; its kernel
is $\{(h_v,v)\colon v\in T_qP\}$. Let $|\lambda|=\gamma$, ${\mathfrak s}\gamma\le1/(2c_1)$, and
consider $(A_{\mathfrak s}-\lambda)z=f$ with $z\in\mathcal W$ and
$f=(f_f,f_c,f_\nu)\in\mathcal H_{\mathfrak s}$. Put $\tilde v=\Pi\nu(\pi)$ and $\tilde
U_c=U_c-h_{\tilde v}$. The harmonic part of the trace of $U_c$ is $\nu(\pi)$, so that of
$\tilde U_c$ is $(1-\Pi)\nu(\pi)\in J\,T_qP$ and $\tilde U_c\in\mathcal W^\perp_q$; since
$\mathcal Bh_{\tilde v}=0$, the collapsing equation $\mathcal BU_c-{\mathfrak s}\lambda U_c={\mathfrak s} f_c$
becomes $(\mathcal B^\perp-{\mathfrak s}\lambda)\tilde U_c={\mathfrak s} f_c+{\mathfrak s}\lambda h_{\tilde v}$. Thus
\[
\tilde U_c={\mathfrak s}(\mathcal B^\perp-{\mathfrak s}\lambda)^{-1}(f_c+\lambda h_{\tilde v}),\qquad
\|\tilde U_c\|_{H^1}\le C{\mathfrak s}(\|f_c\|+|\tilde v|),
\]
and the matching condition at $\theta=\pi$ becomes the perturbed Lagrangian condition
\begin{equation}\label{eq:perturbed-lagrangian}
(1-\Pi)\nu(\pi)={\mathfrak s}\lambda\,\Lambda_{{\mathfrak s}\lambda}\tilde v+{\mathfrak s}\,\mathsf m_f,\qquad
\Lambda_\mu v=\tau\bigl((\mathcal B^\perp-\mu)^{-1}h_v\bigr),\quad
\mathsf m_f=\tau\bigl((\mathcal B^\perp-{\mathfrak s}\lambda)^{-1}f_c\bigr),
\end{equation}
where $\tau$ is the harmonic part of the trace on $\Sigma$. Here $\Lambda_\mu\colon
T_qP\to J\,T_qP$ is bounded uniformly, and $|{\mathfrak s}\mathsf m_f|\le C{\mathfrak s}\|f_c\|\le
C\sqrt{\mathfrak s}\|f\|_{\mathfrak s}$ because ${\mathfrak s}\|f_c\|^2\le\|f\|_{\mathfrak s}^2$. Conversely, a pair
$(U_f,\nu)$ satisfying the fixed-branch and section equations and
\eqref{eq:perturbed-lagrangian}, completed by $U_c=h_{\tilde v}+\tilde U_c$, is a solution.

\emph{Reduction to $A_0$ and the Neumann series.} Fix a cutoff $\chi_\pi$ on $[0,\pi]$,
equal to $1$ near $\pi$ and to $0$ on $[0,\pi/2]$, and for $(U_f,\nu)$ satisfying
\eqref{eq:perturbed-lagrangian} put $\tilde\nu=\nu-\chi_\pi({\mathfrak s}\lambda\Lambda_{{\mathfrak s}\lambda}\tilde
v+{\mathfrak s}\mathsf m_f)$. Then $\tilde\nu(\pi)=\tilde v\in T_qP$ and $\Pi\tilde\nu(\pi)=\tilde v$, so
$y=(U_f,\tilde\nu)\in\mathcal W_0$, and the fixed-branch and section equations become
\[
(A_0-\lambda)y=(f_f,f_\nu)+\Theta_\lambda y+{\mathfrak s} F_\lambda f_c,
\]
with $\Theta_\lambda y=\bigl(0,(J\partial_\theta+\lambda)(\chi_\pi{\mathfrak s}\lambda
\Lambda_{{\mathfrak s}\lambda}\Pi\tilde\nu(\pi))\bigr)$ and ${\mathfrak s} F_\lambda f_c=\bigl(0,(J\partial_\theta
+\lambda)(\chi_\pi{\mathfrak s}\mathsf m_f)\bigr)$. Since the trace at $\pi$ is bounded by the graph
norm, $\|\Theta_\lambda y\|_{\mathcal H_0}\le C{\mathfrak s}\|y\|_{\mathcal W_0}$, and $\|{\mathfrak s}
F_\lambda f_c\|\le C\sqrt{\mathfrak s}\|f\|_{\mathfrak s}$. On $|\lambda|=\gamma$ the resolvent
$(A_0-\lambda)^{-1}\colon\mathcal H_0\to\mathcal W_0$ is bounded by $C_\gamma$ in the graph
norm, so $\|(A_0-\lambda)^{-1}\Theta_\lambda\|_{\mathcal W_0\to\mathcal W_0}\le\frac12$ for
${\mathfrak s}\le{\mathfrak s}_0$, and
\begin{equation}\label{eq:reduced-resolvent}
y=\sum_{j\ge0}\bigl((A_0-\lambda)^{-1}\Theta_\lambda\bigr)^j(A_0-\lambda)^{-1}
\bigl((f_f,f_\nu)+{\mathfrak s} F_\lambda f_c\bigr).
\end{equation}
Conversely \eqref{eq:reduced-resolvent}, followed by
$\nu=\tilde\nu+\chi_\pi({\mathfrak s}\lambda\Lambda_{{\mathfrak s}\lambda}\tilde v+{\mathfrak s}\mathsf m_f)$ and
$U_c=h_{\tilde v}+\tilde U_c$, defines a solution $z\in\mathcal W$; for $f=0$ it gives
$y=0$, then $\tilde v=0$, $\nu=0$ and $U_c=0$. Hence the circle $|\lambda|=\gamma$ lies in
the resolvent set of $A_{\mathfrak s}$. Let $\varpi(U_f,U_c,\nu)=(U_f,\nu)$ and
$\iota(U_f,\nu)=(U_f,0,\nu)$; then $\varpi\iota=1$, $\iota\varpi=1-P_c$, $\iota$ is an
isometry $\mathcal H_0\to\mathcal H_{\mathfrak s}$ and $\varpi=\iota^*$. From
\eqref{eq:reduced-resolvent} and the bounds above,
\[
\|\varpi(A_{\mathfrak s}-\lambda)^{-1}\iota-(A_0-\lambda)^{-1}\|\le C{\mathfrak s},\qquad
\|P_c(A_{\mathfrak s}-\lambda)^{-1}\|+\|(A_{\mathfrak s}-\lambda)^{-1}P_c\|\le C\sqrt{\mathfrak s}
\]
on $|\lambda|=\gamma$: for the first, $f=\iota g$ has $f_c=0$ and the series differs from its
first term by $O({\mathfrak s})$, as does $\nu$ from $\tilde\nu$; for the second, the collapsed
component of $(A_{\mathfrak s}-\lambda)^{-1}f$ has $\mathcal H_{\mathfrak s}$-norm
$\sqrt{\mathfrak s}\|U_c\|\le\sqrt{\mathfrak s}(C|\tilde v|+\|\tilde U_c\|)\le C\sqrt{\mathfrak s}\|f\|_{\mathfrak s}$, and a source
$f=P_cf$ has $\|f\|_{\mathfrak s}=\sqrt{\mathfrak s}\|f_c\|$, so its solution is $O({\mathfrak s}\|f_c\|)=O(\sqrt{\mathfrak s}\|f\|_{\mathfrak s})$.

\emph{Riesz projections.} Let $P_{\mathfrak s}$ and $P_0$ be the Riesz projections of $A_{\mathfrak s}$ and
$A_0$ for the circle $|\lambda|=\gamma$. Integrating the two estimates,
$\|\varpi P_{\mathfrak s}\iota-P_0\|\le C{\mathfrak s}$ and $\|P_cP_{\mathfrak s}\|+\|P_{\mathfrak s} P_c\|\le C\sqrt{\mathfrak s}$.
Since $P_{\mathfrak s}-(1-P_c)P_{\mathfrak s}(1-P_c)=P_cP_{\mathfrak s}+P_{\mathfrak s} P_c-P_cP_{\mathfrak s} P_c$,
\[
\|P_{\mathfrak s}-\iota P_0\varpi\|_{\mathcal H_{\mathfrak s}}\le\|P_cP_{\mathfrak s}\|+\|P_{\mathfrak s} P_c\|+\|P_cP_{\mathfrak s}
P_c\|+\|\iota(\varpi P_{\mathfrak s}\iota-P_0)\varpi\|\le C\sqrt{\mathfrak s} .
\]
The operator $\iota P_0\varpi$ is an orthogonal projection on $\mathcal H_{\mathfrak s}$ of rank
$\operatorname{rank}P_0=n$, and two projections at distance less than one have the same
rank; so $P_{\mathfrak s}$ has rank $n$ for ${\mathfrak s}$ small. Its range contains the $n$-dimensional
kernel found above. As $A_{\mathfrak s}$ is self-adjoint and $\pm\gamma$ are not in its spectrum,
$P_{\mathfrak s}$ is the spectral projection of $(-\gamma,\gamma)$, and the spectrum of $A_{\mathfrak s}$ in
$(-\gamma,\gamma)$ is $\{0\}$. Uniformity over $q$ holds because $P$ is compact and the
gap of $A_0$ is uniform \cite[(4.24)]{DFL}.

\emph{Spatial estimate.} Write $U_c=h_{\Pi\nu(\pi)}+\tilde U_c$ with $\tilde U_c\in\mathcal
W^\perp_q$, as above. By \eqref{eq:complementary-perp},
$\|\tilde U_c\|_{H^1}\le c_1\|\mathcal B\tilde U_c\|=c_1\|\mathcal BU_c\|\le
c_1\sqrt{\mathfrak s}\,\|A_{\mathfrak s} z\|_{\mathfrak s}$, since $\|A_{\mathfrak s} z\|_{\mathfrak s}^2\ge{\mathfrak s}^{-1}\|\mathcal BU_c\|^2$.
\end{proof}

Along the endpoint neck the scale varies, ${\mathfrak s}(\rho)={\mathfrak s}_{\rm b}e^{\pm\rho}$. The
linearized operator is $\mathcal D=\partial_\rho-A_{{\mathfrak s}(\rho)}+2(\log{\mathfrak s})'\Pi_{\rm f}$,
which is $\partial_\rho-A_{{\mathfrak s}(\rho)}-2\Pi_{\rm f}$ for ${\mathfrak s}={\mathfrak s}_{\rm b}e^{-\rho}$; here $\Pi_{\rm f}$
projects to the scalar of the collapsing component, acting on sections of the
bundle of Hilbert spaces $\mathcal H_{{\mathfrak s}(\rho)}$. Since the quadratic form of
$A_{\mathfrak s}$ is independent of ${\mathfrak s}$, the cross term in
$\|(\partial_\rho-A)z\|^2$ is $-2\int\langle z',Az\rangle_{\mathfrak s}=-\int\partial_\rho\langle
z,Az\rangle_{\mathfrak s}\,d\rho=0$ for compactly supported $z$, so
\[
\int\|(\partial_\rho-A)z\|_{\mathfrak s}^2\,d\rho=\int\bigl(\|z'\|_{\mathfrak s}^2+\|Az\|_{\mathfrak s}^2\bigr)d\rho ,
\]
with no derivative of ${\mathfrak s}$. The metric derivative enters only through
$\partial_\rho\|z\|_{\mathfrak s}^2=2\langle z',z\rangle_{\mathfrak s}+(\log{\mathfrak s})'\|P_cz\|_{\mathfrak s}^2$, where $P_c$ projects to
the collapsing component, and hence only in weighted estimates, where it is
absorbed by Proposition~\ref{prop:transition-gap}.

\begin{proposition}\label{prop:transition-inverse}
\begin{enumerate}[label=(\roman*)]
\item (Tail estimate.) Let the scale satisfy $|(\log{\mathfrak s})'|\le1$, and put
  $\|y\|_X^2=\int(\|y\|_{\mathfrak s}^2+\|\partial_\rho y\|_{\mathfrak s}^2+\|A_{\mathfrak s} y\|_{\mathfrak s}^2)d\rho$ and
  $\|y\|_Y^2=\int\|y\|_{\mathfrak s}^2d\rho$. For $0<|\delta|<\gamma/2$ there is
  ${\mathfrak s}_0(\delta)>0$ such that for $z$ supported where ${\mathfrak s}(\rho)\le{\mathfrak s}_0(\delta)$,
  \[
  \|e^{\delta\rho}z\|_X\le C_\delta\|e^{\delta\rho}\mathcal Dz\|_Y,\qquad
  \|e^{\delta\rho}z\|_X\le C_\delta\|e^{\delta\rho}\mathcal D^*z\|_Y,
  \]
  where $\mathcal D^*=-\partial_\rho-A_{\mathfrak s}-(\log{\mathfrak s})'P_c+2(\log{\mathfrak s})'\Pi_{\rm f}$ is the formal
  adjoint.
\item (Conserved pairing.) For $h$ in the kernel of the adjoint and $v\in T_qP$,
  $\langle k_v,h\rangle_{{\mathfrak s}(\rho)}$ is independent of $\rho$.
\item (Transition.) Let $\mathcal T$ be the \emph{complete flat transition}: the
  domain consisting of the rescaled endpoint core (the half-disk of radius
  $c_{\rm core}\eps$ about $z_\pm$, $c_{\rm core}$ as in Definition~\ref{def:insertion}(i), rescaled by $\eps^{-1}$, with its fixed conformal
  structure), continued by the endpoint neck to $\rho=\infty$, with the flat
  pair $(B_q,q)$ and the scale ${\mathfrak s}\to0$ along the neck. With the weight
  $e^{\delta\rho}$, $0<\delta<\gamma/2$, and fixed limit $0$, $\mathcal D$ is an isomorphism; with a
  free limit in $T_qP$ it is onto with kernel $\{k_v\}$ and index $n$. For
  $0<\delta\le\delta_0$ the fixed-limit inverse satisfies
  $\|\mathcal D^{-1}\|\le C/\delta$, uniformly in $q\in P$.
\end{enumerate}
\end{proposition}

\begin{proof}
(i) Let $y=e^{\delta\rho}z$, so that
$e^{\delta\rho}\mathcal Dz=(\partial_\rho-A_{\mathfrak s}-\delta+2(\log{\mathfrak s})'\Pi_{\rm f})y$. Since the
quadratic form of $A_{\mathfrak s}$ does not depend on ${\mathfrak s}$, and
$\partial_\rho\|y\|_{\mathfrak s}^2=2\langle\partial_\rho y,y\rangle_{\mathfrak s}+(\log{\mathfrak s})'\|P_cy\|_{\mathfrak s}^2$,
\[
\int\|(\partial_\rho-A_{\mathfrak s}-\delta)y\|_{\mathfrak s}^2\,d\rho
=\int\bigl(\|\partial_\rho y\|_{\mathfrak s}^2+\|(A_{\mathfrak s}+\delta)y\|_{\mathfrak s}^2
+\delta(\log{\mathfrak s})'\|P_cy\|_{\mathfrak s}^2\bigr)d\rho .
\]
Because $0<|\delta|<\gamma/2$, Proposition~\ref{prop:transition-gap} gives
$\|(A_{\mathfrak s}+\delta)y\|_{\mathfrak s}\ge c_\delta(\|y\|_{\mathfrak s}+\|A_{\mathfrak s} y\|_{\mathfrak s})$ with $c_\delta>0$.
Its spatial estimate, together with $\|h_{\Pi\nu(\pi)}\|\le C|\nu(\pi)|\le
C(\|y\|_{\mathfrak s}+\|A_{\mathfrak s} y\|_{\mathfrak s})$, gives $\|P_cy\|_{\mathfrak s}\le
C\sqrt{\mathfrak s}(\|y\|_{\mathfrak s}+\|A_{\mathfrak s} y\|_{\mathfrak s})$, and $\|\Pi_{\rm f}y\|_{\mathfrak s}\le C{\mathfrak s}\|A_{\mathfrak s} y\|_{\mathfrak s}$
by Lemma~\ref{lem:hodge}(i) and the Poincaré inequality. For ${\mathfrak s}\le{\mathfrak s}_0(\delta)$ the last term and the term $2(\log{\mathfrak s})'\Pi_{\rm f}$ are therefore absorbed. For $\mathcal D^*$ the conjugated
operator is $-\partial_\rho-(A_{\mathfrak s}-\delta)-(\log{\mathfrak s})'P_c+2(\log{\mathfrak s})'\Pi_{\rm f}$, the same
computation with $-\partial_\rho$ gives the identity with $A_{\mathfrak s}-\delta$ in place of
$A_{\mathfrak s}+\delta$, and the extra term $(\log{\mathfrak s})'P_cy$ is absorbed in the same way.

(ii) For $h$ with $\mathcal D^*h=-h'-Ah+P_ch-2\Pi_{\rm f}h=0$, differentiate
$\langle k_v,h\rangle_{\mathfrak s}$: the term $P_ch$ in the equation cancels the metric
derivative, and $Ak_v=\Pi_{\rm f}k_v=0$.

(iii) \emph{Kernel with fixed limit.} If $\mathcal Dz=0$ with decay, the linear
energy identity on the core and neck, whose seam terms cancel by
Convention~\ref{conv:orientation}, gives $z$ flat and constant, hence $0$.
\emph{Cokernel.} A cokernel element $h$ of the free operator pairs trivially
with $\mathcal D(\beta k_v)=\beta'k_v$ for a cutoff $\beta$, so by (ii) its stationary
limit is orthogonal to the kernel and $h$ decays; it then vanishes by the
argument of Proposition~\ref{prop:caps}. 

\emph{The bound.} We use the norms of (i) on the neck, weighted by $e^{\delta\rho}$,
together with the fixed $H^1$ and $L^2$ norms on the core, and write $X_\delta$,
$Y_\delta$. Fix $\bar\delta\in(0,\gamma/2)$. By what has been shown, $\mathcal D\colon
X_{\bar\delta}\to Y_{\bar\delta}$ with fixed limit $0$ is an isomorphism; write
$\mathcal D^{-1}_{\bar\delta}$ for its inverse, which is bounded uniformly in $q$.

\emph{Stationary splitting.} In the component coordinates the map $k\colon
T_qP\to\mathcal H_{\mathfrak s}$, $v\mapsto k_v$, does not depend on $\rho$. Let $P^{\rm st}_\rho$
be the $\mathcal H_{{\mathfrak s}(\rho)}$-orthogonal projection onto $\ker A_{{\mathfrak s}(\rho)}=k(T_qP)$,
let $P^{\rm n}_\rho=1-P^{\rm st}_\rho$, and write $P^{\rm st}_\rho y=k_{\mathsf a(y)}$. With
$I_{\mathfrak s}=\operatorname{diag}(1,{\mathfrak s},1)$ and the Gram matrix $\mathsf G_{\mathfrak s}=k^*I_{\mathfrak s} k$,
which is uniformly positive, $P^{\rm st}=k\mathsf G_{\mathfrak s}^{-1}k^*I_{\mathfrak s}$ and $\mathsf
a=\mathsf G_{\mathfrak s}^{-1}k^*I_{\mathfrak s}$. Since $\partial_\rho I_{\mathfrak s}=(\log{\mathfrak s})'I_{\mathfrak s} P_c$,
\begin{equation}\label{eq:projection-derivative}
\partial_\rho P^{\rm st}=(\log{\mathfrak s})'\,P^{\rm st}P_cP^{\rm n},\qquad
\|\partial_\rho P^{\rm st}\|\le C|(\log{\mathfrak s})'|\sqrt{\mathfrak s} ,
\end{equation}
the bound because the collapsed component of $k_v$ has $\mathcal H_{\mathfrak s}$-norm
$\sqrt{\mathfrak s}\|h_v\|$. Three facts give the structure of the equation: $\mathcal
Dk_a=k_{a'}$ for a function $a$ with values in $T_qP$, since $k$ does not depend on
$\rho$ and $A_{\mathfrak s} k_v=\Pi_{\rm f}k_v=0$; $P^{\rm st}A_{\mathfrak s}=0$ and $P^{\rm st}\Pi_{\rm f}=0$, since
$A_{\mathfrak s}$ is symmetric and annihilates $k_v$, and $k_v$ has no scalar component; and
$P^{\rm st}\partial_\rho y=-(\partial_\rho P^{\rm st})y$ when $P^{\rm st}y=0$.
Consequently, for $z=k_a+z^{\rm n}$ with $a=\mathsf a(z)$ and $z^{\rm n}=P^{\rm n}z$,
the equation $\mathcal Dz=f$ is equivalent to the triangular system
\begin{equation}\label{eq:triangular}
\mathcal D^{\rm n}z^{\rm n}=P^{\rm n}f,\qquad a'=\mathsf a(f)-\mathsf b_{\mathfrak s} z^{\rm n},
\end{equation}
where $\mathcal D^{\rm n}=P^{\rm n}\partial_\rho-A_{\mathfrak s}+2(\log{\mathfrak s})'\Pi_{\rm f}$ acts on
$P^{\rm n}$-valued fields and $\mathsf b_{\mathfrak s} y=-(\log{\mathfrak s})'\mathsf
G_{\mathfrak s}^{-1}k^*I_{\mathfrak s} P_cy$ satisfies $|\mathsf b_{\mathfrak s} y|\le
C|(\log{\mathfrak s})'|\sqrt{\mathfrak s}\,\|y\|_{\mathfrak s}$. The first equation does not involve $a$.

\emph{The normal equation.} Let $\mathcal Y$ be the complete line obtained by
continuing the neck to $\rho\to-\infty$ with the slice
$(C,g_C)\sqcup[0,\pi]\sqcup(-C,{\mathfrak s}_{\mathcal Y}^2g_C)$ and the flat pair $(B_q,q)$, where
${\mathfrak s}_{\mathcal Y}$ is smooth and nonincreasing, equals ${\mathfrak s}(\rho)$ for
$\rho\ge\rho_1+2$, is constant for $\rho\le\rho_1$, and satisfies
$|(\log{\mathfrak s}_{\mathcal Y})'|\le1$; here $\rho_1$ is fixed below, independently of
$\delta$ and $q$. All of the above holds on $\mathcal Y$ with ${\mathfrak s}_{\mathcal Y}$ in place
of ${\mathfrak s}$. For a compactly supported $P^{\rm n}$-valued field $y$ in the domain,
$P^{\rm n}\partial_\rho y=\partial_\rho y+(\partial_\rho P^{\rm st})y$. By the identity
displayed before the proposition, the gap ($\|y\|_{\mathfrak s}\le\gamma^{-1}\|A_{\mathfrak s} y\|_{\mathfrak s}$ for
$y$ orthogonal to $\ker A_{\mathfrak s}$), \eqref{eq:projection-derivative} and $\|\Pi_{\rm f}y\|_{\mathfrak s}\le
C{\mathfrak s}\|A_{\mathfrak s} y\|_{\mathfrak s}$,
\[
\|\mathcal D^{\rm n}y\|_Y\ge\bigl(\|\partial_\rho y\|_Y^2+\|A_{\mathfrak s} y\|_Y^2\bigr)^{1/2}
-C{\mathfrak s}(\rho_1)^{1/2}\|y\|_X\ge c\|y\|_X
\]
once ${\mathfrak s}(\rho_1)$ is small. Conjugation by $e^{\sigma\rho}$ replaces $\mathcal D^{\rm n}$ by
$\mathcal D^{\rm n}-\sigma$, so the same bound holds in the norms weighted by
$e^{\sigma\rho}$ for $|\sigma|\le\sigma_0:=c/2$, with a constant independent of $\sigma$.
The formal adjoint $\mathcal D^*$ maps $P^{\rm n}$-valued fields to $P^{\rm n}$-valued
fields: for $P^{\rm st}y=0$, $P^{\rm st}\mathcal D^*y=(\partial_\rho P^{\rm
st})y-(\log{\mathfrak s})'P^{\rm st}P_cy=0$ by \eqref{eq:projection-derivative}. The same
computation therefore bounds $\|y\|_X$ by the norm of the normal adjoint applied to $y$,
and, weak solutions being strong by the local regularity used above for $A_{\mathfrak s}$,
$\mathcal D^{\rm n}$ is an isomorphism between the $P^{\rm n}$-valued parts of the
weighted spaces on $\mathcal Y$:
\begin{equation}\label{eq:normal-estimate}
\|e^{\sigma\rho}y\|_{X(\mathcal Y)}\le C_{\rm n}\|e^{\sigma\rho}\mathcal D^{\rm
n}y\|_{Y(\mathcal Y)},\qquad|\sigma|\le\sigma_0,
\end{equation}
with $C_{\rm n}$ independent of $\sigma$ and $q$.

\emph{The stationary equation.} For $0<\delta\le\sigma_0$ and $f\in Y_\delta(\mathcal Y)$
let $z^{\rm n}=(\mathcal D^{\rm n})^{-1}P^{\rm n}f$, $g=\mathsf a(f)-\mathsf b_{\mathfrak s}
z^{\rm n}$ and $a(\rho)=-\int_\rho^\infty g$. By \eqref{eq:normal-estimate},
$\|e^{\delta\rho}g\|_{L^2}\le C\|f\|_{Y_\delta}$. Since
$e^{\delta\rho}a(\rho)=-\int_0^\infty e^{-\delta u}(e^{\delta\rho}g)(\rho+u)\,du$ is the
convolution of $e^{\delta\rho}g$ with a one-sided kernel of $L^1$ norm $\delta^{-1}$ and
$L^2$ norm $(2\delta)^{-1/2}$,
\begin{equation}\label{eq:stationary-integral}
\|e^{\delta\rho}a\|_{L^2}\le\delta^{-1}\|e^{\delta\rho}g\|_{L^2},\qquad
\sup_\rho e^{\delta\rho}|a(\rho)|\le(2\delta)^{-1/2}\|e^{\delta\rho}g\|_{L^2},\qquad a'=g .
\end{equation}
As $A_{\mathfrak s} k_a=0$, $\partial_\rho k_a=k_{a'}$ and $\|k_v\|_{\mathfrak s}\le C|v|$, the field
$S_\delta f=k_a+z^{\rm n}$ solves $\mathcal DS_\delta f=f$ on $\mathcal Y$, tends to zero,
and satisfies $\|S_\delta f\|_{X_\delta(\mathcal Y)}\le C\delta^{-1}\|f\|_{Y_\delta(\mathcal
Y)}$. The factor $\delta^{-1}$ arises only in the undifferentiated stationary
coefficient; the normal part is solved with no loss.

\emph{The core.} Fix $\rho_2\ge\rho_1+2$ and cutoffs $\chi_1,\chi_0$ of $\rho$ with
$\chi_1=0$ for $\rho\le\rho_2$, $\chi_1=1$ for $\rho\ge\rho_2+1$, $\chi_0=0$ for
$\rho\le\rho_2+2$ and $\chi_0=1$ for $\rho\ge\rho_2+3$; on the support of $\chi_1$ the
operators of $\mathcal T$ and of $\mathcal Y$ agree. For $0<\delta\le\delta_0:=\min(\sigma_0,
\bar\delta)$ and $f\in Y_\delta$, regard $\chi_0f$ as a source on $\mathcal Y$ and put
\[
z_1=\chi_1S_\delta(\chi_0f),\qquad r=f-\mathcal Dz_1=(1-\chi_0)f-\chi_1'S_\delta(\chi_0f).
\]
The source $r$ is supported in the fixed compact region $\{\rho\le\rho_2+3\}$ of
$\mathcal T$, where $e^{\delta\rho}$ and $e^{\bar\delta\rho}$ are comparable with constants
independent of $\delta$. By \eqref{eq:stationary-integral} the stationary part of
$S_\delta(\chi_0f)$ is bounded there by $C\delta^{-1/2}\|f\|_{Y_\delta}$, so
$\|r\|_{Y_{\bar\delta}}\le C\delta^{-1/2}\|f\|_{Y_\delta}$. Then $z=z_1+\mathcal
D^{-1}_{\bar\delta}r$ solves $\mathcal Dz=f$ with limit zero, and since $\|\cdot\|_{X_\delta}\le
C\|\cdot\|_{X_{\bar\delta}}$ for $\delta\le\bar\delta$,
\[
\|z\|_{X_\delta}\le C\delta^{-1}\|f\|_{Y_\delta},
\]
while on every fixed compact subset of $\mathcal T$ the field $z$ is bounded by
$C\delta^{-1/2}\|f\|_{Y_\delta}$. The fixed-limit kernel is zero at every positive weight,
so $z=\mathcal D^{-1}f$. The constants are uniform in $q$, because $P$ is compact and
$\gamma$, $C_{\rm n}$ and $\mathcal D^{-1}_{\bar\delta}$ are uniform. The bound is sharp:
$z_\delta=k_v\beta(\delta(\rho-\rho_0))$ for a fixed cutoff $\beta$ has
$\|z_\delta\|/\|\mathcal Dz_\delta\|\ge c/\delta$.

\end{proof}

\subsection{Transport along a path in $L(C)$}\label{subsec:transport}

Let $q(t)$ be a smooth path in $P=L(C)$ with bounded derivatives of every order;
in Section~\ref{sec:collapse-coverage} it is the path $q_R(t)$ of section boundary
values of the approximate solution $W_\eps$. Write
$\mathscr L_t=T_{q(t)}P$, $j_tv=(0,h_v)$ for the harmonic extension,
$\ell_t=(j_t^*j_t)^{-1}j_t^*$ and $\Pi_t=j_t\ell_t$, all bulk operators on
$L^2(C)$. A variation $U=(p,b)$, $p={\mathfrak s}\psi$, splits as $U=j_tv+U_\perp$ with
$v=\ell_tU$ and $U_\perp=(1-\Pi_t)U$. On the domain (vanishing normal component,
closed traces) the complementary estimate
\begin{equation}\label{eq:complementary}
\|U_\perp\|_{H^{r+1}(C)}\le C_r\|\mathcal B_tU\|_{H^r(C)}\qquad(r\ge0)
\end{equation}
holds: solve the complementary problem with harmonic trace in $J\mathscr L_t$, use its
elliptic regularity, and remove the finite harmonic kernel by the bulk
projection.

\begin{lemma}[transport]\label{lem:transport}
There is a connection $\nabla_t=\mathsf D_t+C_t$ along $q(t)$ on the bundle of
spatial domains, with $\|C_t\|_{H^r\to H^r}\le C_r$ and the same bounds for its
time derivatives, which preserves the vanishing of the normal component, the
closedness of the traces and the harmonic seam matching, and satisfies
\[
\nabla_t(j_tv)=j_t\nabla^{\mathscr L}_tv,\qquad\nabla^{\mathscr L}_t(\ell_tU)=\ell_t\nabla_tU,\qquad
[\nabla_t,\Pi_t]=0,\qquad[\nabla_t,\Pi_{\rm f}]=0 .
\]
Consequently the iterated commutators $\mathcal B^{[0]}_t=\mathcal B_t$,
$\mathcal B^{[j+1]}_t=[\nabla_t,\mathcal B^{[j]}_t]$ satisfy $\mathcal B^{[j]}_tj_t=0$
and $\|\mathcal B^{[j]}_tU\|_{H^r}\le C_{j,r}\|\mathcal B_tU\|_{H^r}$. The connection is
the derivative of transports $\widetilde T$ from the domain at $q_0$ to the domain at
$q$, defined for $q$ near $q_0$ and bounded on every $H^r(C)$, with
$\widetilde Tj_0=j_qK$, $\ell_q\widetilde T=K\ell_0$ and $\Pi_q\widetilde T=\widetilde T\Pi_0$,
where $K\colon T_{q_0}\Msf\to T_q\Msf$ is a smooth family of isomorphisms carrying
$\mathscr L_{q_0}$ to $\mathscr L_q$, used also on the section collar.
\end{lemma}

\begin{proof}
For $q$ near $q_0$, choose a smooth family of isomorphisms $K\colon T_{q_0}\Msf\to T_q\Msf$
carrying $\mathscr L_{q_0}$ to $\mathscr L_q$; domain elements have arbitrary harmonic
trace, and the same $K$ transports the section values on the collar, so the
harmonic seam matching is preserved. On every
parallel surface of the boundary collars, identify the exact, coexact and
harmonic summands of the Hodge decompositions for the two flat restrictions,
blended to the identity in the interior, leaving scalars and normal
components unchanged; this gives $T^0$, which identifies the domains and their
harmonic boundary values by $K$ and is bounded on every $H^r(C)$. Correct it on
the harmonic fields: $T^1=T^0+(j_qK-T^0j_0)\ell_0$, which changes no boundary
condition and satisfies $T^1j_0=j_qK$. For an interior cutoff $\chi$ close to one
the matrix $M_q=\ell_q\chi j_q$ is invertible; with $e_q=\chi j_qM_q^{-1}$, which has
zero boundary trace, put $\widetilde T=T^1+e_q(K\ell_0-\ell_qT^1)$. Then
$\widetilde Tj_0=j_qK$, $\ell_q\widetilde T=K\ell_0$ and
$\Pi_q\widetilde T=\widetilde T\Pi_0$. Differentiating these transports along
$q(t)$ defines $\nabla_t$, and the identities follow. Since $\mathcal B_tj_t=0$ and
$\nabla_t$ preserves harmonic extension, induction gives
$\mathcal B^{[j]}_tj_t=0$; the commutators have spatial order at most one, so
$\|\mathcal B^{[j]}_tU\|_{H^r}\le C\|U_\perp\|_{H^{r+1}}$, and~\eqref{eq:complementary}
gives the bound.
\end{proof}

\subsection{The inverse on the glued domain}\label{subsec:patching}

Let $V$ be a regular index-zero solution on $X_0$ with limits $q_\pm$ at the
new ends and $\ker D_V=0$ (on the extended weighted spaces), and let $W_\eps$ be the
approximate solution on $X_\eps$ of Section~\ref{sec:collapse-coverage}:
$V$ truncated at radius $R=\eps^{1/4}$ about each endpoint, the flat
transitions at $q_\pm$ inserted, and the flat lift of the section boundary
values over $I\times C$. Put $L=\log(1/\eps)$ and $\mathsf A=1+L$.

\emph{Norms.} On $X_\eps$ use the complete conformal representatives: on the
middle compression $dt^2+(\eps/r)^2g_C$ with $dt=dx/r$; on the endpoint necks
$ds^2+g_C$ and the collapsing branch at scale ${\mathfrak s}=\eps e^s$; on the rescaled
endpoint cores the fixed transition metric; on the section the corresponding
fixed and logarithmic strip metrics. On a small compression with $U=(p,b)$,
$p={\mathfrak s}\psi$, the norms $X_\eps$ and $Y_\eps$ are
\[
\int\bigl({\mathfrak s}|U|^2+{\mathfrak s}|\mathsf D_tU|^2+{\mathfrak s}^{-1}|\mathcal B_tU|^2\bigr)dt,
\qquad\int{\mathfrak s}|f|^2dt ,
\]
on the section and the fixed branch the $H^1$ and $L^2$ norms, on the cores
the fixed mixed $H^1$ and $L^2$ norms, patched by uniformly bounded cutoffs,
and on the original ends the fixed weights with the free limiting vectors. The
norm includes the $L^2$ part of a stationary field on every finite neck: for a
constant stationary vector on a region of length $O(L)$,
\begin{equation}\label{eq:stationary-cost}
\|k_v\|_{X_\eps(\text{region})}\le C\sqrt{\mathsf A}\,|v|.
\end{equation}

\emph{Components.} The two \emph{transitions} are the complete flat
transitions $\mathcal T_\pm$ of Proposition~\ref{prop:transition-inverse} at
$q_\pm$. The \emph{auxiliary exterior} is the following complete domain.

\begin{definition}[auxiliary exterior]\label{def:aux-exterior}
Let $r$ be a smooth positive function on $I$ equal to the distance to the
nearer endpoint where that distance is at most $x_0/2$, and let $t$ be the
complete coordinate on $I$ with $dt=dx/r$, where $x$ is arclength; near
$z_\pm$ the coordinate $t$ agrees, up to sign and an additive constant, with
$s=-\log r$ of Lemma~\ref{lem:conformal-models}(ii). Fix a smooth
$F\colon(0,\infty)\to(0,1]$ with $F(y)=y$ for $y\le\frac12$, $F(y)=1$ for
$y\ge2$ and $0\le yF'(y)/F(y)\le1$, and put
\[
{\mathfrak s}_E=\sqrt\eps\,F\bigl(\sqrt\eps/r\bigr).
\]
The \emph{auxiliary exterior} $E_\eps$ is $X_0$ with the gauge piece $I\times C$,
carrying the metric $dt^2+{\mathfrak s}_E^2g_C$, attached along $I$: the Lagrangian
condition $L(C)$ on $I$ is replaced by the matching condition of
Definition~\ref{def:admissible-config}. Its configuration $W^E_\eps$ is the
truncation $V_R$ of $V$ at $R=\eps^{1/4}$ of
Section~\ref{sec:collapse-coverage} on $X_0$, and on $I\times C$ the flat lift
$B_{q_R(t)}$ of the boundary values $q_R$ of $V_R$ on $I$, with temporal
component in the covariant Neumann slice, $\partial_tB-d_B\Psi=h_{\dot q_R}$.
The operator $D_E$ is the linearized mixed operator at $W^E_\eps$, with the gauge-fixing component of the metric $dt^2+{\mathfrak s}_E^2g_C$ on the compression and of the complete
metrics of $X_0$ elsewhere, as for $D_\eps$ in the conformal representative above, on the
domain given along $I$ by the conditions of
Section~\ref{subsec:model} (vanishing normal component, closed traces, harmonic
part of $b|_\Sigma$ equal to the section value) and by the original conditions
elsewhere, on the extended weighted spaces at the original mixed ends and at the two new ends.
\end{definition}

The scale satisfies $\sup{\mathfrak s}_E\le\sqrt\eps$ and $|\partial_t\log{\mathfrak s}_E|\le H$
with $H$ independent of $\eps$; moreover ${\mathfrak s}_E=\eps/r$ where $r\ge2\sqrt\eps$
and ${\mathfrak s}_E=\sqrt\eps$ where $r\le\sqrt\eps/2$. Near $z_\pm$, in the coordinate
$s$, the domain $E_\eps$ is a half-cylinder with slice
$(C,g_C)\sqcup[0,\pi]\sqcup(-C,{\mathfrak s}_E(s)^2g_C)$, where ${\mathfrak s}_E(s)=\eps e^s$ for
$s\le L/2-\log2$ and ${\mathfrak s}_E(s)=\sqrt\eps$ for $s\ge L/2+\log2$. There is a
fixed $c_R$ such that $W^E_\eps$ coincides with $V$ and its flat lift for
$s\le L/4-c_R$, and with the flat pair $(B_{q_\pm},q_\pm)$ on all three branches
for $s\ge L/4+c_R$. On the region $\{r\ge2\sqrt\eps\}$, which contains the whole
exterior side of $s=L/2-\log2$ and in particular $J_\pm$, the domain, the
configuration and the operator coincide with $X_\eps$, $W_\eps$ and $D_\eps$ in the
conformal representative above.

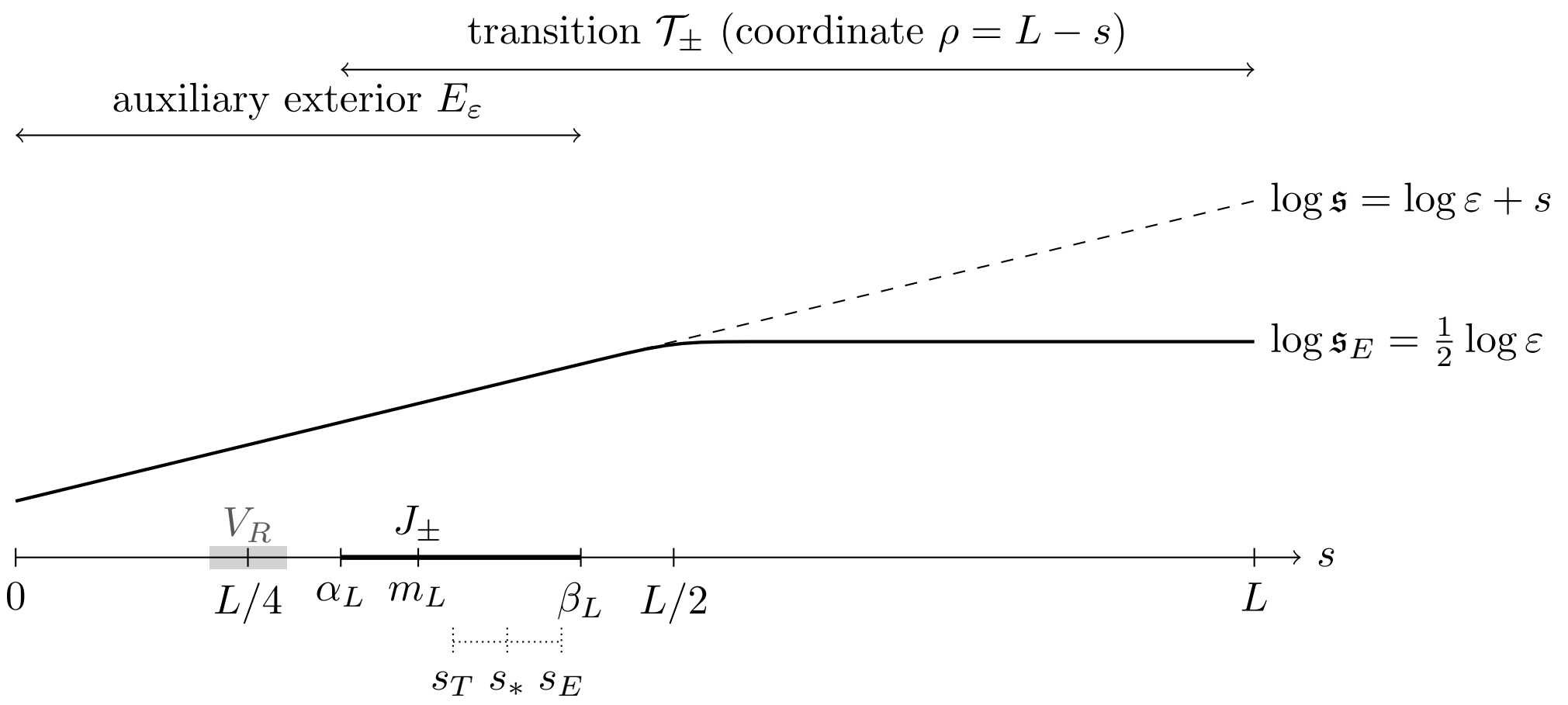
\begin{figure}[ht]
\centering
\begin{tikzpicture}[x=0.66cm,y=0.8cm,font=\small]
\fill[gray!35] (2.5,-0.12) rectangle (3.5,0.12);
\draw[->] (0,0) -- (16.6,0) node[right] {$s$};
\draw (0,0.1) -- (0,-0.1) node[below] {$0$};
\draw (16,0.1) -- (16,-0.1) node[below] {$L$};
\draw[very thick] (4.2,0) -- (7.3,0);
\foreach \x/\l in {3/{L/4},4.2/{\alpha_L},5.2/{m_L},7.3/{\beta_L},8.5/{L/2}}
  \draw (\x,0.1) -- (\x,-0.1) node[below] {$\l$};
\node at (5.2,0.35) {$J_\pm$};
\node[gray!70!black] at (3,0.35) {$V_R$};
\foreach \x/\l in {5.65/{s_T},6.35/{s_*},7.05/{s_E}}
  \draw[densely dotted] (\x,-0.75) -- (\x,-1.05) node[below] {$\l$};
\draw[densely dotted] (5.65,-0.9) -- (7.05,-0.9);
\draw[dashed] (0,0.6) -- (16,3.8);
\draw[thick] (0,0.6) -- (7.3,2.06) .. controls (8.5,2.3) .. (9.7,2.3) -- (16,2.3);
\node[right] at (16,3.8) {$\log{\mathfrak s}=\log\eps+s$};
\node[right] at (16,2.3) {$\log{\mathfrak s}_E=\tfrac12\log\eps$};
\draw[<->] (0,4.5) -- (7.3,4.5) node[midway,above] {auxiliary exterior $E_\eps$};
\draw[<->] (4.2,5.2) -- (16,5.2) node[midway,above] {transition $\mathcal T_\pm$ (coordinate $\rho=L-s$)};
\end{tikzpicture}
\caption{An endpoint neck in the exterior coordinate $s\in[0,L]$, $L=\log(1/\eps)$, with
$r=e^{-s}$ (not to scale). The dashed line is $\log{\mathfrak s}$ for the scale
${\mathfrak s}=\eps e^{s}$ of $X_\eps$; the solid line is the scale ${\mathfrak s}_E$ of the auxiliary
exterior, equal to ${\mathfrak s}$ for $s\le L/2-\log2$ and to $\sqrt\eps$ for $s\ge L/2+\log2$.
The shaded interval about $L/4$ ($r=\eps^{1/4}$) is where the truncation $V_R$ of the limit
is joined to the flat pair. The two components are used over the arrows and overlap the
glued domain on $J_\pm=[\alpha_L,\beta_L]$, with joining point $m_L$
(Section~\ref{subsec:patching}). The dotted marks are the cuts $s_T<s_E$ of the proof of
Proposition~\ref{prop:coverage}, at distance $C_*+D$ and $C_*-D$ from $L/2$, and their common
midpoint $s_*=L/2-C_*$; for large $L$ they lie between $m_L$ and $\beta_L$.}
\label{fig:neck}
\end{figure}

On the new ends of the components, which are not present on $X_\eps$, use the weight $\delta_L=c_0/(1+L)$,
bounded above and below on the parts of the components used in the gluing. Each transition overlaps the
actual domain on $J_\pm=[\alpha_L,\beta_L]$ with $\alpha_L=L/4+b_0$, with a fixed $b_0\ge c_R+1$ so that $W^E_\eps$ is flat on $J_\pm$,
$\beta_L=L/2-\log2$, of length $\mathsf T\sim L/4$, on which the three operators agree
under $\rho=L-s$ (Figure~\ref{fig:neck}).

\emph{Norms and asymptotic values.} Let $\phi_E$ be a smooth function on $E_\eps$, equal to
$\delta_*s$ on the new ends, to $\delta_as$ on the original end $a$ (the fixed
weights of~\ref{hyp:X3}) and to $0$ on a fixed compact core, and depending on the
compression only through $t$. Let $\mathbb V$ be the direct sum of
$T_{q_-}P\oplus T_{q_+}P$ and of the tangent spaces of the critical families at
the original mixed ends. Let $\mathfrak k_0(v)$, $v\in\mathbb V$, be the cut-off stationary fields of $D_V$ of Proposition~\ref{prop:fredholm}: a fixed cutoff of the end
coordinate times the stationary field of the end, in the end chart of $V$; at
the new ends this is $(h_v,v)$ on the slice $C\sqcup[0,\pi]$. Let
$\mathfrak k_0^\eps(v)$ be the same fields built with the end charts of $V_R$;
they coincide with $\mathfrak k_0(v)$ where $s\le L/4-c_R$ and differ from it by
$O(e^{-\beta s}|v|)$ elsewhere. On $E_\eps$ put $\mathfrak k_E(v)=\mathfrak
k^\eps_0(v)$ off the compression and $\mathfrak k_E(v)=(0,h_{\nu(t)})$ on
$I\times C$, where $\nu(t)\in T_{q_R(t)}P$ is the section component of $\mathfrak
k^\eps_0(v)$ on $I$. Then $\mathfrak k_E(v)$ lies in the domain of $D_E$. Every
element of $X_E$ is $z=\bar z+\mathfrak k_E(v)$ with $\bar z$ decaying, and
\begin{gather*}
\|z\|_{X_E}^2=\|e^{\phi_E}\bar z\|^2_{H^1}
+\int_Ie^{2\phi_E}\bigl({\mathfrak s}_E\|\bar U_c\|^2+{\mathfrak s}_E\|\mathsf D_t\bar U_c\|^2
+{\mathfrak s}_E^{-1}\|\mathcal B_t\bar U_c\|^2\bigr)dt+|v|^2,\\
\|f\|_{Y_E}^2=\|e^{\phi_E}f\|^2_{L^2}+\int_Ie^{2\phi_E}{\mathfrak s}_E\|f_c\|^2dt,
\end{gather*}
where the first terms are the norms of the section and of the gauge fields of
$X_0$, and the compression integrals use the fixed metric $g_C$ on $C$. On $J_\pm$
these are the norms $X_\eps$, $Y_\eps$ of this section multiplied by the weight.
Since $\mathcal B_t$ and $\Pi_{\rm f}$ annihilate harmonic fields, the compression
component of $D_E\mathfrak k_E(v)$ is $\mathsf D_t(0,h_{\nu(t)})$; as
${\mathfrak s}_E\le\sqrt\eps$, and $\partial_t\nu$ and $\dot q_R$ decay at the rate
$\beta>\delta_*$ along the ends,
\begin{equation}\label{eq:column-error}
\|(D_E\mathfrak k_E(v))_c\|_{Y_E}\le C\eps^{1/4}|v|,\qquad
\|D_E\mathfrak k_E(v)\|_{Y_E}\le C|v|,
\end{equation}
and $e^{\phi_E}(\mathfrak k_0^\eps(v)-\mathfrak k_0(v))$,
$e^{\phi_E}(D_{V_R}\mathfrak k_0^\eps(v)-D_V\mathfrak k_0(v))$ have $H^1$ and $L^2$
norms at most $C\eps^{(\beta-\delta_*)/4}|v|$.

\begin{proposition}[uniform inverse on the auxiliary exterior]\label{prop:exterior-inverse}
Fix $\delta_*$ with $0<\delta_*<\min(\beta,\gamma/2)$, where $\beta$ is the decay
rate of $V$ and $\gamma$ the gap of Proposition~\ref{prop:transition-gap}. There
are $\eps_1>0$ and $C$ such that for $0<\eps<\eps_1$ the operator
$D_E\colon X_E\to Y_E$ is an isomorphism and
\[
\|z\|_{X_E}\le C\|D_Ez\|_{Y_E}\qquad(z\in X_E).
\]
In particular the limiting values $v_\pm$ of $z$ at the new ends satisfy
$|v_-|+|v_+|\le C\|D_Ez\|_{Y_E}$.
\end{proposition}

At the weight $\delta_*$, which lies below the gaps of the new mixed ends of $X_0$
and below $\beta$, the operator $D_V$ on the extended weighted spaces is an isomorphism by
Hypothesis~\ref{hyp:collapse}\ref{hyp:X4}, $V$ being regular of index zero.

\begin{proof}
\emph{Step 1: a graph estimate with compact remainder.} Fix $s_1$ large, to be
chosen in (b) independently of $\eps$, and put $r_1=e^{-s_1}$. We first derive
a consequence of Lemma~\ref{lem:tie}: at every $t$,
${\mathfrak s}_E\|U_c\|^2\le c\,{\mathfrak s}_E\bigl(\|\mathcal B_tU_c\|^2+|\nu|_I|^2\bigr)\le
c\sup{\mathfrak s}_E\bigl({\mathfrak s}_E^{-1}\|\mathcal B_tU_c\|^2+|\nu|_I|^2\bigr)$, and the
trace inequality on a collar of $I$ in the section gives
\begin{equation}\label{eq:compression-small}
\int_Ie^{2\phi_E}{\mathfrak s}_E\|\bar U_c\|^2dt\le C\eps^{1/2}\|\bar z\|^2_{X_E}.
\end{equation}

(a) \emph{The middle of $I$.} On $\{r\ge r_1/4\}$ use conformal coordinates
$t'+i\theta'$ on a collar of $I$ in $S_0$ with $I=\{\theta'=0\}$, for instance
$t'+i\theta'=\log\bigl((z-z_-)/(z_+-z)\bigr)$ in the half-plane chart, and let
$\theta_1>0$ be such that $\{r\ge r_1/4,\ 0\le\theta'\le\theta_1\}$ lies in the
collar. There $dt'/dt$ is bounded above and below with all its derivatives,
independently of $\eps$. In the time coordinate $t'$ the compression metric is
$(dt/dt')^2(dt'^2+{\mathfrak s}'^2g_C)$ with ${\mathfrak s}'={\mathfrak s}_E\,dt'/dt$. The anti-self-duality
equation is conformally invariant, and in dimension four the gauge-fixing component of a metric changed by a conformal factor depending on $t'$ alone differs from the original one by a
positive factor and by a zeroth-order term in the temporal component. Hence, on fields
supported in $\{r\ge r_1/2\}$ whose section component vanishes for
$\theta'\ge\theta_1/2$, $\Upsilon D_E=D_{{\mathfrak s}'}+2h''\Pi_{\rm f}$, where $\Upsilon$ multiplies the gauge-fixing component by a positive function of $t'$ bounded above and below with its derivatives,
$D_{{\mathfrak s}'}$ is the operator of Section~\ref{subsec:model} with time $t'$, the strip
$\RR\times[0,\theta_1]$, the path $q_R$ and the section of $V_R$, and $h''$ is bounded
independently of $\eps$; a term $2h''\Pi_{\rm f}$ is allowed in the proof of
Proposition~\ref{prop:garding} as part of the endomorphism $Z$. Here
$\sup{\mathfrak s}'\le C\eps/r_1$, while $|\partial_{t'}\log{\mathfrak s}'|$,
$\sup|\partial_{t'}q_R|$ and the $C^1$ norms of the section are bounded
independently of $\eps$. The Lagrangian condition at the far edge of the strip
enters the proof of Proposition~\ref{prop:garding} only through the boundary term
$\omega_0(\xi_1,\partial_t\xi_1)$ and the trace term $|\xi_1|^2$, which vanish for
these fields, and the proof is unchanged for a strip of fixed width $\theta_1$.
The norms of Section~\ref{subsec:model} and those of $X_E$, $Y_E$ are uniformly
equivalent on this region. Hence, for such fields and $\eps$ small,
\[
\|\bar z\|_{X_E}\le C\bigl(\|D_E\bar z\|_{Y_E}+\|\nu\|_{L^2(\{r\ge r_1/2,\
\theta'\le\theta_1/2\})}\bigr).
\]

(b) \emph{The new ends.} For $s\ge s_1-1$ the configuration $W^E_\eps$ is within
$Ce^{-\beta s}$ of the flat pair $(B_{q_\pm},q_\pm)$ in every $C^k$ norm, with $C$
independent of $\eps$, by the exponential convergence of $V$
(Proposition~\ref{prop:weak-limit}(iii)) and the construction of $V_R$. Identify
the slice domains of $W^E_\eps$ with those of the flat pair by the map equal on the
collapsing branch to the transport $\widetilde T$ from $q_\pm$ to $q_R(s)$
of Lemma~\ref{lem:transport}, and on the fixed branch
and the section to an end chart of $V_R$ compatible with the identification $K$
of $\widetilde T$ at the seam; since $\widetilde T$ carries $h_v$ to $h_{Kv}$ and
harmonic boundary values by $K$, this map preserves the domain. In this
identification $D_E=\mathcal D_\pm+\mathfrak E$, where $\mathcal D_\pm$ is the neck
operator of Section~\ref{subsec:transition} at the flat pair with scale
${\mathfrak s}_E(s)$, written in the coordinate $s$ as on $J_\pm$, and
$\|\mathfrak Ey\|_{Y_E}\le Ce^{-\beta s_1}\|y\|_{X_E}$ for $y$ supported in
$\{s\ge s_1-1\}$. Indeed the only coefficient of $\mathfrak E$ carrying ${\mathfrak s}_E^{-1}$
is ${\mathfrak s}_E^{-1}(\mathcal B_{q_R(s)}\widetilde T-\widetilde T\mathcal B_{q_\pm})$; it
vanishes on harmonic fields, and it is a first-order operator whose coefficients
are $O(e^{-\beta s})$, so by Lemma~\ref{lem:hodge}(ii) it is bounded by
$Ce^{-\beta s}{\mathfrak s}_E^{-1}\|\mathcal B_{q_\pm}U_c\|$. By
Proposition~\ref{prop:transition-inverse}(i), which uses of the scale only
$|\partial\log{\mathfrak s}|\le1$ and a small upper bound, and applies here because
$0\le\partial_s\log{\mathfrak s}_E\le1$ and ${\mathfrak s}_E\le\sqrt\eps$, we have
$\|e^{\delta_*s}y\|_X\le C\|e^{\delta_*s}\mathcal D_\pm y\|_Y$ for decaying $y$
supported in $\{s\ge s_1-1\}$, with $C$ independent of $\eps$; by
Lemma~\ref{lem:hodge} these norms are uniformly equivalent to $X_E$ and $Y_E$ on
the ends. Choosing $s_1$ large, independently of $\eps$, to absorb $\mathfrak E$,
\[
\|\bar z\|_{X_E}\le C\|D_E\bar z\|_{Y_E}\qquad\text{for decaying $\bar z$ supported in
$\{s\ge s_1-1\}$.}
\]

(c) \emph{The remaining regions.} On the original ends, beyond a fixed start, the
weighted end estimates of \cite[Section~4.2.1]{DFL} at the fixed weights give
$\|\bar z\|_{X_E}\le C\|D_E\bar z\|_{Y_E}$ for decaying $\bar z$ supported there; on the gauge
ends these are the estimates for the trajectories of $\mathrm{CS}_{\mathfrak h}+f$, whose
Hessian at the critical points is that of $\mathrm{CS}_{\mathfrak h}$ because $f$ vanishes near
the critical set. The
rest of $E_\eps$ is a fixed compact subset of $X_0$ at positive distance from $I$
and from all ends, on which $W^E_\eps=V$ for $\eps$ small; there the local mixed
estimates \cite[Theorem~4.17 and Remark~4.19]{DFL} give
$\|\bar z\|_{H^1}\le C(\|D_E\bar z\|_{L^2}+\|\bar z\|_{L^2})$ for $\bar z$ supported in it.

(d) \emph{Patching.} Choose a finite partition $\sum_\alpha\chi_\alpha^2=1$
subordinate to the regions of (a)--(c), by cutoffs that are functions of $t$ alone
on the compression, take the same value on the two sides of the seam over $I$ (on
the section the cutoff of (a) is $\chi(t')\psi(\theta')$ with $\psi=1$ near $0$ and
$\psi=0$ for $\theta'\ge\theta_1/2$), and have gradients supported in a fixed
compact set. These cutoffs preserve the domain of $D_E$, whose conditions are linear
and pointwise in $t$, and $[D_E,\chi_\alpha]$ is multiplication by derivatives of
$\chi_\alpha$. Its compression component has $Y_E$-norm at most
$C(\int_Ie^{2\phi_E}{\mathfrak s}_E\|\bar U_c\|^2dt)^{1/2}\le C\eps^{1/4}\|\bar z\|_{X_E}$ by
\eqref{eq:compression-small}; its other components are bounded by the $L^2$ norm of
$\bar z$ over a fixed compact set $\Omega_0\subset\overline{X_0}$ containing these
gradients and the section region of (a). Squaring and summing the estimates of
(a)--(c) for $\chi_\alpha\bar z$, and absorbing the compression term for $\eps$
small, gives
\begin{equation}\label{eq:exterior-garding}
\|\bar z\|_{X_E}\le C\bigl(\|D_E\bar z\|_{Y_E}+\|\bar z\|_{L^2(\Omega_0)}\bigr),
\end{equation}
where $\|\cdot\|_{L^2(\Omega_0)}$ is the $L^2$ norm of the section and of the gauge
fields of $X_0$ over $\Omega_0$, and $C$ and $\Omega_0$ do not depend on $\eps$.
The same argument proves \eqref{eq:exterior-garding} for the formal adjoint
$\widehat D_E^*$ of the conjugated operator $\widehat D_E=e^{\phi_E}D_Ee^{-\phi_E}$
in the unweighted spaces: Proposition~\ref{prop:garding} contains its adjoint,
Proposition~\ref{prop:transition-inverse}(i) contains $\mathcal D^*$ with weights of
either sign, and the fixed estimates of (c) hold for adjoints.

\emph{Step 2: the injectivity estimate.} We claim that
$\|z\|_{X_E}\le C\|D_Ez\|_{Y_E}$ for $\eps$ small. Otherwise there are $\eps_i\to0$
and $z_i=\bar z_i+\mathfrak k_{E}(v_i)$ on $E_{\eps_i}$ with $\|z_i\|_{X_E}=1$ and
$\|D_Ez_i\|_{Y_E}\to0$. After passing to a subsequence, $v_i\to v$, and the
restrictions of $e^{\phi_E}\bar z_i$ to the section and to the gauge pieces of $X_0$,
which form a fixed space, converge weakly in $H^1$ and strongly in $L^2$ on compact
subsets of $\overline{X_0}$, with traces converging strongly in $L^2$ on compact
subarcs of $I$; let $\bar z_0$ be the limit and $z_0=\bar z_0+\mathfrak k_0(v)$.
Lemma~\ref{lem:tie} gives $\dist(\nu_i(t),T_{q_{R_i}(t)}P)\le c\|\mathcal
B_tU_{i,c}(t)\|$ at every $t$, hence
$\int_Ie^{2\phi_E}\dist(\nu_i,T_{q_{R_i}}P)^2dt\le c\sup{\mathfrak s}_E\to0$; since
$q_{R_i}\to q_V$ uniformly on compact subarcs, the section component of $z_0$ has
its trace on $I$ in $T_{q_V}L(C)$. The other boundary conditions are closed linear
conditions on traces and pass to the limit. On every compact subset of $X_0$ the
coefficients of $D_E$ off the compression equal those of $D_V$ for large $i$, so
$D_Vz_0=0$ in the interior and up to the boundary. By Fatou's lemma, the
convergence of $\mathfrak k^{\eps_i}_0(v_i)$ to $\mathfrak k_0(v)$ and
\eqref{eq:column-error}, the field $z_0$ lies in the domain of $D_V$ on the extended weighted space at weight $\delta_*$, with limiting vector $v$. Since this operator is injective,
$z_0=0$ and $v=0$. Then $\|\bar z_i\|_{L^2(\Omega_0)}\to0$, and
$\|D_E\bar z_i\|_{Y_E}\le\|D_Ez_i\|_{Y_E}+C|v_i|\to0$ by \eqref{eq:column-error}, so
\eqref{eq:exterior-garding} gives $\|\bar z_i\|_{X_E}\to0$; with $|v_i|\to0$ this
contradicts $\|z_i\|_{X_E}=1$.

\emph{Step 3: surjectivity.} By Step~2, $D_E$ is injective with closed range for
$\eps$ small. In the conjugated coordinates $\hat z=e^{\phi_E}\bar z$ the operator is
$\widehat D_E$ on the decaying part, and the asymptotic values act by $\widehat
Fv=e^{\phi_E}D_E\mathfrak k_E(v)$. Suppose that for $\eps_i\to0$ there are $f_i$ of
unit norm orthogonal to the range. Then $\langle\widehat D_E\hat z,f_i\rangle=0$
for all decaying $\hat z$ in the domain and $\langle\widehat Fv,f_i\rangle=0$ for all
$v\in\mathbb V$. The first identity says that $f_i$ is a weak solution of
$\widehat D_E^*f_i=0$. At the fixed scale $\eps_i$ the local regularity of
\cite[Theorem~4.17 and Remark~4.19]{DFL}, applied with cutoffs, places $f_i$ in the
domain, and \eqref{eq:exterior-garding} for $\widehat D_E^*$ bounds $f_i$ in the
graph norm by $C\|f_i\|_{L^2(\Omega_0)}\le C$. After passing to a subsequence, $f_i$
converges weakly on $X_0$, and strongly in $L^2(\Omega_0)$, to some $f$, while its
compression component satisfies $\int_I{\mathfrak s}_E\|f_{i,c}\|^2dt\le C\eps_i^{1/2}$, by
the argument of \eqref{eq:compression-small}. Let $\zeta$ be a smooth compactly
supported decaying field in the domain of $\widehat D_V$, and lift it to
$E_{\eps_i}$ by $\zeta_i=(\zeta,(0,h_{\nu_\zeta(t)}))$, where $\nu_\zeta(t)\in
T_{q_V(t)}P$ is the section component of $\zeta$ on $I$. For large $i$ the support of
$\zeta$ lies where $q_{R_i}=q_V$, so $\zeta_i$ lies in the domain of $\widehat D_E$;
off the compression $\widehat D_E\zeta_i=\widehat D_V\zeta$, and its compression
component $e^{\phi_E}\mathsf D_t(0,h_{\nu_\zeta})$ has norm $O(\eps_i^{1/4})$. Hence
$0=\langle\widehat D_E\zeta_i,f_i\rangle\to\langle\widehat D_V\zeta,f\rangle$.
Likewise $0=\langle\widehat Fv,f_i\rangle\to\langle e^{\phi_E}D_V\mathfrak
k_0(v),f\rangle$, because $\widehat Fv$ converges strongly by
\eqref{eq:column-error}. Smooth compactly supported fields are dense in the decaying
domain, so $f$ annihilates the range of the conjugated operator $D_V$ on the extended weighted spaces, which is the whole space; hence $f=0$. Then
$\|f_i\|_{L^2(\Omega_0)}\to0$, and the graph bound gives $\|f_i\|\to0$, a
contradiction. Thus $D_E$ is onto for $\eps$ small, and the bound of Step~2 is the
bound of the proposition.
\end{proof}

\begin{lemma}\label{lem:exterior-small-weight}
The auxiliary exterior operator $D_E$ on the extended weighted spaces at the new ends is an
isomorphism at weight $\delta_L$ on the new ends, and
\[
\|z_E\|_{X_E}\le C\mathsf A\|D_Ez_E\|,\qquad |v_-|+|v_+|\le C\sqrt{\mathsf A}\|D_Ez_E\|,
\]
where $v_\pm$ are the limiting values of $z_E$ at the new ends. A solution of a source
supported in $\{s\le s'\}$ on a new end, with $s'\ge L/4+c_R$ so that $W^E_\eps$ is flat
beyond $s'$, differs from its stationary limit, a distance $\mathsf T$ beyond $s'$, by
$O(e^{-\sigma_0\mathsf T})$ times the bound, with $\sigma_0$ the rate
of~\eqref{eq:normal-estimate}.
\end{lemma}

\begin{proof}
At the fixed weight $\delta_*$ the operator is an isomorphism with a bound independent
of $\eps$ by Proposition~\ref{prop:exterior-inverse}. Passing from $\delta_*$ to
$\delta_L$ on the two new ends is the argument of the proof of
Proposition~\ref{prop:transition-inverse}(iii) on each end, with the inverse of
Proposition~\ref{prop:exterior-inverse} in place of the fixed-weight transition
inverse: the normal part is solved with no loss, and the stationary coefficient costs
$\delta_L^{-1}=C\mathsf A$ in the weighted $L^2$ norm, which gives the first bound, and
$\delta_L^{-1/2}=C\sqrt{\mathsf A}$ in the weighted supremum, which bounds the limiting values.
On the part $s\le L/4+c_R$ of an end the configuration $W^E_\eps$ is not flat; in the
identification of Step~1(b) of the proof of Proposition~\ref{prop:exterior-inverse}
its deviation adds to both equations of the triangular system~\eqref{eq:triangular}
terms bounded by $Ce^{-\beta s}(|a|+\|z^{\rm n}\|)$, and the scale adds the term
$\mathsf b z^{\rm n}$ with coefficient $O({\mathfrak s}_E^{1/2}|\partial_s\log{\mathfrak s}_E|)$.
The tail is solved with fixed limit $0$, as $S_\delta$ is in
Proposition~\ref{prop:transition-inverse}(iii), and the limiting vector comes from
the $\delta_*$-inverse applied to the residual supported near $s_1$.
Measuring $a$ by $\sup_se^{\delta_Ls}|a(s)|$ and $z^{\rm n}$ in the weighted graph
norm, these terms enter the normal equation and the stationary integral with
coefficients whose $L^1$ and $L^2$ norms over the end are $O(e^{-\beta s_1}+\eps^{1/4})$,
and they are absorbed with no further loss in $\delta_L$ for $s_1$ fixed large and $\eps$
small; their contribution to the weighted $L^2$ norm of $a$ is at most
$C(1+\eps^{1/4}L^{1/2})$ times the weighted supremum, hence $O(\delta_L^{-1/2})$. The
homogeneous-tail statement is the exponential decay of the normal part and the
integration of the stationary equation from the source-free region.
\end{proof}

\begin{proposition}\label{prop:attachment-inverse}
For $\eps$ small, the linearization $D_\eps$ of the mixed equation at $W_\eps$
is invertible, with $\|D_\eps^{-1}f\|_{X_\eps}\le C\mathsf A\,\|f\|_{Y_\eps}$, and
$\ind D_\eps=0$.
\end{proposition}

\begin{proof}
\emph{Source removal.} On each overlap put $u=(s-\alpha_L)/\mathsf T$ and fix cutoffs
$\theta,\omega$ of $u$ with $\theta\omega=\omega$, $\theta$ supported in $(1/8,7/8)$ and
equal to one on $[1/4,3/4]$, $\omega$ supported in $(1/4,3/4)$ and equal to one on
$[1/3,2/3]$. Let $S$ be the inverse of the stationary line problem on the
overlap, extended to a complete line of constant scale beyond it, at weight
$e^{\delta_L(s-\alpha_L)}$. The line estimate for $S_\delta$ in the proof of
Proposition~\ref{prop:transition-inverse}(iii) holds for every scale profile with
$|(\log{\mathfrak s})'|\le1$ and small $\sup{\mathfrak s}$, in particular for this line of
constant scale on both sides, and gives $\|S(\omega f)\|\le C\mathsf A\|f\|$. Put $z_{\mathrm{loc}}=\sum_\pm\theta S(\omega f)$ and
$f'=f-D_\eps z_{\mathrm{loc}}=(1-\omega)f-[D_\eps,\theta]S(\omega f)$. Since the
derivative of $\theta$ is $O(\mathsf T^{-1})=O(\mathsf A^{-1})$, $\|f'\|\le C\|f\|$, and $f'=0$ on
$1/3\le u\le2/3$.

\emph{Components.} Split $f'$ into sources $f_E$ on the exterior and $f_\pm$ on
the transitions, with the partition changing inside the band $1/3\le u\le2/3$ on which $f'=0$. Solve
$z_E=D_E^{-1}f_E$ by Lemma~\ref{lem:exterior-small-weight}, and on each
transition $z_\pm=k_{v_\pm}+\mathcal D^{-1}f_\pm$ with the limiting values $v_\pm$ of $z_E$.
This imposes two separate equalities between the exterior and each transition;
no equality between $v_-$ and $v_+$, and no rank condition on the exterior
evaluation, is imposed. Both component solutions are homogeneous a distance at
least $\mathsf T/6-O(1)$ before the joining point $m_L=(\alpha_L+\beta_L)/2$. On the
exterior Lemma~\ref{lem:exterior-small-weight} applies; on the transitions the
triangular system~\eqref{eq:triangular} shows that the normal part decays at rate
$\sigma_0$ while the stationary part changes only by $\int\mathsf bz^{\rm n}$. Hence,
on a fixed collar of $m_L$,
$\|z_E-z_\pm\|\le C\mathsf Ae^{-\sigma_0\mathsf T/6}\|f\|$; the stationary parts cancel because
the limiting values agree.

\emph{Patching.} With cutoffs $\beta_E$ and $\beta_\pm=1-\beta_E$ of fixed slope at
$m_L$, put $\mathcal P(f)=z_{\mathrm{loc}}+\beta_Ez_E+\sum\beta_\pm z_\pm$. The boundary
conditions are linear, so $\mathcal P(f)$ lies in the domain, and
$D_\eps\mathcal P(f)=f+\sum[D_\eps,\beta_E](z_E-z_\pm)$, whose error has norm at most
$C\mathsf Ae^{-\sigma_0\mathsf T/6}\|f\|$. Converting the component norms to $X_\eps$ with
\eqref{eq:stationary-cost}, the stationary pieces cost $C\sqrt{\mathsf A}(|v_-|+|v_+|)\le
C\mathsf A\|f\|$, so $\|\mathcal P(f)\|\le C\mathsf A\|f\|$. For large $L$ the Neumann series gives a
right inverse $\mathcal Q_\eps=\mathcal P(1+\mathcal E_{\mathcal P})^{-1}$, where
$\mathcal E_{\mathcal P}=D_\eps\mathcal P-1$ has norm at most $C\mathsf Ae^{-\sigma_0\mathsf T/6}$, with
$\|\mathcal Q_\eps\|\le C\mathsf A$.

\emph{Injectivity.} Let $D_\eps z=0$. On each overlap write $z=k_{a}+w$ with $w$
normal; $w$ decays exponentially toward $m_L$ in terms of its collar norms, and
$|a(m_L)|\le C\|z\|$. Cut $z$ into component fields $z_E^{\rm cut}$, $z_\pm^{\rm cut}$
equal to $z$ on $\operatorname{supp}\beta_E$, respectively
$\operatorname{supp}\beta_\pm$, and completed beyond by the stationary value
$k_{a(m_L)}$, so that $z=\beta_Ez_E^{\rm cut}+\sum\beta_\pm z_\pm^{\rm cut}$; their sources are cutoff derivatives times
$z-k_{a(m_L)}$, of norm $Ce^{-\sigma_0\mathsf T/6}\|z\|$. The component inverses and
\eqref{eq:stationary-cost} give $\|z\|\le C\mathsf Ae^{-\sigma_0\mathsf T/6}\|z\|$, so $z=0$ for large $L$,
since $\mathsf T\sim L/4$ and $\mathsf A=1+L$.

\emph{Index.} $D_\eps$ is invertible, so its index is $0$, in agreement with
Theorem~\ref{thm:collapse}(i), since $W_\eps$ is a configuration in the class of
$V$.
\end{proof}

The loss of one factor $\mathsf A$ is sharp in these norms: a stationary field $k_v$
multiplied by a cutoff of slope $\delta_L$ has norm ratio $c/\delta_L$. It is harmless
because the approximation error of $W_\eps$ is a positive power of $\eps$
(Section~\ref{sec:collapse-coverage}).

\subsection{The scalar Laplacian on the glued domain}\label{subsec:scalar}

Let $g_\eps$ be the complete
conformal metric on the gauge domain of $X_\eps$ used in
Section~\ref{subsec:patching}. Choose a smooth $\lambda_\eps>0$ equal to ${\mathfrak s}$
on the compression, equal to one on the unchanged pieces and away from the
ends of the transition cores, interpolated on fixed buffers and constant in
the normal direction at the boundary, and put $\widehat g_\eps=\lambda_\eps^{-2}g_\eps$,
$d\mu_\eps=\lambda_\eps^2dV_{\widehat g_\eps}$. On the compression
$\widehat g_\eps=dT^2+g_C$ with $dT=dt/{\mathfrak s}$; the rescaled metric has uniformly
bounded local geometry, and its large parts are product cylinders with
compact cross-section joined through finitely many fixed cores. In
dimension four the physical norm of a one-form is
$\|a\|^2_{L^2(g_\eps)}=\int|a|^2_{\widehat g_\eps}d\mu_\eps$, which is
$\int{\mathfrak s}|(p,b)|^2dt\,dV_C$ on the compression, with $p={\mathfrak s}\psi$.

Write $H^j_{\rm uloc}(\widehat g_\eps)$ for the supremum of the
$H^j$ norms over a uniform cover by unit patches in the rescaled
metric, with covariant derivatives taken using the reference
connection.

\begin{lemma}[scalar inverse]\label{lem:scalar-inverse}
Let $\mathcal L_\eps=\lambda_\eps^2d^{*,g_\eps}_{\mathbb A}d_{\mathbb A}$ with Neumann condition, for a
reference $\mathbb A$ within a fixed small $C^k$-distance of $W_\eps$ in the rescaled
metric. Then
$\|\chi\|_{L^2(\mu_\eps)}\le C_P\|d_{\mathbb A}\chi\|_{L^2(\mu_\eps)}$ and
$\|\mathcal L_\eps^{-1}f\|_{H^{m+1}_{\rm uloc}(\widehat g_\eps)}\le C_m\|f\|_{H^{m-1}_{\rm
uloc}(\widehat g_\eps)}$, with constants independent of $\eps$; the same holds
with the weights of the original ends.
\end{lemma}

\begin{proof}
The quadratic form of $\mathcal L_\eps$ is $\int|d_{\mathbb A}\chi|^2d\mu_\eps$. On each
product piece of the compression the Poincaré inequality is the spatial
Neumann gap, uniform over $P$ because every flat connection on $C$ is
irreducible, integrated in time. The transition cores have compact
Neumann inequalities once enlarged to contain an irreducible slice (on each core
$\lambda_\eps$ is comparable with one), the tails of the original ends have their spatial
gaps, and a fixed compact core of the unchanged gauge domain, enlarged to contain an
irreducible slice, has its Neumann inequality. Each of these regions $R_k$ therefore carries an inequality
$\|\chi\|_{L^2(R_k)}\le C\|d_{\mathbb A}\chi\|_{L^2(R_k)}$ with no mean-value term, since
$\mathbb A$ is irreducible on every slice, and the regions cover with overlap at most $N_0$, so
\[
\|\chi\|^2_{L^2(\mu_\eps)}\le\sum_k\|\chi\|^2_{L^2(R_k)}
\le C\sum_k\|d_{\mathbb A}\chi\|^2_{L^2(R_k)}\le CN_0\|d_{\mathbb A}\chi\|^2_{L^2(\mu_\eps)} .
\]
For the inverse,
cover by unit rescaled patches $Q_j$; the number of patches in a distance shell
grows polynomially. Conjugating the form by $e^{\gamma d_j}$, with $d_j$ a
smoothed distance from $Q_j$ and $\gamma$ small, and using the variational
inverse and local elliptic regularity gives
\[
\|\mathcal L_\eps^{-1}f_j\|_{H^{m+1}(Q_\ell)}\le C_me^{-\gamma d(Q_\ell,Q_j)}
\frac{\lambda_\eps(Q_j)}{\lambda_\eps(Q_\ell)}\|f_j\|_{H^{m-1}(Q_j)}
\]
for $f_j$ supported in $Q_j$. Since $\partial_T({\mathfrak s}^{-1})=-h$, the function
$\lambda_\eps^{-1}$ is uniformly Lipschitz in the rescaled distance and bounded
below, so the ratio is at most $C(1+d(Q_\ell,Q_j))$, which is absorbed by half
the exponential. Writing $f=\sum_jf_j$ with $f_j$ supported in $Q_j$ and
$\|f_j\|_{H^{m-1}}\le C\|f\|_{H^{m-1}(Q_j')}$ for a fixed enlargement $Q'_j$,
\begin{align*}
\|\mathcal L_\eps^{-1}f\|_{H^{m+1}(Q_\ell)}
&\le C_m\sum_je^{-\gamma d(Q_\ell,Q_j)/2}\|f_j\|_{H^{m-1}(Q_j)}\\
&\le C_m\Bigl(\sup_{n\ge0}\,\#\{j:n\le d(Q_\ell,Q_j)<n+1\}\,e^{-\gamma n/4}\Bigr)
\sum_{n\ge0}e^{-\gamma n/4}\,\|f\|_{H^{m-1}_{\rm uloc}},
\end{align*}
and the supremum is finite, uniformly in $\ell$ and $\eps$, because the shell counts grow
polynomially in $n$. This is the inverse bound; uniqueness in the
uniformly local space follows from the quadratic identity with a small
decreasing exponential weight. The weights of the original ends, whose logarithmic
gradients are smaller than $\gamma$, are handled in the same way.
\end{proof}

\section{The collapse theorem: gluing and surjectivity}\label{sec:collapse-coverage}

Throughout, $V$ is a regular index-zero solution on $X_0$ with limits
$q_\pm$ at the new ends, $L=\log(1/\eps)$ and $\mathsf A=1+L$. Constants $C,c$
depend on the fixed data and on the finite set $\cM^0(X_0)$ but not on
$\eps$, on neck lengths, or on the minimum of the compression scale.

The section has two parts. Existence (Section~\ref{subsec:charts} and Proposition~\ref{prop:existence}):
the approximate solution $W_\eps$ is corrected by a contraction in strong norms, whose inverse
loses the factor $\mathsf A^{m+1}$, so the solution is unique only in a ball of radius
shrinking like a power of $1/\mathsf A$. Surjectivity (Proposition~\ref{prop:coverage}): an
arbitrary solution near $V$ is first shown, by cutting it along the necks and comparing the
pieces with the components, to lie within a positive power of $\eps$ of $W_\eps$, by estimates
that lose no logarithm; only then does it enter the uniqueness ball.

\subsection{Strong norms and the strong inverse}\label{subsec:charts}

\emph{Strong norms.} On the small compression, with $\mathsf w$ the fixed end weight
(equal to one on the region $N_I\cap X_\eps$, which contains the part where $X_\eps$ differs from $X_0$, Definition~\ref{def:insertion}(iv)), put
\[
\|f\|_{Y^m_\eps}^2=\sum_{i+r\le m}\|\nabla^i_t(\mathsf w{\mathfrak s}^{1/2}f)\|^2_{L^2_tH^r(C)},
\]
\[
\|z\|_{X^m_\eps}^2=\|\mathsf wv\|^2_{H^m_t}+\sum_{i=0}^m\|\mathsf w{\mathfrak s}^{-1/2}\nabla^i_tU_\perp\|^2_{L^2_t
H^{m+1-i}(C)}+\|\mathsf w{\mathfrak s}^{1/2}\nabla_t^{m+1}U\|^2_{L^2_tL^2(C)},
\]
and on the other gauge regions and the section weighted $H^m$ and $H^{m+1}$,
patched by cutoffs of bounded overlap, with the free limiting vectors at the
original ends added. For $m=0$ this is equivalent to $X_\eps$. On the region
$N_I\cap X_\eps$
\begin{equation}\label{eq:strong-control}
\|U\|_{H^m}+\|{\mathfrak s}^{1/2}U\|_{H^{m+1}}+\|{\mathfrak s}^{-1/2}U_\perp\|_{H^m}+\|{\mathfrak s}^{-1/2}\nabla_CU_\perp\|_{H^m}
\le C_m\|z\|_{X^m_\eps},
\end{equation}
and Sobolev embedding on unit patches gives unscaled $C^1$ bounds for the
spatial coefficients, $p$ and the section variation.

\begin{proposition}[strong inverse]\label{prop:strong-inverse}
$D_\eps\colon X^m_\eps\to Y^m_\eps$ is invertible with
$\|D_\eps^{-1}\|\le C_m\mathsf A^{m+1}$.
\end{proposition}

\begin{proof}
Let $\mathcal V$ be a global tangential differential operator equal to $\nabla_t$ where
${\mathfrak s}\le{\mathfrak s}_0/2$ on the compression, extended over the adjacent section collar
by the corresponding boundary derivative, cut off where ${\mathfrak s}\ge{\mathfrak s}_0$ and in
the collar interior, and vanishing on the other regions and on the supports of
all perturbations. By Lemma~\ref{lem:transport}, $\mathcal V$ preserves the linear
domain. We prove by induction on $j\le m$ that
\begin{equation}\label{eq:time-derivatives}
\|\mathcal V^jz\|_{X_\eps}\le C_j\mathsf A^{j+1}\|D_\eps z\|_{Y^j_\eps}.
\end{equation}
For $j=0$ this is Proposition~\ref{prop:attachment-inverse}. At fixed $\eps$ the region
where the geometry changes is compact, so $\mathcal V^jD_\eps^{-1}f$ lies in $X_\eps$ by
fixed-geometry regularity \cite[Theorem~4.17 and Remark~4.19]{DFL}, or by difference
quotients; the estimates below then bound it. Commuting $D_\eps$ through
$\mathcal V^j$ gives
\[
D_\eps(\mathcal V^jz)=\mathcal V^jD_\eps z+
\sum_{i<j}\sum_{k=0}^{j-i}{\mathfrak s}^{-1}a_{ijk}(t)
\mathcal B^{[k]}_t\mathcal V^iU+\sum_{i<j}b_{ij}\mathcal V^iz,
\]
where $\mathcal B^{[0]}_t=\mathcal B_t$, all coefficients are bounded,
and $b_{ij}$ collects the remaining terms: those of order zero, and the
first-order terms in $\nu$ produced by the section-collar cutoff in $\mathcal V$,
which are bounded by $\|\mathcal V^i\nu\|_{H^1}\le\|\mathcal V^iz\|_X$. Differentiating the scale
produces the terms with $k=0$, even when the background is constant. By
Lemma~\ref{lem:transport}, for every $k\ge0$ the operator
$\mathcal B^{[k]}_t$ annihilates the harmonic fields and
$\|\mathcal B^{[k]}_tU\|\le C\|\mathcal B_tU\|$, so no factor ${\mathfrak s}^{-1}$ multiplies a slow
term and
\[
\|{\mathfrak s}^{-1}a_{ijk}\mathcal B^{[k]}_t\mathcal V^iU\|_{Y_\eps}\le C\|{\mathfrak s}^{-1/2}\mathcal B_t\mathcal V^iU\|\le
C\|\mathcal V^iz\|_{X_\eps}.
\]
Applying Proposition~\ref{prop:attachment-inverse} to $\mathcal V^jz$, which costs one factor
$\mathsf A$, and the induction hypothesis to the terms with $i<j$,
\[
\|\mathcal V^jz\|_{X_\eps}\le C\mathsf A\Bigl(\|D_\eps z\|_{Y^j_\eps}+\sum_{i<j}C_i\mathsf A^{i+1}
\|D_\eps z\|_{Y^i_\eps}\Bigr)\le C_j\mathsf A^{j+1}\|D_\eps z\|_{Y^j_\eps},
\]
which is~\eqref{eq:time-derivatives}. The spatial derivatives cost nothing: on the
compression they are recovered from the equation by~\eqref{eq:complementary}; for
the section, the derivative in the direction normal to the seam is recovered
from the Cauchy--Riemann equation and the tangential derivative, and the seam
trace of $\nabla_t^j\nu$ is controlled through the harmonic matching, which
$\nabla_t$ preserves by Lemma~\ref{lem:transport}; both bounds are uniform in
${\mathfrak s}$. In the regions where the geometry is uniformly nondegenerate, ordinary
mixed elliptic estimates give $H^{m+1}$ control by $\|D_\eps z\|_{Y^m}+\|z\|_{X_\eps}$.
Summing gives the bound with $j=m$.
\end{proof}

\emph{The exact matching chart.} Along $q(t)$, let $B_t(v)$ be the flat
connection in the spatial Neumann Coulomb slice of $B_t=B_{q(t)}$ with bulk
harmonic coordinate $v$: $F_{B_t(v)}=0$, $d^*_{B_t}(B_t(v)-B_t)=0$, vanishing normal
component, $\ell_t(0,B_t(v)-B_t)=v$, $D_vB_t(0)=j_t$. Let $e_t$ be a section chart
whose restriction to $\mathscr L_t$ is the boundary class of $B_t(v)$. For
$U=j_tv+U_\perp$ with $U_\perp=(p,w)$ put
\[
\Phi_t(U)=B_t+\bigl(B_t(v)-B_t+T_{t,v}w+\mathcal R_t(v,w)\bigr)+(\Psi_t+p/{\mathfrak s})\,dt,
\]
where $T_{t,v}$ is the transport of the complementary field to $B_t(v)$ by the surface Hodge
decompositions on the parallel surfaces of a fixed collar of $\partial C$, blended to the
identity in the interior and leaving the normal component unchanged, with the harmonic summands
identified so that the class of the trace $(T_{t,v}w)|_\Sigma$ is $de_t(v)\,\tau_tw$; and
$\mathcal R_t(v,w)$ is the surface correction defined below, which is tangential, is supported in
the collar, makes the trace flat with class $e_t(v+\tau_tw)$, and satisfies $\mathcal R_t(v,0)=0$
and $D_w\mathcal R_t(v,0)=0$. On the section use $e_t$ at the seam extended into a
collar. In this chart the matching and flatness conditions hold exactly, and
the mixed operator $\mathfrak B_t(v,U_\perp)$ satisfies
\begin{equation}\label{eq:exact-chart}
\mathfrak B_t(v,0)=0,\qquad D_{U_\perp}\mathfrak B_t(0,0)=\mathcal B_t ,
\end{equation}
so the quadratic part vanishes on the slow directions.

\emph{The surface correction.} Let $[0,r_1)\times\partial C$ be a product collar of $\partial C$
with normal coordinate $r$, write $\Sigma_r=\{r\}\times\partial C$, and let $\chi$ be a cutoff of
$r$ equal to $1$ near $0$ and to $0$ for $r\ge r_1/2$. For given $t$, $v$ and $r$ let
$\alpha=\alpha_{t,v,r}$ be the restriction of $B_t(v)$ to $\Sigma_r$. It is flat and, by (C1),
irreducible on every component, so $\Delta_\alpha=d_\alpha d_\alpha^*$ is invertible on
two-forms and $\ker d_\alpha^*$ on one-forms is the sum of the harmonic space $\mathcal H_\alpha$
and $\operatorname{im}d_\alpha^*$. Let $\mathfrak i_{t,v,r}\colon T_{q(t)}\Msf\to\mathcal H_\alpha$ be the
isomorphism given by $de_t(v)$ and the identification of the cohomology of $\Sigma_r$ with that
of $\Sigma$ by the flat connection $B_t(v)$ on the collar (on the torus faces $\mathcal
H_\alpha=0$). For a tangential one-form $\phi$ on $\Sigma_r$ small in $H^{1/2}$ and
$\vartheta\in T_{q(t)}\Msf$ small, let $\zeta=\mathcal S_\alpha(\phi,\vartheta)$ be the small
solution in $\ker d^*_\alpha$ of the \emph{surface flatness equation}
\begin{equation}\label{eq:surface-flatness}
d_\alpha\zeta+\tfrac12[(\phi+\zeta)\wedge(\phi+\zeta)]=0,\qquad
\text{harmonic part of }\zeta=\mathfrak i_{t,v,r}\vartheta .
\end{equation}
Equivalently $\zeta=\mathfrak i_{t,v,r}\vartheta+\mathsf K_\alpha(\phi+\zeta,\phi+\zeta)$ with
$\mathsf K_\alpha(a,b)=-\frac12d_\alpha^*\Delta_\alpha^{-1}[a\wedge b]$, a bilinear operator of order
$-1$ on $\Sigma_r$; since $[\cdot\wedge\cdot]$ maps $H^{1/2}\times H^{1/2}$ to $L^2$ on a
surface, the contraction principle gives $\mathcal S_\alpha$, smooth in $(\phi,\vartheta)$ and in
the parameters, with values in $H^1$ and $\mathcal S_\alpha(0,0)=0$. By construction
$F_{\alpha+\phi+\zeta}=d_\alpha\phi$; in particular $\alpha+\phi+\mathcal S_\alpha(\phi,\vartheta)$ is
flat whenever $d_\alpha\phi=0$. Let $\phi_{t,r}$ be the tangential part of $T_{t,v}w$ on
$\Sigma_r$. On $\Sigma_0=\partial C$ it is closed, because $w$ has closed traces and $T_{t,v}$
preserves closedness. Let $\vartheta=\vartheta_t(v,w)$ be the small solution of the
finite-dimensional equation
\begin{equation}\label{eq:class-equation}
\bigl[\alpha_{t,v,0}+\phi_{t,0}+\mathcal S_{\alpha_{t,v,0}}(\phi_{t,0},\vartheta)\bigr]=e_t(v+\tau_tw),
\end{equation}
which exists by the implicit function theorem, since $\vartheta\mapsto[\alpha_{t,v,0}+\mathcal
S_{\alpha_{t,v,0}}(0,\vartheta)]$ is the harmonic chart of $\Msf$ at $e_t(v)$, with derivative
$de_t(v)$ at $0$; the class map is smooth on flat connections of class $H^{1/2}$ near
$\alpha_{t,v,0}$. At $w=0$ the derivative of the left side in the class equation is
$de_t(v)(\tau_t\dot w+\dot\vartheta)$, whereas that of the right side
is $de_t(v)\tau_t\dot w$. Thus $\vartheta_t(v,0)=0$ and
$D_w\vartheta_t(v,0)=0$ for every $v$.

\emph{Extension in the normal direction.} The harmonic parameter is a
function of the boundary trace of $w$. Its extension must respect
both the time and the normal derivatives in $X^m_\eps$.
Choose a nonnegative $\rho\in C^\infty_c((-1,1))$ with integral one.
At a fixed scale $s\in(0,1]$, in a smooth frame of $T_{q(t)}\Msf$, set
\[
(\mathcal E_s\vartheta)(t,r)
 =\int\rho(u)\vartheta(t-sru)\,du,\qquad r\ge0.
\]
Its boundary value is $\vartheta(t)$. For a varying scale, choose
unit time intervals, a subordinate partition $\sum\psi_j=1$, and
larger intervals on which ${\mathfrak s}/{\mathfrak s}(j)$ and its
logarithmic derivatives are bounded. Extend $\psi_j\vartheta$ by
$\mathcal E_{{\mathfrak s}(j)}$ and sum, with a time cutoff equal to
one on the slightly enlarged support of each summand. The collar
radius is fixed small enough that all translates used in this formula
remain in the larger interval. Denote the result by
$\mathcal E_{\mathfrak s}\vartheta$. The overlap and cutoff constants
are independent of the minimum scale and of the number of intervals.
Define the exact correction by
\begin{equation}\label{eq:trace-extension-chart}
\Theta(t,r)=(\mathcal E_{\mathfrak s}\vartheta(v,w))(t,r),\qquad
\mathcal R(v,w)(t,r)=\chi(r)\mathcal S_{\alpha_{t,v,r}}
                  (\phi_{t,r},\Theta(t,r)).
\end{equation}
At $r=0$ this is the solution of \eqref{eq:class-equation}. Thus the
trace is flat with the prescribed class. The correction and its
first derivative in $w$ vanish at $w=0$, so \eqref{eq:exact-chart}
and the linearized operator are unchanged.

Here is the estimate for this extension. On a constant-scale patch,
write $\mathcal X_s^m$ for the complementary part of $X^m_\eps$,
with end weight one. The trace norm $\mathcal T_s^m$ is given, in
local time Fourier coordinates and surface spectral coordinates, by
\[
\|f\|_{\mathcal T_s^m}^2=
\int\sum_\lambda s^{-1}(1+\xi^2+\lambda^2)^m
 (1+\lambda^2+s^2\xi^2)^{1/2}|\widehat f(\xi,\lambda)|^2\,d\xi.
\]
Equivalent norms are obtained with a finite surface atlas and smooth
cutoffs. For a finite-dimensional harmonic coefficient omit $\lambda$.
The trace and extension estimates are
\begin{equation}\label{eq:anisotropic-trace}
\|w|_{\partial C}\|_{\mathcal T_s^m}\le C_m\|w\|_{\mathcal X_s^m},
\qquad
\|\chi\mathcal E_s\vartheta\|_{\mathcal X_s^m}
 \le C_m\|\vartheta\|_{\mathcal T_s^m}.
\end{equation}
For the first, apply the normal one-dimensional trace inequality
with parameter $(1+\lambda^2+s^2\xi^2)^{1/2}$ to the time and
surface transform. Tangential differentiation gives the factor
$(1+\xi^2+\lambda^2)^m$. The resulting bulk norm is equivalent to
the displayed complementary norm in $X^m_\eps$: terms containing
$s^2\xi^2$ and a spatial derivative are dominated by its spatial
sum, and the remaining highest time term is
$s\|\partial_t^{m+1}w\|^2$. For the second estimate the multiplier
is $\widehat\rho(sr\xi)$. Its derivatives decrease rapidly for
$sr|\xi|\gg1$. Integrating each normal derivative in $r$ gives
\[
\int_0^{r_1}|\partial_r^j(\chi(r)\widehat\rho(sr\xi))|^2\,dr
 \le C_j(1+s|\xi|)^{2j-1},
\]
with the same assertion for $j=0$. Since $s\le1$, inserting the
time derivatives and summing proves \eqref{eq:anisotropic-trace}.
In particular the top time derivative is estimated through a
half-order boundary norm, not an $L^2$ trace of that derivative.

\begin{lemma}[quadratic estimates in the matching chart]\label{lem:block}
Fix $m\ge4$ and a compact smooth family of reference paths $q$.
In the chart \eqref{eq:trace-extension-chart}, write
$\mathcal F_\eps(W_\eps+z)=f_\eps+D_\eps z+N_\eps(z)$.
For $z,z'$ in a sufficiently small ball of $X^m_\eps$, the \emph{strong ball}, with
constants independent of $\eps$,
\begin{equation}\label{eq:strong-quadratic}
\|N_\eps(z)-N_\eps(z')\|_{Y^m_\eps}
 \le C_m(\|z\|_{X^m_\eps}+\|z'\|_{X^m_\eps})
             \|z-z'\|_{X^m_\eps}.
\end{equation}
Let $\mathsf d(z)$ be the supremum of $\|v\|_{C^0}+
\|U_\perp\|_{C^0}$ on compression patches, of $\|a\|_{C^0}$ on
other gauge patches, and of $\|\nu\|_{C^0}$ on section patches.
On a small smooth coefficient neighbourhood one also has
\begin{equation}\label{eq:low-quadratic}
\|N_\eps(z)\|_{Y_\eps}\le C\mathsf d(z)\|z\|_{X_\eps}.
\end{equation}
If $z_*$ is in the strong ball and $z$ is smooth with small
coefficient size, then
\begin{equation}\label{eq:asymmetric-quadratic}
\|N_\eps(z)-N_\eps(z_*)\|_{Y_\eps}
 \le C(\mathsf d(z)+\|z_*\|_{X^m_\eps})\|z-z_*\|_{X_\eps}.
\end{equation}
The same local assertions hold on the component and exterior domains.
\end{lemma}

\begin{proof}
We first describe the surface estimate that is needed in addition to
\eqref{eq:anisotropic-trace}. Freeze $t,v$ and use the harmonic frame
$\mathfrak i_{t,v,0}$. Thus the prescribed linear class is
$h=\mathfrak i_{t,v,0}^{-1}H_\alpha\phi=\tau_tw$, and $A(v,h)$
is the flat slice representative of $e_t(v+h)$ relative to
$\alpha=B_t(v)|_\Sigma$. We write $A(h)$ while $v$ is fixed.
For differentiation, use the following gauge equation also for
nonclosed $\phi$, with $h=\mathfrak i_{t,v,0}^{-1}H_\alpha\phi$.
It defines an extension of the harmonic class parameter off the
closed-form subspace; on that subspace it equals the parameter in
\eqref{eq:class-equation}. The class equation is
implemented by solving
\[
d_\alpha^*\bigl(g^*A(h)-\alpha-\phi\bigr)=0,
\qquad g=\exp\eta.
\]
For a fixed $p>2$ this is an equation from $W^{1,p}$ to
$W^{-1,p}$, with invertible derivative $d_\alpha^*d_\alpha$ at
zero. Small $C^0$ norm of $\phi$ gives a small $W^{1,p}$ solution
$\eta$; hence $\eta$ is small in $C^0$. The differentiated equation,
viewed from $H^1$ to $H^{-1}$, gives
\[
\|D\eta[\dot\phi]\|_{H^1}\le C\|\dot\phi\|_{L^2}.
\]
Indeed its leading coefficient is close in $C^0$ to that of the
covariant Laplacian, and its remaining coefficients are small in
$L^p$. H\"older's inequality and the surface Sobolev inequality
absorb those terms using the irreducible scalar gap.

Now $\vartheta=H_\alpha(g^*A(h)-\alpha-\phi)$ in the chosen
harmonic frame. Harmonic projection kills the exact part of
$D\eta$, and the linear terms $\dot h-H_\alpha\dot\phi$ cancel.
The remaining terms contain $\phi$, $g-1$, or the nonlinear part
of the flat slice. Pairing with smooth harmonic forms and using
the preceding estimate gives
\begin{equation}\label{eq:class-map-L2}
|D_\phi\vartheta[\dot\phi]|
 \le C\|\phi\|_{C^0}\|\dot\phi\|_{L^2}.
\end{equation}
For example, the commutator term is bounded by
$C\|g^*A(h)-\alpha\|_{L^2}\|D\eta[\dot\phi]\|_{L^2}$.
The quantitative differentiated estimates, in the fixed frame, are
\begin{align}
|D_\phi^2\vartheta[\psi_1,\psi_2]|&\le
 C\|\psi_1\|_2\|\psi_2\|_2,\nonumber\\
|D_vD_\phi\vartheta[b,\psi]|&\le
 C|b|\|\phi\|_2\|\psi\|_2,\nonumber\\
|D_v\vartheta[b]|&\le C|b|\|\phi\|_2^2.
\label{eq:class-map-mixed}
\end{align}
To see these bounds, differentiate the weak gauge equation, retaining
it in divergence form. Each differentiated gauge term has at most
one derivative of a gauge parameter. A first variation lies in
$H^1$; products of an undifferentiated $H^1$ variation and an
$L^2$ differentiated variation lie in $L^r$ for every $r<2$.
The inverse weak Laplacian gives $W^{1,r}$ for the next gauge
variation. Surface Sobolev embedding and H\"older's inequality
then bound its harmonic pairings in $L^1$. Repeating this argument
at any fixed finite order gives bounded multilinear harmonic
pairings, choosing the intermediate exponents sufficiently close
to two. No derivative of $\phi$ is used in these estimates.
Differentiating the smooth finite-dimensional flat slice and
harmonic projection adds only terms with the same bounds. The
identities $\vartheta(v,0)=D_\phi\vartheta(v,0)=0$ eliminate the
constant and linear terms in every $v$-derivative. Integrating the
bounded $D_vD_\phi^2\vartheta$ gives the last two estimates in
\eqref{eq:class-map-mixed}, and integrating these differential
bounds gives their difference versions.

Here is an explicit strong trace bound. On a constant-scale patch
put $Z=\|z\|_{X_s^m}$. The trace theorem and the linear matching
identity give
\[
\|\phi\|_{H^m_tL^2_\Sigma}\le C\sqrt{s}\,Z,\qquad
\|\phi\|_{H^{m+1/2}_tL^2_\Sigma}
 +\|v\|_{H^{m+1/2}_t}\le CZ.
\]
For the second bound on $v$, use
$v=\Pi\nu|_\Sigma-\Pi\tau_tw$: the section has $H^{m+1}$
regularity and its trace has $H^{m+1/2}$ regularity, while the
harmonic part of the complementary trace is bounded by
\eqref{eq:anisotropic-trace}. The transport $T_{t,v}$ obeys the
same estimates by the product rule. The one-dimensional tame
product inequalities, applied to the differentiated harmonic
pairings above, now give
\begin{align*}
\|\vartheta\|_{H^m_t}&\le C\|\phi\|_{H^m_tL^2_\Sigma}^2
 (1+\|v\|_{H^m_t}),\\
\|\vartheta\|_{H^{m+1/2}_t}&\le
 C\|\phi\|_{H^m_tL^2_\Sigma}
 \bigl(\|\phi\|_{H^{m+1/2}_tL^2_\Sigma}
       +\|\phi\|_{H^m_tL^2_\Sigma}\|v\|_{H^{m+1/2}_t}\bigr).
\end{align*}
Since $\|\vartheta\|_{\mathcal T_s^m}$ is equivalent to
$s^{-1/2}\|\vartheta\|_{H^m_t}+\|\vartheta\|_{\dot H^{m+1/2}_t}$,
these prove $\|\vartheta\|_{\mathcal T_s^m}\le C\sqrt{s}\,Z^2$
in a fixed small ball. The differentiated inequalities give
\[
\|\vartheta(z)-\vartheta(z')\|_{\mathcal T_s^m}
 \le C\sqrt{s}(\|z\|_{X_s^m}+\|z'\|_{X_s^m})
                  \|z-z'\|_{X_s^m}.
\]
These are the weighted class-parameter estimates used below.

The spatial part of $\mathcal S_\alpha$ is controlled by its
formula $\mathfrak i\Theta+\mathsf K_\alpha(\phi+\zeta,\phi+\zeta)$.
The operator $\mathsf K_\alpha$ gains one surface derivative.
Time and normal derivatives commute through it with bounded
lower-order commutators; the coefficients depend smoothly on the
fixed family and on $v$. For $m\ge4$, the bulk $H^m$ algebra
inequality and the boundary $H^{m-1/2}$ algebra inequality control
all products below the highest order. At the highest time order,
use \eqref{eq:class-map-L2}, its differentiated equation, and
\eqref{eq:anisotropic-trace}. The half-order time seminorm is
estimated by differences before differentiating the extension.
This gives the tame quadratic bound for $\mathcal R$ in
$\mathcal X_s^m$. The same argument applies to the difference of
two corrections. No extra top-order boundary trace is taken.

At first order this argument uses only small coefficients. For
clarity, the harmonic trace norm with $m=0$ is equivalent to
$s^{-1}\|\vartheta\|_2^2+\|\vartheta\|_{\dot H^{1/2}_t}^2$.
Apply \eqref{eq:class-map-L2} to the differences in the second
term. Terms involving a difference of $v$ are treated by splitting
the difference integral at $|t-u|=s$. On the near part,
$|v(t)-v(u)|^2\le |t-u|\int_u^t|v'|^2$ supplies the weighted
quantity $s\|v'\|_2^2$; on the far part use $2\|v\|_\infty$
and the $s^{-1}\|\phi\|_2^2$ term. The two factors of $\phi$ in
a $v$-derivative ensure that the undifferentiated factor is a
small coefficient. This proves the absolute first-order estimate
for the boundary correction. Expanding about a strong reference
puts its differentiated coefficients on that reference and treats
the unknown remainder quadratically; it gives the asymmetric
version by the same difference-integral argument. More explicitly,
write $y=z-z_*$ and use the integral Taylor formula at $z_*$.
Terms quadratic in its complementary remainder are bounded by
$C\mathsf d(y)\|y\|_{X_s}$. In the linear and mixed terms,
$\|\phi_*\|_{C^1_tC^0_\Sigma}\le C\sqrt{s}\,a$ and
$\|v_*\|_{C^1}\le Ca$, where $a=\|z_*\|_{X_s^m}$.
Multiplication by these coefficients is bounded on the half-order
time norm. A term linear in the slow remainder has coefficient
$O(s a^2)$ by \eqref{eq:class-map-mixed}; this cancels the
$s^{-1/2}$ in the low-frequency trace norm, and the high-frequency
part is controlled by the matching bound
$\|v_y\|_{H^{1/2}}\le C\|y\|_{X_s}$.
The near/far split treats the remaining mixed differences.
Consequently
\[
\|\vartheta(z)\|_{\mathcal T_s^0}\le C\mathsf d(z)\|z\|_{X_s},
\quad
\|\vartheta(z)-\vartheta(z_*)\|_{\mathcal T_s^0}
 \le C(\mathsf d(z)+a)\|z-z_*\|_{X_s}.
\]

The other chart terms have the same estimates directly. In the
surface Hodge transport, the principal symbols do not change with
the reference flat connection. Thus $T_{t,v}-1$ has order at most
$-1$ on the fixed surface and is bounded on $C^0$, with norm
$O(|v|)$; its finite-dimensional harmonic part has the same bound.
The spatial operator vanishes on the exact flat family $B_t(v)$.
Every singular remaining term has a complementary factor, measured
with $s^{-1/2}$ in $X_s$. A differentiated slow factor in the
temporal terms carries $s^{1/2}$: in particular
$s^{1/2}|v'w|$ is bounded by $\|w\|_\infty$ times the weighted
temporal norm of $v$. The terms from $T_{t,v}-1$, the scalar row,
and the local bracket are bounded in this way. The pure section
operator has the usual quadratic estimate in its fixed chart.
These estimates prove all three inequalities at constant scale.
Finally freeze the weights on the unit patches used above and sum.
Bounded logarithmic derivatives, bounded overlap and the support
of $\rho$ give the same constants for the varying-scale extension.
\end{proof}

The chart is locally invertible also in uniformly local rescaled
H\"older norms: use $T=t/s$ on a constant-scale patch. The extension
is then $\int\rho(u)\vartheta(T-ru)du$, bounded in $C^\alpha$;
the surface constructions and harmonic projections are smooth in
these norms, and the correction has zero derivative at zero.
The inverse function theorem gives a uniform small neighbourhood.
Smooth fields with small $C^0$ coefficients and uniformly bounded
higher rescaled derivatives enter it by interpolation. This is the
neighbourhood used for the unknown solutions in the proof of surjectivity.
First differentiate the inverse equation in time and in the surface
directions. The leading operator remains a small perturbation of
the identity on the local H\"older space. Then recover normal
regularity from the explicit extension: each normal derivative of
$\int\rho(u)\vartheta(T-ru)du$ is an integral of a time derivative
of $\vartheta$, already controlled by the first step. This triangular
induction preserves every fixed higher rescaled local bound; a
normal derivative of the extension is not a lower-order commutator. At a fixed positive $\eps$ the changed region
is compact, so these smooth local bounds imply finite first-order
norm there. On the original ends subtract the free stationary
limits first. The extension has uniformly bounded time support and
is bounded on the fixed weighted spaces; the trace estimates and
the small-coefficient inverse give the same weighted bounds for the
remainders. Thus the inverse preserves the end limits and finite
$X_\eps$ norm; this assertion uses the trace estimates as well as
local H\"older invertibility.

\subsection{The approximate solution and existence}

Let $R=\eps^{1/4}$. Near each endpoint, for $x\le R$, replace $V$ by its limiting
flat pair $(B_{q_\pm},q_\pm)$ using a cutoff of fixed slope in the logarithmic
coordinate, performed in the exact matching chart about the flat pair, so that
the matching slices remain flat; since $V$ converges to its limits at rate $\beta>0$
(Proposition~\ref{prop:weak-limit}(iii)), this changes $V$ by $O(R^\beta)$ in every
$C^k$ norm. Over $I\times C$ use the flat lift $B_{q_R(t)}$ in the covariant
Neumann slice of the boundary values $q_R(t)$ of the truncated section, with
temporal component chosen so that $\partial_tB-d_B\Psi=h_{\dot q_R}$. On the necks
and transitions use the flat constant pair. The only errors are the truncation
error and the horizontal velocity of the flat lift, whose weighted size is
${\mathfrak s}^{1/2}|\dot q_R|\le C\eps^{1/2}r^{\beta-1/2}$ on the compression; since
$\nabla_t$ preserves the harmonic extension, the same bound holds for all its
derivatives. Integrating,
\begin{equation}\label{eq:approx-error}
\|\mathcal F_\eps(W_\eps)\|_{Y^m_\eps}\le C_m\,e_\eps,\qquad
e_\eps=R^\beta+\sqrt{\eps K_\beta(R)},\quad
K_\beta(R)=1+\int_R^{r_0}r^{2\beta-2}dr ,
\end{equation}
and $e_\eps\le C\eps^{\sigma_1}$ for some $\sigma_1>0$.

\begin{proposition}[existence]\label{prop:existence}
Let $\varrho\le r_0$ be a radius below which the summed quadratic and difference
estimates hold, and with $\varrho\le1/(4C_q)$, where $C_q$ is the constant of the
quadratic estimate. For $\eps$ small there is a unique solution $\cG_\eps(V)=W_\eps+z^*$ with
$\|z^*\|_{X^m}\le\varrho/(C_m\mathsf A^{m+1})$, and in fact $\|z^*\|_{X^m}\le C_m\mathsf A^{m+1}e_\eps$.
It is regular.
\end{proposition}

\begin{proof}
In the exact matching chart the equation is
$\mathcal F_\eps(W_\eps)+D_\eps z+N_\eps(z)=0$. By Proposition~\ref{prop:strong-inverse},
the quadratic estimate and~\eqref{eq:approx-error}, the map
$z\mapsto-D_\eps^{-1}(\mathcal F_\eps(W_\eps)+N_\eps(z))$ is a contraction of the ball of
radius $\varrho/(C_m\mathsf A^{m+1})$ for small $\eps$: its Lipschitz constant is at most
$2C_q\varrho\le\frac12$, and it maps the ball to itself since $\mathsf A^{2m+2}e_\eps\to0$.
For regularity, the asymmetric quadratic bound gives
$\|DN_\eps(z^*)\|_{X\to Y}\le C\,a$ with $a=\|z^*\|_{X^m}\le C\mathsf A^{m+1}e_\eps$, and
$\|D_\eps^{-1}\|\le C\mathsf A$, so $\mathsf A\,a\le C\mathsf A^{m+2}e_\eps\to0$ and the
linearization at $W_\eps+z^*$ is invertible.
\end{proof}

\subsection{The endpoint necks}\label{subsec:necks}

On each endpoint neck, at a time $s$, the restriction of a configuration
to the slice is a triple $c=(B_f,B_c,\nu)$ of a connection on the fixed
copy of $C$, a connection on the collapsing copy, and a path in the
moduli space of $\Sigma$ with $\nu(0)$ and $\nu(\pi)$ the restrictions of
$B_f$ and $B_c$. The critical family is $\cF=\{(B_q,B_q,q):q\in P\}$.
Gluing the triple gives one connection on the double of $C$ along $\Sigma$
with a collar of length $\pi$ inserted, and the mixed action becomes the
Chern--Simons functional of that three-manifold. For $C=C_1$ it has two torus
boundary components carrying the flat condition; since the flat connection on
the odd bundle over $T^2$ is isolated with $H^1(T^2;\ad)=0$, the functional is
still single valued on a gauge-invariant neighbourhood of $\cF$, and its critical set
there is $\cF$, of Morse--Bott type. Since $\cF$ is connected, the action is constant on
it, and we normalize it to vanish there.

\begin{proposition}[neck decay]\label{prop:neck-decay}
\begin{enumerate}[label=(\roman*)]
\item There are $\delta,C,{\mathfrak s}_2>0$ such that for ${\mathfrak s}\in(0,{\mathfrak s}_2]$
  and all triples within $\delta$ of $\cF$ in the $H^1$ norm of the fixed metric $g_C$ on
  both branches,
  $|a(c)|\le C\|\nabla_{\mathfrak s} a(c)\|_{\mathfrak s}^2$, where the gradient is taken
  in the metric weighting the collapsing component by ${\mathfrak s}$.
\item Let $c(s)$, $s\in[s_0,s_1]$, solve the neck equation with
  ${\mathfrak s}(s)={\mathfrak s}_{\rm b}e^{\pm s}\le{\mathfrak s}_2$, all slices within $\delta$ of
  $\cF$, and energy $e$. Let $p(s)\in P$ be the projection determined by
  the flat projection of the fixed branch $B_f$. Then the energy on
  $[s_0+r,s_1-r]$ is at most $e\,e^{-\gamma r}$ and
  \[
  \operatorname{length}_P\bigl(p|_{[s_0+r,s_1-r]}\bigr)\le C\sqrt e\,e^{-\gamma r/2},
  \]
  with $C,\gamma$ independent of $s_1-s_0$ and of ${\mathfrak s}_{\rm b}$.
\end{enumerate}
\end{proposition}

\begin{proof}
(i) Let ${\mathfrak s}_2$ be the constant ${\mathfrak s}_0$ of Proposition~\ref{prop:transition-gap}.
The action $a(c)$ does not depend on ${\mathfrak s}$, and the squared gradient norm
$\|F_{B_f}\|^2+{\mathfrak s}^{-1}\|F_{B_c}\|^2+\|\partial_\theta\nu\|^2$ (up to the slice terms)
is nonincreasing in ${\mathfrak s}$. The neighbourhood in the statement does not depend on
${\mathfrak s}$, so the inequality at ${\mathfrak s}={\mathfrak s}_2$ implies it for all ${\mathfrak s}\le{\mathfrak s}_2$. At
${\mathfrak s}={\mathfrak s}_2$ write $c$ as the image of $x$ under the exact matching
chart~\eqref{eq:exact-chart} centred at a point $c_0\in\cF$, with $x$ in the Coulomb slice and in
the linear domain of $A_{{\mathfrak s}_2}$, and orthogonal to $T_{c_0}\cF$; in this chart the
matching and flatness conditions hold exactly, and the decomposition exists for $\delta$ small
because $P$ is compact. The Hessian at $c_0$ is $A_{{\mathfrak s}_2}$, whose kernel is
$T_{c_0}\cF$ and whose spectrum meets $(-\gamma_0,\gamma_0)$ only in $0$, with $\gamma_0$ the gap of
Proposition~\ref{prop:transition-gap}; with the elliptic estimate at the fixed scale
${\mathfrak s}_2$ this gives $\|x\|_{H^1}\le C\|A_{{\mathfrak s}_2}x\|$. Taylor expansion, with the
multiplication $H^1\times H^1\to L^2$ in dimension three and the chart bounds of
Lemma~\ref{lem:block}, gives
\[
|a(c)-\tfrac12\langle A_{{\mathfrak s}_2}x,x\rangle|\le C\|x\|_{H^1}^3,\qquad
\|\nabla_{{\mathfrak s}_2}a(c)-A_{{\mathfrak s}_2}x\|\le C\|x\|_{H^1}^2 .
\]
Hence $\|x\|_{H^1}\le C\|\nabla_{{\mathfrak s}_2}a(c)\|+C\|x\|_{H^1}^2$, so
$\|x\|_{H^1}\le2C\|\nabla_{{\mathfrak s}_2}a(c)\|$ once $\|x\|_{H^1}\le 1/(2C)$, and
\[
|a(c)|\le C\|x\|_{H^1}^2\le C'\|\nabla_{{\mathfrak s}_2}a(c)\|^2 .
\]

(ii) Along a solution $a'=-\sigma^2$ with
$\sigma^2=\|\xi_f\|^2+{\mathfrak s}\|\xi_c\|^2+\|\partial_s\nu\|^2$, the kinetic
energy density. By (i), $a'\le-a/C$ from both ends of the neck, which
gives the exponential decay of the energy. The projection $p$ determined
by the fixed branch satisfies $|p'|\le C\sigma$, because the fixed branch
is measured in the unweighted metric; integrating $\sigma\le\sqrt{-a'}$
over dyadic intervals gives the length bound with no factor depending on
the length of the neck.
\end{proof}

The choice of projection matters: with the nearest-point projection of
the collapsing branch in the unscaled metric, the method only gives the
bound multiplied by $(\inf{\mathfrak s})^{-1/2}$.

\subsection{Convergence to the flat family on the necks}

Let $U_\nu$ be index-zero solutions on $X_{\eps_\nu}$ converging to $V$
as in Proposition~\ref{prop:weak-limit}. By
Theorem~\ref{thm:collapse}(ii) there is no energy loss and no breaking,
so the energy $e_\nu$ of $U_\nu$ on the part of each endpoint neck at
distance at least $s_0$ from the exterior tends to the tail energy of $V$,
which is small for $s_0$ large, and on windows $[S,L-S]$ with
$S\to\infty$ it tends to zero.

\begin{lemma}[convergence to the flat family]\label{lem:entry}
Fix $e_0>0$ small. Let $J$ be an interval in the neck and let $J^+$
contain its closed three-unit neighbourhood in the complete neck
coordinate. Include in $J^+$ the full fixed-radius seam half-disks
used below. If the energy $e$ on this larger region is at most $e_0$,
then, for $\nu$ large, the following hold on $J$, uniformly in its
position and in $\nu$:
\begin{enumerate}[label=(\roman*)]
\item the spatial curvatures satisfy
  $\sup_s\|F_{B_j(s)}\|_{L^2(C)}\le C\sqrt{e}$, $j=f,c$, and there are
  gauges in which $B_j(s)=B_{q_j(s)}+\beta_j(s)$ with $\beta_j$ in the
  Neumann slice orthogonal to the harmonic fields and
  $\|\beta_j(s)\|_{H^1}\le C\sqrt e$;
\item on the collapsing seam the conclusions of
  Lemma~\ref{lem:seam-regularity} hold with $E(b,r_0)\le e$; in particular
  the section is uniformly Hölder on a fixed collar, and
  $\sup_\theta d(u(s,\theta),q_f(s))+d_P(q_c(s),q_f(s))\le C\sqrt e$.
\end{enumerate}
\end{lemma}

\begin{proof}
(i) Near a seam center $b$ use the holomorphic coordinate $\Phi_b$ with
$\Phi_b'(s)=1/{\mathfrak s}(s)$; in it the collapsing compression and the section
have exactly the fixed product mixed geometry, and fixed patches
correspond to regions of size $C{\mathfrak s}_b$. Choose their inner
patches to cover $J$ and their larger patches inside $J^+$, decreasing
the fixed patch radius if necessary. The energy on each larger patch
is at most $e$, by conformal invariance. The mixed mean-value inequality
\cite[Theorem~4.2]{DFL-mixed}, the interior estimate
\cite[Lemma~4.50]{DFL-mixed} and the derivative bounds of
\cite[Lemma~4.25]{DFL-mixed} give the curvature bound after integration
over a finite cover; the spatial slice and a Hodge bootstrap give the
gauges.

(ii) The hypothesis of Lemma~\ref{lem:seam-regularity} is the energy bound
$E(b,r_0)\le e\le e_0$. Its part (ii) gives the Hölder bound, and
integrating the derivative bound across the section interval from the
fixed branch, together with its part (iii), gives the last statement.
\end{proof}

Lemma~\ref{lem:entry} is applied with a fixed small bound $e_0$ rather
than with vanishing energy. Apply it to $[s_0,L-s_0]$ using the energy
on $[s_0-3,L-s_0+3]$ and its seam patches, choosing $s_0$ larger first.
It gives a uniform Hölder bound for the
section on unit patches of the complete section geometry, including the
region between the exterior and the windows. With vanishing energy on the
windows $[S,L-S]$ it gives convergence into any prescribed flat neighbourhood
there, and by Proposition~\ref{prop:neck-decay} and
Proposition~\ref{prop:weak-limit}, the rescaled transitions converge to
the flat pairs at $q_\pm$.

\subsection{First-order estimates on the necks and cuts}

On a neck interval $[a,b]$ with midpoint $s_*$, let $q_{\rm m}$ be the
fixed-branch projection at $s_*$ and $B_*=(B_{q_{\rm m}},B_{q_{\rm m}},q_{\rm m})$ the flat pair
at $q_{\rm m}$, with temporal coefficient zero. Let
$\|z\|^2_{\mathfrak s}=\|U_f\|^2+{\mathfrak s}\|U_c\|^2+\|\nu\|^2$ and $A_{\mathfrak s}$ be as in
Section~\ref{subsec:transition}, and define
$\|z\|_{X^0(J)}^2=\int_J(\|z\|^2_{\mathfrak s}+\|z_s\|_{\mathfrak s}^2+\|A_{\mathfrak s}
z\|^2_{\mathfrak s})$. Put $d(s)=\min(s-a,b-s)$ and $\mathsf w=e^{c\,d(s)}$ with
$2c<\gamma$.

\begin{lemma}[first-order neck estimate]\label{lem:neck-first-order}
For $e$ and the neighbourhood of Lemma~\ref{lem:entry} small, there is a
reference Coulomb gauge on each branch of the neck at $B_*$ in which the
exact matching coordinate $z$ (Section~\ref{subsec:charts}) satisfies
\[
\|\mathsf wz\|_{X^0([a+3,b-3])}\le C\sqrt e,\qquad
\|z\|_{X^0([a+r,b-r])}\le C\sqrt e\,e^{-cr}.
\]
\end{lemma}

\begin{proof}
\emph{Gauge.} Start from the gauges of Lemma~\ref{lem:entry}(i), in which
$B_j(s)=B_{q_j(s)}+\beta_j(s)$ with $\beta_j$ in the spatial Neumann slice.
Differentiating the slice condition in $s$ gives
$d^*_Bd_B\psi_j=d^*_B(\partial_sB_j)$ for the temporal coefficients, so the
spatial Neumann gap gives $\|\psi_j(s)\|\le C\sigma(s)$, where $\sigma$ is the
kinetic energy density of Proposition~\ref{prop:neck-decay}(ii). Since the
energy on $[a+r,b-r]$ is at most $e\,e^{-\gamma r}$ and $2c<\gamma$, summation over
unit intervals gives $\|\mathsf w\psi\|\le C\sqrt e$ in the physical measure
$\int(\|\psi_f\|^2+{\mathfrak s}\|\psi_c\|^2)ds$, and likewise
$\sup_s \mathsf w\,d_P(q_j(s),q_{\rm m})\le C\sqrt e$ by the length bound. For the bound on
$\psi_j$: in the Coulomb form of the local section $q\mapsto B_q$,
$d^*_B\partial_sB_j=O(|\beta_j|(|q_j'|+\|\partial_s\beta_j\|))$, and
$|q_j'|+\|\partial_s\beta_j\|\le C(\sigma+\|\psi_j\|_{H^1})$, so the term is absorbed by the
smallness of $\beta_j$. Let $\chi$ solve the reference Coulomb equation
$d^*_{B_*}\bigl((\exp\chi)^*A-B_*\bigr)=0$ on each branch, in the physical metric,
with Neumann conditions on $\partial C$ and at $s=a,b$; it exists and is small on the
finite neck by the Neumann gap and the implicit function theorem. Test it against $\mathsf w^2\chi$ in the physical one-form measure.
The spatial slice parts pair trivially with $d_C\chi$ up to the difference of
the references $B_{q_j(s)}-B_*$, which is bounded by $C\,d_P(q_j(s),q_{\rm m})$, and the
remaining linear source is the temporal coefficient. Since
$|\partial_T\log \mathsf w|\le c$ in the conformal time, absorption gives
$\|\mathsf w\,d\chi\|\le C\sqrt e$.

\emph{Equation.} In the exact matching chart the equation is
$D_{\mathfrak s} z+N_{\mathfrak s}(z)=0$ with $D_{\mathfrak s}=\partial_s-A_{\mathfrak s}+2\varsigma\Pi_{\rm f}$ and
$\|N_{\mathfrak s}(z)\|_{\mathfrak s}\le C\mathsf d(z)(\|z_s\|_{\mathfrak s}+\|A_{\mathfrak s} z\|_{\mathfrak s})$, where
the coefficient size $\mathsf d(z)$ bounds the spatial coefficients and the $C^0$ norm of the section.
For compactly supported $v$, $\|(\partial_s-A_{\mathfrak s})v\|^2=\|v_s\|^2+\|A_{\mathfrak s}
v\|^2$, since $I_{\mathfrak s} A_{\mathfrak s}$ is independent of $s$ and the Green
pairings cancel; hence $\|v\|_{X^0}\le C(\|D_{\mathfrak s} v\|+\|v\|)$, the
$L^2$ term being kept for the stationary kernel. Apply this to $v=\mathsf w\zeta z$
with a cutoff $\zeta$, and use $\|\mathsf wz\|_{L^2}\le C\sqrt e$, which follows by
applying Lemma~\ref{lem:entry} on unit windows with the local energy
$e\,e^{-\gamma d(s)}$ of Proposition~\ref{prop:neck-decay}(ii), since $2c<\gamma$,
together with the gauge estimate.
\end{proof}

\emph{Cut positions} (Figure~\ref{fig:neck})\emph{.} Fix $D$ larger than all cutoff buffers and choose
$C_*>D+1+\log2$. On each neck, of length $L$ in the exterior cylindrical
coordinate $s$, cut the exterior piece near $s_E=L/2-C_*+D$ and the
transition piece near $s_T=L/2-C_*-D$, and use the same midpoint
$s_*=L/2-C_*$ to define the stationary parameter of both pieces. The parts of $U_\nu$ used for the two pieces therefore overlap on a band of fixed width $2D$ on which
both representatives are the actual solution, and both cuts, with their cutoff bands of unit
width, lie in $\{s\le L/2-\log2\}$, where $E_\eps$ coincides with $X_\eps$. By
Lemma~\ref{lem:neck-first-order}, applied on the interval $[s_0,L-2C_*-s_0]$, whose
midpoint is $s_*$, each exact matching cut has first-order trace of size
$C\sqrt e\,e^{-c(L/2-O(1))}$ at the cut, since each cut lies at distance
$L/2-O(1)$ from both ends of that interval, where $L=\log(1/\eps)$.
With component end weights $\delta\le c/4$ on the new ends, the cut error
in the weighted first-order norms is
\begin{equation}\label{eq:cut-error}
\tau_{\rm cut}\le C\sqrt e\,e^{-cL/2}e^{\delta L/2}\le C\eps^{c/4}.
\end{equation}

\subsection{Comparison on the transitions and the exterior}\label{subsec:cut-comparison}

The comparison estimates are proved in Coulomb gauges on the two kinds of
components. These gauges are constructed with the scalar inverse of Lemma~\ref{lem:scalar-inverse}.

Let $U^E_\nu$ be the \emph{cut exterior piece} of $U_\nu$: equal to $U_\nu$ on the
exterior side of the exterior cut $s_E$ on each neck, and completed for $s\ge s_E+1$
by the flat pair at the midpoint value $q_{\rm m}$ of that neck, by a cutoff of fixed
slope in $s$ applied in the gauge of Lemma~\ref{lem:neck-first-order}. Likewise the
\emph{cut transition} is $U_\nu$ on the transition side of $s_T$, completed by the
flat pair at $q_{\rm m}$.

\begin{lemma}[component Coulomb gauges]\label{lem:component-gauge}
Lemma~\ref{lem:scalar-inverse} holds, with the same proof, on the auxiliary
exterior $E_\eps$ with the metric of Definition~\ref{def:aux-exterior}, the weight
$e^{\phi_E}$ and $\lambda$ equal to ${\mathfrak s}_E$ on the compression, and on the complete
flat transition $\mathcal T$ with the weight $e^{\delta\rho}$. For $\nu$ large there is a
gauge transformation of $U^E_\nu$, trivial at the original ends in the limit, after
which $U^E_\nu=W^E_\eps+\alpha_E$ with
\[
\lambda^2d^*_{W^E_\eps}\alpha_E=0,\qquad \alpha_E(n)=0,\qquad
\|\alpha_E\|_{H^N_{\rm uloc}}\to0,
\]
and $\alpha_E-\mathfrak b_E(v_*)$ decaying at the new ends. Here $v_*=(v_-,v_+)$
records the midpoint values in the Neumann Coulomb coordinates at $q_\pm$, and
$\mathfrak b_E(v_*)$ is a cutoff extension of the actual stationary difference:
its connection component on each new end is $B_{q_\pm}(v_\pm)-B_{q_\pm}$,
using the flat family of Section~\ref{subsec:charts}, and its section component
has the corresponding limit at $q_{\rm m}$. In the exact matching coordinate
$z$, the remainder $z-\mathfrak k_E(v_*)$ decays. The same holds for the cut transition relative to
the flat pair at $q_{\rm m}$, with fixed limit zero.
\end{lemma}

\begin{proof}
The rescaled geometry of $E_\eps$ and of $\mathcal T$ is again a union of product
cylinders with compact cross-section joined through finitely many fixed cores,
$\lambda^{-1}$ is Lipschitz in the rescaled distance because
$|\partial_T{\mathfrak s}_E^{-1}|=|\partial_t\log{\mathfrak s}_E|\le H$, and every slice is irreducible, so the
proof of Lemma~\ref{lem:scalar-inverse} applies verbatim. Take the local gauges of Proposition~\ref{prop:weak-limit}(i) on compact subsets
of $X_0$, of Corollary~\ref{cor:seam-convergence} on the middle compression and of
Lemma~\ref{lem:neck-first-order} on the necks. At each new end choose the
midpoint flat representative in the Neumann Coulomb slice of $B_{q_\pm}$
given by Section~\ref{subsec:charts}. On an overlap of two of them the two
representatives are related by a gauge transformation $g$; on each slice both are
close to the same irreducible flat connection, so after a choice of central sign
$g=\exp\chi$ with $\chi$ small, and the gauge equation bounds $\chi$ by the two
connection differences. Interpolate by $\exp(\theta\chi)$ with a fixed cutoff $\theta$
whose derivative in the product time is bounded; the overlaps form a tree, so the
central signs can be chosen consistently. This gives a representative $W^E_\eps+\alpha$
with $\|\alpha\|_{H^N_{\rm uloc}}$ small. The reference Coulomb equation
$\lambda^2d^*_{W^E_\eps}\bigl((\exp\chi)^*(W^E_\eps+\alpha)-W^E_\eps\bigr)=0$, $\nabla_n\chi=0$, has
derivative $\lambda^2d^*d$ at $\chi=0$, which is invertible on uniformly local spaces
by the first part; the implicit function theorem gives a small solution.
On each new end the nonlinear stationary difference
$B_{q_\pm}(v_\pm)-B_{q_\pm}$ is time-independent, coclosed relative to
$B_{q_\pm}$ and normal-free by its Neumann Coulomb-slice construction.
It therefore contributes no Coulomb source, and the gauge parameter decays.
\end{proof}

\begin{lemma}[transition comparison]\label{lem:transition-comparison}
In the Coulomb gauge of Lemma~\ref{lem:component-gauge}, which retains the midpoint
value $q_{\rm m}$ as a free stationary parameter, the cut transition differs from the flat
pair at $q_{\rm m}$ by at most $C\tau_{\rm cut}$ in the weighted first-order norm.
\end{lemma}

\begin{proof}
The flat pair at $q_{\rm m}$ is a solution. Write the cut transition as the flat pair
plus $z$ in the exact matching chart. Then
\[
\mathcal Dz=-N(z)+e_{\rm cut},\qquad \|N(z)\|_{Y_\delta}\le C\mathsf d(z)\|z\|_{X_\delta},\qquad
\|e_{\rm cut}\|_{Y_\delta}\le C\tau_{\rm cut},
\]
where $e_{\rm cut}$ is the error of the cutoff to the flat pair and $\mathsf d(z)$, the
size of the coefficients, tends to zero; the gauge-fixing component carries no
source because $z$ is in Coulomb gauge. The fixed-limit inverse of
Proposition~\ref{prop:transition-inverse}(iii) at weight $\delta$ gives
$\|z\|_{X_\delta}\le C\delta^{-1}(\mathsf d(z)\|z\|_{X_\delta}+C\tau_{\rm cut})$, and for $\mathsf d(z)$ small
the first term is absorbed.
\end{proof}

\begin{lemma}[exterior comparison]\label{lem:exterior-comparison}
Let $z$ be the exact matching coordinate (Section~\ref{subsec:charts}) of
$W^E_\eps+\alpha_E$, the cut exterior piece $U^E_\nu$ (Section~\ref{subsec:cut-comparison}) in the Coulomb gauge of Lemma~\ref{lem:component-gauge}, relative to the
configuration $W^E_\eps$ of the auxiliary exterior (Definition~\ref{def:aux-exterior}),
with a weight $\delta_*\le c/4$ on the new ends.
Then
\[
\|z\|_{X_E}\le C\bigl(\tau_{\rm cut}+\eps^{1/4}+\eps^{(\beta-\delta_*)/4}\bigr),
\]
where $X_E$ is the norm of Proposition~\ref{prop:exterior-inverse}, and the free end
coordinates of $z$, which are the midpoint parameters of the two necks, obey the same
bound.
\end{lemma}

\begin{proof}
Write $\mathcal F_E(z)=f_E+D_Ez+\mathcal N_E(z)$. The reference error is
$\|f_E\|_{Y_E}\le C(\eps^{1/4}+\eps^{(\beta-\delta_*)/4})$: the first term is the
horizontal velocity ${\mathfrak s}_E^{1/2}|\dot q_R|$ of the flat lift, with
${\mathfrak s}_E\le\sqrt\eps$; the second is the truncation of $V$ at $R=\eps^{1/4}$, which is
supported where $s=L/4+O(1)$, has size $O(R^\beta)$ there, and carries the weight
$e^{\delta_*s}=O(R^{-\delta_*})$. The nonlinear term with comparison coordinate zero
satisfies $\|\mathcal N_E(z)\|_{Y_E}\le C\mathsf d(z)\|z\|_{X_E}$, because the spatial
remainder vanishes on the exact flat family, the temporal coefficient derivative is
quadratic, and the stationary fields are exact flat members. Here $\mathsf d(z)\to0$
is the coefficient size, which tends to zero by Proposition~\ref{prop:weak-limit},
Corollary~\ref{cor:seam-convergence} and Lemma~\ref{lem:entry}. The operator $D_E$ is
inverted from $Y_E$ to $X_E$ with a bound independent of $\eps$ by
Proposition~\ref{prop:exterior-inverse}; the gauge-fixing component of $z$ vanishes by the choice of
gauge, and the cut error, which comes from the cutoff to the flat pair at $q_{\rm m}$ and
is bounded by~\eqref{eq:cut-error}, is measured in the same norms. Absorption gives
the bound. The end coordinate of $z$ is the midpoint parameter
because the same midpoint is used for both cuts.
\end{proof}

\subsection{The global gauge and the uniform estimate on the glued domain}\label{subsec:global-gauge}

\emph{Zeroth-order norm.} We use the rescaled metric $\widehat g_\eps$ and measure $\mu_\eps$ of Section~\ref{subsec:scalar}. A variation $x$ on the glued domain is written $x=\bar x+\mathfrak e_\eps(v_{\rm orig})$, where $v_{\rm orig}$ is the vector of limiting values at the original
ends and $\mathfrak e_\eps$ a fixed extension supported in the tails of the original ends, with zero
temporal and normal components and $\|wD_\eps \mathfrak e_\eps(v)\|\le C|v|$. Put
\[
\|x\|_{Z_\eps}=\|w\bar x\|_{L^2(\mu_\eps)}+\|w\bar x_{\rm sec}\|_{L^2}+|v_{\rm orig}|,
\]
the physical zeroth-order norm of gauge field and section, and use the same
expression for nonlinear configurations relative to $W_\eps$. It contains the
stationary fields along the entire finite necks, which on a neck of length
comparable to $L$ contribute up to $\sqrt L$ times their pointwise size.

\begin{lemma}[global gauge]\label{lem:global-gauge}
In the situation of Proposition~\ref{prop:coverage}, for large $\nu$ there is
a gauge transformation of $U_\nu$ on the glued domain, trivial at the
original ends in the limit, after which $U_\nu=W_\eps+x$ with $x$ in the
reference Coulomb gauge of $W_\eps$, $x$ of finite $X_\eps$ norm, and
\[
\|x\|_{Z_\eps}\le C\sqrt{\mathsf A}\,(\tau_{\rm cut}+\eps^{1/4}+R^{\beta-\delta_*}),\qquad
\|x\|_{H^N_{\rm uloc}(\widehat g_\eps)}\to0 .
\]
\end{lemma}

\begin{proof}
\emph{Components.} By Lemmas~\ref{lem:transition-comparison}
and~\ref{lem:exterior-comparison} the cut transitions and the cut exterior
differ from the pieces of $W_\eps$ by at most $C(\tau_{\rm cut}+\eps^{1/4}+R^{\beta-\delta_*})$ in first-order
norms, and the change from $V$ to its truncation costs $CR^\beta$. A stationary
parameter of size $\rho$ on a neck of length $O(L)$ has physical $L^2$ norm
$C\sqrt{\mathsf A}\rho$; the compression contribution is no larger. This gives the
$Z_\eps$ bound on the part of each piece that is used. The qualitative $H^N_{\rm uloc}$ bound
follows from Proposition~\ref{prop:weak-limit}, Corollary~\ref{cor:seam-convergence}
and Lemma~\ref{lem:entry}.

\emph{Joining.} On an overlap band the two representatives represent
the same solution, $A_T=g^*A_E$. On each slice both are close to the same
irreducible flat connection, so $g$ may be chosen near the identity, $g=\exp\chi$,
and the spatial gauge equation with the Neumann gap gives
$\|\chi(T)\|_{H^1(C)}\le C(\|(A_E-B_q)_{\rm sp}\|+\|(A_T-B_q)_{\rm sp}\|)$; the
temporal part of the gauge equation controls $\partial_T\chi$ in the same way.
Interpolate by $\exp(\theta\chi)$ with a fixed cylindrical cutoff $\theta$, whose
product derivative on the compression is ${\mathfrak s}\theta'$ and hence bounded. The
result is an actual gauge-equivalent representative $W_\eps+\alpha$ with the
two stated bounds. The pieces overlap along a tree of bands, so the choices of central sign on the bands can be made consistently with the marked bundle.

\emph{Coulomb gauge.} The reference Coulomb equation
$\lambda_\eps^2d^{*,g_\eps}_{A^0}((\exp\chi)^*(A^0+\alpha)-A^0)=0$, $\nabla_n\chi=0$,
with $A^0$ the connection of $W_\eps$, has derivative $\mathcal L_\eps$ at
$\chi=0$; by Lemma~\ref{lem:scalar-inverse} and the small $H^N_{\rm uloc}$ norm of
$\alpha$ it has a small solution with
$\|\chi\|_{H^{m+1}_{\rm uloc}}\le C\|\alpha\|_{H^m_{\rm uloc}}$, which belongs to the
weighted space at the original ends. To estimate it in the physical norm, write
$(\exp\chi)^*(A^0+\alpha)-A^0=J(\chi)d_0\chi+\operatorname{Ad}_{e^{-\chi}}\alpha$ with
$J(0)=1$, and test the equation against $w^2\chi$:
\[
0=\int\bigl\langle d_0(w^2\chi),J(\chi)d_0\chi+\operatorname{Ad}_{e^{-\chi}}\alpha\bigr\rangle
d\mu_\eps .
\]
The main term is at least $(1-o(1))\|wd_0\chi\|^2$; the weight terms are of
relative size $\delta$, the maximal weight exponent of the original ends. The stationary part
$\mathfrak e_\eps(v_{\rm orig})$ of $\alpha$ is integrated by parts onto its divergence,
bounded by $C|v_{\rm orig}|$, and $\operatorname{Ad}_{e^{-\chi}}\alpha-\alpha$ is small
times $\chi$. Absorption gives
$\|wd_0\chi\|\le C(\|w\bar\alpha\|+|v_{\rm orig}|)$, and the same bound for the
change of the connection. The section is unchanged. Finite $X_\eps$ norm
holds for each fixed $\nu$ since the region $N_I\cap X_\eps$ is compact.
\end{proof}

\begin{proposition}[uniform estimate on the glued domain]\label{prop:attached-garding}
Uniformly for small $\eps$,
$\|x\|_{X_\eps}\le C(\|D_\eps x\|_{Y_\eps}+\|x\|_{Z_\eps})$ for all $x$ in the
domain, with no normalization of the asymptotic values.
\end{proposition}

\begin{proof}
Write $D_\eps=P_\eps+K_\eps$, with $P_\eps$ the differential part and $K_\eps$ the
perturbation term. Let ${\mathfrak s}_1$ be the threshold of Proposition~\ref{prop:garding}. On a unit
patch of $N_I\cap X_\eps$ on which ${\mathfrak s}\le{\mathfrak s}_1$, the configuration $W_\eps$ is the flat lift of
the path $q(t)$ and $P_\eps$ is the operator of Proposition~\ref{prop:garding} (linearization at
the flat family along $q(t)$, scale with $|{\mathfrak s}'/{\mathfrak s}|\le H$); since $\|\nu\|_{L^2(S)}\le\|z\|_{L^2}$,
that proposition gives
\begin{equation}\label{eq:local-glued}
\|z\|_{X_\eps}^2\le C\bigl(\|P_\eps z\|_{Y_\eps}^2+\|z\|_{L^2}^2\bigr)
\qquad(\operatorname{supp}z\text{ in one patch}).
\end{equation}
On unit patches of the endpoint necks the slice is
$C\sqcup[0,\pi]\sqcup(-C,{\mathfrak s}^2g_C)$, not the model of
Proposition~\ref{prop:garding}, whose far edge is Lagrangian; there
\eqref{eq:local-glued} follows instead from the neck identity
$\|(\partial_\rho-A_{\mathfrak s})z\|^2=\|z'\|^2+\|A_{\mathfrak s}z\|^2$ of
Section~\ref{subsec:transition}, the term $2\varsigma\Pi_{\rm f}$ being absorbed into the $L^2$ term. The
identity is exact on the flat part $s\ge L/4+c_R$; on the remaining patches the seam and
$\dot q$ remainders are absorbed by trace interpolation, as in the proof of
Proposition~\ref{prop:garding}.
On patches where ${\mathfrak s}\ge{\mathfrak s}_1/2$ the rescaled geometry lies in a fixed compact family, and
the ordinary local mixed estimates give \eqref{eq:local-glued} there as well.

Use a partition $\sum\chi_\alpha^2=1$ of bounded overlap by cutoffs that depend
on the complete time coordinate only on the compression (constant on each
fibre $C$), share the same value on the two sides of every seam, and are
extended constantly across a fixed section collar; a cutoff in
$\log({\mathfrak s}/{\mathfrak s}_1)$ has derivative bounded by $C|\varsigma|$. A spatial cutoff on $C$
would have commutator ${\mathfrak s}^{-1}d_C\chi$, which is not bounded uniformly; these
cutoffs preserve the domain and have localized commutator bounds
$\|[P_\eps,\chi_\alpha]z\|_Y\le C\|z\|_{L^2(\operatorname{supp}d\chi_\alpha)}$. Applying \eqref{eq:local-glued} to each
$\chi_\alpha z$ and summing,
\[
\|z\|_{X_\eps}^2\le C\sum_\alpha\|\chi_\alpha z\|_{X_\eps}^2
\le C\sum_\alpha\bigl(\|\chi_\alpha P_\eps z\|_{Y_\eps}^2+\|[P_\eps,\chi_\alpha]z\|_{Y_\eps}^2
+\|\chi_\alpha z\|_{L^2}^2\bigr)\le C\bigl(\|P_\eps z\|_{Y_\eps}^2+\|z\|_{L^2}^2\bigr),
\]
the first inequality because the cutoffs have bounded derivatives and are constant on the
fibres, and the last because the supports of the $\chi_\alpha$ and of the $d\chi_\alpha$ have
bounded overlap. This is the claim for the differential part $P_\eps$, with constants
depending on the overlap and not on the number of patches. The perturbation
term has output and dependence region in the unchanged geometry and is bounded
on $L^2$. Conjugation by the weights of the original ends adds bounded zeroth-order terms, and
subtracting $\mathfrak e_\eps(v_{\rm orig})$ reduces to decaying fields.
\end{proof}

\subsection{Surjectivity of the gluing map}

\begin{proposition}[surjectivity of gluing]\label{prop:coverage}
Let $\eps_\nu\to0$ and let $U_\nu$ be index-zero solutions on
$X_{\eps_\nu}$ converging, in the sense of
Proposition~\ref{prop:weak-limit}, to the regular index-zero solution
$V$. Then $U_\nu=\cG_{\eps_\nu}(V)$ for all large $\nu$.
\end{proposition}

\begin{proof}
We write $\eps=\eps_\nu$; the sequence enters only through thresholds.

\emph{Step 1: global gauge.} By Lemma~\ref{lem:global-gauge}, after a gauge
transformation $U_\nu=W_\eps+x$ with $x$ of finite $X_\eps$ norm,
$\|x\|_{Z_\eps}\le C\sqrt{\mathsf A}(\tau_{\rm cut}+\eps^{1/4}+R^{\beta-\delta_*})$ and coefficient size
$\mathsf d(x)\to0$.

\emph{Coordinates.} Write $a_{\rm aff}$ for the affine difference
produced by Lemma~\ref{lem:global-gauge}. Its coefficients tend to
zero, and its higher rescaled local norms are bounded. The rescaled
H\"older inverse of the chart therefore applies. Write $x$ for the
resulting exact matching coordinate; it is smooth, has finite
$X_\eps$ norm, and $\mathsf d(x)\to0$. We need only
\begin{equation}\label{eq:weak-chart-comparison}
\|x\|_{Z_\eps}\le C\|a_{\rm aff}\|_{Z_\eps}
                     +C\mathsf d(x)\|x\|_{X_\eps}.
\end{equation}
Indeed $a_{\rm aff}-x$ is the sum of the quadratic term in the harmonic coordinate,
the transport correction, and $\mathcal R(v,w)$, with the analogous
section correction. The first two are bounded in $Z_\eps$ by
$C\mathsf d(x)\|x\|_{X_\eps}$. For the third use the trace and
extension estimates \eqref{eq:anisotropic-trace} at order zero and
the surface estimate \eqref{eq:class-map-L2}. Its $Z_\eps$ norm
has the factor ${\mathfrak s}^{1/2}$, so the bound follows from the
first-order correction estimate. The section correction has the
ordinary $C^0$ quadratic bound. This proves
\eqref{eq:weak-chart-comparison}.
The gauge-fixing row is the affine Coulomb condition composed with
the chart.

\emph{Step 2: a priori estimate in the first-order norm.} Write
$\mathcal F_\eps(W_\eps+x)=f_\eps+D_\eps x+N_\eps(x)$ with $\|f_\eps\|\le Ce_\eps$
by~\eqref{eq:approx-error} and $\|N_\eps(x)\|_{Y_\eps}\le C\mathsf d(x)\|x\|_{X_\eps}$:
the pure harmonic spatial term vanishes on the exact flat family
by~\eqref{eq:exact-chart}, every remaining singular term contains a
complementary field, which carries the weight ${\mathfrak s}^{-1}$ of the norm, and the section uses only
$C^0$ smallness. Proposition~\ref{prop:attached-garding},
\eqref{eq:weak-chart-comparison}, and \eqref{eq:low-quadratic} give
$\|x\|_{X_\eps}\le C(e_\eps+\|a_{\rm aff}\|_{Z_\eps})+
C\mathsf d(x)\|x\|_{X_\eps}$. Absorbing the last term gives
\[
\|x\|_{X_\eps}\le C\bigl(e_\eps+\sqrt{\mathsf A}(\tau_{\rm cut}+\eps^{1/4}+R^{\beta-\delta_*})\bigr)
\le C\eps^{\sigma_{\rm c}}
\]
for some $\sigma_{\rm c}>0$, since $\tau_{\rm cut}$, $e_\eps$ and $R^{\beta-\delta_*}$ are positive powers of
$\eps$ by~\eqref{eq:cut-error} and~\eqref{eq:approx-error}.

\emph{Step 3: interpolation.} On four-dimensional unit patches of the
conformal geometry, which are copies of one fixed compact manifold,
interpolation between the local $H^N$ bounds, which hold uniformly by the
convergence of Proposition~\ref{prop:weak-limit},
Corollary~\ref{cor:seam-convergence} and Lemma~\ref{lem:entry} and are preserved by the gauges of Step~1, and the
first-order bound of Step~2 gives $\|U_\perp\|_{C^k}\le C\eps^{\sigma_{\rm c}\theta_1}$
for the complementary gauge field, using
$\int{\mathfrak s}^{-1}\|U_\perp(t)\|^2_{H^1}dt=\int\|U_\perp(T)\|^2_{H^1}dT$ with $dt={\mathfrak s}\,dT$,
so that no power of the minimal scale is lost. The slow coefficient is not
interpolated on product patches; its supremum is obtained from the
section through the matching identity $v=\Pi\nu|_{\mathrm{seam}}-
\Pi\tau_t(U_\perp|_\Sigma)$. The section is interpolated on two-dimensional
patches between the uniform Hölder bound of Lemma~\ref{lem:entry} and the
$L^2$ bound. Hence the coefficient size satisfies $\mathsf d(x)\le C\eps^{\sigma}$
for some $\sigma>0$.

\emph{Step 4: uniqueness.} The norm $X_\eps$ on the glued domain carries
only the original end weights; the artificial weights of the component
problems do not enter. Let $z^*=\cG_\eps(V)-W_\eps$, so that
$\|z^*\|_{X^M}\le a:=C_M\mathsf A^{M+1}e_\eps$ by Proposition~\ref{prop:existence}.
The difference $y=x-z^*$ satisfies $D_\eps y=N_\eps(z^*)-N_\eps(x)$. The
difference estimate holds in the asymmetric form
\[
\|N_\eps(x)-N_\eps(z^*)\|_{Y_\eps}\le C\bigl(\mathsf d(x)+a\bigr)\|x-z^*\|_{X_\eps},
\]
where $\mathsf d(x)$ is the unweighted coefficient size of $x$ from Step~3 and
the weighted suprema ${\mathfrak s}^{-1/2}|p_{z^*}|$, ${\mathfrak s}^{-1/2}|\nabla_C(U_\perp)_{z^*}|$ needed by
the singular terms, such as ${\mathfrak s}^{-1}[p_{z^*},h_{v_x-v_{z^*}}]$, are placed on the
factor $z^*$, where they are bounded by $a$ through~\eqref{eq:strong-control}.
The inverse of Proposition~\ref{prop:attachment-inverse} then gives
\[
\|y\|\le C\mathsf A\bigl(\mathsf d(x)+a\bigr)\|y\|.
\]
Both $\mathsf A\,\mathsf d(x)$ and $\mathsf A\,a$ are positive powers of $\eps$ multiplied by powers
of $\mathsf A$, hence tend to zero, and $y=0$ for large $\nu$.
\end{proof}

The unknown solution never has to enter the logarithmically shrinking
ball of Proposition~\ref{prop:existence} directly. It first reaches a
positive power of $\eps$ by estimates free of logarithms, and only then
meets the single factor $\mathsf A$ of the inverse on the glued domain.

\subsection{Proof of Theorem~\ref{thm:collapse}(iii)}

The order of choices is: the data of the pieces and the bound $b_\eta$; the neck
lengths of the ordinary necks, for which
Hypothesis~\ref{hyp:collapse}\ref{hyp:X4} holds, and then $N_I$, which may meet the necks; the finite set
$\cM^0(X_0)$, which fixes the decay rates $\beta_V$ and the weights; and
finally $\eps$. Negative indices: if solutions of negative index existed
for a sequence $\eps_\nu\to0$, a subsequence would converge by
Theorem~\ref{thm:collapse}(ii); the energy bound of part (i) leaves
finitely many negative indices, and the loss formula with
Hypothesis~\ref{hyp:collapse}\ref{hyp:X4} excludes each. For index zero,
$\cG_\eps$ is defined by Proposition~\ref{prop:existence} and is injective for
small $\eps$: $\cM^0(X_0)$ is finite, and $\cG_\eps(V)\to V$ as $\eps\to0$ in the
sense of Proposition~\ref{prop:weak-limit}, since
$\|\cG_\eps(V)-W_\eps\|_{X^m}\to0$ and $W_\eps$ agrees with $V$ away from
shrinking neighbourhoods of the rectangle. If it were not surjective for
a sequence $\eps_\nu\to0$, a sequence of solutions outside its image would
converge to some $V\in\cM^0(X_0)$ by Theorem~\ref{thm:collapse}(ii),
contradicting Proposition~\ref{prop:coverage}. Regularity of the image is
part of Proposition~\ref{prop:existence}. \qed

\section{Two-handles}\label{sec:two-handles}

\subsection{Heegaard triples and the three Lagrangians}\label{subsec:triple-data}

Let $(W_2,E,\nu)\colon Z_0\to Z_1$ be a two-handle cobordism along a framed
link $L$ of $n$ components, with $E$ trivial at both ends
(Proposition~\ref{prop:outer-bundles}), with a filling frame
(Definition~\ref{def:normal-form}) of lifting class $\zeta$. Choose a
subordinate Heegaard triple $(\Sigma;H_0,H_1,H_2)$ of genus $g$ adapted
to the ball (Section~\ref{subsec:triples}), and take the boundary sum with the fixed torus: $C_i=H_i\natural([0,1]\times T^2)$, $h=g+1\ge2$ (after a
stabilization in a small ball of the standard collar, away from the bouquet, if necessary). Put
$r=g-n$. The transition $\rho_\zeta$ of the filling frame on $\Sigma$ induces a
symplectomorphism $\tau_\zeta$ of $\Msf_h$, and the three Lagrangians are
\begin{equation}\label{eq:triple}
(\tau_\zeta L_{\mathfrak h_0}(C_0),\ L(C_1),\ L_{\mathfrak h_2}(C_2))=:(P_0,P_1,P_2),
\end{equation}
where the primary perturbation of $C_1$ is zero, that of $C_2$ is the height
function $t\,h_Q$ of the next subsection, and that of $C_0$ is small and chosen
so that $P_0$ is transverse to $P_1$ and $P_2$ and misses $P_1\cap P_2$.
Since $\rho_\zeta$ is an automorphism of the $\SO(3)$ bundle, it conjugates
every operator in the problem; in the trivialization of the filling frame over $e_0\times H_0$ the boundary problem of $C_0$ is literally the untwisted one, and $\zeta$ enters only
through the global topology of the triangle.

\subsection{The auxiliary pair}

The manifold $Y_{12}=H_1\cup H_2$ is $\#^r(S^1\times S^2)$ and, for vanishing
primary perturbations of $C_1$ and $C_2$, the pair $(L(C_1),L(C_2))$
intersects cleanly in $Q\cong\SU(2)^r$, with $TQ=TL(C_1)\cap TL(C_2)=
H^1(Y_{12}^\#;\ad)$. Here the flat connections of $Y_{12}^\#$ are determined
by the holonomies $b_j\in\SU(2)$ of the $r$ free generators of $\pi_1$, with
lifts fixed by the global frame of $H_1\cup H_2$; the torus summand
contributes the fixed quaternion pair, and the nonzero element of $\phi_M$ identifies its two sign choices. The primary perturbation of $C_2$ is the
height function $t\,h_Q$, $h_Q=\sum_j(1-\tfrac12\tr b_j)$, with $t>0$ small,
and the primary perturbation of $C_1$ is zero; transversality of
$(P_0,P_2)$ and empty triple intersection are obtained by perturbing $C_0$
alone. The critical points are the $2^r$ points $a_S$ with $b_j=-1$ for
$j\in S$ and $b_j=1$ otherwise; $a_S$ has Morse index $m(a_S)=3|S|$, and there
is a unique maximum $q_*=a_{\{1,\dots,r\}}$.

For these auxiliary data we distinguish the relative class from the
total perturbed energy. For a gauge, mixed or strip configuration
from $a$ to $x$, put
\begin{equation}\label{eq:aux-class}
\kappa_{\mathrm{aux}}=
\kappa-\frac{t}{4\pi^2}\bigl(h_Q(a)-h_Q(x)\bigr).
\end{equation}
Here $\kappa$ is the topological energy of
Section~\ref{subsec:energy-convention}. Thus, for a strip,
$4\pi^2\kappa_{\mathrm{aux}}=
\int u^*\omega-t(h_Q(a)-h_Q(x))$.
The relative class vanishes on a reference path in the clean
intersection $Q$ and takes values in $\tfrac12\ZZ$.
For a filling with outgoing limit $a$, we instead normalize
$\kappa_{\mathrm{aux}}$ to vanish on its flat disk-extendible
reference with that limit. The following proof verifies these
normalizations and the associated index formulas.

\begin{lemma}\label{lem:aux-index}
Use the conventions under which the anti-self-duality equation on an
outgoing end $[0,\infty)\times Y$ is the downward gradient flow of the
Chern--Simons functional perturbed by $+t\,h_Q$ with $t>0$.
\begin{enumerate}[label=(\roman*)]
\item For the gauge filling $F_r=\natural^r(S^1\times D^3)\natural(T^2\times D^2)$
  of $Y_{12}^\#$ with outgoing limit $a$, on the component of the flat
  fillings, $\ind=3r-m(a)$; other components have index
  $3r-m(a)+8\kappa_{\mathrm{aux}}$.
\item For a gauge, mixed or ordinary configuration from $a$ to $x$, with
  the relative class just defined,
  $\ind=m(a)-m(x)+8\kappa_{\mathrm{aux}}$, and the full analytic energy is
  $8\pi^2\kappa_{\mathrm{aux}}+2t\,(h_Q(a)-h_Q(x))$
  plus secondary terms.
\item If $2t\max h_Q+2\|f\|_\infty+
  2\sup_A\|\eta(A)\|_2^2<4\pi^2$ on each counted domain,
  every counted configuration has $\kappa_{\mathrm{aux}}\ge0$.
\end{enumerate}
Reversing one of the three conventions (the sign of the flow, the sign of
$t$, outgoing versus incoming end) replaces $3r-m(a)$ by $m(a)$ and the
maximum by the minimum throughout; reversing two restores the statement.
\end{lemma}

\begin{proof}
(i) At a flat connection $A$ on $F_r$ with the unperturbed end, the
operator $-d_A^*\oplus d_A^+$ on decaying spaces is injective and
surjective: a decaying kernel element is closed and represents a class in
the image of $H^1_c(F_r;\ad)\cong H_3(F_r;\ad)=0$, hence is exact, $a=d_Af$,
and $f=0$ by $H^0(Y_{12}^\#;\ad)=0$; a cokernel element has vanishing
function part by $H^0(F_r;\ad)=0$ and a self-dual part that is exact by
$H^2(F_r;\ad)=H^2(T^2;\ad)=0$, hence decays and vanishes by
$\|\omega\|^2=\int\omega\wedge\omega=0$. So the decaying index is $0$. The
asymptotic operator on the tangent space $TQ$ of the critical manifold
moves from $-\delta$ to $t\operatorname{Hess}_a h_Q$, which has $m(a)$ negative
and $3r-m(a)$ positive eigenvalues; for an outgoing end each eigenvalue
crossing from negative to positive increases the index by one. Hence
$\ind=3r-m(a)$. This is a spectral-flow calculation; it does not
identify solutions of the perturbed four-dimensional equation with
finite-dimensional gradient trajectories in $Q$.
Relative to the flat reference, a point insertion changes
$\kappa_{\mathrm{aux}}$ by $1$ and the index by $8$. The possible
nonextending meridional lifting class changes these by a half-integer
and an integer congruent to $4$ modulo $8$, respectively.
Indeed, the filling retracts onto $T^2\vee\bigvee^rS^1$, and
the quotient of the boundary lifting classes by those extending
over it is generated by the torus meridian. With that character
fixed, relative $H^4$ records the point insertion. The
Pontryagin-square congruence gives the stated half-integral lattice
and the formula in (i).

(ii) Join $a$ to $x$ by a path of flat classes in $Q$, constant
near its endpoints. Before the height perturbation, the cylinder
with positive weights at both clean ends has index $-3r$.
The physical incoming and outgoing conditions admit $m(a)$ and
$3r-m(x)$ small modes, respectively. Hence its perturbed index
is $m(a)-m(x)$. Formula (2.52) of \cite{DFL} gives its
topological energy $t(h_Q(a)-h_Q(x))/(4\pi^2)$.
Changing the relative class changes energy and index by
$k/2$ and $4k$, with $k\in\ZZ$.
Gluing this reference to the constant mixed pair at $x$ and
using \cite[Lemma~4.62, Proposition~4.65]{DFL} proves the mixed
formula. Gluing the constant mixed pair at $a$ to a strip and
using \cite[Proposition~4.69]{DFL} proves the strip formula.
In particular, the same integer function $m$ serves on both sides;
it is the relative class $\kappa_{\mathrm{aux}}$, rather than
the total perturbed energy $\kappa$, which occurs in this formula.
The full energy identity follows from~\eqref{eq:aux-class} and
the identity in the proof of Theorem~\ref{thm:mixed-iso}.

(iii) On the outgoing filling collar, let $\beta$ and $\rho$
be increasing cutoffs switching on $t h_Q$ and $f$, respectively.
Integration by parts gives the full energy
\[
8\pi^2\kappa_{\mathrm{aux}}-2t h_Q(a)
 +2t\int\beta'h_Q(A_s)\,ds
 +2\int\rho'f(A_s)\,ds+2\|\eta(A)\|_2^2.
\]
Since $0\le h_Q\le\max h_Q$ and both cutoffs have total
variation one, the absolute value of the height correction is
at most $2t\max h_Q$. The same bound applies to the endpoint
term in (ii). The remaining corrections are bounded by
$2\|f\|_\infty+2\sup_A\|\eta(A)\|_2^2$.
Nonnegative analytic energy and
$\kappa_{\mathrm{aux}}\in\tfrac12\ZZ$ therefore exclude every
negative relative class under the stated inequality.
\end{proof}

With Scaduto's $\ZZ/4$ grading, $\mathrm{gr}=m-r$ modulo four on every
generator \cite[Section~7.3]{Scaduto}; in particular $q_*$ has degree $2r$, the
degree of the relative invariant of $\natural^r(S^1\times D^3)$
\cite[Sections~7.7--7.8]{Scaduto}.

\begin{lemma}[top generator]\label{lem:calibration}
Let $C,C'$ be finite complexes over $\FF$ with bases indexed by the same
set, graded by $m$, and let $N\colon C\to C'$ be an isomorphism whose
matrix entries $N_{a,x}$ vanish unless $m(x)=m(a)+4k$ with $k\ge0$.
Then $N$ maps the unique generator of maximal $m$ to itself.
\end{lemma}

\begin{proof}
If $a$ has maximal $m$, then $N(a)$ is a multiple of $a$, and the
multiple is nonzero because $N$ is injective.
\end{proof}

The hypothesis $k\ge0$ cannot be replaced by a congruence modulo four:
for $r=4$ the isomorphism $N=I+E_{\{1,2,3,4\},\varnothing}$, which adds
the minimum to the image of the maximum, satisfies the congruence and
violates the conclusion.

\begin{definition}
On $\CF(P_1,P_2)$ define the relative grading $m$ by $m(a)-m(x)=\ind(u)-
8\kappa_{\mathrm{aux}}(u)$ for any strip class $u$ from $a$ to $x$; it is independent of $u$
by monotonicity, and for $t$ small it is the Morse index of $h_Q$. Let
$\theta_S\in\CF(P_1,P_2)$ be the unique generator of maximal $m$; it is the
maximum $q_*$ of $h_Q$. Let $\theta_G\in C_G(Y_{12})$ be the relative invariant
of the filling $F_r$: the count of index-zero solutions on $F_r$ with
cylindrical end $Y_{12}^\#$, with lifts by even gauge transformations, projected to the framed
complex.
\end{definition}

\begin{proposition}\label{prop:filling-classes}
$\theta_G=q_*$ and $\theta_S=q_*$ as chains; both are cycles, their
classes are nonzero, and $N_{12}\theta_G=\theta_S$.
\end{proposition}

\begin{proof}
By Lemma~\ref{lem:aux-index} an index-zero filling has
$\kappa_{\mathrm{aux}}=0$ and $m(a)=3r$, so $a=q_*$.
First set $f=\eta=0$. The filling energy identity in that lemma
then gives
$E=-2t\max h_Q+2t\int\beta'h_Q(A_s)\,ds\le0$.
Thus $E=0$. In temporal gauge the end is stationary at the flat
maximum, and the connection on the unperturbed core is flat.
The boundary fundamental group surjects onto that of the filling,
so this boundary class has exactly one flat extension. Its fixed-limit
operator is invertible. Part (i) gives index zero; for injectivity use the
triangular splitting of Section~\ref{subsec:transition} on the end. The
stationary coefficient satisfies $a'=-t\beta\,\mathrm{Hess}(h_Q)a+O(t^2)a$ with
$\mathrm{Hess}(h_Q)\le-c$ at the maximum $q_*$, so
$\frac{d}{ds}|a|\ge t\beta(c|a|-Ct\sup|a|)$, which excludes an interior maximum of $|a|$;
since $a\to0$, this forces $a\equiv0$; the remainder lies in the decaying space, where part (i) gives
an inverse, and the $O(t)$ perturbation is absorbed. For sufficiently small secondary
regularizations this is still the complete count: any additional
index-zero solutions along data tending to zero would have energy
tending to zero, hence no bubbles or nonconstant end trajectories,
and would enter the inverse-function neighbourhood of that solution.
Therefore $\theta_G=q_*$. (This
agrees with \cite[Sections~7.7--7.8]{Scaduto}, where the relative invariant of
$\natural^r(S^1\times D^3)$ is a basis element of the free group
$\Ish(\#^r(S^1\times S^2))$ of rank $2^r$, by the boundary-sum identity
\cite[Section~7.2]{Scaduto} and \cite[Proposition~5.4]{KM-unknot} at the
connected cut $T^3$.) A trajectory from $q_*$ has index
$3r-m(x)+8\kappa_{\mathrm{aux}}$ with $\kappa_{\mathrm{aux}}\ge0$, which is never one, so $q_*$ is a
cycle on both sides. By Lemma~\ref{lem:aux-index} the map $N_{12}$ of
Theorem~\ref{thm:mixed-iso} satisfies the hypothesis of
Lemma~\ref{lem:calibration}, so $N_{12}q_*=q_*$. The filtrations on the
two sides are defined by the same function $m$, so the gauge filling
class and the direction of triangularity of $N_{12}$ cannot come apart
under a change of conventions.
\end{proof}

\subsection{The torus-summed triangle and its lifting group}\label{subsec:torus-summed-triangle}

Let $X^\#$ be the torus-summed triangle: the union of $\Delta\times\Sigma_h$ and
$\Delta\times T^2$ (the two surface sheets) with the edge pieces
$e_i\times C_i$, where each torus compression body $C_i$ is glued to both
sheets along $e_i\times(\Sigma_h\sqcup T^2)$. Its ends are $Y_{ij}^\#=C_i\cup
(-C_j)$. By Proposition~\ref{prop:tube} and the construction of
Section~\ref{sec:category}, removing the tube from $W_2$ and inserting
$T^\circ\times[-1,1]$ gives $W_2^\#=X^\#\cup_{Y_{12}^\#}F_r$ with its bundle. The
classes of $H^1(X^\#;\FF)$ that vanish on both surface sheets restrict to zero
on each compression piece, since $H^1(C_i)\to H^1(\Sigma_h)$ is injective, so
they come from the incidence graph of the torus region, a pair of pants with
three boundary circles: this group is $\Phi=\FF^3/\langle(1,1,1)\rangle$, whose
restriction to the end $Y_{ij}^\#$ is $s_{ij}=t_i+t_j$. Each class is realized by
an automorphism supported on the torus region and equal to the identity near
its joining tree.

\subsection{The gauge two-handle map as a triangle map}

Let $E^\#$ be the bundle on the torus-summed triangle $X^\#$ of
Section~\ref{subsec:torus-summed-triangle}, and
$\Phi=\FF^3/\langle(1,1,1)\rangle$ its lifting group, with $s_{ij}=t_i+t_j$. Let $\mu_G^E$ be the $\Phi$-framed gauge
triangle map $C_G(Y_{01})\otimes C_G(Y_{12})\to C_G(Y_{02})$.

\begin{proposition}\label{prop:gauge-triangle}
Let $F_G(W_2,E)$ be the map of \cite[Section~5]{KM-unknot} of
$(W_2^\#,E^\#)$ with the subgroup $\phi_{W_2}$, between the complexes of
$(Z_i^\#,E^\#|_{Z_i^\#})$ with the markings of the filling frame. Then
$F_G(W_2,E)=[\mu_G^E(\,\cdot\otimes\theta_G)]$ on homology.
\end{proposition}

\begin{proof}
By Proposition~\ref{prop:tube},
$W_2^\#=X^\#\cup_{Y_{12}^\#}F_r$, where $F_r$ is the gauge filling of
Lemma~\ref{lem:aux-index}(i), with the tube in the product region.
Write $W_2^\#=X^\#\circ W'$ with $W'=(Y_{01}^\#\times[0,1])\sqcup F_r$,
cut along the disconnected manifold $Y_{01}^\#\sqcup Y_{12}^\#$. The
restriction $H^0(W')\oplus H^0(X^\#)=\FF^2\oplus\FF\to H^0(Y_{01}^\#\sqcup
Y_{12}^\#)=\FF^2$ is $(a,b,c)\mapsto(a+c,b+c)$, which is onto, so the
connecting image $\psi$ of \cite[Proposition~5.4]{KM-unknot} vanishes.
The subgroup of $W_2^\#$ restricts to the index-two subgroup
$\{s_{12}=0\}\subset\Phi$ on $X^\#$ and to $0$ on $F_r$, since the
torus meridian class does not extend over $T^2\times D^2$. By
\cite[Proposition~5.4]{KM-unknot} and the K\"unneth isomorphism of
\cite[Section~5.2]{KM-unknot},
$F_G(W_2,E)=\mu^{\{s_{12}=0\}}(\,\cdot\otimes\tilde\theta_G)$, where
$\tilde\theta_G$ is the relative invariant of $F_r$ (as in the definition of $\theta_G$),
computed modulo even gauge transformations only. Finally let
$g\in\Phi$ with $s_{12}=1$, let $\iota$ be its action on the critical points of $Y_{12}^\#$, and
let $\pi$ be the projection to the $\Phi$-quotient complex. With $\Phi$-invariant data, $g$ maps
solutions with $Y_{12}$-limit $b$ bijectively to those with limit $\iota b$,
preserving the orbits at the other two ends, whose restriction of
$\{s_{12}=0\}$ is already onto. So $g$ acts freely on the lifts at $Y_{12}$ and trivially on the
orbits at the other ends, and each $\Phi$-orbit of solutions contains exactly one solution for
each lift of $b$. Hence
$\mu^{\{s_{12}=0\}}(a\otimes b)=\mu^\Phi(a\otimes\pi b)$ with no factor
two, and $\mu^{\{s_{12}=0\}}(\,\cdot\otimes\tilde\theta_G)=\mu_G^E(\,\cdot
\otimes\theta_G)$. Nontrivial $E$ is covered by
\cite[Section~5.1]{KM-unknot}, which allows an arbitrary bundle.
\end{proof}

\subsection{The mixed triangle}\label{subsec:mixed-triangle}

\begin{definition}\label{def:mixed-triangle}
The \emph{mixed triangle} is the following domain for the mixed equation.
Let $\Delta$ be a triangle with incoming vertices $v_{01},v_{12}$, outgoing
vertex $v_{02}$ and edges $e_i$, and let $\gamma$ be a proper arc from the
interior of $e_0$ to the interior of $e_2$ separating $v_{02}$ from the incoming
vertices, dividing $\Delta$ into $G$ (containing $v_{01},v_{12}$) and $S$. The
gauge part is the restriction of $X^\#$ over $G$: the two surface sheets over
$G$ and the edge pieces over $e_i\cap G$. The section is a map $S\to\Msf_h$, with
boundary conditions $P_0$ on $e_0\cap S$ and $P_2$ on $e_2\cap S$, matched along
$\gamma$ to both sheets (the torus sheet matching the one-point moduli space).
At $\gamma\cap e_0$ and $\gamma\cap e_2$ the domain has the mixed ends of
$C_0$ and $C_2$. The ends are the gauge inputs $Y_{01}^\#$, $Y_{12}^\#$, the strip
output at $P_0\cap P_2$, and the two mixed ends. Metrics are products near
the seams, and the section almost complex structure is the complex structure $*_\Sigma$ near $\gamma$.
Let $B\subset G$ be the rectangle joining a compact subinterval $I$ of the middle edge $e_1$
to an arc $J\subset\gamma$; the collapse of Theorem~\ref{thm:collapse} takes place over $B$.
\end{definition}

\begin{figure}[ht]
\centering
\begin{tikzpicture}[scale=1.0,font=\small]
\coordinate (A) at (-3,0); \coordinate (B) at (3,0); \coordinate (Cc) at (0,4.2);
\fill[gray!12] (-1.43,2.2) .. controls (-0.5,1.75) and (0.5,1.75) .. (1.43,2.2) -- (Cc) -- cycle;
\fill[orange!18] (A) -- (B) -- (1.43,2.2) .. controls (0.5,1.75) and (-0.5,1.75) .. (-1.43,2.2) -- cycle;
\fill[orange!45] (-0.35,0) -- (0.35,0) -- (0.35,1.88) -- (-0.35,1.88) -- cycle;
\draw[thick] (A) -- (B) -- (Cc) -- cycle;
\draw[very thick,blue!70!black] (-1.43,2.2) .. controls (-0.5,1.75) and (0.5,1.75) .. (1.43,2.2);
\node[blue!70!black] at (0.95,1.68) {$\gamma$};
\node at (0,2.9) {$S$};
\node at (-1.55,0.75) {$G$};
\node at (0,0.95) {$B$};
\node[below] at (0,0) {$e_1$};
\node[left] at (-1.95,1.3) {$e_0$};
\node[right] at (1.95,1.3) {$e_2$};
\node[left] at (-0.85,3.25) {$P_0$};
\node[right] at (0.85,3.25) {$P_2$};
\node[below left] at (A) {$v_{01}$: $Y_{01}^\#$};
\node[below right] at (B) {$v_{12}$: $Y_{12}^\#$};
\node[above] at (Cc) {$v_{02}$: strip output at $P_0\cap P_2$};
\fill (-1.43,2.2) circle (1.6pt); \fill (1.43,2.2) circle (1.6pt);
\end{tikzpicture}
\caption{The mixed triangle of Definition~\ref{def:mixed-triangle}. The gauge part lies
over $G$ and the section over $S$; they are matched along the arc $\gamma$, whose
endpoints are the mixed ends of $C_0$ and $C_2$. The rectangle $B$ joins a
subinterval of the middle edge $e_1$ to $\gamma$. Collapsing $B$
(Theorem~\ref{thm:collapse}) and stretching the ordinary necks gives the symplectic
triangle composed with $N_{01}\otimes N_{12}$; stretching along the closed
three-manifold $Y_{02}^\#$, which by Lemma~\ref{lem:cut-crosses} passes through $B$,
gives $N_{02}$ composed with the gauge triangle.}
\label{fig:mixed-triangle}
\end{figure}
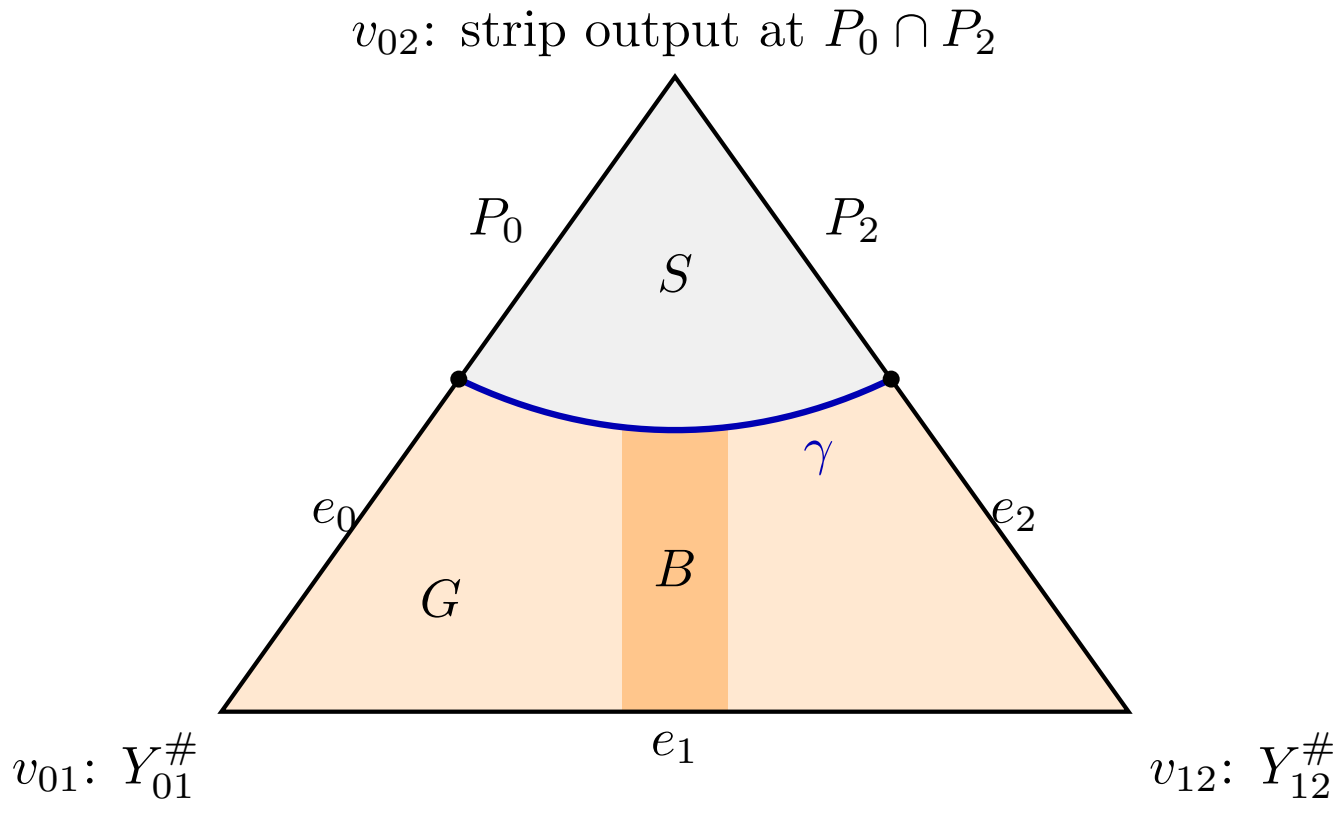

\begin{lemma}[vertex identification]\label{lem:vertex-identification}
Let $D$ be any domain of DFL type for a pair $(P,P')$ covered by
Theorem~\ref{thm:mixed-iso} and Remark~\ref{rem:clean-pairs}, with the same
end types as the special quintuple of \cite[Section~3.1]{DFL}, and regular
invariant data of the following form: near each seam the metric is of the form
used in the proof of Theorem~\ref{thm:family}, conformal near each seam to the product metric
\cite[(1.1)]{DFL-mixed}, that is, with $e^{-f}c\equiv1$ in the notation there, and the
section almost complex structure is $*_\Sigma$; on each mixed end the metric is
conformal to a product $d\theta^2+g_{\rm slice}$. The translation-invariant data on the mixed ends may differ
from those of the quintuple. Then its index-zero count is chain homotopic to the vertex
map of Theorem~\ref{thm:mixed-iso} for that pair.
\end{lemma}

\begin{proof}
First give each mixed end its cylindrical conformal representative.
Conformal invariance preserves the solution space and surjectivity
modulo gauge; the Coulomb row and weighted realization are chosen
for this representative. The space of these domains and data,
with fixed end types, is connected: the
conformal structure of the section region, the placement of the seam, the
metrics near the mixed ends and the neck lengths vary in convex sets, and the
changes are compactly supported. Theorem~\ref{thm:family} gives a chain
homotopy along a path, and Theorem~\ref{thm:mixed-iso}(ii) handles the choices
of $f$, $\chi$ and $\eta$. A change of the translation-invariant data on a mixed
end (the polar geometry of Lemma~\ref{lem:conformal-models}(ii) against the
product geometry of the quintuple, matched along the seam) is a family whose
parameter derivative is independent of the end coordinate there; the flat
references on $C$ and the constant sections do not depend on the end metric,
and the mixed ends have a uniform gap on the complement of the stationary
kernel. Lemma~\ref{lem:family-stationary} then gives the chain homotopy once its
hypothesis of uniform energy decay on the mixed ends is verified. By conformal
invariance of the two equations and their energies (Lemma~\ref{lem:conformal-models})
we compute in the product metric on each mixed end. We choose the path
within the data of the statement, with the complex structure of the section on each
mixed end equal to the standard one pulled back by $(\sigma,\theta)\mapsto(\varphi(\sigma),\theta)$,
where $\varphi$ is a diffeomorphism of the section interval equal to the identity near the
seam and to a translation near the Lagrangian arc; changes of the width of the interval
are made by such maps. Along a compact path of these data the three steps of the proof
of Proposition~\ref{prop:growing-ends}(i) hold with constants independent of the
parameter. First, local compactness and the seam estimates hold with constants
uniform over $C^\infty$-compact families of such data (the compactness paragraph of the
proof of Theorem~\ref{thm:family}), and Gromov compactness is uniform for $J$ in a
compact family; so there is a threshold, independent of the parameter, below which the
kinetic energy on an interval of fixed length forces the slices into a fixed neighbourhood
of the critical family. Second, the local action on that neighbourhood, the Chern--Simons
functional with the primary function plus the symplectic action of the section,
involves neither the metric nor $J$, so its critical set does not depend on the parameter.
The Coulomb--Neumann slices and the charts in which it is expanded depend continuously on
the parameter, so the remainder bounds of the expansion hold with uniform constants; the
gradient is taken in the $L^2$ inner product of the parameter, and these inner products are
uniformly equivalent. With the uniform gap this gives
$|\mathcal F-\mathcal F_{\rm crit}|\le C\|\nabla\mathcal F\|^2$ with $C$ independent of the parameter. Third, the product metric and the form of $j$ on the
end give $\partial_\theta\mathcal F=-\|\nabla\mathcal F\|^2$. Hence the kinetic energy beyond
$\theta$ decays at a rate $c$ independent of the parameter. With the weight $\delta$ of
Lemma~\ref{lem:family-stationary} chosen below $c$, this is its hypothesis.
\end{proof}

The DFL pieces produced by the degenerations below are such domains, so the
maps they count are the mixed maps up to chain homotopy.

\subsection{Critical points on the twisted inputs and the gauge triangle}

\begin{lemma}\label{lem:twisted-critical}
Let $ij\in\{01,12,02\}$, with the $\rho_\zeta$-marking on $C_0$ when
$i=0$, and let $\alpha$ be a critical point of the primary equation
$*F_B+\nabla\mathfrak h_{ij}(B)=0$ on $Y_{ij}^\#$. Then $\alpha$ is flat on a collar of
a copy of $\Sigma_h$ disjoint from the primary dependence regions,
$w_2$ is nonzero on $\Sigma_h$, $\alpha$ is irreducible, an automorphism
fixing $\alpha$ whose lifting class vanishes on $\Sigma_h$ is the identity,
the nonzero element of the lifting group acts freely on the critical set, and the critical
set corresponds bijectively to $P_i\cap P_j$ and is nondegenerate.
Moreover end perturbations as in Proposition~\ref{prop:input-functions}
can be chosen so that the critical set of $\mathrm{CS}_{\mathfrak h_{ij}}+f_{ij}$ equals that
of $\mathrm{CS}_{\mathfrak h_{ij}}$.
\end{lemma}

\begin{proof}
The pair $(C_0,C_j)$ with the $\rho_\zeta$-composed boundary marking is a
pair of the type of \cite[Section~2.3]{DFL}, and its Lagrangian is
$\tau_\zeta L(C_0)$. Flatness near $\Sigma_h$, irreducibility and the
correspondence with $P_i\cap P_j$ with nondegeneracy are
\cite[Lemmas~2.43, 2.45, 2.47]{DFL}, whose proofs use only flatness and
irreducibility on the overlap and Mayer--Vietoris. The restriction of
$E_{ij}^\#$ to a push-off of $\Sigma_h$ into $C_j$ is the odd surface bundle,
unchanged by $\rho_\zeta$. Freeness follows from
Theorem~\ref{thm:regularity}(i). For the last statement, the twisted-marked pairs are not presentations, so
Proposition~\ref{prop:input-functions} is not applied as stated; its proof is
repeated for each pair $ij$, as follows. If a sequence
$f_n$ with vanishing $C^0$ norms of gradients and of their derivatives on
fixed balls had critical points $B_n$ not gauge equivalent to critical
points of $\mathrm{CS}_{\mathfrak h}$, then Uhlenbeck compactness on the closed three-manifold
and elliptic regularity would give a limit critical point $\alpha$ of
$\mathrm{CS}_{\mathfrak h}$, near which the inverse function theorem for the gauge-fixed
critical equation, with invertible extended Hessian, allows only the zero
$\alpha$ itself, since $f_n$ vanishes near $\alpha$.
\end{proof}

\begin{proposition}\label{prop:gauge-triangle-regular}
On $X^\#$ with the bundle $E^\#$ and the three closed ends of
Lemma~\ref{lem:twisted-critical}, arbitrarily small $\Phi$-invariant
secondary perturbations give: no solutions of negative index, a finite
regular index-zero space, and regular one-dimensional spaces whose
boundary strata are an index-zero solution with one index-one trajectory
at one end. The count $\mu_G^E$ is a chain map.
\end{proposition}

\begin{proof}
The ends are nondegenerate and irreducible by
Lemma~\ref{lem:twisted-critical}, so the gauge-fixed operator is Fredholm, the gauge-fixing component having positive gaps. The index-energy law
$\ind(A)=8\kappa(A)+c(\alpha_0,\alpha_1,\alpha_2)$ holds on the pure gauge triangle:
for two connections $A,A'$ with the same limits, glue to both the same
connections on compact fillings of the three ends (every $\SO(3)$ bundle over
a closed three-manifold extends over some compact four-manifold); since the
ends are acyclic, additivity of the index \cite[Proposition~3.9]{Donaldson-book}
reduces $\ind(A)-\ind(A')$ to the difference on the closed manifold, where the
Atiyah--Singer formula gives $8(\kappa(A)-\kappa(A'))$. Compactness is
Uhlenbeck compactness with point concentrations of energy $8\pi^2$ in the
normalization of Section~\ref{subsec:energy-convention} and index loss
eight; a small three-sphere about a concentration point carries no
nonzero lifting class, so the limiting bundle is identified with $E^\#$ without any
unitary lift. Surjectivity in the invariant class is
Theorem~\ref{thm:regularity}(ii), with the loops supplied by an incoming
end: on $Y_{12}^\#$ the bundle admits a $\U(2)$ lift because $H^3$ of a
closed oriented three-manifold is torsion free, the critical limit is
flat on a collar of $\Sigma_h$ with odd bundle, and loops on that surface,
which have zero lifting class in $\Phi$, are transported into the interior along a
path; holonomies are conjugated by parallel transport, which preserves a
basis. Finite selection after compactness at each index and freeness of
$\Phi$ complete the proof.
\end{proof}

\subsection{The triangle identity}

\begin{theorem}\label{thm:triangle}
With the mixed maps $N_{01}$, $N_{12}$, $N_{02}$ of
Theorem~\ref{thm:mixed-iso} for the three pairs of~\eqref{eq:triple} and
the symplectic triangle map $\mu_S^E$, which counts holomorphic triangles only,
\[
N_{02}\circ\mu_G^E\simeq\mu_S^E\circ(N_{01}\otimes N_{12})
\quad\text{over }\FF .
\]
\end{theorem}

The proof compares two degenerations of the mixed triangle of
Definition~\ref{def:mixed-triangle}, with the regions $G$, $S$, the arc $\gamma$ and the
rectangle $B$ fixed there.

\begin{lemma}\label{lem:cut-crosses}
Every proper arc in $G$ from the interior of $e_0$ to the interior of $e_2$
meets $B$. Consequently the closed cut $Y_{02}^\#$ passes through the
collapsing rectangle, and a product collar of the cut exists only when the
compression scale over $I$ is of order one.
\end{lemma}

\begin{proof}
$B$ meets the frontier of $G$ only in $I$ and $J$, and $G\setminus B$ has two
components, one containing each input vertex. The arcs $e_0\cap G$ and
$e_2\cap G$ lie in different components, so an arc joining them meets $B$.
Over $B$ the metric is $dt^2+\lambda_\eps^2(d\rho^2+g_\Sigma)$, whose
surface factor is scaled by $\eps$, so no neighbourhood of the cut is a
product with a fixed closed metric until $\lambda$ is raised to one.
\end{proof}

\begin{lemma}\label{lem:family-path}
There is a smooth path $F_\tau$, $\tau\in[0,1]$, of data on one marked
domain with the bundle $E$, with the following properties.
\begin{enumerate}[label=(\roman*)]
\item Near $\tau=0$ it is the collapse datum of Theorem~\ref{thm:collapse}
  at small scale $\eps$, with the ordinary necks of Theorem~\ref{thm:neck},
  satisfying Hypothesis~\ref{hyp:collapse}. Near $\tau=1$ it is the datum
  with a long product collar of the closed cut $Y_{02}^\#$ inserted.
\item It satisfies the hypotheses of Theorem~\ref{thm:family}:
  \ref{H:lattice} and \ref{H:support}--\ref{H:framing} on every fibre,
  and \ref{H:regular} at the endpoints, as verified in the proof of
  Theorem~\ref{thm:two-handle}. The inserted necks of \ref{H:support}
  are the ordinary necks near $\tau=0$ and the cut collar near $\tau=1$.
\item The scales stay in a compact positive range: the compression
  scale over $I$ is at least $\eps$, and the cut collar has bounded length.
\item The holonomy loop families are fixed. Only their forms vary,
  with support away from the seam collars and the ends; the loops
  themselves may pass through these regions. The derivative
  $\partial_\tau F_\tau$ has compact support disjoint from all ends.
\item Every $F_\tau$ is invariant under $\Phi$.
\end{enumerate}
\end{lemma}

\begin{proof}
The path has five stages, each inside a convex set of such data:
(A) raise the scale over $I$ to one by the convex combination
$\mu_s^2=(1-s)\lambda_\eps^2+s$ of metrics of the form $dt^2+\mu^2g_C$ and of
the matching section metrics; (B) interpolate linearly to a metric which is
a product $a^2du^2+g_{02}$ on a collar of the cut, the two metrics agreeing
with every prescription on every overlap because on the compression pieces
the collar coordinate is the edge coordinate and $g_{02}$ restricts to the
compression metrics, and on $B$ it is the stage (A) metric; (C) lengthen
the collar; (D) switch on the end perturbation of $Y_{02}^\#$ on the
collar, by a cutoff equal to one on the middle of the collar and monotone on its two end
bands, which become the cutoff bands of the adjacent ends; during this stage the term is a
compactly supported nonautonomous perturbation, and it is fixed from then on; (E) interpolate the remaining metrics, section structures and
forms to the data with the collar of the cut inserted. Compatible almost complex structures are
joined by $J_k=(-A_k^2)^{-1/2}A_k$ for the convex combination $k$ of the
associated metrics, where $\omega(u,v)=k(A_ku,v)$; this equals the Hodge
structure wherever both ends do. Forms are replaced by their self-dual
projections in the varying metric. The bundle over the collar is a product
because all transitions are constant along the edges, so stretching the
collar lifts to $E$ with no clutching. Invariance holds because every class
of $\Phi$ vanishes on the surface sheets and hence, by injectivity of
$H^1(C_i)\to H^1(\Sigma_h)$, on the compression pieces, and every loop has
zero lifting class in $\Phi$.
\end{proof}

\begin{proof}[Proof of Theorem~\ref{thm:triangle}]
\emph{Endpoint $1$: the outgoing cut.} Stretch along the closed
three-manifold $Y_{02}^\#$ separating the gauge triangle from the DFL domain
for the $\rho_\zeta$-marked pair $(P_0,P_2)$. We verify the hypotheses of
Theorem~\ref{thm:neck} for these two pieces and this connected cut.
\ref{H:lattice}: on the gauge piece the index-energy law is part of
Proposition~\ref{prop:gauge-triangle-regular}; on the DFL piece it is
\cite[Proposition~3.20, Lemma~3.15]{DFL} for the marked pair, and on the
joined problem it follows by additivity across the cut, the constant pair at
each output generator having index zero, so that the constant of the joined
problem is the sum of the two constants. \ref{H:regular}: $\mu_G^E$ is
regular by Proposition~\ref{prop:gauge-triangle-regular}, and $N_{02}$ is
regular in the invariant class by Theorem~\ref{thm:regularity}, with the
loops supplied from its incoming end $Y_{02}^\#$, whose critical limits are
flat on a collar of an odd copy of $\Sigma_h$ by
Lemma~\ref{lem:twisted-critical}; the section part of a cokernel element is
killed by harmonic traces on the matching arc, and the two
compression ends are conjugate to the untwisted ones by $\tau_\zeta$. The
cross-section at the cut is invertible because the critical points of
$Y_{02}^\#$ are nondegenerate and irreducible
(Lemma~\ref{lem:twisted-critical}), so the extended Hessian, including the gauge-fixing component, has a spectral gap. \ref{H:support}--\ref{H:framing} hold by
construction. Theorem~\ref{thm:neck} gives the count
$N_{02}\circ\mu_G^E$ for long cut collars, the DFL piece being identified with
the mixed map by Lemma~\ref{lem:vertex-identification}, and a finite change of collar
length is a compact family (Theorem~\ref{thm:family}). The lifting groups glue by the fibre
product: restriction $\Phi\to\FF^2$ to the two incoming ends is an
isomorphism, the output group of the DFL piece restricts isomorphically to
the lifting group of the cut, and the connecting $H^0$ term at the connected cut
vanishes.

\emph{Endpoint $0$: the interval collapse.} Collapsing $B$ by
Theorem~\ref{thm:collapse}, the limit domain is a disk with a physical
$P_1$ boundary interval and two new mixed ends; cutting off two DFL disks
along ordinary strip ends leaves a triangle domain with boundary
$(P_0,P_1,P_2)$. The triangle piece is regular for a generic domain-dependent
almost complex structure supported away from the neighbourhood $N_I$ of the
collapsed arc, where Hypothesis~\ref{hyp:collapse}\ref{hyp:X1} fixes the complex
structure $*_\Sigma$: there are no constant triangles because the triple intersection is
empty, so every solution is nonconstant on an open set away from $N_I$ by unique
continuation, and the standard transversality argument applies there.
Choose the ordinary neck lengths first: by
Theorem~\ref{thm:neck} the zero-scale domain is then regular in index at
most zero with index-zero count $\mu_S^E(N_{01}\otimes N_{12})$, the two DFL
pieces being identified with the mixed maps by
Lemma~\ref{lem:vertex-identification}; this is
Hypothesis~\ref{hyp:collapse}\ref{hyp:X4}, and \ref{hyp:X5} holds by
Lemma~\ref{lem:zero-scale-law}. Then choose $\eps$ small; by
Theorem~\ref{thm:collapse}(iii) the count at scale $\eps$ equals the
zero-scale count. Hypotheses \ref{hyp:X1}--\ref{hyp:X3} hold because the
middle primary perturbation vanishes, the end perturbations are the fixed autonomous
perturbations of the original gauge ends, with outputs and loops in those ends and cutoff bands
placed beyond the compact set $N_I$, the secondary perturbations have outputs and loop regions in a compact set
disjoint from $N_I$, and the original ends are regular; at the original
mixed ends this holds because a nonconstant mixed-end trajectory has
$\kappa\ge\frac12$, the action being constant on the connected critical family,
and hence index at least four. The zero-scale gauge region consists of two DFL
pieces (Lemma~\ref{lem:cut-crosses}), with lifting group $\FF\times\FF$; restriction
$\Phi\to\FF^2$, $t\mapsto(s_{01},s_{12})$, is an isomorphism, so the
marking-preserving bijection of Theorem~\ref{thm:collapse}(iii) compares
$\Phi$-framed counts with no factor two.

\emph{The family.} Apply Theorem~\ref{thm:family} to the path of
Lemma~\ref{lem:family-path}. Its hypotheses: the index-energy law on every
fibre with the constant of the twisted endpoint, since a reference
configuration on one fibre is a reference on every fibre and the
topological energy involves no metric (so the index of a fixed reference
is constant in $\tau$ by continuity of the Fredholm family, and the law at
$\tau=0$ is Theorem~\ref{thm:collapse}(i)); an energy bound independent of
$\tau$, the terms $2\int\chi'f$ being bounded by twice the $C^0$ norms of the
end perturbations times the total variation of their cutoffs, which is one for the original
ends and at most two for the switch of stage (D), whose cutoff rises and falls once; family compactness, with weak limits identified
with the prescribed bundle $E$ by modifying bundle isomorphisms only on
small balls around the concentration points joined by arcs to an incoming
slice, which introduces no lifting class; family regularity in the invariant
class by Theorem~\ref{thm:regularity}, with holonomy terms built from local
unitary lifts, which exist on neighbourhoods of the loops and are unique up
to tensoring with a flat real line, so that each term is defined up to sign
and the span is unaffected; and no internal neck breaks, since the scales
stay in a compact positive range. The breaks of the parametrized
one-dimensional space therefore occur only at the original ends: the two
closed inputs and the ordinary output. No break occurs at the two original mixed
ends at $\gamma\cap e_0$ and $\gamma\cap e_2$: a nonconstant mixed-end trajectory has
index at least four, as noted at Endpoint~$0$. At an incoming break the index $-1$
piece may carry the other even lift of that input; a unique class of $\Phi$
aligns it, and by invariance of the whole family the aligned piece is
already counted. The boundary count gives the chain homotopy.
\end{proof}

\subsection{The two-handle square}

Let $p_i^{\phi}$ be the admissible presentations of $Z_0=Y_{01}$ and
$Z_1=Y_{02}$ given by the splittings $(H_0,H_1)$ and $(H_0,H_2)$ with the
global frames $\phi_i^{\mathrm{gl}}$ of Section~\ref{subsec:two-handle-bundles},
whose primary data are those transported by the frames. These presentations are admissible:
$\operatorname{span}F$ contains each trace function and the constants, so both transports act on
it linearly, and the shared perturbation on $C_0$ is taken small enough that both transported
data lie in the open neighbourhood $\mathcal U(F)$.
Let $P^G_{\phi_i}$ and $P^S_{\phi_i}$ be the transport of complexes from the
filling frame to the global frame: on the gauge side the pushforward by the
bundle isomorphism, and on the symplectic side the pushforward by the
symplectomorphism $\tau_{a_1}$ (respectively $\tau_{a_2'}$) of $\Msf_h$. By
naturality of every construction under bundle isomorphisms,
$N_{p_i^\phi}P^G_{\phi_i}=P^S_{\phi_i}N_{0,i+1}$ on chains, and by the definition
of the morphisms of \cite{KM-unknot}, which include the identifications at the
ends, $\Ish(W_2,E;\phi^{\mathrm{gl}})=P^G_{\phi_1}F_G(W_2,E)(P^G_{\phi_0})^{-1}$.

\begin{theorem}\label{thm:two-handle}
Define the symplectic map
\[
S_{W_2,E}=P^S_{\phi_1}\,\mu_S^E(\,\cdot\otimes\theta_S)\,(P^S_{\phi_0})^{-1}
\colon\HF(Z_0,p_0^\phi)\to\HF(Z_1,p_1^\phi).
\]
Then $H(N_{p_1^\phi})\circ\Ish(W_2,E;\phi_0^{\mathrm{gl}},\phi_1^{\mathrm{gl}})
=S_{W_2,E}\circ H(N_{p_0^\phi})$.
\end{theorem}

\begin{proof}
Insert $x\otimes\theta_G$ into the chain homotopy of
Theorem~\ref{thm:triangle}. Since $\theta_G$ is a cycle,
$N_{02}\mu_G^E(x\otimes\theta_G)$ is homologous to
$\mu_S^E(N_{01}x\otimes N_{12}\theta_G)=\mu_S^E(N_{01}x\otimes\theta_S)$ by
Proposition~\ref{prop:filling-classes}. By
Proposition~\ref{prop:gauge-triangle}, $N_{02}F_G(W_2,E)=\mu_S^E(\,\cdot
\otimes\theta_S)N_{01}$ on homology. Conjugating by the transports gives the
claim.
\end{proof}

The map $S_{W_2,E}$ is defined by two pushforwards by symplectomorphisms and
a triangle count; it does not involve gauge theory.

\begin{remark}
For $E$ nontrivial the gauge map can vanish for grading reasons: for
the $(-1)$-framed unknot with the nontrivial bundle, $\chi=1$,
$\sigma=-1$, $P(E)=1$, and the map $\Ish(S^3)\to\Ish(S^3)$ has degree
two by Remark~\ref{rem:grading}, hence is zero
\cite[Section~7.6]{Scaduto}. The theorem then says that the symplectic
triangle map $S_{W_2,E}$ vanishes on homology.
\end{remark}

\section{Excision along a torus for the mixed equation}\label{sec:torus}

The one-handle comparison of Section~\ref{sec:one-three} passes through a
mixed domain whose gauge input is a connected sum along a torus. On the
gauge side such a manifold is split by Floer's excision theorem
\cite{BD, KM-sutures}, in the form of \cite[Theorem~5.6]{KM-unknot}. This section proves the corresponding mixed statement,
Theorem~\ref{thm:excision}.

\subsection{The folded splittings}\label{subsec:folded}

Let $M$ have presentation data with torus compression bodies
$U_0,U_1$ of genus $h\ge2$, and let $C=C_{h,g}$, $g=h+k$, be a compression
body from $\Sigma_h$ to $\Sigma_g$ adding $k$ handles. Put $V_i=C\circ U_i$,
the enlarged torus compression bodies. Then $V_0\cup(-V_1)$ is the torus-summed object $(M^+)^\#$ of $M^+=M\#^k(S^1\times S^2)$: every added meridian disk
occurs on both sides, and the two disks form a nonseparating sphere. Let
$L_C\subset\Msf_h^-\times\Msf_g$ be the Lagrangian correspondence of $C$. For
a Lagrangian $P\subset\Msf_h$ transverse to the lower projection, the
composition $P\circ L_C$ is embedded, because the lower projection of $L_C$ is
a submersion and its upper projection an embedding: $L_C$ is the flat moduli
space of $C$, which fibres over the moduli space of $\Sigma_h$ with fibre
$\SU(2)^k$, and restriction to $\Sigma_g$ is injective since $\pi_1(\Sigma_g)
\to\pi_1(C)$ is onto \cite[Proposition~2.21]{DFL}.

\emph{Height Lagrangians.} For $t>0$ small, perturb the primary data of $V_i$ by
$t$ times a height function whose restriction to the clean intersection
$L(V_i)\cap(L(U_i)\circ L_C)\cong\SU(2)^{g-1}$ has the $2^{g-1}$ nondegenerate
critical points of Section~\ref{sec:two-handles}; call the resulting upper
Lagrangians $R_i\subset\Msf_g$ the \emph{height Lagrangians} at height $t$. The top generator $q_i^*$ of each auxiliary pair
below is defined intrinsically, by the grading $m$ of
Section~\ref{subsec:continuations}; at unperturbed lower data and uncorrected height data it
is the maximum or the minimum of the height function, according to whether the perturbed
Lagrangian is the first or the second entry of the pair, and in general it is the continuation
of that critical point (Section~\ref{subsec:one-three-squares}; Lemma~\ref{lem:aux-index-folded}(i)
for $k\ge1$).

\emph{Lower Lagrangians.} Then choose small primary data $\widetilde P_i$ of $U_i$, in a
neighbourhood of $L(U_i)$ that is small relative to $t$: for $k=0$ lower data in the set
$\mathcal U_F(t)$ of Definition~\ref{def:admissible-nbhd}, whose Lagrangians lie in
$\mathcal N(t)$, and in general the neighbourhoods
$\mathcal N_{C,R}$ of Lemma~\ref{lem:aux-index-folded}(i). In
$\mathbb X=\Msf_h^-\times\Msf_g$ put
\[
A_0=\widetilde P_0^-\times R_0,\qquad A_1=\widetilde P_1^-\times R_1 ,
\]
consider the triple $(A_0,L_C,A_1)$, and require $A_0\pitchfork L_C$, $L_C\pitchfork A_1$, $A_0\pitchfork A_1$ and
$A_0\cap L_C\cap A_1=\varnothing$. These are generic conditions: the triple
intersection has expected dimension $3d-2\cdot2d=-d<0$, where
$d=(3h-3)+(3g-3)$ is the dimension of each of the three Lagrangians in the
manifold $\mathbb X$ of dimension $2d$, and the
universal restriction map for primary perturbations is a submersion
\cite[Section~2.3]{DFL}. (With unperturbed data the triple intersection
would contain the fibres of $L_C$ over $L(U_0)\cap L(U_1)$.) The pairs
$(A_0,L_C)$ and $(L_C,A_1)$ have intersection generators identified
with those of $(R_0,\widetilde P_0\circ L_C)$ and
$(\widetilde P_1\circ L_C,R_1)$, hence, by continuation, with the critical sets of the
height functions. Their strip problems remain the product-target
problems in $\mathbb X$; exact unfolding is used only for the diagonal case
$k=0$.

The gauge pieces are $D_i=(-U_i)\cup_{T^2}V_i$, with $\partial D_i=(-\Sigma_h)
\sqcup\Sigma_g$, and $C$. Flat connections on $D_i$ restrict bijectively to
pairs of flat connections on $U_i$ and $V_i$: across the odd torus the gluing
is unique, since the flat connection on the odd torus has central
determinant-one stabilizer and $H^*(T^2;\ad)=0$; so $L(D_i)=L(U_i)^-\times L(V_i)$,
perturbed to $A_i$. The closed cut of the pair $(A_0,A_1)$ is
\[
Z=D_0\cup(-D_1)\cong(-M)\#M^+\#T^3,
\]
with one torus factor. For a split almost complex structure $J_-\times J_+$ on
$\mathbb X$, folding gives an isomorphism of complexes $\mathcal
K\colon\CF_{\mathbb X}(A_0,A_1)\to\CF(\widetilde P_0,\widetilde P_1)^\vee\otimes\CF(R_0,R_1)$,
and a cycle $\sum x_i^\vee\otimes y_i$ on the right is the chain map
$v\mapsto\sum x_i^\vee(v)y_i$.

The metrics are chosen as follows: $U_0,U_1$ share their metric near
$\Sigma_h$, $V_0,V_1$ near $\Sigma_g$, and all torus collars carry the flat
metric induced from the fixed metric on $T^3$, with the fixed odd bundle.

\subsection{The torus region and the domains}\label{subsec:torus-domains}

Let $T=\RR^2/\Lambda$ be the flat torus, $P_T$ the odd bundle and $\alpha$ its
flat connection, unique up to gauge. The adjoint bundle splits into three flat
real line bundles with nontrivial characters, so $H^*(T;\ad\alpha)=0$, the
twisted Laplacian satisfies $\Delta_\alpha\ge\mu_T^2$ in every degree with
$\mu_T=\pi\min\{|\xi|:\xi\in\Lambda^*\setminus2\Lambda^*\}$, and the
determinant-one stabilizer of $\alpha$ is the centre.

Consider the special quintuple $X_D$ of \cite[Section~3.1]{DFL} for a pair of
torus compression bodies $(Y,Y')$ with $\partial Y=\Sigma'\sqcup T$. Its gauge
part contains $\mathbb U_-\times T$, where $\mathbb U_-$ is a planar region with one gauge
input arm, two mixed-end arms $\theta\to\pm\infty$ and the matching line
$\mathbb U_\partial=\{0\}\times\RR_\theta$; the metric on $\mathbb U_-$ is flat near the matching
line. Hence $X_D$ contains the flat product region
\[
N_X=T\times[-1,0]_s\times\RR_\theta ,
\]
whose face $\{s=0\}$ is the torus component of the matching surface, and for
$|\theta|\ge2$ this region lies in the mixed ends. Since the moduli space of the
odd torus is a point, the matching condition there says only that the
restriction of the connection to each fibre $T\times\{(0,\theta)\}$ is flat; we
call it the \emph{flat-fibre condition}. Its linearization is $d_\alpha\zeta_T=0$
on the face, $\zeta_T$ the component along the torus; together with the Neumann
gauge condition $\zeta(\partial_s)=0$ this gives the boundary conditions used below,
the longitudinal component $\zeta(\partial_\theta)$ being free.

\begin{definition}\label{def:torus-domains}
\begin{enumerate}[label=(\alph*)]
\item For $0\le \mathsf R<\infty$, the \emph{collared domain} is
  $W_\mathsf R=X_D\cup(T\times[0,\mathsf R]_r\times\RR_\theta)$, glued along the torus face,
  with the flat product metric and no perturbation on the collar and the
  flat-fibre condition at $r=\mathsf R$. Its mixed-end slices are
  $Y_\mathsf R=Y\cup([0,\mathsf R]\times T)$ and $Y'_\mathsf R$. It is diffeomorphic to $X_D$ by a map
  $\psi_\mathsf R$ that stretches $s\in[-3/4,-1/4]$ and is the identity outside $N_X$,
  hence on the supports of all perturbations, on the gauge and strip ends,
  near the other matching surface and on the section.
\item The \emph{completed domain} is $W_\infty=X_D\cup(T\times[0,\infty)_r\times
  \RR_\theta)$, with no torus face and mixed slices $Y_\infty$, $Y'_\infty$.
\item For the two sides $X_{D,-}$ (the reflected presentation of the lower end $M$
  of the one-handle cobordism) and
  $X_{D,+}$ (the enlarged presentation of $M^+$) and $\lambda>0$, the
  \emph{joined domain} $W^J_\lambda$ is $W_{-,\lambda}\cup W_{+,\lambda}$ glued along
  their torus faces by the marked bundle isometry reversing the torus
  orientation together with the normal coordinate. It has the bridge
  $T\times[-\lambda,\lambda]\times\RR_\theta$, two gauge inputs $(-M)^\#$ and
  $(M^+)^\#$, a section $\mathbb U_+\to\Msf_h^-\times\Msf_g$ with the split structure
  $J_-\times J_+$ on all of $\mathbb U_+$, strip output at $A_0\cap A_1$, and mixed-end
  slices $D_{0,\lambda}=Y_-\cup_{[-\lambda,\lambda]\times T}Y_+$ and
  $D_{1,\lambda}=Y'_-\cup_{[-\lambda,\lambda]\times T}Y'_+$: on each arm the slices of the
  two sides are joined through the bridge.
\end{enumerate}
Weights are functions $\tau\ge0$ with $|d\tau|\le1$ and $|\nabla d\tau|\le1$, equal to
the end coordinate on the ends and, on the torus region, to a function
$\tau_0(\theta)$ of $\theta$ alone. The domains $W_{\pm,\mathsf R}$ and $W_{\pm,\infty}$ are
called \emph{separate}, as opposed to the \emph{joined} domains $W^J_\lambda$; their torus
collars, and the mixed maps defined on them, are called separate as well.
\end{definition}

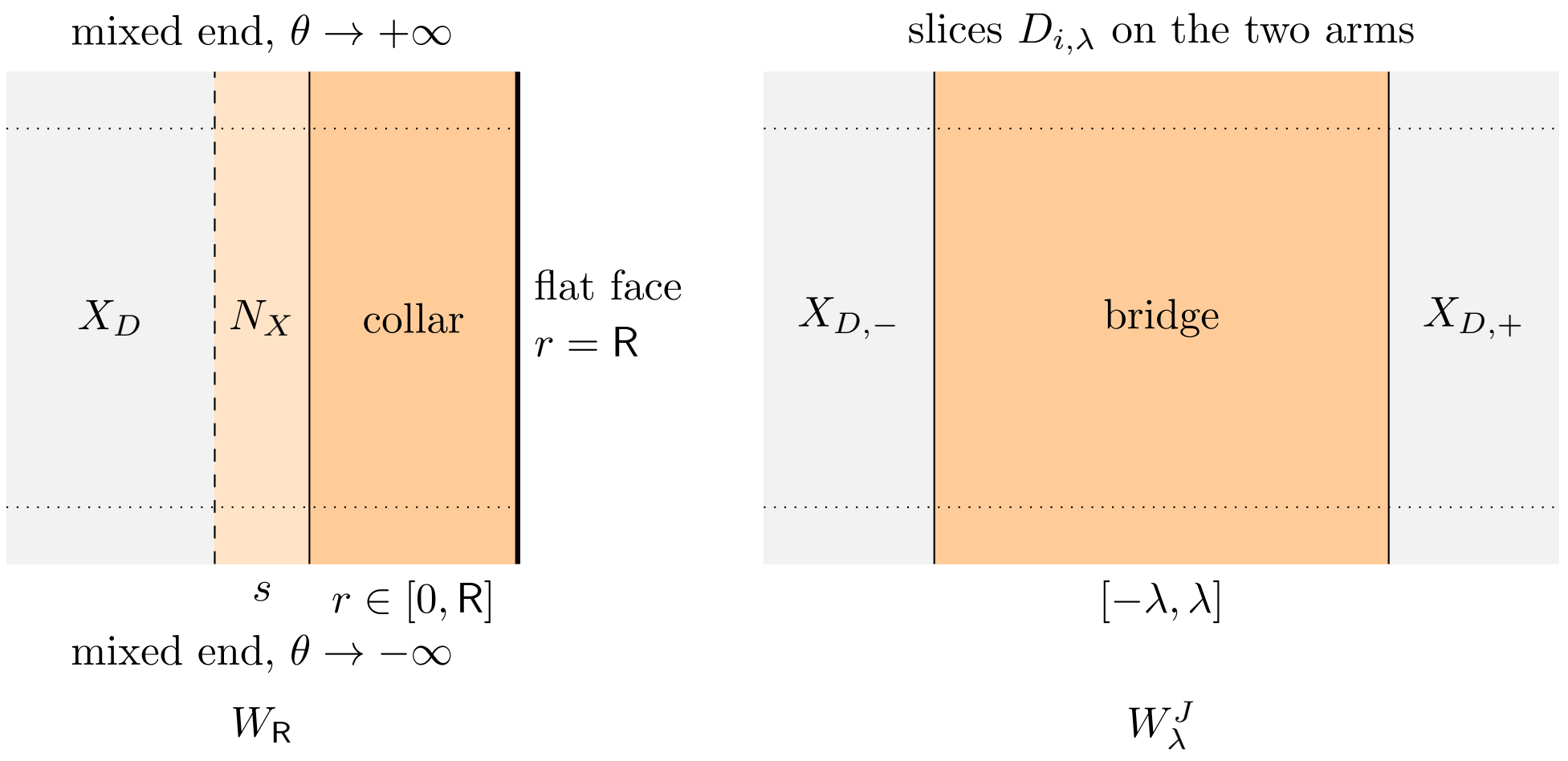
\begin{figure}[ht]
\centering
\begin{tikzpicture}[scale=0.78,font=\small]
\begin{scope}
\fill[gray!10] (-3.2,-2.6) rectangle (-1,2.6);
\fill[orange!22] (-1,-2.6) rectangle (0,2.6);
\fill[orange!40] (0,-2.6) rectangle (2.2,2.6);
\draw[dashed] (-1,-2.6) -- (-1,2.6);
\draw (0,-2.6) -- (0,2.6);
\draw[very thick] (2.2,-2.6) -- (2.2,2.6);
\draw[dotted] (-3.2,2) -- (2.2,2); \draw[dotted] (-3.2,-2) -- (2.2,-2);
\node at (-2.1,0) {$X_D$};
\node at (-0.5,0) {$N_X$};
\node at (1.1,0) {collar};
\node[right] at (2.2,0.35) {flat face};
\node[right] at (2.2,-0.25) {$r=\mathsf R$};
\node[below] at (-0.5,-2.62) {$s$};
\node[below] at (1.1,-2.6) {$r\in[0,\mathsf R]$};
\node at (-0.5,3.0) {mixed end, $\theta\to+\infty$};
\node at (-0.5,-3.55) {mixed end, $\theta\to-\infty$};
\node at (-0.5,-4.3) {$W_\mathsf R$};
\end{scope}
\begin{scope}[xshift=9cm]
\fill[gray!10] (-4.2,-2.6) rectangle (-2.4,2.6);
\fill[gray!10] (2.4,-2.6) rectangle (4.2,2.6);
\fill[orange!40] (-2.4,-2.6) rectangle (2.4,2.6);
\draw (-2.4,-2.6) -- (-2.4,2.6); \draw (2.4,-2.6) -- (2.4,2.6);
\draw[dotted] (-4.2,2) -- (4.2,2); \draw[dotted] (-4.2,-2) -- (4.2,-2);
\node at (-3.3,0) {$X_{D,-}$};
\node at (3.3,0) {$X_{D,+}$};
\node at (0,0) {bridge};
\node[below] at (0,-2.6) {$[-\lambda,\lambda]$};
\node at (0,3.0) {slices $D_{i,\lambda}$ on the two arms};
\node at (0,-4.3) {$W^J_\lambda$};
\end{scope}
\end{tikzpicture}
\caption{The torus region, drawn in the base coordinates (normal coordinate
horizontally, the end coordinate $\theta$ vertically); every point carries a torus
$T$ with the flat odd connection. Left: the collared domain $W_\mathsf R$ of
Definition~\ref{def:torus-domains}(a), with the flat product region $N_X$ of $X_D$, the
collar and the flat face at $r=\mathsf R$; for $|\theta|\ge2$ (dotted lines) the region lies in
the mixed ends. The completed domain $W_\infty$ has $\mathsf R=\infty$. Right: the joined
domain $W^J_\lambda$, in which the collars of the two sides are glued along their flat
faces into a bridge of width $2\lambda$. No closed hypersurface separates the two
sides, since the torus region runs into the mixed ends.}
\label{fig:torus-domains}
\end{figure}

The critical manifold of the joined mixed ends is $P_i\times Q_i$ for the two
separate critical manifolds, because the two flat connections glue uniquely
across the odd torus. The lifting group of $W^J_\lambda$ is $\Phi_J\cong\FF^2$,
restricting isomorphically to $\Phi_-\oplus\Phi_+$, with lifting class $s_Z=s_-+s_+$
at the cut $Z$; its elements are lifting classes that vanish on every piece (they come from the graph with a vertex for each surface sheet and each compression piece and an edge for each contact region), since
restriction of $H^1$ from each compression body to its full boundary is
injective.

Counting index-zero solutions defines: the separate maps $N_\pm^{(\mathsf R)}$ on
$W_\mathsf R$, which for $\mathsf R=0$ are $N_\pm$; the completed maps $N_{\pm,\infty}$ on
$W_\infty$; and the joined maps
$N^J_\lambda\colon C_G((-M)^\#)\otimes C_G((M^+)^\#)\to\CF_{\mathbb X}(A_0,A_1)$.

\subsection{Statement}

Let $N_-\colon C_G((-M)^\#)\to\CF(\widetilde P_0,\widetilde P_1)^\vee$ be the mixed map
of the reflected presentation of the lower end $M$ (Section~\ref{subsec:reflection}),
$N_+\colon C_G((M^+)^\#)\to\CF(R_0,R_1)$ that of the enlarged presentation with the height Lagrangians, and $N_Z\colon C_G(Z)\to\CF_{\mathbb X}(A_0,A_1)$ the mixed map of the
pair $(A_0,A_1)$ (Remark~\ref{rem:clean-pairs}) with the data at width $\lambda_0$ fixed in Section~\ref{subsec:closed-cut}. Let $\mathsf E\colon C_G(Z)\to
C_G((-M)^\#)\otimes C_G((M^+)^\#)$ be the excision map of
\cite[Theorem~5.6]{KM-unknot}, with the end-perturbation data $f_Z$ of the same data on $Z$, and $\overline{\mathsf E}$ the reverse excision map. For other data on $Z$ both
sides change by the same gauge continuation map.

\begin{theorem}[mixed torus excision]\label{thm:excision}
On homology over $\FF$,
\[
[\mathcal K][N_Z]=([N_-]\otimes[N_+])\,[\mathsf E].
\]
\end{theorem}

\emph{The gauge half.} $Z$ contains two parallel tori separating the two
summands; cutting along them and regluing crosswise gives
$Z'=(-M)^\#\sqcup(M^+)^\#$. The bundle is odd on these tori and every component
of $Z$ and $Z'$ meets one. The excision cobordism $\mathcal E\colon Z\to Z'$ is
$T^2\times P$, $P$ a pair of pants with input $c$ and outputs $c_\pm$, with the
products of the punctured summands inserted along two arcs from $c$ to $c_\pm$;
its lifting group is $\Phi_{\mathcal E}=\{(s,s_-,s_+):s=s_-+s_+\}$, dual to the
four hypersurfaces of the construction of \cite[Section~5.3]{KM-unknot}. By
\cite[Theorem~5.6]{KM-unknot}, $[\mathsf E]$ is an isomorphism with inverse
$[\overline{\mathsf E}]$. The theorem is stated over $\ZZ$; the maps $\mathsf E$, $\overline{\mathsf E}$ and
the chain homotopies in its proof are defined over $\ZZ$, so the conclusion
holds after reduction to $\FF$, and the
Künneth identification for $Z'$ is over the field $\FF$.

\emph{Plan of the proof.} Only the mixed statement is new. It is not an
instance of neck-stretching along a closed hypersurface: the torus region runs
into the mixed ends, so no closed three-manifold separates the two sides. The
proof has three steps:
\begin{enumerate}[label=(\alph*)]
\item (closed cut) $\mathcal KN^J_{\lambda_0}\simeq\mathcal KN_Z\overline{\mathsf E}$ at a fixed
  width $\lambda_0$ (Proposition~\ref{prop:closed-cut}, including the data change
  of Lemma~\ref{lem:closed-cut-data});
\item (width) $N^J_b\simeq N^J_a$ for $\lambda_*\le a\le b<\infty$
  (Proposition~\ref{prop:width});
\item (completion) for large $\lambda$ and large $\mathsf R$,
  $\mathcal KN^J_\lambda=N_{-,\infty}\otimes N_{+,\infty}$ and $N_\pm^{(\mathsf R)}=N_{\pm,\infty}$
  as matrices, and $N^{(\mathsf R)}_\pm\simeq N_\pm$
  (Theorem~\ref{thm:torus-completion}, Proposition~\ref{prop:collar-length}).
\end{enumerate}
Then $[\mathcal KN_Z][\overline{\mathsf E}]=[N_-\otimes N_+]$, and since $[\overline{\mathsf E}]=[\mathsf E]^{-1}$
the theorem follows.

\subsection{The linearized operator at the odd flat torus}

The constant $\mu_T>0$ of Section~\ref{subsec:torus-domains}, the square root of the lowest eigenvalue of the twisted Laplacian of $\alpha$, is the spectral gap that makes the following operators invertible.

\begin{proposition}[the operator at the odd flat torus]\label{prop:flat-fibre}
Let $B$ be the Euclidean plane, a half-plane, or a straight strip,
with the product metric on $T\times B$. Let
$D=(d_\alpha^*,\sqrt2\,d_\alpha^+)$ act on $\ad$-valued one-forms,
with $a(\nu)=0$ and $d_\alpha(a|_T)=0$ on the boundary. Both $D$
and its adjoint are invertible and
$\|D^{-1}\|_{L^2\to L^2}\le\mu_T^{-1}$, independently of the strip
width. Conjugating by $e^{\delta\tau}$, where $|d\tau|\le1$ and
$\tau$ is constant on torus fibres, gives the bound
$(\mu_T-|\delta|)^{-1}$ for $|\delta|<\mu_T$.
\end{proposition}

\begin{proof}
\emph{The operator and its boundary conditions.}
Complexify the flat adjoint bundle and decompose into its twisted
Fourier modes. There is no zero mode. On a mode of torus frequency
$k\ge\mu_T$, rotate the orthonormal torus frame so that its covector
is $k\,dx$ and write $a=f\,ds+g\,dt+p\,dx+q\,dy$.
In the isometric identification of self-dual forms with three
components, the four components of $Da$ are
\[
-f_s-g_t-ikp,\quad g_s-f_t+ikq,\quad
p_s-ikf-q_t,\quad q_s+p_t-ikg.
\]
At a straight boundary the conditions are $f=q=0$: the latter is
exactly the closedness of the torus restriction for this mode.
Integration by parts cancels all cross terms and gives
\[
\|Da\|^2=\|\partial_sa\|^2+\|\partial_ta\|^2+k^2\|a\|^2.
\]
The boundary Green form pairs the free components $g,p$ with the
second and third output components. Thus those output components
vanish in the adjoint domain. In invariant notation these conditions
are $\omega(\nu,\tau_B)=0$ and
$d_\alpha^*((\iota_\nu\omega)|_T)=0$; the scalar output is free.
The same calculation for the adjoint gives the identical coercive
bound. Summing over the torus modes proves the estimates on the
plane, half-plane and strip. The closed-range theorem and the
vanishing adjoint kernel give surjectivity. The argument uses no
constant depending on the distance between the two strip boundaries.
Conjugation by the weight adds a zeroth-order principal-symbol term
of norm at most $|\delta|$; a Neumann series proves the weighted claim.
\end{proof}

\begin{corollary}[the stationary operator at the odd flat torus]\label{cor:flat-fibre-line}
Let $I\subset\RR$ be a compact interval, a closed half-line or $\RR$, with
coordinate $r$. On pairs $\zeta=(\varphi,b)$ of an $\ad$-valued function and
one-form on $T\times I$, with the product metric and the pulled-back flat
connection $\alpha$, let
\[
\mathcal A_\alpha(\varphi,b)=\bigl(d_\alpha^*b,\ *d_\alpha b+d_\alpha\varphi\bigr),
\]
with the domain given by $b(\partial_r)=0$ and $d_\alpha\bigl(b|_{T\times\{r\}}\bigr)=0$
at every endpoint $r$ of $I$, and no condition on $\varphi$. Then
$\mathcal A_\alpha$ is self-adjoint and $\|\mathcal A_\alpha\zeta\|\ge\mu_T\|\zeta\|$.
If $\tau=\tau(r)$ satisfies $|\tau'|\le1$ and $|\delta|<\mu_T$, then
$\|e^{\delta\tau}\mathcal A_\alpha\zeta\|\ge(\mu_T-|\delta|)\|e^{\delta\tau}\zeta\|$.
\end{corollary}

\begin{proof}
Let $\zeta=(\varphi,b)$ lie in the domain and have compact support, and let
$\beta\in C^\infty_c(\RR)$ be a nonzero cutoff. On $T\times B$ with
$B=\RR_t\times I$, a plane, a half-plane or a strip, put
$a_N=\beta(t/N)\,(-\varphi\,dt+b)$. It satisfies the boundary conditions of
Proposition~\ref{prop:flat-fibre}: $a_N(\partial_r)=\beta(t/N)b(\partial_r)=0$,
and its restriction to each boundary torus fibre is closed. For a one-form
$a=\psi\,dt+c$ on $T\times B$ one has $d^*a=-\partial_t\psi+d^*c$ and
$\sqrt2\,d^+a$ corresponds, under the identification of $\Lambda^+$ with the
one-forms of the three-dimensional slice used in \cite[(4.16)]{DFL}, to
$\partial_tc+*d_\alpha c-d_\alpha\psi$, isometrically. Hence
$\|Da_N\|=\|(\partial_t+\mathcal A_\alpha)(\beta(t/N)\zeta)\|$, and
\[
\|Da_N\|^2=N\|\beta\|^2\|\mathcal A_\alpha\zeta\|^2+N^{-1}\|\beta'\|^2\|\zeta\|^2,
\]
because the cross term is a multiple of $\int(\beta^2)'=0$. Proposition~\ref{prop:flat-fibre}
gives $\|Da_N\|^2\ge\mu_T^2\|a_N\|^2=\mu_T^2N\|\beta\|^2\|\zeta\|^2$. Dividing by $N$
and letting $N\to\infty$ gives the inequality for compactly supported $\zeta$, and
the general case follows by approximation. The weight $\tau(r)$ is constant on
the torus fibres and independent of $t$, so the weighted version of
Proposition~\ref{prop:flat-fibre} gives the weighted inequality in the same way.
For $\zeta_1,\zeta_2$ in the domain, Green's formula leaves the boundary terms
$\pm\int_{T\times\partial I}\langle b_1\wedge b_2\rangle$ and terms containing
$b_i(\partial_r)$; the latter vanish, and the former vanish because the torus
restrictions are closed, hence exact since $H^1(T;\ad\alpha)=0$, and
$\int_T\langle d_\alpha u\wedge d_\alpha v\rangle=\int_Td\langle u,d_\alpha v\rangle=0$.
The boundary condition is elliptic in the sense of \cite{BB}, the closed and the
exact conditions on $T$ coinciding, so $\mathcal A_\alpha$ is self-adjoint; with
the lower bound it is invertible.
\end{proof}

\begin{theorem}[regularity at the flat torus boundary]\label{thm:flat-boundary-regularity}
For $\rho>0$ let $D^+_\rho=\{(s,t):s\ge0,\ s^2+t^2<\rho^2\}$ and
$Q_\rho=T\times D^+_\rho$. Let $A$ be an $L^p_1$ connection on $Q_2$, $p>2$,
anti-self-dual, with weakly flat restriction to $T\times\{z\}$ for $z$ on the
boundary diameter.
\begin{enumerate}[label=(\alph*)]
\item $A$ is gauge equivalent to a connection smooth up to the boundary.
\item There are $\hbar>0$ and $C_k$, depending only on $T$ and $k$, such that if
  $\epsilon^2=\int_{Q_2}|F_A|^2\le\hbar$ then $\sup_{Q_1}|\nabla^k_AF_A|\le C_k\epsilon$,
  and in a Coulomb--Neumann gauge $A=\alpha+a$ with $\|a\|_{C^k(Q_1)}\le C_k\epsilon$.
\item Suppose that $A$ satisfies the flat-fibre condition and has finite
  energy on the half-space $T\times\{s\ge0\}$. Then
  \[
  \sup_{T\times D_1(z)}|\nabla^k_AF_A|\longrightarrow0\qquad(|z|\to\infty).
  \]
\end{enumerate}
\end{theorem}

\begin{proof}
(a) is the regularity theorem of \cite[Theorem~1]{DFL-mixed} with the torus as
matching surface; its proof uses the Lagrangian boundary condition only through
local charts, and applies to the one-point moduli space. We prove (b) in three steps.

\emph{A fixed Neumann slice.} Choose rounded contractible domains
$\Omega_1\Subset\Omega_2\Subset D^+_2$, with relative closure along
the physical diameter, such that $D^+_1\Subset\Omega_1$.
The qualitative compactness argument of
\cite[Theorem~4.2, Lemma~4.25 and pp.~38--39]{DFL-mixed}, applied
with the one-point torus moduli space, implies that connections
whose energy tends to zero converge smoothly in gauges on
$T\times\Omega_2$ to a flat connection. That connection is gauge
equivalent to $\alpha$: the base is contractible and the odd torus
has a single flat class. On a slightly smaller rounded domain,
still denoted by $\Omega_2$, impose
\[
d_\alpha^*a=0,\qquad a(\nu)=0
       \quad\hbox{on the whole boundary of }T\times\Omega_2.
\]
There is no prescribed gauge on the artificial boundary. The scalar
Neumann linearization is invertible since $H^0(T;\ad\alpha)=0$,
so the slice follows from the implicit function theorem. In this
slice $a$ tends to zero smoothly as the energy tends to zero.
The absolute twisted cohomology in degree one vanishes on
$T\times\Omega_2$. The Hodge estimate and four-dimensional
Sobolev embedding therefore give
\[
\|a\|_{H^1}\le C\|d_\alpha a\|_2
 \le C\bigl(\epsilon+\|a\|_4^2\bigr)
 \le C\epsilon+C'\|a\|_{H^1}^2.
\]
Qualitative smallness absorbs the last term. Thus
$\|a\|_{H^1(T\times\Omega_2)}\le C\epsilon$.

\emph{The boundary estimate.} On the physical face let
$P_{\rm coex}$ denote fibrewise coexact projection and set
$Q_T=d_\alpha^*(d_\alpha d_\alpha^*)^{-1}$ on torus two-forms.
Flatness is equivalent to
\[
P_{\rm coex}a_T=-\tfrac12Q_T[a_T\wedge a_T].
\]
The operator $Q_T$ is bounded on each isotropic Sobolev space
on $T\times\RR_t$; no gain in the $t$ direction is asserted.
For a form $w$ supported in a coordinate half-disk and with zero
normal component, the linear estimate is
\[
\|w\|_{H^{j+1}}\le C_j\bigl(\|Dw\|_{H^j}
 +\|P_{\rm coex}w_T\|_{H^{j+1/2}}+\|w\|_2\bigr).
\]
For zero coexact trace, its $H^1$ version is the Fourier calculation
in Proposition~\ref{prop:flat-fibre}. Tangential derivatives commute
with $D$ and preserve the boundary condition; the equation recovers
normal derivatives. For nonzero trace, subtract its Poisson extension,
which maps $H^{j+1/2}$ to $H^{j+1}$ and has zero normal component.
This proves the displayed estimate at every integer $j\ge0$.

Choose a scalar base cutoff $\chi$ and apply the estimate to
$w=\chi a$. Its equation and boundary value are
\[
Dw=-\sqrt2\chi(a\wedge a)^++[D,\chi]a,
\qquad P_{\rm coex}w_T=-\tfrac12Q_T[a_T\wedge w_T].
\]
At each fixed order, qualitative smooth smallness of $a$ makes
multiplication by $a$, both in the interior and on the boundary,
small in the displayed norms. The trace theorem therefore absorbs
the terms containing the highest derivatives of $w$ into the same
$H^{j+1}$ norm on the left. The cutoff commutator contains only
$j$ derivatives of $a$. Induction on nested half-disks, starting
with the $H^1$ estimate, proves the $H^{j+1}$ bound $C_j\epsilon$.
Interior patches have the same argument without the trace term.
Sobolev embedding then gives the asserted connection and curvature
bounds at each fixed order, initially below an order-dependent
threshold $\hbar_k$.

\emph{A common energy threshold.} Fix $\hbar$ below the concentration
threshold in the cited compactness theorem and small enough to
enter the single low-order Neumann slice on $T\times\Omega_2$.
The family of energy at most $\hbar$ is smoothly precompact modulo
gauge on smaller domains: any failure would give concentration,
whose energy is bounded below by that threshold. The normalized
Neumann representatives are consequently smoothly precompact on
$T\times\Omega_1$; the slice correction is unique near the identity
and elliptic regularity controls it at every higher order. Let
$M_k$ be a bound for their $C^k$ connection and curvature norms.
For energies $E\in[\hbar_k,\hbar]$ use
$M_k\le(M_k/\sqrt{\hbar_k})\sqrt E$, and for $E<\hbar_k$ use the
linear estimate just proved. Increasing $C_k$ gives (b) with one
$\hbar$ for all $k$.

(c) Given $\epsilon>0$, choose a compact set outside which the energy is less than
$\min(\epsilon^2,\hbar)$; for $|z|$ large $T\times D_2(z)$ lies outside it, and (b), or its
interior version (Remark~\ref{rem:interior-boxes}), gives the bound.
\end{proof}

\begin{remark}\label{rem:interior-boxes}
Part (b) of Theorem~\ref{thm:flat-boundary-regularity} holds, with constants
depending only on $T$ and $k$, on interior boxes $T\times D_\rho(z)$ disjoint
from the face and with no boundary condition. Below the threshold the curvature
bounds are Uhlenbeck's interior estimates \cite[Section~2.3]{Donaldson-book},
applied on balls of a fixed radius, and the gauge statement follows from the
contradiction argument of (b), since a flat connection on $T\times D_2(z)$ is
gauge equivalent to $\alpha$. The local gauges may be taken even: the stabilizer of
$\alpha$ is the Klein four-group, whose three nontrivial elements are odd and realize
the three nonzero classes of $H^1(T;\FF)$, so composing with one of them makes any
local gauge even without moving $\alpha$. Below, \emph{small-energy regularity} refers to (b)
and to this interior version.
\end{remark}

\subsection{The stationary problems and uniform decay}

\begin{proposition}[uniform stationary gap]\label{prop:stationary-gap}
Let $q\in P$ and $\mathsf R\in[0,\infty]$. The stationary operator
$\mathcal A_{q,\mathsf R}(\varphi,b,\nu)=(d^*b,(*d_B+\operatorname{Hess}h)b+d\varphi,J\partial_\sigma\nu)$
on $Y_\mathsf R$ and the section interval, with the matching conditions at the surface,
the Lagrangian condition at the end of the interval and the flat-fibre
condition at $r=\mathsf R<\infty$, is self-adjoint with kernel $T_qP$, and its nonzero
spectrum avoids $(-\gamma,\gamma)$ for a $\gamma>0$ independent of $q$ and $\mathsf R$. The same
holds for the joined slices $D_{i,\lambda}$, with kernel $T(P_i\times Q_i)$,
uniformly in $\lambda\ge\lambda_*$, for a width threshold $\lambda_*>0$. The stationary representatives $B_q$ may be taken
equal to $\alpha$ on the collar.
\end{proposition}

\begin{proof}
Self-adjointness uses the Green cancellation at the matching surface
(Convention~\ref{conv:orientation}) and the elliptic boundary theory of
\cite{BB}. The kernel is computed as in Proposition~\ref{prop:transition-gap}.
For uniformity, suppose that there are $q_n\in P$, $\mathsf R_n\in[0,\infty]$ and $\zeta_n$
orthogonal to the kernel with $\|\zeta_n\|=1$ and $\|\mathcal A_{q_n,\mathsf R_n}\zeta_n\|\to0$. Let
$\chi$ be a cutoff of $r$, equal to $0$ for $r\le1$ and to $1$ for $r\ge2$; it commutes with
the primary Hessian, which is supported away from the collar. On the collar the
operator is the stationary operator $\mathcal A_\alpha$ at the odd flat torus on
$T\times[0,\mathsf R_n]$, so Corollary~\ref{cor:flat-fibre-line} gives
\[
\mu_T\|\chi\zeta_n\|\le\|\mathcal A_\alpha(\chi\zeta_n)\|\le\|\mathcal A_{q_n,\mathsf R_n}\zeta_n\|
+\|\chi'\zeta_n\|_{L^2(1\le r\le2)} .
\]
Hence $\|\zeta_n\|_{L^2(K)}\ge1/C$ for the fixed compact region $K=\{r\le2\}$ of the slice
together with the section interval. Local elliptic estimates give a bound in
$H^1(K)$, and after passing to a subsequence $q_n\to q$, $\mathsf R_n\to\mathsf R\in[0,\infty]$
(for a finite limit the slices $Y_{\mathsf R_n}$ are transported to $Y_{\mathsf R}$ by stretching
maps as for $\psi_{\mathsf R}$), and
$\zeta_n$ converges weakly to a limit $\zeta\ne0$, strongly in $L^2(K)$, with
$\mathcal A_{q,\mathsf R}\zeta=0$ and the limiting boundary conditions. The kernel elements
satisfy $\mathcal A_\alpha\zeta=0$ on the collar and therefore decay in weighted $L^2$ like
$e^{-\delta r}$ for every $\delta<\mu_T$, by Corollary~\ref{cor:flat-fibre-line} applied with the
truncated weights $\min(r,M)$, uniformly in $q$; so the kernels converge strongly
and orthogonality to the kernel passes to the limit;
$\zeta$ is a nonzero kernel element orthogonal to the kernel, a contradiction.
\end{proof}

\begin{proposition}[uniform decay on growing mixed ends]\label{prop:growing-ends}
There are $e_0,c>0$ and, for every $m\ge0$, numbers $e'_m\in(0,e_0]$ and $C_m$, all
independent of $\mathsf R\in[0,\infty]$ and of $\lambda\ge\lambda_*$, with the following
properties. Let $U$ be a finite-energy solution of the autonomous mixed equation on a
half-cylinder $[\theta_0,\infty)$ over $Y_\mathsf R$ (or $D_{i,\lambda}$), and let $E(\theta)$ be its
kinetic energy beyond $\theta$.
\begin{enumerate}[label=(\roman*)]
\item If $E(\theta_0)<e_0$, then $E(\theta)\le e^{4c}E(\theta_0)e^{-2c(\theta-\theta_0)}$ for
  $\theta\ge\theta_0$, and $U$ converges to a stationary solution $(B_q,q)$.
\item (A priori estimate.) If moreover $U$ is smooth and $E(\theta_0)<e'_m$, then the limit of
  $U$ in a temporal gauge is smooth, and in the temporal gauge normalized by a
  $\theta$-independent even gauge transformation so that the limit is exactly $\pi^*B_q$,
  \begin{equation}\label{eq:end-decay}
  \|U-\pi^*B_q\|_{H^m([\theta,\theta+1]\times Y_\mathsf R)}\le C_mE(\theta_0)^{1/2}e^{-c(\theta-\theta_0)}
  \qquad(\theta\ge\theta_0+2),
  \end{equation}
  where the section minus the constant $q$ is included in the norm (and $Y_\mathsf R$ is
  replaced by $D_{i,\lambda}$ in the joined case).
\end{enumerate}
\end{proposition}

\begin{proof}
In a temporal gauge the solution is a path $\theta\mapsto(B(\theta),\nu(\theta))$ of slice
configurations, a connection on the slice and a path in $\Msf$ across the section interval,
subject to the matching condition at the surface, the Lagrangian condition at the end of
the interval and, for $\mathsf R<\infty$, the flat-fibre condition at $r=\mathsf R$. It is the downward
gradient flow of the local action $\mathcal F$, the Chern--Simons functional with the primary
function plus the symplectic action of $\nu$. Let $Y'$ be the fixed compact part of the slice
consisting of $Y$ (of the two genus pieces, in the joined case) and the part $r\le3$ of the
collar (of the bridge); it contains the loops of the primary perturbation. On the rest of
the slice no perturbation acts, and anti-self-duality on the product gives
$|F_A|^2=2|\iota_{\partial_\theta}F_A|^2$ there.

\emph{The neighbourhood.} For $\epsilon>0$ let $\mathcal N_\epsilon$ be the set of slice
configurations that, after an even gauge transformation of the slice, differ from a
stationary pair $(B_q,q)$ by a field $\xi$ with $\|\xi\|_{L^2_1}+\|\xi\|_{C^1}\le\epsilon$, the norms
being taken over the whole slice with the covariant derivative of $B_q$. For every $\epsilon>0$
there is $e_0>0$, independent of $\mathsf R$ and $\lambda$, such that a solution whose kinetic
energy on $[\theta-2,\theta+2]$ is less than $e_0$ has its slice at $\theta$ in $\mathcal N_\epsilon$. On
$Y'$ this follows by contradiction from the local compactness theorem
\cite[Theorem~5.2]{DFL}, \cite[Theorem~3]{DFL-mixed}, in its local form: a sequence of
solutions on $[\theta-2,\theta+2]\times Y''$, where
$Y''=Y\cup\{r\le4\}$, whose kinetic energy tends to zero converge,
after gauge transformations, in $C^\infty$ on
$[\theta-1,\theta+1]\times Y'$ to a stationary pair. For
$\mathsf R\in[0,4]$ use $\psi_\mathsf R$ and
Theorem~\ref{thm:flat-boundary-regularity} at the physical face;
compactness of this parameter interval makes the bounds uniform.
Part (i) is gauge invariant. On the unperturbed collar interior
regularity and Theorem~\ref{thm:flat-boundary-regularity}(a) give
a smooth representative, with local gauges patched as in Step~1
of Lemma~\ref{lem:adapted}. Apply the torus slice lemma to this
representative. On the collar,
small-energy regularity (Theorem~\ref{thm:flat-boundary-regularity}(b),
Remark~\ref{rem:interior-boxes}) at orders at most three on the boxes $T\times D_2(z)$, and
Lemma~\ref{lem:torus-slice} below on the base $[1,\mathsf R]\times[\theta-1,\theta+1]$, with
\eqref{eq:torus-slice} for $k\le2$, write the connection as $\alpha+b+\phi_1dr+\phi_2d\theta$ with
$b,\phi$ and their first derivatives bounded pointwise by the fibre norms of $F_A$ and its
covariant derivatives of order at most three. These are bounded by the $L^2$ norms of $F_A$ on
the boxes; summing over a cover of bounded overlap bounds the $L^2_1$ and $C^1$ norms of the
slice field $b+\phi_1dr$ on the collar by $CE(\theta-2)^{1/2}$, with $C$ independent of $\mathsf R$
and $\lambda$. On the overlap $1\le r\le3$ the two gauges differ by a gauge
transformation which maps a configuration close to $\alpha$ to one close to $\alpha$; it is
close to the stabilizer of $\alpha$, the Klein four-group of constant parallel elements. If
it is close to a nontrivial element $v$, replace the collar gauge by its composite with
$v$, which fixes $\alpha=B_q|_{\rm collar}$. The transition is then close to the even
stabilizer, which is trivial, so it is $\exp\varphi$ with $\varphi$ small, and a cutoff of
$\varphi$ in $r$ joins the two gauges. The Coulomb--Neumann implicit function theorem
below is uniform in $\mathsf R$ and $\lambda$ because the Neumann scalar Laplacian has a
uniform gap on $Y_{\mathsf R}$ and $D_{i,\lambda}$, by the cutoff argument of
Proposition~\ref{prop:stationary-gap} and $\Delta_\alpha\ge\mu_T^2$ on the collar.

\emph{Gradient inequality.} Fix $\epsilon$ small. A configuration in $\mathcal N_\epsilon$ is written
as $(B_q,q)+Q_q(\xi)$, where $\xi$ satisfies $d^*_{B_q}\xi=0$, has vanishing normal component on
the boundary of the slice, satisfies the linearized matching, Lagrangian and flat-fibre
conditions and is $L^2$-orthogonal to the tangent space of the critical family at $q$; here $q$
is chosen by the implicit function theorem, the critical family being compact modulo gauge
and Morse--Bott, and $Q_q(\xi)=\xi+O(|\xi|^2)$ is the chart of \cite[Section~4.1]{DFL} (the maps
$E_\tau$ on a collar of the surface, and a metric in which the Lagrangian is totally geodesic
on the section interval), together with the fibrewise correction of
Lemma~\ref{lem:face-chart}(a) at the flat face, used here on the single torus at the
three-dimensional face of the slice in the $L^2_1\cap C^1$ topology: the contraction is
fibrewise, and $\|fg\|_{H^{1/2}}\le C\|f\|_{C^1}\|g\|_{H^{1/2}}$ gives
$\|c(b)\|_{H^{3/2}(T)}\le C\|b\|_{C^1}\|b\|_{H^{1/2}}$. In this chart
\begin{gather*}
\mathcal F-\mathcal F_{\rm crit}=\tfrac12\langle\mathcal A_q\xi,\xi\rangle+R_3,\qquad
\nabla\mathcal F=\mathcal A_q\xi+R_2,\\
|R_3|\le C\|\xi\|_{C^1}\|\xi\|_{L^2_1}^2,\qquad\|R_2\|_{L^2}\le C\|\xi\|_{C^1}\|\xi\|_{L^2_1},
\end{gather*}
where $\mathcal A_q=\mathcal A_{q,\mathsf R}$ is the operator of
Proposition~\ref{prop:stationary-gap} and the gradient is identified with a vector at $\xi$
through the derivative of $Q_q$, which differs from the inclusion by $O(\|\xi\|_{C^1})$. The
remainders are the cubic terms of the Chern--Simons functional and of the symplectic
action, the second-order terms of the chart, and those of the primary function, which is
smooth \cite[Section~6]{DFL} and supported in $Y$; $C$ is independent of $\mathsf R$ and $\lambda$
because the slice is covered by unit patches of finitely many isometry types. The
flat-fibre condition removes the boundary term at $r=\mathsf R$ in the Hessian. By
Proposition~\ref{prop:stationary-gap}, $\|\mathcal A_q\xi\|\ge\gamma\|\xi\|$, and the elliptic
estimate $\|\xi\|_{L^2_1}\le C(\|\mathcal A_q\xi\|+\|\xi\|)$, uniform because it is local and the slice has
bounded geometry (as in model (E) of the proof of Proposition~\ref{prop:torus-inverse}), gives $\|\xi\|_{L^2_1}\le C\|\mathcal A_q\xi\|$. For $\epsilon$
small, $\|\nabla\mathcal F\|\ge\frac12\|\mathcal A_q\xi\|$ and
\[
|\mathcal F-\mathcal F_{\rm crit}|\le\tfrac12\|\mathcal A_q\xi\|\|\xi\|+C\epsilon\|\xi\|^2_{L^2_1}
\le C'\|\mathcal A_q\xi\|^2\le4C'\|\nabla\mathcal F\|^2,
\]
with $C'$ depending only on $\gamma$ and on the constants above.

\emph{Decay.} For $\theta\ge\theta_0+2$ the slice lies in $\mathcal N_\epsilon$, and along the solution
$\partial_\theta\mathcal F=-\|\nabla\mathcal F\|^2\le-(\mathcal F-\mathcal F_{\rm crit})/(4C')$. So the kinetic energy
$E(\theta)=\mathcal F(\theta)-\mathcal F_{\rm crit}$ decays at the rate $2c=1/(4C')$ beyond $\theta_0+2$;
since $E$ is nonincreasing, this gives the bound of (i). The length of the path of slices over $[\theta,\infty)$ is at most the sum over
unit intervals of the square roots of their kinetic energies, which decays at the rate $c$;
so the solution converges to a stationary solution $(B_q,q)$.

\emph{The a priori estimate.} Let $U$ be smooth, and let $\zeta=(\iota_{\partial_\theta}F_A,\partial_\theta u)$
be its gauge-covariant $\theta$-derivative, whose $L^2$ norm on $[\theta,\infty)$ is $E(\theta)^{1/2}$.
On the collar, small-energy regularity at order $k$ gives
$\sup_{T\times D_1(z)}|\nabla^k_A\zeta|\le C_k\|F_A\|_{L^2(T\times D_2(z))}$ once the energy of $T\times D_2(z)$ is
below the threshold at order $k$, and that energy is at most $2E(\theta(z)-2)$. On $Y'$, $\zeta$
satisfies the linearized equation by translation invariance.
In temporal gauge its gauge row vanishes:
$d_A^*\zeta=-d_B^*(*F_B+\nabla h(B))=0$, by Bianchi and gauge
invariance of the cylinder function \cite[Proposition~2.10]{DFL}.
At the seam its normal component is $F_A(\partial_\theta,\partial_s)=0$:
anti-self-duality identifies this component with the curvature on
the flat surface slice. Thus it belongs to the boundary domain
of the gauge-fixed linearization. The local estimates of
\cite[Proposition~5.27 and Remark~4.19]{DFL}, stated there over the whole slice, hold
on the larger region $Y''$, which has an artificial boundary at $r=4$: their proof by
contradiction uses only interior estimates near $r=4$. They give
$\sup_{[\theta,\theta+1]\times Y'}|\nabla_A^k\zeta|
\le C_k\|\zeta\|_{L^2([\theta-1,\theta+2]\times Y'')}$
when the kinetic energy is below a threshold depending on $k$.
Choose $e'_m\le e_0$ below these
thresholds and below half the thresholds of small-energy regularity, for $k\le m+1$. Summing the collar bounds over a cover of bounded overlap, for
$\theta\ge\theta_0+2$ and $k\le m+1$,
\begin{equation}\label{eq:zeta-bounds}
\sup_{[\theta,\theta+1]\times Y_\mathsf R}|\nabla^k_A\zeta|+\|\nabla^k_A\zeta\|_{L^2([\theta,\theta+1]\times Y_\mathsf R)}\le C_kE(\theta-2)^{1/2},
\end{equation}
with $C_k$ independent of $\mathsf R$, $\lambda$ and $U$. In a temporal gauge, which is smooth,
$\partial_\theta A=\iota_{\partial_\theta}F_A$. Expanding the covariant derivative of $A(\theta)$ about that of
$A(\theta_0+2)$ and inducting on $k$, \eqref{eq:zeta-bounds} and integration in $\theta$ bound the
derivatives of $A(\theta)-A(\theta_0+2)$ uniformly in $\theta$, and $A(\theta)$ converges in every
$C^k$ to a limit $B$, which is therefore smooth. By (i), $B$ and $B_q$ represent the same
point of the critical manifold, a quotient by the determinant-one group, so $B=g^*B_q$ for an
even gauge transformation $g$ of the slice, smooth because $dg=gB-B_qg$. In the normalized
temporal gauge $A^N=(g^{-1})^*A$ the limit is $B_q$, and
$A^N(\theta)-B_q=-\int_\theta^\infty\iota_{\partial_\theta}F_{A^N}\,d\theta'$. Now
$\nabla_{B_q}=\nabla_{A^N}-[A^N-B_q,\cdot\,]$, so $\nabla^k_{B_q}\zeta$ is $\nabla^k_A\zeta$ plus terms that are
polynomial in the derivatives of order less than $k$ of $A^N-B_q$ and linear in those of $\zeta$. Inducting on $k$ and
integrating from infinity, \eqref{eq:zeta-bounds} and (i) give
\begin{multline*}
\sup_{Y_\mathsf R}|\nabla^k_{B_q}(A^N(\theta)-B_q)|+\|\nabla^k_{B_q}(A^N(\theta)-B_q)\|_{L^2(Y_\mathsf R)}\\
\le C_k\sum_{n\ge0}E(\theta+n-2)^{1/2}\le C'_kE(\theta_0)^{1/2}e^{-c(\theta-\theta_0)}
\end{multline*}
for $k\le m$, the lower-order coefficients being bounded at each stage by the previous
one. The $\theta$-derivatives of $A^N-\pi^*B_q$ are the covariant $\theta$-derivatives of $\zeta$,
bounded by \eqref{eq:zeta-bounds}, and the section is treated in the same way with
$\partial_\theta u$. This proves \eqref{eq:end-decay}.
\end{proof}

\subsection{Decay along the infinite torus collar}\label{subsec:torus-slice}

A gauge transformation of the odd bundle over a region of $T\times\RR^2$ is
\emph{even} if it lifts to the determinant-one group. The even stabilizer of
$\alpha$ is trivial (Section~\ref{subsec:torus-domains}); the full stabilizer is
the Klein four-group, whose three nontrivial elements are odd.
For a connection $A$ on $T\times\Omega$, $\Omega\subset\RR^2$ with coordinates
$z=(z_1,z_2)$, write $F_{jT}=(\iota_{\partial_{z_j}}F_A)|_{T\times\{z\}}$ and $F_T$ for the
restriction of $F_A$ to a torus fibre.

\begin{lemma}[torus slice]\label{lem:torus-slice}
There are $e_T>0$ and, for each $k\ge0$, constants $C_k$ depending only on $T$,
with the following property. Let $\Omega\subset\RR^2$ be a simply connected
domain and $A$ a smooth connection on $T\times\Omega$ such that every $z\in\Omega$
has a neighbourhood $N$ and an even gauge transformation on $T\times N$ after
which $\|A-\alpha\|_{C^1(T\times N)}\le e_T$. Then there is a unique even gauge
transformation $g$ of $T\times\Omega$ with
\[
g^*A=\alpha+b+\phi_1\,dz_1+\phi_2\,dz_2,\qquad d_\alpha^*b(z)=0,\qquad
\|b(z)\|_{C^0(T)}\le C_0e_T\quad(z\in\Omega),
\]
where $b(z)$ is a one-form and the $\phi_j(z)$ are sections of $\ad$, both on the
fibre $T$. If moreover $\|F_A\|_{C^{k+1}_A(T\times\Omega)}\le e_T$, then for every
$z\in\Omega$
\begin{equation}\label{eq:torus-slice}
\sum_{|\beta|+j\le k+1}\bigl(\|\partial_z^\beta b(z)\|_{H^j(T)}
+\|\partial_z^\beta\phi(z)\|_{H^j(T)}\bigr)
\le C_k\sum_{i\le k+1}\|(\nabla_A^iF_A)|_{T\times\{z\}}\|_{L^2(T)} .
\end{equation}
The fibrewise equations gain derivatives only in the torus directions, so the
top base derivatives of $\phi$ require the curvature derivative of order $k+1$.
We call $g^*A$ the \emph{torus slice} of $A$.
\end{lemma}

\begin{proof}
\emph{The fibre slice.} Consider the map $(\xi,b)\mapsto\exp(\xi)^*(\alpha+b)$
from $\Omega^0(T;\ad)\times\ker d_\alpha^*$ to connections on the odd bundle over
$T$, in $C^{1,\beta}\times C^{0,\beta}$. Its derivative at $(0,0)$ is
$(\xi,b)\mapsto d_\alpha\xi+b$, an isomorphism onto $C^{0,\beta}$ one-forms, because
$\Omega^1(T;\ad\alpha)=\im d_\alpha\oplus\ker d_\alpha^*$ and $d_\alpha$ is injective
($H^0=H^1=0$). By the implicit function theorem every connection $C^1$-close to
$\alpha$ is $\exp(\xi)^*(\alpha+b)$ with $d_\alpha^*b=0$ and $\xi$, $b$ small. The
decomposition is unique among even gauge transformations: if $g$ is even and
$g^*(\alpha+b)=\alpha+b'$ with $b,b'$ small in $C^0$, then $dg=g(\alpha+b')-(\alpha+b)g$
bounds $g$ in $C^1$, and as $b,b'\to0$ every limit of $g$ lies in the stabilizer of
$\alpha$ and is even, hence is the identity; so for $e_T$ small $g$ is close to
the identity and equals the gauge transformation given by the implicit function
theorem.

\emph{The global gauge.} By hypothesis and the fibre slice, each fibre
$T\times\{z\}$ carries a unique even $g(z)$ putting $A|_{T\times\{z\}}$ into the
slice; by uniqueness and the implicit function theorem with parameters, $g(z)$
depends smoothly on $z$. The resulting family is fibrewise even, and it lifts to
the determinant-one group over $T\times\Omega$ because the obstruction lies in
$H^1(T\times\Omega;\FF)=H^1(T;\FF)$ and vanishes on a fibre. Uniqueness of $g$ follows
from the fibrewise uniqueness.

\emph{Estimates.} On a fibre, $F_T=d_\alpha b+\tfrac12[b\wedge b]$. Since
$H^1(T;\ad\alpha)=0$, the elliptic estimate for $d_\alpha+d_\alpha^*$ gives
$\|b\|_{H^{j+1}(T)}\le C\|d_\alpha b\|_{H^j(T)}\le C(\|F_T\|_{H^j(T)}+\|b\|_{C^0}\|b\|_{H^{j+1}(T)})$,
and the last term is absorbed. For the base components,
$F_{jT}=\partial_{z_j}b-d_{\alpha+b}\phi_j$. Applying $d_\alpha^*$ and using
$d_\alpha^*\partial_{z_j}b=\partial_{z_j}d_\alpha^*b=0$ gives
\[
\Delta_\alpha\phi_j+d_\alpha^*[b,\phi_j]=-d_\alpha^*F_{jT},
\]
and $\Delta_\alpha\ge\mu_T^2$ yields $\|\phi_j\|_{H^{j'+1}(T)}\le C\|F_{jT}\|_{H^{j'}(T)}$.
Then $\partial_{z_j}b=F_{jT}+d_{\alpha+b}\phi_j$ controls the first base derivatives
of $b$. Differentiating these identities in $z$, and writing
$\nabla_A=\nabla_\alpha+[b+\phi\,dz,\cdot\,]$, an induction on $|\beta|$ gives
\eqref{eq:torus-slice}: for $|\beta|=k+1$ the equation for $\partial_z^\beta\phi_j$ has
source $d_\alpha^*\partial_z^\beta F_{jT}$, a curvature derivative of order $k+1$, and
the lower-order products are absorbed because $\|F_A\|_{C^{k+1}_A}\le e_T$.
\end{proof}

\begin{proposition}[decay along the completed tail]\label{prop:tail-decay}
Let $V=(A,u)$ be a finite-energy solution on $W_\infty$ (or on $W_{\pm,\infty}$).
Let $\gamma_T=\mu_T/4$ and $0\le\delta\le\delta_T$, where
\[
\delta_T=\min\{\mu_T/4,\ \mu_T^2/8,\ c/2\},
\]
$c$ is the rate of Proposition~\ref{prop:growing-ends}, and the weight on the torus region is
$\tau_0(\theta)$ with $|\tau_0'|\le1$ and $|\tau_0''|\le1$. There is $S_0\ge3$ such that:
\begin{enumerate}[label=(\roman*)]
\item for every $k$ there is $C_k$ with $\sup_{T\times D_1(z)}|\nabla_A^kF_A|\le C_ke^{-\gamma_T(r(z)-S_0)}$
  for every $z$ with $r(z)\ge S_0+2$, uniformly in $\theta(z)$;
\item Lemma~\ref{lem:torus-slice} applies on $\Omega=(S_0,\infty)\times\RR_\theta$, and for every
  $m$ there are $S_m\ge S_0$ and $C_m$ such that in the torus slice $A=\alpha+a_V$
  \[
  \|a_V\|_{H^m_\delta(T\times[S,\infty)\times\RR)}\le C_me^{-\gamma_T(S-S_0)}\qquad(S\ge S_m+2),
  \]
  where $H^m_\delta$ is weighted by $e^{\delta\tau_0(\theta)}$.
\end{enumerate}
\end{proposition}

\begin{proof}
\emph{Smallness.} The energy of $V$ on $T\times[S,\infty)\times\RR$ tends to zero as
$S\to\infty$, since these regions decrease to the empty set and carry no
perturbation. By small-energy regularity (Remark~\ref{rem:interior-boxes}) applied on the
boxes $T\times D_2(z)$ with $r(z)\ge S_0-1$, for $S_0$ large we have
$\sup_{r\ge S_0-1}|\nabla_A^kF_A|\le e_T$ for $k\le2$, together with local even gauges
in which $\|A-\alpha\|_{C^1}\le e_T$. So Lemma~\ref{lem:torus-slice} applies on
$(S_0-1,\infty)\times\RR$, and on each fibre $A=\alpha+b$ with $\|b\|_{C^0}\le Ce_T$. This
fixes $S_0$, independently of $k$ and $m$.

\emph{A differential inequality.} Put $e(z)=\int_{T\times\{z\}}|F_A|^2$. On the flat
product the anti-self-dual curvature satisfies the Weitzenböck identity
$\nabla_A^*\nabla_AF_A=-\{F_A,F_A\}$, with no curvature term and
$|\langle\{F,F\},F\rangle|\le c_W|F|^3$. Summing $\partial_i^2|F|^2$ over the four
coordinate directions and integrating over $T$, where the torus directions
integrate to zero,
\[
(\partial_{z_1}^2+\partial_{z_2}^2)e\ \ge\ 2\int_T|\nabla_{A,T}F_A|^2-2c_W\sup|F_A|\,e .
\]
The components of $F_A$ in the parallel frame of $\Lambda^2(\RR^4)$ are sections of
$\ad$ over the fibre, and $\int_T|\nabla^\alpha_Tw|^2=\langle\Delta_\alpha w,w\rangle\ge\mu_T^2\|w\|^2$
for such sections. With $A|_T=\alpha+b$ and $\|b\|_{C^0}$ small this gives
$\int_T|\nabla_{A,T}F_A|^2\ge\tfrac34\mu_T^2e$ for $e_T$ small, and therefore
\begin{equation}\label{eq:fibre-subharmonic}
(\partial_{z_1}^2+\partial_{z_2}^2)e\ge\mu_T^2e\qquad\text{on }r\ge S_0-1 .
\end{equation}

\emph{Integrability in $\theta$.} Let $E_\pm(\theta)$ be the kinetic energy of $V$ on the mixed
end $\pm\theta\to\infty$ beyond $|\theta|$, and let $\Theta_*\ge4$ be such that $E_\pm(\Theta_*-2)$ is
less than $e_0$ and than half the threshold of small-energy regularity at order two. For
$r(z)\ge S_0-1$ and $|\theta(z)|\ge\Theta_*$ the box $T\times D_2(z)$ lies in the unperturbed
collar of a mixed end, where $|F_A|^2=2|\iota_{\partial_\theta}F_A|^2$ by anti-self-duality; so its
energy is at most $2E_\pm(|\theta(z)|-2)\le Ce^{-2c|\theta(z)|}$, by
Proposition~\ref{prop:growing-ends}(i). Small-energy regularity on the box
(Remark~\ref{rem:interior-boxes}) gives $\sup_{T\times D_1(z)}|\nabla^i_AF_A|\le Ce^{-c|\theta(z)|}$ for
$i\le2$, uniformly in $r$. The function $e$ and its first two base derivatives are integrals
over $T$ of expressions quadratic in $F_A$ and its covariant derivatives of order at most
two, so $e+|\partial e|+|\partial^2e|\le Ce^{-2c|\theta|}$ for $|\theta|\ge\Theta_*$, uniformly in
$r\ge S_0-1$, where $\partial$ is any base derivative; for $|\theta|\le\Theta_*$ these quantities
are bounded by the smallness step. Since $2\delta\le c<2c$, the function
$f(r)=\int_\RR e^{2\delta\tau_0(\theta)}e(r,\theta)\,d\theta$ is finite and of class
$C^2$ for $r\ge S_0-1$, and the integrations by parts below have no boundary terms.
Finally, for $|\theta|\ge\Theta_*$ the same identity gives
$\int_\theta^{\theta\pm1}\int_{S_0}^\infty e\,dr\,d\theta'\le2E_\pm(|\theta|)\le Ce^{-2c|\theta|}$; with
$e^{2\delta\tau_0}\le Ce^{2\delta|\theta|}$ and $2\delta<2c$ the sum over unit intervals converges, and
the range $|\theta|\le\Theta_*$ contributes at most a constant times the energy of $V$. Hence
$\int_{S_0}^\infty f(r)\,dr<\infty$.

\emph{Decay.} By~\eqref{eq:fibre-subharmonic},
\[
f''=\int e^{2\delta\tau_0}\bigl((\partial_r^2+\partial_\theta^2)e\bigr)d\theta
-\int(e^{2\delta\tau_0})''e\,d\theta\ \ge\ (\mu_T^2-4\delta^2-2\delta)f\ \ge\ \tfrac12\mu_T^2f,
\]
using $(e^{2\delta\tau_0})''=(4\delta^2\tau_0'^2+2\delta\tau_0'')e^{2\delta\tau_0}$ and
$\delta\le\delta_T$. So $f$ is convex, nonnegative and integrable on
$[S_0,\infty)$, hence nonincreasing with limit $0$. Put $\gamma_T'=\mu_T/\sqrt2$ and
$h(r)=f(r)-f(S_0)e^{-\gamma_T'(r-S_0)}$. Then $h''\ge\gamma_T'^2h$, $h(S_0)=0$ and
$h\to0$; at a positive interior maximum $h''\le0<\gamma_T'^2h$, which is impossible.
Hence
\[
f(r)\le f(S_0)e^{-\gamma_T'(r-S_0)}\qquad(r\ge S_0).
\]

\emph{Higher norms.} (i) Fix $k$. The energy of $T\times D_2(z)$ is at most
$\int_{r(z)-2}^{r(z)+2}f$, because $e^{2\delta\tau_0}\ge1$; let $S'_k\ge S_0+2$ be such that this is
below the threshold of small-energy regularity at order $k$ for $r(z)\ge S'_k$. For such $z$,
small-energy regularity bounds $\sup_{T\times D_1(z)}|\nabla^k_AF_A|$ by
$C_k(\int_{r(z)-2}^{r(z)+2}f)^{1/2}\le C'_ke^{-\gamma_T(r(z)-S_0)}$, since $\gamma_T'/2>\gamma_T$. On the
strip $S_0+2\le r(z)\le S'_k$ the same supremum is bounded: for $|\theta(z)|$ large the box lies
in the unperturbed collar of a mixed end with energy below that threshold, by
Proposition~\ref{prop:growing-ends}(i), and small-energy regularity applies; the rest of the
strip is compact, and there $V$ is locally gauge equivalent to a smooth connection by interior
regularity of the anti-self-duality equation, no perturbation acting on the torus region.
Enlarging $C'_k$ gives (i). (ii) Fix $m$, and let $S_m\ge S'_m+2$ be such that, by (i),
$\|F_A\|_{C^m_A}\le e_T$ on $r\ge S_m$. By small-energy regularity, which is linear in the square
root of the local energy, cover each unit square $Q$ by a fixed
finite collection of disks $D_1(z)$ whose larger disks $D_2(z)$
lie in $\{r\ge S'_m\}$. Apply the estimate on these larger disks
for every $i\le m$, taking the minimum of the finitely many energy
thresholds. The weight varies by a bounded factor on each disk.
Summing over a cover with bounded overlap,
\[
\sum_{i\le m}\int_{r\ge S}e^{2\delta\tau_0}\|\nabla_A^iF_A\|^2_{L^2(T)}\,dz
\le C\int_{S-2}^\infty f(r)\,dr\le Cf(S_0)e^{-\gamma_T'(S-2-S_0)}\qquad(S\ge S_m+2) .
\]
On $(S_m,\infty)\times\RR$ Lemma~\ref{lem:torus-slice} applies with $\|F_A\|_{C^m_A}\le e_T$, and by
fibrewise uniqueness its torus slice is the restriction of that on $(S_0,\infty)\times\RR$.
Integrating~\eqref{eq:torus-slice} with $k=m-1$ against the weight gives (ii), since
$\gamma_T'/2>\gamma_T$.
\end{proof}

For each mixed end we use the stationary representatives $B_q$ of
Proposition~\ref{prop:stationary-gap}, equal to $\alpha$ on the torus collar. For $q$
near a given point of the critical manifold we fix a smooth local family
$q\mapsto B_q$ with this property, and in temporal gauge on product collars of the
genus boundary surfaces: the flat restrictions to a product collar of $T$
inside the compression body are gauge equivalent to $\alpha$ by gauge
transformations depending smoothly on $q$, and these are extended over the
compression body by cutoffs of their logarithms; the primary perturbation is
supported away from the collars, where the connections are flat, so the temporal
gauge along a collar is available. On $Y_\mathsf R$ and on $D_{i,\lambda}$ the family is
independent of $\mathsf R$ and $\lambda$ away from the collar, and equal to $\alpha$ on the
collar (on the bridge in the joined case). Its tangent vectors are supported
away from the torus collar and have vanishing normal components at the genus
boundary. More is true: every $B_q$ is in temporal gauge on the product collars of the
genus boundary surfaces, so the normal component of $B_q-B_{q'}$ vanishes on these
collars for all $q,q'$ in the domain of the family, not only to first order. Such a family
is defined near a base point. Where finitely many base points occur, as for the limits of
the finite completed index-zero sets in Section~\ref{subsec:torus-completion}, we fix one
family near each distinct base point, on pairwise disjoint neighbourhoods in the critical
manifold, and solutions with the same limit use the same family.

From now on we fix a level $m\ge4$.

\begin{lemma}[adapted representatives]\label{lem:adapted}
Let $V$ be a finite-energy solution on $W_\infty$ with limits $q_e$ at the mixed
ends $e$. There are $S_V\ge S_m+3$, $\Theta_V$ and a representative $V^{\rm ad}=(A^{\rm ad},u)$ of
$V$, obtained from $V$ by an even gauge transformation of class $H^{m+1}_{\rm loc}$, which is
smooth on $W_\infty$ up to the genus matching faces and the Lagrangian arcs, such that:
\begin{enumerate}[label=(\roman*)]
\item on $T\times[S_V,\infty)\times\RR_\theta$, $A^{\rm ad}=\alpha+a_V$ is the torus slice, so
  Proposition~\ref{prop:tail-decay}(ii) applies to $a_V$; moreover $a_V$ has bounded
  derivatives of every order there, and for every $k$ there are $S(k)$ and $C_k$ with
  $\sup_{T\times D_1(z)}|\nabla^k_\alpha a_V|\le C_ke^{-\gamma_T(r(z)-S_0)}$ for $r(z)\ge S(k)$;
\item on each mixed end $e$, with outward end coordinate $\theta$, for every $k$ there is $C_k$
  with
  \[
  \bigl\|(A^{\rm ad}-\pi^*B_{q_e},\ u-q_e)\bigr\|_{H^k([\theta,\theta+1]\times Y_{e,\infty})}\le C_ke^{-c\theta}
  \qquad(\theta\ge\Theta_V),
  \]
  on the whole slice, its torus part included; in particular, since $\delta<c$,
  $V^{\rm ad}-\pi^*B_{q_e}$ (and the section minus the constant $q_e$) lies in
  $H^m_\delta(\{\theta\ge\Theta_V\})$, and its norm on $\{\theta\ge\Theta\}$ tends to $0$
  as $\Theta\to\infty$;
\item on the gauge end the limit is exactly the fixed representative of the critical
  point, and on the gauge and strip ends $V^{\rm ad}$ converges to its limit in every $C^k$ at
  every exponential rate below the spectral gap of the asymptotic operator
  \cite[Section~4.1]{Donaldson-book}, \cite{RS-strips}.
\end{enumerate}
The same holds on $W_{\pm,\infty}$ and for pairs on $W_{-,\infty}\sqcup W_{+,\infty}$. Solutions on the
finite domains are treated in Lemma~\ref{lem:torus-regularity}.
\end{lemma}

\begin{proof}
\emph{Step 1: a smooth representative.} Every point of $W_\infty$ has a neighbourhood on
which $V$ is gauge equivalent, by a gauge transformation of class $H^{m+1}$, to a smooth
solution: near the genus matching faces, on charts $D_-(\rho)\times\Sigma'$ together with the
half-disks $D_+(\rho)$ of the section, by the seam regularity theorem \cite[Theorem~1]{DFL-mixed}
in a local Coulomb gauge, as in the proof of \cite[Proposition~3.6]{DFL}; for the section, by
elliptic regularity for pseudo-holomorphic maps with Lagrangian boundary conditions
\cite[Theorem~B.4.1]{MS}; at the other points, by elliptic regularity of the perturbed
anti-self-duality equation in a local Coulomb gauge. The holonomy terms are not local, so
the last step is carried out level by level on nested locally finite
covers. In each initial ball use radii $\rho_i(1+2^{-k})$ at
level $k$, with the balls of radius $\rho_i$ still covering.
Fix the local Coulomb gauges at level $m$ and shrink one step
for each interior elliptic estimate. If the chart representatives are of class $L^2_k$ and
the transition functions of class $L^2_{k+1}$, then every holonomy term is of class $L^2_k$ on
every chart: the holonomies along the loops, based at the points of the chart, are products
of parallel transports along subarcs contained in single charts and of transition functions
at the switching points, and the terms are obtained from them by smooth functions and by
integration against smooth compactly supported densities
\cite[Proposition~3.1 and Section~3.2]{Kronheimer-higher}.
The elliptic estimate for
$(d^*,d^+)$ with a source of class $L^2_k$ then gives class $L^2_{k+1}$. No holonomy term acts
near the genus matching faces, on the torus region or on the section (the primary loops lie
in the interiors of the compression pieces, and the secondary perturbations are supported in
the core away from the matching lines), so the seam charts do not enter the induction. By
induction the chart representatives are smooth, and so are the transition functions, since
$dt_{ij}=t_{ij}A_j-A_it_{ij}$ with smooth $A_i$, $A_j$. Hence $V$ is a smooth connection on a smooth
bundle $P'$ with a continuous bundle isomorphism $\Phi\colon P'\to P$. A smooth isomorphism $\Psi$
that is $C^0$-close to $\Phi$ exists, and $g=\Psi\Phi^{-1}$ is $C^0$-close to the identity, hence even;
$V^{\rm sm}=g\cdot V$ is smooth.

\emph{Step 2: the ends.} (a) Let $K_e=Y_{e,\infty}\cap\{r\le S_V+3\}$, a compact part of the slice
containing $Y$, and let $\Theta_0$ be such that on each mixed end the kinetic energy of $V$ beyond
$\Theta_0$ is below $e'_0$. By Proposition~\ref{prop:growing-ends}(ii) the temporal gauge of $V^{\rm sm}$
on $\{\theta\ge\Theta_0\}$ has a smooth limit, and after a $\theta$-independent even gauge transformation,
smooth because it relates two smooth stationary connections, the limit is $\pi^*B_{q_e}$; call
the result $V^N$. The normalization is even because the critical manifolds are quotients by
the determinant-one group. For every $k$, the $C^k$ norm of $V^N-\pi^*B_{q_e}$ on
$[\theta,\theta+1]\times K_e$, section minus $q_e$ included, is at most $C_ke^{-c\theta}$ for $\theta\ge\Theta_0$: for
large $\theta$ by~\eqref{eq:end-decay} with $m=k+3$, applied beyond a $\theta$ at which the kinetic
energy is below $e'_{k+3}$, and Sobolev embedding; on the remaining finite range by
compactness.
(b) By the smallness step of the proof of Proposition~\ref{prop:tail-decay},
Lemma~\ref{lem:torus-slice} applies to $V^{\rm sm}$ on $(S_V-1,\infty)\times\RR$; the fibrewise implicit
function theorem with parameters, at smooth data, gives a smooth gauge, so the torus slice
$V^T=\alpha+a_V$ is smooth.
(c) Let $\Theta_1\ge\Theta_0$ be large. On $\{\theta\ge\Theta_1,\ S_V-1\le r\le S_V\}$ the representatives $V^N$
and $V^T$ are related by an even gauge transformation $h$, smooth because it relates smooth
connections. It maps a
connection close to $\alpha$ to one close to $\alpha$ (the first because $B_{q_e}=\alpha$ on the
collar), so it is close to the even stabilizer of $\alpha$, which is trivial, and $h=\exp\varphi$ with
$\varphi$ small. From $\nabla_\alpha h=h\,a_V-(V^N-\alpha)h$ and the injectivity of $d_\alpha$ on
$\ad$-valued functions on each fibre,
\[
\|\varphi\|_{C^{k+1}(T\times D_1(z))}\le C_k\bigl(\|a_V\|_{C^k(T\times D_2(z))}+\|V^N-\alpha\|_{C^k(T\times D_2(z))}\bigr).
\]
Both terms are $O(e^{-c\theta(z)})$: the second by (a); the first by~\eqref{eq:torus-slice} on
$\{r>S_V-1,\ \theta>\Theta\}$ for large $\Theta$, whose torus slice is the restriction of $V^T$ by
fibrewise uniqueness, because small-energy regularity and
Proposition~\ref{prop:growing-ends}(i) give
$\sup_{T\times D_1(z)}|\nabla^i_AF_A|\le C_iE(\theta(z)-2)^{1/2}\le C'_ie^{-c\theta(z)}$, uniformly in $r$.
(d) Put $V^{E}=\exp(\chi(r)\varphi)^*V^N$, with $\chi(r)=0$ for $r\le S_V-1$ and $\chi(r)=1$ for $r\ge S_V$. It
agrees with $V^T$ for $r\ge S_V$ and with $V^N$ for $r\le S_V-1$.
(e) On the gauge and strip ends, $V^{\rm sm}$ converges in a temporal gauge, in every $C^k$ and at
every rate below the spectral gap, to a nondegenerate critical point, respectively to a
transverse intersection point \cite[Section~4.1]{Donaldson-book}, \cite{RS-strips}; a
smooth normalization independent of the end coordinate makes the gauge-end limit the fixed
representative.

\emph{Step 3: assembly.} Let $\Theta_V\ge\Theta_1$. On the outer region $O=\{r\ge S_V\}\cup\{|\theta|\ge\Theta_V\}$, together with
the gauge and strip ends, Step 2 gives representatives related to $V^{\rm sm}$ by smooth even
gauge transformations that agree on the overlaps. They lift to the determinant-one group, so
the obstruction to extending them over the core $K$ lies in $H^4(K,K\cap O;\ZZ)$, which vanishes
because $K$ meets the free genus matching faces. A continuous even extension therefore
exists; it can be taken smooth relative to a neighbourhood of $K\cap O$ within its homotopy
class, hence even. The result is $V^{\rm ad}$, which is smooth.

\emph{Step 4: (i).} Given $k$, let $S(k)\ge S_0+3$ be such that the bounds of
Proposition~\ref{prop:tail-decay}(i) at orders at most $k+2$ are at most $e_T$ on
$\{r\ge S(k)-1\}$. On $(S(k)-1,\infty)\times\RR$
Lemma~\ref{lem:torus-slice} applies with $\|F_A\|_{C^{k+2}_A}\le e_T$, its torus slice is the
restriction of $V^{\rm ad}$ by fibrewise uniqueness, and~\eqref{eq:torus-slice} at order $k+1$,
with $H^2(T)\subset C^0(T)$, bounds $\sup_{T\times D_1(z)}|\nabla^k_\alpha a_V|$ by
$C_ke^{-\gamma_T(r(z)-S_0)}$ for $r(z)\ge S(k)$. On $\{S_V\le r\le S(k)\}$ the derivatives of $a_V$ of
order $k$ are bounded by Step 2(c) for $|\theta|$ large and by compactness otherwise. Since
$S_V\ge S_m+3$, Proposition~\ref{prop:tail-decay}(ii) applies to $a_V$.

\emph{Step 5: (ii).} For $\theta\ge\Theta_V$ split the slice into three parts. On $\{r\le S_V-1\}$,
$V^{\rm ad}=V^N$ and (a) applies. On the band $S_V-1\le r\le S_V$, $V^{\rm ad}-\alpha$ is obtained from
$V^N-\alpha$ by the gauge transformation $\exp(\chi\varphi)$, and it is $O(e^{-c\theta})$ in every $C^k$ by (a)
and (c). On $\{r\ge S_V\}$, $V^{\rm ad}-\pi^*B_{q_e}=a_V$; there, for $\theta$ beyond the threshold at
which $\|F_A\|_{C^{k+1}_A}\le e_T$ on $\{r>S_V-1,\ \theta'>\theta-1\}$, integrating the square
of~\eqref{eq:torus-slice} at order $k$ gives
\begin{align*}
\|a_V\|_{H^{k+1}([\theta,\theta+1]\times\{r\ge S_V\})}&\le C_k\sum_{i\le k+1}\|\nabla^i_AF_A\|_{L^2([\theta,\theta+1]\times\{r\ge S_V\})}\\
&\le C'_k\|F_A\|_{L^2([\theta-1,\theta+2]\times\{r\ge S_V-1\})}\le C''_kE(\theta-1)^{1/2},
\end{align*}
by small-energy regularity summed over unit squares and $|F_A|^2=2|\iota_{\partial_\theta}F_A|^2$ on the
collar, and this is $O(e^{-c\theta})$ by Proposition~\ref{prop:growing-ends}(i). On the finite
range of $\theta$ below that threshold the norm is finite by Step 4. Gauge transformations do
not change the section, and (a) bounds $u-q_e$ on the section interval. Since $\delta<c$, (ii)
follows.

\emph{Step 6: (iii)} is Step 2(e), which Step 3 does not change.
\end{proof}

\begin{lemma}[regularity on the torus domains]\label{lem:torus-regularity}
Let $U$ be a finite-energy solution of class $H^m_{\rm loc}$ on $W_{\pm,\mathsf R}$ with $\mathsf R<\infty$, or on
$W^J_\lambda$ with $\lambda\ge\lambda_*$, for the data of this section (the primary data, end
perturbations and almost complex structures of Section~\ref{subsec:folded}, with secondary
perturbations supported in the cores away from the torus region and the matching lines),
or for data that differ from these by compactly supported changes of the metric, cylinder-function terms,
four-dimensional holonomy terms and deformations of the almost complex structure, all
vanishing near the genus matching lines and, when present, the
flat torus face. Then there is an even gauge transformation of class
$H^{m+1}_{\rm loc}$ after which $U$ is smooth up to the genus matching faces, the flat face and
the Lagrangian arcs; its gauge-end limits are the fixed representatives; it converges on the
gauge and strip ends in every $C^k$ at every rate below the spectral gaps; and on each mixed
end its difference from its stationary limit, section included, has $C^k$ norm on
$[\theta,\theta+1]$ times the compact slice at most $C_ke^{-c\theta}$, for every $k$.
\end{lemma}

\begin{proof}
Steps 1, 2(e), 3 and 6 of the proof of Lemma~\ref{lem:adapted} apply on the domain at hand,
with three changes. At the flat face the local smooth gauges are given by
Theorem~\ref{thm:flat-boundary-regularity}(a); no holonomy term acts there. The mixed-end
slices are compact, so Proposition~\ref{prop:growing-ends}(ii) applies on the whole slice and
no torus slice is needed: the normalized temporal gauge serves on the whole end. The
additional terms are compactly supported and of the kinds handled in Step 1, the
cylinder-function terms living on a compact cylinder over a closed three-manifold
\cite[Sections~3.4--3.5]{KM-unknot}; they vanish near the genus matching lines, where Step 1
uses the seam regularity theorem. This covers the closed-cut data $\mathfrak d_Z(\ell)$ and the
interpolation of Lemma~\ref{lem:closed-cut-data}: the secondary data of $N_Z$ are supported
away from the matching lines (Section~\ref{subsec:mixed-iso} and
Remark~\ref{rem:clean-pairs}); those of $\overline{\mathsf E}$ lie in the piece $\mathrm{(II)}$, which does not
contain the matching lines (Lemma~\ref{lem:closed-cut-embedding}(ii)); the further
perturbation of Proposition~\ref{prop:closed-cut} is chosen away from the matching lines; and
the neck with $f_Z$ and the change of metric lie in a collar of $Z_{\rm cut}$, which meets the
genus parts in $\{|\operatorname{Re}+3|<1\}\cup
\{|\operatorname{Re}-2\lambda-3|<1\}$, at distance at least two
from the matching lines.
\end{proof}

The moduli spaces on $W_{\pm,\mathsf R}$ and $W^J_\lambda$ used in this section therefore consist of
classes of smooth solutions with these asymptotics, as on the domains of
\cite[Proposition~3.6]{DFL}.

\subsection{Energy, index and the loss formula}

\begin{theorem}\label{thm:torus-period}
\begin{enumerate}[label=(\roman*)]
\item For every $\mathsf R\in[0,\infty]$ and every configuration $c$ on $W_\mathsf R$
  with labels $\mathbf a=(a,x)$, the linearized operator on the extended weighted spaces is Fredholm and $\ind D_c=8\kappa(c)+\mathfrak c(x)-\mathfrak c(a)$, with
  $2\kappa(c)-2\kappa_0(\mathbf a)\in\ZZ$, the constants being those of
  \cite[Proposition~3.20]{DFL} for $X_D$. On each mixed end $e$ the stationary kernel is
  represented by the Coulomb fields $h_e(v)=(h_e^A(v),v)$, with $v$ in the tangent space at
  $q_e$ of the critical family of the end $e$ and $h_e^A(v)$ in the kernel of the stationary
  operator $\mathcal A_{q_e}$, with zero normal component; the operator realized with these
  representatives on the extended weighted spaces, the \emph{extended Coulomb realization},
  is onto at every regular solution.
\item For every $\lambda\ge\lambda_*$, the width threshold of
  Proposition~\ref{prop:stationary-gap}, and every configuration $c$ on $W^J_\lambda$ with labels
  $(\mathbf a_-,\mathbf a_+)$, the same holds with the sum of the two constants, and
  $\ind D_c=\ind D_{c_-}+\ind D_{c_+}$, where $c_\pm$ are obtained by first
  deforming $c$, with fixed labels and hence fixed $\kappa$ and index, to a
  configuration equal to $\alpha$ on a slice $T\times\{r_0\}\times\RR$ of the bridge,
  then cutting there and completing by $\alpha$.
\item The constants depend only on the labels, not on $\mathsf R$ or $\lambda$.
\end{enumerate}
\end{theorem}

\begin{proof}
(i) For $\mathsf R<\infty$, pull back by $\psi_\mathsf R$: the pulled-back configuration is a DFL
configuration on $X_D$ with the metric changed on $T\times[-\tfrac34,-\tfrac14]\times\RR_\theta$,
a region that runs into both mixed ends and changes their slices from $Y$ to
$Y_\mathsf R$. Along the family of metrics $g_s$, $s\in[0,\mathsf R]$, the operator on the extended weighted spaces stays Fredholm, because Proposition~\ref{prop:stationary-gap} gives a gap
uniform in the collar length; so the index is constant along the family and
equals the index on $X_D$, given by \cite[Proposition~3.20 and Lemma~3.15]{DFL}.
The topological energy is independent of the metric, and the boundary term in
Stokes' theorem at the flat face vanishes, the fibre restriction being flat and
$H^1(T;\ad)=0$.
For $\mathsf R=\infty$, deform $c$ to equal $\alpha$ far out on the tail and compare
$W_\infty$ with $W_\mathsf R$ by excision across a band of the collar on which all
operators are the operator of Proposition~\ref{prop:flat-fibre}. The excision is run on
the unextended $\delta$-weighted spaces, so that no extension vector crosses the cut, and
the same dimension correction is added on both sides: the operators
$D_{W_\infty}\oplus M_\mathsf R$ and $D_{W_\mathsf R}\oplus M_\infty$, where $M_{\mathsf R}$ is the invertible
model of Proposition~\ref{prop:flat-fibre} on $T\times(-\infty,\mathsf R]\times\RR$ with the flat
face and $M_\infty$ the model on $T\times\RR^2$, have equal index
once the band is longer than a constant depending only on the Fredholm
estimate of $D_{W_\infty}$; that estimate, for $D_{W_\infty}$ and its adjoint, is the
compact-remainder estimate of Step~2 of the proof of
Proposition~\ref{prop:torus-inverse} with an infinite collar. The tail contributes nothing to the topological
energy. The extended Coulomb realization is onto at a regular solution: its kernel and that
of the realization of \cite{DFL} both equal $H^1$ of the deformation complex, and their
indices agree by the finite-rank comparison of Remark~\ref{rem:torus-index-comparison}. (ii) The energy is local and the section splits, so
$\kappa(c)=\kappa(c_-)+\kappa(c_+)$; the index splits by the same excision with two
bands, the plane models having index zero, and the asymptotic values add. (iii)
holds because every constant was fixed on $X_D$.
\end{proof}

\begin{proposition}[loss formula]\label{prop:torus-loss}
Let $U_j$ be solutions of index $d$ on $W_{\mathsf R_j}$ with $\mathsf R_j\to\infty$, or on
$W^J_{\lambda_j}$ with $\lambda_j\to\infty$. After passing to a subsequence there are
solutions $V_\beta$ on the completed domains and broken trajectories $\Gamma_\ell$ on
the original ends with
$d-\sum_\beta\ind(V_\beta)-\sum_\ell\ind(\Gamma_\ell)=4n$, $n\in\ZZ_{\ge0}$. If all $\ind(V_\beta)\ge0$, all $\ind(\Gamma_\ell)\ge1$ and $d\le3$, then $n=0$ and no energy
is lost: after gauge transformations $U_j$ converges to the limit pieces in
$C^\infty$ on compact subsets of each component of the broken limit.
For every $\varsigma>0$ there are compact sets in all its components
such that, for large $j$, the energy outside their images in the
finite domain is less than $\varsigma$; trajectory components are
translated by their breaking parameters. If no trajectory breaks,
these are fixed compact sets in the completed core domains. In
particular this is the case when $d=0$ under the stated index assumptions.
The same statements hold for sequences of solutions on a fixed
completed domain $W_\infty$ or $W_{-,\infty}\sqcup W_{+,\infty}$.
\end{proposition}

\begin{proof}
Energy can be lost through interior instantons, through finite-energy
instantons on $T\times\RR^2$ or on $T$ times a half-plane with flat-fibre boundary,
through concentration on a boundary torus fibre, through seam and section
bubbles, and through escape along the mixed ends, possibly spread along the
lengthening slices. No classification is needed. Extract limits $V_\beta$ on
compact subsets of the completed domains, using fixed-geometry compactness and
Theorem~\ref{thm:flat-boundary-regularity} at the flat face, and chains $\Gamma_\ell$
at the original ends. Build a reference configuration $\mathcal R_j$ on the same finite geometry as $U_j$,
with the labels of $U_j$, as follows. First make each $V_\beta$ and each $\Gamma_\ell$
exactly stationary near the attaching regions on the original ends, by cutoffs
along paths of configurations with fixed limits; this changes neither the
topological energy nor the index. For a separate sequence, pass to the adapted
representative $V^{\rm ad}=\alpha+a_V$ on $r\ge S_V$ (Lemma~\ref{lem:adapted}) and
truncate its torus tail: put $\alpha+\chi_j(r)a_V$ on $T\times[S_V,\mathsf R_j]\times\RR$, with
$\chi_j=1$ for $r\le \mathsf R_j/2$, $\chi_j=0$ for $r\ge \mathsf R_j/2+1$ and $|\chi_j'|\le2$. The result
equals $\alpha$ near $r=\mathsf R_j$, so it satisfies the flat-fibre condition and is a configuration on $W_{\mathsf R_j}$. For a joined sequence truncate the two tails in the same
way at $r_\pm=\lambda_j/2$ and join the flat parts by the bridge identification, which
carries $\alpha$ to $\alpha$. Finally attach the chains $\Gamma_\ell$ at the original ends by
the cutoff construction of \cite[Propositions~4.65 and~4.69]{DFL}. The path
$s\mapsto\alpha+(1-s+s\chi_j)a_V$ consists of configurations on $W_\infty$ with
fixed labels, so by the Stokes argument in the proof of
Theorem~\ref{thm:torus-period} the topological energy is constant along it, and
the truncated configuration has the same energy on $W_{\mathsf R_j}$ because it is flat
beyond $r=\mathsf R_j/2+1$. Hence $\kappa(\mathcal R_j)=\sum\kappa(V_\beta)+\sum\kappa(\Gamma_\ell)$
exactly, and by Theorem~\ref{thm:torus-period}(i)--(iii), with additivity of the
constants $\mathfrak c$ along the original ends, $\ind D_{\mathcal R_j}=I=\sum\ind V_\beta+\sum\ind \Gamma_\ell$. By
Theorem~\ref{thm:torus-period}, $d-I=4n$ and $\kappa(U_j)-\kappa(\mathcal R_j)=
\tfrac12n$. The analytic energy is $8\pi^2\kappa$ plus a correction from the
perturbations that converges, by strong $L^p$ convergence of the holonomy
sources away from the concentration set; lower semicontinuity of the analytic
energy gives $n\ge0$. If $n=0$, the limit of the analytic energies equals the sum of the analytic
energies of the limit pieces. The energy on a compact set converges to that of the
limit on the same set, by the smooth convergence away from the concentration set,
and every concentration point carries a fixed quantum; so there is no concentration
and the energy outside the union of compact sets in the core and
translated trajectory components is small when those sets are large.
When there is no breaking, only the core compact sets are required.
For a fixed completed domain the same
argument applies with $\mathcal R_j$ replaced by a fixed truncation.
\end{proof}

\subsection{Families with stationary ends}

\begin{lemma}\label{lem:family-stationary}
Let $F_\lambda$, $\lambda\in[0,1]$, satisfy the hypotheses of
Theorem~\ref{thm:family} except compact support of $\partial_\lambda F$, and instead:
on each mixed end the critical set is connected, and it and a smooth family of
stationary references $R_q$ are independent of $\lambda$, with $F_\lambda(R_q)=0$; there $\partial_\lambda F_\lambda$
is independent of the end coordinate, and it vanishes on the gauge and ordinary
ends; and the stationary operators have a gap and the energy decay of
Proposition~\ref{prop:growing-ends}(i) holds along the family, with constants $e_0$ and $c$
independent of $\lambda$ and with $c$ larger than the weight $\delta$ of the weighted spaces. Then the
conclusions of Theorem~\ref{thm:family} hold.
\end{lemma}

\begin{proof}
On a mixed end, in the chart $R_q+z$ with $z$ in the weighted domain (weight
$e^{\delta\tau}$ in the end coordinate $\tau$), differentiating $F_\lambda(R_q)=0$ in
$\lambda$ gives $\partial_\lambda F_\lambda(R_q)=0$, since $R_q$ does not depend on $\lambda$. Hence
\begin{gather*}
\partial_\lambda F_\lambda(R_q+z)=\int_0^1D(\partial_\lambda F_\lambda)(R_q+sz)\,z\,ds,\\
\|e^{\delta\tau}\partial_\lambda F_\lambda(R_q+z)\|_{L^2}\le
\sup_{s\in[0,1]}\|D(\partial_\lambda F_\lambda)(R_q+sz)\|\,\|e^{\delta\tau}z\|_{L^2_1},
\end{gather*}
and the operator norm is bounded uniformly along the end because $\partial_\lambda F_\lambda$ is
independent of the end coordinate and $\|z\|_{C^0}$ is small there. So the parameter
derivative lies in the weighted target, $\lambda\mapsto F_\lambda$ is $C^1$ as a map into it, and
the augmented operator $(\dot\lambda,z)\mapsto\dot\lambda\,\partial_\lambda F_\lambda+D_\lambda z$ is a
rank-one extension of a Fredholm operator. Breaks occur only at the gauge and
ordinary ends. A finite-energy trajectory on a mixed end runs from the connected
critical family to itself. Its energy is a period of the action; if it is zero the
trajectory is constant, and otherwise $\kappa>0$, so by the lattice
condition~\ref{H:lattice} its index is at least four. Such a trajectory cannot split off a
family of fibre index at most one whose other components, being regular, have index at
least $-1$. The preglued configuration agrees with the
core on the mixed end, so the parameter derivative vanishes on the joined end and neck.
The solutions are smooth up to gauge (\cite[Proposition~3.6]{DFL},
Lemma~\ref{lem:torus-regularity}), and the argument of the proof of
Proposition~\ref{prop:growing-ends}(ii), applied to each solution at its own parameter, shows
that it decays at the rate $c>\delta$ in every $H^k$ of its mixed-end slices, so it lies in the
weighted domain. Family compactness uses the uniform energy decay and the $\lambda$-independent
constants of Theorem~\ref{thm:torus-period}.
\end{proof}

\begin{proposition}\label{prop:collar-length}
Changing the length of a separate torus collar from $\ell_0$ to $\ell_1$, with the
outer flat-fibre boundary present throughout, changes the regular separate map
by a chain homotopy.
\end{proposition}

\begin{proof}
Apply Lemma~\ref{lem:family-stationary}; its hypothesis that the critical set is
connected holds, the critical sets being the perturbed $L(U_i)$ and $R_i$,
diffeomorphic to $\SU(2)^{h-1}$ and $\SU(2)^{g-1}$. The stationary references are
perturbed-flat and equal to $\alpha$ on the collar; the perturbed flatness equation
involves the metric only on the support of the perturbation, so the references
are independent of the collar length; the pulled-back reference satisfies the
equation on the collar because it is flat there; and the variation of the metric
is supported in a region of the collar independent of the end coordinate.
\end{proof}

\begin{proposition}[width homotopy]\label{prop:width}
For $\lambda_*\le a\le b<\infty$, with $\lambda_*$ the width threshold of
Proposition~\ref{prop:stationary-gap}, $N^J_b\simeq N^J_a$.
\end{proposition}

\begin{proof}
Change the bridge length by a family of metrics supported on the bridge. The
exact stationary references $R_{\lambda,q}$ on the joined mixed ends, flat and equal
to $\alpha$ on the bridge, depend on the width only through the length of the flat
region; differentiating $F_\lambda(R_{\lambda,q})=0$ shows that the nondecaying part of the parameter derivative vanishes. After pulling back by width-changing diffeomorphisms
that are the identity outside the bridge, the hypotheses of
Lemma~\ref{lem:family-stationary} hold; on the compact range $[a,b]$ each slice has a
gap, depending continuously on the width, so the gap is uniform along the family, and so is
the energy decay, by the argument of the proof of Proposition~\ref{prop:growing-ends}(i).
\end{proof}

\subsection{The closed cut}\label{subsec:closed-cut}

We describe the gauge part of $W^J_\lambda$ in the coordinates of the special
quintuple \cite[Section~3.1]{DFL}. For each side, $\mathbb U_-\subset\CC$ is the region of
\cite[Figure~1]{DFL}: it contains the input arm $(-\infty,-3]\times[-1,1]$, the arms
$[-2,0]\times[2,\infty)$ and $[-2,0]\times(-\infty,-2]$ and the square $[-1,0]\times[-2,2]$,
its boundary consists of the matching line $\mathbb U_\partial=i\RR$ and two curves $\eta_-$,
$\eta_-'$, and the gauge domain is
$\RR\times Y_0\cup \mathbb U_-\times\Sigma\cup\RR\times(-Y'_0)$, where $\RR\times\partial Y_0$ is glued to
$\eta_-\times\Sigma$ through the parameter $t$ of $\eta_-$ and the collar coordinate of a
neighbourhood $\RR\times(\tfrac12,1]$ of $\eta_-$, and similarly for $\eta'_-$. Along the
input arm, $t=\operatorname{Re}$, $\eta_-$ is the line $\operatorname{Im}=1$ and
$\eta'_-$ the line $\operatorname{Im}=-1$. For the two sides of $W^J_\lambda$ write
$Y_\pm,Y'_\pm$ for the compression bodies ($Y_-=-U_0$, $Y'_-=-U_1$, $Y_+=V_0$,
$Y'_+=V_1$), $\Sigma_-=\Sigma_h$, $\Sigma_+=\Sigma_g$, and $\mathbb U_-^{(\pm)}$, $\eta^{(\pm)}_-$,
$\eta'^{(\pm)}_-$ for the planar data. The torus part of the gauge domain of
$W^J_\lambda$ is $T\times\widetilde{\mathbb U}_\lambda$ with
\[
\widetilde{\mathbb U}_\lambda=\mathbb U_-^{(-)}\cup\bigl([0,2\lambda]\times\RR\bigr)\cup\varrho_\lambda\bigl(\mathbb U_-^{(+)}\bigr),
\qquad \varrho_\lambda(z)=2\lambda-\bar z ,
\]
where the bridge coordinate of Definition~\ref{def:torus-domains} is
$\operatorname{Re}z-\lambda$, and the torus factors of the two sides are identified by
the marked bundle isometry. The genus parts are $\mathbb U_-^{(-)}\times\Sigma_h$ and
$\mathbb U_-^{(+)}\times\Sigma_g$.

\emph{The cut.} Let $\gamma_0\subset\widetilde{\mathbb U}_\lambda$ be the arc consisting of the segment
$\{-3\}\times[\tfrac12,1]$, the segment $[-3,2\lambda+3]\times\{\tfrac12\}$ and the segment
$\{2\lambda+3\}\times[\tfrac12,1]$, with its two corners rounded, and $\gamma_1$ its mirror
image in the real axis. The endpoints of $\gamma_0$ are the points of $\eta^{(-)}_-$ and
$\varrho_\lambda(\eta^{(+)}_-)$ with parameter $t=-3$, and the arc leaves them along the
collar coordinate; similarly for $\gamma_1$. Let $\mathsf a_\pm=\{-3\}\times[-1,1]$ in
the input arm of $\mathbb U_-^{(\pm)}$. Put
\[
Z_{\rm cut}=\bigcup_\pm\bigl(\{-3\}\times(Y_\pm)_0\ \cup\ \{-3\}\times(-Y'_\pm)_0
\ \cup\ \mathsf a_\pm\times\Sigma_\pm\bigr)\ \cup\ (\gamma_0\sqcup\gamma_1)\times T .
\]
Let $\widetilde{\mathbb U}_{\rm up}$, $\widetilde{\mathbb U}_{\rm mid}$, $\widetilde{\mathbb U}_{\rm down}$ be the
closures of the components of $\widetilde{\mathbb U}_\lambda\setminus(\gamma_0\cup\gamma_1)$ containing
the upper arm, the two input arms, and the lower arm.

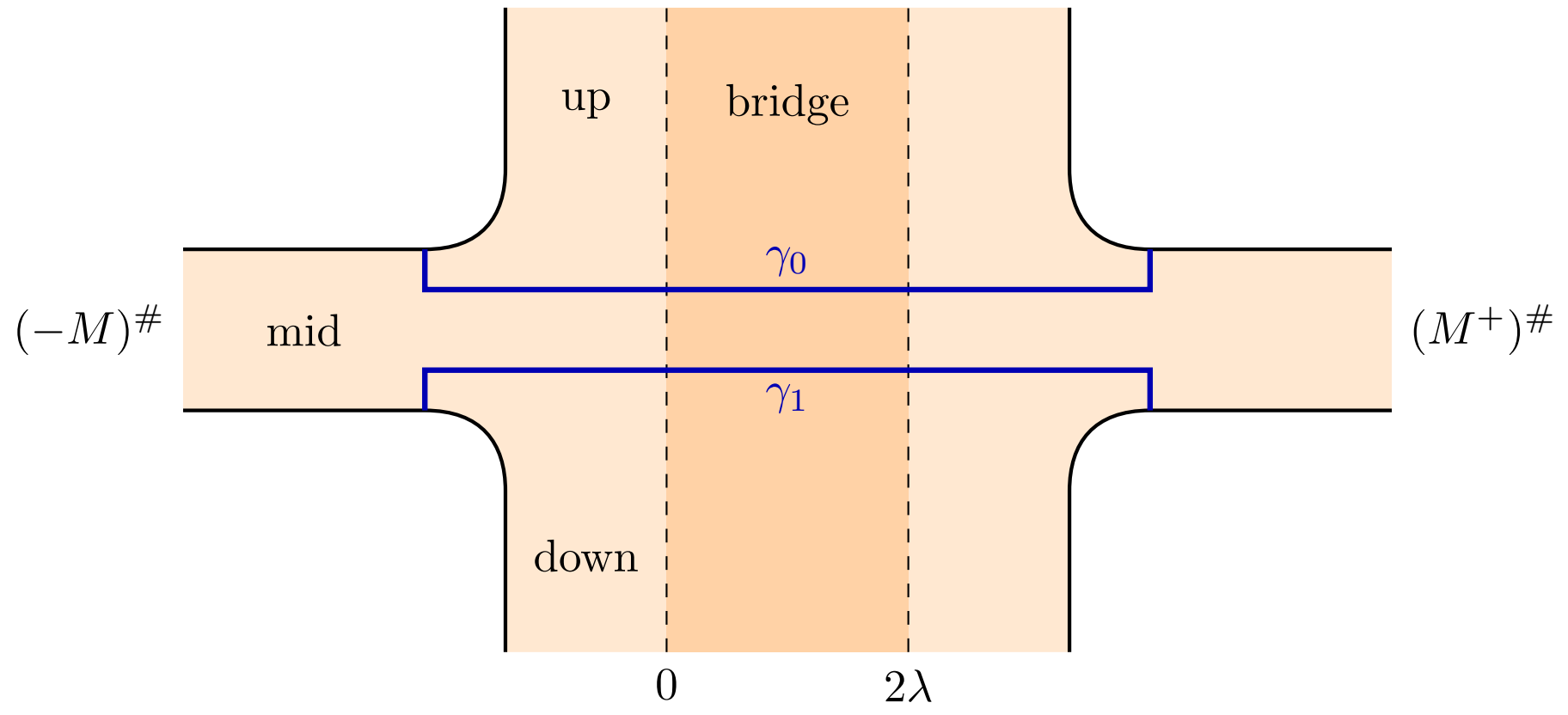
\begin{figure}[ht]
\centering
\begin{tikzpicture}[scale=0.62,font=\small]
\fill[orange!18] (-6,-1) -- (-3,-1) .. controls (-2.3,-1) and (-2,-1.4) .. (-2,-2) -- (-2,-4) -- (5,-4) -- (5,-2) .. controls (5,-1.4) and (5.3,-1) .. (6,-1) -- (9,-1) -- (9,1) -- (6,1) .. controls (5.3,1) and (5,1.4) .. (5,2) -- (5,4) -- (-2,4) -- (-2,2) .. controls (-2,1.4) and (-2.3,1) .. (-3,1) -- (-6,1) -- cycle;
\fill[orange!35] (0,-4) rectangle (3,4);
\draw[thick] (-6,1) -- (-3,1) .. controls (-2.3,1) and (-2,1.4) .. (-2,2) -- (-2,4);
\draw[thick] (-6,-1) -- (-3,-1) .. controls (-2.3,-1) and (-2,-1.4) .. (-2,-2) -- (-2,-4);
\draw[thick] (9,1) -- (6,1) .. controls (5.3,1) and (5,1.4) .. (5,2) -- (5,4);
\draw[thick] (9,-1) -- (6,-1) .. controls (5.3,-1) and (5,-1.4) .. (5,-2) -- (5,-4);
\draw[dashed] (0,-4) -- (0,4); \draw[dashed] (3,-4) -- (3,4);
\draw[very thick,blue!70!black] (-3,1) -- (-3,0.5) -- (6,0.5) -- (6,1);
\draw[very thick,blue!70!black] (-3,-1) -- (-3,-0.5) -- (6,-0.5) -- (6,-1);
\node[blue!70!black] at (1.5,0.85) {$\gamma_0$};
\node[blue!70!black] at (1.5,-0.85) {$\gamma_1$};
\node at (1.5,2.8) {bridge};
\node at (-1,2.8) {up};
\node at (-4.5,0) {mid};
\node at (-1,-2.8) {down};
\node[below] at (0,-4) {$0$}; \node[below] at (3,-4) {$2\lambda$};
\node[left] at (-6,0) {$(-M)^\#$};
\node[right] at (9,0) {$(M^+)^\#$};
\end{tikzpicture}
\caption{The torus base $\widetilde{\mathbb U}_\lambda$ of the joined domain $W^J_\lambda$
(Section~\ref{subsec:closed-cut}): the planar region $\mathbb U_-$ of the lower side on the
left, its reflection for the upper side on the right, and the bridge between the
two matching lines (dashed). The arcs $\gamma_0$ and $\gamma_1$ run from the boundary
curves of one input arm to those of the other. The closed cut $Z_{\rm cut}$ consists of
$(\gamma_0\sqcup\gamma_1)\times T$ together with the slices of the input arms at $t=-3$; it
separates the part containing the matching lines (regions up and down, a domain
for the pair $(A_0,A_1)$) from the part containing the two gauge inputs and the
middle of the bridge (region mid, Kronheimer and Mrowka's reverse excision
cobordism).}
\label{fig:closed-cut}
\end{figure}

\begin{lemma}[the closed cut in the joined domain]\label{lem:closed-cut-embedding}
\begin{enumerate}[label=(\roman*)]
\item $Z_{\rm cut}$ is a closed connected hypersurface in the gauge part of
  $W^J_\lambda$, disjoint from the matching lines, the section and the supports of the
  secondary perturbations, and there is an orientation-preserving diffeomorphism
  $Z_{\rm cut}\cong Z=D_0\cup(-D_1)$, covered by a bundle isomorphism, which is the
  identity on the pieces $(Y_\pm)_0$, $(-Y'_\pm)_0$ and the genus collars, and carries
  $\gamma_i\times T$ onto the torus neck of $D_i$.
\item $Z_{\rm cut}$ separates the gauge part into two pieces. The piece $\mathrm{(I)}$
  containing the matching lines, together with the section $\mathbb U_+$ and a cylindrical
  end $(-\infty,0]\times Z$ attached along $Z_{\rm cut}$, is a domain of DFL type for the
  pair $(A_0,A_1)$ with compression bodies $D_0$, $D_1$, with the end types of the
  special quintuple of $(D_0,D_1)$: its genus parts are
  $(\mathbb U^{(\pm)}_-\cap\{\operatorname{Re}\ge-3\})\times\Sigma_\pm$, its compression pieces
  $[-3,\infty)\times(Y_\pm)_0$ and $[-3,\infty)\times(-Y'_\pm)_0$ together with
  $(\widetilde{\mathbb U}_{\rm up}\sqcup\widetilde{\mathbb U}_{\rm down})\times T$, its gauge input slice is $Z_{\rm cut}$, and
  its mixed-end slices are $D_{0,\lambda}$ and $D_{1,\lambda}$, which are $D_0$ and $D_1$ with
  torus necks of length $2\lambda+4$, the slices including their arm collars $[-2,0]$ on
  both sides of the interval $[-\lambda,\lambda]$ of Definition~\ref{def:torus-domains}(c).
\item The other piece $\mathrm{(II)}$, with a cylindrical end attached along $Z_{\rm cut}$, is
  the reverse excision cobordism $\overline W\colon Z'\to Z$ of
  \cite[Section~5.3]{KM-unknot}, with $Z'=(-M)^\#\sqcup(M^+)^\#$ and cylindrical input
  ends. Here the two tori of Theorem~5.6 of \cite{KM-unknot} are
  $T_a=T\times\{x_a\}$ and $T_b=T\times\{x_b\}$ for points $x_a\in\gamma_0$, $x_b\in\gamma_1$,
  identified through the product structure of $T\times\widetilde{\mathbb U}_\lambda$; the product
  part $T\times U$ of \cite[Figure~4]{KM-unknot} is $T\times\widetilde{\mathbb U}_{\rm mid}$, truncated
  on the input arms; and the bundle is odd on $T_a$ and $T_b$.
\item Let $d_1,d_2,d_3,d_4\in H^1(W^J_\lambda;\FF)$ be dual to the hypersurfaces
  $\eta^{(-)}_-\times T$, $\eta'^{(-)}_-\times T$, $\varrho_\lambda(\eta^{(+)}_-)\times T$,
  $\varrho_\lambda(\eta'^{(+)}_-)\times T$. Then $d_1=d_2$, $d_3=d_4$ and $\Phi_J=\langle d_1,d_3\rangle$.
  Restricted to $\mathrm{(II)}$, $\Phi_J$ is the subgroup $\phi_{\overline W}$ of
  \cite[Section~5.3]{KM-unknot} generated by the duals of the four copies of an
  interval times $T$ over the dotted edges of $U$, and its image at the ends is
  $\Phi_{\overline{\mathsf E}}=\{(s,s_-,s_+):s=s_-+s_+\}$. Restricted to $\mathrm{(I)}$, $\Phi_J$ is the
  lifting group of the domain of the pair $(A_0,A_1)$, generated by the class that
  restricts to $s_Z$ on $Z$; and $\Phi_J$ is the fibre product of the two restrictions
  over $\phi_Z=\langle s_Z\rangle$.
\item After a change of the metric of $\widetilde{\mathbb U}_\lambda$ supported in a compact
  neighbourhood of $\gamma_0\cup\gamma_1$ inside the regions $\mathbb U^{(\pm)}$, disjoint from the bridge
  ends, the arms and the matching lines, $Z_{\rm cut}$ has a neighbourhood
  $(-1,1)\times Z_{\rm cut}$ with the product metric.
\end{enumerate}
\end{lemma}

\begin{proof}
(i) Near the input arms the gauge domain is the product of $t=\operatorname{Re}$
with the slice $\{t\}\times(Y_\pm)_0\cup\{t\}\times[-1,1]\times\Sigma_\pm\cup\{t\}\times(-Y'_\pm)_0$ in
the genus part and with the collar segments $\{t\}\times(\tfrac12,1]$ in the torus part.
Reading $Z_{\rm cut}$ across: the arc $\gamma_0$ joins the torus boundary of $(Y_-)_0$ at
$\eta^{(-)}_-(-3)$ to that of $(Y_+)_0$ at $\varrho_\lambda(\eta^{(+)}_-)(-3)$, so
$(Y_-)_0\cup\gamma_0\times T\cup(Y_+)_0$ is $Y_-\cup_T Y_+=(-U_0)\cup_TV_0=D_0$ with a torus neck of
length equal to the length of $\gamma_0$; likewise
$(-Y'_-)_0\cup\gamma_1\times T\cup(-Y'_+)_0=-D_1$; and $\mathsf a_-\times\Sigma_h$,
$\mathsf a_+\times\Sigma_g$ glue $D_0$ to $-D_1$ along $\Sigma_h$ and $\Sigma_g$. So
$Z_{\rm cut}\cong D_0\cup(-D_1)=Z$. The bundle is the given one on each piece and the
pulled-back odd bundle on $\gamma_i\times T$, so the identification is a bundle
isomorphism. The primary perturbations of the pieces $\RR\times(Y_\pm)_0$ and
$\RR\times(-Y'_\pm)_0$ are independent of $t$ and restrict on $Z_{\rm cut}$ to the primary
perturbations of $D_0$ and $D_1$. Connectedness is clear. The matching lines lie
in $\operatorname{Re}\in\{0,2\lambda\}$, the section lies in $\{t>-3\}$, and the
secondary supports of the data $\mathfrak d_Z(\ell)$ avoid $Z_{\rm cut}$ by their
choice; the closed-cut comparison (Lemma~\ref{lem:closed-cut-data}) does not use this
for the joined data $\mathfrak d_J$. The end-perturbation terms $\chi_\pm f_\pm$ on
the two gauge inputs are chosen with $\operatorname{supp}\chi_\pm\subset\{t\le-4\}$, so
that the transition lies in the piece $\mathrm{(II)}$, whose interior data are
arbitrary, and the neck along $Z_{\rm cut}$ carries only the end perturbation
$f_Z$.

(ii) The component of the complement containing the matching lines consists of
the pieces listed. The region $\widetilde{\mathbb U}_{\rm up}$ is bounded by $\gamma_0$, by the
parts $t\ge-3$ of $\eta^{(-)}_-$ and $\varrho_\lambda(\eta^{(+)}_-)$, and by the end of the
upper arm, whose cross-section is $[-2,2\lambda+2]$; it is diffeomorphic to
$[-3,\infty)\times[0,1]$ by a diffeomorphism that is an isometry onto a product for large
$t$ and that matches the collar coordinates at the two $\eta$-curves. Hence
$[-3,\infty)\times(Y_-)_0\cup\widetilde{\mathbb U}_{\rm up}\times T\cup[-3,\infty)\times(Y_+)_0$ is
diffeomorphic to $[-3,\infty)\times(D_0)_0$, and similarly for $D_1$. With the genus
parts this is the special quintuple of $(D_0,D_1)$ restricted to $t\ge-3$, except
that the metric of the compression pieces is induced from $\widetilde{\mathbb U}_\lambda$: on
$\{t\le t_1\}$ for a fixed $t_1$ it differs from a product in a compact set, and for
$t\ge t_1$ it is the product with the slices $D_{i,\lambda}$, whose torus necks have length
$2\lambda+4$ rather than the length of $\gamma_i$. Its section is $\mathbb U_+$ with
the Lagrangians $A_0,A_1$, the perturbed Lagrangians of $D_0$, $D_1$ being
$\widetilde P_i^-\times R_i=A_i$ (Section~\ref{subsec:folded}).

(iii) The other component is
\begin{multline*}
(-\infty,-3]\times\bigl((Y_-)_0\sqcup(-Y'_-)_0\sqcup(Y_+)_0\sqcup(-Y'_+)_0\bigr)\\
\cup\ \bigl(\{\operatorname{Re}\le-3\}\text{-parts of the input arms}\bigr)\times\Sigma_\pm
\ \cup\ \widetilde{\mathbb U}_{\rm mid}\times T .
\end{multline*}
The first two parts form $(-\infty,-3]\times Z'_{\rm cut}$, where
$Z'_{\rm cut}=\bigl(Y_-\cup_{\Sigma_h}(-Y'_-)\bigr)\sqcup\bigl(Y_+\cup_{\Sigma_g}(-Y'_+)\bigr)$ is
$Z'=(-M)^\#\sqcup(M^+)^\#$ cut along its two tori, with four boundary tori. Truncate
the input arms at $\operatorname{Re}=-T_1$ and $\operatorname{Re}=2\lambda+T_1$. Then
$U:=\widetilde{\mathbb U}_{\rm mid}\cap\{-T_1\le\operatorname{Re}\le2\lambda+T_1\}$ is a disk with eight
corners, whose edges in cyclic order are: $\gamma_0$; the part $t\in[-T_1,-3]$ of
$\varrho_\lambda(\eta^{(+)}_-)$; the end of the right input arm; the corresponding part of
$\varrho_\lambda(\eta'^{(+)}_-)$; $\gamma_1$; the parts of $\eta'^{(-)}_-$ and, after the end of
the left input arm, of $\eta^{(-)}_-$. The four $\eta$-edges (the dotted edges of
\cite[Figure~4]{KM-unknot}) times $T$ are glued to $[-T_1,-3]\times\partial Z'_{\rm cut}$;
the two arm ends times $T$ are two-sided collars of the tori of $Z'$; and
$\gamma_0\times T$, $\gamma_1\times T$ are two-sided collars of $T_a$, $T_b$ in $Z$. This is
the description of the standard cobordism in the proof of
\cite[Theorem~5.6]{KM-unknot}: the complement of $T\times U$ is a product of an interval
with $Z'_{\rm cut}=Z_{\rm cut}\setminus(\text{collars of }T_a\cup T_b)$, and $T\times U$ is a
product. The cobordism runs from $Z'$ at $\operatorname{Re}=-T_1$ to $Z$, so it is the
reverse cobordism $\overline W$. The bundle is the product of the odd bundle with $U$
on $T\times U$ and the product bundle of $Z'_{\rm cut}$ elsewhere; the identification of
$T_a$ with $T_b$ required in \cite[Section~5.3]{KM-unknot} is the product structure.

(iv) The hypersurfaces $\eta^{(-)}_-\times T$ and $\eta'^{(-)}_-\times T$, together with the
matching face $U^{(-)}_\partial\times\Sigma_h$ and parts of the ends, bound the region
$\RR\times(Y_-)_0\cup U^{(-)}_-\times\Sigma_h\cup\RR\times(-Y'_-)_0$; so their duals agree in
$H^1(W^J_\lambda;\FF)$, which is Lefschetz dual to $H_3$ relative to the boundary and
the ends. The same holds on the $+$ side. The class $d_1$ restricts to the dual
of the torus in $(-M)^\#$ and to zero in $(M^+)^\#$, and $d_3$ the other way round, so
$\langle d_1,d_3\rangle$ maps isomorphically onto $\Phi_-\oplus\Phi_+$, which is $\Phi_J$. On
$Z$, both $d_1$ and $d_3$ restrict to the dual of $T_a$, which equals that of $T_b$
since the two tori cobound $Z'_{\rm cut}$-pieces; this is $s_Z$. The restrictions to
$\mathrm{(II)}$ are the duals of the dotted-edge hypersurfaces of (iii), which generate
$\phi_{\overline W}$, with image $\{s=s_-+s_+\}$ at the ends. In $\mathrm{(I)}$ the two
hypersurfaces $\eta^{(-)}_-|_{t\ge-3}\times T$ and $\varrho_\lambda(\eta^{(+)}_-)|_{t\ge-3}\times T$
cobound $\widetilde{\mathbb U}_{\rm up}\times T$ relative to the boundary, so $d_1$ and $d_3$ restrict
to the same class, which vanishes on the matching surfaces and restricts to $s_Z$
on $Z$. Since $Z$ is connected, restriction
$H^1(W^J_\lambda)\to H^1(\mathrm{I})\oplus H^1(\mathrm{II})$ is injective with image the pairs agreeing on
$Z$.

(v) Choose coordinates $(u,s)\in(-\tfrac14,\tfrac14)\times\gamma_i$ near $\gamma_i$, so that the
collars of $\gamma_0$ and $\gamma_1$ are disjoint and lie in the domain, equal to the
Euclidean coordinates $(t,\operatorname{Im})$ near the endpoints and to
$(\operatorname{Im}\mp\tfrac12,\operatorname{Re})$ along the horizontal segment, and
replace the metric by $du^2+ds^2$ near the rounded corners and along the horizontal
segments, where \cite{DFL} do not fix the planar metric to be Euclidean,
interpolating in a compact region. Near the endpoints this is the product structure of the input
arms, so $Z_{\rm cut}$ acquires a product collar.
\end{proof}

\begin{remark}\label{rem:excision-models}
Reversing the direction of the interval factors, the same description gives the
forward cobordism $\mathcal E\colon Z\to Z'$ as $T\times U\cup[0,1]\times Z'_{\rm cut}$. Since
$Z'_{\rm cut}$ is the union of $(T\times I_-)\#(-M)$ and $(T\times I_+)\#M^+$, this is $T\times P$
with $P=U\cup([0,1]\times I_-)\cup([0,1]\times I_+)$, a pair of pants, and with the products of
the punctured summands inserted along two arcs from the input circle to the two
output circles; this is the model used in Lemma~\ref{lem:Bk-topology}(ii).
\end{remark}

\emph{Data at the fixed width.} Let $\mathfrak d_J$ be the regular joined data at width
$\lambda_0$: the primary data of Section~\ref{subsec:folded}, the end
perturbations $f_\pm$ on the two gauge inputs, a split invariant secondary
perturbation $\eta_-+\eta_+$ with each term supported in the core of its own side, and
the split structure $J_-\times J_+$ on $\mathbb U_+$. For $\ell\ge0$ let $\mathfrak d_Z(\ell)$ be
the \emph{closed-cut data}: the metric of Lemma~\ref{lem:closed-cut-embedding}(v)
with $[-\ell,\ell]\times Z$ inserted along $Z_{\rm cut}$; the end perturbation $f_Z$ of the data of $N_Z$ at width $\lambda_0$, multiplied by a cutoff of the neck coordinate
supported in $[-\ell,\ell]$ and equal to one on $[-\ell+1,\ell-1]$; regular invariant
secondary data on $\mathrm{(I)}$, chosen directly on $\mathrm{(I)}$, and those of
$\overline{\mathsf E}$ on $\mathrm{(II)}$, supported away from the neck; the same primary data, end perturbations $f_\pm$ and split structure
$J_-\times J_+$. Here $N_Z$ is the mixed map of Remark~\ref{rem:clean-pairs} for
$(D_0,D_1)$ with the metrics induced on $Z_{\rm cut}$ (the \emph{data at width $\lambda_0$}), and $\overline{\mathsf E}$ is computed at the same end data $f_\pm$ and
$f_Z$.

\begin{lemma}[data change at fixed width]\label{lem:closed-cut-data}
For every $\ell\ge0$ for which $\mathfrak d_Z(\ell)$ is regular in index at most one, the
joined counts at $\mathfrak d_J$ and at $\mathfrak d_Z(\ell)$ are chain homotopic over
$\FF$.
\end{lemma}

\begin{proof}
Let $\psi_\ell$ be a diffeomorphism of $W^J_{\lambda_0}$ onto the domain with the
inserted neck, equal to the identity outside a compact neighbourhood of
$Z_{\rm cut}$ and stretching a collar of $Z_{\rm cut}$. Pull back $\mathfrak d_Z(\ell)$ by
$\psi_\ell$. The two sets of data on $W^J_{\lambda_0}$ then differ only in a compact set $K$
of the gauge domain disjoint from all ends: in the metric (near $Z_{\rm cut}$), in
the term $f_Z$ (supported in the neck) and in the secondary perturbations. On
the gauge inputs, the strip output and the mixed ends they coincide, with the
same critical data. Interpolate linearly in the metrics and in the coefficients
of the holonomy terms. The family is invariant under $\Phi_J$: the metrics and
secondary terms are invariant by choice, and $f_Z$ is invariant under $\phi_Z$ and
under the determinant-one group, while $\Phi_J$ restricts on $Z$ into $\phi_Z$
(Lemma~\ref{lem:closed-cut-embedding}(iv)). The parameter derivative is supported
in $K$. Theorem~\ref{thm:family} applies and gives the chain homotopy.
\end{proof}

\begin{proposition}\label{prop:closed-cut}
For every $\lambda_0\ge\lambda_*$, $\mathcal KN^J_{\lambda_0}\simeq\mathcal KN_Z\overline{\mathsf E}$ over $\FF$,
where $N^J_{\lambda_0}$ is the joined map at the data $\mathfrak d_J$ and $N_Z$ is taken with the data at width $\lambda_0$.
\end{proposition}

\begin{proof}
Fix $\lambda_0\ge\lambda_*$. By Lemma~\ref{lem:closed-cut-embedding}, $W^J_{\lambda_0}$ with the data
$\mathfrak d_Z(\ell)$ is the union, along a neck $[-\ell,\ell]\times Z$, of the DFL domain $\mathrm{(I)}$
for the pair $(A_0,A_1)$ and the gauge cobordism $\mathrm{(II)}=\overline W$, and its lifting group is the fibre product of theirs over $\phi_Z$. The cut $Z$ is connected, its
critical points for the perturbed functional with $f_Z$ are nondegenerate and
irreducible, since the odd torus kills the stabilizers, the connecting $H^0$ term
vanishes, and the framed groups glue by the fibre product
\cite[Proposition~5.4]{KM-unknot}. Hypotheses \ref{H:lattice}--\ref{H:framing} hold,
\ref{H:lattice} by Theorem~\ref{thm:torus-period}(ii). By Theorem~\ref{thm:neck},
for $\ell\ge\ell_0$ the data $\mathfrak d_Z(\ell)$ are regular in index at most zero, and a
further invariant perturbation, at this fixed finite $\ell$,
makes the index-one spaces regular. Its output region includes the
compact neck where the fixed term $\chi f_Z$ has output; it remains
away from the matching lines. Choose it small enough that openness
of regularity and compactness preserve the already regular
index-zero count. Its support may meet that closed neck; the joined count is then the
composite of the count of $\mathrm{(I)}$ and that of $\mathrm{(II)}$. By
Lemma~\ref{lem:vertex-identification}, whose mixed-end clause covers the passage from
the slices $D_{i,\lambda_0}$ of $\mathrm{(I)}$ to the ends of the quintuple of $(D_0,D_1)$ with the
metrics of $Z_{\rm cut}$, the count of $\mathrm{(I)}$ is chain homotopic to $N_Z$. The count of
$\mathrm{(II)}$ is chain homotopic to $\overline{\mathsf E}$, cobordism maps being independent of the
interior data \cite[Section~5.1]{KM-unknot}. With Lemma~\ref{lem:closed-cut-data}
this gives $N^J_{\lambda_0}\simeq N_Z\overline{\mathsf E}$, and composing with $\mathcal K$ gives the claim;
Lemma~\ref{lem:closed-cut-data} is applied to $\mathfrak d_Z(\ell)$ with the further
perturbation above, which is regular in index at most one, and its proof is unchanged.
\end{proof}

\subsection{Torus completion}\label{subsec:torus-completion}

The torus completion follows the steps of the proof of Theorem~\ref{thm:neck}: a uniform
inverse at the preglued configuration, the contraction in a chart, a scalar Neumann inverse for
the gauge slice, and the entry of nearby solutions into the chart. What is new are the inputs: the flat-fibre chart at the face
(Lemma~\ref{lem:face-chart}); the uniform Poincar\'e inequality on the torus domains
(Lemma~\ref{lem:torus-scalar}); the uniform inverse across the growing torus collar
(Proposition~\ref{prop:torus-inverse}), which rests on the invertibility of the odd flat torus
operator (Proposition~\ref{prop:flat-fibre}) and on the gap of
Proposition~\ref{prop:stationary-gap} holding uniformly in the length; and the entry of nearby
solutions into the chart (Proposition~\ref{prop:torus-entry}), through the torus slice of
Lemma~\ref{lem:torus-slice}.

\emph{Conventions.} Throughout this subsection $\mathsf L$ stands for $\mathsf R$ on the separate
collared domains $W_\mathsf L=W_{\pm,\mathsf R}$ and for $\lambda$ on the joined domains $W_\mathsf L=W^J_\lambda$,
and $W_\infty$ stands for $W_{\pm,\infty}$ or $W_{-,\infty}\sqcup W_{+,\infty}$. The index $e$ runs
over the mixed ends; the critical manifold of the end $e$ is $\mathcal C_e$ (a critical
manifold of one side on a separate domain, a product $P_i\times Q_i$ on a joined end) and its slice is $Y_{e,\mathsf L}$ ($Y_\mathsf R$,
$Y'_\mathsf R$ or $D_{i,\lambda}$). The level $m\ge4$ is the one fixed before
Lemma~\ref{lem:adapted}. The weight $\tau$ of
Definition~\ref{def:torus-domains} is chosen with $\tau=\tau_0(\theta)$ on $N_X$ and on the
torus region, $|\tau_0'|\le1$, $|\tau_0''|\le1$, $\tau_0(\theta)=|\theta|$ for $|\theta|\ge2$ and
$\tau\ge0$; it has zero normal derivative on $\partial W_\mathsf L$ and bounded derivatives of
every order used. Fix $\delta$ with
\[
0<\delta\le\delta_T,\qquad\delta<\gamma/2,\qquad\delta<1,\qquad\delta<\mu_o\ \text{for every gauge and strip end }o,
\]
where $\delta_T$ is as in Proposition~\ref{prop:tail-decay}, $\gamma$ is the gap of
Proposition~\ref{prop:stationary-gap}, and $\mu_o>0$ is the spectral gap of the asymptotic
operator at $o$: the extended Hessian, gauge-fixing component included, at the
nondegenerate irreducible critical point of a gauge end, and $J\partial_t$ with the Lagrangian
boundary conditions at the transverse intersection point of a strip end. We also take $\delta$
below the weights of the original mixed ends in \cite[Definition~3.9]{DFL}. Thus $\delta$ lies
below the spectral gaps at all ends: at the mixed ends by $\delta<\gamma/2$, and at the gauge and
strip ends by $\delta<\mu_o$. In particular the completed solutions converge on the gauge and
strip ends at a rate larger than $\delta$ (Lemma~\ref{lem:adapted}(iii)), and model (O) in the
proof of Proposition~\ref{prop:torus-inverse} applies. We further
decrease $\delta$ so that $\delta C_P\le\frac12$, where $C_P$ is the constant of
Lemma~\ref{lem:torus-scalar}, which does not depend on $\delta$. All these conditions are
independent of $\mathsf L$.
Constants $C$, $C_{\rm inv}$, $C_{\rm quad}$ and $r_1$ depend on
the fixed data and on the finite sets of completed solutions, not on $\mathsf L$.

\begin{lemma}[regular completed pieces]\label{lem:completed-regular}
For generic secondary perturbations in the invariant class of
Theorem~\ref{thm:regularity}, one can choose a common compact
output and dependence region in the original core, away from the
inserted torus region and the matching lines, such that the completed domains carry no solutions of negative index, and all
solutions of index at most one are regular. The index-zero sets are finite.
\end{lemma}

\begin{proof}
By Theorem~\ref{thm:torus-period} and the energy identity, the index of a solution
with given labels is bounded below. Proceed upward from the least index $d$ that
occurs, the minimum being taken over all of the finitely many labels, since breaking
changes labels. A sequence of solutions of index $d$ on the fixed completed domain has, by
Proposition~\ref{prop:torus-loss} for a fixed completed domain, a limit of total
index $d-4n$ with $n\ge0$ and nonnegative trajectory contributions; since no lower
index occurs, the limit has index $d$, there are no chains and $n=0$, so the space
is compact. Apply Theorem~\ref{thm:regularity}(ii) with $O$ in the
original pieces, including their original ends, and excluding the
inserted torus collar; every fixed primary or secondary perturbation has
support there. Each side keeps an original end on which the limiting
connection is flat on an odd surface collar; in the joined problem use the
two sets $O$ of the two sides. Thus the hypotheses of
Theorem~\ref{thm:regularity}(ii) hold without loops crossing the bridge. Part (iii) selects finitely many directions and
a generic small choice making index $d$ regular; a negative index
then does not occur. At the subsequent steps remove neighbourhoods
of already regular broken trajectories before taking a finite cover.
Ordinary end gluing remains valid with the completed torus tail:
the cutoff and error are on the original end, while the fixed tail
is controlled by its positive torus gap and the compact-remainder
Fredholm estimate; only these two facts about the tail enter, and no invertibility of an
unregularized completed operator is used.
Repeat through index one. The selected loop images and outputs have
compact union. Enlarge the common core and original-end truncations
to contain this union, and fix them before the finite collar or
bridge lengths vary. Beyond those truncations the original
asymptotic equations are unchanged. Finiteness at index zero follows from compactness and
regularity, and the lemma is proved.
\end{proof}

We fix such data on $W_{\pm,\infty}$ and use the same secondary perturbations on
$W_{\pm,\mathsf R}$ and $W^J_\lambda$ (split in the joined case). The primary data, end perturbations and almost complex structures are those of Section~\ref{subsec:folded}.

\emph{Asymptotic values and their tangent lifts.}
Let $P$ be a smooth configuration on $W_\mathsf L$
($\mathsf L\le\infty$) converging on each mixed end $e$ to
$R_e=(\pi^*B_{q_e},q_e)$ in all the weighted norms used below.
In the separate case suppose that $P=\alpha$ on a collar of the
flat face. Choose local stationary representatives equal to $\alpha$
on the torus collar, in product gauge near the genus face, and
put
\[
k_e(v)=\Bigl(\pi^*\tfrac{d}{ds}\Big|_{s=0}
 B_{\exp_{q_e}(sv)},v\Bigr),\qquad v\in T_{q_e}\mathcal C_e.
\]
The adapted exponential is fixed on the whole far strip of each
mixed end \cite[Lemma~4.4 and Section~4.1]{DFL}. Choose a cutoff
$\beta_e(\theta)$ with a fixed transition band, after a common
truncation on which $P$ is close to $R_e$.

A cutoff of $k_e(v)$ need not satisfy the linearized matching
condition at $P$. We use the following corrected lift
$K_{P,e}(v)$. Write $a_e=B_{q_e}|_{\Sigma'}$. On the far genus
trace express
\[
 A_P|_{\Sigma'}=E_{a_e}(b_P),\qquad
 u_P|_{\mathbb U_\partial}=e_{q_e}(v_P),\qquad [b_P]=v_P,
\]
using the flat surface chart of \cite[Lemmas~4.1 and~4.4]{DFL}.
Transport the surface trace of $k_e(v)$ by $DE_{a_e}(b_P)$ and
its section vector by $De_{q_e}(v_P)$, with the corresponding
adapted exponential along the strip. Differentiating flatness
and the identity $[E_{a_e}(b)]=e_{q_e}([b])$ proves the linearized
matching condition at $P$ and tangency on the Lagrangian arcs.
Multiply by $\beta_e$. Extend the difference between this genus
trace and that of $\beta_e k_e^A(v)$ into the fixed genus collar,
with zero normal component, and add it to $\beta_e k_e^A(v)$.
A weighted right inverse of the trace of fixed normal support is
used for this extension. The difference has zero stationary limit.
Thus $K_{P,e}$ has zero normal component, satisfies the actual
linearized boundary conditions, vanishes on the torus region, and
\begin{equation}\label{eq:corrected-profile}
 K_{P,e}(v)-\beta_e k_e(v)\in H^m_\delta,\qquad
 \|K_{P,e}(v)-\beta_e k_e(v)\|_{H^m_\delta}\le C|v|.
\end{equation}
Section differences here are read in the fixed end trivializations.
The constants are uniform for the preglues below: the genus collars
are fixed and their reference errors obey uniform weighted bounds.
The extension may have longitudinal tails; only its trace identity
and fixed normal support are used. Choose the local stationary
parameter frames consistently on the fixed cores. In the joined
case all these choices are products of the choices on the two sides.

\begin{definition}[gluing spaces]\label{def:torus-norms}
Let $P$ be a configuration on $W_\mathsf L$ as in the preceding paragraphs, for instance
the preglued configuration $P_{\mathsf L}(V)$ of Lemma~\ref{lem:preglued} below. Let $\mathcal X^m_\mathsf L$ be the space of pairs $w=(\zeta,\nu)$ with
$\zeta\in H^m_\delta(W_\mathsf L;\Lambda^1\otimes\ad)$ and $\nu\in H^m_\delta(\mathbb U_+;u_P^*T\Msf)$, where
$\|\zeta\|_{H^m_\delta}=\|e^{\delta\tau}\zeta\|_{H^m}$ with derivatives taken with respect to the
Levi-Civita connection and the connection of $P$, and similarly for $\nu$, subject
to the linearized boundary conditions: on the genus matching faces
$\mathbb U_\partial\times\Sigma'$ the normal component of $\zeta$ vanishes, its restriction $\zeta_\theta$
to each surface is $d_{A_\theta}$-closed and $[\zeta_\theta]=\nu(0,\theta)$
\cite[Definition~3.9(iii)]{DFL}; on the Lagrangian arcs $\nu$ is tangent to the
Lagrangians; and on the flat face (separate case) $\zeta(\partial_r)=0$ and
$d_\alpha(\zeta|_{T\times\{(\mathsf L,\theta)\}})=0$. Put
\[
X_\mathsf L=\mathcal X^m_\mathsf L\oplus\bigoplus_eT_{q_e}\mathcal C_e,\qquad\|(w,p)\|_{X_\mathsf L}=\|w\|_{\mathcal X^m_\mathsf L}+|p|,
\]
\[
Y_\mathsf L=H^{m-1}_\delta\bigl(W_\mathsf L;(\Lambda^0\oplus\Lambda^+)\otimes\ad\bigr)\oplus
H^{m-1}_\delta\bigl(\mathbb U_+;\Lambda^{0,1}\otimes u_P^*T\Msf\bigr),
\]
\[
D_{\mathsf L}(w,p)=\Bigl(d_{A_P}^*\zeta,\ \mathfrak L_P\bigl(w+\textstyle\sum_e K_{P,e}(p_e)\bigr)\Bigr),
\]
where $\mathfrak L_P$ is the linearization at $P$ of the mixed operator
$\mathfrak M$ (the perturbed anti-self-duality operator and the section operator).
\end{definition}

\begin{remark}\label{rem:torus-index-comparison}
The image of a profile is bounded in the weighted target. Indeed,
$\mathfrak L_P(\beta_e k_e(v))$ is the sum of
$\sigma(d\beta_e)k_e(v)$ and
$\beta_e(\mathfrak L_P-\mathfrak L_{R_e})k_e(v)$;
the second term decays with $P-R_e$, and
\eqref{eq:corrected-profile} handles the added correction. Hence the
part of $D_\mathsf L$ acting on asymptotic values has norm at most
$C|p|$. To compare indices with Theorem~\ref{thm:torus-period},
represent its extended domain by the same decaying space and the
same finite stationary parameter space, using Coulomb stationary
representatives corrected at the actual matching trace. The two
operators agree on the decaying summand. Their bounded stationary
columns, including their scalar rows, differ by a finite-rank map.
Thus $D_\mathsf L$ is Fredholm with the index of that theorem.
This argument compares indices, not kernels at a nonsolution.
The gauge row of $D_\mathsf L$ acts on $\zeta$ alone; the scalar
row of another realization must therefore also be retained in this
comparison. The nonlinear slice below imposes Coulomb gauge on
this decaying coordinate, rather than on the entire stationary
increment as in \cite[Proposition~4.6]{DFL}.
\end{remark}

\begin{lemma}[preglued references]\label{lem:preglued}
For $V\in\cM^0(W_\infty)$ let $V^{\rm ad}=\alpha+a_V$ on $r\ge S_V$ be as in
Lemma~\ref{lem:adapted}. For $\mathsf L\ge2S_V+2$ define $P_{\mathsf L}(V)$ on $W_\mathsf L$ to be $V^{\rm ad}$
outside $T\times[S_V,\mathsf L]\times\RR$ and $\alpha+\chi_{\mathsf L}(r)a_V$ on $T\times[S_V,\mathsf L]\times\RR$, where
$\chi_\mathsf L=1$ for $r\le \mathsf L/2$, $\chi_\mathsf L=0$ for $r\ge \mathsf L/2+1$ and $|\chi_\mathsf L'|\le2$. For a pair
$(V_-,V_+)$ define $P_\lambda(V_-,V_+)$ on $W^J_\lambda$ by the same truncation of both tails at
$r_\pm=\lambda/2$, with $\alpha$ on the middle of the bridge, joined by the bridge
identification. Then, with $P=P_\mathsf L$:
\begin{enumerate}[label=(\roman*)]
\item $P_\mathsf L$ is a smooth configuration with the labels of $V$ (of the pair), equals $\alpha$ near
  the face, and on each mixed end converges to $\pi^*B_{q_e}$ on $Y_{e,\mathsf L}$ (to
  $B_{p}\cup\alpha\cup B_q$ on $D_{i,\lambda}$);
\item $\mathfrak M(P_\mathsf L)$ is supported in $T\times[\mathsf L/2,\mathsf L/2+1]\times\RR$ (on both sides of the
  bridge) and $\|\mathfrak M(P_\mathsf L)\|_{Y_\mathsf L}\le Ce^{-\gamma_T \mathsf L/2}$;
\item $\ind D_\mathsf L=0$;
\item $\|P_\mathsf L-\alpha\|_{C^m(T\times[S,\mathsf L]\times\RR)}\le Ce^{-\gamma_T(S-S_0)}$ for $S\ge S_V$; on each end
  $\sup_\mathsf L\|P_\mathsf L-\pi^*B_{q_e}\|_{H^m_\delta(\{\theta\ge\Theta\})}\to0$ as $\Theta\to\infty$; and $P_\mathsf L=V^{\rm ad}$
  on $\{r\le \mathsf L/2\}$.
\end{enumerate}
\end{lemma}

\begin{proof}
(i) is immediate from Lemma~\ref{lem:adapted}, whose representative $V^{\rm ad}$ is smooth, with
$a_V$ smooth and $\chi_\mathsf L$ a smooth cutoff. (ii) Only the truncation band
carries an error. There the equation is the unperturbed anti-self-duality
equation, and since $F^+_{\alpha+a_V}=0$,
\[
F^+_{\alpha+\chi a_V}=(d\chi\wedge a_V)^++(\chi^2-\chi)(a_V\wedge a_V)^+ ,
\]
whose $H^{m-1}_\delta$ norm is at most $C\|a_V\|_{H^m_\delta(r\ge \mathsf L/2)}$, bounded by
Proposition~\ref{prop:tail-decay}(ii); the section is unchanged. (iii) The path
$\alpha+(1-s+s\chi_\mathsf L)a_V$ shows, as in the proof of Proposition~\ref{prop:torus-loss}, that $\kappa(P_\mathsf L)=\kappa(V)$ (the sum
$\kappa(V_-)+\kappa(V_+)$ in the joined case). By Theorem~\ref{thm:torus-period}(i)--(iii)
the index is $8\kappa+\mathfrak c(x)-\mathfrak c(a)$ with constants independent of $\mathsf L$, so
$\ind D_{P_\mathsf L}=\ind D_V=0$ (the sum of the two indices, by
Theorem~\ref{thm:torus-period}(ii)). (iv) follows from
Propositions~\ref{prop:tail-decay}(i) and~\ref{prop:tail-decay}(ii) and
Lemma~\ref{lem:adapted}(i)--(ii), using (i) for the pointwise
$C^m$ bound.
\end{proof}

\begin{lemma}[uniform scalar inverse]\label{lem:torus-scalar}
There are $r_0$, $\mathsf L_0$, $C_P$ and $C$ such that for
$\mathsf L\in[\mathsf L_0,\infty]$ and every connection $A$ with
$\|A-A_{P_\mathsf L(V)}\|_{C^m}\le r_0$:
\begin{enumerate}[label=(\roman*)]
\item $\|\xi\|_{L^2}\le C_P\|d_A\xi\|_{L^2}$ for every $\xi\in L^2_1(W_\mathsf L;\ad)$, with no
  boundary condition on $\xi$;
\item the Neumann map
$\xi\mapsto\bigl(d_A^*d_A\xi,\nabla_n\xi|_{\partial W_\mathsf L}\bigr)$
is an isomorphism from $H^{m+1}_\delta(W_\mathsf L;\ad)$ to
$H^{m-1}_\delta(W_\mathsf L;\ad)\oplus
 H^{m-1/2}_\delta(\partial W_\mathsf L;\ad)$, with inverse norm
at most $C$.
\end{enumerate}
At infinity $P_\infty(V)=V^{\rm ad}$. The boundary
consists of the genus matching faces and, when present, the flat face.
\end{lemma}

\begin{proof}
(i) On the far torus region the fibre gap and
Lemma~\ref{lem:preglued}(iv) give
$\|\xi\|_{L^2}\le 2\mu_T^{-1}\|d_A\xi\|_{L^2}$ by integration
over the planar base. There is no boundary term or factor depending
on its area. The stationary spatial compressions have a uniform
scalar gap as well. Otherwise normalized sections with covariant
derivatives tending to zero would have vanishing mass on their
torus collars, and a parallel limit on each fixed compact core.
That limit vanishes on the overlap with the torus collar and hence
on the connected core, a contradiction. The same argument covers
the two cores of a joined slice. Integrating over the mixed-end
coordinate and absorbing the uniformly small reference error gives
the mixed-end estimate. The original gauge ends use their fixed
irreducible cross-sections.

On the whole domain a normalized counterexample now has mass only
on the fixed compact cores. Rellich compactness gives a parallel
limit there, zero on the controlled overlap bands, again a
contradiction. Thus
\[
 \|\xi\|_{L^2}\le C_P\|d_A\xi\|_{L^2}
\]
uniformly at finite and infinite lengths. A sufficiently small
$C^0$ change of $A$ preserves this estimate by absorption; take
$C_P$ for that neighbourhood.

(ii) This is the Lax--Milgram argument of the proof of Theorem~\ref{thm:neck}. For
$U=e^{\delta\tau}\xi$ the weak Neumann form becomes
$\langle d_AU-\delta\,d\tau\,U,\ d_AV+\delta\,d\tau\,V\rangle$, whose diagonal
$\|d_AU\|^2-\delta^2\|d\tau\,U\|^2\ge\frac34\|d_AU\|^2$ is coercive by (i), $|d\tau|\le1$ and
$\delta C_P\le\frac12$. Since $\partial_n\tau=0$, the conjugation does not change the Neumann
condition. Local Neumann regularity on the unit patches of Step~3 of the proof of
Proposition~\ref{prop:torus-inverse}, on which the coefficients of $A$ are bounded in $C^m$ by
Lemma~\ref{lem:preglued}(iv) and the choice of $r_0$, gives the estimate in $H^{m+1}_\delta$ and
surjectivity, directly on the completed domain.
\end{proof}

\begin{proposition}[uniform inverse]\label{prop:torus-inverse}
There are $\mathsf L_0$ and $C_{\rm inv}$ such that for $\mathsf L\ge \mathsf L_0$ and every $V$ (every pair) in the
finite completed index-zero sets, $D_\mathsf L\colon X_\mathsf L\to Y_\mathsf L$ at $P=P_{\mathsf L}(V)$ is an
isomorphism with $\|D_{\mathsf L}^{-1}\|\le C_{\rm inv}$.
\end{proposition}

\begin{proof}
By Lemma~\ref{lem:preglued}(iii) and the remark after
Definition~\ref{def:torus-norms}, $D_\mathsf L$ is Fredholm of index zero, so it suffices to
prove $\|(w,p)\|_{X_\mathsf L}\le C_{\rm inv}\|D_{\mathsf L}(w,p)\|_{Y_\mathsf L}$. Write $D^{\rm dec}_{\mathsf L}w=D_{\mathsf L}(w,0)$.

\emph{Step 1: models.} Fix parameters $S$, $\Theta$, $T_0$, $\ell$, chosen below in this
order.

(T) The torus region $\{r\ge S\}$ (the part of the bridge at distance at least $S$
from both pieces in the joined case), including its corners inside the mixed ends.
There $D^{\rm dec}_\mathsf L=D_\alpha+\mathfrak a$, where $D_\alpha$ is the operator of
Proposition~\ref{prop:flat-fibre} and $\mathfrak a$ is of order zero with
$|\mathfrak a|\le Ce^{-\gamma_T(S-S_0)}$ by Lemma~\ref{lem:preglued}(iv); no perturbation acts
there, and the corners of the torus region inside the mixed ends carry only the
gauge field. The base is a strip $[S,\mathsf L]\times\RR$ with the flat-fibre condition at $r=\mathsf L$,
or a strip in the bridge; on it Proposition~\ref{prop:flat-fibre} with the weight
$e^{\delta\tau_0(\theta)}$ gives $\|e^{\delta\tau}\zeta\|\le(\mu_T-\delta)^{-1}\|e^{\delta\tau}D_\alpha\zeta\|$ for $\zeta$
vanishing near $r=S$. For a straight boundary the Bochner identity for $d+d^*$ on
the base with absolute boundary conditions has no boundary term, so the identity
in the proof of Proposition~\ref{prop:flat-fibre} also bounds $\|\nabla\zeta\|$ by
$\|D_\alpha\zeta\|$. Hence $\|\zeta\|_{H^1_\delta}\le C_T\|D_\alpha\zeta\|_{L^2_\delta}$ with $C_T$ depending only on
$\mu_T$ and $\delta$.

(E) The mixed end $e$, beyond $\Theta$, restricted to $r\le S+\ell$.
First identify the actual boundary domain with the stationary one.
The surface Hodge and harmonic-block maps of
\cite[Proposition~5.27]{DFL-mixed}, applied at paired depths in
the fixed gauge and section collars, give an isomorphism
$Q(\theta)$ from the stationary domain to the actual domain.
It is the identity off those collars and extends by the identity
over the added torus region. It and its inverse are bounded on
$L^2$ and $H^1$, and their differences from the identity and their
longitudinal derivatives tend uniformly to zero as
$\Theta\to\infty$ \cite[Lemma~5.67]{DFL-mixed}. These results are stated for a fixed
compact slice, locally in $\theta$; they apply uniformly on the infinite end with a
growing torus collar because $Q(\theta)$ is supported near the genus surface, where the
geometry is fixed, and is the identity on the torus collar.
Use the corresponding fixed trivializations of the target.
The transformed operator is
$\sigma(\partial_\theta+\mathcal A_{q_e,\mathsf L})+\mathfrak b$, where $\mathcal A_{q_e,\mathsf L}$ is the
stationary operator of Proposition~\ref{prop:stationary-gap} on the whole slice
$Y_{e,\mathsf L}$ \cite[Section~4.2.1, (4.16)]{DFL} and $\|\mathfrak b\zeta\|_{L^2}\le\eta(\Theta)\|\zeta\|_{H^1}$ with
$\eta(\Theta)\to0$ uniformly in $\mathsf L$ by
Lemma~\ref{lem:preglued}(iv) and these transport bounds.
For a cutoff field vanishing near $\theta=\Theta$, apply
$Q^{-1}$ before extending by zero. The resulting field belongs
to the stationary domain on $\RR_\theta\times Y_{e,\mathsf L}$.
Since $\mathcal A=\mathcal A_{q_e,\mathsf L}$ is self-adjoint with
spectrum meeting $(-\gamma,\gamma)$ only in $0$, and $0<\delta<\gamma$, Fourier transform in
$\theta$ gives $\|(\partial_\theta+\mathcal A-\delta)^{-1}\|\le\min(\delta,\gamma-\delta)^{-1}$; the identity
$\|\partial_\theta\xi\|^2+\|(\mathcal A-\delta)\xi\|^2=\|(\partial_\theta+\mathcal A-\delta)\xi\|^2$ and the elliptic
estimate $\|\xi\|_{H^1(Y_{e,\mathsf L})}\le C(\|\mathcal A\xi\|+\|\xi\|)$, which is uniform in $\mathsf L$
because it is local, $Y_{e,\mathsf L}$ has bounded geometry and on the collar
$\mathcal A=\mathcal A_\alpha$ (Corollary~\ref{cor:flat-fibre-line}), give
$\|e^{\delta\theta}\xi\|_{H^1}\le C_E\|e^{\delta\theta}(\partial_\theta+\mathcal A)\xi\|$ with $C_E$ independent of
$\mathsf L$. Transporting back and absorbing the
$H^1$-to-$L^2$ error gives the estimate on the actual mixed end.

(O) The original gauge and ordinary ends beyond $T_0$: translation-invariant
operators with invertible asymptotic operators, as in
Theorem~\ref{thm:neck}, perturbed by terms that are small for $T_0$ large.

(K) A fixed compact region $K'$ containing the rest: the local elliptic estimates
of \cite[Theorem~4.17 and Remark~4.19]{DFL} at the matching faces and interior
estimates give $\|\chi w\|_{H^1}\le C(\|D^{\rm dec}_{\mathsf L}w\|_{L^2(K'')}+\|w\|_{L^2(K'')})$ for a
cutoff $\chi$ supported in $K'$ and a slightly larger compact $K''$.

\emph{Step 2: the compact-remainder estimate.} Let $\chi$ be a smooth function on $\RR$
with $\chi=0$ on $(-\infty,0]$ and $\chi=1$ on $[1,\infty)$. Patch the models with the partition
of unity $\psi_T+\sum_e\psi_{E,e}+\psi_G+\psi_O+\psi_K=1$, where $\psi_T=\chi((r-S)/\ell)$ on the
torus region (in the joined case, the product of the two such functions of $r_-$ and $r_+$),
$\psi_{E,e}=(1-\psi_T)\chi(\theta-\Theta)$ on the end $e$, $\psi_G,\psi_O$ are cutoffs of fixed
slope on the original ends beyond $T_0$, and $\psi_K$ has compact support. The only noncompact
transition region is the corner band $\{\theta\ge\Theta,\ S\le r\le S+\ell\}$, where the torus
tail enters the mixed ends; there both (T) and (E) apply, and the cutoffs vary in $r$ with slope
at most $C/\ell$. The gradient of $\psi_T$ is at most $C/\ell$ on all of $\{S\le r\le S+\ell\}$;
the other gradients are supported in the compact set
$K=\{\Theta\le\theta\le\Theta+1,\ r\le S+\ell\}$, in the fixed parts of the slices and in
compact transition regions on the original ends. Choose $S$ and $\Theta$ with
$C_T\sup|\mathfrak a|\le\tfrac12$ and $C_E\eta(\Theta)\le\tfrac12$, which does not depend on
$\ell$ since $C_T$ and $C_E$ do not, then $T_0$ similarly. Applying the model estimates to
$\psi_iw$ and summing,
\[
\|w\|_{H^1_\delta}\le C\bigl(\|D^{\rm dec}_{\mathsf L}w\|_{L^2_\delta}
+\ell^{-1}\|w\|_{L^2_\delta}+\|w\|_{L^2(K_1)}\bigr),
\]
where $K_1$ is a compact set containing $K$, $K''$ and the compact transition regions. The
coefficient of $\ell^{-1}\|w\|_{L^2_\delta}$ comes from the cutoff gradients through the
models (T) and (E), whose constants do not depend on $\ell$. A choice of $\ell$ with
$C\ell^{-1}\le\tfrac12$ therefore gives
\begin{equation}\label{eq:torus-remainder}
\|w\|_{H^1_\delta}\le C\bigl(\|D^{\rm dec}_{\mathsf L}w\|_{L^2_\delta}+\|w\|_{L^2(K_1)}\bigr)\qquad(\mathsf L\ge \mathsf L_0),
\end{equation}
with $\mathsf L_0=2(S+\ell)+2$, so that $K_1\subset\{r\le \mathsf L/2\}$ and $P_\mathsf L=V^{\rm ad}$ on $K_1$. In the
joined case $K_1=K_{1,-}\sqcup K_{1,+}$, one compact set in the core region of each side,
at a fixed position relative to that side.

\emph{Step 3: higher order and asymptotic values.} Cover $W_\mathsf L$ by unit patches of finitely
many isometry types: flat patches $T\times Q$ in the torus region, flat half-patches at
the face, translates of fixed patches along the mixed ends and the original ends,
and fixed patches of the core. The coefficients of $D_\mathsf L$ are bounded in $C^m$
uniformly in $\mathsf L$ (Lemma~\ref{lem:preglued}(iv)). On the flat patches the local
estimate is that of $D_\alpha+\mathfrak a$ with the flat-fibre domain at the face: the
tangential derivatives commute with $D_\alpha$ and preserve the domain, and the
normal derivatives are recovered from the equation, as in the proof of
Theorem~\ref{thm:flat-boundary-regularity}(b); on the other patches it is
\cite[Theorem~4.17 and Remark~4.19]{DFL}. The weight varies by a bounded factor on
each patch. Summing with bounded overlap and using~\eqref{eq:torus-remainder},
\[
\|w\|_{\mathcal X^m_\mathsf L}\le C\bigl(\|D^{\rm dec}_{\mathsf L}w\|_{H^{m-1}_\delta}+\|w\|_{L^2(K_1)}\bigr).
\]
Since the part of $D_\mathsf L$ acting on asymptotic values is bounded by $C|p|$,
\begin{equation}\label{eq:torus-remainder-full}
\|(w,p)\|_{X_\mathsf L}\le C\bigl(\|D_{\mathsf L}(w,p)\|_{Y_\mathsf L}+\|w\|_{L^2(K_1)}+|p|\bigr).
\end{equation}

\emph{Step 4: the limit operator is injective.}
Let $D_V$ be the operator of Definition~\ref{def:torus-norms} at
$P=V^{\rm ad}$ on $W_\infty$, and suppose $D_V(w,p)=0$.
Write $K_{V,e}=(K^A_{V,e},K^S_{V,e})$ and let
$h_e(v)=(h_e^A(v),v)$ be the Coulomb stationary kernel
representative of Theorem~\ref{thm:torus-period}. Since $k_e^A(v)$
and $h_e^A(v)$ represent the same primary deformation and have
zero normal component, their difference is
$d_{B_{q_e}}\psi_e(v)$, with $\psi_e$ satisfying the spatial
Neumann condition. Put
\[
 \zeta_1=\zeta+\sum_eK^A_{V,e}(p_e)
             -d_V\Bigl(\sum_e\beta_e\psi_e(p_e)\Bigr),\qquad
 \nu_1=\nu+\sum_eK^S_{V,e}(p_e).
\]
This variation satisfies the geometric linearized matching condition
and $\mathfrak L_V(\zeta_1,\nu_1)=0$: the original variation
uses the corrected lifts, and the subtracted term is an
infinitesimal gauge transformation at the solution $V$.
Its asymptotic representatives are $h_e(p_e)$. The remaining field
is decaying by \eqref{eq:corrected-profile}; the other errors contain
$d\beta_e$ or $V-R_e$ and its derivatives.

The divergence of $\zeta_1$ and its normal trace are decaying.
The latter need not vanish before the final gauge correction:
the spatial Neumann condition on $\psi_e$ uses $B_{q_e}$,
whereas $d_V$ uses the actual connection. Apply
Lemma~\ref{lem:torus-scalar} with both these data to solve
\[
 d_V^*d_V\xi=-d_V^*\zeta_1,\qquad
 \nabla_{n,V}\xi=-\iota_n\zeta_1.
\]
Then $(\zeta_1+d_V\xi,\nu_1)$ belongs to the kernel of the
extended Coulomb realization in Theorem~\ref{thm:torus-period}.
At $V$ that operator is onto and has index zero, so its kernel
vanishes. Hence $p=0$, $\nu=0$ and $\zeta=-d_V\xi$.
The original scalar and normal conditions give the homogeneous
Neumann equation for $\xi$, and the scalar lemma implies
$\xi=0$ and $\zeta=0$.

In the joined case the section, the Lagrangians and all the
profile corrections were chosen as products. The limiting
operator is therefore the direct sum of the two separate
operators, with their independent stationary coordinates, and
is injective by the same argument.

\emph{Step 5: contradiction.} If the bound failed, there would be $\mathsf L_j\to\infty$ and
$(w_j,p_j)$ of norm one with $D_{\mathsf L_j}(w_j,p_j)\to0$, for a fixed $V$ (the sets are finite),
and $\|w_j\|_{L^2(K_1)}+|p_j|\ge1/C$ by \eqref{eq:torus-remainder-full}. Every compact subset of
$W_\infty$ embeds isometrically in $W_{\mathsf L_j}$ for large $j$, with $P_{\mathsf L_j}=V^{\rm ad}$
on it, where $D_{\mathsf L_j}$ coincides with $D_V$. A diagonal subsequence therefore converges,
weakly in $H^m_{\rm loc}$ and strongly in $L^2_{\rm loc}$, to some $(w_\infty,p_\infty)$ of norm at
most one that satisfies the linear boundary conditions and $D_V(w_\infty,p_\infty)=0$. Step~4
gives $(w_\infty,p_\infty)=0$, contradicting the lower bound on $K_1$. In the joined case the exhaustions are taken separately on
the two sides: the part of $W^J_{\lambda_j}$ within distance $\lambda_j/2$ of the $\pm$ core is
identified with the corresponding part of $W_{\pm,\infty}$, the bridge between them
escapes to infinity, and the limit is a pair $(w_{\pm,\infty},p_{\pm,\infty})$ on
$W_{-,\infty}\sqcup W_{+,\infty}$, which vanishes by Step 4, while $K_1=K_{1,-}\sqcup K_{1,+}$ retains
mass. Invertibility follows from index zero.
\end{proof}

\begin{lemma}[the flat-fibre chart]\label{lem:face-chart}
Let $s=m-\frac12$ and let $H^s_\delta$ denote the weighted Sobolev space of
$\ad$-valued torus one-forms on the face $T\times\RR_\theta$, with weight $e^{\delta\tau_0(\theta)}$.
Let $\Pi_{\rm cl}$ be the fibrewise projection onto $d_\alpha$-closed forms and
$\mathsf Q=d_\alpha^*(d_\alpha d_\alpha^*)^{-1}$ on fibrewise two-forms.
\begin{enumerate}[label=(\alph*)]
\item There are $r_3$, $C$ such that for every closed $b$ with $\|b\|_{H^s_\delta}\le r_3$ the
  equation $c=-\tfrac12\mathsf Q[(b+c)\wedge(b+c)]$ has a unique coexact solution $c(b)$
  with $\|c\|_{H^s_\delta}\le r_3$. Then $\alpha+b+c(b)$ is flat on every fibre, $c(0)=0$,
  $Dc(0)=0$, and $\|c(b)-c(b')\|_{H^s_\delta}\le C(\|b\|+\|b'\|)\|b-b'\|$. Conversely, if
  $\alpha+\tilde b$ is flat on every fibre and $\|\tilde b\|_{H^s_\delta}\le r_3/C$, then
  $\tilde b=b+c(b)$ with $b=\Pi_{\rm cl}\tilde b$.
\item There is a linear map $\operatorname{Ext}_\delta\colon H^s_\delta(T\times\RR)\to H^m_\delta(W_\mathsf L)$,
  with norm bounded independently of $\mathsf L$, whose values are torus one-forms
  supported in the collar $T\times[\mathsf L-1,\mathsf L]\times\RR$ with vanishing $dr$ and $d\theta$
  components, and whose restriction to the face is the identity.
\item The map $\Psi^{\rm face}(\zeta)=\zeta+\operatorname{Ext}_\delta\bigl(c(\zeta_T|_{r=\mathsf L})\bigr)$,
  defined on the ball of radius $r_3/C$ of fields satisfying the linear flat-fibre
  condition, takes values in fields $x$ for which $\alpha+x$ satisfies the nonlinear
  flat-fibre condition and $x(\partial_r)=\zeta(\partial_r)$ at the face;
  in particular it preserves the zero-normal subspace. It is inverted by
  $x\mapsto x-\operatorname{Ext}_\delta\bigl(c(\Pi_{\rm cl}x_T|_{r=\mathsf L})\bigr)$, and both maps have
  first and second derivatives bounded independently of $\mathsf L$.
\end{enumerate}
\end{lemma}

\begin{proof}
(a) The operators $\Pi_{\rm cl}$ and $\mathsf Q$ act fibrewise and are Fourier
multipliers in the eigenmodes of the flat bundles $\ad\alpha$ on $T$; they commute
with $\partial_\theta$ and with multiplication by functions of $\theta$, in particular with
$e^{\delta\tau_0}$. Since $\Delta_\alpha\ge\mu_T^2$ on each fibre, $\mathsf Q$ is bounded on the
isotropic space $H^s(T\times\RR)$ (its symbol is bounded by $\mu_T^{-1}$ in the torus
frequencies and does not involve the $\theta$-frequency). As $s>\frac32$, $H^s(T\times\RR)$ is an algebra, and
$\|e^{\delta\tau_0}fg\|_{H^s}\le C\|e^{\delta\tau_0}f\|_{H^s}\|e^{\delta\tau_0}g\|_{H^s}$ because $\tau_0\ge0$ and
$e^{-\delta\tau_0}$ has bounded derivatives. Hence the right side of the equation is a
contraction of the ball of coexact forms of radius $r_3$ for $r_3$ small, and the
contraction estimates give the stated bounds. Flatness:
$F_{\alpha+b+c}=d_\alpha c+\tfrac12[(b+c)\wedge(b+c)]=0$, because $d_\alpha b=0$ and $d_\alpha\mathsf Q$
is the identity on two-forms, $H^2(T;\ad\alpha)$ being zero. Conversely, write
$\tilde b=b+\tilde c$ with $b=\Pi_{\rm cl}\tilde b$ and $\tilde c$ coexact; flatness gives
$d_\alpha\tilde c=-\tfrac12[\tilde b\wedge\tilde b]$, and $\mathsf Qd_\alpha$ is the identity on
coexact forms, so $\tilde c$ solves the equation and equals $c(b)$ by uniqueness.
(b) Let $\Lambda=(1+\Delta_{T\times\RR})^{1/2}$, acting on torus one-forms on $T\times\RR$ with the
connection $\alpha$, and $E_0(c)(\rho)=\chi(\rho)e^{-\rho\Lambda}c$, $\rho=\mathsf L-r$, with a fixed cutoff
$\chi$ supported in $[0,1)$ and equal to one near $0$. Mode by mode
$\int_0^\infty|\partial_\rho^je^{-\rho\mu}|^2d\rho=\mu^{2j-1}/2$, so $E_0$ is bounded from $H^{k-1/2}$ to
$H^k$ of the half-space for integers $k$, and by interpolation from $H^s$ to $H^m$.
Put $\operatorname{Ext}_\delta(c)=e^{-\delta\tau_0}E_0(e^{\delta\tau_0}c)$. The collar is a translate of a
fixed product collar, so the norm does not depend on $\mathsf L$. (c) By (b) and the trace
theorem, $\zeta_T|_{r=\mathsf L}$ lies in $H^s_\delta$ with norm at most $C\|\zeta\|_{H^m_\delta}$; the torus
restriction of $\zeta+\operatorname{Ext}_\delta(c)$ at the face is $\zeta_T|_{r=\mathsf L}+c(\zeta_T|_{r=\mathsf L})$, flat by
(a) since $\zeta_T|_{r=\mathsf L}$ is closed, and the normal component is unchanged. For the
inverse, if $\alpha+x$ satisfies the nonlinear condition, (a) gives
$x_T|_{r=\mathsf L}=b+c(b)$ with $b=\Pi_{\rm cl}x_T|_{r=\mathsf L}$, so
$x-\operatorname{Ext}_\delta(c(b))$ has closed restriction $b$. The derivative bounds follow
from those of $c$ and the boundedness of $\operatorname{Ext}_\delta$, the trace and
$\Pi_{\rm cl}$.
\end{proof}

\emph{The chart.}
We specify the genus trace extension, including its stationary
part. For a genus trace $f$ with stationary limits $f_e$, put
\[
 \mathcal T^{\rm ext}f=
 \mathcal T_\delta\Bigl(f-\sum_e\beta_e f_e\Bigr)
                  +\sum_e\beta_e\mathcal T_e f_e.
\]
Here $\mathcal T_\delta$ is a weighted right inverse of trace on
fixed genus collars, and $\mathcal T_e$ is its spatial stationary
counterpart. Both lift surface-tangential forms with zero normal
and longitudinal components. Their norms are uniform, as for the
Poisson extension in Lemma~\ref{lem:face-chart}(b). Thus no
constant longitudinal mode is assigned a weighted Sobolev norm.

For $z=(w,p)$, $w=(\zeta,\nu)$, set
\[
 (\zeta',\nu')=w+\sum_e K_{P,e}(p_e),\qquad
 Q_a(b)=E_a(b)-a-b,
\]
where $E_a$ is the flat surface chart compatible with the adapted
section exponential. Apply $Q$ fibrewise to
$b=\zeta'|_{\Sigma'}$, using $a=A_P|_{\Sigma'}$. At a stationary
reference the resulting connection is
\[
 \widehat B_{e,p}=B_{q_e}+k_e^A(p)
             +\mathcal T_e Q_{a_e}(k_e^A(p)|_{\Sigma'}).
\]
Choose $B_{\exp_{q_e}p}$ to have the same genus trace as
$\widehat B_{e,p}$, and take this re-chosen family as the fixed family $B_q$ of
stationary references used throughout, in particular in
Proposition~\ref{prop:torus-entry}; it is in product gauge on the genus collar, equals
$\alpha$ near the torus, and any $\theta$-independent normalization is admissible in
Proposition~\ref{prop:growing-ends}(ii). The two traces represent the same class, so
a smooth near-identity surface gauge provides this choice; extend
it through the fixed genus collar and as the identity near the
torus. Their first variations agree, so the adjustment is quadratic
in $p$. In product gauge the difference
$C_e(p)=B_{\exp_{q_e}p}-\widehat B_{e,p}$ has zero genus
tangential and normal trace, vanishes on the torus collar, and
satisfies $C_e(0)=DC_e(0)=0$.

For fixed end cutoffs $\varphi_e$, equal to one sufficiently far
out, define $\Psi_\mathsf L(z)$ to have connection
\[
 A_P+\zeta'
 +\mathcal T^{\rm ext}\bigl(Q_{A_P|_{\Sigma'}}
                                      (\zeta'|_{\Sigma'})\bigr)
 +\sum_e\varphi_e\pi^*C_e(p_e)
 +\operatorname{Ext}_\delta
                    \bigl(c(\zeta_T|_{r=\mathsf L})\bigr)
\]
and section $e_{u_P}(\nu')$, with the adapted exponentials of
\cite[Lemma~4.4]{DFL}. The face term is absent in the joined and
completed cases. The genus trace is flat and has precisely the
class of the section. The $C_e$ term preserves that trace and
makes the limiting connection exactly $B_{\exp_{q_e}p_e}$.
All profile and genus corrections vanish on the torus region.
Write $P^{(p)}=\Psi_\mathsf L(0,p)$; it equals $P$ there.
The same formula gives an unsliced chart if the normal component
of the decaying field is allowed to vary. We use that version
before imposing the Neumann gauge condition.

Define
\[
 \mathcal F_\mathsf L(z)=
 \bigl(d_{A_P}^*\zeta,
                   \Pi_P\mathfrak M(\Psi_\mathsf L(z))\bigr)
 \in Y_\mathsf L.
\]
The anti-self-duality row is read in $\Lambda^+\otimes\ad$ for
the given metric. For the section row use its tangent-vector
representation $\partial_su+J(u)\partial_\theta u$ from
\eqref{eq:mixed}. Embed $\Msf$ isometrically in $\mathbb R^N$
and project this vector orthogonally to $T_{u_P}\Msf$.
When the target is written as a $(0,1)$-form, reconstruct that
form from its projected value using the complex structure at
$u_P$. This defines $\Pi_P$ on the section row; it is the identity
on the other rows. The reference
section has uniformly bounded derivatives, so this projection is
bounded on the weighted target. For $z$ small it is an isomorphism
on the tangent spaces along $u_{\Psi_\mathsf L(z)}$, with uniformly
bounded inverse.

\begin{lemma}[chart and quadratic estimate]\label{lem:torus-quadratic}
There are $\mathsf L_0$, $r_1$ and $C_{\rm quad}$ such that for $\mathsf L\ge \mathsf L_0$ and $P=P_{\mathsf L}(V)$:
\begin{enumerate}[label=(\alph*)]
\item $\Psi_{\mathsf L}(z)$ satisfies the nonlinear matching, Lagrangian and flat-fibre
  conditions, and zeros of $\mathcal F_\mathsf L$ in the ball of radius $r_1$ are solutions of the mixed
  equation.
  Conversely, if $U$ satisfies these conditions and $\|U-P^{(p)}\|_{H^m_\delta}+|p|\le r_1/C_{\rm quad}$,
  then there are $\xi$ with $\|\xi\|_{H^{m+1}_\delta}\le C_{\rm quad}\|U-P^{(p)}\|_{H^m_\delta}$ and $w$ with
  \[
  \exp(\xi)^*U=\Psi_{\mathsf L}(w,p),\qquad d^*_{A_P}\zeta=0,\qquad\|w\|_{\mathcal X^m_\mathsf L}\le C_{\rm quad}\bigl(\|U-P^{(p)}\|_{H^m_\delta}+|p|\bigr).
  \]
\item $\mathcal F_{\mathsf L}(0)=(0,\mathfrak M(P_\mathsf L))$ and $D\mathcal F_{\mathsf L}(0)=D_\mathsf L$;
\item $\mathcal F_\mathsf L$ is of class $C^2$ on the ball of radius $r_1$ in $X_\mathsf L$, with
  $\|D^2\mathcal F_\mathsf L\|\le C_{\rm quad}$. Consequently $\mathcal F_{\mathsf L}(z)=\mathcal F_{\mathsf L}(0)+D_{\mathsf L}z+N_{\mathsf L}(z)$ with
  \[
  \|N_{\mathsf L}(z)-N_{\mathsf L}(z')\|_{Y_\mathsf L}\le C_{\rm quad}\bigl(\|z\|_{X_\mathsf L}+\|z'\|_{X_\mathsf L}\bigr)\|z-z'\|_{X_\mathsf L}.
  \]
\end{enumerate}
\end{lemma}

\begin{proof}
(a) The trace identities in the construction give the exact genus
matching and Lagrangian conditions. Lemma~\ref{lem:face-chart}
gives the flat-face condition; the collars are disjoint and their
corrections have zero normal component. The injectivity of $\Pi_P$
on the actual section target proves the assertion about zeros.

For the converse first use the unsliced chart. Extract $p$ from the
stationary limits and subtract $\sum_e\varphi_e\pi^*C_e(p_e)$.
The inverse surface exponential recovers the closed genus trace,
and $\Pi_{\rm cl}$ recovers the closed torus trace. Subtract their
extended nonlinear remainders and the corrected profiles $K_{P,e}$.
Together with the inverse section exponential, this gives the
unique decaying coordinate $\widetilde w$. The trace and extension
bounds give
$\|\widetilde w\|\le C\|U-P^{(p)}\|_{H^m_\delta}$ uniformly.

Let $Z(\exp(\xi)^*U)$ be the connection coordinate of $\exp(\xi)^*U$ in the unsliced chart, and
solve $\bigl(d_{A_P}^*Z,\iota_nZ\bigr)(\exp(\xi)^*U)=0$ for a small decaying $\xi$; the unsliced
chart is needed since a gauge transform need not have zero normal coordinate. As in the entry
step of the proof of Theorem~\ref{thm:neck}, the derivative in $\xi$ at $\xi=0$, $U=P$ is, up to
sign, the Neumann map of Lemma~\ref{lem:torus-scalar}(ii), the second derivatives are bounded
uniformly for $m\ge4$, and the quantitative implicit function theorem gives the fixed radius and
the bounds; at $U=P^{(p)}$ the coordinate vanishes, so the bound for $\xi$ involves only
$\|U-P^{(p)}\|_{H^m_\delta}$.

(b) Each nonlinear surface correction has zero first derivative
at zero, as does $C_e(p)$. Hence
\[
 D\Psi_\mathsf L(0)(w,p)=w+\sum_e K_{P,e}(p_e).
\]
The reference section solves its equation, so differentiating its
target projection adds no term. The gauge row acts on $\zeta$
alone. This gives exactly $D_\mathsf L$.

(c) The genus chart and its inverse have uniform $C^2$ bounds by the argument of
Lemma~\ref{lem:neck-chart}(b): they are assembled from the surface Hodge decomposition, the
exponential of the exact potential and a finite-rank harmonic part along the fixed genus
collars, the surface Hodge inverses have a common gap, the trace space $H^{m-1/2}$ on each
three-dimensional face is an algebra, and the trace extensions are uniformly bounded. The
flat-face chart has the bounds of Lemma~\ref{lem:face-chart}.
On the torus the equation is quadratic curvature; bounded geometry
and bounded overlap give multiplication constants independent of
its length. Fixed holonomy families obey the corresponding local
bounds by differentiating parallel transport; for the primary
terms see \cite[Proposition~6.1]{DFL}. Secondary dependence and
output sets stay in the fixed compact cores.

Finally, on each far mixed end subtract the same chart expression
at its stationary reference. Every remaining term contains a
decaying reference error, a decaying coordinate, or a derivative
of a cutoff on a fixed band. The separated extension
$\mathcal T^{\rm ext}$ respects this decomposition. The stationary
pair with limit $\exp_{q_e}p_e$ solves the equation exactly, so
it contributes no nondecaying term after any differentiation in
$p$. The weighted product estimates therefore bound the second
derivative of $\mathcal F_\mathsf L$ uniformly. The displayed
quadratic estimate follows by the mean value theorem.
\end{proof}

\begin{proposition}[existence and local uniqueness]\label{prop:torus-existence}
Let $r=\min\{r_1,(4C_{\rm inv}C_{\rm quad})^{-1}\}$. There is $\mathsf L_1$ such that for $\mathsf L\ge \mathsf L_1$ and every
$V$ (every pair) in the completed index-zero sets there is a unique $z_{\mathsf L}(V)$ in the
ball of radius $r$ in $X_\mathsf L$ with $\mathcal F_{\mathsf L}(z_{\mathsf L}(V))=0$. It satisfies
$\|z_{\mathsf L}(V)\|\le2C_{\rm inv}\|\mathfrak M(P_\mathsf L)\|_{Y_\mathsf L}\le Ce^{-\gamma_T \mathsf L/2}$, and
$\mathcal G_{\mathsf L}(V):=\Psi_{\mathsf L}(z_{\mathsf L}(V))$ is a regular solution.
\end{proposition}

\begin{proof}
The map $\mathcal F_\mathsf L$ takes values in $Y_\mathsf L$ because its target is projected explicitly: the
anti-self-duality component onto the self-dual forms of $W_\mathsf L$, and the Cauchy--Riemann
component pointwise onto the tangent spaces along the preglued section $u_P$. By
Proposition~\ref{prop:torus-inverse} and Lemma~\ref{lem:torus-quadratic}, the map
$z\mapsto-D_{\mathsf L}^{-1}(\mathcal F_{\mathsf L}(0)+N_{\mathsf L}(z))$ has Lipschitz constant at most
$2C_{\rm inv}C_{\rm quad}r\le\tfrac12$ on the ball of radius $r$, and maps it into itself once
$C_{\rm inv}\|\mathcal F_{\mathsf L}(0)\|\le r/2$, which holds for $\mathsf L\ge \mathsf L_1$ by
Lemma~\ref{lem:preglued}(ii). As in the proof of Theorem~\ref{thm:neck}, this gives $z_\mathsf L(V)$,
the only zero in that ball, with the stated bound; since $r\le r_1$, $\mathcal G_\mathsf L(V)$ is a
solution (Lemma~\ref{lem:torus-quadratic}(a)). The derivative
$D\mathcal F_{\mathsf L}(z_\mathsf L)=D_\mathsf L+DN_{\mathsf L}(z_\mathsf L)$, with
$\|DN_{\mathsf L}(z_\mathsf L)\|\le2C_{\rm quad}r\le(2C_{\rm inv})^{-1}$, is invertible. It is the linearization of the mixed equation at $\mathcal G_{\mathsf L}(V)$ in the chart, with
the Coulomb condition on the decaying coordinate and the corrected
stationary lifts $K_{P,e}$, followed by
$\Pi_P$; at the solution $\mathcal G_\mathsf L(V)$ the Cauchy--Riemann part of the linearization takes values
in $\Lambda^{0,1}\otimes u^*_{\mathcal G_\mathsf L(V)}T\Msf$, on which $\Pi_P$ is an isomorphism onto
$\Lambda^{0,1}\otimes u_P^*T\Msf$, bounded with bounded inverse. Since this
slice and the Coulomb slice of the solution itself are both transverse to the
gauge orbits, surjectivity of the linearization does not depend on the slice, and
$\mathcal G_{\mathsf L}(V)$ is regular.
\end{proof}

A gauge transformation is \emph{end-even} if it is even on the torus region and on all ends.

\begin{proposition}[convergence into the gluing chart]\label{prop:torus-entry}
Let $\mathsf L_j\to\infty$ and let $U_j$ be index-zero solutions on $W_{\mathsf L_j}$ which, in the
sense of Proposition~\ref{prop:torus-loss}, converge without loss to $V$ (to a pair
$(V_-,V_+)$) in the completed index-zero sets. Then for large $j$ there are end-even gauge transformations $g_j$ with $g_j^*U_j=\Psi_{\mathsf L_j}(z_j)$, $\mathcal F_{\mathsf L_j}(z_j)=0$ and
$\|z_j\|_{X_{\mathsf L_j}}\to0$. In particular $[U_j]=[\mathcal G_{\mathsf L_j}(V)]$ for large $j$.
\end{proposition}

\begin{proof}
We treat the separate case; the joined case is identical, with two pieces, the
bridge in place of the tail and the slices $D_{i,\lambda}$. Let $V^{\rm ad}$, $S_V$,
$\Theta_V$ be as in Lemma~\ref{lem:adapted}, and let $q_{j,e}$ be the limits of $U_j$ at the
mixed ends. Fix $\Theta\ge\Theta_V+2$ and $S\ge S_V+2$, to be taken large, in this order.
Consider the regions
\[
R_E^e=\{\theta\ge\Theta\}\times Y_{e,\mathsf L_j}\ \text{(with the section strip)},\qquad
R_T=T\times[S,\mathsf L_j]\times[-\Theta-1,\Theta+1],
\]
the compact core region $R_K$, consisting of the part of $W_{\mathsf L_j}$ with $r\le S+1$,
$|\theta|\le\Theta+1$ and the original ends up to a fixed $T_0$, and the original ends
$R_O$ beyond $T_0-1$. These regions cover $W_{\mathsf L_j}$; the overlaps
$R_E^e\cap R_T=\{\Theta\le|\theta|\le\Theta+1,\ r\ge S\}$ are bands of unit width and growing
length, and all other overlaps are compact.

\emph{Energy exhaustion.} By Proposition~\ref{prop:torus-loss}, for every
$\varsigma>0$ there are a compact set $K_\varsigma$ and $j_\varsigma$ with the energy of $U_j$ outside
$K_\varsigma$ at most $\varsigma$ for $j\ge j_\varsigma$. The energy density outside the fixed secondary
supports and the regions where the cutoffs $\chi_\pm$ are nonconstant is nonnegative: the curvature density on the torus
region and the kinetic energy on the ends. We write $E_j(\Omega)$ for the energy of
$U_j$ on a region $\Omega$ and $o_j(1)$ for quantities tending to zero as $j\to\infty$ with
$\Theta$ and $S$ fixed.

\emph{The mixed ends.} Replacing $U_j$ by an even gauge transform
(Lemma~\ref{lem:torus-regularity}), we may assume that $U_j$ is smooth. The kinetic energy of
$U_j$ beyond $\Theta-2$ converges to that of $V$, which is at most $Ce^{-2c(\Theta-2-\Theta_V)}$ by
Proposition~\ref{prop:growing-ends}(i); for $\Theta$ large it is below the threshold $e'_m$ of
Proposition~\ref{prop:growing-ends}(ii), which is uniform in $\mathsf L$. In a temporal gauge on
$R_E^e$, normalized by a gauge transformation independent of $\theta$ so that the limit is
exactly $\pi^*B_{q_{j,e}}$, the representative $U_j^E$ satisfies, by~\eqref{eq:end-decay} with
$\theta_0=\Theta-2$, summed over unit intervals with the weight $e^{\delta\theta}$ and $\delta<c$,
\[
\|U^E_j-\pi^*B_{q_{j,e}}\|_{H^m_\delta(R^e_E)}\le Ce^{\delta\Theta}E_j(\{\theta\ge\Theta-2\})^{1/2}
\le Ce^{-(c-\delta)\Theta}+o_j(1).
\]
The slices $Y_{e,\mathsf L_j}$ include the whole torus collar, so this controls the corners
$\{\theta\ge\Theta,\ r\ge S\}$ of the torus region. The same estimate for $V$, the
comparison of the sections at $\theta=\Theta$ and the smooth convergence on compact sets
give $d(q_{j,e},q_e)\le Ce^{-c\Theta}+o_j(1)$, so the coordinates
$p_{j,e}=\exp_{q_e}^{-1}q_{j,e}$ satisfy $|p_j|\to0$ along the diagonal choice below.
On $R^e_E$, where $\beta_e=\varphi_e=1$, the moved reference is
$P^{(p_j)}_{\mathsf L_j}=\pi^*B_{q_{j,e}}+(P_{\mathsf L_j}-\pi^*B_{q_e})+O\bigl(|p_j|\,|P_{\mathsf L_j}-\pi^*B_{q_e}|\bigr)$, the last
term coming from the $E$-corrections, so
\[
\|U^E_j-P^{(p_j)}_{\mathsf L_j}\|_{H^m_\delta(R^e_E)}\le\|U^E_j-\pi^*B_{q_{j,e}}\|_{H^m_\delta(R^e_E)}
+C\sup_\mathsf L\|P_\mathsf L-\pi^*B_{q_e}\|_{H^m_\delta(\{\theta\ge\Theta\})},
\]
and both terms are small for $\Theta$ large, uniformly in $j$ large, by the display and
Lemma~\ref{lem:preglued}(iv).

\emph{The torus strip.} Let $R'_T=T\times[S-2,\mathsf L_j]\times[-\Theta-3,\Theta+3]$. For $\varsigma>0$ let
$K_\varsigma$, $j_\varsigma$ be as in the energy exhaustion; if $S-2$ exceeds the largest value of
$r$ on $K_\varsigma$, then $R'_T$ is disjoint from $K_\varsigma$ for every $\Theta$, so
$E_j(R'_T)\le\varsigma$ for $j\ge j_\varsigma$. We write $\varsigma(S)$ for the least such $\varsigma$; it tends to
zero as $S\to\infty$, independently of $\Theta$. By small-energy regularity at the face and in
the interior (Theorem~\ref{thm:flat-boundary-regularity}(b),
Remark~\ref{rem:interior-boxes}), for $S$ and $j$ large the curvature of $U_j$ and its
covariant derivatives up to order $m$ are small on $R_T$ and there are local gauges
with $\|U_j-\alpha\|_{C^1}\le e_T$. Lemma~\ref{lem:torus-slice} on the simply connected base
$[S-1,\mathsf L_j]\times[-\Theta-2,\Theta+2]$ gives the torus slice $U^T_j=\alpha+a_j$, and
\eqref{eq:torus-slice}, integrated over unit patches with small-energy regularity, gives
\[
\|a_j\|_{H^m(R_T)}\le CE_j(R'_T)^{1/2}\le C\varsigma(S)^{1/2}\qquad(j\ge j_{\varsigma(S)}).
\]
The weight on $R_T$ is at most $e^{\delta(\Theta+1)}$. Since $P_{\mathsf L_j}=\alpha+\chi_{\mathsf L_j}a_V$ on $R_T$,
\[
\|U^T_j-P^{(p_j)}_{\mathsf L_j}\|_{H^m_\delta(R_T)}\le e^{\delta(\Theta+1)}\bigl(C\varsigma(S)^{1/2}
+Ce^{-\gamma_T(S-S_0)}\bigr),
\]
which is small once $S$ is chosen large after $\Theta$.

\emph{The core and the original ends.} On $R_K$, which is compact and contained in
$\{r\le \mathsf L_j/2\}$ for large $j$, there are end-even gauge transformations with
$\|U^K_j-V^{\rm ad}\|_{C^{m+1}(R_K)}=o_j(1)$, by the convergence of
Proposition~\ref{prop:torus-loss}; and $P^{(p_j)}_{\mathsf L_j}-V^{\rm ad}=O(|p_j|)$ there. On $R_O$
the standard exponential decay at nondegenerate ends, with constants uniform for
small end energy \cite[Section~4.1]{Donaldson-book}, \cite{RS-strips}, gives
$\|U^O_j-P_{\mathsf L_j}\|_{H^m_\delta(R_O)}\le\varsigma(T_0)+o_j(1)$ in temporal gauges normalized at the
critical points.

\emph{Overlap logarithms.} On $R^e_E\cap R_T$ both $U^E_j$ and $U^T_j$ are $\alpha$ plus a
term that is small in $C^1$, the first because $B_{q_{j,e}}=\alpha$ on the collar. The
relative gauge transformation is even, so it is close to the even stabilizer of
$\alpha$, which is trivial, and it equals $\exp\varphi_j$ with $\varphi_j$ small. On each fibre
$d_\alpha\varphi_j$ equals the difference of the two torus components up to quadratic
terms, and the base derivatives of $\varphi_j$ equal the difference of the base
components up to quadratic terms; as in the proof of Lemma~\ref{lem:torus-slice},
\[
\|\varphi_j\|_{H^{m+1}_\delta(R^e_E\cap R_T)}\le C\bigl(\|U^E_j-\alpha\|_{H^m_\delta(R^e_E\cap R_T)}
+\|U^T_j-\alpha\|_{H^m_\delta(R^e_E\cap R_T)}\bigr),
\]
with $C$ independent of the length of the band. On the compact overlaps with
$R_K$ and $R_O$ the representatives are close to $V^{\rm ad}$, respectively to the
same stationary or critical configuration, and $V^{\rm ad}$ is irreducible on each
overlap, which contains torus fibres; so the relative gauge transformations are
close to the identity after a choice of central sign, and their logarithms are
small in $H^{m+1}$.

\emph{Patching.} Join $U^T_j$ to the two representatives $U^E_j$ across the two bands by
$\exp(\chi(\pm\theta-\Theta)\varphi_j)$ with a fixed cutoff $\chi$ (the overlap graph is a tree),
then to $U^K_j$ and $U^O_j$ across the compact connected overlaps, each time with one choice of
central sign and a cutoff of fixed slope. The result is an end-even gauge transformation $g_j$,
and since all cutoffs have bounded derivatives and all logarithms are small,
\[
\|x_j\|_{H^m_\delta(W_{\mathsf L_j})}\le\eta(\Theta,S,T_0)+o_j(1),\qquad x_j=g_j^*U_j-P^{(p_j)}_{\mathsf L_j},
\]
where $\eta$ tends to zero when $\Theta$, then $S$, then $T_0$ tend to infinity. A
diagonal choice of $\Theta,S,T_0$ along the sequence gives $\|x_j\|_{H^m_\delta}+|p_j|\to0$.

\emph{Slice and chart.} Gauge transformations preserve the flat-fibre, matching
and Lagrangian conditions, so $g_j^*U_j$ satisfies them. By
Lemma~\ref{lem:torus-quadratic}(a) there are $\xi_j$ with
$\|\xi_j\|_{H^{m+1}_\delta}\le C\|x_j\|_{H^m_\delta}$ and $w_j$ with
$\exp(\xi_j)^*g_j^*U_j=\Psi_{\mathsf L_j}(w_j,p_j)$, $d^*_{A_P}\zeta_j=0$ and
$\|z_j\|\le C(\|x_j\|+|p_j|)\to0$, where $z_j=(w_j,p_j)$; since $U_j$ is a solution,
$\mathcal F_{\mathsf L_j}(z_j)=0$. For large $j$, $\|z_j\|<r$, and Proposition~\ref{prop:torus-existence}
gives $z_j=z_{\mathsf L_j}(V)$.
\end{proof}

\begin{theorem}[torus completion]\label{thm:torus-completion}
Choose the secondary data of the completed domains by
Lemma~\ref{lem:completed-regular} and use them on $W_{\pm,\mathsf R}$ and $W^J_\lambda$. There are
$\mathsf R_0$ and $\lambda_1\ge\lambda_*$ such that for $\mathsf R\ge \mathsf R_0$ and $\lambda\ge\lambda_1$:
\begin{enumerate}[label=(\roman*)]
\item there are no solutions of negative index on $W_{\pm,\mathsf R}$ and $W^J_\lambda$;
\item gluing defines bijections
  \[
  \mathcal G_\mathsf R\colon\cM^0(W_\infty)\to\cM^0(W_\mathsf R),\qquad
  \mathcal G_\lambda\colon\cM^0(W_{-,\infty})\times\cM^0(W_{+,\infty})\to\cM^0(W^J_\lambda)
  \]
  onto sets of regular points, preserving labels, relative classes and markings.
\end{enumerate}
In particular, the index-zero counts satisfy $N_\pm^{(\mathsf R)}=N_{\pm,\infty}$ and
$\mathcal KN^J_\lambda=N_{-,\infty}\otimes N_{+,\infty}$ as matrices.
\end{theorem}

\begin{proof}
\emph{Negative indices.} By Theorem~\ref{thm:torus-period} and the energy identity
there are finitely many negative indices at which solutions can occur for given
labels. If solutions $U_j$ of a fixed negative index $d$ existed on $W_{\mathsf L_j}$ with
$\mathsf L_j\to\infty$, Proposition~\ref{prop:torus-loss} would give
$d=\sum\ind V_\beta+\sum\ind \Gamma_\ell+4n$ with $n\ge0$, completed pieces of nonnegative index
(Lemma~\ref{lem:completed-regular}) and chains of index at least one, which is
impossible.

\emph{The map.} For $\mathsf L\ge\max\{\mathsf L_0,\mathsf L_1\}$, Proposition~\ref{prop:torus-existence}
defines $\mathcal G_{\mathsf L}(V)$ for each representative $V$. It is well defined on classes
and equivariant under the lifting groups. If $V'=h^*V$ for an end-even gauge transformation $h$ of $W_\infty$, trivial at the original ends in the limit, then $V'^{\rm ad}$ and $V^{\rm ad}$ differ by an end-even gauge transformation $h'$ which on $\{r\ge S_V\}$ preserves the torus slice, hence is
constant there with value in the stabilizer of $\alpha$; it extends over the
truncated tail and the collar as the same constant, which fixes $\alpha$, so
$P_{\mathsf L}(V')=h'^*P_{\mathsf L}(V)$, and the chart, the equation and the slice condition are
equivariant, so $\mathcal G_{\mathsf L}(V')=h'^*\mathcal G_{\mathsf L}(V)$. The labels are those of $V$, and
the relative class and bundle marking are those of $P_{\mathsf L}(V)$, which equal those of
$V$ by Lemma~\ref{lem:preglued}(iii) and the construction. In the joined case the
lifting group is $\Phi_J\cong\Phi_-\oplus\Phi_+$, the even stabilizer of $\alpha$ on the bridge is
trivial, and the joined reference of a pair is determined up to $\Phi_J$ by the
classes of $V_\pm$, so no gluing parameter occurs.

\emph{Injectivity and surjectivity.} These are the last two steps of the proof of
Theorem~\ref{thm:neck}. Since $\|z_{\mathsf L}(V)\|\to0$ and $P_\mathsf L(V)=V^{\rm ad}$ on
$\{r\le\mathsf L/2\}$, $\mathcal G_{\mathsf L}(V)$ converges to $V$ on compact subsets; if
$\mathcal G_{\mathsf L_j}(V)$ and $\mathcal G_{\mathsf L_j}(V')$ were gauge equivalent for
$\mathsf L_j\to\infty$, the gauge transformations would be bounded on compact sets and a limit
would give $[V]=[V']$, and the completed sets are finite. If for $\mathsf L_j\to\infty$ there were
index-zero solutions $U_j$ on $W_{\mathsf L_j}$ outside the image, Proposition~\ref{prop:torus-loss},
the index balance and Lemma~\ref{lem:completed-regular} would give a subsequence converging
without loss to some $V$ (a pair) of index zero, and Proposition~\ref{prop:torus-entry} would give
$U_j=\mathcal G_{\mathsf L_j}(V)$ for large $j$.

\emph{Counts.} The coefficient of $N^{(\mathsf R)}_\pm$ at labels $(a,x)$ is the number of
index-zero solutions with these labels, which by (ii) equals that of $N_{\pm,\infty}$.
In the joined case the strip output is at $A_0\cap A_1=(\widetilde P_0\cap\widetilde P_1)\times(R_0\cap R_1)$
and the labels of a joined solution are the pairs of labels of its two
components, so the coefficient of $N^J_\lambda$ at $(a_-\otimes a_+,(x_-,x_+))$ is the
product of the coefficients of $N_{-,\infty}$ and $N_{+,\infty}$; composing with the
identification $\mathcal K$ of generators gives the last identity.
\end{proof}

\subsection{Proof of Theorem~\ref{thm:excision}}

\emph{Order of choices.} (1) Fix the separate mixed maps $N_\pm$ on $X_{D,\pm}$.
(2) Fix $\lambda_0\ge\lambda_*$, the data of $N_Z$ at width $\lambda_0$ and the data of
$\overline{\mathsf E}$ (Proposition~\ref{prop:closed-cut}), and regular joined data $\mathfrak d_J$
at width $\lambda_0$. (3) Choose regular data on the completed domains
(Lemma~\ref{lem:completed-regular}). (4) Choose $\mathsf R\ge \mathsf R_0$ and $\lambda\ge\max(\lambda_1,\lambda_0)$
(Theorem~\ref{thm:torus-completion}). (5) At these finite lengths, add an invariant
secondary perturbation supported in the cores, regular in index one and small
enough that the finite regular index-zero sets and the absence of negative indices
persist; by compactness of the index-zero sets and the implicit function theorem
the index-zero counts are unchanged.

By Proposition~\ref{prop:closed-cut}, which contains Lemma~\ref{lem:closed-cut-data};
by Proposition~\ref{prop:width} with $a=\lambda_0$ and $b=\lambda$, applied to the
data $\mathfrak d_J$ changed by a small generic invariant secondary perturbation that
is regular at both widths (Theorem~\ref{thm:regularity}), followed, at width
$\lambda$, by Theorem~\ref{thm:family} for the compactly supported change from these
data to those of steps (3) and (5), both regular; and by
Theorem~\ref{thm:torus-completion},
\[
[\mathcal KN_Z][\overline{\mathsf E}]=[\mathcal KN^J_{\lambda_0}]=[\mathcal KN^J_\lambda]=[N_{-,\infty}]\otimes[N_{+,\infty}]
=[N_-^{(\mathsf R)}]\otimes[N_+^{(\mathsf R)}]=[N_-]\otimes[N_+],
\]
the last step by Proposition~\ref{prop:collar-length}, applied to the data changed by a
small generic invariant secondary perturbation that is regular at both lengths $0$ and
$\mathsf R$ (Theorem~\ref{thm:regularity}), and Theorem~\ref{thm:family} for the compactly
supported changes of secondary data between $N_\pm$, these data and $N^{(\mathsf R)}_\pm$.
All structures on the section are split, so no change of the section structure
occurs. Since $[\overline{\mathsf E}]=[\mathsf E]^{-1}$, the theorem follows. \qed

\section{One- and three-handles}\label{sec:one-three}

We keep the notation of Section~\ref{subsec:folded}.

\subsection{The relative invariant of the one-handle cobordism}\label{subsec:one-handle-invariant}

Let $W_k\colon M\to M^+$ be the trace of $k$ one-handles attached away from the
tube, and let $B_k=(W_k\setminus\operatorname{int}\nu(\gamma))\natural(T^2\times D^2)$, where
$\nu(\gamma)$ is the tube; the two end balls of the tube turn the boundary of the
first factor into $(-M)\#M^+$, and its boundary sum with the torus filling gives
$Z$.

\begin{lemma}\label{lem:Bk-topology}
\begin{enumerate}[label=(\roman*)]
\item $\mathscr W\cup_{Z_0}F_0\cup_{Z_1}F_1\cong B_k$ by an orientation-preserving
  diffeomorphism covered by a bundle isomorphism, where $\mathscr W$ is the gauge
  triangle assembled from $D_0$, $C$, $D_1$ over the two surface sheets
  $\Delta\times\Sigma_h$, $\Delta\times\Sigma_g$ and $F_i=V_i\times[0,1]$ fill its auxiliary boundaries
  $Z_i=V_i\cup(-V_i)$. The diffeomorphism restricts on the boundary to the
  identification $Z=D_0\cup(-D_1)\cong(-M)\#M^+\#T^3$ of Section~\ref{subsec:folded},
  carries the complement of the torus parts onto the complement of the tube
  $\nu(\gamma)$ in $W_k$, and the lifting group of the union is zero, as is
  that of $B_k$.
\item $B_k\cup_Z\mathcal E\cong W_k^\#$, with the incoming end $-M$ read as outgoing,
  and the lifting group of the union is $\phi_{W_k}$.
\end{enumerate}
\end{lemma}

\begin{proof}
\begin{figure}[ht]
\centering
\begin{tikzpicture}[scale=2.1,font=\small]
\fill[gray!12] (0,0) circle (1);
\fill[orange!18] (170:1) arc (170:250:1) -- cycle;
\fill[orange!30] (50:1) arc (50:130:1) -- cycle;
\fill[orange!18] (290:1) arc (290:370:1) -- cycle;
\fill[blue!15,opacity=0.8] (100:0.6) arc (100:190:0.6) -- (190:0.92) arc (190:100:0.92) -- cycle;
\fill[blue!15,opacity=0.8] (-10:0.6) arc (-10:80:0.6) -- (80:0.92) arc (80:-10:0.92) -- cycle;
\draw[thick] (0,0) circle (1);
\draw[thick] (170:1) -- (250:1);
\draw[thick] (50:1) -- (130:1);
\draw[thick] (290:1) -- (10:1);
\node at (210:0.88) {$G_0$};
\node at (90:0.9) {$G_m$};
\node at (330:0.88) {$G_1$};
\node at (0,-0.28) {$H$};
\node at (225:0.6) {$c_0$};
\node at (270:0.08) {};
\node at (72:0.62) {$c_m$};
\node at (318:0.62) {$c_1$};
\node at (145:0.76) {$\Omega_0$};
\node at (35:0.76) {$\Omega_1$};
\node at (150:1.13) {$a_{0m}$};
\node at (30:1.13) {$a_{m1}$};
\node at (270:1.13) {$a_{10}$};
\end{tikzpicture}
\caption{The planar base of the proof of Lemma~\ref{lem:Bk-topology}(i). The fibre is
$\Sigma_h$ over the hexagon $H$, $\Sigma_g$ over $G_m$ and the torus $T$ over $G_0$ and $G_1$;
over the chords $c_0$, $c_m$, $c_1$ the pieces $U_0$, $C$, $U_1$ are inserted. Over the
bands $\Omega_0$ and $\Omega_1$ the four-manifold is $F_0=[0,1]\times V_0$ and
$F_1=[0,1]\times V_1$; the boundary over $\partial\mathbb D$ is $Z$.}
\label{fig:planar-base}
\end{figure}
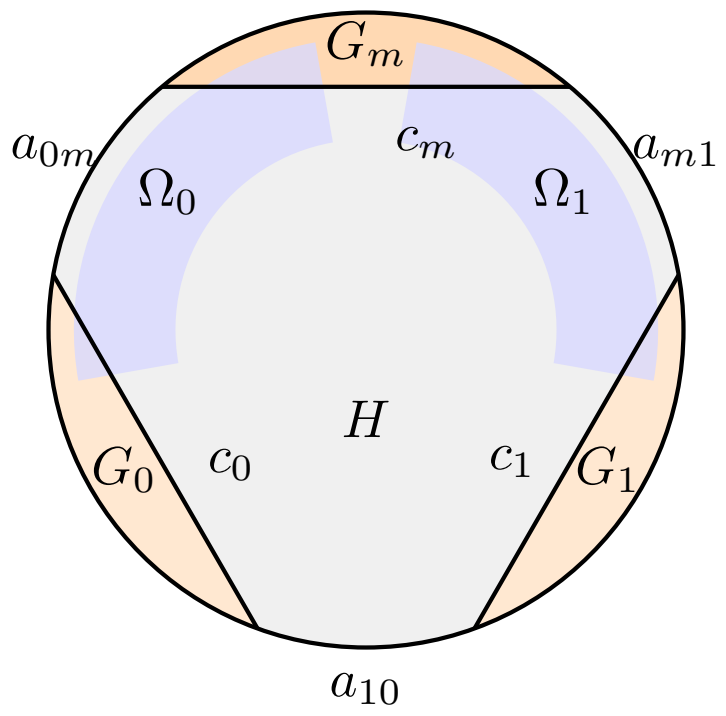

(i) We use four-manifolds built from a planar base. Let $\mathbb D$ be a closed disk
and $c_0,c_m,c_1\subset\mathbb D$ three disjoint chords; they cut $\mathbb D$ into three outer
half-disks $G_0,G_m,G_1$ (with $\partial G_e\setminus\partial\mathbb D=c_e$) and a hexagon $H$, whose boundary
consists alternately of $c_0$, an arc $a_{0m}\subset\partial\mathbb D$, $c_m$, an arc $a_{m1}$, $c_1$ and an
arc $a_{10}$. Put
\[
\mathfrak X=H\times\Sigma_h\ \cup\ G_m\times\Sigma_g\ \cup\ G_0\times T\ \cup\ G_1\times T\ \cup\
c_0\times U_0\ \cup\ c_m\times C\ \cup\ c_1\times U_1 ,
\]
where $c_e\times Y$ is glued to the two adjacent regions along $c_e\times\partial Y$, with
$\partial U_i=\Sigma_h\sqcup T$ and $\partial C=\Sigma_h\sqcup\Sigma_g$, corners rounded in product collars, and
the bundle given by the bundles of $U_i$, $C$ and their boundary restrictions. The
boundary $\partial\mathfrak X$ lies over $\partial\mathbb D$; reading it along $\partial\mathbb D$ from $G_0$ through $a_{0m}$,
$G_m$, $a_{m1}$, $G_1$ and $a_{10}$ gives
$(T\times I)\cup U_0\cup_{\Sigma_h}C\cup_{\Sigma_g}(-C)\cup_{\Sigma_h}(-U_1)\cup_T(T\times I)\cup U_1\cup_{\Sigma_h}(-U_0)$, that is
$V_0\cup(-V_1)$ and $(-U_0)\cup U_1$ glued along two tori, which is
$Z=D_0\cup(-D_1)$.

\emph{Step 1: $\mathscr W\cup F_0\cup F_1\cong\mathfrak X$.} Let $\Omega_0\subset\operatorname{int}\mathbb D$ be a
closed rectangle $[0,1]_y\times[0,1]_x$ near $a_{0m}$ that crosses, in the $x$-direction and
in this order, $G_0$, $c_0$, $H$, $c_m$ and $G_m$, with $c_0\cap\Omega_0=[0,1]\times\{\tfrac14\}$
and $c_m\cap\Omega_0=[0,1]\times\{\tfrac34\}$; let $\Omega_1$ be the analogous rectangle near $a_{m1}$
crossing $G_m$, $c_m$, $H$, $c_1$, $G_1$. Over $\Omega_0$ the fibre across $x$ is
$T$-collar, $U_0$, $\Sigma_h$-collar, $C$, $\Sigma_g$-collar, so
$\mathfrak X|_{\Omega_0}\cong[0,1]\times(U_0\cup_{\Sigma_h}C)=[0,1]\times V_0=F_0$, and likewise
$\mathfrak X|_{\Omega_1}\cong F_1$; the boundary $\partial\Omega_i$ carries $Z_i$. On
$\mathbb P=\mathbb D\setminus\operatorname{int}(\Omega_0\cup\Omega_1)$, a disk with two holes, the chords are cut
into the following arcs: two arcs of $c_0$ from $\partial\mathbb D$ to $\partial\Omega_0$, a short one near
$a_{0m}$ and a long one; three arcs of $c_m$, from $\partial\mathbb D$ to $\partial\Omega_0$, from $\partial\Omega_0$ to
$\partial\Omega_1$, and from $\partial\Omega_1$ to $\partial\mathbb D$; and two arcs of $c_1$ from $\partial\Omega_1$ to $\partial\mathbb D$.
Let $\mathcal N_0\subset\mathbb P$ be a band from $\partial\mathbb D$ to $\partial\Omega_0$ containing the long and the short
arc of $c_0$, the region $G_0\setminus\Omega_0$ between them, the thin region of $H$ between
the short arc of $c_0$ and the first arc of $c_m$, and that arc. Across $\mathcal N_0$ the fibre
is a $\Sigma_h$-collar, $U_0$ traversed from $\Sigma_h$ to $T$, a $T$-collar, $U_0$ traversed from
$T$ to $\Sigma_h$, a $\Sigma_h$-collar, $C$ and a $\Sigma_g$-collar, that is
$(-U_0)\cup_TU_0\cup_{\Sigma_h}C=(-U_0)\cup_TV_0=D_0$ with collars of its boundary; so
$\mathfrak X|_{\mathcal N_0}\cong e_0\times D_0$. Likewise a band $\mathcal N_1$ from $\partial\Omega_1$ to $\partial\mathbb D$ gives
$e_1\times D_1$, and a neighbourhood $\mathcal N_m$ of the middle arc of $c_m$ gives $e_m\times C$. The
complement of $\mathcal N_0\cup \mathcal N_m\cup \mathcal N_1$ in $\mathbb P$ consists of two hexagons, the rest of $H$
with fibre $\Sigma_h$ and the rest of $G_m$ with fibre $\Sigma_g$, each with its three edges on
$\mathcal N_0$, $\mathcal N_m$, $\mathcal N_1$ and its three vertex arcs on $\partial\mathbb D$, $\partial\Omega_0$, $\partial\Omega_1$. This is the
assembly of $\mathscr W$ over the two sheets, with the output vertex at $\partial\mathbb D$ and the
inputs at $\partial\Omega_0$, $\partial\Omega_1$. The identifications are the identity on the pieces
$U_i$, $C$ and on the surfaces, so the bundles correspond.

\emph{Step 2: the Heegaard pieces.} Write $M=H_0\cup_\Sigma H_1$ for the Heegaard
splitting of genus $h-1$ of the presentation, and $U_i=H_i\cup_{D}([0,1]\times T)$, where
the sum disk $D\subset\partial H_i$ is identified with a disk $D\subset\{0\}\times T$; by
Definition~\ref{def:presentation}(ii), $D$ is the equatorial disk of the ball
$b(B^3)$, and we write $o'$ for its centre. So
$\Sigma_h=(\Sigma\setminus D)\cup_c T^\circ$, $T^\circ=T\setminus D$, $c=\partial D$. The compression body $C$ is
$C''\cup_{c\times[0,1]}(T^\circ\times[0,1])$, where $C''$ is obtained from $(\Sigma\setminus D)\times[0,1]$
by attaching $k$ three-dimensional one-handles to $(\Sigma\setminus D)\times\{1\}$ away from
$\partial D$, and $\Sigma_g=(\Sigma^+\setminus D)\cup_cT^\circ$ with $\Sigma^+$ the surgered surface. Accordingly $\mathfrak X=\mathfrak X_M\cup\mathfrak X_T$, where $\mathfrak X_T$
collects the torus parts: $T^\circ$ over $H$ and $G_m$, $T^\circ\times[0,1]$ over $c_m$, the
factor $[0,1]\times T$ of $U_i$ over $c_i$, and $T$ over $G_i$; and $\mathfrak X_M$ collects the
rest: $\Sigma\setminus D$ over $H$, $\Sigma^+\setminus D$ over $G_m$, $C''$ over $c_m$ and $H_i$ over $c_i$. They
meet along $H'\times c\cup(c_0\sqcup c_1)\times D$, where $H'=H\cup G_m\cup(\text{band of }c_m)$.

\emph{Step 3: $\mathfrak X_T$.} Thicken each $c_i$ into a band $b_i=c_i\times[0,1]$, the
factor $[0,1]$ being the interval factor of $[0,1]\times T\subset U_i$, and similarly $c_m$.
Then $\mathbb D'=H'\cup b_0\cup G_0\cup b_1\cup G_1$ is a disk, $H'$ is a band across it separating
$Q_i=b_i\cup G_i$, and the fibres are $T^\circ$ over $H'$ and $T$ over $Q_0\cup Q_1$. Hence
\[
\mathfrak X_T=(\mathbb D'\times T)\setminus\operatorname{int}(H'\times D),
\]
and $H'\times D$ is a tubular neighbourhood of the arc $\beta'=\beta_0\times\{o'\}$, where $\beta_0$
crosses the band $H'$ from $a_{10}$ to the opposite boundary arc and $o'\in D$ is the
centre; the interface with $\mathfrak X_M$ is its lateral boundary
$H'\times c\cup(c_0\sqcup c_1)\times D\cong S^2\times[0,1]$. The arc $\beta'$ is a chord of the disk
$\mathbb D'\times\{o'\}$, so it cobounds a half-disk $\delta$ with an arc of $\partial\mathbb D'\times\{o'\}$. A
regular neighbourhood $B$ of $\delta$ is a four-ball meeting $\partial(\mathbb D'\times T)$ in a
three-ball, and we may take $H'\times D\subset B$ with $\nu(\beta')\cap\partial B=\nu(\beta')\cap\partial(\mathbb D'\times T)$.
Then $(\mathbb D'\times T)\setminus\operatorname{int}B\cong D^2\times T^2$, and $B\setminus\operatorname{int}\nu(\beta')$ is
the complement of a trivial arc in a four-ball: in the model
$B=D^3\times[0,1]$, $\beta'=\{0\}\times[0,1]$, $\nu(\beta')=D^3_{1/2}\times[0,1]$, it is
$S^2\times[\tfrac12,1]\times[0,1]$, with the lateral boundary $S^2\times\{\tfrac12\}\times[0,1]$. Hence
\[
\mathfrak X_T\cong(T^2\times D^2)\ \natural\ \bigl(S^2\times[\tfrac12,1]\times[0,1]\bigr),
\]
the boundary sum being taken along the three-ball $\partial B\cap\operatorname{int}(\mathbb D'\times T)$,
and the interface with $\mathfrak X_M$ is $S^2\times\{\tfrac12\}\times[0,1]$.

\emph{Step 4: $\mathfrak X_M$.} Refill the disk: $\mathfrak X'_M=\mathfrak X_M\cup(H'\times D)\cup
(c_0\times \mathcal V_0\cup c_1\times \mathcal V_1)$, where $\mathcal V_i\cong D\times[0,1)$ is a collar of $D$ in $H_i$ removed
from $\mathfrak X_M$ beforehand (which does not change its diffeomorphism type). Then
$\mathfrak X'_M=H\times\Sigma\cup c_0\times H_0\cup c_1\times H_1\cup c_m\times C'\cup G_m\times\Sigma^+$, with
$C'=C''\cup D\times[0,1]$ the compression body from $\Sigma$ to $\Sigma^+$, and
$\nu=H'\times D\cup c_0\times \mathcal V_0\cup c_1\times \mathcal V_1$ is a tubular neighbourhood of $\beta'$. Identify
$H$ with $[-1,1]_u\times[0,1]_t$ so that $c_0=\{u=-1\}$, $c_1=\{u=1\}$, $a_{10}=\{t=0\}$ and
$c_m=\{t=1\}\times[u_1,u_2]$ with $0<u_1<u_2<1$. Then $H\times\Sigma\cup c_0\times H_0\cup c_1\times H_1=[0,1]_t\times(H_0\cup[-1,1]_u\times\Sigma
\cup H_1)=[0,1]\times M$. The part over $c_m$ and $G_m$ is, after absorbing the collar
$G_m\times\Sigma^+$ of the upper boundary, $[u_1,u_2]\times C'$ attached along $[u_1,u_2]\times\Sigma=\{1\}\times[u_1,u_2]\times\Sigma\subset\{1\}\times M$. As
$C'=\Sigma\times[0,1]\cup h_1\cup\dots\cup h_k$, this is a collar on $\{1\}\times M$ together with
the pieces $[u_1,u_2]\times h_j=[u_1,u_2]\times D^1\times D^2\cong D^1\times D^3$ attached along
$\partial D^1\times(D^2\times[u_1,u_2])$, two three-balls in the top boundary straddling
$\Sigma$: these are four-dimensional one-handles. Hence $\mathfrak X'_M\cong W_k$, with
incoming end $-M$ over $a_{10}$ and outgoing end $M^+$ over the opposite arc. Choosing
$\beta_0=\{0\}\times[0,1]$, which ends on $a_{0m}$, the arc $\beta'$ becomes $\gamma=[0,1]\times\{q\}$
with $q=(0,o')\in[-1,1]\times\Sigma\subset M$, the centre of $b(B^3)$, which is the product
path, and $\nu$ becomes the tube $\nu(\gamma)$ of the decorated cobordism; the
one-handles are attached away from it. So
$\mathfrak X_M\cong W_k\setminus\operatorname{int}\nu(\gamma)$, and the interface with $\mathfrak X_T$ is the lateral
boundary $S^2\times[0,1]$ of the tube.

\emph{Step 5: assembly.} By Steps 3 and 4, $\mathfrak X$ is obtained from
$W_k\setminus\operatorname{int}\nu(\gamma)$ by attaching $S^2\times[\tfrac12,1]\times[0,1]$ along
$S^2\times\{\tfrac12\}\times[0,1]$, which is a collar and does not change the diffeomorphism
type, boundary summed with $T^2\times D^2$ at a ball in the lateral part of the
boundary. Since $\partial(W_k\setminus\operatorname{int}\nu(\gamma))=(-M)\#M^+$ is connected, the boundary
sum does not depend on the ball, and $\mathfrak X\cong B_k$. With Step 1 this proves the
diffeomorphism. The bundle of $\mathfrak X$ is trivial on $\mathfrak X_M$ and the pullback of
the odd bundle $P_T$ under $\mathfrak X_T\subset\mathbb D'\times T\to T$ on $\mathfrak X_T$, glued along the interface
by the fixed trivialization of $P_T$ over $D$; the summand $S^2\times[\tfrac12,1]\times[0,1]$ lies
over $D$, where $P_T$ is trivialized. So the diffeomorphism is covered by an
isomorphism with the bundle of $B_k$, trivial on $W_k\setminus\nu(\gamma)$ and pulled back
from $T^2$ on $T^2\times D^2$, boundary summed along the trivialization over the ball.
On the boundary, the pieces over $a_{10}$ and over the opposite arc are
$-M\setminus B$ and $M^+\setminus B$, and the torus part over $\partial\mathbb D'$ is $T^3$ minus the two
balls at the ends of $\nu(\beta')$, which are the connected-sum balls, in the
identification of Section~\ref{subsec:folded}; the cut tori $T_a$, $T_b$ of
Lemma~\ref{lem:closed-cut-embedding} are the fibres over points of the arcs
$\partial G_0\cap\partial\mathbb D$ and $\partial G_1\cap\partial\mathbb D$.

For the lifting groups: restriction
$H^1(\mathscr W\cup F_0\cup F_1;\FF)\to H^1(\mathscr W)\oplus H^1(F_0)\oplus H^1(F_1)$ is injective, since
$H^0(\mathscr W)\oplus H^0(F_0)\oplus H^0(F_1)\to H^0(Z_0)\oplus H^0(Z_1)$, $(a,b,c)\mapsto(a+b,a+c)$, is
onto. The caps carry the determinant-one group, so a lifting class of the
union vanishes on $Z_0$ and $Z_1$. A lifting class of $\mathscr W$ vanishing on the surface
sheets restricts to zero on $C$ and on each $D_i$: $H^1$ of $C$, of $U_i$ and of $V_i$
injects into $H^1$ of $\Sigma_g$, of $\Sigma_h$ and of $\Sigma_g$ respectively, and since $T$
separates $D_i=(-U_i)\cup_TV_i$, $H^1(D_i)$ injects into $H^1(U_i)\oplus H^1(V_i)$. In addition, the
incidence graph of $\mathscr W$ (two sheets, three edge pieces and six contact
regions) has first Betti number two, and the two corresponding classes restrict
to zero on every piece. These two classes restrict nontrivially and
independently to $Z_0$ and $Z_1$, each being the dual of the surface sheet in the
vertex slice. A lifting class vanishing on the sheets is therefore determined by its
restrictions to $Z_0$ and $Z_1$, where the neck-torus duals and the incidence
classes are nonzero, and so it vanishes. The
lifting group of the union is therefore zero, which is the group of $B_k$.
(ii) Model $\mathcal E$ as $T^2\times P$ with the summand products inserted along two
arcs. Attaching the disk base of $T^2\times D^2$ in $B_k$ caps the input circle $c$
of the pair of pants and restores the product annulus, with the identity
parametrizations at both exterior circles and no Dehn twist; across the cap the
two sum collars join to the lateral collar of the original tube, with the
prescribed normal frame. The complement of the tube and its bundle are unchanged.
The cut $Z$ is connected, so the lifting group of the union consists of the
classes vanishing on $B_k$ and restricting into $\Phi_{\mathcal E}$; compatibility on $Z$
forces $s=0$ and $s_-=s_+$, which is the torus lifting group $\phi_{W_k}$.
\end{proof}

\begin{lemma}\label{lem:handle-state}
$[\mathsf E](\Theta_{B_k})=\mathrm{ten}(\Ish(W_k))$ on homology over $\FF$, where
$\Theta_{B_k}=[\mu_{\mathscr W}](\theta_0,\theta_1)$ is the relative invariant of $B_k$,
$\mu_{\mathscr W}$ is the gauge triangle map of $\mathscr W$ modulo the gauge group enlarged
by $\Phi$,
$\theta_i\in C_G(Z_i)$ is the relative invariant of the filling $F_i$, defined as $\theta_G$ in
Section~\ref{sec:two-handles}, and $\mathrm{ten}(F)=\sum_ia_i\otimes b_i$
for $F(v)=\sum_i\langle a_i,v\rangle b_i$, the pairing being the gauge Kronecker
pairing between $-M$ and $M$.
\end{lemma}

\begin{proof}
\emph{The relative invariant.} By Lemma~\ref{lem:Bk-topology}(i) and
\cite[Proposition~5.4]{KM-unknot} with two caps, whose connecting images vanish
because the auxiliary cuts are connected and meet the triangle, together with the
projection argument of Proposition~\ref{prop:gauge-triangle} for each input,
$\Theta_{B_k}=[\mu_{\mathscr W}](\theta_0,\theta_1)$.

\emph{Gluing.} For a long neck along $Z$ with nondegenerate data, the index-zero
moduli spaces of $B_k\cup_Z\mathcal E$ are, by Theorem~\ref{thm:neck}, the fibre
products over the critical points $z_\sigma$ of $Z$ (with their two even lifts
$\sigma$) of those of $B_k$ and of $\mathcal E$; the gluing parameter is a point since the
determinant-one stabilizer $\{\pm1\}$ extends over both pieces and the connecting
image vanishes.

\emph{Counting.} Let $e_{\sigma\tau}$ be the number of index-zero solutions on $\mathcal E$
from $z_\sigma$ to the output lifts $(a_0,b_\tau)$ for fixed labels. Restriction of
$\Phi_{\mathcal E}$ to the lifts of $z$ and of the $-M$ output is an isomorphism, so the
coefficient of $E([z])$ at $[a]\otimes[b]$ is $e_{00}+e_{01}$. The element of
$\Phi_{\mathcal E}$ acting on the $z$ and $M^+$ lifts gives a bijection between
solutions from $z_1$ to $b_\tau$ and from $z_0$ to $b_{\tau+1}$, so
$e_{1\tau}=e_{0,\tau+1}$. If $\beta_\sigma$ are the coefficients of $\Theta_{B_k}$ at $z_\sigma$,
the glued count at $(a_0,\cdot)$, summed over the output lifts, is
$\sum_\tau\sum_{z,\sigma}\beta_\sigma e_{\sigma\tau}=\sum_z(\beta_0+\beta_1)(e_{00}+e_{01})$, the coefficient of $\mathsf E(\Theta_{B_k})$. Each orbit of lifts is counted once, because the lifting group acts freely and the counts are taken in the quotient.

\emph{Reversal.} By Lemma~\ref{lem:Bk-topology}(ii) the union is $W_k^\#$ with the
incoming end read as outgoing; the same moduli spaces, with the incoming lift
fixed and the output lifts summed, compute the framed coefficient
$\langle\Ish(W_k)[a],[b]\rangle$ \cite[Section~7.2]{Scaduto}, the pairing between
$C_G((-M)^\#)$ and $C_G(M^\#)$ being the reversal pairing
\cite[Section~7.4]{Scaduto}. Independence of data holds on homology by
functoriality under continuation.
\end{proof}

\subsection{Reflection}\label{subsec:reflection}

For an object $(M,b)$ with equatorial disk $D=\{x_3=0\}$ in the ball, put
$\bar b=b\circ r_B$ with $r_B(x_1,x_2,x_3)=(x_1,x_2,-x_3)$; this is an
orientation-preserving parametrization of a ball in $-M$ that fixes $D$ and
exchanges the two half-balls. On the torus $T^3=T^2\times S^1_c$ use
$\rho_T(x,y,c)=(x,y,-c)$, which reverses orientation, preserves the odd two-torus
with its bundle, the class $\alpha_T$ and the flat metric. Transporting the torus
attachment by $r_B$ and $\rho_T$ defines the attachment of the \emph{reflected
object} $(-M,\bar b)$, with $(-M)^\#=-(M^\#)$ as marked bundles; reflecting twice
is the identity. For a presentation $p=(C_0,C_1)$ put $\bar p=(C_1,C_0)$, with all
spatial data transported; ordinary strip data are reflected by
$u(s,t)\mapsto u(-s,1-t)$, $J\mapsto J(-s,1-t)$, $H\mapsto-H(-s,1-t)$, with the
full one-form transformation.

Reflection gives pairings, which are identifications of complexes for the reflected
end perturbation and hold up to gauge continuation, hence on homology, for
independently chosen end and secondary data; only the homology statement is used:
\[
C_G((-M)^\#,\bar p)=C_G(M^\#,p)^\vee,\qquad
\CF(\bar p)=\CF(p)^\vee,
\]
under which the gauge and ordinary continuation maps are transposed:
$\mathrm{cont}_G(\bar P,\bar Q)=\mathrm{cont}_G(Q,P)^\vee$ and $\mathrm{cont}_S(\bar P,\bar Q)=\mathrm{cont}_S(Q,P)^\vee$. The mixed
completion of a reflected presentation is \emph{not} reflected: mixed solutions
are not carried to mixed solutions. We choose independent regular invariant data (end perturbation and secondary perturbation) defining $N_{\bar p}$; by Theorem~\ref{thm:mixed-iso}(ii) its homology
class does not depend on the choice, and in particular
$[N_{\bar{\bar P}}]=[N_P]$.

\subsection{Continuations and composition}\label{subsec:continuations}

The Lagrangians $\widetilde P_i$ and $R_i$ are simply connected with $\pi_2=0$.
The correspondence $L_C$ is an $\SU(2)^k$-bundle over $\Msf_h$, simply connected
with $\pi_2(L_C)=\ZZ$; spheres in $L_C$ have zero area and zero Chern number, so
every pair and triangle used below is monotone with minimal Maslov number $4$.

\emph{Continuation between primary data.} For two admissible choices of primary data
with the same splitting and metric data, the Lagrangians differ by Hamiltonian isotopies.
Choose Hamiltonian isotopies on the two boundary components and extend them
over disjoint boundary collars. Their interpolation in the strip interior
gives a domain-dependent Hamiltonian change of variables which straightens
both labels, with fixed data on the far ends. The transformed equation
has a Hamiltonian one-form $K=K_s\,ds+K_t\,dt$; its tangential boundary
restriction need not vanish. With $\iota_{X_H}\omega=-dH$, the equation is
\[
v_s-X_{K_s}+J(v)(v_t-X_{K_t})=0.
\]
On a truncated strip its full analytic energy is
\[
E_K(v)=2\int v^*\omega-2\int_{\partial S}v^*K
 +2\int_S\bigl(\partial_sK_t-\partial_tK_s
                   +\omega(X_{K_s},X_{K_t})\bigr)(v)\,ds\,dt.
\]
Combining the cut-end terms with fixed reference areas gives the prescribed
chord-action contributions. The physical
boundary work and curvature term have uniform bounds on compact families
of these data. Monotonicity therefore bounds energy at fixed index.
On an inserted autonomous neck the curvature vanishes and a stationary
chord has zero covariant energy, so the bound is independent of neck length.

One may instead add a Hamiltonian boundary connection making the moving
labels parallel before straightening. The resulting fixed-boundary
perturbation has zero tangential boundary restriction. Interpolating the
connections with fixed ends retains the preceding energy bound. Generic
relative perturbations, the exclusion of bubbles (each costs index at least four), and
ordinary strip gluing give a chain homotopy between the two constructions. This
identifies $\mathrm{cont}_S(P,P')$ with Wehrheim and Woodward's data-change
map on fixed Lagrangians. It is a chain map, equals
the identity for constant data, is compatible with composition and is
invariant under homotopies of the transformed data with fixed ends
\cite[Theorems~3.9, 3.13, 4.1, 4.2 and Example~2.6]{WW-relative}. Both boundary
Lagrangians may move. Paths of Lagrangians and their homotopies are taken among graphs of
differentials of functions over the unperturbed Lagrangians, with the functions in a small
convex set; so the continuations used here depend only on their
endpoints. The functions of the Lagrangians in the neighbourhoods $\mathcal N(t)$ of
Definition~\ref{def:admissible-nbhd} form such a set, although the admissible sets of primary
data defined there are not convex; the interpolating paths of Lagrangians may leave
$\mathcal N(t)$, which does not affect the continuation, and the gauge continuations run in the
convex space $\operatorname{span}F$.

\emph{The auxiliary pairs.} We use the actual product-target pairs
$(A_0,L_C)$ and $(L_C,A_1)$ in $\mathbb X$. Their intersection points identify
with those of $(R_0,\widetilde P_0\circ L_C)$ and
$(\widetilde P_1\circ L_C,R_1)$, respectively, and hence, by continuation, with the critical
sets of the height functions of Section~\ref{subsec:folded}. This is
an identification of generators. For a general compression it is not
an unfolding identification of strip moduli spaces; exact diagonal
unfolding is used only when $k=0$. For each ordered pair,
let $\mathcal A_i$ be its critical action: the values at intersections
of its local exact-graph primitive, with the sign determined by the
order of the pair. Set
\[
4\pi^2\kappa_{\mathrm{aux}}(u)
 =\int u^*\omega_{\mathbb X}-\bigl(\mathcal A_i(a)-\mathcal A_i(x)\bigr).
\]
At the flat height data $\mathcal A_i$ is the corresponding signed
height function $t h_i$. Under small primary changes it is the
actual graph primitive, including its higher-order terms. Continue the
normalization from the flat reference. Along a critical branch the graph
action and the perturbed gauge action then have the same derivative,
given by the varying perturbation evaluated on that branch. Their
differences at the generators agree; a common additive constant in the
primitive has no effect on these differences.
On $\CF_{\mathbb X}(A_0,L_C)$ and $\CF_{\mathbb X}(L_C,A_1)$ define the relative grading $m$ by
$m(a)-m(x)=\ind(u)-8\kappa_{\mathrm{aux}}(u)$ for any strip class $u$ from $a$ to $x$; it is
independent of $u$ by monotonicity. For $t$ small it coincides with the Morse
index of the height function (for $i=0$) or with the Morse index of its
negative (for $i=1$), and there is a unique generator $q_i^*$ of maximal $m$.

\begin{lemma}[Regular auxiliary strips and triangles]
\label{lem:auxiliary-split-regularity}
Fix a lower almost complex path, regular for the ordinary strips, and an upper path near
the two strip boundaries, leaving a nonempty free interval between them. There are split
almost complex structures with zero Hamiltonian, extending these choices, that are regular
for the strips of $(A_0,L_C)$, $(L_C,A_1)$ and $(A_0,A_1)$ and for the triangle with these
boundary conditions. The relative parameter families used below are regular for generic
such choices, with only the usual breaking and parameter-end strata in dimensions zero and
one. For these choices the generators $q_i^*$ are cycles.
\end{lemma}
\begin{proof}
Put $M=\Msf_h$, $N=\Msf_g$, and write
$J_{\mathbb X}=(-J_M)\oplus J_N$. The upper projection $L_C\to N$ is an
embedding: restriction to the upper surface determines the
representation of the compression group. In particular,
\begin{equation}\label{eq:upper-trace-injective}
 TL_C\cap(TM\oplus\{0\})=0.
\end{equation}
For a strip $u=(u_M,u_N)$, constancy of $u_N$ forces the lower
trace along its $L_C$ boundary to be constant. Extend $u_M$ by
that constant across a small interval of this boundary. Equal
traces give a weak $W^{1,p}$ solution, $p>2$, for smoothly extended
coefficients. Interior regularity and unique continuation make
$u_M$ constant. Hence a nonconstant strip has nonconstant upper
component, even if its lower component is constant.

For the upper component, the ordinary strip injective-point
argument gives a dense open set of points where
$\partial_su_N\ne0$ and the value $u_N(s,t)$ is neither an end
value nor another value $u_N(s',t)$ with $s'\ne s$. The
componentwise application is given in the proof of
\cite[Theorem~3.2, Step~1a]{WW-constant}. It uses finite energy,
the entire $s$-axis and convergence to the fixed end values;
we do not invoke that theorem's later generic-Hamiltonian
argument for arbitrary constant patches.

Let $\eta=(\eta_M,\eta_N)$ annihilate the linearization together
with the allowed $J_N(t)$ variations. If $\eta_N$ is nonzero,
adjoint unique continuation gives a nonzero point in the free
$t$-interval, and one may choose it in the preceding injective
set. Compatible variations $Y$ send the nonzero vector $u_{N,t}$
onto the full target tangent space. A variation supported near
$(t,u_N(s,t))$ consequently has nonzero pairing with $\eta_N$,
contrary to the annihilator condition. Thus $\eta_N=0$.
Green's formula identifies the remaining adjoint boundary trace
with an element of $TL_C$ whose upper component is zero.
Equation~\eqref{eq:upper-trace-injective} makes the lower trace
zero. Extending it by zero across that boundary interval and
using adjoint regularity and unique continuation gives
$\eta_M=0$. This proves full universal surjectivity without
varying $J_M$. Reversing the strip boundary order gives the
other auxiliary pair.

At a constant strip $x$, the operator is $\partial_s+B_x$.
The self-adjoint strip operator $B_x$ has kernel
$T_xA_i\cap T_xL_C=0$. Its spectral gap and Fourier transform
in $s$ give an inverse, also for small end weights. Constants
are thus regular of index zero. Sard--Smale and a countable
exhaustion now give regular upper paths for every required
auxiliary pair. Intersect the residual sets with those for the
ordinary upper strips, having fixed regular ordinary lower data first.
The output operator for $(A_0,A_1)$ is their direct sum.
Different triangle ends permit independent almost-complex paths;
these paths extend to split structures on the triangle.

For the triangle the same constant-upper argument would force
a constant full map and hence a triple intersection, which is
empty. Domain-dependent upper variations in a free interior
disc therefore kill the upper adjoint; no global injectivity
condition is required for these domain-dependent variations.
Equation~\eqref{eq:upper-trace-injective} and boundary unique
continuation kill the lower adjoint. This proves triangle
regularity while retaining zero Hamiltonian, the fixed strip
ends and the prescribed collar of the collapse interval. For
$k=0$ choose the collars along which the configuration is doubled at the outset; the
variations in the free interior discs remain available.

For a parameter family multiply these variations by parameter bumps, zero near the
prescribed regular endpoint collars. The parametrized operator includes the data parameter
once; if several ends share that parameter, their independently supported variations enter
the one parametrized operator. The Sard--Smale theorem, applied relative to the endpoint
collars, makes the parametrized moduli spaces regular, which is all that is used.

The area-index relation above gives a fixed-index energy bound.
Every positive sphere or disk bubble removes at least four
from the core index. A regular nonconstant autonomous strip
has its nonzero translation vector, so its full index is at
least one. A constant strip has index zero. Thus in full strip
indices one and two no bubble can occur, and the index-two
translation quotient has exactly two index-one strips at its
boundary. The same index count applies to triangles of dimensions zero and one and to
the corresponding regular parametrized families, with the parameter counted once. When
several strips share a parameter, the quotient by translation is taken after forming the
parametrized operator. Ordinary
nondegenerate strip gluing gives the boundary collars. Its
converse follows from zero-loss compactness and exponential
decay; a sequence outside the finitely many local gluing
neighbourhoods would converge strongly into one of them.

Finally the energy and index of an auxiliary strip from the
maximal generator are
\[
 E(u)=8\pi^2\kappa_{\rm aux}
          +2\bigl(\mathcal A_i(q_i^*)-\mathcal A_i(x)\bigr),
 \qquad
 \ind(u)=m(q_i^*)-m(x)+8\kappa_{\rm aux}.
\]
The small action oscillation is less than $2\pi^2$ and
$2\kappa_{\rm aux}$ is integral, so energy excludes negative
classes. The Morse difference is a nonnegative multiple of
three, and $8\kappa_{\rm aux}$ a nonnegative multiple of four.
Their sum is not one. No differential trajectory leaves $q_i^*$.
Changing $J$ has altered neither these ordered actions nor their
relative grading. This proves the cycle assertion directly on
the product-target problem.
\end{proof}

\begin{definition}\label{def:cup}
The \emph{compression cup} of a $k$-fold one-handle, at height Lagrangians, is
\[
U_C=\mathcal B\,\mathcal K\bigl(\mu_{\mathbb X}(q_0^*,q_1^*)\bigr)\colon
\HF(\widetilde P_0,\widetilde P_1)\to\HF(R_0,R_1),
\]
where $\mu_{\mathbb X}\colon\CF_{\mathbb X}(A_0,L_C)\otimes\CF_{\mathbb X}(L_C,A_1)\to\CF_{\mathbb X}(A_0,A_1)$ is the
triangle map in $\mathbb X$, $\mathcal K$ is the folding isomorphism and $\mathcal B$
turns a cycle $\sum x_i^\vee\otimes y_i$ into the map $v\mapsto\sum x_i^\vee(v)
y_i$. At other admissible presentations it is conjugated by continuation maps. The
\emph{compression cap} $V_C\colon\HF(R_0,R_1)\to\HF(\widetilde P_0,\widetilde
P_1)$ of a three-handle is the transpose of the cup of the reflected
one-handle, $V_C=(U_{\bar C})^\vee$ under the pairings of
Section~\ref{subsec:reflection}.
\end{definition}

\begin{figure}[ht]
\centering
\begin{tikzpicture}[scale=1.0,font=\small]
\coordinate (a) at (-2.4,0); \coordinate (b) at (2.4,0); \coordinate (c) at (0,3.2);
\fill[gray!12] (a) -- (b) -- (c) -- cycle;
\draw[thick] (a) -- (b) node[midway,below] {$L_C$};
\draw[thick] (a) -- (c) node[midway,left] {$A_0$};
\draw[thick] (b) -- (c) node[midway,right] {$A_1$};
\fill (a) circle (1.8pt); \fill (b) circle (1.8pt); \fill (c) circle (1.8pt);
\node[below left] at (a) {$q_0^*\in A_0\cap L_C$};
\node[below right] at (b) {$q_1^*\in L_C\cap A_1$};
\node[above] at (c) {output in $A_0\cap A_1$};
\node at (0,1.1) {$\Msf_h^-\times\Msf_{h+k}$};
\end{tikzpicture}
\caption{The folded triangle defining the compression cup
(Definition~\ref{def:cup}). The inputs are the top generators $q_0^*$ and $q_1^*$ of the
auxiliary pairs, and the output lies in $\CF_{\mathbb X}(A_0,A_1)$, which the folding isomorphism
$\mathcal K$ identifies with $\CF(\widetilde P_0,\widetilde P_1)^\vee\otimes\CF(R_0,R_1)$, that is,
with maps from $\CF(\widetilde P_0,\widetilde P_1)$ to $\CF(R_0,R_1)$. Unfolded, the triangle
is a quilted surface whose seam carries the correspondence $L_C$
(Remark~\ref{rem:quilted-cup}).}
\label{fig:folded-triangle}
\end{figure}
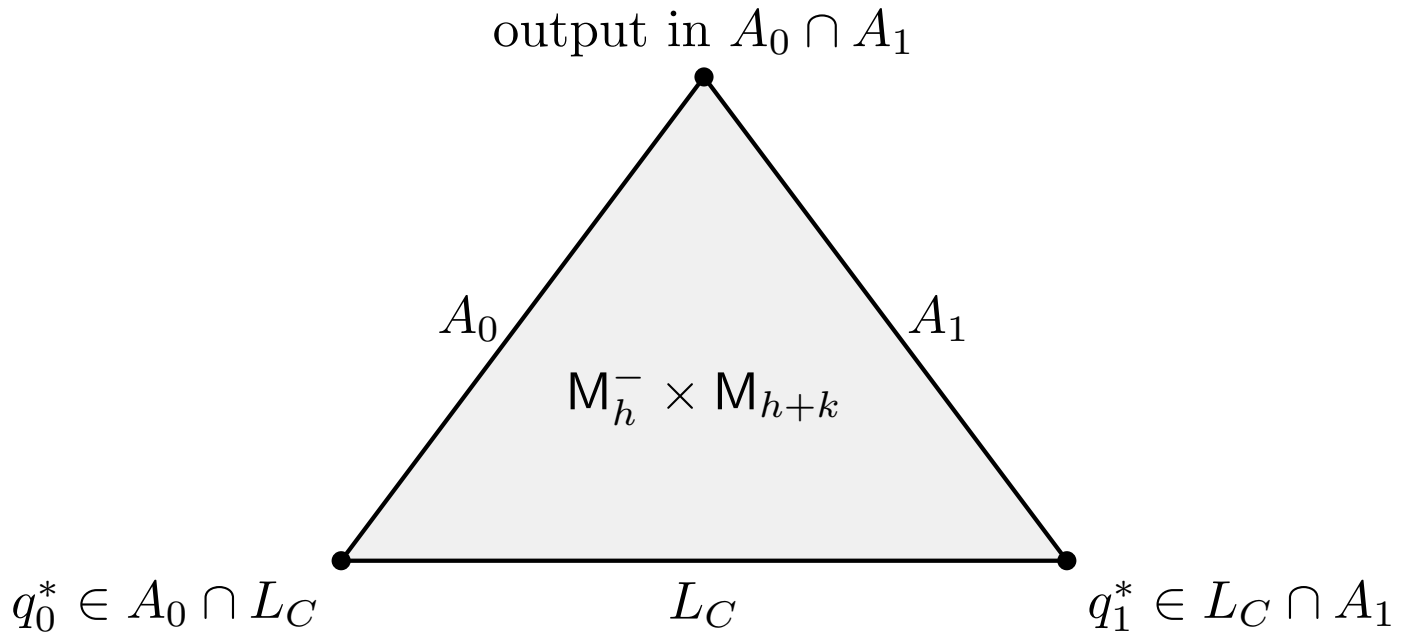

The cup is a count of holomorphic triangles, parallel to the two-handle map
$\mu_S(\,\cdot\otimes\theta_S)$; its definition involves no gauge theory and no
geometric composition.

\begin{lemma}\label{lem:diagonal-cup}
For $k=0$, $L_C=\Delta\subset\Msf_h^-\times\Msf_h$ and $U_C$ equals the continuation
map $\mathrm{cont}_S(\widetilde P,R)$ on homology.
\end{lemma}

\begin{proof}
Let $S$ be the triangle with boundary arcs labelled $A_0$, $\Delta$, $A_1$ and
$D$ the disk obtained by doubling $S$ along its $\Delta$-arc. For data
$(J_D,K_D)$ on $D$ with $K_D$ vanishing on the doubling arc, put
$J_{\mathbb X}(z)=(-J_D(\bar z))\times J_D(z)$ on $S$. The map $v\mapsto(v\circ\bar{\ },v)|_S$ is a
bijection from solutions on $D$ to solutions on $S$ with boundary on
$(A_0,\Delta,A_1)$: a solution on $S$ glues to a Lipschitz map on $D$ that solves
the equation weakly across the arc, hence is smooth; linearized operators
correspond. The two ends of $S$ adjacent to the $\Delta$-arc double to strip
ends of $D$ labelled $(R_0,\widetilde P_0)$ and
$(\widetilde P_1,R_1)$, respectively, and the output end doubles to one
incoming and one outgoing end, which is the folding $\mathcal K$ followed by
$\mathcal B$. So $U_C(x)$ is the count of disks with boundary arcs labelled
$R_1,\widetilde P_1,\widetilde P_0,R_0$ with $x$ at the incoming end and $q_0^*$,
$q_1^*$ at the two auxiliary ends. We construct the unit cap for
each of these ordered pairs. In a Weinstein neighbourhood of
$\widetilde P_i$, write $R_i=\operatorname{graph}(dg_i)$.
For the first pair, fibre translation by $-dg_0$ takes
$(\operatorname{graph}(dg_0),0)$ to
$(0,\operatorname{graph}(-dg_0))$. Realize this translation by a
Hamiltonian isotopy with cutoff equal to one near the whole graph
path. The second pair is already $(0,\operatorname{graph}(dg_1))$.
Thus the ordered primitives are
$f_0=-g_0-\min(-g_0)$ and $f_1=g_1-\min g_1$.
Put $M_i=\max f_i$. For the small primary data in use, $f_i$
has its unique maximum at $q_i^*$. Extend $f_i$ constantly in the
cotangent fibres near the zero section and cut it off farther away
by a function in $[0,1]$. The resulting $K_i$ satisfies
$0\le K_i\le M_i$ and $K_i(q_i^*)=M_i$. Choose the cutoff
equal to one on the whole graph path. With
$\iota_{X_{K_i}}\omega=-dK_i$, the inverse time-one flow
carries the zero section to $\operatorname{graph}(df_i)$.

Use a disk with one outgoing strip end, with
boundary on that zero section and Hamiltonian one-form
$\beta(s)\chi(\theta)(K_i-M_i)\,d\theta$.
Here $\beta$ increases from zero to one, while $\chi\ge0$ is
supported away from the strip sides and has integral one.
On the end, conjugating by the inverse Hamiltonian flow gives
the ordinary strip for $(0,\operatorname{graph}(df_i))$ and fixes
its critical chords. Transport the first cap back by the actual
translation, including its almost-complex and Hamiltonian data.
Its ordered labels are then $(R_0,\widetilde P_0)$; those of the
second cap are $(\widetilde P_1,R_1)$. A cap with output $x$
has index $m(q_i^*)-m(x)+8\kappa_{\mathrm{aux}}$: by additivity of the index, gluing the
constant cap at $q_i^*$ to a strip from $q_i^*$ to $x$ gives the grading
difference of Section~\ref{subsec:continuations}, and the constant cap has
index zero: by spectral flow of $\beta\,\mathrm{Hess}(f_i)$ along the end, the index of the
constant cap at a critical point $c$ is the number of positive eigenvalues of
$\mathrm{Hess}(f_i)(c)$, which is zero exactly at the maximum $q_i^*$. The full energy identity for this cap is
\[
E_{\mathrm{cap}}=8\pi^2\kappa_{\mathrm{aux}}
 -2(K_i(x)-M_i)
 +2\int\beta'\chi(K_i(u)-M_i)\,ds\,d\theta.
\]
The integral is nonpositive and the endpoint term is at most
$2M_i$. Thus $2M_i<4\pi^2$ excludes negative relative classes.
The index formula then forces an index-zero cap to have
$\kappa_{\mathrm{aux}}=0$ and output $q_i^*$. Its energy is
zero, so it is the constant cap. At this solution the linear
energy identity, with the nonpositive Hessian of $K_i-M_i$,
forces a decaying kernel field to be horizontal with zero
outgoing value, hence zero. Index zero gives regularity.
Gluing
the two unit caps into the auxiliary ends \cite[Theorem~3.13]{WW-relative}
produces a strip with moving boundary from
$(\widetilde P_0,\widetilde P_1)$ to $(R_0,R_1)$, whose map is the continuation
\cite[Theorem~4.1, Example~2.6]{WW-relative}.
\end{proof}

\begin{remark}\label{rem:quilted-cup}
By the unit-insertion formula of Wehrheim and Woodward
\cite[Corollary~5.5]{WW-relative}, $U_C$ agrees with the morphism map of the
quilted correspondence $L_C$ \cite{WW-quilt} followed by the geometric
composition isomorphisms $\widetilde P_i\circ L_C$ and continuation to $R_i$.
That formula rests on the strip-shrinking theorem
\cite[Theorem~5.1]{WW-relative}; for the analysis of strip shrinking see
\cite{WW-composition, BW-squiggly}. The
composition isomorphism itself is also established by Lekili and Lipyanskiy
without strip shrinking \cite[Theorem~2 and Corrigendum,
Remark~1(1)]{LL-composition}, whose hypotheses hold here: the targets are
monotone with the same constant, the Lagrangians $\widetilde P_i$ and $L_C$ are
compact, connected, orientable and simply connected with minimal Maslov
number $4$, and the compositions are embedded. Neither statement is used in
this paper.
\end{remark}

\begin{lemma}\label{lem:aux-index-folded}
The statements of Lemma~\ref{lem:aux-index} hold at the height data for the fillings
$F_i=V_i\times[0,1]$ of $Z_i=V_i\cup(-V_i)=(\#^{g-1}S^1\times S^2)\#T^3$ and the
pairs $(A_0,L_C)$ and $(L_C,A_1)$, with $r=g-1$, the height functions defining $R_i$, and the
top generators $q_i^*$. Fix the height $t$ and a finite family $F\supseteq F_*$ as in
Definition~\ref{def:presentation}. The \emph{end data} of an auxiliary pair are its lower
primary datum, in $\operatorname{span}F$, and the primary datum of its height Lagrangian, in
$\operatorname{span}F_*$; the \emph{height data} are the end data $(0,th_i)$. For $k\ge1$ the
height datum, the family $F_*$ and the function $\delta$ of Section~\ref{subsec:one-three-squares}
are those of the splitting $(V_0,V_1)$ of $M^+$, whose height functions are those of
Section~\ref{subsec:folded} on $V_i$, so that $R_i$ is a height Lagrangian in the sense of that
section with $k=0$.
\begin{enumerate}[label=(\roman*)]
\item Let $R$ be a height Lagrangian at height $t$ with primary datum $th_i+g$,
  $g\in\operatorname{span}F_*$, $|g|\le\delta(t)$ (Section~\ref{subsec:one-three-squares}).
  There is a ball $\mathcal N_{C,R}$ about zero in $\operatorname{span}F$, depending
  on $C$, on $R$ and on $F$, such that, for end and secondary perturbations below the threshold
  of (ii) with $c=2\pi^2$, under changes of the lower primary data within
  $\mathcal N_{C,R}$ the relative classes and integer gradings are transported from the height
  data, the endpoint height differences in the energy identities of the gauge trajectories, mixed
  solutions and strips are replaced by the differences
  of the actual critical actions $\mathcal A_i$, and the coefficient of $q_i^*$ in $N_i\theta_i$,
  and in the filling invariant, remains $1$; the representing solutions need not remain flat.
\item Let $(\mathfrak h_s,\mathfrak r_s)$, $0\le s\le1$, be a piecewise smooth path of end data
  starting at the height data, such that at every $s$ the generators of the auxiliary pairs are
  nondegenerate and continue those at the height data, with their gradings $m$ and relative
  classes, and $2\max_{a,x}|\mathcal A_i(a)-\mathcal A_i(x)|<c$ for a constant $c<4\pi^2$
  independent of $s$, where $\mathcal A_i$ are the critical actions of the auxiliary pair at
  $(\mathfrak h_s,\mathfrak r_s)$, lifted continuously along the branches of generators. Then at
  $(\mathfrak h_1,\mathfrak r_1)$, for end perturbations $f$ along the path and secondary
  perturbations $\eta$ at the endpoint with
  $4\sup\|f\|_\infty+2\sup_A\|\eta(A)\|^2_{L^2}<4\pi^2-c$, below the thresholds of
  Theorem~\ref{thm:mixed-iso}: every counted gauge trajectory, mixed solution and strip has
  $\kappa_{\mathrm{aux}}\ge0$; the differential of the auxiliary gauge complex vanishes; the
  filling invariant is the generator $\theta_i$ of maximal $m$; and $N_i\theta_i=q_i^*$, both on
  chains.
\end{enumerate}
\end{lemma}

\begin{proof}
The proof of Lemma~\ref{lem:aux-index} uses only that the filling retracts
onto a wedge of $r$ circles and a torus, so that $H^1_c=H^2=0$ with twisted
coefficients, and that the critical manifold is $\SU(2)^r$ with the height
function; both hold for $F_i$. The mixed maps of the pairs $(A_0,L_C)$ and $(L_C,A_1)$ are those
of Remark~\ref{rem:clean-pairs}. We prove (ii) first.

\emph{Energy.} The gradings and relative classes are transported along the path by hypothesis,
and the critical actions are the perturbed Chern--Simons values at the generators, lifted
continuously along the branches from the height data; their differences at one parameter are
the differences of the graph primitives (Section~\ref{subsec:continuations}). On the gauge
pieces of the domains of the gauge trajectories, mixed solutions and strips the end data are
autonomous, so they enter the analytic energy only through these critical values, as in the
proof of Theorem~\ref{thm:mixed-iso}: the energy identity of Lemma~\ref{lem:aux-index}(ii)
holds at every $s$ with the endpoint height differences replaced by
$2(\mathcal A_i(a)-\mathcal A_i(x))$, and the bounds on $c$, $f$ and $\eta$, with nonnegative
analytic energy and $\kappa_{\mathrm{aux}}\in\frac12\ZZ$, give $\kappa_{\mathrm{aux}}\ge0$ as in
Lemma~\ref{lem:aux-index}(iii). The fillings are not treated in this way: their energy identity
contains the profile term of the lower datum, bounded only by its $C^0$ norm, which the
hypotheses do not control.

\emph{Vanishing differentials.} At every parameter at which the data are regular the auxiliary
gauge complex has $2^r$ generators. It is the framed complex of $\#^r(S^1\times S^2)$, with the
lifting-group quotient of Section~\ref{subsec:complexes}, so its homology is
$\Ish(\#^r(S^1\times S^2);\FF)$, of dimension $2^r$: the integral group is free of rank $2^r$
\cite[Sections~7.7--7.8]{Scaduto}, as recalled in the proof of
Proposition~\ref{prop:filling-classes}, and universal coefficients give the dimension over
$\FF$. Hence its differential vanishes, and a chain map between two such complexes that is a
homotopy equivalence is an isomorphism.

\emph{The mixed coefficient.} At the endpoint $N_i$ is a chain isomorphism by
Theorem~\ref{thm:mixed-iso}(i), and by the index formula of Lemma~\ref{lem:aux-index}(ii) its
matrix entry from $a$ to $x$ vanishes unless $m(x)=m(a)+8\kappa_{\mathrm{aux}}$ with
$\kappa_{\mathrm{aux}}\ge0$. Lemma~\ref{lem:calibration} gives $N_i\theta_i=q_i^*$ on chains
for the generator $\theta_i$ of maximal $m$. This step uses only the data at the endpoint.

\emph{The filling invariant.} At the height data, for end and secondary regularizations below
the threshold of Proposition~\ref{prop:filling-classes}, the filling invariant is $\theta_i$.
The filling with the end data
$(\mathfrak h_1,\mathfrak r_1)$ is joined, through a path of interior data with fixed ends, to
the filling with the height data glued to the continuation cylinder $\RR\times Z_i$ carrying the
path $(\mathfrak h_s,\mathfrak r_s)$ and an interpolation of the end perturbations; by
Theorems~\ref{thm:family} and~\ref{thm:neck}, the filling invariant at the endpoint is, on
homology and hence on chains, $\mathrm{cont}_G$ of the filling invariant at the height data,
where $\mathrm{cont}_G$ is the gauge continuation map along the path, the map of the product
cobordism with varying end data. It remains to show $\mathrm{cont}_G\theta_i=\theta_i$.
Choose a finite partition $s_j$ of $[0,1]$, including the corners of the path,
fine enough for the work and action bounds below. At its finitely many vertices
choose regular end perturbations $f_j$ within the stated bound by
Proposition~\ref{prop:input-functions}, retaining the prescribed regular data
at the two endpoints. Then regularize the finitely many continuation cylinders.
For a continuation trajectory from $a$ at $s_j$ to $x$ at $s_{j+1}$ the relative class
$\kappa_{\mathrm{aux}}$ is taken with respect to the reference at $s_j$ concatenated with the
branch of $x$ from $s_j$ to $s_{j+1}$; it is quantized, and the index formula holds for it, the
spectral flow along a nondegenerate branch being zero. The continuation map along the
path is the composite of the maps along the intervals $[s_j,s_{j+1}]$, by functoriality of the
cobordism maps for the composite of the product cobordisms, up to chain homotopy, hence
exactly, the differentials vanishing. A continuation trajectory of index zero from $\theta_i$ at
$s_j$ to $x$ at $s_{j+1}$ has $m(x)=m(\theta_i)+8\kappa_{\mathrm{aux}}$, the gradings being
continued, and its analytic energy is $8\pi^2\kappa_{\mathrm{aux}}$ plus
$2(\mathcal A_i(\theta_i)-\mathcal A_i(x))$, the actions taken at $s_j$ and at $s_{j+1}$, plus
the end terms, bounded by $2\|f_j\|_\infty+2\|f_{j+1}\|_\infty$, plus the work term
$2\int\partial_\sigma(\mathfrak h,\mathfrak r)_{s(\sigma)}(A_\sigma)\,d\sigma$ of the
interpolation $\sigma\mapsto s(\sigma)$, plus the term of the small compactly supported
holonomy perturbation that makes the continuation cylinder regular
(Theorem~\ref{thm:regularity}(iii)), bounded by twice its squared norm, which is as small as we
please. The work term is bounded by
$2\int_{s_j}^{s_{j+1}}\sup|\partial_s(\mathfrak h_s,\mathfrak r_s)|\,ds$, and the action term
differs from the difference at the single parameter $s_j$ by a quantity tending to zero with the
length of the interval, the actions being continuous along the branches.
The partition depends on this fixed path in the finite-dimensional span:
the sup norm of its derivative is integrable, and the finitely many action
branches are uniformly continuous. Choose the partition first and the
cylinder regularizers afterward, so that these corrections are less than
the margin $4\pi^2-c-4\sup\|f\|_\infty$, which is
positive by hypothesis, so $\kappa_{\mathrm{aux}}\ge0$; hence the continuation map along
$[s_j,s_{j+1}]$ has matrix entries vanishing unless $m(x)=m(a)+8\kappa_{\mathrm{aux}}$ with
$\kappa_{\mathrm{aux}}\ge0$. It is an isomorphism, being a homotopy equivalence between
complexes with vanishing differentials, so Lemma~\ref{lem:calibration} gives
$\theta_i\mapsto\theta_i$. Composing, $\mathrm{cont}_G\theta_i=\theta_i$, and the filling
invariant at the endpoint is $\theta_i$.

(i) follows from (ii) along the path that first moves the height datum from $th_i$ to $th_i+g$
along $\sigma\mapsto th_i+\sigma g$, at lower datum zero, and then the lower datum from zero to
$\mathfrak h$ along $\sigma\mapsto\sigma\mathfrak h$. On the first segment the hypotheses of
(ii) hold with $c=2\pi^2$ by the choice of $\delta(t)$ in
Section~\ref{subsec:one-three-squares}, which is made so that at lower datum zero and for
$|g|\le\delta(t)$ the auxiliary pairs have $2^r$ transverse generators continuing those at $g=0$
with their gradings and relative classes, and critical actions of oscillation less than
$\pi^2$. On the second segment they
hold for $\mathfrak h$ in a sufficiently small ball of $\operatorname{span}F$ about zero, which
is convex: nondegeneracy of the finitely many generators and the bound on the critical actions
are open conditions in the lower datum, satisfied at $\mathfrak h=0$. Under these changes the
relative classes and gradings are transported by construction, and the representing solutions
need not remain flat.
\end{proof}

\subsection{The folded triangle}

\begin{theorem}\label{thm:folded-triangle}
Let $N_0$, $N_1$ be the mixed maps of the product-target pairs
$(A_0,L_C)$ and $(L_C,A_1)$,
with gauge inputs $Z_0=D_0\cup(-C)$ and $Z_1=C\cup(-D_1)$, both diffeomorphic to
$V_i\cup(-V_i)=(\#^{g-1}S^1\times S^2)\#T^3$, and $N_Z$ the mixed map of
$(A_0,A_1)$ with the data at width $\lambda_0$ of Theorem~\ref{thm:excision}. Then
\[
N_Z\circ\mu_{\mathscr W}\simeq\mu_{\mathbb X}\circ(N_0\otimes N_1)\quad\text{over }\FF,
\]
where $\mu_{\mathscr W}$ is the gauge triangle of the pieces $D_0$, $C$, $D_1$.
Moreover $N_i\theta_i=q_i^*$.
\end{theorem}

\begin{proof}
The proof is that of Theorem~\ref{thm:triangle} with the following inputs.
\emph{Collapse.} The middle compression body $C=C_{h,g}$ is an odd
compression body (Section~\ref{subsec:clean-data}), with $L(C)=L_C\subset
\Msf_h^-\times\Msf_g$ and both boundary surfaces matched. Hypotheses
\ref{hyp:X1}--\ref{hyp:X3} hold as in the proof of Theorem~\ref{thm:two-handle}: the middle
primary perturbation vanishes, the end perturbations are the fixed autonomous perturbations of
the gauge input ends $Z_0$, $Z_1$, with outputs and loops in those ends and cutoff bands placed
beyond the compact set $N_I$, and all secondary perturbations have outputs and loop regions in a
compact subset of the interiors of $D_0$ and $D_1$, disjoint from $N_I$. The triangle piece in $\mathbb X$ is regular for a generic domain-dependent almost
complex structure supported away from $N_I$, as in the proof of
Theorem~\ref{thm:triangle}. For \ref{hyp:X4}, stretch the two ordinary necks cutting off the
DFL disks of the pairs $(A_0,L_C)$ and $(L_C,A_1)$ from the triangle domain of $(A_0,L_C,A_1)$;
Theorem~\ref{thm:neck} gives regularity in index at most zero with count
$\mu_{\mathbb X}(N_0\otimes N_1)$. The linear inputs of Section~\ref{sec:collapse-linear}
hold in this instance with the two-trace slice space
$H^1(\Sigma_h;\ad)\oplus H^1(\Sigma_g;\ad)$: the gap is uniform in the scale
because no harmonic field on $C$ has vanishing harmonic trace
($H^1(C,\partial C;\ad)=0$), and the reversed lower sheet enters the Green identity
with the sign of Convention~\ref{conv:orientation}. The period lattice of capped
configurations is the image of $\ZZ^2/\ZZ(1,1)$, of index step four, as in
Section~\ref{subsec:energy-convention}.
\emph{Outgoing cut.} Stretch along the closed cut $Z\cong(-M)\#M^+\#T^3$, which is
connected. Its critical points are nondegenerate and irreducible, the connecting
$H^0$ term vanishes, and the framed groups glue by the fibre product;
Theorem~\ref{thm:neck} gives $N_Z\circ\mu_{\mathscr W}$.
\emph{Family.} The path of data has the five stages of
Lemma~\ref{lem:family-path}: (A) raise the scale of $C=C_{h,g}$ from $\eps$ to
one; (B) interpolate to a metric that is a product $a^2du^2+g_Z$ on a collar of
the cut $Z$, where the cut crosses the edge pieces $D_0$ and $D_1$, on which the
collar coordinate is the edge coordinate, and the two scaled sheets over $B$
(Lemma~\ref{lem:cut-crosses}), and $g_Z$ restricts to the metrics of $D_0$ and $D_1$; (C) lengthen the
collar; (D) switch on the end perturbation of $Z$ on the collar; (E)
interpolate the remaining data to the data with the collar of the cut inserted. Each stage lies in a convex set of data satisfying the hypotheses of Theorem~\ref{thm:family}, the parameter derivative is compactly
supported away from the ends, and invariance holds because the classes of the
lifting group vanish on $\Sigma_h\sqcup\Sigma_g$. Theorem~\ref{thm:family}
applies.
\emph{Classes.} The pair $(R_i,\widetilde P_i\circ L_C)$ is the auxiliary pair of a height Lagrangian with unique top generator $q_i^*$, and
Lemma~\ref{lem:aux-index-folded} gives $N_i\theta_i=q_i^*$: part (i) for lower data in
$\mathcal N_{C,R}$, and part (ii), for $k=0$, for lower data in the set $\mathcal U_F(t)$ of
Definition~\ref{def:admissible-nbhd}, along the paths provided there.
\end{proof}

\begin{proposition}[cup formula]\label{prop:cup-formula}
Put $N^\dagger_P=N_-^\vee\colon\CF(\widetilde P_0,\widetilde P_1)\to C_G(M^\#)$, the transpose of the mixed
map $N_-=N_{\bar P}$ of the reflected presentation, so that $N^\dagger_P=N_{\bar P}^\vee$, defined by
$\langle a,N^\dagger_Px\rangle_G=(N_-a)(x)$. Then on homology
\[
U_C=N_+\circ\Ish(W_k)\circ N^\dagger_P .
\]
For $k=0$ this reads $\mathrm{cont}_S(\widetilde P,R)=N_R\circ \mathrm{cont}_G(\widetilde P,R)\circ N^\dagger_P$,
with $\mathrm{cont}_G$ the gauge continuation.
\end{proposition}

Here and below the mixed map and its transpose are indexed by the primary data when these
vary: for a presentation $p$ with primary data $P$ we write $N_P=N_p$ and $N^\dagger_P=N^\dagger_p$.

\begin{proof}
Insert $\theta_0\otimes\theta_1$ into Theorem~\ref{thm:folded-triangle}. Here
$N_i\theta_i=q_i^*$ at the data of Definition~\ref{def:cup}, and also at the generic data
of Theorem~\ref{thm:folded-triangle}: a continuation strip from $q_i^*$ of index zero has
$\kappa_{\mathrm{aux}}=0$ and ends at $q_i^*$, and a strip of zero energy is constant, so
the continuation maps $q_i^*$ to $q_i^*$ on chains; and $\mu_{\mathbb X}(q_0^*,q_1^*)$ is natural
under continuation with fixed ends. By Definition~\ref{def:cup} and
Theorem~\ref{thm:excision},
\[
U_C=\mathcal B\mathcal K\mu_{\mathbb X}(q_0^*,q_1^*)=\mathcal B\mathcal K N_Z
\mu_{\mathscr W}(\theta_0,\theta_1)=\mathcal B(N_-\otimes N_+)\mathsf E(\Theta_{B_k}),
\]
and $\mathsf E(\Theta_{B_k})=\mathrm{ten}(\Ish(W_k))=\sum a_i\otimes b_i$ by
Lemma~\ref{lem:handle-state}. Hence
\[
U_C(x)=\sum_i(N_-a_i)(x)\,N_+b_i=\sum_i\langle a_i,N^\dagger_Px\rangle
N_+b_i=N_+\Ish(W_k)N^\dagger_Px. \qedhere
\]
\end{proof}

The factor $N^\dagger_P$ is not yet known to be $N_P^{-1}$; that is proved below, after
naturality, in this order.

\subsection{Order of choices}\label{subsec:order}

In each application of
Proposition~\ref{prop:cup-formula}: first the upper height $t$ and the upper Lagrangians $R_i$; then the lower Lagrangians $\widetilde P_i$ in a neighbourhood small
relative to $t$, subject to the generic incidence conditions of
Section~\ref{subsec:folded} (expected dimension $-(3h+3g-6)<0$); then the data at width $\lambda_0$ of Theorem~\ref{thm:excision}. A finite change of torus
collar length returns the torus factors to the fixed data by
Proposition~\ref{prop:collar-length}, and the secondary data are then replaced
by the prescribed ones by Theorem~\ref{thm:mixed-iso}(ii). Hence the maps $N$
appearing in one application are the mixed maps on homology. Squares proved with
the data of one application are moved to other data at the vertices by
Proposition~\ref{prop:small-primary}, whose own quantifiers are ordered in the
same way: height first, then the neighbourhood $\mathcal N(t)$ of
Definition~\ref{def:admissible-nbhd}, which determines $\mathcal U_F(t)$, then the lower data,
and last a generic height Lagrangian.

\subsection{The one- and three-handle squares}\label{subsec:one-three-squares}

A \emph{height Lagrangian at height $t$} for $k=0$ is an upper Lagrangian $R$ of
Section~\ref{subsec:folded} for the product compression body: the primary datum
of $U$ perturbed by $t$ times the height function $h_i$ of Section~\ref{subsec:folded}, with
its unique maximal generator $q^*$. A further perturbation of its primary datum by
$g\in\operatorname{span}F_*$ with $|g|\le\delta(t)$ keeps it a height Lagrangian at height $t$,
for a function $\delta>0$ depending only on the splitting, the metric data and the height
function, with $\delta(t)\to0$ as $t\to0$: $\delta(t)$ is so small that for $|g|\le\delta(t)$
and at unperturbed lower data the generators of the auxiliary pairs, the intersection points of
$R$ with the unperturbed Lagrangian, are $2^r$ transverse points continuing those at $g=0$ with
their gradings $m$ and relative classes, and that their critical actions have oscillation less
than $\pi^2$. Such $\delta(t)$ exists, for $t$ below a threshold depending only on the splitting,
the metric data and the height function: at $g=0$ the generators are the nondegenerate critical
points of $th_i$ on the unperturbed Lagrangian, with critical actions of oscillation
$t\operatorname{osc}h_i$, and both conditions are open in $g$, by the implicit function theorem
and compactness. Then $q^*$ continues to the unique generator of maximal $m$; the generators
need not remain at the critical points of the height function. The generic perturbations used below are taken from
this ball. They are not required to vanish near the
generators: a perturbation vanishing near a generator cannot move it, and the incidence
conditions may require that it move, for instance at a common flat point of the two
compression Lagrangians that is critical for both height functions.

\begin{proposition}[small-primary naturality]\label{prop:small-primary}
Fix a splitting, metric data and a finite family $F\supseteq F_*$ of primary perturbations.
For every height $t>0$ below a threshold depending only on the splitting, the metric data and
the height functions, and any two regular data $P,Q$ in the neighbourhood
$\mathcal U_F(t)$ of zero in $\operatorname{span}F$ of Definition~\ref{def:admissible-nbhd},
$N_Q\circ \mathrm{cont}_G(P,Q)=\mathrm{cont}_S(P,Q)\circ N_P$ on homology.
\end{proposition}

\begin{proof}
The quantifiers are chosen in this order: first the height $t$, then the neighbourhood
$\mathcal N(t)$ of Definition~\ref{def:admissible-nbhd}, which determines $\mathcal U_F(t)$,
then the lower data, and last a generic height Lagrangian. The conditions defining
$\mathcal U_F(t)$ are symmetric under reflection, which exchanges the two Lagrangians of the
pair, and the choices for the reflected splitting are the reflections of those for the
splitting (Definition~\ref{def:admissible-nbhd}), so $\bar P$ and $\bar Q$ lie in the
corresponding set for the reflected family. Given
regular $P,Q\in\mathcal U_F(t)$, choose one height Lagrangian $A$ at height $t$ for the
reflected presentation which, by a small generic perturbation $g\in\operatorname{span}F_*$,
$|g|\le\delta(t)$, satisfies the incidence conditions of
Section~\ref{subsec:folded} with respect to both $\bar P$ and $\bar Q$; this is
possible because the universal restriction map of $\operatorname{span}F_*$ is a
submersion at the perturbed flat connections of $th_i+g$, for $t$ and $\delta(t)$ below a
threshold fixed with the splitting, by openness as in the proof of
Lemma~\ref{lem:app-separated}(ii), and the conditions are generic for each of the two lower
pairs. By
Lemma~\ref{lem:aux-index-folded}(ii), applied along the path that first moves the height datum
to that of $A$ and then the lower datum from zero to $\bar P$, respectively to $\bar Q$, along
the paths provided by Definition~\ref{def:admissible-nbhd}, its conclusions hold for the
auxiliary pairs of $(\bar P,A)$ and $(\bar Q,A)$, so the product case of
Proposition~\ref{prop:cup-formula} applies to $(\bar P,A)$ and $(\bar Q,A)$ at once and gives
$\mathrm{cont}_S(\bar P,A)=N_A\mathrm{cont}_G(\bar P,A)N_{\bar{\bar P}}^\vee$ and similarly for $\bar Q$;
with $[N_{\bar{\bar P}}]=[N_P]$ (Section~\ref{subsec:reflection}) this is
$B_{\bar P}=N_AG_{\bar P}N_P^\vee$, $B_{\bar Q}=N_AG_{\bar Q}N_Q^\vee$, where $G$ and $B$
are the gauge and symplectic continuations to $A$. Eliminating the invertible
$N_A$ gives $G_{\bar Q}^{-1}G_{\bar P}N_P^\vee=N_Q^\vee B_{\bar Q}^{-1}B_{\bar P}$. The
composites are the transposes of $\mathrm{cont}_G(Q,P)$ and $\mathrm{cont}_S(Q,P)$ under reflection;
transposing gives the claim.
\end{proof}

\begin{proposition}[mixed duality]\label{prop:duality}
Let $t$ and $F$ be as in Proposition~\ref{prop:small-primary}. There is a height $t'\le t$,
depending only on $t$, the splitting, the metric data and the height functions, such that
$N_{\bar P}=(N_P^{-1})^\vee$ on homology for every regular $P$ in the path component
$\mathcal U'_F(t)$ of zero in $\mathcal U_F(t)\cap\mathcal U_F(t')$
(Definition~\ref{def:admissible-nbhd}).
\end{proposition}

\begin{proof}
Choose $t'\le t$ so small that every height Lagrangian at height at most $t'$, with any
$g\in\operatorname{span}F_*$, $|g|\le\delta(t')$, and its reflection lie in $\mathcal N(t)$;
this is possible because these Lagrangians form a compact family converging to the unperturbed
pair as the height tends to zero. For a height Lagrangian $R$ at height $t'$ with primary datum
$t'h_i+g$, the segment from zero to $t'h_i+g$ in $\operatorname{span}F_*\subseteq\operatorname{span}F$
then lies in the preimage of $\mathcal N(t)$, so $R\in\mathcal U_F(t)$.
Given regular $P\in\mathcal U'_F(t)$, let $R$ be a height Lagrangian at height $t'$ for
the presentation of $P$, generic with respect to $P$ as in the proof of
Proposition~\ref{prop:small-primary}. The product case of
Proposition~\ref{prop:cup-formula} at height $t'$, available at $P\in\mathcal U_F(t')$ by
Lemma~\ref{lem:aux-index-folded}(ii), gives
$\mathrm{cont}_S(P,R)=N_R\mathrm{cont}_G(P,R)N_{\bar P}^\vee$, and Proposition~\ref{prop:small-primary}
at height $t$, applied to the primary data $P,R\in\mathcal U_F(t)$, gives
$N_R\mathrm{cont}_G(P,R)=\mathrm{cont}_S(P,R)N_P$. Hence $\mathrm{cont}_S(P,R)=\mathrm{cont}_S(P,R)N_PN_{\bar P}^\vee$;
cancelling the isomorphism $\mathrm{cont}_S(P,R)$ gives $N_PN_{\bar P}^\vee=1$.
\end{proof}

\begin{definition}[admissible neighbourhoods]\label{def:admissible-nbhd}
Fix the splitting and the metric data; the height functions, the family $F_*$ of
Section~\ref{subsec:compressions} and the function $\delta$ above are fixed with them, and so is
a Weinstein neighbourhood of each $L(U_i)$ and of its reflection, in which the Lagrangians of
small primary data are graphs of differentials of functions over $L(U_i)$
(Section~\ref{subsec:compressions}). All choices for the reflected splitting are the
reflections of those for the splitting. For a height $t$ let $\mathcal R(t)$ be the family of
height Lagrangians at height $t$ for $k=0$, of the presentation and of its reflection, with
$g\in\operatorname{span}F_*$ and $|g|\le\delta(t)$; it is $C^\infty$-compact and independent of
any further family. Let $\mathcal N(t)$ be an open $C^1$-neighbourhood of the pair
$(L(U_0),L(U_1))$, together with the corresponding neighbourhood of the reflected pair, each
invariant under the symplectomorphisms $\tau_c$ for the lifting classes $c$ that extend over
both compression bodies (Proposition~\ref{prop:frame-square}), the reflected part under the
reflected $\tau_c$, and so small that for every
lower pair in it and every $R\in\mathcal R(t)$: the unfolding statements of
Section~\ref{subsec:continuations} hold; the lower Lagrangians are graphs of differentials of
functions whose $C^2$-norms are so small relative to $t$ that the generators of the auxiliary
pairs are the $2^r$ nondegenerate critical points of a Morse function $C^2$-close to $t$ times
the height function, with the same indices; and
$2\max_{a,x}|\mathcal A_i(a)-\mathcal A_i(x)|<2\pi^2$ for the critical actions of the auxiliary
pairs, the values at the generators of the difference of the primitives of the two graphs
(Section~\ref{subsec:continuations}). For a finite family $F\supseteq F_*$ let $\mathcal U_F(t)$
be the path component of zero in the set of $\mathfrak h\in\operatorname{span}F$ whose
Lagrangians, and those of their reflections, lie in $\mathcal N(t)$, and let $\mathcal U'_F(t)$
be the path component of zero in $\mathcal U_F(t)\cap\mathcal U_F(t')$, where $t'\le t$ is the
height chosen in the proof of Proposition~\ref{prop:duality}. Fix once and for all a height
$t_*$ below the thresholds of Propositions~\ref{prop:small-primary} and~\ref{prop:duality}, and
for a finite family $F\supseteq F_*$ as in Definition~\ref{def:presentation} put
$\mathcal U(F)=\mathcal U'_F(t_*)$.
\end{definition}

The neighbourhood $\mathcal N(t)$ exists and depends only on the splitting, the metric data and
the height functions. The conditions imposed on a lower pair are open, and they are uniform
over the compact family $\mathcal R(t)$, which is independent of $F$ because the generic
perturbations $g$ are taken from the fixed span of $F_*$; and they concern the Lagrangians
alone: the critical actions of an auxiliary pair agree with the perturbed Chern--Simons values
at its generators (Section~\ref{subsec:continuations}), so that their differences are bounded by
the oscillation of the difference of the primitives, which is controlled by the $C^1$-size of
the Lagrangians, and at the unperturbed lower pair they are less than $2\pi^2$ by the choice of
$\delta(t)$. Intersecting with the
finitely many images under the $\tau_c$ gives the invariance. Along any path of lower data
within the preimage of $\mathcal N(t)$ the generators continue: at every point they are $2^r$
transverse intersection points, so they form a trivial covering of the path, and their relative
classes, quantized and continuous, and their indices, invariant under homotopy, are constant.
Hence the hypotheses of Lemma~\ref{lem:aux-index-folded}(ii) hold, with $c=2\pi^2$, along the
path that first moves the height datum to that of a given $R\in\mathcal R(t)$ and then follows
a smooth path of lower data in the preimage, and the conclusions of that lemma hold at every
point of $\mathcal U_F(t)$ for end and secondary perturbations below the thresholds of
Theorem~\ref{thm:mixed-iso} with $4\sup\|f\|_\infty+2\sup_A\|\eta(A)\|^2_{L^2}<2\pi^2$; by
Theorem~\ref{thm:mixed-iso}(ii) and Section~\ref{subsec:order} they then hold on homology for
the prescribed data. Propositions~\ref{prop:small-primary} and~\ref{prop:duality} use
Theorem~\ref{thm:mixed-iso} only for presentation data, so the definition is not circular.

The sets $\mathcal U_F(t)$, $\mathcal U'_F(t)$ and $\mathcal U(F)$ are open neighbourhoods of
zero in $\operatorname{span}F$: the map $\mathfrak h\mapsto(L_{\mathfrak h},L_{\bar{\mathfrak h}})$
is defined and continuous on the set of primary perturbations in $\operatorname{span}F$ in the
sense of Section~\ref{subsec:compressions}, which is open, transversality and embeddedness being
open conditions by compactness; the preimage of the open set $\mathcal N(t)$ in it is open in the
finite-dimensional space $\operatorname{span}F$, hence locally path connected, and so are its
path components; any two points of a path component are joined by a smooth path in it. They are
monotone in the family: $\mathcal U(F)\subseteq\mathcal U(F')$ for $F\subseteq F'$, since a
path in $\operatorname{span}F$ is a path in $\operatorname{span}F'$. In particular, if
$P\in\mathcal U(F)$ and $G$ is a finite family, then $P+g\in\mathcal U(F\cup G)$ for all
sufficiently small $g\in\operatorname{span}G$; and two admissible choices of primary data in
$\mathcal U(F_1)$ and $\mathcal U(F_2)$ both lie in $\mathcal U(F_1\cup F_2)$, where
Proposition~\ref{prop:small-primary} applies to them. Admissibility is thus a condition on the
primary data as functions, through the path joining them to zero within the family, and not on
the position of their Lagrangians alone: a datum whose Lagrangians lie in $\mathcal N(t_*)$ but
which is not joined to zero within the preimage is not admissible.

\begin{theorem}[one-handle square]\label{thm:one-handle}
Let $W_k\colon M\to M^+$ be a $k$-fold one-handle cobordism with the product bundle,
$p$ an admissible presentation adapted to its attaching balls, and $C\circ p$ the
enlarged presentation with height Lagrangians, with primary data in the admissible
neighbourhoods $\mathcal U(F)$ of Definition~\ref{def:admissible-nbhd} at $p$ and at
$C\circ p$. Then
$N_{C\circ p}\circ\Ish(W_k)=U_C\circ N_p$ on homology.
\end{theorem}

\begin{proof}
Choose auxiliary lower data for $p$ in $\mathcal N_{C,R}\cap\mathcal U(F)$, with
$\mathcal N_{C,R}$ the neighbourhood of Lemma~\ref{lem:aux-index-folded}(i) for the given height
Lagrangians $R$ of $C\circ p$; they exist because both sets are neighbourhoods of zero in
$\operatorname{span}F$ and small generic data are regular. At
these data Propositions~\ref{prop:cup-formula} and~\ref{prop:duality} give
$U_C=N_{C\circ p}\Ish(W_k)N^\dagger_p=N_{C\circ p}\Ish(W_k)N_p^{-1}$. At the given data,
conjugate by the continuations and apply Proposition~\ref{prop:small-primary} at
both ends.
\end{proof}

\begin{theorem}[three-handle square]\label{thm:three-handle}
Let $W_3\colon M^+\to M$ be a three-handle cobordism, read as the one-handle cobordism
$\overline W_3\colon-M\to-M^+$ of the reflected presentation. Then
$N_p\circ\Ish(W_3)=V_C\circ N_{C\circ p}$ on homology, with primary data in the admissible
neighbourhoods of Definition~\ref{def:admissible-nbhd}.
\end{theorem}

\begin{proof}
Reading the same oriented four-manifold from the opposite ends gives
$\Ish(W_3)=\Ish(\overline W_3)^\vee$: the anti-self-duality equation on $W_3$ with ends
read in reflected charts has the same index-zero moduli spaces with input and
output exchanged. By Definition~\ref{def:cup}, $V_C=(U_{\bar C})^\vee$. Transposing
Theorem~\ref{thm:one-handle} for the reflected presentation and applying
Proposition~\ref{prop:duality} at both ends,
\[
N_p\,\Ish(W_3)=(N_{\bar p}^{-1})^\vee\,\Ish(\overline W_3)^\vee
=(U_{\bar C})^\vee\,(N_{\overline{C\circ p}}^{-1})^\vee=V_C\,N_{C\circ p}. \qedhere
\]
\end{proof}

\begin{remark}
A one-handle followed by a three-handle does not cancel: the composite is
$(M\times[0,1])\#(S^1\times S^3)$. Nor do two-handle data determine the
one-handle map: if $a$ and $b$ are maps with $b\circ a=0$ and $b\circ f=1$, then
$b\circ(f+a\circ g)=1$ for every $g$. The one-handle square therefore needs its
own degeneration.
\end{remark}

\section{Presentation independence and frame changes}\label{sec:presentations}

\subsection{The presentation graph}

Let $(M,b)$ be an object. Let $\cP(M,b)$ be the graph whose vertices are
the admissible presentations of Definition~\ref{def:presentation} and whose
edges are:
\begin{enumerate}[label=(P\arabic*)]
\item\label{P:disk} changes of the meridian disk systems in a fixed
  embedded splitting (Definition~\ref{def:presentation} records no disk
  systems, so these edges are loops at a vertex; they are listed to match
  the edges of \cite{JTZ});
\item\label{P:iso} ambient isotopies of the splitting fixed near
  $b(B^3)$;
\item\label{P:stab} stabilization and destabilization away from
  $b(B^3)$;
\item\label{P:model} changes of the metric data, of the primary
  data and of the regularizing data (Definition~\ref{def:presentation}(v)),
  including the almost complex structure on the strip end, of a fixed
  splitting.
\end{enumerate}

\begin{proposition}\label{prop:pres-graph}
$\cP(M,b)$ is connected.
\end{proposition}

\begin{proof}
Remove the interior of $b(B^3)$ and put an annular suture around the
equator $\partial D$ of its boundary sphere. The result is a balanced
sutured manifold $M(b)$ whose positive and negative regions are disks.
Removing $D$ from an embedded Heegaard surface adapted to $b$ gives a
sutured Heegaard diagram with one boundary circle, and conversely. By
\cite[Definitions~2.34--2.35, Proposition~2.36]{JTZ}, any two vertices of the graph of embedded
sutured diagrams of $M(b)$ are joined by a path of edges given by
compression-equivalent changes of attaching sets, stabilizations, and
diffeomorphisms isotopic to the identity in $M(b)$ relative to the
suture \cite[Definition~2.34]{JTZ}. The first are~\ref{P:disk}. The last are
diffeomorphisms isotopic to the identity relative to the suture, which need
not preserve the standard collar of the boundary sphere; since the image
surface meets the boundary sphere transversely in the same circle, the collar
neighbourhood theorem gives an ambient isotopy, fixed on a smaller collar, that
makes it standard on the collar again, and the composite is an edge~\ref{P:iso}. The
stabilizations of \cite[Definition~2.18]{JTZ} need not be local: the
punctured torus $T$ replacing the disk $D$ is not required to lie in a
small ball. Nevertheless every such edge $H_1\to H_2$ factors in
$\cP(M,b)$ as a local stabilization~\ref{P:stab} in a ball $Q$ with $Q\cap
\Sigma_1\subset\operatorname{int}D$, followed by a change of disk system and
an isotopy fixed near $b(B^3)$ (Lemma~\ref{lem:stab-factor}). An additional stabilization, in a small ball of the standard collar, makes the genus at least one throughout: fix a small ball $Q_0$
in the common standard collar, where all diagrams of the path are standard and
all isotopies are the identity, and stabilize every vertex of the path in
$Q_0$. The stabilized vertices form a path of the same kind, and the two ends
are joined to their unstabilized versions by edges~\ref{P:stab}. Hence the
torus genus is at least two throughout. For a fixed splitting,
shrink the admissible collars so that the fixed height solid tori
remain in the compression interiors. Cutting along the framing
sphere $S$ gives neat collars on the two sides, matched along the
collar of $\partial D$. Uniqueness of neat collars up to relative
isotopy, followed by isotopy extension fixing $\Sigma_h$ and
preserving $S$, joins the collar choices away from the fixed torus
collar. Product metrics on a fixed collar form a convex set. The admissible
neighbourhoods $\mathcal U(F)$ of Definition~\ref{def:admissible-nbhd} are path connected,
being path components. An edge
\ref{P:model} joins admissible endpoints within these neighbourhoods;
its intermediate primary pairs need not stay transverse. Thus the
choices are joined by edges~\ref{P:model}.
\end{proof}

\begin{lemma}\label{lem:stab-factor}
Let $H_1=(\Sigma_1,\alpha_1,\beta_1)$ and $H_2=(\Sigma_2,\alpha_2,\beta_2)$ be
embedded diagrams of $M(b)$, standard on a common collar of the boundary,
such that $H_2$ is obtained from $H_1$ by a stabilization with data
$D\subset\Sigma_1$, $T\subset\Sigma_2$, $a=\alpha_2\cap T$, $b=\beta_2\cap T$. Then
there are a ball $Q$ disjoint from the curves of $H_1$ and from $b(B^3)$,
with $Q\cap\Sigma_1\subset\operatorname{int}D$, and an isotopy $\psi_t$ of $M(b)$
with compact support in the interior, such that $\psi_1(\Sigma_2)$ is the
local stabilization of $\Sigma_1$ in $Q$, with the same sides.
\end{lemma}

\begin{proof}
The curves $a$ and $b$ bound disks $D_a$, $D_b$ on the two sides of $\Sigma_2$,
meeting in one point. Destabilizing $\Sigma_2$ along a regular neighbourhood $N$ of
$D_a\cup D_b$, chosen to meet $\Sigma_2$ exactly in $T$, gives a surface $\Sigma''$ that agrees with $\Sigma_1$ outside $D$. Then
$\Sigma''\setminus(\Sigma_1\setminus D)$ and $D$ are two disks with the same boundary. They
lie in the union of the two sides of $\Sigma_1$ glued along $D$, a boundary sum of
irreducible compression bodies, which is irreducible. After an innermost-disk argument
removing the intersections of their interiors, the sphere they form bounds a ball, and the
two disks are isotopic relative to their boundary, by an isotopy
supported in the interior and away from the curves of $H_1$ and from $b(B^3)$. This
isotopy carries $N$, which contains $T$, into a ball $Q$, which a final small isotopy near
$D$ places with $Q\cap\Sigma_1\subset\operatorname{int}D$, and carries $\Sigma_2$ to the local stabilization of $\Sigma_1$ in $Q$,
with the same sides. The framing ball and the torus are never moved.
\end{proof}

\begin{lemma}\label{lem:output-identification}
Let $Q\subset M\setminus b(B^3)$ be a ball, $W_1\colon M\to M_m$ a one-handle
and $W_2\colon M_m\to M'$ a two-handle attached in $Q\times[0,1]$ whose
attaching circle meets the belt sphere of $W_1$ transversely once, with product
bundles and vertical tube. For every diffeomorphism $\iota\colon M'\to M$ equal
to the identity outside $Q$, $(W_2,\iota)\circ W_1\cong1_M$ in $\bB$.
\end{lemma}

\begin{proof}
Give $W_2\circ W_1$ a Morse function with the two critical points and a
gradient-like field equal to $\partial_t$ on $(M\setminus Q)\times[0,1]$. By the
first cancellation theorem \cite[Theorem~5.4]{Milnor-h-cob} it can be replaced
by one without critical points, the change supported near the connecting
trajectory between the two critical points, which misses $(M\setminus Q')\times
[0,1]$ for a slightly larger ball $Q'$. Integrating the new field gives a
diffeomorphism to $M\times[0,1]$ equal to the identity on the incoming end and
on $(M\setminus Q')\times[0,1]$. The induced self-diffeomorphism of $M$ at the
outgoing end, composed with $\iota^{-1}$, is supported in $Q'$, hence isotopic
to the identity relative to the boundary of $Q'$ by Cerf's theorem
\cite{Cerf}, and the mapping cylinder of the isotopy is an isomorphism of
decorated cobordisms fixed on the tube.
\end{proof}

\subsection{The edge squares}

For each edge $e\colon p\to q$ we define an independent symplectic map
$\Gamma_e\colon\HF(M,p)\to\HF(M,q)$ and prove
\begin{equation}\label{eq:edge-square}
H(N_q)\circ G(q,p)=\Gamma_e\circ H(N_p),
\end{equation}
where $G(q,p)$ is the gauge continuation map, i.e.\ $\Ish$ of the product
cobordism with vertical tube and product bundle. In the notation of
Section~\ref{subsec:continuations}, $G(q,p)=\mathrm{cont}_G(p,q)$; here the arguments are written
in the order of composition, so that $G(r,q)G(q,p)=G(r,p)$.

\begin{proposition}\label{prop:pres-squares}
Every edge of $\cP(M,b)$ has an independently defined map $\Gamma_e$
satisfying~\eqref{eq:edge-square}.
\end{proposition}

\begin{proof}
\ref{P:disk}: The Lagrangians $L(C_i)$ depend on the embedded compression
bodies, their metrics and primary perturbations, and not on the disk
systems, so $\Gamma_e=1$ and~\eqref{eq:edge-square} is trivial.

\ref{P:model}, primary data: $\Gamma_e$ is the symplectic continuation of
Section~\ref{subsec:continuations}, and~\eqref{eq:edge-square} is
Proposition~\ref{prop:small-primary} inside one finite family; two
finite families are bridged through their union, which is again a finite
family containing both.

\ref{P:model}, regularizing data: for changes of $f$, $\chi$, $\eta$ the map
$\Gamma_e$ is the identity and the square is Theorem~\ref{thm:mixed-iso}(ii). For a
change of the almost complex structure $J$ on the strip end, $\Gamma_e$ is the
symplectic continuation. Insert a strip carrying the continuation from $J$ to
$J'$ at the ordinary end of the mixed domain; by Theorem~\ref{thm:neck} at the
ordinary neck, for a long neck the count on the glued domain is $\Gamma_e\circ N_p$,
and the glued domain is a compactly supported deformation of the domain for
$J'$, so Theorem~\ref{thm:family} identifies the count with $N_{p'}$ up to
chain homotopy.

\ref{P:stab}: Let $r$ be the local stabilization of $p$ in a ball $Q$. Let
$W_1\colon M\to M\#(S^1\times S^2)$ be the one-handle and $W_2$ the cancelling
two-handle, both attached in $Q\times[0,1]$. By
Lemma~\ref{lem:output-identification}, $W_2\circ W_1\cong1_M$ for any
identification of the outgoing end that is the identity outside $Q$, so
$G(r,p)=\Ish(W_2)\Ish(W_1)$. The one-handle square
(Theorem~\ref{thm:one-handle}) ends at the enlarged presentation $C\circ p$,
with compression bodies $(V_0,V_1)$ carrying height Lagrangians $(R_0,R_1)$ at a
height $t$, while the two-handle square (Theorem~\ref{thm:two-handle}) starts
at a vertex whose second body is unperturbed. The two are joined by
primary-data changes, as follows. Give $C\circ p$ the metric data that
agrees on $V_0$ with the metric of $r$ on its first body; the output of the
cancelling two-handle has bodies $(V_0,C_2)$, where $C_2$ is the second body of
$r$ with the metric of $r$. Take first the heights $t_*$ of Definition~\ref{def:admissible-nbhd} for the splittings and metric data of
$C\circ p$ and of $r$, then $t$ and the two-handle height $t_2$ so small that
all data below lie in the corresponding neighbourhoods and in the admissible
neighbourhoods, and then the lower Lagrangians of $p$ as in
Theorem~\ref{thm:one-handle}. Let $R_0'$ be a height Lagrangian obtained from $R_0$
by a small generic perturbation from $\operatorname{span}F_*$ as in
Section~\ref{subsec:one-three-squares}, transverse to
$L(V_1)$ and to $L_{t_2h_Q}(C_2)$ and missing their intersection; this is a
generic condition, and the perturbation may move the generators of $R_0$, whose labels are
continued. Put
\begin{itemize}[leftmargin=2em]
\item $(C\circ p)'$: the presentation $C\circ p$ with primary data $(R_0',0)$;
\item $r'$: the presentation $r$ with primary data $(R_0',t_2h_Q)$.
\end{itemize}
Both are admissible. The subordinate triple of the cancelling two-handle is
$(C_0,C_1,C_2)=(V_0,V_1,C_2)$ with primary data $(R_0',0,t_2h_Q)$, which
satisfies the hypotheses of Section~\ref{subsec:triple-data}, so
Theorem~\ref{thm:two-handle} holds from $(C\circ p)'$ to $r'$. Define
\[
S_{r,p}=\mathrm{cont}_S(r',r)\circ S_{W_2}\circ \mathrm{cont}_S(C\circ p,(C\circ p)')\circ U_C .
\]
The squares of Theorem~\ref{thm:one-handle}, of
Proposition~\ref{prop:small-primary} for the two changes of primary data, and of
Theorem~\ref{thm:two-handle} compose to $H(N_r)G(r,p)=S_{r,p}H(N_p)$.

The reverse composite is a two-handle $W_2'\colon M\to M\#(S^1\times S^2)$
followed by a three-handle $W_3$, both in $Q\times[0,1]$. Reversing the
orientation of the cobordism exchanges them with a one-handle followed by a
cancelling two-handle, so Lemma~\ref{lem:output-identification} applied to the
reversed cobordism gives $W_3\circ W_2'\cong1_M$ and $G(p,r)=\Ish(W_3)\Ish(W_2')$.
The triple of $W_2'$ is $(V_0,C_2,V_1)$ with primary data $(R_0'',0,t_2'h_Q')$, where $h_Q'$ is the corresponding height function on $V_1$,
chosen in the same way, and
\[
D_{p,r}=V_C\circ \mathrm{cont}_S\bigl((C\circ p)'',C\circ p\bigr)\circ S_{W_2'}\circ
\mathrm{cont}_S(r,r''),
\]
where $r''$ carries $(R_0'',0)$ and $(C\circ p)''$ carries $(R_0'',t_2'h_Q')$.
Theorems~\ref{thm:two-handle} and~\ref{thm:three-handle} and
Proposition~\ref{prop:small-primary} give $H(N_p)G(p,r)=D_{p,r}H(N_r)$.
Then $D_{p,r}S_{r,p}H(N_p)=H(N_p)G(p,r)G(r,p)=
H(N_p)$, and since $H(N_p)$ is invertible, $D_{p,r}S_{r,p}=1$; similarly
$S_{r,p}D_{p,r}=1$.

\ref{P:model}, metric data: for two presentations $p,q$ of one
splitting that differ in their metric data, choose a common stabilization $r$ and put
$\Gamma_e=D_{q,r}S_{r,p}$. This is possible because the upper metric data in
Theorems~\ref{thm:two-handle} and~\ref{thm:one-handle} is free within the class
of Section~\ref{subsec:folded}: the flat-connection Lagrangians do not depend on
the metric, so one upper choice serves both $p$ and $q$. The square follows from the stabilization squares
and $G(q,r)G(r,p)=G(q,p)$; independence of $r$ follows by invertibility of
$N$.

\ref{P:iso}: An ambient isotopy $\phi_t$ fixed near $b(B^3)$ transports all
data. Transporting data conjugates every equation, operator and
perturbation, so the counts agree exactly, and the mapping cylinder of
the isotopy is the decorated product. The isotopy transports the primary data as well as the
metrics, so the transported vertex differs from the prescribed one in both.
Hence $\Gamma_e$ is the map induced by the transported data followed by the
metric-data change and then the primary-data change~\ref{P:model}, and the square
follows from the previous cases.
\end{proof}

\begin{corollary}\label{cor:transitive}
For paths $\gamma$ in $\cP(M,b)$ from $p$ to $q$, the composite
$\Gamma_\gamma$ depends only on $p$ and $q$; the maps $\Gamma(q,p)$ satisfy
$\Gamma(r,q)\Gamma(q,p)=\Gamma(r,p)$ and $\Gamma(p,p)=1$.
\end{corollary}

\begin{proof}
By~\eqref{eq:edge-square} and functoriality of gauge continuation,
$\Gamma_\gamma H(N_p)=H(N_q)G(q,p)$, and $H(N_p)$ is invertible.
\end{proof}

\subsection{Frame changes}

\begin{proposition}\label{prop:frame-square}
Let $u\colon M\to\SO(3)$ equal $1$ near $b(B^3)$, and let $p$, $q$ be admissible
presentations. The restriction $u|_\Sigma$, extended by $1$ over the torus
summand, is an automorphism $\hat u$ of the odd surface bundle that extends
over both compression bodies; its lifting class $c(\hat u)$ therefore lies in
$\Lambda_0\cap \Lambda_1$, and the induced symplectomorphism $\tau_{c(\hat u)}$ of $\Msf_h$
carries the Lagrangians of $p$ to those of the transported presentation
$u_*p$. Let $u_S$ be the induced pushforward chain isomorphism and
$S_{[u]}=\Gamma(q,u_*p)\circ u_S$. Then $H(N_q)\circ\Ish([u])=S_{[u]}\circ H(N_p)$.
\end{proposition}

\begin{proof}
On chains $N_{u_*p}\circ u_G=u_S\circ N_p$, where $u_G$ is the gauge
pushforward: applying the bundle automorphism $u$ to all data transports
the equations of the mixed domain, the compression bodies, the torus
marking and the tube (on which $u=1$), and preserves the lifting group,
whose support is disjoint from that of $u$. The transported presentation
$u_*p$ is admissible for the family $F\cup u_*F$, with the transported mixed map: the path
joining the primary data of $p$ to zero in $\operatorname{span}F$ is transported to a path in
$\operatorname{span}u_*F$ whose Lagrangians, the images under $\tau_c$ of the original ones, and
their reflections lie in $\mathcal N(t_*)\cap\mathcal N(t_*')$, these neighbourhoods being
invariant under the symplectomorphisms $\tau_c$ for these classes $c$
(Definition~\ref{def:admissible-nbhd}). The gauge map of the
frame change is $\Ish([u])=G(q,u_*p)\circ u_G$. Hence, by
Corollary~\ref{cor:transitive},
$H(N_q)\Ish([u])=H(N_q)G(q,u_*p)u_G=\Gamma(q,u_*p)H(N_{u_*p})u_G=\Gamma(q,u_*p)u_SH(N_p)$.
\end{proof}

A lifting class in $\Lambda_0+\Lambda_1$ would not suffice: if $a\in \Lambda_0\setminus \Lambda_1$, then
$\tau_a$ changes the sign of the lifted holonomy of a meridian of $H_1$ and
moves $L(C_1)$ off itself. The condition $\zeta\in \Lambda_0+\Lambda_1$ on the filling-frame class $\zeta$ is the triviality
condition for the bundle on $H_0\cup H_1$ and is used in the construction of
the global frames of Section~\ref{subsec:two-handle-bundles}.

\section{Proof of Theorem~\ref{thm:A}}\label{sec:assembly}

\subsection{A descent lemma}

Let $V(M)$ denote the set of admissible presentations of an object $M$
(Section~\ref{sec:models}). Consider the category $\bB^V$ whose objects
are pairs $(M,p)$ with $p\in V(M)$ and whose morphisms
$(M,p)\to(M',q)$ are the morphisms $M\to M'$ of $\bB$. The projection
$\bB^V\to\bB$ is an equivalence. For each $(M,p)$ the gauge complex
$C_G(M,p)$ computes $\Ish(M)$: the gauge continuation maps provide
canonical isomorphisms $H(C_G(M,p))\cong\Ish(M)$ compatible with the
maps $\Ish(W,E,\nu)$ (Proposition~\ref{prop:gauge-functor}). We
therefore regard $\Ish$ as a functor on $\bB^V$ with
$\Ish(M,p)=H(C_G(M,p))$, and write $N_{(M,p)}=H(N_p)\colon
\Ish(M,p)\to\HF(M,p)$.

An \emph{elementary arrow} $e\colon(M,p)\to(M',q)$ is a triple consisting
of an underlying morphism $[e]\in\bB(M,M')$, source and target
presentations, and a linear map $S_e\colon\HF(M,p)\to\HF(M',q)$ defined
by holomorphic curve counts. The arrows used are listed in
Section~\ref{subsec:generation}.

\begin{lemma}\label{lem:descent}
Let $\mathcal A$ be a set of elementary arrows such that:
\begin{enumerate}[label=(\alph*)]
\item each $N_{(M,p)}$ is an isomorphism;
\item every $e\in\mathcal A$ satisfies
  $N_{(M',q)}\circ\Ish([e])=S_e\circ N_{(M,p)}$;
\item for all $(M,p)$, $(M',q)$ and every $W\in\bB(M,M')$ there is a
  word $w=(e_1,\dots,e_m)$ of elements of $\mathcal A$ from $(M,p)$ to
  $(M',q)$ with $[e_m]\circ\dots\circ[e_1]=W$.
\end{enumerate}
For a word $w$ put $S_w=S_{e_m}\circ\dots\circ S_{e_1}$. Then
\begin{enumerate}[label=(\roman*)]
\item $N_{(M',q)}\circ\Ish(W)=S_w\circ N_{(M,p)}$ for every word $w$
  representing $W$;
\item $S_w$ depends only on $W$, $p$ and $q$; write it $S(W)_{qp}$;
\item $W\mapsto S(W)_{qp}$ is a functor $\bB^V\to\mathrm{Vect}_{\FF}$,
  and $N$ is a natural isomorphism $\Ish\Rightarrow S$;
\item the maps $\Gamma(q,p)=S(1_M)_{qp}$ satisfy $\Gamma(r,q)\Gamma(q,p)=\Gamma(r,p)$ and
  $\Gamma(p,p)=1$. With $\HF^\#(M)=\operatorname{colim}_p\HF(M,p)$ along
  $\Gamma$, the functor $S$ and the transformation $N$ descend to a functor
  $\HF^\#$ on $\bB$ and a natural isomorphism $\Ish\Rightarrow\HF^\#$.
\end{enumerate}
\end{lemma}

\begin{proof}
(i) Induction on $m$. For $m=0$ the word is empty, $W=1_M$ and $p=q$.
For $m\ge1$ let $w'=(e_1,\dots,e_{m-1})$ end at $(M'',r)$. By
functoriality of $\Ish$, (b) for $e_m$, and the induction hypothesis,
\[
N_q\Ish(W)=N_q\Ish([e_m])\Ish([w'])=S_{e_m}N_r\Ish([w'])
=S_{e_m}S_{w'}N_p.
\]
(ii) By (i) and (a), $S_w=N_q\Ish(W)N_p^{-1}$, which does not involve
$w$. (iii) $S(W'\circ W)=N\Ish(W')\Ish(W)N^{-1}=S(W')S(W)$ and
$S(1)=1$; naturality is (i). (iv) is (iii) applied to identity
morphisms; a system of isomorphisms satisfying the cocycle law has a
colimit with isomorphisms $i_p\colon\HF(M,p)\to\HF^\#(M)$, and
$N_M=i_pN_{(M,p)}$ is independent of $p$ by (i).
\end{proof}

Lemma~\ref{lem:descent} is formal. All of the content of
Theorem~\ref{thm:A} lies in (a), (b) and (c). In (b) the maps
$N_{(M,p)}$ must be the ones fixed at the vertices: a square proved with
an auxiliary mixed map at the source or target is usable only after
that map has been identified, on homology, with the mixed map. This is
the role of Theorem~\ref{thm:mixed-iso}(ii) and of the data-change squares of Section~\ref{sec:presentations}.

\subsection{Generation}\label{subsec:generation}

The set $\mathcal A$ consists of the following arrows.
\begin{enumerate}[label=(A\arabic*)]
\item\label{arr:pres} \emph{Presentation moves} $(M,p)\to(M,q)$ with
  $[e]=1_M$: changes of disk system, stabilization and destabilization
  away from the ball, isotopies fixed near the ball, and changes of the metric and primary data (Section~\ref{sec:presentations}).
\item\label{arr:frame} \emph{Frame changes} $[u]$
  (Section~\ref{subsec:frame-change}), with $S_{[u]}$ the map induced by
  the symplectomorphism of $\Msf_h$ given by $u|_\Sigma$
  (Proposition~\ref{prop:frame-square}).
\item\label{arr:one} \emph{One-handles} $(M,p)\to(M^+,C\circ p)$, where $C\circ p$
  carries height Lagrangians at some height $t$ and $p$, $C\circ p$ satisfy the
  neighbourhood conditions of Theorem~\ref{thm:one-handle}, with the
  compression cup.
\item\label{arr:two} \emph{Two-handle cobordisms} $(W_2,E;\phi_0^{\mathrm{gl}},
  \phi_1^{\mathrm{gl}})$ with bundles trivial at both ends, between the
  presentations $p_0^\phi$, $p_1^\phi$ of a subordinate triple with the
  global frames of Section~\ref{subsec:two-handle-bundles}, carrying the
  transported data of Section~\ref{subsec:triple-data} with $t$ below the
  threshold of Lemma~\ref{lem:aux-index}(iii) and so small that the primary data are joined
  to zero within the admissible neighbourhoods; these vertices are admissible
  presentations (Definition~\ref{def:admissible-nbhd}). The map is $S_{W_2,E}$ (Theorem~\ref{thm:two-handle}).
\item\label{arr:three} \emph{Three-handles} $(M^+,C\circ p)\to(M,p)$, with the
  same vertex conditions as in~\ref{arr:one}, with the compression cap
  (Theorem~\ref{thm:three-handle}).
\end{enumerate}

\begin{proposition}\label{prop:generation}
$\mathcal A$ satisfies condition~(c) of Lemma~\ref{lem:descent}.
\end{proposition}

\begin{proof}
Let $(W,E,\nu)\colon M\to M'$ and presentations $p,q$ be given. By
Propositions~\ref{prop:handles} and~\ref{prop:outer-bundles},
$(W,E,\nu)=W_3\circ(W_2,E)\circ W_1$ with product bundles on $W_1$ and
$W_3$ and connected levels $Z_0,Z_1$. Write $W_1$ and $W_3$ as composites
of single one- and three-handles. By Lemma~\ref{lem:two-handle-frames},
$(W_2,E)=[u_1]\circ(W_2,E;\phi^{\mathrm{gl}})\circ[u_0]^{-1}$, a two-handle
arrow between frame changes, or a single frame change if there are no
two-handles. Each handle can be
attached in a presentation adapted to it
(Section~\ref{subsec:triples} and the paragraph preceding it), and the
two-handle cobordism in the presentation given by a subordinate triple
(Proposition~\ref{prop:tube}). By Proposition~\ref{prop:pres-graph},
any two presentations of the same object, and in particular consecutive
presentations at each level as well as $p$, $q$ at the ends, are joined
by presentation moves; in particular each handle arrow can be taken between
vertices satisfying the conditions of~\ref{arr:one}--\ref{arr:three}, since
such vertices exist for every splitting (small heights and generic small data) and are reached by presentation moves. Concatenation gives a word of arrows in
$\mathcal A$ from $(M,p)$ to $(M',q)$ whose underlying composite is
$(W,E,\nu)$, since presentation moves have underlying morphism the
identity.
\end{proof}

\subsection{Proof of Theorem~\ref{thm:A}}

Condition~(a) of Lemma~\ref{lem:descent} is
Theorem~\ref{thm:mixed-iso}. Condition~(b) is
Proposition~\ref{prop:pres-squares} for~\ref{arr:pres},
Proposition~\ref{prop:frame-square} for~\ref{arr:frame},
Theorem~\ref{thm:one-handle} for~\ref{arr:one}, at presentations with primary data in the
admissible neighbourhoods of Definition~\ref{def:admissible-nbhd} (arrows may be taken at any
convenient vertices, since presentation moves connect them),
Theorem~\ref{thm:two-handle} for~\ref{arr:two} and
Theorem~\ref{thm:three-handle} for~\ref{arr:three}. Condition~(c) is
Proposition~\ref{prop:generation}. Lemma~\ref{lem:descent} gives the
equation of Theorem~\ref{thm:A} for every word, the independence of
$S_d$, the functor $\HF^\#$ and the natural isomorphism. The statement
in the homotopy category of finite complexes follows from
Lemma~\ref{lem:field}. \qed

\begin{remark}
No relation among symplectic words is proved directly. Independence of
the decomposition follows from the squares, the connectivity of the
presentation graph and the invertibility of $N$. The argument does not
give integral statements: over $\ZZ$ the squares would be needed with
signs, and Lemma~\ref{lem:field} fails.
\end{remark}

\subsection{The isomorphism of Daemi, Fukaya and Lipyanskiy}

\begin{corollary}\label{cor:dfl-natural}
Every object has admissible presentations with DFL-regular primary data, at every Heegaard
splitting and all metric data. Let $p_0$, $p_1$ be admissible presentations of $M_0$, $M_1$ with DFL-regular primary
data, and let $N^{\rm DFL}_{p_i}$ be the isomorphisms of \cite[Theorem~3.29]{DFL} at the
data of Proposition~\ref{prop:dfl-comparison}. Then for every morphism $(W,E,\nu)$ and
every admissible decomposition $d$,
\[
H(N^{\rm DFL}_{p_1}\otimes\FF)\circ\Ish(W,E,\nu)=S_d(W,E,\nu)\circ H(N^{\rm DFL}_{p_0}\otimes\FF).
\]
In particular, for two such choices of data at presentations $p,q$ of one object,
$H(N^{\rm DFL}_q\otimes\FF)\circ G(q,p)=\Gamma(q,p)\circ H(N^{\rm DFL}_p\otimes\FF)$ as maps from
the gauge homology at $p$, where $G(q,p)$ is the gauge continuation and $\Gamma(q,p)$ the
symplectic transitive system of Theorem~\ref{thm:A}, defined by holomorphic curve counts.
So within this class of data, including changes of the Heegaard splitting, the reduction
modulo two of the isomorphism of \cite{DFL} on homology does not depend on the choices.
\end{corollary}

\begin{proof}
The first sentence is Theorem~\ref{thm:dfl-regular}, applied with any separated family as in
Definition~\ref{def:presentation} and generic primary data near zero in $\mathcal U(F)$.
By Proposition~\ref{prop:dfl-comparison}, $N^{\rm DFL}_{p_i}\otimes\FF=N_{p_i}$ at these data,
and Theorem~\ref{thm:A} applies to $N_{p_i}$. The statement for two choices of data at one
object is the case of the identity morphism: on the gauge homology at $p$ and $q$ the identity of $\Ish(M)$ is the
gauge continuation $G(q,p)$, and $S(1_M)_{qp}=\Gamma(q,p)$ by Lemma~\ref{lem:descent}.
\end{proof}

Let $p$ be an admissible presentation of $M$, with family $F_p$, and let $q$ be an
admissible presentation with the same splitting and metric data and DFL-regular primary data,
admissible for a family $F_q$ (Theorem~\ref{thm:dfl-regular}). The case $W=1_M$ of
Theorem~\ref{thm:A}, with Proposition~\ref{prop:dfl-comparison} at $q$, gives
\[
H(N_p)=\Gamma(p,q)\circ H(N^{\rm DFL}_q\otimes\FF)\circ G(q,p)
\]
on homology, where $G(q,p)$ is the gauge continuation map. Both primary data lie in
$\mathcal U(F_p\cup F_q)$, by the monotonicity of the admissible neighbourhoods
(Definition~\ref{def:admissible-nbhd} and the paragraph following it), so
Proposition~\ref{prop:small-primary} applies to them in
the family $F_p\cup F_q$; with the invertibility of $H(N_p)$ it identifies $\Gamma(q,p)$ with the
symplectic continuation map. So the mixed map at
any admissible presentation is, on homology, the reduction modulo two of the isomorphism of
\cite{DFL} at an admissible presentation with the same splitting and metric data, composed with
gauge and symplectic continuation maps.

\appendix
\section{Results used from the literature}\label{app:inputs}

We list the published results used in the proof, where they enter, and
the hypotheses that have to be checked.

\begin{enumerate}[label=(\arabic*),leftmargin=2em]
\item \cite[Theorem~1, Propositions~2.21, 3.19, 3.20, Lemma~3.15,
  Theorem~3.29]{DFL}: monotonicity with minimal Maslov number four,
  embeddedness of $L(C)$, the index of constant pairs, the index-energy
  law and triangularity of the mixed map. Used in
  Sections~\ref{sec:models} and~\ref{sec:collapse}. Hypotheses: the
  bundles are odd on the matching surface; the Lagrangians are embedded,
  which holds for compression bodies.
\item \cite[Theorem~4.17, Remark~4.19, Section~4.2.1, (4.20)--(4.24),
  Theorem~5.2, Lemma~5.23, Proposition~5.29, Section~5.4, Section~6]{DFL}:
  local elliptic estimates for the mixed operator, theory of mixed ends,
  small-energy regularity, exponential decay, compactness and
  perturbations. Used throughout. Hypotheses: product metrics and the complex structure $*_\Sigma$ near the seams, which Convention~\ref{conv:orientation}
  and Hypothesis~\ref{hyp:collapse} arrange.
\item \cite[Section~2.3, Proposition~2.53, Section~3.1 and Figure~1, Definitions~3.8--3.9,
  (3.24), Proposition~3.18]{DFL}: primary perturbations and the submersion property of
  the universal restriction map, the special quintuple and its planar region, the
  domains of DFL type with their weighted spaces, and the Fredholm property of the
  mixed operator. Used in Sections~\ref{sec:models}, \ref{sec:collapse}
  and~\ref{sec:two-handles}. Hypotheses: those of~(1).
\item \cite[Section~4.1 with (4.3), (4.5), (4.7), (4.8), Lemmas~4.1, 4.2 and~4.4,
  Proposition~4.6, Section~4.2.1 with (4.16) and (4.24), Propositions~2.10, 4.65, 4.69,
  5.27, 5.28 and~5.29]{DFL}: the chart $P_p$ and the cutoff maps $E_\tau$, the flat
  extension of \cite[Lemma~4.2]{DFL}, a metric in which the Lagrangian is totally
  geodesic, the Coulomb slice and its variant, the Coulomb gauge on the three-manifold,
  the slice operator on a cylinder quintuple, additivity of energy and index under
  cutoff gluing, and exponential decay on mixed ends. Used in
  Sections~\ref{sec:families} and~\ref{sec:torus}. Lemma~4.1 is used with its radius
  measured in $L^q$ for some $q>2$ (Remark~\ref{rem:dfl-radius}). Hypotheses:
  nondegenerate or Morse--Bott ends as stated there.
\item \cite[Theorems~1 and~3, Theorem~4.2, Lemmas~4.25 and~4.50, Theorem~4.44]{DFL-mixed}:
  seam regularity, compactness, seam estimates and removal of seam
  singularities for the mixed equation.
  Hypotheses: fixed product collars after conformal rescaling; the
  disconnected odd matching surface is allowed.
\item \cite[Sections~5.1--5.2, Definition~5.1, Propositions~5.2
  and~5.4, Theorem~5.6]{KM-unknot}: the enlarged gauge groups, the
  quotient by $\phi$, the composition law with the connecting-map
  condition, and excision. Used in Sections~\ref{sec:category},
  \ref{sec:two-handles} and~\ref{sec:torus}. Hypotheses: the
  $\phi$-non-integral condition at every end, and the connecting image
  $\psi$ contained in $\phi$; both are checked where used.
\item \cite[Sections~7.1--7.4, 7.6--7.8, Proposition~7.1]{Scaduto}: the
  framed gauge transformations, the $\ZZ/4$ grading and degree formula,
  and $\Ish(\#^r(S^1\times S^2))$ with its relative invariants.
\item \cite[Definition~4.2, Proposition~4.3, Lemma~4.5]{OS-triangles}:
  subordinate Heegaard triples and the identification of the filled
  triangle with the two-handle trace; the tube is controlled by
  Proposition~\ref{prop:tube}.
\item \cite[Definitions~2.34--2.35, Proposition~2.36]{JTZ}: connectivity of the
  graph of embedded sutured diagrams; applied to the sutured manifold
  $M(b)$ in Proposition~\ref{prop:pres-graph}.
\item \cite[Theorems~3.9, 3.13, 4.1, 4.2, Example~2.6 and Corollary~5.5]{WW-relative}:
  quilted Floer maps, gluing, data-change isomorphisms, the strip example
  (Example~2.6) and unit insertion (Corollary~5.5). Hypotheses: monotonicity with a common constant and minimal
  Maslov number at least three for all Lagrangians and seams, and
  Hamiltonian one-forms vanishing on the boundary; the continuation of
  Section~\ref{subsec:continuations} is arranged in this form.
\item \cite[Section~14.4]{Oh-book}: continuation with moving Lagrangian
  boundary conditions.
\item \cite[Propositions~3.8--3.10, 4.3, 4.4 and~4.16, Theorems~4.17 and~5.4,
  Sections~5.3 and~5.5]{Donaldson-book}: additivity of the index under gluing along
  acyclic ends, weights and kernels, exponential decay of solutions of small energy
  on long tubes, also for the holonomy-perturbed equation, neck gluing with
  surjectivity, and gluing at index $-1$ in families;
  the addition property of \cite[Section~3.3]{Donaldson-book}, on which the excision argument of Lemma~\ref{lem:excision} is modelled.
\item Gluing along ordinary strip necks with transverse intersections is
  proved in Theorem~\ref{thm:neck} by the argument of
  \cite[Section~3]{Donaldson-book} with the local estimates of
  \cite[Section~4]{DFL-mixed} and the exponential decay of holomorphic strips
  of small energy \cite[Proposition~2.3.1]{Simcevic}, whose hypotheses (H) and (H1),
  closed Lagrangians with compact clean intersection, hold here; monotonicity enters
  only through compactness. For exact Lagrangians the corresponding statement is
  \cite[Proposition~2.22]{Seidel-LES}.
\item \cite[Theorem~B.4.1]{MS}: boundary regularity for $J$-holomorphic
  curves with Lagrangian boundary conditions, from $W^{1,p}$, $p>2$.
\item \cite[Proposition~8.6, Theorem~C.1]{SW-lagrangian}: unique
  continuation for the trajectory equation.
\item \cite{BB}, \cite[Corollary~2.1.6, Theorem~2.1.7]{Schwarz}, \cite{Robbin-Salamon},
  \cite{Aronszajn}: elliptic boundary value problems for first-order
  operators, the Gaffney--Friedrichs inequality, regularity for
  families of self-adjoint operators, and unique continuation.
\end{enumerate}
The quilted description of the compression cup in
Remark~\ref{rem:quilted-cup} uses in addition \cite[Corollary~5.5]{WW-relative},
\cite{WW-quilt}, \cite{WW-composition} and \cite[Theorem~2 and Corrigendum,
Remark~1(1)]{LL-composition}; none of these enters the proofs of
Theorems~\ref{thm:A} and~\ref{thm:B}.

\section{Homological algebra over a field}\label{app:algebra}

\begin{lemma}\label{lem:field}
Let $(C,d_C)$ and $(D,d_D)$ be finite-dimensional chain complexes over a
field $k$, and $f,g\colon C\to D$ chain maps with $H(f)=H(g)$. Then $f$
and $g$ are chain homotopic. Consequently, the homology functor from the
homotopy category of finite-dimensional complexes over $k$ to graded
vector spaces is fully faithful, and an equation between composites of
chain maps that holds on homology holds up to chain homotopy.
\end{lemma}

\begin{proof}
It suffices to treat $g=0$. Over a field every finite complex splits: there
are maps $i_C\colon H(C)\to C$, $p_C\colon C\to H(C)$ and
$h_C\colon C\to C$ of degree one with $p_Ci_C=1$,
$1-i_Cp_C=d_Ch_C+h_Cd_C$, and $d_Ci_C=0$, $p_Cd_C=0$; similarly for $D$.
Put
\[
H=h_Df+i_Dp_Dfh_C .
\]
Using $d_Di_D=0$ and $d_Df=fd_C$,
\[
d_DH+Hd_C=(d_Dh_D+h_Dd_D)f+i_Dp_Dfh_Cd_C
=(1-i_Dp_D)f+i_Dp_Df\,(1-i_Cp_C-d_Ch_C).
\]
Since $p_Dfd_C=p_Dd_Df=0$ and $p_Dfi_C=H(f)=0$, the last term equals
$i_Dp_Df$, and $d_DH+Hd_C=f$.
\end{proof}

The splitting maps depend on choices, so the homotopy $H$ is not
canonical, and no higher coherence is asserted. Over $\ZZ$ the lemma
fails, which is one reason Theorem~\ref{thm:A} is stated over $\FF$.

\section{Existence of DFL-regular primary data}\label{app:dfl-regular}

We prove Lemma~\ref{lem:dfl-detect} and Theorem~\ref{thm:dfl-regular}. The scheme of the
induction is that of \cite[proof of Proposition~6.8]{DFL} and
\cite[Section~5.5.1]{Donaldson-book}; the new input is Lemma~\ref{lem:dfl-detect}, whose proof
occupies Sections~\ref{subsec:app-setting}--\ref{subsec:app-proof-lemma}.

\subsection{Setting}\label{subsec:app-setting}

Fix a Heegaard splitting and metric data as in Definition~\ref{def:presentation}(i)--(iii),
and write $\Sigma=\Sigma_h\sqcup T^2$. Following \cite[Section~2.1, (2.1)]{DFL}, $M^\#$
contains a product region $(-2,2)\times\Sigma$ with the metric $dt^2+g_\Sigma$, and
\[
M^\#=C_0^\circ\cup N\cup C_1^\circ,\qquad N=[-1,1]\times\Sigma,\qquad
C_0^\circ\cap N=\{1\}\times\Sigma,\qquad C_1^\circ\cap N=\{-1\}\times\Sigma ;
\]
in the notation of \cite{DFL}, $C_0^\circ=Y_0$ and $C_1^\circ=Y_0'$. All gauge groups are
determinant-one \cite[Section~2.2]{DFL}, connections are of class $L^2_l$ with $l\ge2$, and
$L^2$ inner products use $\langle a,b\rangle=-\tr(ab)$. Gauge transformations of $C_i^\circ$ carry
no boundary condition. At an irreducible connection $B$ on $C_i^\circ$ the operator
$d_B\colon L^2_1(C_i^\circ)\to L^2(C_i^\circ)$ is injective and satisfies
$\|\zeta\|_{L^2_1}\le C(\|d_B\zeta\|+\|\zeta\|)$, so by the compactness of $L^2_1\subset L^2$ it
has closed range. The projection $\pi^i_B$ of Section~\ref{subsec:complexes} is
$\pi^i_Bb=b-d_B\zeta$, where $\zeta$ solves the Neumann problem $d_B^*d_B\zeta=d_B^*b$,
$*d_B\zeta=*b$ on $\partial C_i^\circ$, and it depends continuously on $B$ in operator norm near
irreducible connections.

Let $\mathfrak h$ be primary data and $g$ a split addition, and put $h=\mathfrak h+g=h_0+h_1$,
where $h_i=\mathfrak h_i+g_i$ depends only on the restriction to $C_i^\circ$. Put
\[
V(A)=*F_A+\nabla h(A),\qquad K_A=DV_A=*d_A+\operatorname{Hess}_Ah .
\]
For a cylinder function $h$, the gradient $\nabla h(B)$ lies in $\ker d_B^*$ and is supported in
the images of the immersions of $h$; the map $B\mapsto\nabla h(B)$ is smooth, with
$\|\nabla h(B)-\nabla h(B')\|_{L^p}\le C\|B-B'\|_{L^p}$, and $\operatorname{Hess}_Bh$ is bounded on
$L^2_k$ for $0\le k\le l$ and depends smoothly on $B$ \cite[Proposition~6.1]{DFL}; it is
symmetric \cite[Lemma~2.13]{DFL}.

The images of the immersions $\gamma_\alpha$ of the separated family and of the immersions
used by $g$ are compact subsets of
$\operatorname{int}C_0^\circ\sqcup\operatorname{int}C_1^\circ$. Fix $\eps\in(0,1)$ such that
$U=(-1-\eps,1+\eps)\times\Sigma$ is disjoint from all of them, and put
\[
V_0=(1,1+\eps)\times\Sigma\subset\operatorname{int}C_0^\circ,\qquad
V_1=(-1-\eps,-1)\times\Sigma\subset\operatorname{int}C_1^\circ .
\]
On $U$ we have $V(A)=*F_A$. As in Section~\ref{subsec:complexes}, $X_i\subset\cB(C_i^\circ)$ is the
set of restrictions of the critical points of $\mathrm{CS}_{\mathfrak h}$.

\begin{lemma}[locality]\label{lem:app-locality}
For $x\in C_i^\circ$ the value $V(A)(x)$ depends only on $A|_{C_i^\circ}$, and $(K_Ab)(x)$ only
on $b|_{C_i^\circ}$; write $V^{(i)}$ and $K^{(i)}$ for the resulting maps on $C_i^\circ$. For every
$\zeta\in\Omega^0(C_i^\circ;\ad)$, with no boundary condition,
\[
K^{(i)}_A(d_A\zeta)=[V^{(i)}(A),\zeta].
\]
The same identity holds on $M^\#$ with $V$ and $K$.
\end{lemma}

\begin{proof}
The gradient of $h_{1-i}$ vanishes on $C_i^\circ$, and $*F_A$ is local. The function $h_i$ is
invariant under all determinant-one gauge transformations of $C_i^\circ$, because its loops lie
in the interior and traces of holonomy are invariant under conjugation; so
$V^{(i)}(u\cdot A)=uV^{(i)}(A)u^{-1}$ for such $u$. Differentiating at $u=1$ in the direction
$\zeta$, with $u\cdot A=A-(d_Au)u^{-1}$, gives the identity. The argument on $M^\#$ is the same.
\end{proof}

\subsection{Trajectories and cokernels}

Trajectories of $\mathrm{CS}_h$ are written in temporal gauge as solutions
$\gamma=(A_s)_{s\in\RR}$ of $\partial_sA_s=-V(A_s)$ between determinant-one lifts of critical
points, in any path class. Reversing $s$ gives the opposite convention and changes the signs of
all eigenvalues below; the argument uses only that two of them have opposite signs. The
\emph{extended Hessian}
\[
L_A(b,c)=(K_Ab-d_Ac,\,-d_A^*b)
\]
acts on $H=L^2(M^\#;(\Lambda^1\oplus\Lambda^0)\otimes\ad)$ with domain $L^2_1$
(compare \cite[(2.23)]{Donaldson-book}) and is symmetric. Let $\tilde a$ be a lift of a
nondegenerate critical point; it is irreducible (Lemma~\ref{lem:app-irreducible}). Then $L_{\tilde a}$ is self-adjoint with
compact resolvent, and $\ker L_{\tilde a}=0$: if $K_{\tilde a}b=d_{\tilde a}c$ and
$d_{\tilde a}^*b=0$, then $\|d_{\tilde a}c\|^2=\langle b,K_{\tilde a}d_{\tilde a}c\rangle
=\langle b,[V(\tilde a),c]\rangle=0$ by Lemma~\ref{lem:app-locality}, so $c=0$ by
irreducibility and $b=0$ by nondegeneracy. The linearization at a trajectory $\gamma$, with the
Coulomb condition, is $D_\gamma=\partial_s+L_{A_s}$ (compare \cite[Section~3.1]{Donaldson-book}), and its
cokernel consists of the $L^2$ solutions of $\partial_sy=L_{A_s}y$, which are smooth and decay
exponentially at both ends.

\begin{lemma}\label{lem:app-trajectories}
Let $\gamma=(A_s)$ be a trajectory of $\mathrm{CS}_h$.
\begin{enumerate}[label=(\roman*)]
\item $x(s)=(\partial_sA_s,0)$ satisfies $\partial_sx=-L_{A_s}x$.
\item If $y=(\phi,c)$ is in the cokernel of $D_\gamma$, then $c\equiv0$; hence
  $d^*_{A_s}\phi_s=0$ and $\partial_s\phi=K_{A_s}\phi$.
\item Two solutions of the trajectory equation of $\mathrm{CS}_h$ on an interval that agree
  at one time agree on the whole interval; the proof is local in time. If $\gamma$
  is nonconstant, then $\partial_sA_s\ne0$ for every $s$, and if $y\ne0$ in (ii), then
  $\phi_s\ne0$ for every $s$.
\end{enumerate}
\end{lemma}

\begin{proof}
(i) Differentiating the flow gives $\partial_s\dot A=-K_A\dot A$, where $\dot A=\partial_sA$, and
$d_A^*\dot A=-d_A^*V(A)=0$ by the Bianchi identity and the gauge invariance of $h$.

(ii) The equation reads $\partial_s\phi=K_A\phi-d_Ac$, $\partial_sc=-d_A^*\phi$. Define $\{a,b\}$
by $\langle\zeta,\{a,b\}\rangle=\langle[a,\zeta],b\rangle$ for all sections $\zeta$. Then
$\partial_s(d^*_{A_s})b=\{\dot A_s,b\}$, and taking adjoints in Lemma~\ref{lem:app-locality} on
$M^\#$ gives $d_A^*K_Ab=\{V(A),b\}$. Hence
\[
\partial_s(d_A^*\phi)=\{\dot A+V(A),\phi\}-d_A^*d_Ac=-d_A^*d_Ac ,
\]
so $\partial_s^2c=d_A^*d_Ac$ and
$\tfrac{d^2}{ds^2}\|c\|^2=2\|\partial_sc\|^2+2\|d_Ac\|^2\ge0$. The function $\|c\|^2$ is convex
and tends to zero at both ends, so $c\equiv0$.

(iii) We use the theorem of Agmon and Nirenberg in the form \cite[Theorem~C.1]{SW-lagrangian},
in $L^2(M^\#)$, as in the proof of \cite[Proposition~8.6]{SW-lagrangian}; since $M^\#$ is closed
there are no boundary terms. For two trajectories $A$, $B$ put $x=A-B$. Then
\[
\partial_sx+*d_{B}x=-\tfrac12*[x\wedge x]-\bigl(\nabla h(A)-\nabla h(B)\bigr),
\]
whose $L^2$ norm is at most $c_1(s)\|x\|$ with $c_1$ locally bounded, by
\cite[Proposition~6.1(ii)]{DFL}. The operators $*d_{B(s)}$ are symmetric, and
$\tfrac{d}{ds}\langle x,*d_Bx\rangle-2\langle\partial_sx,*d_Bx\rangle=\langle
x,*[\partial_sB\wedge x]\rangle\le c_3(s)\|x\|^2$. So \cite[Theorem~C.1(i)]{SW-lagrangian}
gives $x\equiv0$ forward in time from a time at which $x$ vanishes, and backward after
reversing time, which replaces the operators by their negatives. The linear equations of (i)
and (ii) are treated in the same way, with the symmetric operators $\pm L_{A_s}$, for which
$c_1=0$ and $c_3=\|\partial_sL_{A_s}\|$ is locally bounded by \cite[Proposition~6.1(iii)]{DFL}.
\end{proof}

\subsection{Asymptotics at $-\infty$}

\begin{lemma}\label{lem:app-asymptotics}
Let $H$ be a real Hilbert space, $L$ a self-adjoint operator on $H$ with compact resolvent and
$\ker L=0$, and $B\colon(-\infty,0]\to\mathcal L(H)$ continuous with
$\beta(s):=\|B(s)\|\le C_0e^{\delta s}$ for some $\delta>0$. Let
$x\in C^1((-\infty,0];\operatorname{dom}L)$, with $\operatorname{dom}L$ carrying the graph norm,
solve $x'=-(L+B)x$, with $x(s)\ne0$ for all $s$, $x(s)\to0$ as $s\to-\infty$, and
\begin{equation}\label{eq:app-R}
\|Lx(s)\|\le C_1\sup_{|\sigma-s|\le1}\|x(\sigma)\|\qquad(s\le-1).
\end{equation}
Then there are $\lambda>0$ and $v\ne0$ with $Lv=-\lambda v$ such that $e^{-\lambda s}x(s)\to v$ as
$s\to-\infty$. If instead $x'=(L+B)x$, the same holds with $Lv=\lambda v$.
\end{lemma}

\begin{proof}
Put $\hat x=x/\|x\|$ and $N=\langle L\hat x,\hat x\rangle$.

(a) We have $(\log\|x\|)'=-N-\langle B\hat x,\hat x\rangle$ and
$\hat x'=-(L-N)\hat x-B\hat x+\langle B\hat x,\hat x\rangle\hat x$. Since
$\langle(L-N)\hat x,\hat x\rangle=0$,
\[
N'=2\langle(L-N)\hat x,\hat x'\rangle=-2\|(L-N)\hat x\|^2-2\langle(L-N)\hat x,B\hat x\rangle
\le-\|(L-N)\hat x\|^2+\beta^2 .
\]

(b) Hence $N-\int_{-\infty}^s\beta^2$ is nonincreasing, $N$ is bounded below on $(-\infty,0]$,
and $N_{-\infty}=\lim_{s\to-\infty}N(s)\in(-\infty,\infty]$ exists. If $N_{-\infty}>0$, then by
(a) $(\log\|x\|)'\le-c<0$ near $-\infty$, and $\|x(s)\|\to\infty$ as $s\to-\infty$, which is
excluded; so $N_{-\infty}\le0$. By (a), $\int_{-\infty}^0\|(L-N)\hat x\|^2<\infty$, so
$\|(L-N_{-\infty})\hat x(s_k)\|\to0$ along a sequence $s_k\to-\infty$. Hence
$N_{-\infty}\in\operatorname{spec}L$, and $N_{-\infty}=-\lambda$ with $\lambda>0$ since
$0\notin\operatorname{spec}L$.

(c) Let $P$ be the orthogonal projection onto the finite-dimensional space $\ker(L+\lambda)$,
$Q=1-P$, and $Q_\pm$ the spectral projections of $L$ onto the parts of the spectrum above and
below $-\lambda$; let $g>0$ be the distance from $-\lambda$ to the rest of the spectrum. Choose
$s_1$ with $|N+\lambda|\le g/4$ on $(-\infty,s_1]$. For $y=Q_+\hat x$ and $z=Q_-\hat x$, which
have norm at most one, (a) gives on $(-\infty,s_1]$
\[
(\|y\|^2)'\le-\tfrac32g\|y\|^2+4\beta\|y\|\le-g\|y\|^2+\tfrac8g\beta^2,\qquad
(\|z\|^2)'\ge\tfrac32g\|z\|^2-4\beta\|z\|\ge g\|z\|^2-\tfrac8g\beta^2 .
\]
Integrating the first inequality from $-\infty$, and the second backward from $s_1$, gives
$\|Q\hat x(s)\|\le Ce^{\delta's}$ for any fixed $\delta'<\min(\delta,g/2)$.

(d) The function $N$ is bounded on $(-\infty,0]$, so $(\log\|x\|)'$ is bounded, and
\eqref{eq:app-R} bounds $\|L\hat x\|$. Hence
$|N+\lambda|=|\langle(L+\lambda)Q\hat x,Q\hat x\rangle|\le Ce^{\delta's}$ is integrable, and by (a)
$e^{-\lambda s}\|x(s)\|$ converges to some $c\in(0,\infty)$.

(e) $(e^{-\lambda s}Px)'=-e^{-\lambda s}PBx$ has norm at most $C\beta$, so $e^{-\lambda s}Px$
converges to some $v$, while $e^{-\lambda s}Qx=e^{-\lambda s}\|x\|\,Q\hat x\to0$. So
$e^{-\lambda s}x\to v$, with $\|v\|=c>0$ and $Lv=-\lambda v$. The second statement follows by
applying the first to $-L$ and $-B$.
\end{proof}

\begin{corollary}\label{cor:app-leading}
Let $\gamma=(A_s)$ be a nonconstant trajectory of $\mathrm{CS}_h$, and $y=(\phi,0)\ne0$ in the
cokernel of $D_\gamma$. Let $\tilde a$ be the limit of $A_s$ as $s\to-\infty$. There are
$\lambda,\mu>0$ and nonzero $v,w\in\ker d^*_{\tilde a}$ with $K_{\tilde a}v=-\lambda v$ and
$K_{\tilde a}w=\mu w$ such that, in $L^2(M^\#)$,
\begin{equation}\label{eq:app-leading}
\partial_sA_s=e^{\lambda s}(v+o(1)),\qquad \phi_s=e^{\mu s}(w+o(1))\qquad(s\to-\infty).
\end{equation}
\end{corollary}

\begin{proof}
In a temporal gauge, $A_s$ converges to a lift $\tilde a$ of a critical point exponentially fast
in every $C^k$ norm \cite[Theorem~4.2 and Proposition~4.3]{Donaldson-book}, whose proofs apply
to the holonomy-perturbed equation with nondegenerate critical points
\cite[Section~5.5]{Donaldson-book}. Apply Lemma~\ref{lem:app-asymptotics} with
$H=L^2(M^\#;(\Lambda^1\oplus\Lambda^0)\otimes\ad)$, $L=L_{\tilde a}$ and
$B(s)=L_{A_s}-L_{\tilde a}$ to $x=(\partial_sA_s,0)$ and to $y$
(Lemma~\ref{lem:app-trajectories}(i),(ii)). The operator $B(s)$ consists of brackets with
$A_s-\tilde a$ and of $\operatorname{Hess}_{A_s}h-\operatorname{Hess}_{\tilde a}h$, so
$\|B(s)\|\le C\|A_s-\tilde a\|_{L^2_l}\le C_0e^{\delta s}$ by \cite[Proposition~6.1(iii)]{DFL}.
The solutions vanish nowhere by Lemma~\ref{lem:app-trajectories}(iii), and \eqref{eq:app-R} is
the interior elliptic estimate for $\partial_s+L_{A_s}$ on the slabs $[s-1,s+1]\times M^\#$, whose
coefficients are bounded uniformly in $s$, since $\operatorname{Hess}_{A_s}h$ is bounded on $L^2_k$
for $k\le2\le l$; the estimate at a fixed time uses the slab estimate one derivative higher,
for which $s\mapsto\operatorname{Hess}_{A_s}h$ is $C^1$ into the bounded operators on $L^2_1$
\cite[Proposition~6.1(iii)]{DFL}. The limits have vanishing $\Lambda^0$-components, because $x$ and $y$ have, so
they are eigenvectors of $K_{\tilde a}$ in $\ker d^*_{\tilde a}$.
\end{proof}

The two eigenvalues in \eqref{eq:app-leading}, $-\lambda<0<\mu$, have opposite signs; this is
what replaces the $L^2$-orthogonality of $\partial_s\gamma$ and $\psi$ used in
\cite[p.~79]{DFL}.

\subsection{Eigenvectors on one piece}

\begin{lemma}\label{lem:app-irreducible}
The restriction to $C_i^\circ$ of a lift $\tilde a$ of a critical point of $\mathrm{CS}_h$ is
irreducible. If $A_s\to\tilde a$ as $s\to-\infty$, then $A_s|_{C_i^\circ}$ is irreducible for all
$s$ near $-\infty$.
\end{lemma}

\begin{proof}
On $V_i\cap((-2,2)\times T^2)$, an interval times $T^2$, the connection $\tilde a$ satisfies
$*F_{\tilde a}=V(\tilde a)=0$. The bundle is odd on $T^2$, so the holonomy of this flat
connection lifts to the quaternion group in $\SU(2)$, which fixes no nonzero vector of
$\mathfrak{su}(2)$. A parallel section of $\ad$ over the connected $C_i^\circ$ therefore vanishes
on this open set, hence everywhere. Irreducibility is an open condition, and
$A_s\to\tilde a$ in $C^\infty$.
\end{proof}

\begin{lemma}\label{lem:app-eigen}
Let $a$ be a nondegenerate critical point of $\mathrm{CS}_h$, $i\in\{0,1\}$, and
$v_1,v_2\in\ker d_a^*$ eigenvectors of $K_a$ with eigenvalues $\nu_1\ne\nu_2$.
\begin{enumerate}[label=(\roman*)]
\item If $v_j|_{C_i^\circ}$ lies in the closure of $d_aL^2_1(C_i^\circ)$, then
  $v_j|_{C_i^\circ}=0$.
\item If $v_1|_{C_i^\circ}\ne0\ne v_2|_{C_i^\circ}$, then $\pi^i_av_1$ and $\pi^i_av_2$ are
  linearly independent.
\end{enumerate}
\end{lemma}

\begin{proof}
Since $V(a)=0$, Lemma~\ref{lem:app-locality} gives $K^{(i)}_a(d_a\zeta)=0$ for smooth $\zeta$, and
hence, by continuity of $K^{(i)}_a$ from $L^2$ to distributions on
$\operatorname{int}C_i^\circ$, for every $b$ in the closure of $d_aL^2_1(C_i^\circ)$. On
$C_i^\circ$ the eigenvector equation reads $K^{(i)}_av_j=\nu_jv_j$, and $\nu_j\ne0$ by
nondegeneracy; this gives (i). For (ii) suppose $\pi^i_av_2=c\,\pi^i_av_1$, so that
$v_2=cv_1+b$ on $C_i^\circ$ with $b$ in the closure of $d_aL^2_1$. If $c=0$ this contradicts
(i). Otherwise applying $K^{(i)}_a$ gives $\nu_2(cv_1+b)=c\nu_1v_1$, so
$c(\nu_2-\nu_1)v_1=-\nu_2b$ on $C_i^\circ$, and $v_1|_{C_i^\circ}$ lies in the closure of
$d_aL^2_1$, which contradicts (i).
\end{proof}

\subsection{Non-trapping}

The following Carleman estimate is proved by the method of \cite[Lemmas~5 and~6]{Booss-ucp},
after a change of variables that makes the inner product on the slices independent of the
normal variable; we give the proof in full because the slices below are tori and the operator
contains a nonlocal term acting slice by slice.

\begin{lemma}\label{lem:app-carleman}
Let $N_1$ be a compact region in a Riemannian three-manifold and $t\in C^\infty(N_1)$ with
$dt\ne0$ and $N_1=t^{-1}([\tau_0,\tau_1])$, whose level sets $S_t$ are closed surfaces. Let
$E\to N_1$ be a Hermitian bundle with a Clifford multiplication
$\mathbf c\colon T^*N_1\to\operatorname{End}(E)$, $\mathbf c(\xi)^*=-\mathbf c(\xi)$,
$\mathbf c(\xi)^2=-|\xi|^2$, let $\nabla$ be a Hermitian connection on $E$ and $\Theta$ a smooth
endomorphism, and put $D=\sum_j\mathbf c(e^j)\nabla_{e_j}+\Theta$ for local orthonormal frames.
Let $K$ act on $L^2(N_1;E)$ slice by slice: there are operators $K_t$ on $L^2(S_t;E)$, with the
area measure, such that $\|K_t\|\le b$ and $(Ku)|_{S_t}=K_t(u|_{S_t})$. Then there are
$T_0,C,R_0>0$, depending on $N_1$, $D$ and $b$, such that for every
$[t_0,t_1]\subset[\tau_0,\tau_1]$ with $T:=t_1-t_0\le T_0$:
\begin{enumerate}[label=(\roman*)]
\item for every $R\ge R_0$ and every $v\in L^2_1$ with compact support in $t^{-1}((t_0,t_1))$,
  \begin{equation}\label{eq:app-carleman}
  R\int e^{R(t_1-t)^2}|v|^2\le C\int e^{R(t_1-t)^2}|(D+K)v|^2 ;
  \end{equation}
\item if $u\in L^2_1(t^{-1}([t_0,t_1]);E)$ satisfies $(D+K)u=0$ and vanishes on
  $\{t<t_0+\eta\}$ for some $0<\eta<T$, then $u=0$ on $\{t\le t_0+T/2\}$.
\end{enumerate}
\end{lemma}

\begin{proof}
It suffices to prove \eqref{eq:app-carleman} for smooth $v$, since both sides are continuous in
$L^2_1$.

(a) \emph{Coordinates.} Let $X=\nabla t/|\nabla t|^2$. Its flow gives a diffeomorphism
$\Phi\colon[\tau_0,\tau_1]\times S\to N_1$, $S=S_{\tau_0}$, with $t(\Phi(t,y))=t$. Since $X$ is
orthogonal to the level sets, $\Phi^*g=|\nabla t|^{-2}dt^2+h_t$ with no cross terms. Write
$\Phi^*d\mathrm{vol}=\kappa\,dt\,d\mathrm{vol}_S$ with $\kappa>0$, for a fixed area form
$d\mathrm{vol}_S$ on $S$.

(b) \emph{A unitary change of variables.} Let $P_t$ be parallel transport for $\nabla$ along the
flow lines, and put $(\mathcal Uu)(t,y)=\kappa^{1/2}P_t^{-1}u(\Phi(t,y))$. Then $\mathcal U$ is a
unitary map from $L^2(N_1;E)$ to $L^2([\tau_0,\tau_1]\times S,\,dt\,d\mathrm{vol}_S;E|_S)$, and
the inner product on each slice is now $L^2(S,d\mathrm{vol}_S)$, independent of $t$.

(c) \emph{Normal form.} Let $e_0=\nabla t/|\nabla t|$, and let $e_1,e_2$ be tangent to the level
sets. Then $D=|\nabla t|\,\mathbf c(e^0)\nabla_X+D_S+\Theta$ with
$D_S=\sum_{k=1,2}\mathbf c(e^k)\nabla_{e_k}$, and $\mathcal U\nabla_X\mathcal U^{-1}=\partial_t
-\tfrac12\partial_t\log\kappa$. Hence
\[
\mathcal UD\,\mathcal U^{-1}=\sigma(\partial_t+\mathcal B_t),\qquad
\sigma=|\nabla t|\,P_t^{-1}\mathbf c(e^0)P_t,\qquad
\mathcal B_t=-\tfrac12\partial_t\log\kappa+\sigma^{-1}\,\mathcal U(D_S+\Theta)\,\mathcal U^{-1},
\]
where $\sigma$ is invertible, $\sigma^{-1}=-|\nabla t|^{-1}P_t^{-1}\mathbf c(e^0)P_t$, and
$\mathcal B_t$ is a first-order operator on $S$. For $\xi\in T^*_yS$ let $\eta$ be the covector at
$\Phi(t,y)$ with $\eta(X)=0$ whose pullback to $S$ is $\xi$; since $X$ is normal to the level
set, $\eta$ is orthogonal to $e^0$. The principal symbol of $\mathcal B_t$ at $\xi$ is
$-i|\nabla t|^{-1}P_t^{-1}\mathbf c(e^0)\mathbf c(\eta)P_t$. It is self-adjoint, because
$\mathbf c(\eta)\mathbf c(e^0)=-\mathbf c(e^0)\mathbf c(\eta)-2\langle e^0,\eta\rangle
=-\mathbf c(e^0)\mathbf c(\eta)$, and invertible for $\xi\ne0$, because
$(\mathbf c(e^0)\mathbf c(\eta))^2=-|\eta|^2$. Put $B_t=\tfrac12(\mathcal B_t+\mathcal B_t^*)$ and
$C_t=\tfrac12(\mathcal B_t-\mathcal B_t^*)$, adjoints in $L^2(S,d\mathrm{vol}_S)$. Then $B_t$ is
symmetric and elliptic of first order, $C_t$ is a skew-adjoint endomorphism, since the
principal symbols cancel, and both depend smoothly on $t$. Fix constants $c,c_1,s_\sigma$,
uniform over $[\tau_0,\tau_1]$, with $\|f\|_{H^1(S)}\le c(\|f\|+\|B_tf\|)$,
$|\langle f,(-\partial_tB_t+[B_t,C_t])f\rangle|\le c_1\|f\|\,\|f\|_{H^1(S)}$ (the operator is of
first order) and $\|\sigma^{-1}\|\le s_\sigma$.

(d) \emph{The estimate for $\partial_t+B_t+C_t$.} Let $w$ be smooth on $[t_0,t_1]\times S$ and
vanish near $t_0$ and $t_1$. Put $\mathcal S_t=B_t+R(t_1-t)$, which is symmetric, and
$J=\int_{t_0}^{t_1}\|\partial_tw+C_tw+\mathcal S_tw\|^2dt$. Since the slice inner product does not depend
on $t$,
\[
2\,\mathrm{Re}\int\langle\partial_tw,\mathcal S_tw\rangle\,dt=-\int\langle w,(\partial_t\mathcal S_t)w\rangle\,dt
=R\int\|w\|^2dt-\int\langle w,(\partial_tB_t)w\rangle\,dt ,
\]
and since $C_t$ is skew, $2\,\mathrm{Re}\langle C_tw,\mathcal S_tw\rangle=\langle w,[B_t,C_t]w\rangle$.
Expanding the square,
\[
J\ge\int\|\mathcal S_tw\|^2+R\int\|w\|^2+J_3,\qquad
J_3=\int\langle w,(-\partial_tB_t+[B_t,C_t])w\rangle .
\]
With $\|B_tw\|\le\|\mathcal S_tw\|+RT\|w\|$,
\[
|J_3|\le cc_1\int\|w\|\bigl(\|\mathcal S_tw\|+(1+RT)\|w\|\bigr)
\le\tfrac12\int\|\mathcal S_tw\|^2+\bigl(cc_1(1+RT)+\tfrac12(cc_1)^2\bigr)\int\|w\|^2 .
\]
For $T\le T_0:=1/(4cc_1)$ and $R\ge4cc_1+2(cc_1)^2$ this gives $J\ge\frac R2\int\|w\|^2$. For
smooth $v$ supported in $t^{-1}((t_0,t_1))$ put $w=e^{R(t_1-t)^2/2}\,\mathcal Uv$. Then
$e^{R(t_1-t)^2/2}(\partial_t+B_t+C_t)\mathcal Uv=(\partial_t+B_t+C_t+R(t_1-t))w$ and
$(\partial_t+B_t+C_t)\mathcal Uv=\sigma^{-1}\mathcal UDv$. Since $\mathcal U$ is unitary and the
weight depends only on $t$,
\[
R\int e^{R(t_1-t)^2}|v|^2\le 2s_\sigma^2\int e^{R(t_1-t)^2}|Dv|^2 .
\]

(e) \emph{Absorption of $K$.} By the coarea formula $d\mathrm{vol}=|\nabla t|^{-1}dA_t\,dt$, and
the weight is constant on each slice, so
$\int e^{R(t_1-t)^2}|Kv|^2\le b'^2\int e^{R(t_1-t)^2}|v|^2$ with
$b'^2=b^2\sup|\nabla t|/\inf|\nabla t|$. With $C_D=2s_\sigma^2$,
\[
R\int e^{R(t_1-t)^2}|v|^2\le2C_D\int e^{R(t_1-t)^2}|(D+K)v|^2+2C_Db'^2\int e^{R(t_1-t)^2}|v|^2 ,
\]
and for $R\ge R_0:=\max\bigl(4cc_1+2(cc_1)^2,\,4C_Db'^2\bigr)$ this is \eqref{eq:app-carleman}
with $C=4C_D$.

(f) \emph{Continuation.} Let $\varphi(t)=1$ for $t\le t_0+0.8T$ and $\varphi(t)=0$ for
$t\ge t_0+0.9T$, and put $v=\varphi u$. Then $v\in L^2_1$ has compact support in
$t^{-1}((t_0,t_1))$, because $u=0$ on $\{t<t_0+\eta\}$. Since $K$ commutes with multiplication by
functions of $t$, $(D+K)v=\varphi'(t)\,\mathbf c(dt)u$, which is supported where
$(t_1-t)^2\le T^2/25$, while $(t_1-t)^2\ge T^2/4$ on $\{t\le t_0+T/2\}$, where $v=u$. By
\eqref{eq:app-carleman},
\[
\int_{\{t\le t_0+T/2\}}|u|^2\le\frac CR\,e^{-\frac{21}{100}RT^2}\int|\varphi'\,\mathbf c(dt)u|^2 ,
\]
which tends to zero as $R\to\infty$.
\end{proof}

\begin{proposition}\label{prop:nontrapping}
Let $\mathfrak h$ be primary data in the span of a separated family with embeddings
$\gamma_\alpha$, and choose $r<1$ such that every two-form $\mu$ in the averaged trace functions
of $\mathfrak h$ is supported in the disk $D_r$ of radius $r$. Put
$\Omega'=M^\#\setminus\bigcup_\alpha\gamma_\alpha(S^1\times\bar D_r)$. Let $a$ be a critical point
of $\mathrm{CS}_{\mathfrak h}$ and $v\in\ker d_a^*$ an eigenvector of $K_a$. If $v$ vanishes on a
nonempty open subset of $\Omega'$, then $v=0$.
\end{proposition}

\begin{proof}
A critical point of $\mathrm{CS}_{\mathfrak h}$ is gauge equivalent to a smooth connection, by
elliptic bootstrapping with $\nabla\mathfrak h\colon L^2_k\to L^2_k$ for every $k$
\cite[Proposition~6.1]{DFL}; we work in such a gauge. Let $\nu$ be the eigenvalue.
\emph{Outside the solid tori.} On $\Omega'$ the Hessian of
$\mathfrak h$ vanishes, so $a$ is flat there and $*d_av=\nu v$, $d_a^*v=0$. Applying $*d_a$ gives
$d_a^*d_av=\nu^2v$, so $\Delta_av=\nu^2v$ on $\Omega'$, a second-order system whose principal part
is the scalar Laplacian in every component. The set $\Omega'$ is connected: it is the complement
in the connected $M^\#$ of finitely many disjoint closed solid tori, each of which has a connected
collar $\gamma_\alpha(S^1\times\{r<|z|<1\})$ in $\Omega'$, so a decomposition of $\Omega'$ into two
open sets would give one of $M^\#$. By elliptic regularity $v$ is smooth on $\Omega'$. By the
unique continuation theorem of \cite{Aronszajn}, which applies to such systems
\cite[Remark~3]{Aronszajn} (see also \cite[Remark~8(a)]{Booss-ucp}), applied on an exhaustion of
$\Omega'$ by connected relatively compact open subsets, $v=0$ on $\Omega'$.

\emph{Inside a solid torus.} Fix $\alpha$ and use the coordinates $(\theta,z)\in S^1\times D^2$ of
$\gamma_\alpha$. Since the images are disjoint and $\mathfrak h$ is a linear combination of averaged
trace functions, on $\gamma_\alpha(S^1\times D^2)$ the operator $\operatorname{Hess}_a\mathfrak h$ is a
linear combination of Hessians of averaged trace functions of $\gamma_\alpha$. The gradient of such
a function at $\gamma_\alpha(\theta,z)$ is determined by $\mu(z)$ and by the traceless part of the
holonomy around the loop $\gamma_\alpha(S^1\times\{z\})$ based at that point
\cite[Lemma~5.12]{Donaldson-book}; its derivative in a direction $b$ is therefore an integral
along that loop of $b$, transported by $a$, with bounded kernel. So the operator
$\widetilde B(c,b)=(0,\operatorname{Hess}_a\mathfrak h\,(b))$ maps sections over a union of loops
$\gamma_\alpha(S^1\times\{z\})$ to sections over the same union, with a uniform bound on $L^2$ of
such unions. With $\mathcal D_a(c,b)=(d_a^*b,\,d_ac+*d_ab)$ on
$E=(\Lambda^0\oplus\Lambda^1)\otimes\ad\otimes\CC$, the section $u=(0,v)$ satisfies
\[
(\mathcal D_a-\nu)u+\widetilde Bu=0 .
\]
The operator $\mathcal D_a-\nu$ has the form of Lemma~\ref{lem:app-carleman}, with
$\mathbf c(\xi)(c,b)=(-\iota_{\xi^\sharp}b,\,\xi c+*(\xi\wedge b))$, which is skew and satisfies
$\mathbf c(\xi)^2=-|\xi|^2$, with $\nabla$ the Levi-Civita connection coupled to $a$, and
$\Theta=-\nu$.

Let $\omega'$ be the set of $z\in D^\circ$ such that $u$ vanishes on $\gamma_\alpha(S^1\times W)$
for a neighbourhood $W$ of $z$. It is open, and it contains $\{r<|z|<1\}$ by the first step.
Suppose $\omega'\ne D^\circ$. Its complement in $D^\circ$ is compact, being closed in $D^\circ$ and
contained in $\bar D_r$, and $D^\circ$ is connected; so there is $z_1\in\omega'$ whose distance
$\rho$ to $D^\circ\setminus\omega'$ is less than $(1-r)/2$, hence less than
$\dist(z_1,\partial D)$. Let $z_0\in D^\circ\setminus\omega'$ attain this distance. Then
$D(z_1,\rho)\subset\omega'$. Fix $0<\eta_1<\rho$ and $T_1>0$ with
$\rho+T_1<\dist(z_1,\partial D)$, put $t=|z-z_1|-\rho$, and let
$N_1=\gamma_\alpha(S^1\times\{-\eta_1\le t\le T_1\})$. On $N_1$ we have $dt\ne0$, and the level
sets $S_t=\gamma_\alpha(S^1\times\{|z-z_1|=\rho+t\})$ are closed tori, each a union of loops
$\gamma_\alpha(S^1\times\{z\})$, so $K=\widetilde B$ acts slice by slice with a uniform bound.
Let $T_0$ be given by Lemma~\ref{lem:app-carleman} for $N_1$, and choose $0<\eta<T'$ with
$T'\le T_1$, $\eta\le\eta_1$ and $T'+\eta\le T_0$. Lemma~\ref{lem:app-carleman}(ii) with
$[t_0,t_1]=[-\eta,T']$, of length $T'+\eta$, where $u=0$ for $t<0=t_0+\eta$, gives $u=0$ on
$\{t\le(T'-\eta)/2\}$. Together with $u=0$ over $D(z_1,\rho)$, this is a neighbourhood of
$\gamma_\alpha(S^1\times\{z_0\})$, so $z_0\in\omega'$, a contradiction. Hence $v$ vanishes on every
solid torus, and $v=0$.
\end{proof}

\begin{corollary}\label{cor:app-nontrapping}
Let $\mathfrak h$ and $g$ be as in Lemma~\ref{lem:dfl-detect}, with $\eps$ as in
Section~\ref{subsec:app-setting}. For every critical point $a$ of $\mathrm{CS}_{\mathfrak h+g}$,
every eigenvector $v\ne0$ of $K_a$ in $\ker d_a^*$ and each $i$, we have $v|_{V_i}\ne0$, and
hence $\pi^i_av\ne0$.
\end{corollary}

\begin{proof}
The critical points of $\mathrm{CS}_{\mathfrak h+g}$ are those of $\mathrm{CS}_{\mathfrak h}$, and
$g$ vanishes on a neighbourhood of each of them, so $K_a$ is the Hessian of
$\mathrm{CS}_{\mathfrak h}$ at $a$. The open set $V_i$ lies in $\Omega'$ by the choice of $\eps$,
so Proposition~\ref{prop:nontrapping} gives $v|_{V_i}\ne0$, and
Lemma~\ref{lem:app-eigen}(i) gives $\pi^i_av\ne0$.
\end{proof}

\subsection{Restricted paths}

For a trajectory $\gamma=(A_s)$ of $\mathrm{CS}_h$ put
$\Gamma_i(s)=[A_s|_{C_i^\circ}]\in\cB(C_i^\circ)$, extended continuously to $\pm\infty$ by the
restrictions of the limits.

\begin{lemma}\label{lem:app-hybrid}
Suppose that $\mathrm{CS}_h$ and $\mathrm{CS}_{\mathfrak h}$ have the same critical points, and
let $\gamma$ be a trajectory of $\mathrm{CS}_h$ and $i\in\{0,1\}$. Either $\Gamma_i$ is constant,
or it is injective on $\RR$ and $\Gamma_i(\RR)\cap X_i=\varnothing$.
\end{lemma}

\begin{proof}
Take $i=0$; the case $i=1$ is symmetric. Put $O_1=\RR\times\operatorname{int}C_0^\circ$ and
$O_2=\RR\times(U\cup C_1^\circ)$. Then $O_1\cup O_2=\RR\times M^\#$ and
$O_1\cap O_2=\RR\times V_0$.

\emph{Returns.} Suppose $\Gamma_0(s_1)=\Gamma_0(s_0)$ with $\tau=s_1-s_0\ne0$. Choose a gauge
transformation $k$ of $C_0^\circ$ with $k\cdot A_{s_1}=A_{s_0}$ on $C_0^\circ$; it is smooth. It
extends to a determinant-one gauge transformation of $M^\#$, since the obstructions lie in
$H^j(M^\#,C_0^\circ;\pi_{j-1}\SU(2))=0$ for $j\le3$. Put $\tilde\gamma(s)=k\cdot\gamma(s+\tau)$, a
trajectory of $\mathrm{CS}_h$ in temporal gauge with $\tilde\gamma(s_0)=\gamma(s_0)$ on $C_0^\circ$.
On the open set $V_0$ both $\tilde\gamma$ and $\gamma$ solve the unperturbed equation
$\partial_sA+*F_A=0$, so they agree on $\RR\times V_0$ by \cite[Lemma~8.7(i)]{SW-lagrangian}.
Define $\Xi$ to be $\tilde\gamma$ on $O_1$ and $\gamma$ on $O_2$. It is smooth, and it is a
trajectory of $\mathrm{CS}_h$: on $O_1$ the equation involves only the restriction to $C_0^\circ$
(Lemma~\ref{lem:app-locality}), where $\Xi=\tilde\gamma$; on $O_2$ it involves the unperturbed
region $U$ and the restriction to $C_1^\circ$, where $\Xi=\gamma$. Since $\Xi(s_0)=\gamma(s_0)$ on
$M^\#$, Lemma~\ref{lem:app-trajectories}(iii) gives $\Xi=\gamma$, so $\tilde\gamma=\gamma$ on
$\RR\times C_0^\circ$. Hence $\Gamma_0$ is periodic with period $\tau$; it converges as
$s\to-\infty$, so it is constant.

\emph{Critical passages.} Suppose $\Gamma_0(s_1)$ is the restriction of a critical point, with
lift $\tilde c$. The same argument, with $\tilde\gamma$ the constant trajectory at $k\cdot\tilde c$
for a gauge transformation $k$ with $k\cdot\tilde c=A_{s_1}$ on $C_0^\circ$, applies: $\tilde c$ is
flat on $V_0$ and $V(\tilde c)=0$ on $M^\#$. It gives $\Gamma_0\equiv\Gamma_0(s_1)$.
\end{proof}

\subsection{Holonomy detection and localization}

\begin{lemma}\label{lem:app-detection}
Let $B$ be a smooth irreducible connection on $C_i^\circ$. All trace functions below are averaged
trace functions of embeddings into $\operatorname{int}C_i^\circ$.
\begin{enumerate}[label=(\roman*)]
\item If a smooth one-form $b$ on $C_i^\circ$ satisfies $D\tau_B(b)=0$ for every trace function
  $\tau$, then $b\in d_B\Omega^0(C_i^\circ)$.
\item If $S$ is a finite-dimensional space of smooth one-forms on $C_i^\circ$ with
  $S\cap d_B\Omega^0(C_i^\circ)=0$, there are finitely many trace functions whose derivatives at
  $B$ are jointly injective on $S$.
\item If $\mathcal K\subset\cB(C_i^\circ)$ is compact and $[B]\notin\mathcal K$, finitely many
  trace functions separate $[B]$ from $\mathcal K$.
\end{enumerate}
\end{lemma}

\begin{proof}
(i) Letting $\mu$ concentrate at a point shows that the derivative in the direction $b$ of the
trace of the holonomy around every embedded loop in $\operatorname{int}C_i^\circ$ vanishes, and by
approximation \cite[p.~141]{Donaldson-book} the same holds for every piecewise smooth loop. Fix
$p\in\operatorname{int}C_i^\circ$. Since $B$ is irreducible, there are loops $\delta_1,\dots,\delta_m$
based at $p$ whose holonomies do not commute. Fix a unitary frame of the fibre at $p$, and write
the holonomies of $\delta_1,\dots,\delta_m$ and of a further loop $\delta$ based at $p$ as
$\lambda_jg_j$ and $\lambda g$, with $g_j,g\in\SU(2)$ and $\lambda_j,\lambda\in U(1)$ independent of
$B$, as in Section~\ref{subsec:compressions}. Identify $\SU(2)$ with the unit quaternions, let
$\rho$ be the homomorphism from the free group on $m+1$ letters to $\SU(2)$ sending the letters
to $g_1,\dots,g_m,g$, and extend it to an algebra map $\rho\colon\RR[F_{m+1}]\to\mathbb H$. For a
word $W$ the holonomy around the corresponding loop is $\lambda_W\rho(W)$ with $\lambda_W\in U(1)$
independent of $B$, so the derivative of $\tr\rho(W)$ in the direction $b$ vanishes. Let
$\dot\rho(W)\in\mathbb H$ be the derivative of $\rho(W)$ in the direction $b$, extended linearly.
Then $\dot\rho(UV)=\dot\rho(U)\rho(V)+\rho(U)\dot\rho(V)$ and $\tr\dot\rho=0$ on
$\RR[F_{m+1}]$. The image of $\rho$ is the real span of a noncommutative subgroup of the unit
quaternions; the proper subalgebras of $\mathbb H$ are $\RR$ and copies of $\CC$, which are
commutative, so the image is $\mathbb H$. If $\rho(x)=0$, then
$\tr(\dot\rho(x)\rho(y))=\tr\dot\rho(xy)-\tr(\rho(x)\dot\rho(y))=0$ for all $y$, and since the
trace form of $\mathbb H$ is nondegenerate, $\dot\rho(x)=0$. Hence $\dot\rho=D\circ\rho$ for a
derivation $D$ of $\mathbb H$. Every derivation of $\mathbb H$ is inner, so
$D=[\eta,\cdot\,]$ with $\eta\in\operatorname{Im}\mathbb H=\mathfrak{su}(2)$, and $\eta$ is
unique because $g_1,\dots,g_m$ do not commute. So the derivative of the holonomy in the direction
$b$ is the infinitesimal conjugation by $\eta$ for every loop based at $p$, and one $\eta$ serves
for all of them. Let $\zeta_0$ be a smooth section whose value at $p$ induces this infinitesimal
conjugation. Then the derivative of the holonomy in the direction $b'=b-d_B\zeta_0$ vanishes for
every loop based at $p$, and $\zeta_1(q)=\int_\ell b'$, computed in a parallel frame along a path
$\ell$ from $p$ to $q$, is independent of $\ell$ and satisfies $d_B\zeta_1=b'$
\cite[(6.4)]{DFL}. So $b=d_B(\zeta_0+\zeta_1)$.

(ii) By (i), no nonzero element of $S$ is annihilated by the derivatives of all trace functions,
so these derivatives span the dual of $S$, and $\dim S$ of them suffice.

(iii) A connection is determined up to gauge equivalence by its holonomy representation up to
conjugacy, and this representation by the traces of all words \cite[Lemma~5.13 and
p.~141]{Donaldson-book}. So for each $[B']\in\mathcal K$ some loop, which may be taken embedded,
has trace of holonomy at $B'$ different from that at $B$, and so has an averaged trace function
with $\mu$ concentrated near the centre. This is an open condition in $[B']$, and $\mathcal K$ is
compact.
\end{proof}

\begin{lemma}\label{lem:app-localize}
Let $\nu\colon\RR\to\RR^n$ be $C^2$ and $W\colon\RR\to\RR^n$ be $C^1$ and integrable. Let
$s_0\in\RR$ and $\delta>0$ be such that $\nu$ is an injective immersion on $(s_0-\delta,s_0+\delta)$,
$\nu(s_0)$ is not in the closure of $\nu(\RR\setminus(s_0-\delta,s_0+\delta))$, and
$W(s_0)\notin\RR\,\nu'(s_0)$. Then every neighbourhood of $\nu(s_0)$ contains the support of a
function $G\in C^\infty_c(\RR^n)$ with $\int_\RR\langle\nabla G(\nu(s)),W(s)\rangle\,ds\ne0$.
\end{lemma}

\begin{proof}
Write $\nu_0=\nu(s_0)$. Near $s_0$ the distance $|\nu(s)-\nu_0|$ is strictly monotone on each
side of $s_0$, because $\nu'(s_0)\ne0$, and the rest of $\nu(\RR)$ stays at positive distance
from $\nu_0$, by injectivity on $(s_0-\delta,s_0+\delta)$ and the closure hypothesis. So for small
$r$ the preimage of the ball of radius $r$ about $\nu_0$ is an interval $J\ni s_0$ of length
$O(r)$. Choose $n\perp\nu'(s_0)$ with $\langle n,W(s_0)\rangle=1$, a function
$\chi\in C^\infty_c([0,1))$ with $\chi\ge0$ and $\chi=1$ near $0$, and put
$G(x)=\chi(|x-\nu_0|/r)\,\langle n,x-\nu_0\rangle$. Along $\nu$ we have
$\langle n,\nu(s)-\nu_0\rangle=O((s-s_0)^2)$, so the term of
$\int\langle\nabla G(\nu),W\rangle$ containing the derivative of $\chi$ is $O(r^2)$, while the
term $\int\chi(|\nu-\nu_0|/r)\langle n,W\rangle$ equals $c_\chi r/|\nu'(s_0)|+O(r^2)$ with
$c_\chi=\int_\RR\chi(|x|)\,dx>0$. For small $r$ the integral is nonzero.
\end{proof}

\begin{lemma}\label{lem:app-reduction}
Let $\mathfrak h$, $g$, $\gamma=(A_s)$ and $\psi=(\phi,0)$ be as in Lemma~\ref{lem:dfl-detect},
$i\in\{0,1\}$ and $s_0\in\RR$. Suppose that $B=A_{s_0}|_{C_i^\circ}$ is irreducible, that
$\pi^i_B\partial_sA_{s_0}$ and $\pi^i_B\phi_{s_0}$ are linearly independent, and that
$\Gamma_i(s_0)\notin\Gamma_i([-\infty,\infty]\setminus\{s_0\})\cup X_i$. Then the conclusion (iii)
of Lemma~\ref{lem:dfl-detect} holds.
\end{lemma}

\begin{proof}
Let $S$ be the span of $\pi^i_B\partial_sA_{s_0}$ and $\pi^i_B\phi_{s_0}$. These are smooth, by
elliptic regularity for the Neumann problem that defines $\pi^i_B$, and orthogonal to
$d_BL^2_1(C_i^\circ)$, so $S\cap d_B\Omega^0=0$. Trace functions are invariant under
all gauge transformations of $C_i^\circ$, so their derivatives at $B$ vanish on $d_BL^2_1$. By
Lemma~\ref{lem:app-detection}(ii) there are trace functions $\tau'$ with $D\tau'_B$ injective on
$S$; then $D\tau'_B(\partial_sA_{s_0})=D\tau'_B(\pi^i_B\partial_sA_{s_0})\ne0$, so $\tau'\circ\Gamma_i$
is an injective immersion on some interval $(s_0-\delta,s_0+\delta)$. The set
$\mathcal K=\Gamma_i([-\infty,s_0-\delta]\cup[s_0+\delta,\infty])\cup X_i$ is compact and does not
contain $\Gamma_i(s_0)$, so Lemma~\ref{lem:app-detection}(iii) gives trace functions $\tau''$
separating $\Gamma_i(s_0)$ from $\mathcal K$. Put $\tau=(\tau',\tau'')$, $\nu=\tau\circ\Gamma_i$ and
$W(s)=D\tau_{A_s}(\phi_s)$, which is integrable because $\phi$ decays exponentially and the
gradients of trace functions are bounded. The hypotheses of Lemma~\ref{lem:app-localize} hold:
$\nu(\RR\setminus(s_0-\delta,s_0+\delta))\subset\tau(\mathcal K)$, a compact set not containing
$\nu(s_0)$, and $\nu'(s_0)=D\tau_B(\pi^i_B\partial_sA_{s_0})$ and
$W(s_0)=D\tau_B(\pi^i_B\phi_{s_0})$ are linearly independent. It gives $G$ supported in a ball about
$\nu(s_0)$ whose closure misses $\tau(\mathcal K)\supset\tau(X_i)$. The cylinder function
$g_i'=G\circ\tau$ has loops in $\operatorname{int}C_i^\circ$ and vanishes on the neighbourhood
$\tau^{-1}(\RR^n\setminus\operatorname{supp}G)$ of $X_i$, and by the chain rule
$\int\langle\nabla g_i'(A_s),\phi_s\rangle\,ds=\int\langle\nabla G(\nu(s)),W(s)\rangle\,ds\ne0$.
\end{proof}

\subsection{Proof of Lemma~\ref{lem:dfl-detect}}\label{subsec:app-proof-lemma}

By Lemma~\ref{lem:app-trajectories}(ii) the cokernel element is $\psi=(\phi,0)$. Fix $i$ and let
$\tilde a$ be the limit of $A_s$ as $s\to-\infty$. By Corollary~\ref{cor:app-leading},
\eqref{eq:app-leading} holds with eigenvalues $-\lambda<0<\mu$ of $K_{\tilde a}$. By
Corollary~\ref{cor:app-nontrapping}, $\pi^i_{\tilde a}v\ne0\ne\pi^i_{\tilde a}w$, and since
$-\lambda\ne\mu$, Lemma~\ref{lem:app-eigen}(ii) shows that $\pi^i_{\tilde a}v$ and
$\pi^i_{\tilde a}w$ are linearly independent. The restriction of $\tilde a$ to $C_i^\circ$ is
irreducible (Lemma~\ref{lem:app-irreducible}), so $\pi^i_{A_s}\to\pi^i_{\tilde a}$ in operator norm,
and by \eqref{eq:app-leading}
\[
e^{-\lambda s}\pi^i_{A_s}\partial_sA_s\to\pi^i_{\tilde a}v,\qquad
e^{-\mu s}\pi^i_{A_s}\phi_s\to\pi^i_{\tilde a}w .
\]
Linear independence of a pair of vectors is an open condition, so there is $T$ such that (ii)
holds for all $s_0<-T$. In particular $\pi^i\partial_sA_{s_0}\ne0$, so $\Gamma_i$ is not constant,
and Lemma~\ref{lem:app-hybrid} gives (i). The limits $\Gamma_i(\pm\infty)$ lie in $X_i$, so by (i)
$\Gamma_i(s_0)\notin\Gamma_i([-\infty,\infty]\setminus\{s_0\})\cup X_i$ for every $s_0$, and
Lemma~\ref{lem:app-reduction} at any $s_0<-T$ gives (iii). \qed

\subsection{Proof of Theorem~\ref{thm:dfl-regular}}

\begin{lemma}\label{lem:app-separated}
\begin{enumerate}[label=(\roman*)]
\item For every Heegaard splitting and metric data there are finite separated families of
  primary perturbations containing the height functions of Section~\ref{sec:torus} with the
  submersion property of Section~\ref{subsec:compressions}: the
  universal restriction map of their span is a submersion at every flat connection of $C_0$ and
  of $C_1$. The family $F_*$ fixed there may be taken to be one of them.
\item For every finite family $F$ containing such a family, in particular for every
  $F\supseteq F_*$, generic primary data in a
  neighbourhood of zero in $\operatorname{span}F$ have transverse Lagrangians.
\end{enumerate}
\end{lemma}

\begin{proof}
(i) Fix $i$. For a flat connection $B$ on $C_i$ let
$\mathcal H_B=\{b:\ d_Bb=0,\ d_B^*b=0,\ *b|_{\partial C_i}=0\}$. It maps isomorphically onto
$H^1(C_i;\ad B)$. By Definition~\ref{def:clean-datum}(C2), $B$ is irreducible and
$H^2(C_i;\ad B)=0$, so $\dim\mathcal H_B=-3\chi(C_i)$ is independent of $B$, and $\mathcal H_B$
depends continuously on $B$ in a local slice of the flat connections. Let $0\ne b\in\mathcal H_B$.
If $b|_{C_i^\circ}\in d_B\Omega^0(C_i^\circ)$, then the class of $b$ vanishes, because
$C_i^\circ\hookrightarrow C_i$ is a homotopy equivalence; so $b=d_B\zeta$ on $C_i$ and
$\|b\|^2=\langle\zeta,d_B^*b\rangle+\int_{\partial C_i}\langle\zeta,*b\rangle=0$, a contradiction.
Hence the restrictions to $C_i^\circ$ of the elements of $\mathcal H_B$ form a space $S$ with
$S\cap d_B\Omega^0(C_i^\circ)=0$, and Lemma~\ref{lem:app-detection}(ii) gives finitely many averaged
trace functions of embeddings into $\operatorname{int}C_i^\circ$ whose derivatives at $B$ are
jointly injective on $\mathcal H_B$. At a flat connection $B'$ and for $d_{B'}b=0$, the derivative in the direction $b$ of
$\tr g_\ell$ is unchanged by a homotopy of $\ell$, with $\lambda$ transported along it: the first
variation of the holonomy between homotopic loops differs by the integral of $d_{B'}b$ over the swept
surface. So on $\mathcal H_{B'}$ the derivative of an averaged trace function equals that of
$\tr g_\ell$ for its core $\ell$, and up to sign it depends only on the free homotopy class of
$\ell$. Injectivity on $\mathcal H_{B'}$ is open in $B'$, and the space of gauge equivalence classes
of flat connections on $C_i$ is compact, so the cores $\ell_1,\dots,\ell_n$ of finitely many of the
functions above serve every flat connection on $C_i$. Let $W$ be the complement in
$\operatorname{int}C_i^\circ$ of the solid tori of the height functions. It is connected, and every
loop in $C_i^\circ$ is freely homotopic into $W$: push it off the cores of those solid tori by
general position, then radially out of each of them. By general position there are pairwise
disjoint embedded loops $\ell_k'\subset W$ with $\ell_k'$ freely homotopic to $\ell_k$. Take pairwise
disjoint thin tubular neighbourhoods of them in $W$, and any $\mu$. At every flat $B'$ the
derivatives on $\mathcal H_{B'}$ of the resulting averaged trace functions agree up to sign with those
of the $\tr g_{\ell_k}$, so they are jointly injective. Let $F$ consist of the height functions and
these averaged trace functions; it is separated.

By the proofs of \cite[Propositions~6.5 and~6.7]{DFL}, with the Banach space of perturbations
replaced by $\operatorname{span}F$, the universal restriction map is a submersion at a flat
connection $B$ on $C_i$ provided two cokernels vanish. The first,
$H^1(C_i,\partial C_i;\ad B)\cong H^2(C_i;\ad B)^*$, is zero by
Definition~\ref{def:clean-datum}(C2). The second consists of the $b\in\mathcal H_B$ annihilated by
$D\tau_B$ for every $\tau\in F$, and it is zero by construction. So $\operatorname{span}F$ has the submersion property of
Section~\ref{subsec:compressions}.

(ii) The submersion property holds for $\operatorname{span}F$ at every flat connection, since it
holds for a subspace, and is an open condition, so it holds over a neighbourhood of zero in
$\operatorname{span}F$, and
transversality of the two Lagrangians for generic primary data there follows as in
\cite[proof of Lemma~2.56, p.~78]{DFL}.
\end{proof}

\begin{lemma}\label{lem:app-critical}
Let $\mathfrak h$ have transverse Lagrangians. For every split addition $g$ sufficiently small
in the norm
$\sup_B\|\nabla g(B)\|_{L^\infty}+\sup_{B\ne B'}\|\nabla g(B)-\nabla g(B')\|_{L^2}/\|B-B'\|_{L^2}$,
which is finite by \cite[Proposition~6.1]{DFL},
$\mathrm{CS}_{\mathfrak h+g}$ and $\mathrm{CS}_{\mathfrak h}$ have the same critical points, with
the same Hessians.
\end{lemma}

\begin{proof}
Critical points correspond to intersection points of the Lagrangians
\cite[Lemmas~2.45 and~2.47]{DFL}. For $g$ small the perturbed flat connections on $C_i$ are
$C^1$-close to those for $\mathfrak h$, by the implicit function theorem. Those whose restrictions
to $\Sigma_h$ are near a point of $L_{\mathfrak h_0}(C_0)\cap L_{\mathfrak h_1}(C_1)$ are near the
corresponding perturbed flat connection, where $g_i$ vanishes; so each intersection point
remains an intersection point with unchanged tangent spaces, and transversality together with
$C^1$-closeness of the Lagrangians, uniform in that norm, gives uniqueness of the intersection
near it. Away from a neighbourhood of the intersection the two Lagrangians
for $\mathfrak h$ are at positive distance, and so are those for $\mathfrak h+g$ when $g$ is
small. The Hessians agree because $g$ vanishes near the critical points.
\end{proof}

\begin{proof}[Proof of Theorem~\ref{thm:dfl-regular}]
The first statement is Lemma~\ref{lem:app-separated}(i) and~(ii). Let $\mathfrak h\in\mathcal U(F)$ have
transverse Lagrangians; its critical points are nondegenerate \cite[Lemmas~2.45
and~2.47]{DFL}. The functional $\mathrm{CS}_{\mathfrak h+g}$ is invariant under $\alpha_M$, which is dual to the $T^2$
component of $\Sigma$; that torus misses $\operatorname{int}C_i^\circ$, so every trace function with loops in
$\operatorname{int}C_i^\circ$ is $\alpha_M$-invariant; the gauge complex modulo $\phi_M$ is
therefore defined. We follow \cite[proof of Proposition~6.8]{DFL}, by induction over the energy
as in \cite[p.~145]{Donaldson-book}. Nonconstant trajectories have positive energy, and the positive energies of path classes of
index at most seven in the determinant-one theory form a finite set $e_1<\dots<e_N$,
independent of $g$ since $g$ vanishes at the critical points: there are finitely many critical
points, and within each class the index determines the energy. Suppose a split addition $g_{k-1}$ has been
found, as small as Lemma~\ref{lem:app-critical} requires, such that all trajectories of
$\mathrm{CS}_{\mathfrak h+g_{k-1}}$ of index at most seven and energy less than $e_k$ are regular.
The non-regular trajectories of index at most seven and energy $e_k$ form a compact set. The
components of a broken limit have smaller energy. None has index at least eight, since the
others would then have total index at most $-1$ and one of them negative index; so each has
index at most seven, hence is regular by the induction hypothesis, and has index at least one.
Bubbling is excluded by the energy and index bounds \cite[p.~79]{DFL}, and trajectories near a
broken trajectory with regular components are themselves regular, by gluing
\cite[p.~144]{Donaldson-book}; so a limit of non-regular trajectories is unbroken. So does the set of pairs $(\gamma,\psi)$ with $\gamma$ such a trajectory, up
to translation, and $\psi$ a unit cokernel element \cite[p.~144]{Donaldson-book}. For each such pair
Lemma~\ref{lem:dfl-detect}(iii), applied with $i=0$ to $\mathfrak h$ and $g_{k-1}$, gives a
cylinder function with loops in $\operatorname{int}C_0^\circ$, vanishing near $X_0$, whose
gradient pairs nontrivially with $\psi$. The pairing is continuous in $(\gamma,\psi)$ and
invariant under translation, so finitely many such functions $g^{(1)},\dots,g^{(m)}$ suffice for
all pairs. The linearization of the family of equations parametrized by
$P_\#=\operatorname{span}\{g^{(j)}\}$ is then onto at every non-regular trajectory of energy $e_k$,
and the Sard--Smale theorem gives arbitrarily small $\rho\in P_\#$ for which all trajectories of
index at most seven and energy $e_k$ of $\mathrm{CS}_{\mathfrak h+g_{k-1}+\rho}$ are regular \cite[p.~80]{DFL},
\cite[pp.~143--144]{Donaldson-book}. For $\rho$ small the trajectories of smaller energy stay regular,
by compactness and openness, and $g_k=g_{k-1}+\rho$ is a split addition to which
Lemma~\ref{lem:app-critical} applies. After the stage $k=N$ the primary data $\mathfrak h+g$,
$g=g_N$, are DFL-regular and have the critical points of $\mathfrak h$. Using
Lemma~\ref{lem:dfl-detect}(iii) with $i=1$ as well gives additions in both pieces. The additions
have loops in $\operatorname{int}C_0^\circ\sqcup\operatorname{int}C_1^\circ$, and the metric data are
never changed.

Fix the required smallness of the total addition before the finite induction, choosing each
increment small enough to meet this requirement. Let $G$ be the family of all cylinder
functions added; this family may depend on the prescribed smallness bound.
Admissibility is preserved stage by stage: if $\mathfrak h+g_{k-1}\in\mathcal U(F\cup G_{k-1})$,
where $G_{k-1}$ is the family added so far, then
$\mathfrak h+g_{k-1}\in\mathcal U(F\cup G_k)$, which is open in $\operatorname{span}(F\cup G_k)$
(Definition~\ref{def:admissible-nbhd}), so the arbitrarily small increment $\rho_k$ may be chosen
with $\mathfrak h+g_k\in\mathcal U(F\cup G_k)$; after the last stage
$\mathfrak h+g\in\mathcal U(F\cup G)$.
\end{proof}

\bigskip
\noindent{\sc Bernd Johannes Wuebben, New York, NY,}
\texttt{wuebben@gmail.com}

\end{document}